\documentclass[a4paper,11pt]{amsart}
\usepackage{amssymb, amsmath, amsthm}
\usepackage{mdwlist}

\def\res{\hbox{ {\vrule height .3cm}{\leaders\hrule\hskip.3cm}}\hskip5.0\mu}

\newcommand\beqn{\begin{equation}}
\newcommand\eeqn{\end{equation}}
\newcommand\beqny{\begin{eqnarray}}
\newcommand\eeqny{\end{eqnarray}}
\newcommand\beqnyn{\begin{eqnarray*}}
\newcommand\eeqnyn{\end{eqnarray*}}

\usepackage{amsmath,amsfonts,amssymb,amsthm}
\newtheorem{theorem}{Theorem}[section]
\newtheorem{lemma}[theorem]{Lemma}
\newtheorem{proposition}[theorem]{Proposition}
\newtheorem{corollary}[theorem]{Corollary}

\def\a{\alpha}           
\def\b{\beta}          
\def\g{\gamma}           
\def\d{\delta}         
\def\e{\epsilon}         
\def\z{\zeta}            
\def\th{\theta}          
\def\k{\kappa}           
        
\def\m{\mu}       
\def\n{\nu}
\def\t{\tau}       
       
\def\p{\pi}       
\def\r{\rho} 

\def\s{\sigma}

\def\E{\mathcal E}
\def\G{\Gamma}    
\def\D{\Delta}  

\def\Th{\Theta}

\begin{document}
\setlength\parskip{5pt}
\title[stable codimension 1  integral varifolds]{A general regularity theory for stable codimension 1  integral varifolds}
\author{Neshan Wickramasekera} 
\thanks{Department of Pure Mathematics and Mathematical Statistics\\
\indent University of Cambridge\\ \indent Cambridge, CB3 0WB, United Kingdom.\\ \indent N.Wickramasekera@dpmms.cam.ac.uk}

%\date{}

\begin{abstract}
We give a necessary and sufficient geometric structural condition, which we call the $\a$-Structural~Hypothesis, for a stable codimension 1 integral varifold on a smooth Riemannian manifold to correspond to an embedded smooth hypersurface away from a small set of generally unavoidable singularities. The $\alpha$-Structural Hypothesis says that no point of the support of the varifold has a neighborhood in which the support is the union of 3 or more embedded $C^{1, \alpha}$ hypersurfaces-with-boundary meeting (only) along their common boundary. We establish that whenever a stable integral $n$-varifold on a smooth $(n+1)$-dimensional Riemannian manifold satisfies the $\alpha$-Structural Hypothesis for some $\alpha \in (0, 1/2)$, its singular set is empty if $n \leq 6$, discrete if $n =7$ and has Hausdorff dimension $\leq n-7$ if $n \geq 8$; in view of well known examples, this is the best possible general dimension estimate on the singular set of a varifold satisfying our hypotheses. We also establish compactness of mass-bounded subsets of the class of stable codimension 1 integral varifolds satisfying the $\a$-Structural Hypothesis for some $\alpha \in (0, 1/2)$.

The $\alpha$-Structural Hypothesis on an $n$-varifold for any $\alpha \in (0, 1/2)$ is readily implied by either of the following two hypotheses: (i) the varifold corresponds to an absolutely area minimizing rectifiable current with no boundary, (ii) the singular set of the varifold has vanishing $(n-1)$-dimensional Hausdorff measure. Thus, our theory subsumes the well known regularity theory for codimension 1 area minimizing rectifiable currents, and settles the long standing question as to which weakest size hypothesis on the singular set of a stable minimal hypersurface guarantees the validity of the above regularity conclusions.

An optimal strong maximum principle for stationary codimension 1 integral varifolds follows from  our regularity and compactness theorems.     
\end{abstract}
%%%%%%%%%%%%%%%%%%%%%%%%%%%%%%%
\maketitle
\newtheorem{hypotheses}{Hypotheses}
\renewcommand{\thehypotheses}{\arabic{section}.\arabic{hypotheses}}

\tableofcontents

\section{Introduction}\label{3-1}
\setcounter{equation}{0}
Here we study regularity properties of stable critical points of the $n$-dimensional area functional in a smooth $(n+1)$-dimensional Riemannian manifold, addressing, among a number of other things, the following basic question:

\textit{When is a stable critical point $V$ of the $n$-dimensional area functional in a smooth $(n+1)$-dimensional Riemannian manifold made-up of pairwise disjoint, smooth, embedded, connected hypersurfaces each of which is itself a critical point of area?} 

Without further hypothesis, $V$ need not satisfy the stated property; this is illustrated by (a sufficiently small region of) any stationary union of three or more hypersurfaces-with-boundary meeting along a common $(n-1)$-dimensional submanifold (e.g.\ a pair of transverse hyperplanes in a Euclidean space). In each of these examples, the connected components of the regular part of the union are not individually critical points of area (in the sense of having vanishing first variation with respect to area for deformations by compactly supported smooth vector fields of the ambient space; see precise definition in Section~\ref{maintheorems} below). 

We give a geometrically optimal answer to the above question by establishing a precise version (given as Corollary 1 below) of the following assertion: 

\textit{Presence of a region of $V$ where three or more hypersurfaces-with-boundary meet along their common boundary is the only obstruction for $V$ to correspond to a locally finite union of pairwise disjoint, smooth, embedded, connected hypersurfaces each of which is itself a critical point of area.} 

This follows directly from our main theorem (the Regularity and Compactness Theorem below) which establishes a precise version of the following regularity statement: 

\textit{Presence of a region of $V$ where three or more hypersurfaces-with-boundary meet along their common boundary is the only obstruction to complete regularity of $V$ in low dimensions, and to regularity of $V$ away from a small, quantifiable, set of generally unavoidable singularities in general dimensions.}

In proving these results, we shall first work in the context where the ambient manifold is an open subset of ${\mathbf R}^{n+1}$ with the Euclidean metric. The differences that arise in the proof in replacing Euclidean ambient space by  a general smooth $(n+1)$-dimensional  Riemannian manifold amount to ``error terms'' in various identities and inequalities that are valid in the Euclidean setting, and can be handled in a straightforward manner. We shall discuss this further in the penultimate section of the paper.
 
Here, a {\em critical point of the $n$-dimensional area} means a \emph{stationary} integral $n$-varifold; i.e.\ an integral $n$-varifold having zero first variation with respect to area under deformation by the flow generated by any compactly supported $C^{1}$ vector field of the ambient space (see hypothesis (${\mathcal S}{\emph 1}$) in  Section~\ref{maintheorems}).

For a varifold $V$, let ${\rm reg} \, V$ denote its regular part, i.e.\  the smoothly embedded part (of the support of the weight measure $\|V\|$ associated with $V$) and 
${\rm sing} \, V$ denote its singular set, i.e.\ the complement of ${\rm reg} \, V$ (in the support of 
$\|V\|$); see Section~\ref{notation} for the precise definitions of these terms. 

A stationary integral varifold  $V$ is {\em stable} if ${\rm reg} \, V$ is stable in the sense that $V$ has non-negative second variation  with respect to area under deformation by the flow generated by any $C^{1}$ ambient vector field which is compactly supported away from ${\rm sing} \, V$ and which, on ${\rm reg} \, V$, is normal to ${\rm reg} \, V$.  In our codimension 1 setting and for Euclidean ambient space,  stability of $V$ whenever ${\rm reg} \, V$ is orientable is equivalent  to requiring that  ${\rm reg} \, V$ satisfies the following  \emph{stability inequality} (\cite{S1}, Section 9):

\begin{equation*}
\int_{{\rm reg} \, V} |A|^{2}\zeta^{2} \,d{\mathcal H}^{n} \leq \int_{{\rm reg} \, V} |\nabla \, \zeta|^{2} \, d{\mathcal H}^{n} \;\;\;\; \forall \;\zeta \in C^{1}_{c}({\rm reg} \, V);
\end{equation*}

\noindent
here $A$ denotes the second fundamental form of ${\rm reg} \, V$, $|A|$ the length of $A$, $\nabla$ the gradient operator on ${\rm reg} \, V$ and ${\mathcal H}^{n}$ is the $n$-dimensional Hausdorff measure on ${\mathbf R}^{n+1}.$ (In fact a slightly weaker form of the stability hypothesis suffices for the proofs of all of  our theorems here, and as a result of that, orientability of ${\rm reg} \,V$ for the varifolds $V$ considered here is a conclusion rather than a hypothesis;  see hypothesis (${\mathcal S}{\emph 2}$) in Section~\ref{maintheorems}, and Corollary~\ref{orientable}.) 

By a \emph{stable integral $n$-varifold} we mean a stable, stationary integral $n$-varifold.

For $\a \in (0, 1)$ and $V$ an integral varifold on a smooth Riemannian manifold, we state the following condition (hypothesis (${\mathcal S}{\emph 3}$) in Section~\ref{maintheorems}) 
which we shall refer to often throughout the rest of the introduction:

%\smallskip

\noindent
{\sc $\a$-Structural Hypothesis}. \emph{No singular point of $V$ has a neighborhood in which  $V$ corresponds to a union of embedded $C^{1, \a}$ hypersurfaces-with-boundary meeting (only) along their common $C^{1, \a}$ boundary (and with multiplicity a constant positive integer on each of the constituent hypersurfaces-with-boundary).}

%\smallskip

Our main theorem (Theorem~\ref{compactness-N}; for Euclidean ambient space, Theorem~\ref{compactness}) can now be stated as follows:

%\smallskip

\noindent
{\sc Regularity and Compactness Theorem}. \emph{A stable integral $n$-varifold $V$ on a smooth $(n+1)$-dimensional Riemannian manifold corresponds to an embedded hypersurface with no singularities when $1 \leq n \leq 6;$ to one with at most a discrete set of singularities when $n=7;$ and to one with a closed set of singularities having Hausdorff dimension at most $n-7$ when $n \geq 8$, (and with multiplicity, in each case, a constant positive integer on each connected component of the hypersurface), provided $V$ satisfies the $\a$-Structural Hypothesis above for some $\a \in (0, 1/2).$}

\emph{Furthermore, for any given $\a \in (0, 1/2)$, each mass-bounded subset of the class of stable codimension 1 integral varifolds satisfying the $\a$-Structural Hypothesis is compact in the topology of varifold convergence.}

In case $V$ corresponds to an absolutely area minimizing codimension 1 rectifiable current, the regularity conclusion of this theorem is well known, and is the result of
combined work of E.~De~Giorgi (\cite{DG}), R.~Reifenberg (\cite{RR}), W.~Fleming (\cite{FW}), F.~Almgren (\cite{A1}), J.~Simons (\cite{SJ}) and H.~Federer (\cite{F1}). While our work uses  ideas and results from some of these pioneering work,  it does not rely upon the fact that the conclusions hold in the area minimizing case;  it is interesting to note that the above theorem indeed subsumes the regularity theory for codimension 1 area minimizing rectifiable currents for the following simple reason: If  
$T$ is a rectifiable current on an open ball and if $T$ has no boundary in the interior of the ball and is supported on a union of 3 or more embedded hypersurfaces-with-boundary meeting only along their common boundary,  then $T$ cannot be area minimizing. 

 Let $V$ be a stationary integral $n$-varifold on a Riemannian manifold $N$. Once we know that the singular set of $V$ is sufficiently small---in fact, as small as having vanishing $(n-1)$-dimensional  Hausdorff measure---it is not difficult to check that the multiplicity 1 varifold associated with each connected component of the regular part of $V$ is itself stationary in $N$. Thus we deduce from the Regularity and Compactness Theorem the following:

%\smallskip

\noindent
{\sc Corollary 1}. \emph{The $\a$-Structural Hypothesis (see above)  for some $\a \in (0, 1/2)$ is necessary and sufficient
 for a stable codimension 1 integral varifold $V$ on a smooth Riemannian manifold $N$ to have the following ``local decomposability property''}: 

\emph{For each open $\Omega \subset N$ with compact closure in $N$, there exist a finite number of pairwise disjoint, smooth, embedded, connected hypersurfaces $M_{1}, M_{2}, \ldots, M_{k}$ of $\Omega$ (possibly with 
a non-empty interior singular set ${\rm sing} \, M_{j} = (\overline{M}_{j} \setminus M_{j}) \cap \Omega$ for each $j=1, 2, \ldots, k$) and positive integers $q_{1}, q_{2}, \ldots, q_{k}$ 
 such that the multiplicity 1 varifold $|M_{j}|$ defined by $M_{j}$ is stationary in $\Omega$ for each $j=1, 2, \ldots, k$ and $V \res \Omega = \sum_{j=1}^{k} q_{j}|M_{j}|.$} 

%\smallskip

In 1981, R. Schoen and L. Simon (\cite{SS}) proved that the conclusions of the Regularity and Compactness Theorem hold for the $n$-dimensional stable minimal hypersurfaces (viz.\ embedded hypersurfaces which are stationary and stable as multiplicity 1 varifolds) satisfying, in place of the $\a$-Structural Hypothesis, the (much more restrictive) property that the singular sets have locally finite $(n-2)$-dimensional Hausdorff measure. Since then, it has remained an open question as to what the \emph{weakest} size hypothesis (in terms of Hausdorff measure) on the singular sets is that would guarantee the validity of the same conclusions. Since vanishing of the $(n-1)$-dimensional Hausdorff measure of the singular set trivially implies the $\a$-Structural Hypothesis, we have the following immediate corollary of the Regularity and Compactness Theorem, which settles this question:

%\smallskip

\noindent
{\sc Corollary 2.} \emph{The conclusions of the Regularity and Compactness Theorem hold for the $n$-dimensional stable minimal hypersurfaces with singular sets of vanishing $(n-1)$-dimensional Hausdorff measure. In fact, a stable codimension 1 integral $n$-varifold $V$ satisfies the $\a$-Structural Hypothesis for some $\a \in (0, 1/2)$ if and only if its singular set has vanishing $(n-1)$-dimensional Hausdorff measure.} 

%\smallskip

A union of two transversely intersecting hyperplanes in a Euclidean space shows that for no $\g > 0$ can the singular set hypothesis in Corollary 2 be weakened to vanishing of the $(n-1+\g)$-dimensional Hausdorff measure. 

In contrast to our $\a$-Structural Hypothesis, the singular set hypothesis of \cite{SS} (i.e.\ the hypothesis that ${\mathcal H}^{n-2} \, ({\rm sing} \, V \cap K) < \infty$ for each compact subset $K$ of the ambient space), together with stability away from the singular set, \emph{a priori} implies, by a straightforward argument, that the singularities are ``removable for the stability inequality''--- that is to say, the above stability inequality is valid for the larger class of test functions $\z$ which are the restrictions to the hypersurface of compactly supported smooth functions of the ambient space (that are not required to vanish near the singular set). The techniques employed in \cite{SS} in the proof of the regularity theorems therein relied on this fact in an essential way. Interestingly, the $\a$-Structural Hypothesis, or, for that matter, vanishing of the $(n-1)$-dimensional Hausdorff measure of the singular set, does not seem to imply \emph{a priori} even local finiteness of total curvature, viz.\ $\int_{{\rm reg} \, V \cap K} |A|^{2} < \infty$ for each compact subset $K$ of the ambient space (whereas the singular set hypothesis of \cite{SS} does, in view of the strengthening of the stability inequality just mentioned). This means that we cannot in our proof use the stability inequality in a direct way over arbitrary regions of the varifolds. (Of course \emph{a posteriori} we can strengthen the stability inequality in the manner described above, so in particular it is true under our hypotheses that
$\int_{{\rm reg} \, V \cap K} |A|^{2} < \infty$ for each compact subset $K$ of the ambient space.) Our proof nevertheless at several stages makes indispensable use of the work of Schoen and Simon---specifically, Theorem~\ref{SS} below; indeed, application of Theorem~\ref{SS} in regions where we have sufficient control over the singular set is a principal way in which the stability hypothesis enters our proof. 
 
The Regularity and Compactness Theorem  is optimal in several ways. A key aspect of the theorem is that it \emph{requires no hypothesis concerning the size of the singular sets}; nor does it require any hypothesis concerning the generally-difficult-to-control set of points where some tangent cone is a plane of multiplicity 2 or higher. What suffices is the  $\a$-Structural Hypothesis, which is easier to check in principle. As mentioned before, stationary unions of half-hyperplanes of a Euclidean space meeting along common axes illustrate that the $\a$-Structural Hypothesis is a sharp condition for the regularity conclusions of the theorem.

In view of well known examples of 7-dimensional stable hypercones with isolated singularities (e.g. the cone over ${\mathbf S}^{3}(1/\sqrt{2}) \times {\mathbf S}^{3}(1/\sqrt{2}) \subset {\mathbf R}^{8}$), the Regularity and Compactness Theorem is also optimal with regard to its conclusions in the sense that it gives, in dimensions $\geq 7,$ the optimal general estimate on the Hausdorff dimension of the singular sets.

It remains an open question as to what one can say about the size of the singular sets if the stability hypothesis in the theorem is removed. Obviously one cannot in this case draw the same conclusions in view of the fact that there are embedded non-equatorial minimal surfaces of ${\mathbf S}^{3}$ (e.g.\ $\mathbf{S}^{1}(1/\sqrt{2}) \times \mathbf{S}^{1}(1/\sqrt{2}) \subset {\mathbf S}^{3}$), the cones over which provide examples of stationary (unstable) hypercones in ${\mathbf R}^{4}$ with isolated singularities. There are no 2-dimensional singular stationary hypercones satisfying the $\a$-Structural Hypothesis; it is however not known whether there is a singular 2-dimensional stationary integral varifold $V$ in ${\mathbf R}^{3}$ such that $V$ either satisfies the $\a$-Structural Hypothesis or has a singular set of vanishing 1-dimensional Hausdorff measure or has an isolated singularity. It also remains largely open what one can say concerning stable integral varifolds of codimension $>1$. Again, the same conclusions as in our theorem cannot be made in this case due to the presence of branch point singularities, as illustrated by 2-dimensional holomorphic varieties with isolated branch points. See the remark following the statement of the Regularity and Compactness Theorem in Section 3 (Theorem~\ref{compactness}) for a further discussion on optimality of our results here.

For a general stationary integral varifold, a point where some tangent cone is a plane of multiplicity 2 or higher may or may not be a regular point. Our ``Sheeting Theorem'' (Theorem~\ref{main} below) implies that if the varifold satisfies the hypotheses of the Regularity and Compactness Theorem, then such a point is a regular point. (As is well known, a point where there is a multiplicity 1 tangent plane is always a regular point, for any stationary integral varifold, by the regularity theorem of W.~K.~Allard ([\cite{AW}, Section 8]; also [\cite{S1}, Theorem 23.1])). Indeed, the Sheeting Theorem is one of the two principal ingredients of the proof of the Regularity and Compactness Theorem; the other is the ``Minimum Distance Theorem'' (Theorem~\ref{no-transverse}) which implies that no tangent cone to a varifold satisfying the hypotheses of the Regularity and Compactness Theorem can be supported by a union of three or more half-hyperplanes meeting along a common $(n-1)$-dimensional axis.

A direct consequence of Allard's regularity theorem is that the regular part of a stationary integral varifold is a non-empty---in fact a dense---subset of its support [\cite{AW}, Section 8.1]. Thus, given stationarity of the varifold,  our stability hypothesis, which concerns only the regular part of the varifold, is never vacuously true. However, an open, dense subset could have arbitrarily small (positive) measure, and in fact, as mentioned above, under the stationarity hypothesis alone no general result whatsoever is known concerning the Hausdorff measure of the singular sets. Closely related to this is the point made before that from the hypotheses of the Regularity and Compactness Theorem, not even local finiteness of total curvature seems to follow \emph{a priori}. In light of these considerations which indicate that our hypotheses are rather mild, it is somewhat surprising that they imply optimal regularity of the hypersurfaces.

We may summarise as follows all of the various regularity results discussed above and established in subsequent sections of the paper:

%\smallskip
\noindent
{\sc Theorem.} \emph{Let $V$ be a stable integral $n$-varifold on a smooth $(n+1)$-dimensional Riemannian manifold $N$. The following statements concerning $V$ are equivalent:
\begin{itemize}
\vspace{-.15in}
\item[(a)] For some $\a \in (0, 1/2)$, $V$ satisfies the $\a$-Structural Hypothesis, viz. no singular point of $V$ has a neighborhood in which  $V$ corresponds to  a union of $C^{1, \a}$ embedded hypersurfaces-with-boundary meeting (only) along their common boundary, with multiplicity a constant positive integer on each constituent hypersurface-with-boundary.
\item[(b)] ${\rm sing} \, V = \emptyset$ if $1 \leq n \leq 6$, ${\rm sing} \, V$ is discrete if $n=7$ and 
${\mathcal H}^{n-7+\g} \, ({\rm sing} \, V) = 0$ for each $\g >0$ if $n \geq 8.$
\item[(c)] ${\mathcal H}^{n-1} \, ({\rm sing} \, V) = 0.$
\item[(d)] $V$ has the local decomposability property (defined in Corollary 1 above), viz. for each open $\Omega \subset N$ with compact closure in $N$, there exist a finite number of pairwise disjoint, smooth, embedded, connected hypersurfaces $M_{1}, M_{2}, \ldots, M_{k}$ of $\Omega$ (possibly with $(\overline{M}_{j} \setminus M_{j}) \cap \Omega$ non-empty for each $j=1, 2, \ldots, k$) and positive integers $q_{1}, q_{2}, \ldots, q_{k}$ 
 such that the multiplicity 1 varifold $|M_{j}|$ defined by $M_{j}$ is stationary in $\Omega$ for each $j=1, 2, \ldots, k$ and $V \res \Omega = \sum_{j=1}^{k} q_{j}|M_{j}|.$
\item[(e)] No tangent cone of $V$ corresponds to a union of three or more half-hyperplanes meeting along a common $(n-1)$-dimensional subspace, with multiplicity a constant positive integer on each constituent half-hyperplane.
\item[(f)] $V$ satisfies the $\a$-Structural Hypothesis for each $\a \in (0, 1/2).$
\end{itemize}}

%\smallskip

Finally, we mention another direct implication of the Regularity and Compactness Theorem, namely, the following optimal strong maximum principle (Theorem~\ref{maxm}) for codimension 1 stationary integral varifolds:

\smallskip

\noindent
{\sc Varifold Maximum principle}. \emph{Let $N$ be a smooth $(n+1)$-dimensional Riemannian manifold and 
%\begin{itemize}
%\vspace{-.15in}
%\item[(a)] If $V_{1}$, $V_{2}$ are stationary integral $n$-varifolds on $N$ such that ${\mathcal H}^{n-1} \, ({\rm spt} \, \|V_{1}\| \cap {\rm spt} \, \|V_{2}\|) = 0,$ then ${\rm spt} \, \|V_{1}\| \cap {\rm spt} \, \|V_{2}\| = \emptyset.$
%\item[(b)] 
let $\Omega_{1}$, $\Omega_{2}$ be open subsets of $N$ such that $\Omega_{1} \subset \Omega_{2}$. Let
$M_{i} = \partial \, \Omega_{i}$ for $i=1, 2.$ If for $i=1, 2$, $M_{i}$ is connected,  ${\mathcal H}^{n-1} \, ({\rm sing} \, M_{i} ) = 0$ and $V_{i} \equiv |M_{i}|$ is stationary in $N,$ then either ${\rm spt} \, \|V_{1}\|  = {\rm spt} \, \|V_{2}\| $ or ${\rm spt} \, \|V_{1}\| \cap {\rm spt} \, \|V_{2}\|  = \emptyset.$ 
Here ${\rm sing} \, M_{i} = M_{i} \setminus {\rm reg} \, M_{i}$ where ${\rm reg} \, M_{i}$ is the set of points $X \in M_{i}$ such that $M_{i}$ is a smooth, embedded submanifold near $X$.} 
%\end{itemize}}

%\smallskip

\noindent
See Section~\ref{notation} for explanation of notation used here. If the varifolds $V_{1}$ and $V_{2}$ are both free of singularities, the theorem is easily seen to follow from the Hopf maximum principle. B.~Solomon and B.~White (\cite{SoW}) proved the theorem assuming only that one of $V_{1}$ or $V_{2}$ is free of singularities (allowing the other to be arbitrary with no restriction on its singular set). L.~Simon (\cite{S4}) and independently M.~Moschen (\cite{Mo}) established the result in case $V_{1}$ and $V_{2}$ correspond to area minimizing  integral currents, both possibly singular. Using  the Schoen--Simon regularity theory (\cite{SS}), some key ideas from \cite{S4} as well as 
the Solomon--White theorem, T.~Ilmanen (\cite{I}) established the theorem (for stationary $V_{1}$, $V_{2}$) subject to the stronger condition ${\mathcal H}^{n-2} \, ({\rm sing} \, M_{i}) < \infty$ for $i=1, 2$. The version above follows directly from the argument in \cite{I}, in view of the fact that we may use Corollary 2 above in places where the argument in \cite{I} depended on the Schoen-Simon theory. This version is optimal in the sense that larger singular sets cannot generally be allowed.

\noindent
{\bf Outline of the method:} Below we give a brief description of the proof of the Regularity and Compactness Theorem.

Fix any $\a \in (0, 1)$ and let ${\mathcal S}_{\a}$ denote the family of stable integral $n$-varifolds of the open ball $B_{2}^{n+1}(0) \subset {\mathbf R}^{n+1}$ satisfying the $\a$-Structural Hypothesis. The proof of the Regularity and Compactness Theorem  is based on establishing the fact that no tangent cone at a singular point of a varifold belonging to the varifold closure of ${\mathcal S}_{\a}$ can be  supported by (a) a hyperplane or (b) a union of half-hyperplanes meeting along an $(n-1)$-dimensional subspace. Once this is established, it is not difficult to reach the conclusions of the theorem with standard arguments. 

The assertion in case  (a) is implied by the following regularity result (Theorem~\ref{main}):

\noindent
{\sc Sheeting Theorem}. \emph{Whenever a varifold in ${\mathcal S}_{\a}$ is weakly close to 
a given hyperplane ${\mathbf P}_{0}$ of constant positive integer multiplicity, it must break up in the interior into disjoint, embedded smooth graphs (``sheets'') of small curvature over ${\mathbf P}_{0}$}.

The assertion in case (b) is a consequence of the following (Theorem~\ref{no-transverse}):

\noindent
{\sc Minimum Distance 
Theorem}. \emph{No varifold in ${\mathcal S}_{\a}$ can be weakly close to a given 
stationary integral hypercone ${\mathbf C}_{0}$ corresponding to a union of three or more half-hyperplanes meeting along an $(n-1)$-dimensional subspace (and with constant positive integer multiplicity on each  half-hyperplane).} 

Our strategy is to prove both the Sheeting Theorem and the Minimum Distance Theorem simultaneously by an inductive argument, inducting on the multiplicity $q$ of ${\mathbf P}_{0}$ for the Sheeting Theorem and on the density $\Theta \, (\|{\mathbf C}_{0}\|, 0)$ (= $q$ or $q + 1/2$) of  ${\mathbf C}_{0}$ at the origin for the Minimum Distance Theorem, where $q$ is an integer $\geq 1.$ Approaching both theorems inductively and \emph{simultaneously} in this manner makes it possible to establish, for varifolds in ${\mathcal S}_{\a}$ (satisfying appropriate ``small excess'' hypotheses in accordance with the theorems) and for their ``blow-ups,'' many of the necessary a priori estimates which seem inaccessible via an approach (inductive or otherwise) aimed at proving the two theorems separately. 

The main general idea in the argument is the following: Let $q$ be an integer $\geq 2$ and assume by induction the validity of the Sheeting Theorem when ${\mathbf P}_{0}$ has multiplicity $\in \{1,  \ldots, q-1\}$ and of the Minimum Distance Theorem when $\Theta \, (\|{\mathbf C}_{0}\|, 0) \in \{3/2, \ldots, q-1/2, q\}.$ Then, in a region of a varifold in ${\mathcal S}_{\a}$ where no singular point has density $\geq q$, we may apply the induction hypotheses together with a theorem of J. Simons (\cite{SJ}; see also \cite{S1}, Appendix B) and the ``generalized stratification of stationary integral varifolds'' due to F. J. Almgren Jr. (\cite{A}, Theorem 2.26 and Remark 2.28) to reduce the dimension of the singular set to a low value. This permits effective usage of the stability hypothesis, including applicability of the Schoen-Simon version (\cite{SS}, Theorem 2; also Theorem~\ref{SS} below) of the Sheeting Theorem, in such a region. On the other hand, in the presence of singularities of density $\geq q$ (and whenever the density ratio of the varifold at scale 1 is close to  $q$), it is possible to make good use of the monotonicity formula; most notable among its consequences in the present context are versions (Theorem~\ref{L2-est-1} and Corollary~\ref{L2-est-2}), for a varifold in ${\mathcal S}_{\a}$ with small ``height excess'' relative to a hyperplane and lower order height excess relative to certain cones, of L.~Simon's (\cite{S}) a priori $L^{2}$-estimates, and an analogous, new, ``non-concentration-of-tilt-excess'' estimate (Theorem~\ref{non-concentration}(b)) giving control of the amount of its ``tilt-excess'' relative to the hyperplane in regions where there is a high concentration of points of density $\geq q$. 

Combining these techniques, we are able to fully analyse, under the induction hypotheses, the  ``coarse blow-ups,''  namely, the compact class ${\mathcal B}_{q} \subset W^{1, 2}_{\rm loc} \, (B_{1}; {\mathbf R}^{q}) \cap L^{2} \ (B_{1}; {\mathbf R}^{q})$ ($B_{1} =$ the open unit ball in ${\mathbf R}^{n}$) consisting of ordered $q$-tuples of functions produced by blowing up sequences of varifolds in ${\mathcal S}_{\a}$ converging weakly to a multiplicity $q$ hyperplane. (See precise definition of ${\mathcal B}_{q}$ at the end of Section~\ref{blow-up}.) One of the key properties that need to be established for ${\mathcal B}_{q}$ is that it does not contain an element $H$ 
whose graph is the union of $q$ half-hyperplanes in one half-space of ${\mathbf R}^{n+1}$ and $q$ half-hyperplanes in the complementary half-space, with all half-hyperplanes meeting along a common $(n-1)$-dimensional subspace and at least two of them distinct on one side or the other. (This is a Minimum Distance Theorem for ${\mathcal B}_{q}$, analogous to the Minimum Distance Theorem for ${\mathcal S}_{\a}$.) Establishing this property takes considerable effort and occupies a significant part (Sections~\ref{step3} through ~\ref{propertiesIII}) of our work, and is achieved as follows:

First we rule out (in Section~\ref{step3}), by a first variation argument utilizing the non-concentration-of-tilt-excess estimate of Theorem~\ref{non-concentration}(b), the possibility that there is such $H \in {\mathcal B}_{q}$ with its graph having all $q$ half-hyperplanes on one side coinciding (but not on the other). 

The second, more involved step is to rule out the existence of such an element in ${\mathcal B}_{q}$ (call it $H^{\prime}$) with its graph having at least two distinct half-hyperplanes on each side. To this end we assume, arguing by contradiction, that there is such $H^{\prime} \in {\mathcal B}_{q}$ and use the induction hypotheses to implement   
a ``fine blow-up''  procedure (see definition at the end of Section~\ref{fineexcessblowup}), where certain sequences of varifolds in ${\mathcal S}_{\a}$ are blown-up by their height excess (the ``fine excess'') relative to appropriate unions of half-hyperplanes (corresponding to ``vertical'' scalings of $H^{\prime}$ by the coarse excess of the varifolds giving rise to $H^{\prime}$). We use first variation arguments (in particular, Simon's $L^{2}$-estimates and the non-concentration-of-tilt-excess estimate of Theorem~\ref{non-concentration}(b)) and standard $C^{1,\b}$ boundary regularity theory for harmonic functions to prove a uniform interior continuity estimate (Theorem~\ref{c1alpha}) for the first derivatives of the fine blow-ups, and use it, via an excess improvement argument, to show that our assumption $H^{\prime} \in {\mathcal B}_{q}$ must contradict one of the induction hypotheses, namely, that the Minimum Distance Theorem is valid when $\Theta \, (\|{\mathbf C}_{0}\|, 0) = q$. This enables us to conclude that the coarse blow-up class ${\mathcal B}_{q}$ has the asserted property, viz. that the only elements in ${\mathcal B}_{q}$ which are given by linear functions on either of two complementary half-spaces are the ones given by $q$ copies of a single linear function everywhere.

Equipped with this fact and a number of other key properties that we establish for the coarse blow-ups (see items (${\mathcal B}{\emph 1}$)-(${\mathcal B}{\emph 7}$) of Section~\ref{proper-blow-up} for a complete list), we ultimately obtain (in Theorems~\ref{blowups} and ~\ref{blowup-reg}), subject to the induction hypotheses,  interior $C^{1}$ regularity of coarse blow-ups and consequently, that any coarse blow-up is an ordered set of $q$ harmonic functions (a Sheeting Theorem for ${\mathcal B}_{q}$, analogous to the Sheeting Theorem for ${\mathcal S}_{\a}$); furthermore, we show that these harmonic functions all agree if  infinitely many members of a sequence of varifolds giving rise to the blow-up contain, in the interior, points of density $\geq q$.

The preceding result is the key to completion of the induction step for the Sheeting Theorem. Together with the Schoen-Simon version of the Sheeting Theorem, it enables us to prove a De-Giorgi type lemma (Lemma~\ref{excess-s}), the iterative application of which leads us to the following conclusion: Let ${\mathbf P}_{0}$ be a hyperplane with multiplicity $q,$ and suppose that $V$ is a  varifold in ${\mathcal S}_{\a}$ lying weakly close to ${\mathbf P}_{0}$ in a unit cylinder over ${\mathbf P}_{0}$. Let $D$ be the region of ${\mathbf P}_{0}$ inside a cylinder slightly smaller than the unit cylinder. Then (i) there is a closed subset of $D$ over each point of which the support of $V$ consists of a single point; furthermore, at this point, $V$ has a unique multiplicity $q$ tangent hyperplane almost parallel to ${\mathbf P}_{0},$ and relative to this tangent hyperplane, the height excess of $V$ satisfies a uniform decay estimate; and  (ii) over the complementary open set, $V$ corresponds to embedded  graphs of $q$ ordered, analytic functions of small gradient solving the minimal surface equation. Facts (i), (ii) and elliptic estimates imply, by an elementary general argument (Lemma~\ref{general}), that the varifold corresponds to $q$ ordered graphs over all of $D$ and that each graph satisfies a uniform $C^{1, \b}$ estimate (Theorem~\ref{sheetingthm}) for some fixed $\b \in (0, 1),$ completing the induction step for the Sheeting Theorem.

The final step of the argument is to complete induction for the Minimum Distance Theorem, which requires showing that the Minimum Distance Theorem holds whenever  
$\Theta \, (\|{\mathbf C}_{0}\|, 0) \in \{q + 1/2, q+1\},$ where ${\mathbf C}_{0}$ is a stationary cone as in the theorem. Since we may now assume the validity of the Sheeting Theorem for multiplicity up to and including $q,$ we have all the necessary ingredients to establish (in Theorem~\ref{no-transverse-q}) that given such ${\mathbf C}_{0}$, if there is a varifold  $V \in {\mathcal S}_{\a}$ weakly close to ${\mathbf C}_{0},$ then it must in the interior be made up of $C^{1, \a}$ embedded hypersurfaces-with-boundary meeting along their common boundary; this directly contradicts the $\a$-Structural Hypothesis and proves the Minimum Distance Theorem, subject to the induction hypotheses, when $\Theta \, (\|{\mathbf C}_{0}\|, 0) \in \{q + 1/2, q+1\}.$ Our argument also establishes the Minimum Distance Theorem when $\Theta \, (\|{\mathbf C}_{0}\|, 0) \in \{3/2, 2\},$ since in this case we have, in place of the induction hypotheses, Allard's Regularity Theorem which implies the Sheeting Theorem when $q=1$. 

This completes the outline of the proof.

%\medskip

\noindent
{\bf Acknowledgements:} 
Several technical aspects of this work owe a large mathematical debt to the ideas of Leon~Simon contained in his seminal work \cite{S} on asymptotics for minimal submanifolds near singularities. A crucial first step needed for the theory developed here is an a priori estimate, due to Rick~Schoen and Leon~Simon (\cite{SS}), for sufficiently regular stable minimal hypersurfaces.  Other key 
previously known results used here include: an integral estimate  (the Hardt--Simon inequality in Section~\ref{proper-blow-up}) for (coarse) blow-ups of stationary varifolds, first established in the work of Bob~Hardt and Leon~Simon (\cite{HS}); Almgren's multiple-valued approximate graph decomposition theorem (\cite{A}); Almgren's generalised stratification principle for stationary integral varifolds (\cite{A}) and J.~Simons' theorem concerning non-existence of low dimensional singular stable minimal hypercones (\cite{SJ}). The basic ``harmonic approximation''  idea, used at several places in this work in various forms and levels of sophistication, has its origin in the work of  De~Giorgi (\cite{DG}), and is found in a form closer to the way it is used here in 
the work of Bill~Allard (\cite{AW}).

I am very grateful to the two referees for taking the time to read the manuscript so carefully and pointing out various corrections---one report in particular was extremely detailed and contained also many valuable suggestions for clarifying some arguments and improving  the exposition. I am also grateful to Yoshi Tonegawa for valuable comments on a first draft---in particular for not being entirely happy (!) with the (then) level of details in Section 10 of the paper. I  thank Brian~White for a conversation concerning orientability of the regular part of a stationary codimension 1 varifold on a ball,  which lead to a slightly weaker and more natural form of the stability hypothesis than in an earlier version of the paper, and to Corollary~\ref{orientable}. At an early stage of this work, details of some preliminary ideas were worked out during the 2008 Calculus of Variations workshop at Oberwolfach, and I thank MFO for its hospitality.

\section{Notation}\label{notation}
The following notation will be used throughout the paper:

$n$ is a fixed positive integer $\geq 2,$  ${\mathbf R}^{n+1}$ denotes the $(n+1)$-dimensional Euclidean space and 
$(x^{1}, x^{2}, y^{1}, y^{2}, \ldots, y^{n-1}),$ which we shall sometimes abbreviate  as 
$(x^{1}, x^{2}, y),$ denotes a general point in ${\mathbf R}^{n+1}$. We shall identify 
${\mathbf R}^{n}$ with the hyperplane $\{x^{1} = 0\}$ of ${\mathbf R}^{n+1},$ and ${\mathbf R}^{n-1}$ with the subspace 
$\{x^{1} = x^{2} = 0\}.$

For $Y \in {\mathbf R}^{n+1}$ and $\r >0$, $B_{\r}^{n+1}(Y) = \{X \in {\mathbf R}^{n+1} \, : \, |X - Y| < \r\}.$

For $Y \in {\mathbf R}^{n}$ and $\r >0$, $B_{\r}(Y) = \{X \in {\mathbf R}^{n} \, : \, |X -Y| < \r\}.$ We shall often abbreviate $B_{\r}(0)$ as $B_{\r}.$

For $Y \in {\mathbf R}^{n+1}$ and $\r>0$, $\eta_{Y, \r} \, : \, {\mathbf R}^{n+1} \to {\mathbf R}^{n+1}$ is the map defined by 
$\eta_{Y, \r}(X) = \r^{-1}(X - Y)$ and $\eta_{\r}$ abbreviates $\eta_{0, \r}$.

${\mathcal H}^{k}$ denotes the $k$-dimensional Hausdorff measure in ${\mathbf R}^{n+1}$, and 
$\omega_{n} = {\mathcal H}^{n} \, (B_{1}(0)).$

For $A, B \subset {\mathbf R}^{n+1}$, ${\rm dist}_{\mathcal H} \, (A, B)$ denotes the Hausdorff distance between $A$ and $B$.  

For $X \in {\mathbf R}^{n+1}$ and $A \subset {\mathbf R}^{n+1}$, ${\rm dist} \, (X, A) = \inf_{Y \in A} \, |X - Y|.$

For $A \subset {\mathbf R}^{n+1}$, $\overline{A}$ denotes the closure of $A$.

$G_{n}$ denotes the space of hyperplanes of ${\mathbf R}^{n+1}.$

 For an $n$-varifold $V$ (\cite{AW}; see also [\cite{S1}, Chapter 8]) on an open subset $\Omega$ of ${\mathbf R}^{n+1}$, an open subset $\widetilde{\Omega}$ of $\Omega$, a Lipschitz mapping $f \, : \, \Omega \to {\mathbf R}^{n+1}$ and a countably $n$-rectifiable subset $M$ of $\Omega$ with locally finite ${\mathcal H}^{n}$-measure, we use the following notation:

$V \, \res \,\widetilde{\Omega}$ abbreviates the restriction $V \, \res \, (\widetilde{\Omega} \times G_{n})$ of $V$ to $\widetilde{\Omega} \times G_{n}.$

$\|V\|$ denotes the weight measure on $\Omega$ associated with $V$. 

${\rm spt} \, \|V\|$ denotes the support of $\|V\|.$ 

$f_{\#} \, V$ denotes the image varifold under the mapping $f.$ 

$|M|$ denotes the multiplicity 1 varifold on $\Omega$ associated with $M$.

For $Z \in {\rm spt} \, \|V\| \cap \Omega$, ${\rm VarTan} \, (V, Z)$ denotes the set of tangent cones to $V$ at $Z$.

 ${\rm reg} \, V$ denotes the (interior) regular part of 
${\rm spt} \, \|V\|.$ Thus, $X \in {\rm reg} \,V$ if and only if $X \in {\rm spt} \, \|V\| \cap \Omega$ and there exists $\r >0$ such that $\overline{B^{n+1}_{\r}(X)} \cap {\rm spt} \, \|V\|$ is a smooth, compact, connected, embedded hypersurface-with-boundary, with its boundary contained in $\partial \, B^{n+1}_{\r}(X).$

${\rm sing} \, V$ denotes the interior singular set of ${\rm spt} \, \|V\|.$ Thus, 
${\rm sing} \, V = ({\rm spt} \, \|V\| \setminus {\rm reg} \, V) \cap \Omega.$

%\medskip

\section{Statement of the main theorems}\label{maintheorems}

\noindent
{\bf The Class ${\mathcal S}_{\a}$:} Fix any $\a \in (0, 1).$ Denote by ${\mathcal S}_{\a}$ the collection of all  integral 
$n$-varifolds $V$ on  $B_{2}^{n+1}(0)$  with 
$0 \in {\rm spt} \, \|V\|$, $\|V\|(B_{2}^{n+1}(0)) < \infty$ and satisfying the following conditions:
\begin{itemize}
\item[(${\mathcal S{\emph 1}}$)] {\sc Stationarity:} $V$ has zero first variation with respect to the area functional in the following sense: For any given vector field 
$\psi \in C^{1}_{c}(B_{2}^{n+1}(0); {\mathbf R}^{n+1}),$ $\e > 0$ and  $C^{2}$ map 
$\varphi \, : \, (-\epsilon, \epsilon) \times B_{2}^{n+1}(0) \to B_{2}^{n+1}(0)$ such that 
\begin{itemize} 
\item[(i)] $\varphi(t, \cdot) \, : \, B_{2}^{n+1}(0) \to B_{2}^{n+1}(0)$ is a $C^{2}$ diffeomorphism for each $t \in (-\epsilon, \epsilon)$ with $\varphi(0, \cdot)$ equal to the identity map on $B_{2}^{n+1}(0),$
\item[(ii)] $\varphi(t, x) = x$ for each $(t, x)   \in (-\epsilon, \epsilon) \times \left(B_{2}^{n+1}(0) \setminus {\rm spt} \, \psi\right)$ and
\item[(iii)] $\left.\partial \, \varphi(t, \cdot)/\partial \, t\right|_{t=0} = \psi$ 
\end{itemize}
(the flow generated by $\psi$ for instance gives rise to such a family $\varphi(t, \cdot)$), we have that 
$$\left.\frac{d}{dt}\right|_{t=0} \, \|\varphi(t, \cdot)_{\#} \, V\|(B_{2}^{n+1}(0)) = 0;$$  
equivalently (see [\cite{S1}, Section 39]),
\begin{equation}\label{firstvariation}
\int_{B_{2}^{n+1}(0) \times G_{n}} {\rm div}_{S} \,\psi(X) \, dV(X, S) = 0
\end{equation}
for every vector field 
$\psi \in C^{1}_{c}(B_{2}^{n+1}(0); {\mathbf R}^{n+1}).$\\
\item[(${\mathcal S{\emph 2}}$)] {\sc Stability:} For each open set $\Omega \subset B_{2}^{n+1}(0)$ such that ${\rm sing} \, V \cap \Omega = \emptyset$ in case $2 \leq n \leq 6$ or 
${\mathcal H}^{n-7+\g}({\rm sing} \, V \cap \Omega) = 0$ for every $\g > 0$ in case $ n\geq 7$, we have that
\begin{equation}\label{stability}
\int_{{\rm reg} \, V \cap \Omega} |A|^{2}\zeta^{2} \,d{\mathcal H}^{n} \leq \int_{{\rm reg} \, V \cap \Omega} |\nabla \, \zeta|^{2} \, d{\mathcal H}^{n} \;\;\;\; \forall \; \zeta \in C^{1}_{c}({\rm reg} \, V \cap \Omega),
\end{equation}
\noindent
where $A$ denotes the second fundamental form of ${\rm reg} \, V$, $|A|$ the length of $A$ and $\nabla$ denotes the gradient operator on ${\rm reg} \, V;$  equivalently (see \cite{S1}, Section 9),  for each such $\Omega,$ 
$V$ has non-negative second variation with respect to area for normal deformations compactly supported in $\Omega \setminus {\rm sing} \, V$, in the following sense:  for any given vector field 
$\psi \in C^{1}_{c}(\Omega \setminus {\rm sing} \, V; {\mathbf R}^{n+1})$ with 
$\psi(X) \perp T_{X}{\rm reg} \, V$ for each $X \in {\rm reg} \, V \cap \Omega$, 
$$\left.\frac{d^{2}}{dt^{2}}\right|_{t=0} \, \|\varphi(t, \cdot)_{\#} \, V\|(B_{2}^{n+1}(0)) \geq 0$$
where $\varphi(t, \cdot),$ $t \in (-\epsilon, \epsilon)$, are the $C^{2}$ diffeomorphisms of $B_{2}^{n+1}(0)$ associated with $\psi,$ described in (${\mathcal S{\emph 1}}$) above.\\ 
\item[(${\mathcal S{\emph 3}}$)] {\sc $\a$-Structural Hypothesis:} For each given $Z \in {\rm sing} \, V,$ there exists {\em no} $\r >0$ such that ${\rm spt} \, \|V\| \cap B_{\r}^{n+1}(Z)$ is equal to 
the union of a finite number of embedded $C^{1,\a}$ hypersurfaces-with-boundary of $B^{n+1}_{\r}(Z),$ all having a common $C^{1, \a}$ boundary in $B_{\r}^{n+1}(Z)$ containing $Z$ and no two intersecting except along their common boundary.
\end{itemize}

\noindent
{\bf Remarks}:  {\bf (1)} Note that the stability hypothesis $({\mathcal S{\emph 2}})$ concerns only 
the regular part ${\rm reg} \, V$, and by Allard's regularity theorem, ${\rm reg} \, V \neq \emptyset$---in fact ${\rm reg} \, V$ is an open, dense subset of ${\rm spt} \, \|V\|$---whenever $V$ is stationary (\cite{AW}, Section 8.1). Thus given hypothesis 
$({\mathcal S}\emph{1})$, hypothesis $({\mathcal S{\emph 2}})$ is never vacuously true. However an open, dense subset can have arbitrarily small positive measure, so it is not at all obvious whether hypothesis $({\mathcal S{\emph 2}})$ is sufficiently strong to give any control over the singular set. By our main theorem (Theorem~\ref{compactness} below) however, we conclude that for $V\in {\mathcal S}_{\a},$ ${\rm sing} \, V$ must in fact be very low dimensional. 

\noindent
{\bf (2)} The hypothesis ${\mathcal H}^{n-1} \, ({\rm sing} \, V) = 0$ trivially implies $({\mathcal S{\emph 3}}),$ so all of our theorems concerning the class ${\mathcal S}_{\a}$ in particular apply to the class of stable minimal hypersurfaces $M$ of $B_{2}^{n+1}(0)$ (that is, smooth embedded hypersurfaces $M$ of $B_{2}^{n+1}(0)$ with their associated multiplicity 1 varifolds $V=|M|$ satisfying (${\mathcal S}{\emph 1}$) and (${\mathcal S}{\emph 2}$)) with no removable singularities 
(thus, if $X \in \overline{M} \cap B_{2}^{n+1}(0)$ and $\overline{M}$ is a smooth, embedded hypersurface near $X$, then $X \in M$) and with $${\mathcal H}^{n-1} \, ({\rm sing} \, M) = 0,$$
\noindent
where ${\rm sing} \, M = ({\overline M} \setminus M) \cap B_{2}^{n+1}(0).$ In fact, by Theorem~\ref{compactness}, these two classes are the same.

\noindent
{\bf (3)} Hypothesis $({\mathcal S{\emph 3}})$ will of course be satisfied if no tangent cone to $V$ at a singular point is supported by a union of three or more distinct $n$-dimensional half-hyperplanes meeting along an $(n-1)$-dimensional subspace. By Theorem~\ref{no-transverse} below, for stable codimension 1 integral varifolds, this condition on the tangent cones is in fact equivalent to hypothesis $({\mathcal S}{\emph 3}).$ 

%\medskip

Our main theorem concerning the varifolds in ${\mathcal S}_{\a}$ is the following:

%\medskip

\begin{theorem}[{\sc Regularity and Compactness Theorem}]\label{compactness}
Let $\a \in (0, 1).$ Let $\{V_{k}\} \subset {\mathcal S}_{\a}$ be a sequence with 
$$\limsup_{k \to \infty} \, \|V_{k}\|(B_{2}^{n+1}(0))< \infty.$$ 
There exist a subsequence $\{k^{\prime}\}$ of $\{k\}$ and a varifold $V \in {\mathcal S}_{\a}$  with
${\mathcal H}^{n-7+\g} \,({\rm sing} \ V \cap B_{2}^{n+1}(0)) = 0$ for each $\g >0$ if $n \geq 7$, ${\rm sing} \, V \cap B_{2}^{n+1}(0)$ discrete if $n= 7$ and ${\rm sing} \,V \cap B_{2}^{n+1}(0) = \emptyset$ if $2 \leq n \leq 6$  
such that $V_{k^{\prime}} \to V$ as varifolds of $B_{2}^{n+1}(0)$ and smoothly (i.e. 
in the $C^{m}$ topology for every $m$) locally in $B_{2}^{n+1}(0) \setminus {\rm sing} \, V$.

In particular, if $W \in {\mathcal S}_{\a}$, then ${\mathcal H}^{n-7+\g} \,({\rm sing} \, W \cap B_{2}^{n+1}(0)) = 0$ for each $\g >0$ if $n \geq 7$, ${\rm sing} \, W \cap B_{2}^{n+1}(0)$ is discrete if $n= 7$ and ${\rm sing} \,W \cap B_{2}^{n+1}(0) = \emptyset$ if $2 \leq n \leq 6.$  
\end{theorem} 

%\medskip

Note that we do not a priori assume orientability of ${\rm reg} \, V$ for $V \in {\mathcal S}_{\a}$; indeed, by virtue of low dimensionality of ${\rm sing} \, V$ guaranteed by 
Theorem~\ref{compactness}, orientability of ${\rm reg} \, V$ \emph{follows} if $V \in {\mathcal S}_{\a}$:

\begin{corollary}\label{orientable}
If $V \in {\mathcal S}_{\alpha},$ then ${\rm reg} \, V$ is orientable.
\end{corollary} 

%\medskip

Our proof of Theorem~\ref{compactness} will be based on the following two theorems:

\begin{theorem}[{\sc Sheeting Theorem}] \label{main} Let $\a \in (0, 1).$ Corresponding to each $\Lambda\in [1, \infty)$ and $\th \in (0, 1)$, there exists a number $\e_{0} = \e_{0}(n, \Lambda,\a, \th) \in (0, 1)$ such that if $V \in {\mathcal S}_{\a}$, 
$(\omega_{n}2^{n})^{-1}\|V\|(B_{2}^{n+1}(0)) \leq \Lambda$ and 
${\rm dist}_{\mathcal H} \, ({\rm spt} \, \|V\| \cap ({\mathbf R} \times B_{1}),\{0\} \times  B_{1}) < \e_{0}$ then 
$$V \,\res \, ({\mathbf R} \times B_{\th}) = \sum_{j=1}^{q} \, |{\rm graph} \, u_{j}|$$
for some integer $q,$ where $u_{j} \in C^{1, \b}(B_{\th})$ for each $j=1, 2, \ldots, q$; $u_{1} \leq u_{2} \leq \ldots \leq u_{q};$ 
$$\sup_{B_{\th}} \, \left(|u_{j}| + |Du_{j}|\right) + \sup_{X_{1}, X_{2} \in B_{\th}, \, X_{1} \neq X_{2}} \, \frac{|Du_{j}(X_{1}) - Du_{j}(X_{2})|}{|X_{1} - X_{2}|^{\b}}\leq C\left(\int_{{\mathbf R} \times B_{1}}|x^{1}|^{2} \, d\|V\|(X)\right)^{1/2};$$
Furthermore,  $u_{j}$ solves the minimal surface equation weakly on $B_{\th}$  and hence in fact $u_{j} \in C^{\infty}(B_{\th})$ for each $j=1, 2, \ldots, q$. Here $C = C(n, \Lambda, \a, \th) \in (0, \infty)$ and 
$\b  = \b(n, \Lambda, \a, \th) \in (0, 1)$.
\end{theorem}

%\medskip

\begin{theorem}[{\sc Minimum Distance Theorem}]\label{no-transverse} 
Let $\a \in (0, 1)$. Let $\d \in (0, 1/2)$ and ${\mathbf C}_{0}$ be an $n$-dimensional stationary cone in ${\mathbf R}^{n+1}$ such that 
${\rm spt} \, \|{\mathbf C}_{0}\|$ is equal to a finite union of at 
least three distinct $n$-dimensional half-hyperplanes of ${\mathbf R}^{n+1}$ meeting along an $(n-1)$-dimensional subspace. There exists $\e = \e(n, \a, \d, {\mathbf C}_{0}) \in (0, 1)$  such that if $V \in {\mathcal S}_{\a},$ $\Theta \, (\|V\|, 0) \geq \Theta \, (\|{\mathbf C}_{0}\|, 0)$ and
$(\omega_{n}2^{n})^{-1}\|V\|(B_{2}^{n+1}(0)) \leq \Theta_{{\mathbf C}_{0}}(0) +\d$, then 
$${\rm dist}_{\mathcal H} \, ({\rm spt} \, \|V\| \cap B_{1}^{n+1}(0), {\rm spt} \, \|{\mathbf C}_{0}\| \cap B_{1}^{n+1}(0)) \geq \e.$$
\end{theorem}

The proofs of Theorems~\ref{compactness},~\ref{main}~and~\ref{no-transverse} will be given in Sections~\ref{reg-compactness},~\ref{sheeting} and~\ref{MD} respectively.

\noindent
{\bf Remark}:  Theorems ~\ref{compactness}, ~\ref{main} and ~\ref{no-transverse} are optimal in several ways:

\noindent
{\bf(a)} Examples such as pairs of transverse hyperplanes or a union of three half-hyperplanes meeting at $120^{o}$ angles along a common axis show that  Theorems~\ref{main}, \ref{no-transverse} and \ref{compactness} do not hold  if the structural hypothesis (${\mathcal S{\emph 3}}$)  is removed (or replaced by  the condition ${\mathcal H}^{n-1+\g} \, ({\rm sing} \, V) = 0$ for any $\g >0$). Stable \emph{branched} minimal hypersurfaces (e.g. those constructed in \cite{SW1} or
in \cite{R}) show that in the absence of hypothesis (${\mathcal S{\emph 3}}$), even when $n=2$, there is no hope of proving regularity of stable codimension 1 integral varifolds away from the set of points where (${\mathcal S{\emph 3}}$) fails. Thus hypothesis (${\mathcal S}{\emph 3}$) can in particular be viewed as a geometric condition that implies non-existence of branch points in stable codimension 1 integral varifolds.

\noindent 
{\bf(b)} Appropriate rescalings of a standard 2-dimensional Catenoid  in ${\mathbf R}^{3}$ show that Theorem~\ref{main} does not hold without the stability hypothesis (${\mathcal S}{\emph 2}$). Similarly, rescalings of a Scherk's 2nd surface show that Theorem~\ref{no-transverse} does not hold without (${\mathcal S}{\emph 2}$). It is however an open question, even when $n=2$, 
whether some form of Theorem~\ref{compactness} giving a bound on the singular set holds without (${\mathcal S}{\emph 2}$). In fact it remains open whether  2 dimensional stationary integral varifolds in ${\mathbf R}^{3}$ must be regular almost everywhere, even subject to 
a condition such as (${\mathcal S}{\emph 3}$).

\noindent 
{\bf(c)} There are many examples provided by 
complex algebraic varieties demonstrating that Theorems~\ref{main} and ~\ref{compactness} do not hold in codimension $>1$ even if the stability hypothesis (${\mathcal S}{\emph2}$) (where the corresponding higher codimensional  stability inequality takes a different form from (\ref{stability}); see \cite{S1}, Section 9) is replaced by the (stronger) absolutely area minimizing hypothesis. For instance, the  holomorphic varieties $V_{t} = \{(z, w) \, : \, z^{2} = tw^{3} + tw\} \cap B_{1}^{4}(0)  \subset {\mathbf C} \times {\mathbf C} \equiv {\mathbf R}^{4}$,  $t \in {\mathbf R}$, which are smooth, embedded area minimizing submanifolds lying close to the plane $\{z=0\} \cap B_{1}^{4}(0)$  for small $|t| \neq 0$, show that Theorem~\ref{main} does not hold in codimension $>1$. Those holomorphic varieties with branch point singularities such as $V = \{(z, w) \, : \, z^{2} = w^{3}\} \cap B_{1}^{4}(0)
\subset {\mathbf C} \times {\mathbf C}$ show that even in 2 dimension, $C^{2}$ regularity, and hence Theorem~\ref{compactness}, is false if codimension $> 1$. (For area minimizing currents of dimension $n$ and arbitrary codimension, Almgren's theorem (\cite{A}) gives the optimal bound on the Hausdorff dimension of the interior singular sets, namely, $n-2.$) Since the cone ${\mathbf C}_{0}$ in 
Theorem~\ref{no-transverse} is not area minimizing, there are no area minimizing examples nearby, but 
a given transverse pair of planes in ${\mathbf R}^{3} \times\{0\} \subset {\mathbf R}^{4},$ for instance, can be perturbed in ${\mathbf R}^{4}$ into a union of two planes intersecting only at the origin, and the latter union is of course stable and satisfies (${\mathcal S}{\emph 3}$), showing that Theorem~\ref{no-transverse} is false in codimension $> 1.$

Our theorems generalize the regularity and compactness theory of R. Schoen and L. Simon \cite{SS}, which established Theorems~\ref{main} and ~\ref{compactness} for stable codimension 1 integral varifolds $V$ on $B_{2}^{n+1}(0)$ under the hypothesis 
${\mathcal H}^{n-2} \, ({\rm sing} \, V \cap K) < \infty$ for each compact $K \subset B_{2}^{n+1}(0)$ in place of our hypothesis (${\mathcal S{\emph 3}}$). (Under this more stringent hypothesis on the singular set, Theorem~\ref{no-transverse} is a straightforward consequence of Theorem~\ref{main} and inequality (\ref{stability}).) Our proofs of Theorem~\ref{main} and Theorem~\ref{no-transverse} however rely on the Schoen-Simon version of Theorem~\ref{main} in an essential way; in fact, what we need is the following slightly weaker version of their theorem:

\begin{theorem}[\cite{SS}, special case of Theorem 2]\label{SS}
Let $V$ be an integral $n$-varifold on $B_{2}^{n+1}(0)$ and assume in place of (${\mathcal S{\emph 3}}$) the (stronger) condition that ${\mathcal H}^{n-7 +\g} \, ({\rm sing} \, V) = 0$ for every $\g >0$ 
in case $n \geq 7$ and ${\rm sing} \, V = \emptyset$ in case 
$2 \leq n \leq 6.$ Let all other hypotheses be as in Theorem~\ref{main}. Then the conclusions of Theorem~\ref{main} hold.  
\end{theorem}

\noindent
{\bf Remark:} It suffices to prove Theorem~\ref{main} for $\th = 1/8$ and arbitrary $\Lambda \in [1, \infty).$ To see this, suppose that case $\th = 1/8$ of the theorem is true, with $\e^{\prime} = \e^{\prime}(n, \a, \Lambda) \in (0, 1)$ corresponding to $\e_{0}.$ Let $\th \in (1/8, 1)$ and let the hypotheses be as in the theorem with $\e_{0} = \e_{0}(n, \a, \Lambda, \th) \in (0, 1)$ satisfying $\e_{0} < \left(\frac{1-\th}{8}\right)\e^{\prime}(n, \a, 3^{n}\Lambda).$  We may then apply the case $\th=1/8$ of the theorem with $3^{n}\Lambda$ in place of $\Lambda$ and with $\widetilde{V} = \left(\eta_{Z, (1-\th)/2}\right)_{\#} \, V \in {\mathcal S}_{\a}$ in place of $V,$ where $Z \in {\rm spt} \, \|V\| \cap ({\mathbf R} \times B_{\th})$ is arbitrary; since we may cover ${\rm spt} \, \|V\| \cap ({\mathbf R} \times B_{\th})$ by a collection of balls $B_{(1-\th)/2}^{n+1}(Z_{j})$, $j=1, 2, \ldots, N,$ with $Z_{j} \in {\rm spt} \, \|V\| \cap ({\mathbf R} \times B_{\th})$ and $N  = N(n, \Lambda, \th)$, the required estimate  
follows. 

So assume $\th = 1/8$ and let the hypotheses be as in Theorem~\ref{main}. It follows from Allard's integral varifold compactness theorem (\cite{AW}, Theorem 6.4) and the Constancy Theorem for stationary integral varifolds (\cite{S1}, Theorem 41.1) that if $\e_{0}= \e_{0}(n, \Lambda) \in (0, 1)$ is sufficiently small, then  
$q - 1/2 \leq \left(\omega_{n}R^{n}\right)^{-1}\|V\|({\mathbf R} \times B_{R}) < q + 1/2$ for some integer $q \in [1, \Lambda+1)$ and 
$R \in \{1/3, 2/3\}$. Then $V_{1} \equiv \eta_{0, 1/3 \, \#} \, V$ satisfies $(\omega_{n}2^{n})^{-1}\|V_{1}\|(B_{2}^{n+1}(0)) < q + 1/2$ and $q - 1/2 \leq \omega_{n}^{-1}\|V_{1}\|({\mathbf R} \times B_{1}) < q + 1/2.$ Thus in order to prove the special case $\th=1/8$ of Theorem~\ref{main} (and therefore the general version), it suffices to establish the following:

\noindent
{\bf Theorem ${\bf 3.3^{\prime}}$\; {\sc (Sheeting Theorem)}.}\label{main-special}
{\em Let $\a \in (0, 1)$. Let $q$ be any integer $\geq 1.$ There exists a number 
$\e_{0} = \e_{0}(n, \a, q) \in (0, 1)$ such that if $V \in {\mathcal S}_{\a}$, $(\omega_{n}2^{n})^{-1}\|V\|(B_{2}^{n+1}(0)) < q + 1/2,$ 
$q-1/2 \leq \omega_{n}^{-1}\|V\| ({\mathbf R} \times B_{1}) < q + 1/2$ and 
${\rm dist}_{\mathcal H} \, ({\rm spt} \, \|V\| \cap ({\mathbf R} \times B_{1}), \{0\} \times B_{1}) < \e_{0}$ then 
$$V \res ({\mathbf R} \times B_{3/8}) = \sum_{j=1}^{q} \, |{\rm graph} \, u_{j}|$$
where $u_{j} \in C^{1, \b}(B_{3/8})$ for each $j=1, 2, \ldots, q;$ $u_{1} \leq u_{2} \leq \ldots \leq u_{q};$ 
$$\sup_{B_{3/8}} \, \left(|u_{j}| + |Du_{j}|\right) + \sup_{X_{1}, X_{2} \in B_{3/8},\, X_{1} \neq X_{2}} \, \frac{|Du_{j}(X_{1}) - Du_{j}(X_{2})|}{|X_{1} - X_{2}|^{\b}}\leq C\left(\int_{{\mathbf R} \times B_{1}}|x^{1}|^{2} \, d\|V\|(X)\right)^{1/2};$$
\noindent
and $u_{j}$ solves the minimal surface equation (weakly) on $B_{3/8}$. Here $C = C(n, q, \a) \in (0, \infty)$ and $\b = \b(n, q, \a) \in (0, 1)$.}

%\medskip

Finally, we note that in the absence of the $\a$-Structural Hypothesis (${\mathcal S{\emph 3}}$), Theorems~\ref{compactness},  ~\ref{main} and the upper semi-continuity of density of stationary integral varifolds readily imply the following:

\begin{corollary}
Let $V$ be a stable integral $n$-varifold on $B_{2}^{n+1}(0)$ (in the sense that $V$ satisfies
(\ref{firstvariation}) and (\ref{stability})). 

If $Z \in {\rm sing} \, V$ and one of the tangent cones to $V$ at $Z$ is 
(the varifold associated with) a hyperplane with multiplicity $q \in \{2, 3, \ldots\}$, then for any $\a \in (0, 1)$, there exist a sequence of points  $Z_{j} \in {\rm sing} \, V$ with $Z_{j} \neq Z$, $Z_{j} \to Z$ and a sequence of numbers 
$\s_{j}$ with $0< \s_{j} <  |Z_{j} - Z|$ such that for each $j=1, 2, 3, \ldots$,  
${\rm spt} \, \|V\| \cap B_{\s_{j}}^{n+1}(Z_{j})$ is the union of at least 3 and at most $2q$ embedded $C^{1,\a}$ hypersurfaces-with-boundary meeting only along an $(n-1)$-dimensional $C^{1, \a}$ submanifold of $B_{\s_{j}}^{n+1}(Z_{j})$ containing $Z_{j}.$

In fact, if $Z \in {\rm sing} \, V$ is such that one tangent cone ${\mathbf C}$ to $V$ at $Z$ has the form, after a rotation, 
${\mathbf C} = {\mathbf C}^{\prime} \times {\mathbf R}^{n-k}$ for some $k \in \{0, 1, \ldots, \min \, \{6, n\}\}$, then for any $\a \in (0, 1)$, there exist a sequence of points  $Z_{j} \in {\rm sing} \, V$ with $Z_{j} \neq Z$, $Z_{j} \to Z$ and a sequence of numbers $\s_{j}$ with $0< \s_{j} < |Z_{j} - Z|$ such that for each $j=1, 2, 3, \ldots$,  
${\rm spt} \, \|V\| \cap B_{\s_{j}}^{n+1}(Z_{j})$ is the union of at least 3 and at most $2\Theta \, (\|V\|, Z)$ embedded $C^{1,\a}$ hypersurfaces-with-boundary meeting only along an $(n-1)$-dimensional $C^{1, \a}$ submanifold of $B_{\s_{j}}^{n+1}(Z_{j})$ containing $Z_{j}.$
\end{corollary}

\section{Proper blow-up classes}\label{proper-blow-up}
\setcounter{equation}{0}

Fix an integer $q \geq 1$  and a constant $C \in (0, \infty).$ Consider a
family ${\mathcal B}$ of functions $v = (v^{1}, v^{2}, \ldots, v^{q}) \, : \, B_{1} \to {\mathbf R}^{q}$ satisfying the following properties:

\begin{itemize}
\item[(${\mathcal B{\emph 1}}$)] ${\mathcal B} \subset W^{1, 2}_{\rm loc} \, (B_{1}; {\mathbf R}^{q}) \cap L^{2} \, (B_{1}; {\mathbf R}^{q}).$
\item[(${\mathcal B{\emph 2}}$)] If $v \in {\mathcal B}$, then $v^{1} \leq v^{2} \leq \ldots \leq v^{q}$.
\item[(${\mathcal B{\emph 3}}$)] If $v \in {\mathcal B},$ then $\Delta \, v_{a} = 0$ in $B_{1}$ where $v_{a} = q^{-1}\sum_{j=1}^{q} v^{j}.$
\item[(${\mathcal B{\emph 4}}$)] For each $v \in {\mathcal B}$ and each $z \in B_{1}$, either (${\mathcal B}{\emph 4 \, I}$) or (${\mathcal B}{\emph 4 \, II}$) below is true: 
\begin{itemize}
\item[(${\mathcal B{\emph 4 \,I}}$)] The {\em Hardt-Simon inequality} $$\sum_{j=1}^{q} \int_{B_{\r/2}(z)} R_{z}^{2-n} \left(\frac{\partial \, \left((v^{j} - v_{a}(z))/R_{z}\right)}{\partial \, R_{z}}\right)^{2} \leq C \, \r^{-n-2}\int_{B_{\r}(z)} |v - \ell_{v,\, z}|^{2}$$ 
holds for each $\r \in (0, \frac{3}{8}(1- |z|)]$, where $R_{z}(x) = |x - z|,$ $\ell_{v, \, z}(x) = v_{a}(z) + Dv_{a}(z) \cdot (x -z)$ 
and $v - \ell_{v, \, z} = (v^{1} - \ell_{v,\,z}, v^{2} - \ell_{v, \, z}, \ldots, v^{q} - \ell_{v, \, z}).$
\item[(${\mathcal B{\emph 4 \, II}}$)] There exists $\s  = \s(z) \in (0, 1 - |z|]$ such that $\Delta \, v = 0$ in $B_{\s}(z).$
\end{itemize}
\item[(${\mathcal B{\emph 5}}$)] If $v \in {\mathcal B}$, then
\begin{itemize}
\item[(${\mathcal B{\emph 5 \, I}}$)]${\widetilde v}_{z, \s}(\cdot) \equiv \|v(z + \s(\cdot))\|_{L^{2}(B_{1}(0))}^{-1}v(z + \s(\cdot)) \in {\mathcal B}$ for each $z \in B_{1}$ and $\s \in (0, \frac{3}{8}(1 - |z|)]$ whenever $v \not\equiv 0$ in $B_{\s}(z);$
\item[(${\mathcal B{\emph 5 \, II}}$)] $v \circ \g \in {\mathcal B}$ for each orthogonal rotation $\g$ of ${\mathbf R}^{n}$ and
\item[(${\mathcal B{\emph 5 \, III}}$)] $\|v - \ell_{v}\|_{L^{2}(B_{1}(0))}^{-1}\left(v- \ell_{v}\right) \in {\mathcal B}$ whenever $v  - \ell_{v}\not\equiv 0$ in $B_{1}$, where  $\ell_{v}(x) = v_{a}(0) + Dv_{a}(0) \cdot x$ for $x \in {\mathbf R}^{n}$ and $v- \ell_{v} =(v^{1} - \ell_{v},v^{2} -  \ell_{v}, \ldots, v^{q} - \ell_{v}).$ 
\end{itemize}
\item[(${\mathcal B{\emph 6}}$)] If $\{v_{k}\}_{k=1}^{\infty} \subset {\mathcal B}$ then there exists a subsequence $\{k^{\prime}\}$ 
of $\{k\}$ and a function $v \in {\mathcal B}$ such that $v_{k^{\prime}} \to v$ locally in $L^{2}(B_{1})$ and locally weakly in $W^{1, 2}(B_{1}).$
\item[(${\mathcal B{\emph 7}}$)] If $v \in {\mathcal B}$ is such that for each $j=1, 2, \ldots, q$, there exist linear functions 
$L^{j}_{1}, L^{j}_{2} \; : \; {\mathbf R}^{n} \to {\mathbf R}$ with $v^{j}(x^{2}, y) = L^{j}_{1}(x^{2}, y)$ if $x^{2} > 0$, $v^{j}(x^{2}, y) = 
L^{j}_{2}(x^{2}, y)$ if $x^{2} \leq 0$ and $L^{j}_{1}(0, y) = L^{k}_{2}(0, y)$ for $1 \leq j, k \leq q,$ 
$y \in  {\mathbf R}^{n-1}$     
then 
$v^{1} = v^{2} = \ldots = v^{q} = L$ for some linear function $L \; : \; {\mathbf R}^{n} \to {\mathbf R}.$ 
\end{itemize}

We shall refer to any such class ${\mathcal B}$ as a {\em proper blow-up class}.

Our main result in this section (Theorem~\ref{blowup-reg} below) is that \emph{functions in any proper blow-up class are harmonic}. Subsequently, we shall prove that the collection of functions arising as ``coarse blow-ups''  (see Section~\ref{blow-up}
for the definition) of mass-bounded sequences of varifolds in ${\mathcal S}_{\a}$  converging weakly to a hyperplane is a proper blow-up class  for a suitable constant $C$ depending only on $n$ and the mass bound.

\noindent
{\bf Remark:} The first use of the inequality in (${\mathcal B}{\emph 4 \,I}$) in the context of regularity theory for minimal submanifolds is due to R. Hardt and L. Simon (\cite{HS}).

Let ${\mathcal B}$ be a proper blow-up class. There exists a   constant $\t = \t({\mathcal B}) \in (0, 1/4)$ such that 
if $v \in {\mathcal B},$ $v_{a}(0) = 0$ and property (${\mathcal B}{\emph 4 \, I}$) holds with $z =0$, then 
\begin{equation}\label{blow-up-non-conc}
\int_{B_{1} \setminus B_{\t}}|v|^{2} \geq \frac{1}{2}\int_{B_{1}}|v|^{2}.
\end{equation}  
To see this, note that since every weakly convergent sequence in $W^{1, 2} \, (B_{2/3})$ is bounded in $W^{1,2} \, (B_{2/3})$, it follows from the compactness property (${\mathcal B}{\emph 6}$) and property (${\mathcal B}{\emph 5 \, I}$) that there exists a constant $C_{1}= C_{1}({\mathcal B}) \in (0, \infty)$ 
such that $\int_{B_{1/4}}|Dv|^{2} \leq C_{1}\int_{B_{1}}|v|^{2}$ for every $v \in {\mathcal B}.$ Hence by property (${\mathcal B}{\emph 4 \, I}$) with $z = 0$ and $\r  = 3/8$, we see that if $v_{a}(0) = 0$ then 
$$\int_{B_{3/16}} \frac{|v|^{2}}{R^{2}} \leq 2(C_{2} + C_{1})\int_{B_{1}}|v|^{2},$$ where $C_{2} = C_{2}(C, n)$ and we have used 
the fact that, since $v_{a}$ is harmonic, $|\ell_{v, 0}(x)|^{2}= |Dv_{a}(0)|^{2}|x|^{2} \leq C_{3}\int_{B_{1/4}}|Dv|^{2} \leq C_{3}C_{1}\int_{B_{1}}|v|^{2},$ with $C_{3} = C_{3}(n, q).$  
This readily implies that for each $\t \in (0, 3/16)$, $\int_{B_{\t}}|v|^{2} \leq 2(C _{2}+ C_{1}) \t^{2}\int_{B_{1}}|v|^{2},$ and choosing 
$\t = \t({\mathcal B}) \in (0, 3/16)$ such that  $2(C_{2} + C_{1})\t^{2} < 1/2$, we deduce (\ref{blow-up-non-conc}).

\begin{theorem}\label{blowup-reg}
If ${\mathcal B}$ is a proper blow-up class for some $C \in (0, \infty)$, then 
each $v \in {\mathcal B}$ is harmonic in $B_{1}$. Furthermore, if $v \in {\mathcal B}$ and there is a point $z \in B_{1}$ such that 
(${\mathcal B{\it{4} \, I}}$) is satisfied, then $v^{1} = v^{2} = \ldots = v^{q}.$
\end{theorem}

The proof of this Theorem will be based on the following Proposition:

\begin{proposition}\label{homog}
Let ${\mathcal B}$ be a proper blow-up class, and let $\t = \t({\mathcal B}) \in (0, 1/4)$ be the constant as in (\ref{blow-up-non-conc}). If $v \in {\mathcal B}$  satisfies property (${\mathcal B}{\it{4} \, I}$) with $z = 0$ and if $v$ is homogeneous of degree 1
in the annulus $B_{1} \setminus B_{\t}$, viz. $\frac{\partial \, \left(v/R\right)}{\partial \, R} = 0$ a.e. in $B_{1} \setminus B_{\t},$ then $v^{j} = L$ in $B_{1}$ for some linear function $L$ and all $j \in \{1, 2, \ldots, q\}.$ 
\end{proposition}

For the proofs of Theorem~\ref{blowup-reg}, Proposition~\ref{homog} and subsequently, we shall need the following general principle:

\begin{lemma}\label{general}
Let $w \in L^{2}(B_{1}; {\mathbf R}^{q}).$ Suppose there is a closed subset $\G \subset B_{1}$ and numbers $ \b, \b_{1}, \b_{2} \in (0, \infty)$,  $\m  \in (0, 1)$ and $\e \in (0, 1/4)$ such that the following hold: For each $z \in \G \cap B_{3/4},$ there is an affine function $\ell_{z} \, : \, {\mathbf R}^{n} \to {\mathbf R}^{q}$ with $\sup_{B_{1}} \, |\ell_{z}| \leq \b$ such that  
$$\s^{-n-2}\int_{B_{\s}(z)} |w - \ell_{z}|^{2} \leq \beta_{1} \left(\frac{\s}{\r}\right)^{\m}\r^{-n-2}\int_{B_{\r}(z)}
|w - \ell_{z}|^{2}$$ 
\noindent
for all $0< \s \leq \r/2 \leq \e/2$ and  for each $z \in B_{3/4} \setminus \G$, there is an affine function $\ell_{z} \, : \, {\mathbf R}^{n} \to {\mathbf R}^{q}$ such that 
$$\s^{-n-2}\int_{B_{\s}(z)} |w - \ell_{z}|^{2} \leq \beta_{2} \left(\frac{\s}{\r}\right)^{\m}\r^{-n-2}\int_{B_{\r}(z)}
|w - \ell|^{2}$$ 
for each affine function $\ell \, : \, {\mathbf R}^{n} \to {\mathbf R}^{q}$ and all $0 < \s \leq \r/2 < \frac{1}{2}\min \, \{1/4, {\rm dist} \, (z, \G)\}.$ Then $w \in C^{1, \lambda}(B_{1/2})$ for some $\lambda = \lambda(n, q, \b_{1}, \b_{2}, \e, \m) \in (0, 1)$ with 
$$\sup_{B_{1/2}} \, \left(|w| + |Dw|\right) + \sup_{x, y \in B_{1/2}, x \neq y} \, \frac{|Dw(x) - Dw(y)|}{|x - y|^{\lambda}}
\leq C\left(\b^{2} + \int_{B_{1}}|w|^{2}\right)^{1/2}$$

\noindent
where $C = C(n, q, \b_{1}, \b_{2}, \e) \in (0, \infty).$
\end{lemma}

\noindent
{\bf Remark:} In our applications of the lemma, the component functions of $w$, in $B_{1} \setminus \G$, will either be harmonic or smooth functions with small gradient solving the minimal surface equation; the second estimate in the hypotheses, with $\ell_{z}(x) = w(z) + Dw(z) \cdot (x-z)$ and $\b_{2}$ depending only on $n,$ follows in these cases from standard interior estimates for second derivatives of harmonic functions and solutions to uniformly elliptic equations.

\begin{proof} Consider an arbitrary point 
$y \in B_{3/4}$ and a number $\r \in (0, \e).$ With $\g = \g(n, \b_{1}, \e, \m) \in (0, 1/8)$ to be chosen, if there is a point $z \in \G \cap \overline{B_{\g\r}(y)}$, then by the given condition
with  $\r - |z - y|$ in place of $\r$ and $\s = \g\r + |z-y|$, 
\begin{eqnarray*}
&&(\g\r)^{-n-2}\int_{B_{\g\r}(y)} |w - \ell_{z}|^{2} \leq \left(1 + \frac{|z-y|}{\g\r}\right)^{n+2} (\g\r + |z-y|)^{-n-2} \int_{B_{\g\r + |z-y|}(z)} |w - \ell_{z}|^{2} \nonumber\\ 
&&\hspace{1in}\leq 2^{n+2}\b_{1}\left(\frac{\g\r + |z-y|}{\r - |z-y|}\right)^{\m} (\r - |z - y|)^{-n-2}\int_{B_{\r - |z-y|}(z)} |w - \ell_{z}|^{2}\nonumber\\ 
&&\hspace{1in}\leq 4^{n+2}\b_{1} \left(\frac{2\g}{1 - \g}\right)^{\m} \r^{-n-2}\int_{B_{\r}(y)} |w - \ell_{z}|^{2}.
\end{eqnarray*}
Choosing $\g = \g(n, \b_{1}, \e, \m) \in (0, \e)$ such that $4^{n+2}\b_{1} \left(\frac{2\g}{1 - \g}\right)^{\m} < 1/4$, 
we see from this that 
\begin{equation*}\label{general-1}
(\g\r)^{-n-2}\int_{B_{\g\r}(y)} |w - \ell_{z}|^{2}  
\leq 4^{-1}\r^{-n-2}\int_{B_{\r}(y)} |w - \ell_{z}|^{2}
\end{equation*} 
for any $y \in B_{3/4}$ and $\r \in (0, \e)$ provided there is a point $z \in \G \cap \overline{B_{\g\r}(y)}.$ In particular, 
if $z^{\star} \in \G$ is such that $|y - z^{\star}| = {\rm dist} \, (y, \G),$ then 
\begin{equation}\label{general-1-1}
(\g\r)^{-n-2}\int_{B_{\g\r}(y)} |w - \ell_{z^{\star}}|^{2}  
\leq 4^{-1}\r^{-n-2}\int_{B_{\r}(y)} |w - \ell_{z^{\star}}|^{2}
\end{equation} 
for each $\r \in (0, \e)$ such that $\g\r \geq|y - z^{\star}|.$ If on the other hand $\G \cap \overline{B_{\g\r}(y)}  = \emptyset$, then again by the given condition we know that for any affine function $\ell$
 \begin{equation}\label{general-2}
(\s\g\r)^{-n-2}\int_{B_{\s\g\r}(y)}|w - \ell_{y}|^{2} \leq \b_{2}\s^{\m}(\g\r)^{-n-2}\int_{B_{\g\r}(y)}|w - \ell|^{2}
\end{equation}
for all $\s \in (0, 1/2].$ Iterating  inequality (\ref{general-1-1}) with $\r = \g^{j}$, $j=1, 2, \ldots$ and using inequality (\ref{general-2}), we see that for each $y \in B_{3/4} \setminus \G$, there is an integer $j^{\star} \geq 1,$ an affine function $\ell_{\star}$ ($= \ell_{z^{\star}}$) with $\sup_{B_{1}} \, |\ell_{\star}| \leq \b$ and an affine 
function $\ell_{y}$ such that 
\begin{equation}\label{general-4}
(\s\g^{j^{\star}+1})^{-n-2}\int_{B_{\s\g^{j^{\star} +1}}(y)}|w - \ell_{y}|^{2} \leq \b_{2}\s^{\m}(\g^{j^{\star}+1})^{-n-2}\int_{B_{\g^{j^{\star}+1}}(y)}|w - \ell|^{2}
\end{equation}
for each affine function $\ell$ and each $\s \in (0, 1/2];$ and 
\begin{eqnarray}\label{general-5}
(\g^{j})^{-n-2}\int_{B_{\g^{j}}(y)}|w - \ell_{\star}|^{2} &\leq& 4^{-1}(\g^{j-1})^{-n-2} \int_{B_{\g^{j-1}}(y)}|w - \ell_{\star}|^{2} \nonumber\\
&\leq& 4^{-(j-1)} \g^{-n-2}\int_{B_{\g}(y)}|w - \ell_{\star}|^{2} \;\;\; \mbox{for each $j= 1, 2, \ldots, j^{\star}.$} 
\end{eqnarray}

By taking $\ell = \ell_{\star}$, $\s =1/2$ in (\ref{general-4}) and $j = j^{\star}$ in (\ref{general-5}), and using the triangle inequality,  we see  that
 $$\left(\frac{1}{2}\g^{j^{\star} + 1}\right)^{-n-2}\int_{B_{\frac{1}{2}\g^{j^{\star}+1}}(y)} |\ell_{y} - \ell_{\star}|^{2} \leq C 4^{-(j^{\star} - 1)}\int_{B_{\g}(y)} |w - \ell_{\star}|^{2}$$
which in particular implies 
\begin{equation}\label{general-6}
\left(\g^{j}\right)^{-n-2}\int_{B_{\g^{j}}(y)}|\ell_{y} - \ell_{\star}|^{2} \leq C4^{-j} \int_{B_{\g}(y)}|w - \ell_{\star}|^{2}
\end{equation}
for $j=1, 2, \ldots, j^{\star},$ where $C = C(n, \b_{2}, \mu, \g) \in (0, \infty)$.

By (\ref{general-5}) and (\ref{general-6}), we conclude that
\begin{equation}\label{general-8}
(\g^{j})^{-n-2}\int_{B_{\g^{j}}(y)} |w- \ell_{y}|^{2} \leq C4^{-(j-1)}\int_{B_{\g}(y)} |w - \ell_{\star}|^{2}
\end{equation}
for each $j=1, 2, \ldots, j^{\star}.$ Thus if $y \in B_{3/4} \setminus \G$, we deduce that
\begin{equation}\label{general-9}
\r^{-n-2}\int_{B_{\r}(y)}|w - \ell_{y}|^{2} \leq C\r^{\lambda} 
\int_{B_{\g}(y)}|w-\ell_{\star}|^{2}
\end{equation}
for all $\r \in (0, \g/2]$, by considering, for any given $\r \in (0, \g/2]$, the two alternatives: (i) $2\r \leq \g^{j^{\star} +1}$, in which case 
$\r = \s\g^{j^{\star} +1}$ for some $\s \in (0, 1/2]$ and we use (\ref{general-4}) provided  
$\g = \g(n, q, \b_{1}, \b_{2}, \m, \e)$ is chosen to satisfy $\g^{\m} < 1/4$ also, or (ii) $\g^{j+1} < 2\r \leq \g^{j}$ for some $j \in \{1, 2, \ldots, j^{\star}\}$, in which case we use (\ref{general-8}). 

 In view of (\ref{general-9}) (in case $y \in B_{3/4} \setminus \G$) and the given condition (in case $y \in B_{3/4} \cap \G$), we conclude that for each $y \in B_{3/4}$, there exists an affine function $\ell_{y}$ such that 
\begin{equation}\label{general-10}
\r^{-n-2}\int_{B_{\r}(y)} |w - \ell_{y}|^{2} \leq C\r^{\lambda} \left(\b^{2} + \int_{B_{1}}|w|^{2}\right)
\end{equation}
for all $\r \in (0, \g/2]$, where $C = C(n, q, \b_{1}, \b_{2}, \m, \e) \in (0, \infty)$ and 
$\lambda = \lambda(n, q, \b_{1}, \b_{2}, \m, \e) \in (0, 1).$ It is standard that from this the assertions of the lemma follow. 
\end{proof}

In the proofs of Proposition~\ref{homog}, Theorem~\ref{blowup-reg} and subsequently, we let, for $v \in {\mathcal B},$ 
$$\G_{v} = \left\{z \in B_{1} \setminus \Omega_{v} \, : \, \mbox{(${\mathcal B{\emph 4 \, I}}$) holds}\right\}$$
where
$$\Omega_{v} = \left\{z \in B_{1} \, : \, \mbox{$\exists \r \in (0, 1-|z|]$ such that}\right.\hspace{2in}$$
$$\hspace{1in}\left.\mbox{$v^{1}(x) =v^{2}(x) = \ldots =v^{q}(x) (= v_{a}(x))$ for a.e. $x \in B_{\r}(z)$}\right\}.$$ 

\noindent
{\bf Remark:} Note that it follows directly from  property (${\mathcal B{\emph 4}}$) that $\G_{v}$ is a relatively closed subset of $B_{1}$ and on $B_{1} \setminus \G_{v}$, $v^{j}$ is a.e. equal to a harmonic function  for each $j=1, 2, \ldots, q.$

\begin{proof}[Proof of Proposition~\ref{homog}]
Let $\t = \t({\mathcal B}) \in (0, 1/4)$ be as in (\ref{blow-up-non-conc}). Note first that if $v \in {\mathcal B}$  is homogeneous of degree 1 in any annulus $B_{1}\setminus B_{\t^{\prime}},$ $\t^{\prime} \in (0, 1)$, viz. $v$ satisfies $\frac{\partial \, \left(v/R\right)}{\partial \, R} = 0$ a.e. in $B_{1} \setminus B_{\t^{\prime}},$ then, since $v_{a} = q^{-1} \sum_{j=1}^{q} v^{j}$ is harmonic  in $B_{1}$ by property (${\mathcal B{\emph 3}}$), it follows that $v_{a}$ is a linear function in $B_{1}.$

Let ${\mathcal H}$ denote the collection of all homogeneous of degree 1 functions $\widetilde{v} \, : \, {\mathbf R}^{n} \to {\mathbf R}^{q}$ such that $\left.\widetilde{v}\right|_{B_{1} \setminus B_{\t}} \equiv \left.v\right|_{B_{1} \setminus B_{\t}}$ for some $v \in {\mathcal B}$ satisfying property (${\mathcal B}{\emph 4 \, I}$) with $z = 0.$ For any given $\widetilde{v} \in {\mathcal H}$, 
let $T(\widetilde{v}) = \{z \in {\mathbf R}^{n} \, : \, \widetilde{v}(x + z) = \widetilde{v}(x) \; \mbox{for a.e.} \; x \in {\mathbf R}^{n}\}.$ It is standard to verify using homogeneity of $\widetilde{v}$ that $T(\widetilde{v})$ is a linear subspace of ${\mathbf R}^{n}.$

For $k=0, 1, 2, \ldots, n$, let ${\mathcal H}_{k} = \{ \widetilde{v} \in {\mathcal H} \, : \, {\rm dim} \, T(\widetilde{v}) = n-k\}$ so that ${\mathcal H} = \cup_{k=0}^{n} {\mathcal H}_{k}.$ Clearly ${\mathcal H}_{0} = \{0\}.$ Let $\widetilde{v} \in {\mathcal H}_{1},$ and let $v$ be any element $\in {\mathcal B}$ which is homogeneous of degree 1 in $B_{1} \setminus B_{\t}$ 
such that $v$ satisfies property (${\mathcal B}{\emph 4 \, I}$) with $z = 0$ and $v$ agrees with $\widetilde{v}$ on $B_{1} \setminus B_{\t}$. We wish to show that there exists a linear function $L$ such that $v^{j} = L$ in $B_{1}$ for each $j \in \{1, \ldots, q\}.$ This is true if $v^{j} = v_{a}$ on $B_{1}$ for each $j \in \{1, \ldots, q\}$, so suppose $v - v_{a} = (v^{1} - v_{a}, \ldots, v^{q} - v_{a}) \not\equiv 0$ in $B_{1}$ and let  $w = \|v -v_{a}\|^{-1}(v - v_{a}).$ Then $w \in {\mathcal B}$ by property (${\mathcal B}{\emph 5 \, III}$), $w \not\equiv 0$, $w_{a} \equiv 0$ and property (${\mathcal B}{\emph 4 \, I}$) is satisfied with $w$ in place of $v$ and $z = 0$, and   hence by (\ref{blow-up-non-conc}), $w \not\equiv 0$ in $B_{1} \setminus B_{\t}.$ By the definition of ${\mathcal H}_{1}$ and property (${\mathcal B{\emph 5 \, II}}$) (of $v$), we may assume that $T(\widetilde{v}) = \{0\} \times {\mathbf R}^{n-1}$, and by homogeneity of $w$ in $B_{1} \setminus B_{\t},$ it then follows that there exist constants $\lambda_{1} \geq \lambda_{2} \geq \ldots \geq  \lambda_{q}$, $\m_{1} \leq \m_{2} \leq  \ldots \leq\m_{q}$, with $\sum_{j=1}^{q} \lambda_{j}= \sum_{j=1}^{q}\m_{j} = 0$ such that, for each $j \in \{1, \ldots, q\}$,   
$w^{j}(x^{2}, y) = \lambda_{j}x^{2}$ for each $(x^{2}, y) \in (B_{1} \setminus B_{\t}) \cap \{x^{2} < 0\}$ and 
$w^{j}(x^{2}, y) = \m_{j}x^{2}$ for each $(x^{2}, y) \in (B_{1} \setminus B_{\t}) \cap \{x^{2} > 0\}.$ Moreover, since $w \not\equiv 0$ in $B_{1} \setminus B_{\t}$, we must have some $j_{0} \in \{1, \ldots, q-1\}$ such that either $\lambda_{j_{0}}  > \lambda_{j_{0}+1}$ or 
$\m_{j_{0}} < \m_{j_{0}+1}$. Thus, taking any point $(0, y_{1}) \in (B_{1} \setminus B_{\t}) \cap (\{0\} \times {\mathbf R}^{n-1})$ and
any number $\s_{1}$ with $0 < \s_{1} < {\rm min} \, \{1 - |y_{1}|, |y_{1}| - \t\}$ and setting $\widetilde{w} = \|w((0, y_{1}) + \s_{1}(\cdot))\|_{L^{2}(B_{1})}^{-1}w((0, y_{1}) + \s_{1}(\cdot))$, we produce an element $\widetilde{w} \in {\mathcal B}$ whose 
existence contradicts the fact that ${\mathcal B}$ satisfies property (${\mathcal B}{\emph 7}$). Hence it must be that $v- v_{a} = 0$ in $B_{1},$ and ${\mathcal H}_{1}$ consists of linear functions. 

Now let $k_{1}$ be the smallest integer $\in \{2, 3, \ldots, n\}$ such that ${\mathcal  H}_{k_{1}} \neq \emptyset.$
Consider any $\widetilde{v} \in {\mathcal H}_{k_{1}},$ and let $v$ be any element $\in {\mathcal B}$ such that $v$ satisfies property (${\mathcal B}{\emph 4 \, I}$) with $z = 0$ and $v$ agrees with $\widetilde{v}$ on $B_{1} \setminus B_{\t}$. By property (${\mathcal B{\emph 5 \, II}}$) (of $v$), we may assume that $T(\widetilde{v}) = \{0\} \times {\mathbf R}^{n-k_{1}}.$ If $\G_{v} \cap (B_{1} \setminus \overline{B_{\t}}) \subseteq \{0\} \times {\mathbf R}^{n-k_{1}}$, then by the remark immediately following the definition of $\G_{v}$, $v^{j}$ is harmonic in $(B_{1} \setminus \overline{B_{\t}}) \setminus \left(\{0\} \times {\mathbf R}^{n-k_{1}}\right)$ for each $j \in \{1, 2, \ldots, q\}$, whence by homogeneity, $\widetilde{v}$ is harmonic on $B_{1} \setminus \left(\{0\} \times {\mathbf R}^{n-k_{1}}\right).$ Since $\widetilde{v}^{j} \in W^{1, 2}_{\rm loc} \, ({\mathbf R}^{n})$ and independent of the last $(n-k_{1})$ variables, it follows that $\widetilde{v}^{j}$ is harmonic in all of ${\mathbf R}^{n}.$ By homogeneity of $\widetilde{v}$ again and property $({\mathcal B{\emph 2}})$ of $v$, it follows that $\widetilde{v}^{1} = \widetilde{v}^{2} = \ldots = \widetilde{v}^{q} = L$ for some linear function $L,$  contrary to the assumption that $\widetilde{v} \in {\mathcal H}_{k_{1}}$ for $k_{1} \geq 2$. So we must have that $\G_{v} \cap (B_{1} \setminus \overline{B_{\t}}) \setminus \left(\{0\} \times {\mathbf R}^{n-k_{1}}\right) \neq \emptyset.$ We shall contradict this also. 

Let $K$ be any compact subset 
of $(B_{1} \setminus \overline{B_{\t}}) \setminus \left(\{0\} \times {\mathbf R}^{n-k_{1}}\right).$ We claim that there exists $\e = \e(v, K, {\mathcal B}) \in (0, 1)$ such that for each $z \in K \cap \G_{v}$ and each $\r$ with $0 < \r  \leq \e$, 
\begin{equation}\label{homog-1}
\sum_{j=1}^{q}\int_{B_{\r}(z) \setminus B_{\t\r}(z)} R_{z}^{2-n} \left(\frac{\partial\, \left((v^{j} - v_{a})/R_{z}\right)}{\partial\, R_{z}}\right)^{2} 
\geq \e \r^{-n-2} \sum_{j=1}^{q}\int_{B_{\r}(z)}|v^{j} - v_{a}|^{2}.
\end{equation}
(Recall that $v_{a}$ is a linear function.) If this were false, then there would exist points $z, z_{i} \in K \cap \G_{v},$ $i=1, 2, 3, \ldots,$ 
with $z_{i} \to z$, and radii $\r_{i} \to 0$ such that $v - v_{a}  \not\equiv 0$ in $B_{\r_{i}}(z_{i})$ for each $i=1, 2, 3, \ldots$ and 
\begin{equation}\label{homog-2}
 \sum_{j=1}^{q}\int_{B_{\r_{i}}(z_{i}) \setminus B_{\t\r_{i}}(z_{i})} R_{z_{i}}^{2-n} \left(\frac{\partial\, \left((v^{j}-v_{a})/R_{z_{i}}\right)}{\partial\, R_{z_{i}}}\right)^{2} 
< \e_{i} \r_{i}^{-n-2}\sum_{j=1}^{q}\int_{B_{\r_{i}}(z_{i})}|v^{j} -v_{a}|^{2}
\end{equation}
where $\e_{i} \to 0^{+}.$ By property (${\mathcal B{\emph 5 \, III}}$), we have that 
$w \equiv \|v - v_{a}\|_{L^{2}(B_{1})}^{-1}(v-v_{a}) \in {\mathcal B}$, so that, by property (${\mathcal B{\emph 5 \, I}}$),
$w_{i}  \equiv w_{z_{i}, \r_{i}}  = \|w(z_{i} + \r_{i}(\cdot))\|_{L^{2}(B_{1})}^{-1}w(z_{i} + \r_{i}(\cdot))$ also belongs to ${\mathcal B}$  for each sufficiently large $i$, and hence, by property 
(${\mathcal B{\emph 6}}$), there exists $w_{\star} \in {\mathcal B}$ such that after passing to 
a subsequence, $w_{i} \to w_{\star}$ locally in $L^{2}(B_{1})$ and locally weakly 
in $W^{1, 2}(B_{1}).$ Since $\|w_{i}\|_{L^{2}(B_{1})} = 1$, it follows from (\ref{homog-2}) that $\|w_{i}\|_{L^{2}(B_{3/4})} > c$ for sufficiently large $i$, where $c = c(n) >0$. Hence $w_{\star} \not\equiv 0$ in $B_{1}.$ In view of the strong convergence $w_{i} \to w_{\star}$  locally in $L^{2}(B_{1})$ and the weak convergence $Dw_{i} \to Dw_{\star}$
locally in $L^{2}(B_{1})$ (which in particular implies that $\int_{B_{1-\e}(0) \setminus B_{\e^{\prime}}(0)} |x|^{-n-2}(Dw_{\star} \cdot x)^{2} \leq \liminf_{i \to \infty} \, \int_{B_{1-\e}(0) \setminus B_{\e^{\prime}}(0)} |x|^{-n-2}(Dw_{i} \cdot x)^{2}$ for any $\e, \e^{\prime} \in (0, 1/4)$), it follows from (\ref{homog-2}) that $w_{\star}$ is homogeneous of degree 1 in $B_{1} \setminus B_{\t},$ and since property (${\mathcal B}{\emph 4 \, I}$) is satisfied with $w_{i}$ in place of $v$ and $z= 0,$ that it is also satisfied with $w_{\star}$ in place of $v$ and $z =0.$ Thus
if $\widetilde{w}_{\star}$ denotes the homogeneous of degree 1 extension of $\left.w_{\star}\right|_{B_{1} \setminus B_{\t}}$ to all of ${\mathbf R}^{n}$, then $\widetilde{w}_{\star} \in {\mathcal H}.$ Note also that $\{0\} \times {\mathbf R}^{n-k_{1}} \subseteq T(\widetilde{w}_{\star}).$ 

Now by homogeneity of $v$ in $B_{1} \setminus B_{\t}$, we have that for each $y \in B_{1}$, sufficiently small $\s>0$ and sufficiently large $i$, 
\begin{eqnarray*}
\s^{-n}\int_{B_{\s}(y)} w_{i}(x + z) \,dx = \e_{i}^{-1}\s^{-n}\int_{B_{\s}(y)} w(z_{i} + \r_{i}(x + z)) \, dx&&\nonumber\\
&&\hspace{-3in}=(1 + \r_{i})\e_{i}^{-1}\s^{-n}\int_{B_{\s}(y)} w(z_{i} +(1 + \r_{i})^{-1}\r_{i}(z - z_{i}) + (1+ \r_{i})^{-1}\r_{i}x) \, dx\nonumber\\
&&\hspace{-3in}=(1 + \r_{i})^{n+1} \e_{i}^{-1}\s^{-n}\int_{B_{(1+\r_{i})^{-1}\s}((1+\r_{i})^{-1}(z-z_{i}+y))}w(z_{i} + \r_{i}x)\, dx\nonumber\\
&&\hspace{-3in}= (1+\r_{i})^{n+1}\s^{-n}\int_{B_{(1+\r_{i})^{-1}\s}((1+\r_{i})^{-1}(z-z_{i}+y))} w_{i}(x) \, dx
\end{eqnarray*}
where $\e_{i} = \|w(z_{i} + \r_{i}(\cdot))\|_{L^{2}(B_{1})},$ so first letting $i \to \infty$ in this (noting that 
$z_{i} \to z$) and then letting $\s \to 0$, 
we conclude that $\widetilde{w}_{\star}(y + z) = \widetilde{w}_{\star}(y)$ for a.e. $y$. i.e. 
$z \in T(\widetilde{w}_{\star}).$ But $z \in B_{1} \setminus \left(\{0\} \times {\mathbf R}^{n-k_{1}}\right)$ (since $z \in K$), and therefore we must have  
${\rm dim} \, T(\widetilde{w}_{\star}) > n-k_{1}$.  Note on the other hand that by the definition of $k_{1}$, either $k_{1} = 2$ or (in case $k_{1} \geq 3$) 
${\mathcal H}_{k}  =\emptyset$ for all $k=2, \ldots, (k_{1}-1),$ so that, in either case, whenever ${\rm dim} \, T(\widetilde{v}) > n- k_{1}$ for some $\widetilde{v} \in {\mathcal H}$, it follows that $\widetilde{v} \in {\mathcal H}_{1}.$ Thus we have shown that 
$\widetilde{w}_{\star} \in {\mathcal H}_{1}$ and hence that
$\widetilde{w}_{\star}^{1} = \widetilde{w}_{\star}^{2} =\ldots = \widetilde{w}_{\star}^{q} = L$ for some linear function $L$. But since $(w_{i})_{a} \equiv 0$ for each $j=1, 2, \ldots,$ it follows that 
$L = (\widetilde{w}_{\star})_{a} = 0$, which is a contradiction. Thus (\ref{homog-1}) must hold 
for some $\e = \e(v, K, {\mathcal B}) \in (0, 1)$ and all $z \in K$, $\r \in (0, \e]$ as claimed.

Combining (\ref{homog-1}) with property (${\mathcal B {\emph 4 \, I}}$), we then have that 
$$\sum_{j=1}^{q}\int_{B_{\r}(z) \setminus B_{\t\r}(z)} R_{z}^{2-n} \left(\frac{\partial\, \left((v^{j} - v_{a})/R_{z}\right)}{\partial\, R_{z}}\right)^{2} \geq \frac{\e}{C} \sum_{j=1}^{q}\int_{B_{\t\r}(z)} R_{z}^{2-n} \left(\frac{\partial\, \left(v^{j} - v_{a})/R_{z}\right)}{\partial\, R_{z}}\right)^{2}$$ 
which implies that 
\begin{equation}\label{homog-3}
\sum_{j=1}^{q}\int_{B_{\t\r}(z)} R_{z}^{2-n} \left(\frac{\partial\, \left((v^{j} - v_{a})/R_{z}\right)}{\partial\, R_{z}}\right)^{2} \leq \th \sum_{j=1}^{q}\int_{B_{\r}(z)} R_{z}^{2-n} \left(\frac{\partial\, \left(v^{j} - v_{a})/R_{z}\right)}{\partial\, R_{z}}\right)^{2} 
\end{equation} 
for all $z \in K \cap \G_{v}$  and $\r \in (0,\e]$, where $\th = \th(v, K, {\mathcal B}) \in (0, 1).$ Iterating this (for fixed $z \in K \cap \G_{v}$) with $\t^{i}\r$, $i=1, 2, 3, \ldots$ in place of $\r$, we see that 
$$\sum_{j=1}^{q}\int_{B_{\t^{i}\r}(z)} R_{z}^{2-n} \left(\frac{\partial\, \left((v^{j} - v_{a})/R_{z}\right)}{\partial\, R_{z}}\right)^{2} \leq \th^{i}\sum_{j=1}^{q}\int_{B_{\r}(z)} R_{z}^{2-n} \left(\frac{\partial\, \left(v^{j} - v_{a})/R_{z}\right)}{\partial\, R_{z}}\right)^{2}$$ 
for $i=0, 1, 2, 3 \ldots$, which readily implies that 
$$\sum_{j=1}^{q}\int_{B_{\s}(z)} R_{z}^{2-n} \left(\frac{\partial\, \left((v^{j} - v_{a})/R_{z}\right)}{\partial\, R_{z}}\right)^{2} \leq \b\left(\frac{\s}{\r}\right)^{\m}\sum_{j=1}^{q}\int_{B_{\r}(z)} R_{z}^{2-n} \left(\frac{\partial\, \left(v^{j} - v_{a})/R_{z}\right)}{\partial\, R_{z}}\right)^{2}$$ 
for any $z \in K \cap \G_{v}$ and all $0 <\s \leq \r/2 \leq \e/2,$ where the constants $\b = \b(v, K, {\mathcal B}) \in (0, \infty)$ and $\m = \m(v, K, {\mathcal B}) \in (0, 1)$ are independent of $z$. By property (${\mathcal B{\emph 4 \, I}}$) and inequality 
(\ref{homog-1}), this yields the estimate 
\begin{equation}\label{homog-4}
\sum_{j=1}^{q}\s^{-n-2}\int_{B_{\s}(z)}|v^{j} -v_{a}|^{2} \leq 2^{-n-2}\e^{-1}C\b \left(\frac{\s}{\r}\right)^{\m}\r^{-n-2}\sum_{j=1}^{q}\int_{B_{\r}(z)}|v^{j} - v_{a}|^{2}
\end{equation}
for each $z \in K \cap \G_{v}$ and $0 < \s \leq  \r/2 \leq \e/4.$ Since property (${\mathcal B{\emph 4 \, II}}$) and the definition of $\G_{v}$ imply that $v$ is harmonic in ${\mathbf R}^{n} \setminus \G_{v}$,  we deduce from Lemma~\ref{general}, the remark immediately following Lemma~\ref{general} and the arbitrariness of $K$ that 
$v \in C^{1}\left((B_{1} \setminus \overline{B_{\t}}) \setminus \left(\{0\} \times {\mathbf R}^{n-k_{1}}\right)\right).$ 

Now by property (${\mathcal B{\emph 4 \, I}}$), 
$\G_{v} \cap (B_{1} \setminus \overline{B_{\t}}) \setminus \left(\{0\} \times {\mathbf R}^{n-k_{1}}\right) \subset$ the zero set of 
$u^{j} \equiv \left.\left(v^{j} - v^{j-1}\right)\right|_{B_{1} \setminus \overline{B_{\t}}}$ for each $j=2, \ldots, q.$ Since $u^{j}$ is non-negative and $C^{1}$ in 
$(B_{1} \setminus \overline{B_{\t}}) \setminus \left(\{0\} \times {\mathbf R}^{n-k_{1}}\right)$, it follows that $Du^{j}(z) = 0$ for any 
$z \in \G_{v} \cap (B_{1} \setminus \overline{B_{\t}}) \setminus \left(\{0\} \times {\mathbf R}^{n-k_{1}}\right)$. Also, by property (${\mathcal B{\emph 4\, II}}$) and the definition of $\G_{v}$, $u^{j}$ is harmonic in $(B_{1} \setminus \overline{B_{\t}}) \setminus \left(\G_{v} \cup \left(\{0\} \times {\mathbf R}^{n-k_{1}}\right)\right)$. In order to derive a contradiction, pick any point $z_{1} \in \G_{v} \cap (B_{1} \setminus \overline{B_{\t}}) \setminus \left(\{0\} \times {\mathbf R}^{n-k_{1}}\right)$ and let $\r_{1} = \frac{1}{4} \, 
{\rm dist} \, \left(z_{1}, \partial \, B_{1} \cup \partial B_{\t} \cup \{0\} \times {\mathbf R}^{n-k_{1}}\right).$ 
If $u^{j}(z) > 0$ for some $z \in B_{\r_{1}}(z_{1})$, then there exists $\r  \in (0, \r_{1})$ such that $u^{j} >0$ in $B_{\r}(z)$ and 
$\partial \, B_{\r}(z) \cap \left(\G_{v} \cap (B_{1} \setminus B_{\t}) \setminus \left(\{0\} \times {\mathbf R}^{n-k_{1}}\right)\right) \neq \emptyset$, contradicting the Hopf boundary point lemma. It follows that $u^{j} \equiv 0$ in $B_{\r_{1}}(z_{1})$ for each $j=2, \ldots, q.$ But since $z_{1} \in \G_{v}$, this is impossible by the definition of $\G_{v},$ so we see that the assumption 
$\G_{v} \cap (B_{1} \setminus \overline{B_{\t}}) \setminus \left(\{0\} \times {\mathbf R}^{n-k_{1}}\right) \neq \emptyset$ leads to a contradiction. Thus ${\mathcal H}_{k} = \emptyset$ for each $k = 2, \ldots, n$, and the Proposition is proved.
\end{proof}

%\medskip
\begin{proof}[Proof of Theorem~\ref{blowup-reg}] The main point is to prove that ${\mathcal B} \subseteq  C^{1}(B_{1}).$ For if this is true, 
then, by exactly the same argument as in the last paragraph of the proof of Proposition~\ref{homog}, we see that 
$\G_{v}  = \emptyset$ for each $v \in {\mathcal B},$ from which the first assertion of the theorem follows immediately. 

In view of Lemma~\ref{general}, property $({\mathcal B{\emph 4\, II}})$ and property $({\mathcal B{\emph 5\, I}})$, to prove that ${\mathcal B}  \subseteq C^{1}(B_{1})$, it suffices to establish that there are fixed constants $\b  = \b({\mathcal B}) \in (0, \infty)$ and 
$\m  = \m({\mathcal B}) \in (0, 1)$ such that for each $v \in {\mathcal B}$, $z \in \G_{v} \cap B_{3/4}$ and $0 < \s \leq \r/2 \leq1/8$,
\begin{equation}\label{blowup-1}
\s^{-n-2}\sum_{j=1}^{q}\int_{B_{\s}(z)} |v^{j} - \ell_{z}|^{2} \leq \b\left(\frac{\s}{\r}\right)^{\m} \r^{-n-2}
\sum_{j=1}^{q}\int_{B_{\r}(z)}|v^{j} - \ell_{z}|^{2}
\end{equation}
where $\ell_{z}$ is the affine function given by $\ell_{z}(x) = v_{a}(z) + Dv_{a}(z) \cdot (x- z)$, $x \in {\mathbf R}^{n}.$ This estimate follows by exactly the same hole-filling
argument used in the proof of Proposition~\ref{homog}. Specifically, we may first prove, by arguing by contradiction and using Proposition~\ref{homog},  that 
there exists a fixed constant $\e = \e({\mathcal B}) >0$ such that if $v \in {\mathcal B}$, $0 \in \G_{v},$ $v_{a}(0) = 0$ and 
$Dv_{a}(0) = 0$, then 
$$\sum_{j=1}^{q}\int_{B_{1/4}(0) \setminus B_{\t/4}(0)} R^{2-n} \left(\frac{\partial\, \left(v^{j}/R\right)}{\partial\, R}\right)^{2} 
\geq \e \sum_{j=1}^{q}\int_{B_{1/4}(0)}|v^{j}|^{2}$$
where $\t = \t({\mathcal B}) \in (0, 1/4)$ is the constant as in (\ref{blow-up-non-conc}). It follows from this and property (${\mathcal B{\emph 4 \, I}}$) (by arguing as in the proof of (\ref{homog-4})) that if $v \in {\mathcal B}$, $0 \in \G_{v},$ $v_{a}(0) = 0$ and $Dv_{a}(0) = 0$, then
$$\r^{-n-2}\int_{B_{\r}(0)}|v|^{2} \leq \b\r^{\m}\int_{B_{1/2}}|v|^{2} \;\;\;\;\; \forall \r \in (0, 1/2]$$
where $\b = \b({\mathcal B}) \in (0, \infty)$ and $\m  = \m({\mathcal B}) \in (0, 1)$. In view of properties (${\mathcal B{\emph 5 \, I}}$) and (${\mathcal B{\emph 5 \, III}}$), the estimate (\ref{blowup-1}) follows from this.

Since finiteness of $\sum_{j=1}^{q}\int_{B_{\r}(z)}R_{z}^{2-n}\left(\frac{\partial((v^{j} - v_{a}(z))/R_{z})}{\partial R_{z}}\right)^{2}$ implies that $v^{1}(z) = v^{2}(z) =\ldots =v^{q}(z)$ $(= v_{a}(z))$, the second assertion of the theorem follows from the first, property (${\mathcal B{\emph 2}}$) and 
the maximum principle. 
\end{proof}

\section{Lipschitz approximation and coarse blow-ups}\label{blow-up}
\setcounter{equation}{0}

We here recall (in Theorem~\ref{flat-varifolds} below) some facts concerning approximation of a stationary integral varifold weakly close to a hyperplane by the graph of a Lipschitz function over the hyperplane. These results were established by Almgren (\cite{A}), adapting, for higher multiplicity setting, the corresponding result of Allard (\cite{AW}) for multiplicity 1 varifolds. We shall use these facts to blow up mass-bounded sequences of varifolds weakly converging to a hyperplane.

First note the following elementary fact, which we shall need here and subsequently: If $V$ is a stationary integral $n$-varifold on 
$B_{2}^{n+1}(0)$, then 
\begin{equation}\label{tilt-ht-0}
\int_{B_{2}^{n+1}(0)} |\nabla^{V} \, x^{1}|^{2} \widetilde{\z}^{2} \,  d\|V\|(X) \leq 4\int_{B_{2}^{n+1}(0)}|x^{1}|^{2}|\nabla^{V} \, \widetilde{\z}|^{2} \, d\|V\|(X)
\end{equation}
for each $\widetilde{\z} \in C^{1}_{c}(B_{2}^{n+1}(0)).$ This is derived simply by taking $\psi(X) = x^{1}{\widetilde \z}^{2}(X)e^{1}$ in the first variation formula (\ref{firstvariation}). 

Let $\r \in (0, 1)$ and suppose that ${\rm spt} \, \|V\| \cap ({\mathbf R} \times B_{(1+\r)/2}) \subset \{|x^{1}| < 1\}.$ Choosing $\widetilde{\z}$ in (\ref{tilt-ht-0}) such that $\widetilde{\z}(x^{1}, x^{\prime})  = \z(x^{\prime})$ in a neighborhood of ${\rm spt} \, \|V\| \cap ({\mathbf R} \times B_{1})$, where 
$\z \in C^{1}_{c}(B_{(1+\r)/2})$ is such that $\z \equiv 1$ on $B_{\r}$, $0 \leq \z \leq 1$ and $|D\z| \leq C$ for some constant $C = C(\r)$ (e.g.\ $\widetilde{\z}(x^{1}, x^{\prime})  = \eta(x^{1})\z(x^{\prime})$ where $\eta \in C^{1}_{c}(-3/2, 3/2)$ with $\eta \equiv 1$ on $[-1, 1]$), we deduce from (\ref{tilt-ht-0}) that for each $\r \in (0, 1)$, 
\begin{equation}\label{tilt-ht}
\int_{{\mathbf R} \times B_{\r}} |\nabla^{V} \, x^{1}|^{2} d\|V\|(X) \leq C{\hat E}_{V}^{2}
\end{equation}
where $C = C(n, \r) \in (0, \infty),$ and ${\hat E}_{V} = \sqrt{\int_{{\mathbf R} \times B_{1}} |x^{1}|^{2} \, d\|V\|(X).}$

\begin{theorem}[\cite{A}, Corollary 3.11] \label{flat-varifolds}
Let $q$ be a positive integer and $\s \in (0, 1)$. There exist numbers $\e_{0} = \e_{0}(n, q, \s) \in (0, 1/2)$ and $\xi = \xi(n, q) \in (0, 1/2)$ such that the following holds: Let $V$ be a stationary integral $n$-varifold on 
$B_{2}^{n+1}(0)$ with $$(\omega_{n}2^{n})^{-1}\|V\|(B_{2}^{n+1}(0)) < q + 1/2, \;\;\;\;\;\;\; q - 1/2 \leq \omega_{n}^{-1}\|V\|({\mathbf R} \times B_{1}) < q + 1/2 \; {\rm and}$$  
$${\hat E}_{V}^{2} \equiv \int_{{\mathbf R} \times B_{1}} |x^{1}|^{2} d\|V\|(X) \leq \e_{0}; \;\; let$$
$$\Sigma = \pi \, \widetilde{\Sigma}_{1} \cup \pi \, \widetilde{\Sigma}_{2} \cup \pi \, \widetilde{\Sigma}_{3} \cup \Sigma^{\prime}$$
where $\pi \,: \, {\mathbf R}^{n+1} \to \{0\} \times {\mathbf R}^{n}$ is the orthogonal projection, 
$$\widetilde{\Sigma}_{1} = \left\{ Y \in {\rm spt} \, \|V\| \cap ({\mathbf R} \times B_{\s}) \, :\, \r^{-n}\int_{{\mathbf R} \times B_{\r}(\pi \, Y)} |\nabla^{V} \, x^{1}|^{2} \, d\|V\|(X) \geq \xi \;\; \mbox{\rm{for some}} \;\; \r \in (0, 1-\s)\right\},$$
\begin{eqnarray*}
\widetilde{\Sigma}_{2} = \left\{Y \in {\rm spt} \, \|V\| \cap ({\mathbf R} \times B_{\s}) \,: \, {\rm either} \;\; {\rm Tan} \, ({\rm spt} \, \|V\|, Y) \neq {\rm Tan}^{n} \, (\|V\|, Y)\right.&&\\
&&\hspace{-3.5in}\left.{\rm or} \;\; {\rm Tan} \, ({\rm spt} \, \|V\|, Y) \not\in G_{n} \;\; {\rm or} \;\; \Theta \, (\|V\|, Y) \; \mbox{\rm{is not a positive integer}}\right\}
\end{eqnarray*}
where ${\rm Tan} \, ({\rm spt} \, \|V\|, Y)$ denotes the tangent cone of ${\rm spt} \, \|V\|$ at $Y$ 
${\rm (}$\cite{F}, 3.1.21$\rm )$ and ${\rm Tan}^{n} \, (\|V\|, Y)$ denotes the $(\|V\|, n)$ approximate tangent cone
of $\|V\|$ at $Y$ ${\rm (}$\cite{F}, 3.2.16${\rm )}$,
$$\widetilde{\Sigma}_{3} = \left\{Y \in {\rm spt} \, \|V\| \cap ({\mathbf R} \times B_{\s}) \setminus \widetilde{\Sigma}_{2}\, : \, 1 -  (e_{1} \cdot \nu(Y))^{2}\geq 1/4\right\}$$
where $\nu(Y)$ is the unit normal to ${\rm Tan} \, ({\rm spt} \, \|V\|, Y)$, and 
$$\Sigma^{\prime} = \left\{Y \in B_{\s} \setminus (\pi \, \widetilde{\Sigma}_{1} \cup \pi \, \widetilde{\Sigma}_{2} \cup \pi \, \widetilde{\Sigma}_{3}) \, : \, \Theta \, (\|\pi_{\#} \, V\|, Y) \leq q-1\right\}.$$
Then  
\begin{itemize}
\item[(a)] ${\mathcal H}^{n} \,(\Sigma) + \|V\|({\mathbf R} \times \Sigma) \leq C{\hat E}_{V}^{2}$ where 
$C = C(n , q, \s) \in (0, \infty).$ 
\item[(b)] There are Lipschitz functions $u^{j} \; : \; B_{\s} \to {\mathbf R},$ with 
${\rm Lip} \, u^{j} \leq 1/2$ for each $j \in \{1, 2, \ldots, q\}$ such that 
$u^{1} \leq u^{2} \leq \ldots \leq u^{q}$  and
$${\rm spt} \, \|V\| \cap ({\mathbf R} \times (B_{\s} \setminus \Sigma)) = \cup_{j=1}^{q} {\rm graph} \, u^{j} \cap ({\mathbf R} \times (B_{\s} \setminus \Sigma)).$$ 
\item[(c)] For each $x \in B_{\s} \setminus \Sigma$ and each $Y \in 
{\rm spt} \, \|V\| \cap \pi^{-1}(x)$, $\Theta \, (\|V\|, Y)$ is a positive integer and  
$$\sum_{Y \in {\rm spt} \, \|V\| \cap \pi^{-1}(x)} \Theta \,(\|V\|, Y) = q.$$ 
\end{itemize}
\end{theorem}

\begin{proof}  
In view of the Constancy Theorem (\cite{S1}, Theorem 41.1), the estimate (\ref{tilt-ht}), and the easily verifiable fact that in the present codimension 1 setting, the ``unordered distance'' is the same as the ``ordered distance'' (that is, if $a_{j}, b_{j} \in {\mathbf R}$ are such that $a_{1} \leq a_{2} \leq \ldots \leq a_{q}$ and $b_{1} \leq b_{2} \leq \ldots \leq b_{q}$, then 
${\mathcal G}(\{a_{1}, \ldots, a_{q}\}, \{b_{1}, \ldots, b_{q}\}) \equiv \inf\left\{\sqrt{\sum_{j=1}^{q} (a_{j} - b_{\s(j)})^{2}} \, : \, \s \; \mbox{is a permutation  of} \; \{1, \ldots, q\}\right\} = \sqrt{\sum_{j=1}^{q} (a_{j} - b_{j})^{2}}),$ the theorem follows immediately from [\cite{A}, Corollary 3.11], which in turn is a fairly straightforward adaptation of the corresponding argument in \cite{AW} for the case $q= 1.$ 
\end{proof}

%\medskip
\noindent
{\bf Remark:} It is an easy consequence of the monotonicity of mass ratio (\cite{S1}, Section 17.5) 
that for each $\s \in (0, 1)$, there exists $\e = \e(n, \s) \in (0, 1)$ such that if $V$ is a stationary integral $n$-varifold on ${\mathbf R} \times B_{1}$ with ${\hat E}^{2}_{V}  = \int_{{\mathbf R} \times B_{1}} |x^{1}|^{2} \, d\|V\|(X) < \e$ then 
$$\sup_{X \in \left({\mathbf R} \times B_{\s}\right) \cap {\rm spt} \, \|V\|} |x^{1}| \leq C{\hat E}_{V}^{1/n}$$
where $C = C(n) \in (0, \infty).$ In particular, under the hypotheses of Theorem~\ref{flat-varifolds}, we have that 
$$\sup_{x \in B_{\s}} |u(x)|  \leq C{\hat E}_{V}^{1/n}$$
where $C = C(n) \in (0, \infty).$

\medskip

Let $q$ be a positive integer. Let $\{V_{k}\}$ be a sequence of $n$-dimensional stationary integral varifolds of $B_{2}^{n+1}(0)$ such that 
\begin{equation}\label{blow-up-1}
(\omega_{n}2^{n})^{-1}\|V_{k}\|(B_{2}^{n+1}(0)) < q + 1/2; \;\;\; q-1/2 \leq \omega_{n}^{-1}\|V_{k}\|({\mathbf R} \times B_{1}) < q + 1/2
\end{equation}
for each $k=1, 2, 3, \ldots$, and ${\hat E}_{k} \to 0,$ where
\begin{equation}\label{blow-up-1-1}
{\hat E}_{k}^{2} \equiv {\hat E}_{V_{k}}^{2} = \int_{{\mathbf R} \times B_{1}} |x^{1}|^{2}\, d\|V_{k}\|(X).
\end{equation}
Let $\s \in (0, 1).$ By Theorem~\ref{flat-varifolds}, for all sufficiently large $k$, there exist Lipschitz functions 
$u_{k}^{j} \, : \, B_{\s} 
\to {\mathbf R},$ $j=1, 2, \ldots, q$, with $u_{k}^{1} \leq u_{k}^{2} \leq \ldots \leq u_{k}^{q}$ and 
\begin{equation}\label{blow-up-2}
{\rm Lip} \, u_{k}^{j} \leq 1/2 \;\;\; \mbox{for each} \;\;\;  j \in \{1, 2, \ldots, q\}
\end{equation}
such that 
\begin{equation}\label{blow-up-4}
{\rm spt} \, \|V_{k}\| \cap ({\mathbf R} \times (B_{\s} \setminus \Sigma_{k})) = 
\cup_{j=1}^{q} {\rm graph} \, u_{k}^{j} \cap ({\mathbf R} \times (B_{\s} \setminus \Sigma_{k}))
\end{equation}
where $\Sigma_{k}$ is the measurable subset of $B_{\s}$ that corresponds to $\Sigma$ in Theorem~\ref{flat-varifolds} 
when $V$ is replaced by $V_{k};$ thus by Theorem~\ref{flat-varifolds},  
\begin{equation}\label{blow-up-3}
\|V_{k}\|({\mathbf R} \times \Sigma_{k}) + {\mathcal H}^{n} \, (\Sigma_{k}) \leq C{\hat E}_{k}^{2}
\end{equation}
where $C  = C(n, q, \s) \in (0, \infty).$ Set 
$v_{k}^{j}(x) = {\hat E}_{k}^{-1}u_{k}^{j}(x)$ for $x \in B_{\s},$ and write $v_{k} = (v_{k}^{1}, v_{k}^{2}, \ldots, v_{k}^{q}).$ 
Then  $v_{k}$ is Lipschitz on $B_{\s};$ and by (\ref{blow-up-3}) and (\ref{blow-up-4}),
\begin{equation}\label{blow-up-5}
\int_{B_{\s}} |v_{k}|^{2} \leq C, \;\;\; C = C(n, q, \s) \in (0, \infty).
\end{equation} 
Furthermore, 
\begin{eqnarray*}
\int_{B_{\s}}(1 + |Du_{k}|^{2})^{-1/2} |Du_{k}|^{2} = \int_{B_{\s} \setminus \Sigma_{k}} (1 + |Du_{k}|^{2})^{-1/2}|Du_{k}|^{2}&&\nonumber\\ 
&&\hspace{-2in} + \; \int_{B_{\s} \cap \Sigma_{k}} (1 + |Du_{k}|^{2})^{-1/2}|Du_{k}|^{2}\nonumber\\
&&\hspace{-2in} \leq \int_{{\mathbf R} \times B_{\s}}  |\nabla^{V_{k}} \, x^{1}|^{2} \, d\|V_{k}\|(X) + C_{1}{\hat E}_{k}^{2} \leq C_{2}{\hat E}_{k}^{2}
\end{eqnarray*}
where $C_{1} = C_{1}(n, q, \s) \in (0, \infty),$ $C_{2} = C_{2}(n, q, \s) \in (0, \infty)$ 
and we have used (\ref{blow-up-2}) in the first inequality and (\ref{tilt-ht}) in the second. By (\ref{blow-up-2}) again, this implies that
\begin{equation}\label{blow-up-6}
\int_{B_{\s}}|Dv_{k}|^{2} \leq C, \;\;\; C = C(n, q,\s) \in (0, \infty).
\end{equation}
In view of the arbitrariness of $\s \in (0, 1),$ by (\ref{blow-up-5}), (\ref{blow-up-6}), Rellich's theorem and a diagonal sequence argument, we obtain a function $v \in W^{1, 2}_{\rm loc} \, (B_{1}; {\mathbf R}^{q}) \cap L^{2} \, (B_{1}; {\mathbf R}^{q})$ and a subsequence $\{k_{j}\}$ of $\{k\}$ such that 
$v_{k_{j}} \to v$ as $j \to \infty$ in $L^{2} \, (B_{\s}; {\mathbf R}^{q})$ and weakly in $W^{1, 2} \, (B_{\s}; {\mathbf R}^{q})$ for every $\s \in (0, 1).$ 

\noindent
{\bf Definitions}: {\bf (1) Coarse blow-ups}:  Let $v \in W^{1,2}_{\rm loc} \, (B_{1}; {\mathbf R}^{q}) \cap L^{2} \, (B_{1}; {\mathbf R}^{q})$ correspond, in the manner described above, to (a subsequence of) a sequence $\{V_{k}\}$ of stationary integral $n$-varifolds of $B_{2}^{n+1}(0)$ satisfying (\ref{blow-up-1}) and with ${\hat E}_{k} \to 0$, where  ${\hat E}_{k}$ is as in (\ref{blow-up-1-1}). We shall call $v$ a {\em coarse blow-up} of the sequence $\{V_{k}\}$.

\noindent
{\bf(2)  The Class ${\mathcal B}_{q}$}: Denote by ${\mathcal B}_{q}$ the collection of all coarse blow-ups of sequences 
of varifolds $\{V_{k}\} \subset {\mathcal S}_{\a}$ satisfying  (\ref{blow-up-1}) and for which 
${\hat E}_{k} \to 0,$ where ${\hat E}_{k}$ is as in (\ref{blow-up-1-1}). 

\section{An outline of the proof of the main theorems}\label{outline}
\setcounter{equation}{0}

Note that if ${\mathbf C}_{0}$ is a stationary cone as in Theorem~\ref{no-transverse}, then 
$\Theta_{{\mathbf C}_{0}}(0) = q-1/2$ or $\Theta_{{\mathbf C}_{0}}(0) = q$ for some integer $q \geq 2.$ We prove both 
Theorem~$3.3^{\prime}$ and Theorem~\ref{no-transverse} simultaneously by induction on $q$. The case $q=1$ of Theorem~$3.3^{\prime}$ is a consequence of Allard's Regularity Theorem. (Note however that setting $q=1$ in the proofs of Lemma~\ref{excess-s} and Theorem~\ref{sheetingthm} 
given below reproduces Allard's argument proving Theorem~$3.3^{\prime}$ in case $q=1$). Validity of the cases 
$\Theta \, (\|{\mathbf C}_{0}\|, 0) = 3/2$ and $\Theta \, (\|{\mathbf C}_{0}\|, 0) = 2$ of Theorem~\ref{no-transverse} will be justified at the end of Section~\ref{MD}. 

Let $q$ be an integer $\geq 2$ and consider the following: 

\noindent
{\sc Induction Hypotheses:}
\begin{itemize}
\item[($H1$)] Theorem~${3.2}^{\prime}$ holds with $1, \ldots, (q-1)$ in place of 
$q.$
\item[($H2$)] Theorem~\ref{no-transverse} holds whenever $\Theta_{{\mathbf C}_{0}}(0) \in \{3/2, 2, 5/2, \ldots, q\}$.
\end{itemize}

\noindent
The inductive proof of Theorem~$3.3^{\prime}$ and Theorem~\ref{no-transverse} is obtained by completing, 
assuming $(H1)$, $(H2)$, the steps below in the order they are listed:
\begin{itemize}
\item[]{\sc Step 1:} Prove that ${\mathcal B}_{q}$ is a proper blow-up class. (Sections~\ref{non-con-tilt}-\ref{propertiesIII})
\item[]{\sc Step 2:} Prove Theorem $3.3^{\prime}$. (Section~\ref{sheeting})
\item[]{\sc Step 3:} Prove Theorem~\ref{no-transverse} when $\Theta \, (\|{\mathbf C}_{0}\|, 0) = q +1/2.$ (Section~\ref{MD})
\item[]{\sc Step 4:} Prove Theorem~\ref{no-transverse} when $\Theta \, (\|{\mathbf C}_{0}\|, 0) = q+1.$ (Section~\ref{MD})
\end{itemize}

\noindent
{\bf Remarks:} {\bf (1)} Let $m \in \{1, 2, \ldots, n\}$ and suppose that ${\mathbf C}$ is an $m$-dimensional stationary integral cone in ${\mathbf R}^{n+1}.$ Let  
$L_{\mathbf C} = \{Y \in {\rm spt} \, \|{\mathbf C}\| \, : \, \Theta \, (\|{\mathbf C}\|, Y) = \Theta \, (\|{\mathbf C}\|, 0)\}.$ It is  a well known consequence of the monotonicity formula that $L_{\mathbf C}$ is a linear subspace of ${\mathbf R}^{n+1}$ of  dimension $\leq m$ and that $Y \in L_{\mathbf C}$ if and only if $T_{Y \, \#} \, {\mathbf C} = {\mathbf C}$, where $T_{Y}\, : \, {\mathbf R}^{n+1} \to {\mathbf R}^{n+1}$ is the translation $T_{Y}(X) = X - Y.$ Let $d_{\mathbf C} = {\rm dim} \, L_{\mathbf C}.$ Then, if $\G_{\mathbf C}$ is a rotation of ${\mathbf R}^{n+1}$ such that $\G_{\mathbf C}(L_{\mathbf C}) = \{0\} \times {\mathbf R}^{d_{\mathbf C}}$, we have that 
$\G_{{\mathbf C} \, \#} \, {\mathbf C} = {\mathbf C}^{\prime} \times {\mathbf R}^{d_{\mathbf C}}$, where ${\mathbf C}^{\prime}$ is a stationary integral cone in ${\mathbf R}^{n+1 - d_{\mathbf C}}.$ Here, given an integer $d \in \{0, 1, 2, \ldots, n\}$ and a rectifiable varifold $V^{\prime}$ of  ${\mathbf R}^{n+1 -d}$, we use the notation $V^{\prime} \times {\mathbf R}^{d}$ to denote the rectifiable varifold $V$ of ${\mathbf R}^{n+1}$ with ${\rm spt} \, \|V\| = {\rm spt} \, \|V^{\prime}\| \times {\mathbf R}^{d}$ and the multiplicity function $\th_{V}$  defined by $\th_{V}(x, y) = \th_{V^{\prime}}(x)$ for $(x, y) \in {\rm spt} \, \|V^{\prime}\| \times {\mathbf R}^{d}$, where $\th_{V^{\prime}}$ is the multiplicity function of $V^{\prime}.$

\noindent
{\bf (2)} \emph{Let $q$ be an integer $\geq 2$ and suppose that the induction hypotheses $(H1)$, $(H2)$ hold. Let $V \in {\mathcal S}_{\a}.$ Then we have the following:
\begin{itemize}
\item[(a)] If $2 \leq n \leq 6$, then ${\rm sing} \, V \cap \{Z \in {\rm spt}  \, \|V\|\, : \, \Th \, (\|V\|, Z) < q\} = \emptyset.$
\item[(b)] If $n \geq 7,$ $Z \in {\rm sing} \, V$ and $\Theta \, (\|V\|, Z) < q$, then $d_{\mathbf C} \leq n-7$ for any 
${\mathbf C} \in \mbox{Var Tan} \, (V, Z).$ 
 \end{itemize}}
 
To see this, suppose either (a) or (b) is false. Then we have either 
\begin{itemize}
\item[(a$^{\prime}$)] $n \in \{2, 3, \ldots, 6\}$ and there exist a varifold $V \in {\mathcal S}_{\a}$ and a point $Z \in {\rm sing} \, V$ such that  $\Th \, (\|V\|, Z) < q$  or
\item[(b$^{\prime}$)] $n \geq 7$ and there exist a varifold $V \in {\mathcal S}_{\a}$ and a point $Z \in {\rm sing} \, V$ with $\Theta \, (\|V\|,Z) < q$ such that 
$d_{\mathbf C} > n-7$ for some ${\mathbf C} \in \mbox{Var Tan} \, (V, Z).$ 
\end{itemize}

If (a$^{\prime}$) holds, fix any ${\mathbf C} \in \mbox{Var Tan} \, (V, Z).$ 

In either case (a$^{\prime}$) or (b$^{\prime}$), the induction hypothesis $(H1)$ implies that $d_{\mathbf C} \neq n;$ for if $d_{\mathbf C} = n$, then ${\mathbf C} = q^{\prime}|P|$ for some integer 
$q^{\prime} \in \{1, 2, \ldots, q-1\}$ and some hyperplane $P$ which we may take without loss of generality to be $\{0\} \times {\mathbf R}^{n},$ whence by the definition of tangent cone and the fact that weak convergence of stationary integral varifolds implies convergence of mass and convergence in Hausdorff distance of the supports of the associated weight measures, for any given $\e >0$, there exists  $\s \in (0, 1 - |Z|/2)$ such that ${\rm dist}_{\mathcal H} \, ({\rm spt} \, \|\eta_{Z, \s \, \#} \, V\| \cap ({\mathbf R} \times B_{1}), \{0\} \times B_{1}) <\e,$ $q^{\prime} - 1/2 \leq \omega_{n}^{-1}\|\eta_{Z, \s \, \#} \, V\|({\mathbf R} \times B_{1}) < q^{\prime} + 1/2$ and $(\omega_{n}2^{n})^{-1}\|\eta_{Z, \s \, \#} \, V\|(B_{2}^{n+1}(0)) < q^{\prime} + 1/2.$ Choosing 
$\e = \e_{0}(n, \a, q^{\prime})$ where $\e_{0}$ is as in Theorem~$3.3^{\prime}$, by $(H1),$ we may apply Theorem~$3.3^{\prime}$ to deduce that near $Z$, $V$ corresponds to an embedded graph of a
$C^{1, \a}$ function over $P$ solving the minimal surface equation, and hence ${\rm spt} \, \|V\|$ near $Z$ is an embedded analytic hypersurface, contradicting our assumption that $Z \in {\rm sing} \,V.$ Thus $d_{\mathbf C} < n.$

Again in either case, the induction hypothesis $(H2)$ implies that $d_{\mathbf C} \neq n-1;$ for if $d_{\mathbf C} = n-1$, then ${\rm spt} \, \|{\mathbf C}\|$ is the union of at least three 
half-hyperplanes meeting along an $(n-1)$-dimensional subspace, and since $\Theta \, (\|{\mathbf C}\|, 0) < q$, we must have 
that $\Theta \, (\|{\mathbf C}\|, 0) \in \{3/2, 2, 5/2, \ldots, q-1/2\}.$ Again by the definition of tangent cone we have that for any given $\e_{1}>0,$ a number $\s \in (0, 1 - |Z|/2)$ such that 
$|(\omega_{n}2^{n})^{-1}\|\eta_{Z, \s \, \#} \, V\|(B_{2}^{n+1}(0)) - \Theta \, (\|{\mathbf C}\|, 0)| < 1/8$ and 
${\rm dist}_{\mathcal H} \, ({\rm spt} \, \|\eta_{Z, \s \, \#} \, V\| \cap B_{1}^{n+1}(0), {\rm spt} \, \|{\mathbf C}\| \cap B_{1}^{n+1}(0)) < \e_{1},$ so choosing $\e_{1} = \frac{1}{2}\e(\a, \frac{1}{8}, {\mathbf C})$ where $\e$ is as in Theorem~\ref{no-transverse}, we 
see by hypothesis $(H2)$ we have a contradiction to Theorem~\ref{no-transverse}. 

Thus $d_{\mathbf C} \leq n-2.$ Assume now without loss of generality that 
$L_{\mathbf C} = \{0\} \times {\mathbf R}^{d_{\mathbf C}}.$ Then ${\mathbf C} = {\mathbf C}^{\prime} \times {\mathbf R}^{d_{\mathbf C}}$, where ${\mathbf C}^{\prime}$ is an  
$(n - d_{\mathbf C})$-dimensional stationary integral cone of ${\mathbf R}^{n+1 - d_{\mathbf C}}$ with $0 \in {\rm sing} \, {\mathbf C}^{\prime}.$ Note that since 
$\Theta \, (\|{\mathbf C}\|, Y) < q$ for each $Y \in {\rm spt} \, \|{\mathbf C}\|$, in view of hypothesis $(H1)$, it follows from 
Theorem~$3.3^{\prime}$ that ${\rm reg} \, {\mathbf C}$ satisfies the stability inequality; viz., 
$\int_{{\rm reg} \, {\mathbf C}} |A_{\mathbf C}|^{2}\z^{2} \leq \int_{{\rm reg} \, {\mathbf C}} |\nabla^{{\mathbf C}} \, \z|^{2}$
for each $\z \in C^{1}_{c}({\rm reg} \, {\mathbf C}),$ where $A_{\mathbf C}$ denotes the second fundamental form of 
${\rm reg} \, {\mathbf C}.$  

Now by a theorem of J. Simons \cite{SJ} (see \cite{S1}, Appendix B for a shorter proof), we know that if $2 \leq n \leq 6$, there does not exist, in ${\mathbf R}^{n+1},$ a minimal hypercone with an isolated singularity and satisfying the stability inequality. Applying this to ${\mathbf C}^{\prime}$, we conclude that if 
${\rm sing} \, {\mathbf C} = \{0\} \times {\mathbf R}^{d_{\mathbf C}}$, then, in either of the cases (a$^{\prime}$) or (b$^{\prime}$), we have a contradiction. Hence there is a point $Z_{1} \in {\rm sing} \, {\mathbf C} \setminus \{0\} \times {\mathbf R}^{d_{\mathbf C}}.$ 

Let ${\mathbf C}_{1} \in \mbox{Var Tan} \, ({\mathbf C}, Z_{1}).$ Then $\{tZ_{1} \, : \, t \in {\mathbf R}\} \times {\mathbf R}^{d_{\mathbf C}} \subseteq L_{{\mathbf C}_{1}}$ so that $d_{{\mathbf C}_{1}} \geq d_{\mathbf C} + 1.$  Since ${\mathbf C}_{1} \, \res \, B_{1}^{n+1}(0) = \lim_{k \to \infty} \, V_{k}$ for some sequence of varifolds $\{V_{k}\} \subset {\mathcal S}_{\a}$ (indeed, $V_{k} = \eta_{\widetilde{Z}_{k}, \s_{k} \, \#} \, V$ for some sequence of points $\widetilde{Z}_{k}$ and a sequence of positive numbers $\s_{k}$ converging to $0$) and $\Theta \, (\|{\mathbf C}_{1}\|, 0) = \Theta \, (\|{\mathbf C}\|, Z_{1}) < q,$ by reasoning as above, we see that $d_{{\mathbf C}_{1}} \leq n-2$ and that ${\rm reg} \, {\mathbf C}_{1}$ satisfies the stability inequality. Thus $d_{\mathbf C} \leq n-3,$ and hence in particular $n \geq 3$. 

By Simons' theorem again, there exists a point $Z_{2} \in {\rm sing} \, {\mathbf C}_{1} \setminus L_{{\mathbf C}_{1}},$ which implies (by reasoning as above considering a cone 
${\mathbf C}_{2} \in \mbox{Var Tan} \, ({\mathbf C}_{1}, Z_{2})$) that $d_{\mathbf C} \leq n-4$ and $n \geq 4.$ Repeating this argument twice more in case (a$^{\prime}$), we produce a cone contradicting  Simons' theorem, and three times more in case (b$^{\prime}$), we reach 
the conclusion $d_{\mathbf C} \leq n-7$ contrary to the assumption. Thus both claims (a) and (b) must hold.  
 
\noindent
{\bf (3)} By Remark (2) above and, in case $n\geq 7$, Almgren's generalized stratification of stationary integral varifolds (\cite{A}, p. 224, Theorem 2.26 and Remark 2.28; see \cite{S3}, Section 3.4 for a concise presentation of the argument in the context of energy minimizing maps), we have the following:

\emph{Let $q$ be an integer $\geq 2$. If the induction hypotheses $(H1)$, $(H2)$ hold, $V \in {\mathcal S}_{\a},$ $\Omega \subseteq B_{2}^{n+1}(0)$ is open and $\Theta \, (\|V\|, Z) < q$ for each $Z \in {\rm spt} \, \|V\| \cap \Omega$, then 
${\mathcal H}^{n-7 + \g} \, ({\rm sing} \, V \, \res \, \Omega) = 0$ for each $\g > 0$ if $n \geq 7$ (with 
${\rm sing} \, V \, \res \, \Omega$ discrete if $n=7$) and 
${\rm sing} \, V \, \res \, \Omega = \emptyset$ if $2 \leq n \leq 6.$}\\

We shall now begin, and end in Section~\ref{propertiesIII}, the central part of our work, namely, the proof that  \emph{for any integer $q \geq 1$, the class of functions ${\mathcal B}_{q}$ (as defined at the  end of Section~\ref{blow-up})  
is a proper blow-up class (as defined in Section~\ref{proper-blow-up})}.

\section{Non-concentration of tilt-excess}\label{non-con-tilt}
\setcounter{equation}{0}

The main result of this section is the estimates of Theorem~\ref{non-concentration}(b), which says that for a stationary integral $n$-varifold on an open ball in ${\mathbf R}^{n+1}$ having small height excess relative to a hyperplane, concentration of points of ``top density'' near an $(n-1)$-dimensional subspace $L$ implies non-concentration, near $L$, of the tilt-excess of the varifold relative to the hyperplane. This estimate will play a crucial role in the proof that ${\mathcal B}_{q}$  (see definition at the end of Section~\ref{blow-up}) is a proper blow-up class---specifically, in establishing property (${\mathcal B{\emph 7}}$) (see Section~\ref{proper-blow-up}) for ${\mathcal B}_{q}$. No stability hypothesis is required for the results of this section.
 
\begin{theorem}\label{non-concentration}
Let $q$ be a positive integer, $\t \in (0, 1/16)$ and $\m \in (0, 1).$There exists a number $\e_{1} = \e_{1}(n, q, \t, \m)~\in~(0,1/2)$ such that if $V$ is a stationary integral $n$-varifold of $B_{2}^{n+1}(0)$ with 
$$(\omega_{n}2^{n})^{-1}\|V\|(B_{2}^{n+1}(0)) < q + 1/2, \;\;\;  q-1/2 \leq \omega_{n}^{-1}\|V\|({\mathbf R} \times B_{1}) < q + 1/2 \;\; {\rm and}$$  
$$\int_{{\mathbf R} \times B_{1}} |x^{1}|^{2}d\|V\|(X)\leq \e_{1},$$
then the following hold:
\begin{itemize}
\item[(a)] For each point $Z = (z^{1}, z^{\prime}) \in {\rm spt} \, \|V\| \cap({\mathbf R} \times B_{9/16})$ with $\Theta \, (\|V\|, Z) \geq q$, $$|z^{1}|^{2} \leq C\int_{{\mathbf R} \times B_{1}} |x^{1}|^{2} \, d\|V\|(X)$$ where $C = C(n, q) \in (0, \infty).$
\item[(b)] If $L$ is an $(n-1)$-dimensional subspace of $\{0\} \times {\mathbf R}^{n}$ such that 
$$L \cap B_{1/2} \subset \left(\{Z  \in {\rm spt} \, \|V\|\, : \, \Theta \, (\|V\|, Z) \geq q\}\right)_{\t}, \;\; {\rm then}$$ 
$$\int_{(L)_{\t} \cap ({\mathbf R} \times B_{1/2})} |\nabla^{V} \, x^{1}|^{2} d\|V\|(X) \leq C\t^{1-\m}\int_{{\mathbf R} \times B_{1}} |x^{1}|^{2}d\|V\|(X)$$
\noindent
where  $C = C(n, q, \m) \in (0, \infty)$. Here for a subset $A$ of 
${\mathbf R}^{n+1}$, we use the notation $(A)_{\t} = \{X \in {\mathbf R}^{n+1} \, : \, {\rm dist} \, (X, A) \leq \t\}.$
\end{itemize} 
\end{theorem}

\noindent
{\bf Remarks:} {\bf (1)} Since Theorem~\ref{flat-varifolds} holds with tilt-excess $\int_{{\mathbf R} \times B_{1}}|\nabla^{V} \, x^{1}|^{2} \, d\|V\|(X)$ in place of the height excess ${\hat E}_{V}^{2}$ (see \cite{A}, Corollary 3.11), an examination of the proof below in fact shows that for any $\m \in (0, 1)$, the more refined estimate 
$$\int_{(L)_{\t} \cap ({\mathbf R} \times B_{1/2})} |\nabla^{V} \, x^{1}|^{2} d\|V\|(X) \leq C\t^{1-\m}\int_{{\mathbf R} \times B_{1}} |\nabla^{V} \, x^{1}|^{2}d\|V\|(X), \;\; C = C(n, q, \m) \in (0, \infty),$$
holds under the hypotheses of Theorem~\ref{non-concentration}(b). We do not however need it here. 

\noindent
{\bf (2)} A similar estimate for height excess relative to certain minimal cones was established in a ``multiplicity 1 setting'' in \cite{S}. Indeed, we shall later need a version of that as well (see Corollaries~\ref{no-gaps} and ~\ref{MD-no-gaps}).

\begin{proof}
The proof is based on the monotonicity formula [\cite{S1}, 17.5], which implies that, for any $Z \in {\rm spt} \, \|V\| \cap ({\mathbf R} \times B_{9/16}),$
\begin{equation}\label{monotonicity-1}
\frac{1}{\omega_{n}}\int_{B^{n+1}_{3/8}(Z)} \frac{|(X - Z)^{\perp}|^{2}}{|X - Z|^{n+2}} d\|V\|(X) = 
\frac{\|V\|(B_{3/8}^{n+1}(Z))}{\omega_{n}(3/8)^{n}} - \Theta \, (\|V\|, Z).
\end{equation}

Write ${\hat E}_{V} = \sqrt{\int_{{\mathbf R} \times B_{1}} |x^{1}|^{2} \, d\|V\|(X)}.$ Assuming $\e_{1} = \e_{1}(n, q) \in (0, \infty)$ is sufficiently small to guarantee the validity of its conclusions, Theorem~\ref{flat-varifolds} with $\s = 15/16$ implies  that 
\begin{eqnarray*}
\|V\|(B_{3/8}^{n+1}(Z)) \leq \|V\|({\mathbf R} \times B_{3/8}(z^{\prime})) = \|V\|({\mathbf R} \times (B_{3/8}(z^{\prime}) \setminus \Sigma)) + \|V\|({\mathbf R} \times (B_{3/8}(z^{\prime}) \cap \Sigma))&&\nonumber\\
&&\hspace{-6in}\leq \sum_{j=1}^{q}\int_{B_{3/8}(z^{\prime}) \setminus \Sigma} \sqrt{1 + |Du^{j}|^{2}} d{\mathcal H}^{n} 
+ \|V\|({\mathbf R} \times \Sigma) \leq \sum_{j=1}^{q}\int_{B_{3/8}(z^{\prime})} \sqrt{1 + |Du^{j}|^{2}} d{\mathcal H}^{n} + C{\hat E}_{V}^{2}
\end{eqnarray*}
where $C = C(n, q) \in (0, \infty),$ and $u^{j}$, $j=1, 2, \ldots, q$, $\Sigma$ are as in Theorem~\ref{flat-varifolds}; if, additionally, $\Theta \, (\|V\|, Z) \geq q,$ it follows that  
\begin{eqnarray}\label{monotonicity-3}
\frac{\|V\|(B_{3/8}^{n+1}(Z))}{\omega_{n}(3/8)^{n}} - \Theta \, (\|V\|, Z) \leq \frac{\|V\|(B_{3/8}^{n+1}(Z))}{\omega_{n}(3/8)^{n}} - q&&\nonumber\\
&&\hspace{-3.5in}\leq \sum_{j=1}^{q}\frac{1}{\omega_{n}(3/8)^{n}}
\int_{B_{3/8}(z^{\prime})} \left(\sqrt{1 + |Du^{j}|^{2}} - 1\right) d{\mathcal H}^{n} + C{\hat E}_{V}^{2} \leq C\sum_{j=1}^{q}\int_{B_{3/8}(z^{\prime})} |Du^{j}|^{2} \, d{\mathcal H}^{n} + C{\hat E}_{V}^{2}\nonumber\\
&&\hspace{-3.5in} \leq C \, \sum_{j=1}^{q}\int_{B_{3/8}(z^{\prime}) \setminus \Sigma} |Du^{j}|^{2} \, d{\mathcal H}^{n} + 
C\sum_{j=1}^{q}\int_{B_{3/8}(z^{\prime}) \cap \Sigma} |Du^{j}|^{2} \, d{\mathcal H}^{n} + C{\hat E}_{V}^{2}\nonumber\\
&&\hspace{-3in} \leq C\, \sum_{j=1}^{q}\int_{B_{3/8}(z^{\prime}) \setminus \Sigma} |Du^{j}|^{2} \, d{\mathcal H}^{n} + C{\hat E}_{V}^{2} \leq C \, \int_{{\mathbf R} \times B_{3/8}(z^{\prime})} |\nabla^{V} \,x^{1}|^{2} \, d\|V\|(X) + C{\hat E}_{V}^{2}
\leq C{\hat E}_{V}^{2}
\end{eqnarray}
where $C = C(n, q) \in (0, \infty),$ and in the last inequality we have used $(\ref{tilt-ht}).$ Thus we deduce from 
(\ref{monotonicity-1}) that 
\begin{equation}\label{monotonicity-2}
\int_{B^{n+1}_{3/8}(Z)} \frac{|(X - Z)^{\perp}|^{2}}{|X - Z|^{n+2}} d\|V\|(X) \leq C{\hat E}_{V}^{2}
\end{equation}
for each $Z \in {\rm spt} \, \|V\| \cap ({\mathbf R} \times B_{9/16})$ with $\Theta \, (\|V\|, Z) \geq q,$ where $C = C(n, q) \in (0, \infty).$ 

To prove the assertion of part (a) of the Theorem, we estimate the left hand side of (\ref{monotonicity-1}) from below as follows:
\begin{eqnarray}\label{monotonicity-4}
\int_{B_{1/4}^{n+1}(Z)} \frac{|(X - Z)^{\perp}|^{2}}{|X - Z|^{n+2}} d\|V\|(X) \geq 4^{n+2}\int_{B_{1/4}^{n+1}(Z)} \left|\sum_{j=2}^{n+1}((x^{j} - z^{j})e_{j}^{\perp} + 
(x^{1} - z^{1})e_{1}^{\perp}\right|^{2} d\|V\|(X)\nonumber\\
&&\hspace{-5.5in}\geq \frac{1}{2}4^{n+2}\int_{B_{1/4}^{n+1}(Z)} |x^{1} - z^{1}|^{2}|e_{1}^{\perp}|^{2} d\|V\|(X) - 4^{n}\int_{B_{1/4}^{n+1}(Z)} \sum_{j=2}^{n+1}|e_{j}^{\perp}|^{2} d\|V\|(X)\nonumber\\
&&\hspace{-5.5in}=   \frac{1}{2}4^{n+2}\int_{B_{1/4}^{n+1}(Z)} |x^{1} - z^{1}|^{2}|e_{1}^{\perp}|^{2} d\|V\|(X) - 4^{n} \int_{B_{1/4}^{n+1}(Z)} |\nabla^{V} \, x^{1}|^{2}d\|V\|(X)\nonumber\\
&&\hspace{-5.5in}\geq  \frac{1}{2}4^{n+2}\int_{B_{1/4}^{n+1}(Z)} |x^{1} - z^{1}|^{2}|e_{1}^{\perp}|^{2}d\|V\|(X) - C{\hat E}_{V}^{2}\nonumber\\
&&\hspace{-5.5in}\geq4^{n+1}|z^{1}|^{2}\int_{B_{1/4}^{n+1}(Z)}|e_{1}^{\perp}|^{2}d\|V\|(X) - C{\hat E}_{V}^{2}\nonumber\\
&&\hspace{-5.5in}\geq 4^{n+1}|z^{1}|^{2}\sum_{j=1}^{q}\int_{B_{1/8}(z^{\prime}) \setminus \Sigma}(1 + |Du^{j}|^{2})^{-1}d{\mathcal H}^{n} - C{\hat E}_{V}^{2} \geq C|z^{1}|^{2}  - C{\hat E}_{V}^{2}
\end{eqnarray}
where, for $\|V\|$ a.e. $X \in {\rm spt} \, \|V\|$, $e_{j}^{\perp}(X)$ is the orthogonal projection of $e_{j}$ onto the orthogonal complement of the approximate tangent plane ${\rm Tan} \, (\|V\|, X)$ and $C = C(n, q) \in (0, \infty).$ Note that we have used the fact that $|Du^{j}| \leq 1/2$ a.e. and 
${\mathcal H}^{n}(B_{1/8}(z^{\prime}) \setminus \Sigma) \geq \frac{1}{2}{\mathcal H}^{n}(B_{1/8}(z^{\prime})) = 
\frac{1}{2}\omega_{n}(\frac{1}{8})^{n},$ which hold by Theorem~\ref{flat-varifolds} provided $\e_{1} = \e_{1}(n, q) \in (0, 1/2)$ is sufficiently small. The estimate of (a) readily follows from this and (\ref{monotonicity-2}).

To see (b), let $Z = (z^{1},z^{\prime}) \in {\rm spt} \, \|V\| \cap ({\mathbf R} \times B_{9/16})$ be an arbitrary point and choose $\z \in C^{1}_{c}({\mathbf R}^{n+1})$ such that
$\z \equiv 1$ on $B_{3/8}^{n+1}(0)$, $\z \equiv 0$ in ${\mathbf R}^{n+1} \setminus B_{1/2}^{n+1}(0),$ $0 \leq \z \leq 1$ and $|D\z| \leq 16$ everywhere. Taking, for $\m \in (0, 1)$, $\psi(X) = \z^{2}(X - Z)|X - Z|^{-n-2+\m}|x^{1} - z^{1}|^{2}(X-Z)$ in the first variation formula (\ref{firstvariation})  (a valid choice as shown by an easy cut-off function argument) and computing and estimating as in \cite{S}, p. 616, we deduce that 
\begin{eqnarray*}
\int_{B_{3/8}^{n+1}(Z)} \frac{|x^{1} - z^{1}|^{2}}{|X - Z|^{n+2-\m}}d\|V\|(X)\nonumber\\ 
&&\hspace{-2in}\leq C\int \left(\z^{2}(X - Z)
\frac{|(X - Z)^{\perp}|^{2}}{|X - Z|^{n+2-\m}} + \frac{|x^{1} - z^{1}|^{2}}{|X - Z|^{n - \m}} |\nabla^{V} \, \z(X - Z)|^{2}\right) d\|V\|(X)
\end{eqnarray*}
where $C = C(n, \m) \in (0, \infty).$ Since ${\rm spt} \, D\z \subset B^{n+1}_{1/2}(0) \setminus B^{n+1}_{3/8}(0)$, this together with
(\ref{monotonicity-2}) and part (a)  implies that 
\begin{equation*}
\int_{B_{3/8}^{n+1}(Z)} \frac{|x^{1} - z^{1}|^{2}}{|X - Z|^{n+2 - \m}} d\|V\|(X) \leq C\int_{{\mathbf R} \times B_{1}} 
|x^{1}|^{2}d\|V\|(X)
\end{equation*}
for every $Z = (z^{1}, z^{\prime}) \in {\rm spt} \, \|V\| \cap ({\mathbf R} \times B_{9/16})$ with 
$\Theta \, (\|V\|, Z) \geq q,$ where $C = C(n, q, \m) \in (0, \infty); $ in particular, 
\begin{equation*}
\int_{B_{4\t}^{n+1}(Z)} |x^{1} -z^{1} |^{2} d\|V\|(X)
\leq C\t^{n+2 - \m}\int_{{\mathbf R} \times B_{1}} |x^{1}|^{2} d\|V\|(X)
\end{equation*}
for each $Z = (z^{1},z^{\prime}) \in {\rm spt} \, \|V\| \cap ({\mathbf R} \times B_{9/16})$ with $\Theta \, (\|V\|, Z) \geq q$ and each $\t \in (0, 1/16).$ In view of the hypothesis 
$$L \cap B_{1/2} \subset \left(\{Z  \in {\rm spt} \, \|V\|\, : \, 
\Theta \, (\|V\|, Z) \geq q\}\right)_{\t},$$
the preceding estimate implies that for each $Y \in L \cap B_{1/2}$, there exists $z^{1} \in {\mathbf R}$ such that
$$\int_{B_{2\t}^{n+1}(Y)} |x^{1} - z^{1}|^{2}d\|V\|(X)\leq C\t^{n+2 -\m}\int_{{\mathbf R} \times B_{1}} |x^{1}|^{2} d\|V\|(X).$$
This in turn implies by (\ref{tilt-ht-0}) (applied with $\eta_{(z^{1}, 0), 1 \, \#} V$ in place of $V$ and a choice of appropriate test function $\widetilde\z$) that for each $Y \in L \cap B_{1/2},$
$$\int_{B_{3\t/2}^{n+1}(Y)} |\nabla^{V} \,x^{1}|^{2} d\|V\|(X) \leq C\t^{n-\m}\int_{{\mathbf R} \times B_{1}} |x^{1}|^{2} d\|V\|(X).$$
Since we may cover the set $(L)_{\t} \cap ({\mathbf R} \times B_{1/2})$ by $N$ balls 
$B_{3\t/2}^{n+1}(Y_{j})$ with $Y_{j} \in L \cap B_{1/2}$ for $j=1, 2, \ldots, N$ 
and with $N \leq C\t^{1-n}$, $C = C(n)$, it follows that 
$$\int_{(L)_{\t} \cap ({\mathbf R} \times B_{1/2})}|\nabla^{V} \,x^{1}|^{2}
d\|V\|(X) \leq C \t^{1-\m}\int_{{\mathbf R} \times B_{1}} |x^{1}|^{2} d\|V\|(X)$$
with $C = C(n , q, \m) \in (0, \infty)$, as required.
\end{proof}

\section{{Properties of coarse blow-ups: Part I}}\label{step2}
\setcounter{equation}{0}

Recall from Section~\ref{proper-blow-up} the defining properties $({\mathcal B{\emph 1}})-({\mathcal B{\emph 7}})$ of a proper blow-up class ${\mathcal B}$, and note that it follows from the discussion in Section~\ref{blow-up} that the class ${\mathcal B} = {\mathcal B}_{q}$ satisfies properties $({\mathcal B{\emph 1}})$ and $({\mathcal B{\emph 2}}).$ In this section, we verify that ${\mathcal B}_{q}$ also satisfies properties $({\mathcal B{\emph 3}})-({\mathcal B{\emph 6}}).$ 

Let $v \in {\mathcal B}_{q}$ be arbitrary. By the definition of ${\mathcal B}_{q}$, there exists, for each $k = 1, 2, 3, \ldots$,  a stationary integral varifold $V_{k} \in {\mathcal S}_{\a}$ such that the following are true: $(\omega_{n}2^{n})^{-1}\|V_{k}\|(B_{2}^{n+1}(0)) < q + 1/2$; $q-1/2 \leq \omega_{n}^{-1} \|V_{k}\|({\mathbf R} \times B_{1}) < q + 1/2;$ ${\hat E}_{k}^{2} \equiv 
\int_{{\mathbf R} \times B_{1}} |x^{1}|^{2}d\|V_{k}\|(X) \to 0$ as $k \to \infty;$ for each 
$\s \in (0, 1)$ and each sufficiently large $k$ depending on $\s$,  if $u_{k}^{j} \, : \, B_{\s} 
\to {\mathbf R}$ are the functions corresponding to $u^{j}$, $j= 1, 2, \ldots, q,$ and  $\Sigma_{k} \subset B_{\s}$ is the measurable set corresponding to $\Sigma$ in Theorem~\ref{flat-varifolds} taken with $V_{k}$ in place of $V,$
then, $u_{k}^{1} \leq u_{k}^{2} \leq \ldots \leq u_{k}^{q};$ $u_{k}^{j}$ is Lipschitz with 
\begin{equation}\label{step2-0}
{\rm Lip} \, u_{k}^{j} \leq 1/2 \;\;\; \mbox{for each $j \in \{1, 2, \ldots, q\}$};
\end{equation}
\begin{equation*}
{\rm spt} \, \|V_{k}\| \cap ({\mathbf R} \times (B_{\s} \setminus \Sigma_{k})) = \cup_{j=1}^{q} {\rm graph} \, u_{k}^{j} \cap ({\mathbf R} \times (B_{\s} \setminus \Sigma_{k})); 
\end{equation*}
\begin{equation}\label{step2-1}
\|V_{k}\|({\mathbf R} \times \Sigma_{k}) + {\mathcal H}^{n} \, (\Sigma_{k}) \leq C{\hat E}_{k}^{2}
\end{equation}
where $C = (n, q, \s) \in (0, \infty);$ and ${\hat E}_{k}^{-1}u_{k}^{j} \to v^{j}$ for each $j = 1, 2, \ldots, q$, where 
the convergence is in $L^{2}(B_{\s})$ and weakly in $W^{1, 2}(B_{\s}).$

To verify that $v$ satisfies property $({\mathcal B{\emph 3}})$, note that by (\ref{firstvariation}), for each $k$ and each function $\z \in C^{1}_{c}(B_{\s})$, we have that 
\begin{equation}\label{step2-2}
\int \nabla^{V_{k}} \, x^{1} \cdot \nabla^{V_{k}} \, {\widetilde \z} \,d\|V_{k}\|(X) = 0
\end{equation}
where ${\widetilde \z}$ is any function in  $C^{1}_{c}({\mathbf R} \times B_{\s})$ such that 
${\widetilde \z} \equiv \z_{1}$ in a neighborhood of ${\rm spt} \, \|V_{k}\| \cap ({\mathbf R} \times B_{\s})$, where $\z_{1}(X)$ is defined for $X = (x^{1}, x^{\prime}) \in {\mathbf R} \times B_{\s}$ by 
$\z_{1}(x^{1}, x^{\prime}) = \z(x^{\prime})$.  Since $x^{1} = {\widetilde u}_{k}^{j}(X)$ 
for $\|V_{k}\|$ a.e. $X = (x^{1}, x^{\prime}) \in {\rm graph} \, u_{k}^{j} \cap {\rm spt} \, \|V_{k}\|$, 
where ${\widetilde u}_{k}^{j}(x^{1}, x^{\prime}) = u_{k}^{j}(x)$ for $(x^{1}, x^{\prime}) \in {\mathbf R} \times B_{\s}$, we deduce from (\ref{step2-2}) that 
\begin{eqnarray*}
\sum_{j=1}^{q}\int_{B_{\s}}(1 + |Du_{k}^{j}|^{2})^{-1/2} \, Du_{k}^{j} \cdot D\z &=& 
-\int_{{\mathbf R} \times (B_{\s} \cap \Sigma_{k})} \nabla^{V_{k}} \, x^{1} \cdot \nabla^{V_{k}} \, 
{\widetilde \z} \, d\|V_{k}\|(X)\nonumber\\ 
&&+ \sum_{j=1}^{q}\int_{B_{\s} \cap \Sigma_{k}} (1 + |Du_{k}^{j}|^{2})^{-1/2} \, Du_{k}^{j} \cdot D\z 
\end{eqnarray*}
which can be rewritten as 
\begin{eqnarray}\label{step2-3}
\sum_{j=1}^{q}\int_{B_{\s}}Du_{k}^{j} \cdot D\z =-\int_{{\mathbf R} \times (B_{\s} \cap \Sigma_{k})} \nabla^{V_{k}} \, x^{1} \cdot \nabla^{V_{k}} \, 
{\widetilde \z} \, d\|V_{k}\|(X)&&\nonumber\\ 
&&\hspace{-3in}+ \; \sum_{j=1}^{q}\int_{B_{\s} \cap \Sigma_{k}} (1 + |Du_{k}^{j}|^{2})^{-1/2} \, Du_{k}^{j} \cdot D\z \; + \; F_{k} \;\;\; {\rm where}
\end{eqnarray}
\begin{eqnarray}\label{step2-4}
|F_{k}| &=& \left|\sum_{j=1}^{q}\int_{B_{\s}} (1 +|Du_{k}^{j}|^{2})^{-1/2}(1 + (1 + |Du_{k}^{j}|^{2})^{1/2})^{-1}|Du_{k}^{j}|^{2} Du_{k}^{j} \cdot D\z\right|\nonumber\\
&&\leq \sup \, |D\z|\int_{{\mathbf R} \times B_{\s}} |\nabla^{V_{k}} \, x^{1}|^{2} \, d\|V_{k}\|(X)\nonumber\\ 
&&\hspace{.5in}+ \; \sup \, |D\z|\sum_{j=1}^{q}\int_{B_{\s} \cap \Sigma_{k}} (1 +|Du_{k}^{j}|^{2})^{-1/2}(1 + (1 + |Du_{k}^{j}|^{2})^{1/2})^{-1}|Du_{k}^{j}|^{3}\nonumber\\
&&\leq \sup \, |D\z| \left(C{\hat E}_{k}^{2} + q{\mathcal H}^{n} \, (\Sigma_{k})\right).
\end{eqnarray}
The last inequality in (\ref{step2-4}),  where $C = C(n, \s) \in (0, \infty),$ follows from (\ref{tilt-ht}) and 
(\ref{step2-0}).

Dividing both sides of (\ref{step2-3}) by ${\hat E}_{k}$ and letting $k \to \infty$, we deduce, using 
(\ref{step2-0}), (\ref{step2-1}) and (\ref{step2-4}), that
$$\sum_{j=1}^{q}\int_{B_{\s}} Dv^{j} \cdot D\z =0$$
for any $\z \in C^{1}_{c}(B_{\s}).$ Since $\s \in (0, 1)$ is arbitrary, this implies that 
$\Delta \, v_{a} = 0$ in $B_{1}$, establishing property $({\mathcal B{\emph 3}})$ for ${\mathcal B}_{q}.$

Next we verify that ${\mathcal B}_{q}$ satisfies properties $({\mathcal B}{\emph 5 \, I}),$ $({\mathcal B}{\emph 5 \, II})$, (${\mathcal B}{\emph 6}$) and (${\mathcal B}{\emph 5 \, III})$, in that order:

Let $z \in B_{1}$, $\s \in (0, (1 - |z|)]$ and $\g$ be an orthogonal rotation of ${\mathbf R}^{n}$, and note that $\widetilde{v}_{z, \s} \equiv  \|v(z + \s(\cdot))\|_{L^{2}(B_{1})}^{-1}v(z + \s(\cdot))$ is the coarse blow-up of the sequence $\{\eta_{(0, z), \s \, \#} \, V_{k}\},$ and $v\circ \g$ is the coarse blow-up of the sequence $\{\widetilde{\g}_{\#} \, V_{k}\}$, where $\widetilde{\g} \, : \, {\mathbf R}^{n+1} \to {\mathbf R}^{n+1}$ is 
the orthogonal rotation defined by $\widetilde{\g}(x^{1}, x^{\prime}) = (x^{1}, \g(x^{\prime})).$ Thus ${\mathcal B}_{q}$ 
satisfies properties $({\mathcal B{\emph 5 \, I}})$ and $({\mathcal B{\emph 5 \, II}})$. 

To verify that ${\mathcal B}_{q}$ satisfies property $({\mathcal B{\emph 6}})$, let $\{v_{\ell}\}_{\ell=1}^{\infty}$ be a sequence of elements in ${\mathcal B}_{q},$ and for each $\ell=1, 2, \ldots,$ let $\{V^{\ell}_{k}\}_{k=1}^{\infty} \subset {\mathcal S}_{\a}$ be a sequence whose coarse blow-up is $v_{\ell}.$ Choose, for each $\ell = 1, 2, \ldots$, a positive integer $k_{\ell}$ such that $k_{1} < k_{2} < k_{3} < \ldots$, ${\hat E}_{V^{\ell}_{k_{\ell}}} < \min\{\ell^{-1}, \e_{0}(n, q, 1 - \ell^{-1})\},$ where 
$\e_{0}$ is as in Theorem~\ref{flat-varifolds}, and 
$\|{\hat E}_{V^{\ell}_{k_{\ell}}}^{-1}u_{\ell, k_{\ell}} - v_{\ell}\|_{L^{2}(B_{1 - \ell^{-1}})} < \ell^{-1}$, where $u_{\ell, k_{\ell}} = (u_{\ell, k_{\ell}}^{1}, u_{\ell, k_{\ell}}^{2}, \ldots, u_{\ell, k_{\ell}}^{q}) \, : \, B_{1 - \ell^{-1}} \to {\mathbf R}^{q}$ is the Lipschitz function (with Lipschitz constant of each component function $\leq 1/2$) corresponding to $u = (u^{1}, u^{2}, \ldots, u^{q})$ 
of Theorem~\ref{flat-varifolds} taken with $V^{\ell}_{k_{\ell}}$ in place of $V$ and with $\s = 1 - \ell^{-1}.$That such a choice 
exists follows from the definition of coarse blow-up. Note also that it follows from (\ref{blow-up-5}) and (\ref{blow-up-6}) that for each $\s \in (0, 1)$ and all sufficiently large $\ell$, $\int_{B_{\s}} |v_{\ell}|^{2} +  |Dv_{\ell}|^{2} < C$, where $C = C(n, q, \s) \in (0, \infty)$ is independent of $\ell.$ Let $v \in {\mathcal B}_{q}$ be the coarse blow-up of 
an appropriate subsequence $\{V^{\ell^{\prime}}_{k_{\ell^{\prime}}}\}$ of the sequence $\{V^{\ell}_{k_{\ell}}\}$. It is then straightforward to check, after passing to a subsequence of $\{\ell^{\prime}\}$ without changing notation, that for each $\s \in (0, 1)$, $v_{\ell^{\prime}} \to v$ in $L^{2}(B_{\s})$ and weakly in $W^{1, 2}(B_{\s}).$ 

In order to verify that ${\mathcal B}_{q}$ satisfies property $({\mathcal B{\emph 5 \, III}}),$ note first that if $y \in {\mathbf R}$ is a constant and $v - y \not\equiv 0$ in $B_{1}$, then $\|v - y\|_{L^{2}(B_{1})}^{-1}(v-y)\in {\mathcal B}_{q}$,  where we have used the notation $v - y = (v^{1} - y, v^{2} - y, \ldots, v^{q} - y).$ To check this, note that $v(\s(\cdot)) - y \not\equiv 0$ for 
all sufficiently large $\s \in (0, 1),$ and  that for any such $\s$, $\|v(\s(\cdot)) - y\|_{L^{2}(B_{1})}^{-1}\left((v(\s(\cdot)) - y\right)$ is the coarse blow-up of the sequence $\{\t_{k \, \#} \, \eta_{\s \, \#} \, V_{k}\}$ where $\t_{k} \, : \, {\mathbf R}^{n+1} \to {\mathbf R}^{n+1}$ is the translation $X \mapsto X - ({\hat E}_{k}y, 0).$ Thus $\|v(\s(\cdot)) - y\|_{L^{2}(B_{1})}^{-1}\left(v(\s(\cdot)) - y\right)\in {\mathcal B}_{q}$ for all sufficiently large $\s \in (0, 1)$, and hence it follows from property $({\mathcal B{\emph 6}})$ that $\|v - y\|_{L^{2}(B_{1})}^{-1}(v-y) \in {\mathcal B}_{q}$ as claimed. Next note that if $L \, : \, {\mathbf R}^{n} \to {\mathbf R}$ is a linear function and 
$v - L \not\equiv 0$ in $B_{1}$, then $\|v - L\|_{L^{2}(B_{1})}^{-1}(v-L) \in {\mathcal B}_{q},$ where, $v - L = (v^{1} - L, v^{2} - L, \ldots, v^{q} - L).$ To check this, assume without loss of generality (in view of property $({\mathcal B{\emph 5 \, II}})$) that $L(x) = \lambda x^{2}$ for some $\lambda \in {\mathbf R},$ and note that for sufficiently large $\s \in (0, 1)$, 
$\|v(\s(\cdot)) - \s L\|_{L^{2}(B_{1})}^{-1}\left(v(\s(\cdot)) - \s L\right)$ is the coarse blow-up of the sequence 
$\{\G_{k \, \#} \, \eta_{\s \, \#} \, V_{k}\}$, where $\G_{k}  \, : \, {\mathbf R}^{n+1} \to {\mathbf R}^{n+1}$ is the rotation 
fixing $\{0\} \times {\mathbf R}^{n-1}$ pointwise and mapping the unit normal $\nu_{k} = \left(1 + {\hat E}_{k}^{2}\lambda^{2}\right)^{-1/2}\left(1, -{\hat E}_{k}\lambda, 0\right)$ to the hyperplane $P_{k} \equiv {\rm graph} \, {\hat E}_{k}L$ to $e^{1}.$  Thus $\|v(\s(\cdot)) - \s L\|_{L^{2}(B_{1})}^{-1}\left(v(\s(\cdot)) - \s L\right) \in {\mathcal B}_{q}$  for all sufficiently large $\s \in (0, 1)$, and it follows from property $({\mathcal B{\emph 6}})$ that $\|v - L\|_{L^{2}(B_{1})}(v - L) \in {\mathcal B}_{q}$ as claimed. We deduce that ${\mathcal B}_{q}$ satisfies property $({\mathcal B{\emph 5 \, III}})$ by applying the above facts with $y = v_{a}(0)$ and with the linear function $L$ defined by 
$L(x) = \|v - v_{a}(0)\|_{L^{2}(B_{1})}^{-1}Dv_{a}(0) \cdot x$ for $x \in {\mathbf R}^{n}.$ (Note that $v  - v_{a}(0) \not\equiv 0$ in $B_{1}$ or else  $v - \ell_{v} \equiv 0$ in $B_{1}$, contrary to the hypothesis of $({\mathcal B{\emph 5 \, III}})$, where $\ell_{v}$ is as in the statement of $({\mathcal B{\emph 5 \, III}})$.) Note that our argument shows more generally that 
\begin{equation}\label{step2-4-1-0}
v \in {\mathcal B}_{q}, \;\; v - \ell_{v, \, z} \not\equiv 0 \;\; \mbox{in} \;\;B_{1} \;\; \implies \;\; \|v  - \ell_{v, \, z}\|_{L^{2}(B_{1})}^{-1}\left(v - \ell_{v, \, z}\right) \in {\mathcal B}_{q}
\end{equation}
for each $z \in B_{1},$ where $\ell_{v, \, z}(x) = v_{a}(z) + Dv_{a}(z) \cdot (x - z)$ and $v - \ell_{v, \, z} = (v^{1} - \ell_{v, \, z}, \ldots, v^{q} - \ell_{v, \, z}).$
 
Finally in this section, we verify that ${\mathcal B}_{q}$ satisfies property $({\mathcal B{\emph 4}})$ with a constant $C = C(n, q) \in (0, \infty)$ to be specified momentarily. First note that for any stationary integral $n$-varifold $V$ on $B_{2}^{n+1}(0)$ with ${\hat E}_{V}$ sufficiently small and satisfying the hypotheses of Theorem~\ref{flat-varifolds} taken with $\s = 15/16$ and for any $Z = (z^{1}, z^{\prime}) \in {\rm spt} \, \|V\| \cap B_{1/8}^{n+1}(0)$ with $\Theta \, (\|V\|, Z) \geq q$, we have that  
\begin{equation}\label{step2-4-1}
\sum_{j=1}^{q}\int_{B_{1/2}(z^{\prime}) \setminus \Sigma}
\left(\frac{R_{z}^{2}}{(u^{j} - z^{1})^{2} + R_{z}^{2}}\right)^{\frac{n+2}{2}} R_{z}^{2-n} \left(\frac{\partial \, \left((u^{j} - z^{1})/R_{z}\right)}{\partial \, R_{z}}\right)^{2} d{\mathcal H}^{n}(x) \leq C_{2}{\hat E}_{V}^{2}
\end{equation}
where $R_{z}(x) = |x - z|$ for $x \in {\mathbf R}^{n}$ and $C_{2} = C_{2}(n, q) \in (0, \infty);$ the set $\Sigma \subset B_{15/16}$ here and the functions 
$u_{j}$, $j=1, 2, \ldots, q$ are as in Theorem~\ref{flat-varifolds} taken with $\s = 15/16.$ To see this, note that by estimating as in (\ref{monotonicity-2}), it follows that
\begin{eqnarray*}
\int_{B^{n+1}_{3/4}(Z)} \frac{|(X - Z)^{\perp}|^{2}}{|X - Z|^{n+2}} d\|V\|(X) \leq C_{2}{\hat E}_{V}^{2}, \;\; C_{2} = C_{2}(n, q) \in (0, \infty),
\end{eqnarray*}
while
\begin{eqnarray*}
\int_{B^{n+1}_{3/4}(Z)} \frac{|(X - Z)^{\perp}|^{2}}{|X - Z|^{n+2}} d\|V\|(X) \geq \int_{{\mathbf R} \times (B_{1/2}(z^{\prime}) \setminus \Sigma)} \frac{|(X - Z)^{\perp}|^{2}}{|X - Z|^{n+2}} d\|V\|(X)&&\nonumber\\
&&\hspace{-4in}\geq\frac{1}{2}\sum_{j=1}^{q} \int_{B_{1/2}(z^{\prime}) \setminus \Sigma} 
\frac{\left((x^{\prime} - z^{\prime}) \cdot Du^{j}(x^{\prime}) - (u^{j}(x^{\prime}) - z^{1})\right)^{2}}
{\left((u^{j}(x^{\prime}) - z^{1})^{2} + |x^{\prime} - z^{\prime}|^{2}\right)^{\frac{n+2}{2}}}d{\mathcal H}^{n}(x^{\prime})\nonumber\\
&&\hspace{-5in}=\frac{1}{2}\sum_{j=1}^{q}\int_{B_{1/2}(z^{\prime}) \setminus \Sigma}
\left(\frac{R_{z^{\prime}}^{2}}{(u^{j} - z^{1})^{2} + R_{z^{\prime}}^{2}}\right)^{\frac{n+2}{2}} R_{z^{\prime}}^{2-n} \left(\frac{\partial \, \left((u^{j} - z^{1})/R_{z^{\prime}}\right)}{\partial \, R_{z^{\prime}}}\right)^{2} d{\mathcal H}^{n}(x^{\prime}).
\end{eqnarray*}

Now let $v \in {\mathcal B}_{q}$ and let $z \in B_{1}$ be such that $({\mathcal B{\emph 4 \, I}})$ with $C = C_{2}$, where $C_{2} = C_{2}(n, q)$ is as 
 in (\ref{step2-4-1}), fails. By (\ref{step2-4-1-0}), $\widetilde{v} \equiv \|v - \ell_{v, \, z}\|^{-1}_{L^{2}(B_{1})}(v - \ell_{v, \, z}) \in {\mathcal B}_{q}.$ Let $V_{k} \in {\mathcal S}_{\a}$ be such that $\widetilde{v}$ is the coarse blow-up of $\{V_{k}\}$ . We claim that then there exists $\s_{1} >0$ such that for all sufficiently large $k$,
\begin{equation}\label{step2-5}
Z \in {\rm spt} \, \|V_{k}\| \cap ({\mathbf R} \times B_{\s_{1}}(z)) \implies \Theta \, (\|V_{k}\|, Z) < q.
\end{equation}
If not, then there would exist, for each positive integer $\ell$, a positive integer  $\{k_{\ell}\}$ 
with $k_{1} < k_{2} < k_{3} < \ldots,$ a point 
$Z_{\ell} = (z^{1}_{\ell}, z^{\prime}_{\ell}) \in {\rm spt} \, \|V_{k_{\ell}}\| \cap ({\mathbf R} \times B_{1/\ell}(z))$ such that  
$\Theta \, (\|V_{k_{\ell}}\|, Z_{\ell}) \geq q.$ Fix any $\r \in (0, \frac{3}{8}(1 - |z|)].$ Applying (\ref{step2-4-1}) with $\eta_{Z_{\ell}, \r \, \#} \, V_{k_{\ell}}$ in place of $V$ and $0$ in place of $Z$, 
we then have, after changing variables, that for all sufficiently large $\ell$, 
\begin{eqnarray}\label{step2-7}
\sum_{j=1}^{q}\int_{B_{\r/2}(z_{\ell}^{\prime}) \setminus \Sigma_{k_{\ell}}}
\left(\frac{R_{z_{\ell}^{\prime}}^{2}}{(u_{k_{\ell}}^{j} - z_{\ell}^{1})^{2} + R_{z_{\ell}^{\prime}}^{2}}\right)^{\frac{n+2}{2}} R_{z_{\ell}^{\prime}}^{2-n} \left(\frac{\partial \, 
\left((u_{k_{\ell}}^{j} - z_{\ell}^{1})/R_{z_{\ell}^{\prime}}\right)}{\partial \, R_{z_{\ell}^{\prime}}}\right)^{2} d{\mathcal H}^{n}(x)&&\nonumber\\ 
&&\hspace{-2in}\leq C_{2} \, \r^{-n-2}\int_{{\mathbf R} \times B_{\r}(z^{\prime}_{\ell})}|x^{1}|^{2}d\|V_{k_{\ell}}\|(X)
\end{eqnarray}

Now for all sufficiently large $\ell$ depending on $\r$, 
$\|V_{k_{\ell}}\|({\mathbf R} \times B_{\r/16}(z^{\prime}_{\ell})) \geq 
C\r^{n}$ for a suitable constant $C = C(n) \in (0, \infty)$, so there exists a  
point $Y_{\ell} = (y_{\ell}^{1}, y_{\ell}^{\prime})  \in {\rm spt} \, \|V_{k_{\ell}}\| 
\cap ({\mathbf R} \times (B_{\r/16}(z^{\prime}_{\ell}))$ such that 
\begin{equation}\label{step2-8}
|y_{\ell}^{1}|^{2} \leq C\r^{-n} \int_{{\mathbf R} \times B_{\r/16}(z^{\prime}_{\ell})} |x^{1}|^{2} \, d\|V_{k_{\ell}}\|(X)
\end{equation}
where $C = C(n) \in (0, \infty).$ Applying Theorem~\ref{non-concentration}(a)
with ${\widetilde V} = \eta_{Y_{\ell}, \r/2 \, \#} \, V_{k_{\ell}}$ in place of $V$ and ${\widetilde Z} = (\r/2)^{-1}(Z_{\ell} - Y_{\ell})$ in place of $Z$ (noting that ${\widetilde Z} \in {\rm spt} \, \|{\widetilde V}\| \cap 
({\mathbf R} \times B_{1/8})$ with $\Theta \, (\|{\widetilde V}\|, {\widetilde Z}) \geq q$), we deduce, 
using also (\ref{step2-8}), that 
\begin{equation}\label{step2-9}
|z^{1}_{\ell}|^{2} \leq C\r^{-n}\int_{{\mathbf R} \times B_{3\r/4}(z_{\ell}^{\prime})}|x^{1}|^{2}d\|V_{k_{\ell}}\|(X)
\end{equation}
for all sufficiently large $\ell$, where $C = C(n, q) \in (0, \infty).$ Dividing both sides of 
(\ref{step2-7}) by ${\hat E}_{k_{\ell}}^{2}$, and letting $\ell \to \infty,$ we conclude, using 
(\ref{step2-9}) and the fact that 
$\sup_{X \in {\rm spt} \, \|V_{k_{\ell}}\| \cap ({\mathbf R} \times B_{3/4}))} \, |x^{1}| \to 0$ as $\ell \to \infty$, that 
\begin{equation}\label{step2-9-1}
\sum_{j=1}^{q}\int_{B_{\r/2}(z)} R_{z}^{2-n} \left(\frac{\partial \, 
\left((\widetilde{v}^{j} - y)/R_{z}\right)}{\partial \, R_{z}}\right)^{2} \leq C_{2} \, \r^{-n-2}\int_{B_{\r}(z)}|\widetilde{v}|^{2}
\end{equation}
for some $y \in {\mathbf R}$ and each $\r \in (0, \frac{3}{8}(1 - |z|].$ Since by the triangle inequality this implies that 
$\int_{B_{\r/2}(z)} R_{z}^{2-n} \left(\frac{\partial \, 
\left((\widetilde{v}_{a} - y)/R_{z}\right)}{\partial \, R_{z}}\right)^{2} < \infty$, it follows that $y = \widetilde{v}_{a}(z) = 0.$ But this contradicts our assumption that property $({\mathcal B{\emph 4 \, I}})$ fails for $v$, leading us to the conclusion that
(\ref{step2-5}) must hold for all sufficiently large $k$.
 
By Remark 3 of Section~\ref{outline} and (\ref{step2-5}), it follows that for all sufficiently large $k$, 
${\mathcal H}^{n-7+\g} \, ({\rm sing} \, V_{k} \cap ({\mathbf R} \times B_{\s_{1}}(z))) = 0$ for every $\g >0$ if 
$n \geq 7$ and ${\rm sing} \, V_{k} \cap ({\mathbf R} \times B_{\s_{1}}(z)) = \emptyset$ if 
$2 \leq n \leq 6$, so we may apply Theorem~\ref{SS} and standard elliptic theory to conclude that 
$$V_{k} \res ({\mathbf R} \times B_{\s_{1}/2}(z)) = \sum_{j=1}^{q} |{\rm graph} \, u_{k}^{j}|$$
where $u_{k}^{j} \, : \, B_{\s_{1}/2}(z) \to {\mathbf R}$ are $C^{2}$ functions satisfying 
$$\sup_{B_{\s_{1}/2}(z)} \, \sum_{j=1}^{q} |Du_{k}^{j}| + |D^{2} \, u_{k}^{j}| \leq C{\hat E}_{k}$$
and solving the minimal surface equation on $B_{\s_{1}/2}(z),$ where $C = C(n, q, \s) \in (0, \infty).$
This readily shows that $\Delta \, \widetilde{v}^{j} = 0$ on $B_{\s_{1}/2}(z)$ for each $j=1, 2, \ldots, q,$
establishing property $({\mathcal B{\emph 4}})$ for ${\mathcal B}_{q}.$

\noindent
{\bf Remarks:} {\bf (1)} The argument leading to (\ref{step2-9-1}) proves the following: 

\emph{Let $\Omega$ be an open subset of $B_{3/4}$. If $v \in {\mathcal B}_{q}$ and $\{V_{k}\} \subset {\mathcal S}_{\a}$ is a sequence whose coarse blow-up is $v$ (in the sense described in Section~\ref{blow-up}), and if for infinitely many $k$, there are points $Z_{k} \in {\rm spt} \, \|V_{k}\| \cap ({\mathbf R} \times \Omega)$ with $\Theta \, (\|V_{k}\|, Z_{k}) \geq q$, then there exists a point $z \in \overline{\Omega}$ such that
\begin{equation*}
\sum_{j=1}^{q}\int_{B_{\r/2}(z)} R_{z}^{2-n} \left(\frac{\partial \, 
\left(({v}^{j} - v_{a}(z))/R_{z}\right)}{\partial \, R_{z}}\right)^{2} \leq C_{2} \, \r^{-n-2}\int_{B_{\r}(z)}|{v}|^{2}
\end{equation*}
for each $\r \in (0, \frac{3}{8}(1- |z|)].$}
 
\noindent
{\bf (2)} \emph{Let $q$ be an integer $\geq 2$. There exist constants $\eta^{\prime} = \eta^{\prime}(n, q, \a) \in (0, 1)$ and 
$\d^{\prime} = \d^{\prime}(n, q, \a) \in (0, 1)$ such that the following is true: If the induction hypotheses $(H1)$, $(H2)$ hold,
$V \in {\mathcal S}_{\a},$ $(\omega_{n}2^{n})^{-1}\|V\|(B_{2}^{n+1}(0)) < q + 1/2,$ $\omega_{n}^{-1}\|V\|({\mathbf R} \times B_{1}) < q + 1/2$, $\int_{{\mathbf R} \times B_{1}} {\rm dist}^{2}(X, {\mathbf P}) \, d\|V\|(X) < \d^{\prime}$ for some union ${\mathbf P} \subset {\mathbf R}^{n+1}$ of finitely many (distinct) affine hyperplanes disjoint in ${\mathbf R} \times B_{1}$ with ${\rm dist}_{\mathcal H} \, ({\mathbf P} \cap ({\mathbf R} \times B_{1}), \{0\} \times B_{1}) < \d^{\prime}$ and, writing ${\mathcal A}$ for the set of affine hyperplanes of ${\mathbf R}^{n+1}$, if
$$\int_{{\mathbf R} \times B_{1}} {\rm dist}^{2}(X, {\mathbf P}) \, d\|V\|(X) < \eta^{\prime} \inf_{L \in {\mathcal A}} \int_{{\mathbf R} \times B_{1}} {\rm dist}^{2} \, (X, L) \, d\|V\|(X),$$
then 
\begin{itemize}
\item[(a)] ${\mathbf P}$ consists of at least two affine hyperplanes; 
\item[(b)] $\{Z \in {\rm spt} \, \|V\| \cap ({\mathbf R} \times B_{3/4}) \, : \, \Theta \, (\|V\|, Z) \geq q\}  = \emptyset$; 
\item[(c)] there exist an integer $p$ with $2 \leq p \leq q,$ positive integers $a_{j} \leq q-1,$  affine hyperplanes $P_{j}^{i}  \subset {\mathbf P},$ $C^{2}$ functions $u_{j}^{i} \, : \, P_{j}^{1} \cap ({\mathbf R} \times B_{3/4}) \to (P_{j}^{1})^{\perp}$ with 
$u_{j}^{1}\cdot e_{1} \leq \ldots \leq u_{j}^{a_{j}} \cdot e_{1}$ for $1 \leq j \leq p$, $1 \leq i \leq a_{j}$ and  
$u_{j-1}^{a_{j-1}}  \cdot e_{1}< u_{j}^{1}\cdot e_{1}$ for $2 \leq j \leq p$ such that 
$\|u_{j}^{i}\|^{2}_{C^{2}(P_{j}^{1} \cap ({\mathbf R} \times B_{3/4}))} < C \int_{{\mathbf R} \times B_{1}}{\rm dits}^{2} \, (X, {\mathbf P}) \, d\|V\|(X),$ $V\res ({\mathbf R} \times B_{5/8}) = \sum_{j=1}^{p} V_{j}$ where 
$V_{j}  = \sum_{i=1}^{a_{j}}|{\rm graph} \, u_{j}^{i} \cap ({\mathbf R} \times B_{5/8})|$   
 and 
$$\int_{{\mathbf R} \times B_{5/8}} {\rm dist}^{2}(X, {\mathbf P}) \, d\|V\|(X)  = \sum_{j=1}^{p}   
\int_{{\mathbf R} \times B_{5/8}} {\rm dist}^{2}(X, {\mathbf P}_{j}) \, d\|V_{j}\|(X)$$
where ${\mathbf P}_{j} = \cup_{i=1}^{a_{j}}P_{j}^{i}$.  Here ${\rm graph} \, u_{j}^{i} = \{X + u_{j}^{i}(X) \, : \, X \in P_{j} \cap ({\mathbf R} \times B_{3/4})\}$.
\end{itemize}}
To see this, argue by contradiction: Were the assertion false, we can find a sequence $V_{k} \in {\mathcal S}_{\a}$ with $(\omega_{n}2^{n})^{-1}\|V_{k}\|(B_{2}^{n+1}(0)) < q + 1/2,$ $\omega_{n}^{-1}\|V_{k}\|({\mathbf R} \times B_{1}) < q + 1/2$  and for each $k$, affine hyperplanes $P_{k}^{1}, \ldots, P_{k}^{n_{k}}$ with $P_{k}^{i} \cap P_{k}^{j} \cap ({\mathbf R} \times B_{1}) = \emptyset$ for $1 \leq i < j \leq n_{k}$ and ${\rm dist}_{\mathcal H} \, ({\mathbf P}_{k} \cap ({\mathbf R} \times B_{1}), \{0\} \times B_{1}) \to 0$ as $k \to \infty$ where ${\mathbf P}_{k} = \cup_{j=1}^{n_{k}} P_{k}^{j},$ such that $\int_{{\mathbf R} \times B_{1}} {\rm dist}^{2}(X, {\mathbf P}_{k}) \, d\|V_{k}\|(X)  \to 0$ and 
\begin{equation}\label{separation'}
\left(\inf_{L \in {\mathcal A}} \int_{{\mathbf R} \times B_{1}} {\rm dist}^{2} \, (X, L) \, d\|V_{k}\|(X)\right)^{-1} \int_{{\mathbf R} \times B_{1}} {\rm dist}^{2}(X, {\mathbf P}_{k}) \, d\|V_{k}\|(X)  \to 0
\end{equation}
and yet, at least one of the conclusions (a)-(c) with $V_{k}$ in place of $V$ and ${\mathbf P}_{k}$  in place of ${\mathbf P}$ fails. Note that $\inf_{L \in {\mathcal A}} \int_{{\mathbf R} \times B_{1}} {\rm dist}^{2} \, (X, L) \, d\|V_{k}\|(X) \to 0$, and choose $L_{k} \in {\mathcal A}$ such that 
$$\int_{{\mathbf R} \times B_{1}} {\rm dist}^{2} \, (X, L_{k}) \, d\|V_{k}\|(X) < \frac{3}{2} \inf_{L \in {\mathcal A}} \int_{{\mathbf R} \times B_{1}} {\rm dist}^{2} \, (X, L) \, d\|V_{k}\|(X).$$ Noting then that 
$L_{k} \to \{0\} \times {\mathbf R}^{n},$ choose rigid motions $\G_{k} \, : \, {\mathbf R}^{n+1} \to {\mathbf R}^{n+1}$ such that $\G_{k} \to {\rm Identity}$ and $\G_{k}(L_{k}) = \{0\} \times {\mathbf R}^{n}$ and let $v = (v^{1}, \ldots, v^{\ell}) \in W^{1, 2}_{\rm loc}(B_{1}; {\mathbf R}^{p}) \cap L^{2}(B_{1}; {\mathbf R}^{p}),$ with $v^{1} \leq v^{2} \ldots \leq v^{\ell}$, be the coarse blow-up, as described in Section~\ref{blow-up}, of (a suitable subsequence of) the sequence $\{{\widetilde V}_{k} = \eta_{0, 13/16 \, \#} \, \G_{k \, \#} \, V_{k}\}$ relative to $\{0\} \times {\mathbf R}^{n}$, where $\ell$ is a positive integer $\leq q$. Let $p \leq \ell$ be the number of distinct functions in the set $\{v^{1}, \ldots, v^{\ell}\},$ denoted $\widetilde{v}^{1}, \ldots, \widetilde{v}^{p}$ with the labelling so chosen that $\widetilde{v}^{1} \leq \ldots \leq \widetilde{v}^{p}.$ Then by (\ref{separation'}), for each $k$, there exists $\{{\widetilde P}_{k}^{1}, {\widetilde P}_{k}^{2}, \ldots, {\widetilde P}_{k}^{p}\} \subset \{P_{k}^{1}, P_{k}^{2}, \ldots, P_{k}^{n_{k}}\}$ such that, writing $\G_{k}{\widetilde P}_{k}^{i} = {\rm graph} \, \widetilde{p}_{k}^{i}$ for an affine function $\widetilde{p}_{k}^{i} \, : \, {\mathbf R}^{n} \to {\mathbf R}$ with labelling so chosen that $\widetilde{p}_{k}^{1} < \ldots < \widetilde{p}_{k}^{p}$ in ${\mathbf R} \times B_{1}$, we have that $\widetilde{v}^{j} = \lim_{k \to \infty} \left({\hat E}_{k}\right)^{-1}\widetilde{p}_{k}^{j}$ for $ 1 \leq j \leq p.$ Thus each $v^{j}$ is affine, and by (\ref{separation'}) again, $p \geq 2$ and $\widetilde{v}^{p} > \widetilde{v}^{1}$ in $B_{1}.$ It then follows from Remark (1) above (taken with $\ell$ in place of $q$) that $\{Z \in {\rm spt} \, \|V_{k}\| \cap ({\mathbf R} \times B_{3/4}) \, : \, \Theta \, (\|V_{k}\|, Z) \geq \ell\}  = \emptyset$ for sufficiently large $k$. The rest of the conclusions with $V_{k}$ in place of $V$ and ${\mathbf P}_{k}$ in place of ${\mathbf P}$ now follow, for all sufficiently large $k$, from 
Remark 3 of Section~\ref{outline}, Theorem~\ref{SS} and standard elliptic estimates, contrary to the assumption that at least one of those conclusions must fail for each $k$.

\noindent
{\bf (3)} \emph{Let $q$ be an integer $\geq 2$. There exists a constant $\d = \d(n, q, \a) \in (0, 1)$ such that the following is true: If the induction hypotheses $(H1)$, $(H2)$ hold,
$V \in {\mathcal S}_{\a},$ $(\omega_{n}2^{n})^{-1}\|V\|(B_{2}^{n+1}(0)) < q + 1/2,$ $\omega_{n}^{-1}\|V\|({\mathbf R} \times B_{1}) < q + 1/2$ and 
$$\int_{{\mathbf R} \times B_{1}} {\rm dist}^{2}(X, {\mathbf P}) \, d\|V\|(X) < \d$$
for some union ${\mathbf P} \subset {\mathbf R}^{n+1}$ of at most $q$ affine hyperplanes disjoint in ${\mathbf R} \times B_{1}$ with ${\rm dist}_{\mathcal H} \, ({\mathbf P} \cap ({\mathbf R} \times B_{1}), \{0\} \times B_{1}) < \d,$ then either 
\begin{itemize}
\item[(a)] $\{Z \in {\rm spt} \, \|V\| \cap ({\mathbf R} \times B_{7/8}) \, : \, \Theta \, (\|V\|, Z) \geq q\}  = \emptyset$ and there exist a positive integer  $\ell$ with $1 \leq \ell \leq q$, distinct affine hyperplanes $P_{1}, P_{2}, \ldots, P_{\ell} \subset {\mathbf P},$  positive integers $q_{1}, q_{2}, \ldots, q_{\ell}$ with $\sum_{k=1}^{\ell} q_{k} \leq q$ and  $C^{2}$ functions $u_{k}^{j} \, : \, P_{k} \cap ({\mathbf R} \times B_{3/4}) \to P_{k}^{\perp}$ with 
$$\sup_{P_{k} \cap ({\mathbf R} \times B_{3/4})} \,|u_{k}^{j}|^{2} + |Du_{k}^{j}|^{2} \leq C\int_{{\mathbf R} \times B_{1}} {\rm dist}^{2} \, (X, {\mathbf P}) \, d\|V\|(X)$$
for $1 \leq k \leq \ell$, $1 \leq j \leq q_{k}$ 
where $C = C(n)$, such that $$V \res ({\mathbf R} \times B_{1/2}) = \sum_{k=1}^{\ell}\sum_{j=1}^{q_{k}}|{\rm graph} \, u_{k}^{j} \cap ({\mathbf R} \times B_{1/2})|, \;\;\;{\rm or}$$
\item[(b)] $\{Z \in {\rm spt} \, \|V\| \cap ({\mathbf R} \times B_{7/8}) \, : \, \Theta \, (\|V\|, Z) \geq q\}  \neq \emptyset,$ $\omega_{n}^{-1}\|V\|({\mathbf R} \times B_{1}) \geq q - 1/2$ and there exist an affine hyperplane $P \subset {\mathbf P},$  a measurable subset $\Sigma \subset P \cap ({\mathbf R} \times B_{13/28})$   Lipschitz functions $u_{1}, u_{2}, \ldots, u_{q}  \, : \, P \cap ({\mathbf R} \times B_{13/28}) \to P^{\perp}$ with ${\rm Lip} \, (u_{j}) \leq 9/16$ for each $j \in \{1, 2, \ldots, q\}$ such that 
$${\mathcal H}^{n} \, (\Sigma) + \|V\|({\mathcal C}_{P}(\Sigma)) + \sum_{j=1}^{q}\int_{P \cap ({\mathbf R} \times B_{13/28}) \setminus \Sigma} |u_{j}|^{2} + |Du_{j}|^{2} \leq C \int_{{\mathbf R} \times B_{1}} {\rm dist}^{2} \, (X, {\mathbf P}) \, d\|V\|(X)$$
and   
$$V \res (({\mathbf R} \times B_{13/28}) \setminus {\mathcal C}_{P}(\Sigma)) = \sum_{j=1}^{q} |{\rm graph} \, u_{j}\cap (({\mathbf R} \times B_{13/28}) \setminus {
\mathcal C}_{P}(\Sigma))|$$
where ${\mathcal C}_{P}(\Sigma) = \{X \in {\mathbf R}^{n+1} \, : \, \pi_{P}(X) \in \Sigma\}$ with $\pi_{P}$ denoting the orthogonal projection of ${\mathbf R}^{n+1}$ onto $P$; furthermore, in this case we have that for each $j \in \{1, 2, \ldots, q\}$, 
$$\sup_{B_{13/28}} |u_{j}| \leq C \d^{1/2n}$$
where $C = C(n) \in (0, \infty).$ 
\end{itemize}}

To see this, let $\eta^{\prime} = \eta^{\prime}(n, q, \a) \in (0, 1)$ and $\d^{\prime} = \d^{\prime}(n, q, \a) \in (0, 1)$ be the constants as in Remark (2) above, let $\e_{0} = \e_{0}(n, q, \a, 3/4) \in (0, 1)$ be the constant as in Theorem~\ref{flat-varifolds}. Let the hypotheses of the assertion of Remark (3) be satisfied for sufficiently small $\d \in (0, \eta^{\prime}\d^{\prime}\e_{0}],$ and note that it follows from the Constancy Theorem (\cite{S1}, Theorem 41.1) that if $\d = \d(n, q, \a) \in (0, 1)$ is sufficiently small, then there exists an integer $m$ with $1 \leq m \leq q$ such that 
$\omega_{n}^{-1}\|V\|(B_{1}^{n+1}(0)) < m + 1/2$ and $m-1/2 \leq \omega_{n}^{-1}2^{n}\|V\|({\mathbf R} \times B_{1/2}) < m + 1/2$. Consider the two alternatives 
\begin{itemize}
\item[(A)] $\int_{{\mathbf R} \times B_{1}} {\rm dist}^{2}(X, {\mathbf P}) \, d\|V\|(X) < \eta^{\prime} \inf_{L \in {\mathcal A}} \int_{{\mathbf R} \times B_{1}} {\rm dist}^{2} \, (X, L) \, d\|V\|(X)$ and
\item[(B)] $\int_{{\mathbf R} \times B_{1}} {\rm dist}^{2}(X, {\mathbf P}) \, d\|V\|(X) \geq \eta^{\prime} \inf_{L \in {\mathcal A}} \int_{{\mathbf R} \times B_{1}} {\rm dist}^{2} \, (X, L) \, d\|V\|(X)$. 
\end{itemize}
In case of alternative (B), choose $\widetilde{L} \in {\mathcal A}$ such that  
$$\int_{{\mathbf R} \times B_{1}} {\rm dist}^{2} \, (X, \widetilde{L}) \, d\|V\|(X) < \frac{3}{2} \inf_{L \in {\mathcal A}} \int_{{\mathbf R} \times B_{1}} {\rm dist}^{2} \, (X, L) \, d\|V\|(X)$$ 
and note, by Theorem~\ref{flat-varifolds}, that if $\d = \d(n, q, \a) \in (0, 1)$ is sufficiently small, then ${\rm dist}_{\mathcal H}^{2} \, (\widetilde{L} \cap ({\mathbf R} \times B_{1}), P \cap ({\mathbf R} \times B_{1}))  \leq C\int_{{\mathbf R} \times B_{1}} {\rm dist}^{2}(X, {\mathbf P}) \, d\|V\|(X)$ for some affine hyperplane $P \subset {\mathbf P},$ where $C = C(n) \in (0, \infty).$ If  now $m=q$ and $\{Z \in {\rm spt} \, \|V\| \cap ({\mathbf R} \times B_{3/4}) \, : \, \Theta \, (\|V\|, Z) \geq q\}  \neq \emptyset$ (in case (B)),  
the assertion with conclusion (b) follows, for sufficiently small $\d = \d(n, q, \a) \in (0, 1)$, by applying 
Theorem~\ref{flat-varifolds} (with $\eta_{1/2 \, \#} \, V$ in place of $V$) and using the estimate (\ref{tilt-ht}) as well as the estimate of the Remark following Theorem~\ref{flat-varifolds}, whereas if $m=q$ and $\{Z \in {\rm spt} \, \|V\| \cap ({\mathbf R} \times B_{3/4}) \, : \, \Theta \, (\|V\|, Z) \geq q\}   = \emptyset$, the assertion with conclusion (a) with $\ell=1$ and $q_{1} = q$ follows from Remark 3 of Section~\ref{outline}, Theorem~\ref{SS} and standard elliptic estimates; if $m \leq q-1,$ 
 hypothesis $(H1)$ implies that conclusion (a) holds. 

In case of alternative (A), we argue by induction on $q$ to see that the assertion with conclusion (a) holds: If $q=2$, the desired conclusion follows directly from Remark (2)(c) above. For general $q,$ let $V_{j}$, ${\mathbf P}_{j}, $ $a_{j}$ be as in Remark 2(c) and note that $a_{j} \leq q-1.$ For each fixed 
$j$, consider the same two alternatives (A) and (B) as above but with $V_{j},$ ${\mathbf P}_{j}$ in place of $V$, ${\mathbf P}$. In case alternative (B) holds (with $V_{j}$, ${\mathbf P}_{j}$ in place of $V$, ${\mathbf P}$), we see by elliptic estimates that conclusion (a) (with $V_{j}$ in place of $V$ and $\ell = 1$) must hold, whereas in case of alternative (A), we may assume by induction the validity of conclusion (a) (with $V_{j}$ in place of $V$ and suitable $\ell_{j}$ in place of $\ell$).

\section{Properties of coarse blow-ups: Part II}\label{step3}
\setcounter{equation}{0}
Fix an integer $q \geq 2$ and suppose that the induction hypotheses $(H1)$, $(H2)$ hold. We begin in this section the proof that the coarse blow-up class ${\mathcal B}_{q}$ satisfies property $({\mathcal B{\emph 7}});$  we shall completed the proof in Section~\ref{propertiesIII}. 

Suppose
\begin{itemize}
\item[(\dag)] $v_{\star}  = (v_{\star}^{1}, v_{\star}^{2}, \ldots, v_{\star}^{q}) \in {\mathcal B}_{q}$ is such that for each $j=1, 2, \ldots, q$, there exist two linear functions $L_{1}^{j}, L_{2}^{j} \, : \, {\mathbf R}^{n} \to {\mathbf R}$ with $L_{1}^{j}(0, y) = L_{2}^{j}(0, y) = 0$ for each $y \in {\mathbf R}^{n-1},$ $v_{\star}^{j}(x^{2}, y) = L_{1}^{j}(x^{2}, y)$ if $x^{2} < 0$ and $v_{\star}^{j}(x^{2}, y) = L_{2}^{j}(x^{2}, y)$ if $x^{2} \geq 0.$ 
\end{itemize}

In order to show that ${\mathcal B}_{q}$ satisfies property $({\mathcal B{\emph 7}})$, we need to prove that 
$v_{\star}^{1} = v_{\star}^{2} = \ldots = v_{\star}^{q} = L$ for some linear function $L \, : \,{\mathbf R}^{n} \to {\mathbf R}.$ We shall do this by establishing the assertions in each of the following two cases:
\begin{itemize}
\item[]{\bf Case 1:} There exists no $v_{\star} \in {\mathcal B}_{q}$ as in ($\dag$) above such that $L_{1}^{1} = L_{1}^{2} = \ldots = L_{1}^{q}$ but $L_{2}^{j} \neq L_{2}^{j +1}$ for some  $j\in \{1, 2, \ldots, q-1\}$.
\item[]{\bf Case 2:} There exists no $v_{\star} \in {\mathcal B}_{q}$ as in ($\dag$) above such that $L_{1}^{i} \neq L_{1}^{i+1}$ for some $i \in \{1, 2, \ldots, q-1\}$ and 
$L_{2}^{j} \neq L_{2}^{j +1}$ for some  $j\in \{1, 2, \ldots, q-1\}$.
\end{itemize}

We prove the assertion of {\bf Case 1} in Lemma~\ref{no-overlap} below and complete the proof that ${\mathcal B}_{q}$ satisfies property $({\mathcal B{\emph 7}})$ (by proving the assertion of {\bf Case 2}) in Corollary~\ref{F7}; the latter requires a number of preliminary results which we shall establish in Sections \ref{fineblowup}-\ref{propertiesIII}.
 
\begin{lemma}\label{no-overlap}
Let $v_{\star}$ and $L_{i}^{j}$, $i \in \{1, 2\}$, $j \in \{1,2, \ldots, q\}$, be as in {\rm (}$\dag${\rm )} above. If $L_{1}^{1} = L_{1}^{2} = \ldots  =L_{1}^{q}$, then {\rm (i)} $L_{2}^{1} = L_{2}^{2} = \ldots = L_{2}^{q}$ and {\rm (ii)} $v_{\star}^{j} = L$ for some linear function $L$ and all $j=1, 2, \ldots, q.$
\end{lemma}

\begin{proof}
The assertion of (ii) follows from that of (i) since the average $(v_{\star})_{a}
 = q^{-1}\sum_{j=1}^{q}v_{\star}^{j}$ is harmonic and hence is a linear function under the hypotheses of the lemma.
 
Suppose, contrary to the assertion of (i), that $L_{2}^{j} \neq L_{2}^{j+1}$ for some $j \in \{1, 2, \ldots, q-1\}.$ By property $({\mathcal B5III}),$ $\frac{v_{\star} - (v_{\star})_{a}}{\|v_{\star} - (v_{\star})_{a}\|} \in {\mathcal B}_{q}$, so we may assume without loss of generality that $L_{1}^{j} =0$ for each $j=1, 2, \ldots, q.$ 
 For $k=1, 2, \ldots,$ let  $V_{k} \in {\mathcal S}_{\a}$ with $(\omega_{n}2^{n})^{-1}\|V_{k}\|(B_{2}^{n+1}(0)) < q + 1/2$, 
$q-1/2 \leq \omega_{n}^{-1}\|V_{k}\|({\mathbf R} \times B_{1}) < q + 1/2$ and ${\hat E}_{k}^{2}  = \int_{{\mathbf R} \times B_{1}} |x^{1}|^{2} d\|V_{k}\|(X)\to 0$ be such that the coarse blow-up of the sequence $V_{k}$, obtained as described in Section~\ref{blow-up}, is $v_{\star}.$ Let the notation be as in Section~\ref{blow-up}. Thus for each $\s \in (0, 1)$ and each sufficiently large $k$ (depending on $\s$), there exists Lipschitz functions $u_{k}^{j} \, : \, B_{\s} \to {\mathbf R},$
$j=1, 2, \ldots, q,$ with ${\rm Lip} \, u_{k}^{j} \leq 1/2$ for each $j \in \{1,2, \ldots, q\}$, such that 
$$v_{\star}^{j} = \lim_{k \to \infty} \, {\hat E}_{k}^{-1}{u}_{k}^{j}$$ 
\noindent
where the convergence is in $L^{2}(B_{\s})$ and weakly in $W^{1,2}(B_{\s})$, and 
\begin{equation}\label{no-overlap-1-1}
{\rm spt} \, \|V_{k}\| \cap \pi^{-1}(B_{\s} \setminus \Sigma_{k}) = \cup_{j=1}^{q} {\rm graph} \, u^{j}_{k} \cap \pi^{-1}(B_{\s} \setminus \Sigma_{k})
\end{equation}
where $\Sigma_{k} \subset B_{\s}$ is the set corresponding to $\Sigma$ in Theorem~\ref{flat-varifolds} when $V$ is replaced with $V_{k}$, so that in particular
\begin{equation}\label{no-overlap-1}
\|V_{k}\|({\mathbf R} \times \Sigma_{k}) + {\mathcal H}^{n} \, (\Sigma_{k}) \leq C{\hat E}^{2}_{k}
\end{equation}
where $C = C(n,q, \s) \in (0, \infty).$ 

In what follows, we take $\s \in [15/16, 1)$ to be fixed. Fix any $\t \in (0, 1/16).$ Since
\begin{eqnarray*}
\int_{({\mathbf R} \times B_{9/16}) \cap \{x^{2} \leq -\t/2\}} |x^{1}|^{2}d\|V_{k}\|(X) &=&
\sum_{j=1}^{q}\int_{(B_{9/16} \setminus \Sigma_{k}) \cap \{x^{2} \leq -\t/2\}} \sqrt{1 + |Du^{j}_{k}|^{2}} \, |u^{j}_{k}|^{2} d{\mathcal H}^{n}\nonumber\\
&&\hspace{.3in}+\int_{({\mathbf R} \times (B_{9/16} \cap \Sigma_{k})) \cap \{x^{2} \leq -\t/2\}} |x^{1}|^{2} d\|V_{k}\|(X),
\end{eqnarray*}
${\hat E}_{k}^{-1}{u}_{k} \to 0$ in $L^{2}$ on $B_{9/16} \cap \{x^{2} \leq -\t/2\}$ and
$$\sup_{X = (x^{1}, x^{\prime}) \in {\rm spt} \, \|V_{k}\| \cap ({\mathbf R} \times B_{9/16})} \, |x^{1}| \to 0,$$
it follows from (\ref{no-overlap-1}) that 
\begin{equation*}
{\hat E}_{k}^{-2}\int_{({\mathbf R} \times B_{9/16}) \cap \{x^{2} \leq -\t/2\}} |x^{1}|^{2}d\|V_{k}\|(X) \to 0
\end{equation*}
and consequently, by (\ref{tilt-ht}), that 
\begin{equation}\label{no-overlap-2}
{\hat E_{k}}^{-2}\int_{({\mathbf R} \times B_{1/2}) \cap \{x^{2} \leq -\t\}} |\nabla^{V_{k}} \, x^{1}|^{2}d\|V_{k}\|(X) \to 0.
\end{equation}

We claim that for all sufficiently large $k$,
\begin{equation}\label{no-overlap-2-1}
\Theta \, (\|V_{k}\|, Z) < q \;\; \mbox{for all} \;\;  Z \in {\rm spt} \, \|V_{k}\|\cap ({\mathbf R} \times B_{5/8}) \cap \{x^{2} > \t/8\}. 
\end{equation}
If this were false, then there would exist a subsequence $\{k^{\prime}\}$ of $\{k\}$ and for each 
$k^{\prime}$, a  point 
$Z_{k^{\prime}} = (z_{k^{\prime}}^{1},z_{k^{\prime}}^{\prime})  \in {\rm spt} \, \|V_{k^{\prime}}\| \cap ({\mathbf R} \times B_{5/8}) \cap \{x^{2} > \t/8\}$
with $\Theta \, (\|V_{k^{\prime}}\|, Z_{k^{\prime}}) \geq q;$ by the reasoning as in the Remark at the end of Section~\ref{step2}, this fact yields
$$\sum_{j=1}^{q}\int_{B_{1/4}(z^{\prime})} R_{z^{\prime}}^{2-n} \left(\frac{\partial \, \left((v_{\star}^{j} - y)/R_{z^{\prime}}\right)}{\partial \, R_{z^{\prime}}}\right)^{2} d{\mathcal H}^{n} \leq C$$
for some $z^{\prime} \in {\overline B}_{5/8} \cap \{x^{2} \geq \t/8\}$ and some $y \in {\mathbf R},$ which  
implies that $v_{\star}^{j}(z^{\prime}) = y$ for all $j=1, 2, \ldots, q.$ But this contradicts our hypothesis
that $L_{2}^{j} \neq L_{2}^{j+1}$ for some $j \in \{1, 2, \ldots, q-1\},$ so (\ref{no-overlap-2-1}) must 
hold for all sufficiently large $k$.

With the help of Remark 3 of Section~\ref{outline}, we deduce from (\ref{no-overlap-2-1}) that for all sufficiently large $k$,  ${\mathcal H}^{n-7 +\g} \, ({\rm sing} \, V_{k} \cap ({\mathbf R} \times B_{5/8}) \cap 
\{x^{2} > \t/8\}) = 0$ for each $\g >0$ if $n \geq 7$ and ${\rm sing} \, V_{k} \cap ({\mathbf R} \times B_{5/8}) \cap \{x^{2} > \t/8\} = \emptyset$ if $2 \leq n \leq 6.$ We may therefore apply Theorem~\ref{SS} and elliptic theory to deduce that, for all sufficiently large $k$, $\Sigma_{k} \cap B_{9/16} \cap \{x^{2} > \t/4\} = \emptyset$; 
\begin{equation}\label{no-overlap-3}
 V_{k} \res(({\mathbf R} \times B_{9/16}) \cap \{x^{2} > \t/4\}) =   \sum_{j=1}^{q} |{\rm graph} \, u_{k}^{j}| \res 
(({\mathbf R} \times B_{9/16}) \cap \{x^{2} > \t/4\});
\end{equation}
and that $u_{k}^{j}$ are $C^{2}$ on $B_{9/16} \cap \{x^{2} > \t/4\}$, solve the minimal surface equation there and satisfy
\begin{equation}\label{no-overlap-4}
\sup_{B_{1/2}\cap \{x^{2} >\t/4\}} |D^{\ell} \, u_{k}|^{2} \leq C_{\t}{\hat E}_{k}^{2}
\end{equation}
for $\ell = 0, 1, 2$, where $C_{\t}$ is a constant depending only on $n$ and $\t$, 
and $D^{\ell}$ denotes the order $\ell$ differentiation.

We next claim that for all sufficiently large $k$,
\begin{equation}\label{no-overlap-4-1-0}
(\{0\} \times {\mathbf R}^{n-1}) \cap B_{1/2}  \subset \left(\{Z  \in {\rm spt} \, \|V_{k}\|\, : \, 
\Theta \, (\|V_{k}\|, Z) \geq q\}\right)_{\t}.
\end{equation}
If this were false, then there would exist a point $(0, y) \in \{0\} \times {\mathbf R}^{n-1} \cap {\overline B}_{1/2}$ and a subsequence $\{k^{\prime}\}$ of $\{k\}$ such that for each $k^{\prime}$, 
$$B^{n+1}_{3\t/4}((0, y)) \cap \{Z  \in {\rm spt} \, \|V_{k^{\prime}}\|\, : \, 
\Theta \, (\|V_{k^{\prime}}\|, Z) \geq q\} = \emptyset.$$
Since ${\rm spt} \, \|V_{k}\| \cap ({\mathbf R} \times B_{3/4}) \to \{0\} \times B_{3/4}$ in Hausdorff distance, it follows that for each $k^{\prime}$ and each 
$Z \in {\rm spt} \, \|V_{k^{\prime}}\| \cap ({\mathbf R} \times B_{\t/2}((0, y)))$ we must have 
$\Theta \, (\|V_{k^{\prime}}\|, Z) < q$. Arguing exactly as for (\ref{no-overlap-3}) and (\ref{no-overlap-4}), we conclude that for all sufficiently large $k^{\prime}$, $\Sigma_{k^{\prime}} \cap B_{\t/4}(0, y) = \emptyset$;
$${\rm spt} \, \|V_{k^{\prime}}\| \cap ({\mathbf R} \times B_{\t/4}(0, y)) = 
\cup_{j=1}^{q} {\rm graph} \, \left.u_{k^{\prime}}^{j}\right|_{B_{\t/4}(0, y)};$$ 
and that $u_{k^{\prime}}^{j}$ are $C^{2}$ functions on $B_{\t/4}(0, y)$, satisfy 
$$\sum_{j=1}^{q} \sup_{B_{\t/4}(0, y)} \, |Du_{k^{\prime}}^{j}| + |D^{2}u_{k^{\prime}}^{j}| \leq C{\hat E}_{k^{\prime}}, \;\;\; C = C(n, \t) \in (0, \infty)$$ 
and solve the minimal surface equation on $B_{\t/4}(0, y).$ Consequently, 
$v_{\star}^{j}|_{B_{\t/4}(0,y)}$ must be harmonic for each $j = 1, 2, \ldots, q$, which is however impossible since by hypothesis, 
$L_{1}^{j} = 0$ for each $j = 1, 2, \ldots, q$ while $L_{2}^{j} \neq L_{2}^{j+1}$ for some $j \in \{1, 2, \ldots, q-1\}.$ This contradiction establishes (\ref{no-overlap-4-1-0}) for all sufficiently large $k$.

We now proceed to derive the contradiction needed for the proof of the lemma. By taking $\psi(X) = \widetilde{\z}(X)e^{2}$ in the first variation formula (\ref{firstvariation}), we deduce that 
\begin{equation}\label{no-overlap-main}
\int \nabla^{V_{k}} \, x^{2} \cdot \nabla^{V_{k}} \,{\widetilde \z}(X) d\|V_{k}\|(X)= 0
\end{equation}
for each $k=1, 2, \ldots$ and each ${\widetilde \z} \in C^{1}_{c}({\mathbf R} \times B_{1}).$ Choosing ${\widetilde \z}$ 
to agree with $\z^{\prime}(x^{1}, x^{\prime}) = \z(x^{\prime})$ in a neighborhood of ${\rm spt} \,\|V\| \cap ({\mathbf R} \times B_{1/4})$, where 
$\z  \in C^{1}_{c}(B_{1/4})$ is arbitrary, we deduce from this that 
\begin{equation}\label{no-overlap-4-1}
\sum_{j=1}^{q}\int_{B_{1/4}}\sqrt{1 + |D{u}_{k}^{j}|^{2}}\left(D_{2}\z - \frac{D_{2}{u}_{k}^{j} (D\z \cdot D{u}_{k}^{j})}{1 + |D{u}_{k}^{j}|^{2}}\right) = F_{k}, \;\;\; \mbox{where}
\end{equation}
\begin{eqnarray*}
F_{k} = -\int_{{\mathbf R} \times (B_{1/4} \cap \Sigma_{k})} \nabla^{V_{k}} \, x^{2} \cdot \nabla^{V_{k}} \,{\widetilde \z}(X) d\|V_{k}\|(X) &&\\
&&\hspace{-3in}+ 
\sum_{j=1}^{q}\int_{B_{1/4} \cap \Sigma_{k}}\sqrt{1 + |D{u}_{k}^{j}|^{2}}\left(D_{2}\z - \frac{D_{2}{u}_{k}^{j} (D\z \cdot D{u}_{k}^{j})}{1 + |D{u}_{k}^{j}|^{2}}\right) 
\end{eqnarray*}
Since $\int_{B_{1/4}} D_{2} \z = 0$, it follows from (\ref{no-overlap-4-1}) that 
\begin{equation}\label{no-overlap-4-2}
\sum_{j=1}^{q}\int_{B_{1/4}}\frac{|D{u}_{k}^{j}|^{2}}{1 + \sqrt{1 + |D{u}_{k}^{j}|^{2}}} D_{2}\z - \frac{D_{2}{u}_{k}^{j} (D\z \cdot D{u}_{k}^{j})}{\sqrt{1 + |D{u}_{k}^{j}|^{2}}} = F_{k}.
\end{equation}
In view of  (\ref{no-overlap-3}) and (\ref{no-overlap-4}), it follows from  the definition of $\Sigma_{k}$ (see Theorem~\ref{flat-varifolds}) that
\begin{equation}\label{no-overlap-4-2-1}
B_{1/4} \cap \Sigma_{k} \subset B_{1/4} \cap \{x^{2} < \t/2\}.
\end{equation}
We claim also that for all sufficiently large $k$,
\begin{equation}\label{no-overlap-4-3}
\|V_{k}\|({\mathbf R} \times (B_{1/4} \cap \Sigma_{k})) + {\mathcal H}^{n}(B_{1/4} \cap  \Sigma_{k}) 
\leq C \int_{({\mathbf R} \times B_{1/2}) \cap \{x^{2} < \t\}} |\nabla^{V_{k}} \, x^{1}|^{2} d\|V_{k}\|(X)
\end{equation}
where $C \in (0, \infty)$ is a fixed constant depending only on $n$ and $q$. To see this, let $\widetilde{\Sigma}_{k}^{(j)},$ $j=1, 2, 3$, correspond to the set $\widetilde{\Sigma}_{j}$ in Theorem~\ref{flat-varifolds} when $V$ is replaced by $V_{k}$, and let $\Sigma_{k}^{\prime}$ correspond to $\Sigma^{\prime}.$ Since for
each $k$, $\r \in (\t/4, 1/16)$ and $Y \in {\rm spt} \, \|V_{k}\| \cap 
({\mathbf R} \times B_{1/2})$, we trivially have that 
\begin{equation}
\r^{-n}\int_{{\mathbf R} \times B_{\r}(\pi \, Y)} |\nabla^{V_{k}} \, x^{1}|^{2} d\|V_{k}\|(X) \leq 4^{n}\t^{-n} \int_{{\mathbf R} \times B_{3/4}} |\nabla^{V_{k}} \, x^{1}|^{2} d\|V_{k}\|(X) \leq 4^{n}\t^{-n}{\hat E}_{k}^{2},
\end{equation}
and since by definition,  
\begin{eqnarray*}
{\widetilde \Sigma}_{k}^{(1)} = \{Y \in {\rm spt} \, \|V_{k}\| \cap ({\mathbf R} \times B_{\s}) \, : \, 
\r^{-n}\int_{{\mathbf R} \times B_{\r}(\pi \, Y)} |\nabla^{V_{k}} \, x^{1}|^{2} d\|V_{k}\|(X) \geq \xi&&\\
&&\hspace{-1in} \;\; \mbox{for some} \;\;
\r \in (0, (1-\s))\}
\end{eqnarray*}
where $\xi = \xi(n, q) \in (0, 1/2)$ is as in Theorem~\ref{flat-varifolds}, it follows that for all sufficiently large $k$ (depending on $\t$), 
$Y \in \widetilde{\Sigma}_{k}^{(1)}$ if and only if $Y \in {\rm spt} \, \|V_{k}\| \cap ({\mathbf R} \times B_{\s})$ and $\r^{-n}\int_{{\mathbf R} \times B_{\r}(\pi \, Y)} |\nabla^{V_{k}} \, x^{1}|^{2} d\|V_{k}\|(X) \geq \xi$ for some 
$\r \in (0, \t/4].$ Also, by Part 3 of the proof of [\cite{A}, Theorem~3.8], 
we have that for each $x \in B_{\s}$ and each $k$, 
\begin{equation}\label{no-overlap-4-3-1}
\sum_{Y \in {\rm spt} \, \|V_{k}\| \cap \pi^{-1}(x) \setminus \left({\widetilde \Sigma}_{k}^{(1)} \cup {\widetilde \Sigma}_{k}^{(2)}\right)}\Theta\left(\|V_{k}\|, Y\right) \leq q. 
\end{equation}
In view of (\ref{no-overlap-4-2-1}), it follows from the Besicovitch covering lemma and (\ref{no-overlap-4-3-1}) that
\begin{equation*}
\|V_{k}\|({\mathbf R} \times (B_{1/4} \cap \pi \, \widetilde{\Sigma}_{k}^{(j)})) + {\mathcal H}^{n}(B_{1/4} \cap  \pi \, \widetilde{\Sigma}_{k}^{(j)}) 
\leq C \int_{(B_{1/2} \times {\mathbf R}) \cap \{x^{2} < \t\}} |\nabla^{V_{k}} \, x^{1}|^{2} d\|V_{k}\|(X)
\end{equation*}
for $j=1$, where $C = C(n, q) \in (0, \infty).$ Since $\|V_{k}\|(\widetilde{\Sigma}_{k}^{(2)}) = 0$ (see Part 2 of the proof of [\cite{A}, Theorem 3.8]), this estimate also follows for $j=2$ in view of (\ref{no-overlap-4-3-1}); it follows for $j=3$, directly from the definition of $\widetilde{\Sigma}_{k}^{(3)}$, (\ref{no-overlap-4-2-1}) and (\ref{no-overlap-4-3-1}); it also holds with $\Sigma_{k}^{\prime}$ in place of $\pi \, \widetilde{\Sigma}_{k}^{(j)}$, by (\ref{no-overlap-4-2-1}) and Part 5 of the proof of 
[\cite{A}, Theorem~3.8].  Thus the estimate (\ref{no-overlap-4-3}), with 
the constant $C$ depending only on $n$ and $q$ (in particular independent of $\t$), holds.

By (\ref{no-overlap-4-3}), Theorem~\ref{non-concentration}(b) (with $\m  = 1/2$) 
and (\ref{no-overlap-2}) we deduce that, since the integrands in both integral
expressions in  $F_{k}$ are bounded, 
\begin{equation}\label{no-overlap-5}
{\hat E}_{k}^{-2}|F_{k}| \leq C \sup \, |D\z| \t^{1/2}
\end{equation}
for all sufficiently large $k$, where $C = C(n, q) \in (0, \infty).$

Abbreviating $w_{k} = \sum_{j=1}^{q}\frac{|D{u}_{k}^{j}|^{2}}{1 + \sqrt{1 + |D{u}_{k}^{j}|^{2}}} D_{2}\z - \frac{D_{2}{u}_{k}^{j} (D\z \cdot D{u}_{k}^{j})}{\sqrt{1 + |D{u}_{k}^{j}|^{2}}}$, note that  
$$\int_{B_{1/4} \setminus \Sigma_{k} \cap \{x^{2} \leq \t\}} |w_{k}| \leq C \sup \, |D\z| \int_{({\mathbf R} \times B_{1/2}) \cap \{x^{2} \leq \t\}} |\nabla^{V_{k}} \, x^{1}|^{2} d\|V_{k}\|(X),$$
and by (\ref{no-overlap-4-3}),
$$ \int_{B_{1/4} \cap \Sigma_{k}} |w_{k}|  \leq C \sup \, |D\z| \int_{({\mathbf R} \times B_{1/2}) \cap \{x^{2} \leq \t\}} |\nabla^{V_{k}} \, x^{1}|^{2} d\|V_{k}\|(X)$$
where $C = C(n)$, so that again by Theorem~\ref{non-concentration}(b) with $\m=1/2$ and (\ref{no-overlap-2}), 
\begin{equation}\label{no-overlap-6}
{\hat E_{k}}^{-2} \left(\int_{B_{1/4} \setminus \Sigma_{k} \cap \{x^{2} \leq \t\}} |w_{k}| + \int_{B_{1/4} \cap  \Sigma_{k}} |w_{k}|\right) \leq C\sup \, |D\z|\t^{1/2}
\end{equation}
for all sufficiently large $k$, where $C = C(n)$. Finally, by (\ref{no-overlap-4}), 
\begin{equation}\label{no-overlap-7}
\lim_{k \to \infty} \, {\hat E_{k}}^{-2} \int_{B_{1/4} \cap \{x^{2} \geq \t\}} w_{k} = -\frac{1}{2}\sum_{j=1}^{q}\int_{B_{1/4} \cap \{x^{2} \geq \t\}} |D_{2}v_{\star}^{j}|^{2}D_{2} \z
\end{equation}
where we have used the fact that $D_{i} v_{\star}^{j} \equiv 0$ for $i=3, \ldots, (n+1)$ and 
$j= 1, 2, \ldots, q.$ Dividing (\ref{no-overlap-4-2}) by ${\hat E}_{k}^{2}$ and first letting $k \to \infty$ and 
then letting $\t \to 0$, we conclude from (\ref{no-overlap-5}), (\ref{no-overlap-6}) and 
(\ref{no-overlap-7}) that 
$$\sum_{j=1}^{q}\int_{B_{1/4} \cap \{x^{2} \geq 0\}} |D_{2}v^{j}_{\star}|^{2} D_{2}\z = 0$$
for any $\z \in C^{1}_{c}(B_{1/4}).$ Since $v^{j}_{\star} = L_{2}^{j}$ on $\{x^{2} \geq 0\}$, this contradicts (for any choice of $\z \in C^{1}_{c}(B_{1/4})$ with 
$\int_{B_{1/4} \cap \{x^{2} \geq 0\}} D_{2}\z \neq 0$) our assumption that  $L_{2}^{j} \neq L_{2}^{j+1}$ 
for some $j \in \{1, 2, \ldots, q-1\}$.
\end{proof}

\noindent
{\bf Remark}: It follows from Lemma~\ref{no-overlap} and the compactness property $({\mathcal B{\emph 6}})$ that there exists a constant $c = c(n, q) \in (0, \infty)$ with the following property:   
If  $v \in {\mathcal B}_{q}$ is such that, for each $j= 1, 2, \ldots, q$, 
$v^{j}(x^{2}, y) =  \ell_{j}x^{2}$ for $x^{2} < 0$; 
$v^{j}(x^{2}, y) = m_{j}x^{2}$ for $x^{2} \geq 0,$ where $\ell_{j}$, $m_{j}$ are constants; and $v^{j} \not\equiv v_{a}$ for some $j \in \{1, 2, \ldots, q\}$, where $v_{a}  \equiv q^{-1}\sum_{j=1}^{q}v^{j}$, then 
$|\ell_{1} - \ell_{q}|^{2} \geq c\sum_{j=1}^{q}\|v^{j} - v_{a}\|^{2}_{L^{2}(B_{1})}$ and $|m_{1} - m_{q}|^{2} \geq c\sum_{j=1}^{q}\|v^{j} - v_{a}\|^{2}_{L^{2}(B_{1})}.$ (Of course once we have completed the proof that ${\mathcal B}_{q}$ satisfies property 
$({\mathcal B{\emph 7}})$, we will have ruled out the existence of such $v \in {\mathcal B}_{q}.$)

\section{Parametric $L^{2}$-estimates in terms of  fine excess}\label{fineblowup}
\setcounter{equation}{0}

This section and all of the subsequent sections up to and including Section~\ref{propertiesIII} will be devoted to the proof of the assertion of {\bf Case 2} set forth at the beginning of Section~\ref{step3}. Crucial to our proof are the $L^{2}$-estimates, given in Theorem~\ref{L2-est-1} and Corollary~\ref{L2-est-2} below, for a varifold $V \in {\mathcal S}_{\a}$ with small coarse excess (relative to a hyperplane) and  lower order ``fine excess'' relative to an appropriate union of half-hyperplanes meeting along  an $(n-1)$-dimensional axis (see Hypotheses~\ref{hyp}(5) below). These results are adaptations to the present ``higher multiplicity'' setting of those proved in (\cite{S}) in the context of ``multiplicity 1  classes''  of minimal submanifolds.

\noindent
{\bf Notation}: {\bf (1)} Let ${\mathcal C}_{q}$ denote the set of hypercones ${\mathbf C}$ of ${\mathbf R}^{n+1}$ such that 
${\mathbf C} = \sum_{j=1}^{q} |H_{j}| + |G_{j}|,$ where for each $j \in \{1, 2, \ldots, q\}$, $H_{j}$ is the half-hyperplane defined by 
$$H_{j} = \{(x^{1}, x^{2}, y) \in {\mathbf R}^{n+1} \, : \, x^{2} < 0 \;\; \mbox{and} \;\; x^{1} = \lambda_{j}x^{2}\},$$ 
$G_{j}$ the half-hyperplane defined by 
$$G_{j} = \{(x^{1},x^{2}, y) \in {\mathbf R}^{n+1}\, : \, x^{2} > 0 \;\; \mbox{and} \;\; x^{1} = \mu_{j}x^{2}\},$$ with $\lambda_{j}, \mu_{j}$ constants,  
$\lambda_{1} \geq \lambda_{2} \geq \ldots \geq \lambda_{q}$ and $\mu_{1} \leq \mu_{2} \leq \ldots \leq \mu_{q}.$ 
Note that we do \emph{not} assume cones in ${\mathcal C}_{q}$ are stationary in ${\mathbf R}^{n+1}.$

\noindent
{\bf (2)} For $p \in \{2, 3, \ldots, 2q\}$, let ${\mathcal C}_{q}(p)$ denote the set of  hypercones ${\mathbf C}  = \sum_{j=1}^{q} |H_{j}| + |G_{j}| \in {\mathcal C}_{q}$ as defined above such that the number of {\em distinct} half-hyperplanes in the set $\{H_{1}, \ldots, H_{q}, G_{1}, \ldots, G_{q}\}$ is $p$. Then ${\mathcal C}_{q} = \cup_{p=2}^{2q} \, {\mathcal C}_{q}(p).$

\noindent
{\bf (3)} For $V \in {\mathcal S}_{\a}$ and ${\mathbf C} \in {\mathcal C}_{q}$ define a height excess (``fine excess'') $Q_{V}({\mathbf C})$ of $V$ relative to ${\mathbf C}$ by
\begin{eqnarray*}
Q_{V}({\mathbf C})  = \left(\int_{{\mathbf R} \times (B_{1/2} \setminus \{|x^{2}| < 1/16\})} {\rm dist}^{2}(X, {\rm spt} \, \|V\|) \,d\|{\mathbf C}\|(X)\right.&&\nonumber\\ 
&&\hspace{-2in}+ \; \left.\int_{{\mathbf R} \times B_{1}} {\rm dist}^{2} \, (X, {\rm spt}\, \|{\mathbf C}\|) \, d\|V\|(X)\right)^{1/2}. 
\end{eqnarray*}

\noindent
{\bf (4)} For $q \geq 2$ and $p \in \{4, \ldots, 2q\}$, let 
$$Q_{V}^{\star}(p) = \inf_{{\mathbf C} \in \cup_{k=4}^{p} {\mathcal C}_{q}(k)} \, Q_{V}({\mathbf C}).$$

Let $\a \in (0, 1)$ and $q$ be an integer $\geq 2.$ In Theorem~\ref{L2-est-1}, Corollary~\ref{L2-est-2} and Lemma~\ref{no-gaps} below and subsequently, we shall consider the following set of hypotheses for appropriately small $\e, \g \in (0, 1)$ to be determined depending only on $n$, $q$ and $\a$:

\begin{hypotheses}\label{hyp}
\begin{itemize}
\item[]
\item[(1)] $V \in {\mathcal S}_{\a}$, \; $\Theta \, (\|V\|, 0) \geq q$, \;$(\omega_{n}2^{n})^{-1}\|V\|(B_{2}^{n+1}(0)) < q + 1/2$, \; $\omega_{n}^{-1}\|V\|({\mathbf R} \times B_{1}) < q + 1/2.$ 
\item[(2)] ${\mathbf C} = \sum_{j=1}^{q} |H_{j}| + |G_{j}| \in {\mathcal C}_{q},$ where for each $j \in \{1, 2, \ldots, q\}$, $H_{j}$ is the half-hyperplane defined by
$H_{j} = \{(x^{1}, x^{2}, y) \in {\mathbf R}^{n+1} \, : \, x^{2} < 0 \;\; \mbox{and} \;\; x^{1} = \lambda_{j}x^{2}\},$ 
$G_{j}$ the half-hyperplane defined by $G_{j} = \{(x^{1},x^{2}, y) \in {\mathbf R}^{n+1}\, : \, x^{2} > 0 \;\; \mbox{and} \;\; x^{1} = \mu_{j}x^{2}\},$ with $\lambda_{j}, \mu_{j}$ constants,  
$\lambda_{1} \geq \lambda_{2} \geq \ldots \geq \lambda_{q}$ and $\mu_{1} \leq \mu_{2} \leq \ldots \leq \mu_{q}$. 
\item[(3)] ${\hat E}_{V}^{2} \equiv \int_{{\mathbf R} \times B_{1}} |x^{1}|^{2} d\|V\|(X) < \e.$
\item[(4)] $\{Z \, : \, \Theta \, (\|V\|, Z) \geq q\} \cap \left({\mathbf R} \times (B_{1/2} \setminus \{|x^{2}| < 1/16\}) \right) = \emptyset.$
\item[(5)] $Q_{V}^{2}({\mathbf C}) < \g {\hat E}_{V}^{2}$.
\end{itemize} 
\end{hypotheses}

\noindent
{\bf Remark:} \emph{There exists $\e = \e(n, q) \in (0, 1)$ such that  if Hypotheses~\ref{hyp} above hold with any $\g \in (0, 1),$ and the induction hypotheses $(H1)$, $(H2)$ hold, then 
\begin{equation}\label{slope-bounds-1}
{\rm max} \, \{|\lambda_{1}|, |\lambda_{q}|\} \leq c_{1} {\hat E}_{V} \;\;\; {\rm and} \;\;\; {\rm max} \, \{|\mu_{1}|, |\mu_{q}|\} \leq c_{1} {\hat E}_{V}
\end{equation}
where $c_{1} = c_{1}(n) \in (0, \infty).$}  These bounds follow from Hypotheses~\ref{hyp}(5) in view of the fact that (by Hypotheses~\ref{hyp}(4), Remark 3 of Section~\ref{outline} and Theorem~\ref{SS}), under Hypotheses~\ref{hyp}, 
$V \res ({\mathbf R} \times (B_{1/4} \setminus \{|x^{2}| < 1/8\})) = \sum_{j=1}^{q} |{\rm graph} \, \widetilde{u}_{j}| + |{\rm graph} \, \widetilde{w}_{j}|$ where, for $j=1, 2, \ldots, q$, 
$\widetilde{u}_{j} \in C^{2}(B_{1/4} \cap \{x^{2} < -1/8\})$, $\widetilde{w}_{j} \in C^{2}(B_{1/4} \cap \{x^{2} > 1/8\})$ with $\sup_{B_{1/4} \cap \{x^{2} < -1/8\}} \, |\widetilde{u}_{j}| \leq C {\hat E}_{V}$ and $\sup_{B_{1/4} \cap \{x^{2} > 1/8\}} \, |\widetilde{w}_{j}| \leq C {\hat E}_{V},$ $C = C(n) \in (0, \infty).$

Let $c_{1} = c_{1}(n)$ be the constant as in (\ref{slope-bounds-1}) above and define a constant $M_{0} = M_{0}(n, q) \in (0, \infty)$ by 
$$M_{0} = \max \, \left\{\frac{3}{2}, \frac{2^{2n+8}\omega_{n}^{2}(2q+1)^{2}c_{1}^{2}}{\overline{C}_{1}}, \frac{2^{2n+8}\omega_{n}(2q+1)}{\overline{C}_{1}}\right\}$$ where $\overline{C}_{1} = \int_{B_{1/2} \cap \{x^{2} >1/16\}} |x^{2}|^{2} \, d{\mathcal H}^{n}(x^{2}, y)$. We shall use this constant 
at several places below.

For $V$ as in Hypotheses~\ref{hyp}, we shall also assume the following for suitable values of $M > 1$:

\noindent
{\bf Hypothesis ($\star$)}: 
\begin{equation*}
{\hat E}_{V}^{2} < M\inf_{\{P = \{x^{1} = \lambda x^{2}\} \in G_{n} \, : \, \lambda \in {\mathbf R}\}} \, \int_{{\mathbf R} \times B_{1}} {\rm dist}^{2} \, (X, P) , d\|V\|(X).
\end{equation*}

\noindent
{\bf Remarks}: {\bf (1)} \emph{If Hypotheses~\ref{hyp} and Hypothesis ($\star$) hold with sufficiently small $\e = \e(n, q) \in (0, 1)$, $\g =\g(n, q) \in (0, 1)$ and with $M = \frac{3}{2}M_{0}^{4}$, then 
\begin{eqnarray}\label{slope-bounds-2}
c {\hat E}_{V}  \leq {\rm max} \, \{|\lambda_{1}|, |\lambda_{q}|\}, \;\;\;\; c{\hat E}_{V}  \leq {\rm max} \, \{|\mu_{1}|, |\mu_{q}|\} \;\;\; {\rm and}&&\nonumber\\
&&\hspace{-3.25in} {\rm min} \, \{|\lambda_{1} - \lambda_{q}|, |\mu_{1} - \mu_{q}|\} \geq 2c {\hat E}_{V}
\end{eqnarray}
for some constant $c = c(n, q) \in (0, \infty)$}. Indeed, the triangle inequality (in the form 
${\rm dist}^{2} \, (X, P) \leq 2{\rm dist}^{2} \, (X, {\rm spt} \, \|{\mathbf C}\|) + 2{\rm dist}^{2}_{\mathcal H} \, (P \cap ({\mathbf R} \times B_{1}), {\rm spt} \, \|{\mathbf C}\| \cap ({\mathbf R} \times B_{1}))$ for $X \in {\mathbf R} \times B_{1},$ applied with $P = \{x^{1}  = \frac{1}{2}(\lambda_{1}+\lambda_{q})x^{2}\}$ or $P = \{x^{1} = \frac{1}{2}(\m_{1}+\m_{2})x^{2}\}$), Hypothesis~($\star$) (with $M = \frac{3}{2}M_{0}^{4}$) and Hypotheses \ref{hyp} (with sufficiently small $\e = \e(n, q) \in (0, 1)$ and $\g = \g(n, q) \in (0, 1)$) imply that $|\lambda_{1} - \lambda_{q}| + |\m_{1} - \m_{q}| \geq \widetilde{c}{\hat E}_{V}$ for some $\widetilde{c} = \widetilde{c}(n, q) \in (0, \infty).$ Lemma~\ref{no-overlap} then implies that 
${\rm min} \, \{|\lambda_{1} - \lambda_{q}|, |\mu_{1} - \mu_{q}|\} \geq 2c {\hat E}_{V}$, $c = c(n, q) \in (0, 1)$; the first two inequalities of (\ref{slope-bounds-2}) follow readily from this.
 
\noindent
{\bf (2)} It follows from the last inequality of (\ref{slope-bounds-2}) that if Hypotheses~\ref{hyp} and 
Hypothesis~($\star$) hold with $\e = \e(n, q),$ $\g = \g(n, q) \in (0, 1)$ sufficiently small and 
$M= \frac{3}{2}M_{0}^{4}$, then ${\mathbf C} \in {\mathcal C}_{q}(p)$ for some $p \in \{4, 5, \ldots, 2q\}.$

Finally, for ${\mathbf C}$, $V$  as in Hypotheses~\ref{hyp} and appropriately small $\b \in (0, 1/2)$ (to be determined depending only on $n$, $q$ and $\a$), we will also need to consider the following:

\noindent
{\bf Hypothesis ($\star\star$)}: \emph{Either
\begin{itemize}
\item[(i)] ${\mathbf C} \in {\mathcal C}_{q}(4)$ or
\item[(ii)] $q \geq 3$, ${\mathbf C} \in {\mathcal C}_{q}(p)$ for some $p \in \{5, \ldots, 2q\}$ and 
$Q_{V}^{2}({\mathbf C}) < \b \left(Q_{V}^{\star}(p-1)\right)^{2}.$
\end{itemize}}

\noindent
{\bf Remarks:} {\bf (1)} Let ${\mathbf C}$ be as in Hypothesis~\ref{hyp}(2). If $V \in {\mathcal S}_{\a},$ ${\mathbf C}$ satisfy Hypothesis~\ref{hyp}(1), Hypothesis~($\star\star$)(ii) with $\b \in (0, 1/4)$ and if 
$\lambda_{1} = \lambda_{1}^{\prime}> \lambda_{2}^{\prime}> \ldots >\lambda_{p_{1}}^{\prime} = \lambda_{q}$ are the distinct elements of the set
$\{\lambda_{1}, \ldots, \lambda_{q}\}$ and $\mu_{1} = \mu_{1}^{\prime} < \mu_{2}^{\prime} < \ldots < \mu_{p_{2}}^{\prime} = \mu_{q}$ are the distinct elements of $\{\mu_{1}, \ldots, \mu_{q}\}$ (notation as in Hypothesis~\ref{hyp}(2)),
then it follows from Hypothesis~($\star\star$) and the triangle inequality that  
\begin{equation}\label{separation}
\lambda_{i+1}^{\prime} - \lambda_{i}^{\prime} \geq 2c^{\prime}  Q_{V}^{\star}(p-1), \;\;\; \mu_{j+1}^{\prime} - \mu_{j}^{\prime} \geq 2c^{\prime} Q_{V}^{\star}(p-1)
\end{equation}
for some constant $c^{\prime} = c^{\prime}(n, q) \in (0, \infty)$ and all $i=1, 2, \ldots, p_{1}-1$ and $j=1, 2, \ldots, p_{2}-1.$

\noindent
{\bf (2)} Suppose $V \in {\mathcal S}_{\a}$, ${\mathbf C} \in {\mathcal C}_{q}$ satisfy Hypotheses~\ref{hyp}, Hypothesis~(${\star}$) and Hypothesis (${\star\star}$) for some $\e, \g, \b \in (0, 1/2).$ If ${\mathbf C}^{\prime} \in {\mathcal C}_{q}$ is any other cone with 
${\rm spt} \, \|{\mathbf C}^{\prime}\| = {\rm spt} \, \|{\mathbf C}\|,$ then Hypotheses~\ref{hyp}, Hypothesis~($\star$) and Hypothesis~($\star\star$) will continue to be satisfied with ${\mathbf C}^{\prime}$ in place of ${\mathbf C}$ provided $\g$, $\b$ are replaced by $2q\g$, $2q\b$ respectively.

\begin{theorem}\label{L2-est-1}
Let $q$ be an integer $\geq 2,$ $\a \in (0, 1)$, $\t \in (0, 1/8)$ and $\m \in (0, 1)$. There exist numbers $\e_{0} = \e_{0}(n , q, \a,\t) \in (0,1),$ 
$\g_{0} = \g_{0}(n, q,\a, \t) \in (0, 1)$ and $\b_{0} = \b_{0}(n, q, \a,\t) \in (0, 1)$ such that the following is true:
Let $V \in {\mathcal S}_{\a}$, ${\mathbf C} \in {\mathcal C}_{q}$ satisfy
Hypotheses~\ref{hyp}, Hypothesis~($\star$) and Hypothesis~($\star\star$) with $M = \frac{3}{2}M_{0}^{4}$ and $\e_{0}$,$\g_{0}$, $\b_{0}$ in place of $\e$, $\g$, $\b$ respectively. Suppose also that the induction hypotheses $(H1)$, $(H2)$ hold. 
Write ${\mathbf C} = \sum_{j=1}^{q}|H_{j}| + |G_{j}|$ where for each $j \in \{1, 2, \ldots, q\}$, $H_{j}$ is the half-space defined by
$H_{j} = \{(x^{1}, x^{2}, y) \in {\mathbf R}^{n+1} \, : \, x^{2} < 0 \;\; \mbox{and} \;\; x^{1} = \lambda_{j}x^{2}\},$ 
$G_{j}$ the half-space defined by $G_{j} = \{(x^{1},x^{2}, y) \in {\mathbf R}^{n+1}\, : \, x^{2} > 0 \;\; \mbox{and} \;\; x^{1} = \mu_{j}x^{2}\},$ with $\lambda_{j}, \mu_{j}$ constants,  
$\lambda_{1} \geq \lambda_{2} \geq \ldots \geq \lambda_{q}$ and $\mu_{1} \leq \mu_{2} \leq \ldots \leq \mu_{q}$; for 
$(x^{2}, y) \in {\mathbf R}^{n}$ and $j = 1, 2, \ldots, q$, define $h_{j}(x^{2}, y) = \lambda_{j}x^{2}$ and 
$g_{j}(x^{2}, y) = \m_{j}x^{2}.$ Then, 
after possibly replacing ${\mathbf C}$ with another cone ${\mathbf C}^{\prime} \in {\mathcal C}_{q}$ with 
${\rm spt} \, \|{\mathbf C}^{\prime}\| = {\rm spt} \, \|{\mathbf C}\|$ and relabelling ${\mathbf C}^{\prime}$ as ${\mathbf C}$ (see the preceding Remark (2)), the following must hold: 
\begin{eqnarray*}
&&\hspace{-.2in}({\rm a}) \;\;\;V \res ({\mathbf R} \times (B_{3/4} \setminus \{|x^{2}| < \t\})) = \sum_{j=1}^{q} |{\rm graph} \, (h_{j} + u_{j})| + |{\rm graph} \, (g_{j} + w_{j})|\nonumber\\ 
&&\hspace{-.2in}\mbox{where, for each $j=1, 2, \ldots, q,$ $u_{j} \in C^{2} \, (B_{3/4} \,\cap \,\{x^{2} < -\t\});$ $w_{j} \in C^{2} \, (B_{3/4} \, \cap \, \{x^{2} > \t\});$}\nonumber\\  
&&\hspace{-.2in}\mbox{$h_{j} + u_{j}$ and $g_{j} + w_{j}$ solve the minimal surface equation on their respective domains;}\nonumber\\  
&&\hspace{-.2in}h_{1} + u_{1} \leq h_{2} + u_{2} \leq \ldots \leq h_{q} + u_{q}; \;\; g_{1} + w_{1} \leq g_{2} + w_{2} \leq \ldots \leq g_{q} + w_{q};\nonumber\\
&&\hspace{-.2in}{\rm dist} \, ((h_{j}(x^{2}, y) + u_{j}(x^{2},y), x^{2}, y), {\rm spt} \,\|{\mathbf C}\|) = 
(1 + \lambda_{j}^{2})^{-1/2}|u_{j}(x^{2}, y)|, \;(x^{2}, y) \in B_{3/4} \cap \{x^{2} < -\t\};\nonumber\\
&&\hspace{-.2in}{\rm dist} \, ((g_{j}(x^{2}, y) + w_{j}(x^{2},y), x^{2}, y), {\rm spt} \,\|{\mathbf C}\|) = 
(1 + \mu_{j}^{2})^{-1/2}|w_{j}(x^{2}, y)|, \; (x^{2}, y) \in B_{3/4} \cap \{x^{2} > \t\}.\nonumber\\
&&\hspace{-.2in}({\rm b})\;\;\;\int_{B_{5/8}^{n+1}(0)} \frac{|X^{\perp}|^{2}}{|X|^{n+2}} \, d\|V\|(X) \leq C\int_{{\mathbf R} \times B_{1}} {\rm dist}^{2} \, (X, {\rm spt} \, \|{\mathbf C}\|) \, d\|V\|(X).\nonumber\\
&&\hspace{-.2in}({\rm c})\;\;\; \sum_{j=3}^{n+1}\int_{B_{5/8}^{n+1}(0)} |e_{j}^{\perp}|^{2} \, d\|V\|(X) \leq C\int_{{\mathbf R} \times B_{1}} {\rm dist}^{2} \, (X, {\rm spt} \, \|{\mathbf C}\|) \, d\|V\|(X).\nonumber\\
&&\hspace{-.2in}({\rm d})\;\;\;\int_{B_{5/8}^{n+1}(0)} \frac{{\rm dist}^{2} \, (X, {\rm spt} \, \|{\mathbf C}\|)}{|X|^{n+2-\m}} \, d\|V\|(X) \leq \widetilde{C}\int_{{\mathbf R} \times B_{1}} {\rm dist}^{2} \, (X, {\rm spt} \, \|{\mathbf C}\|) \, d\|V\|(X).
\end{eqnarray*}
Here $e_{j}^{\perp}(X)$ denotes the orthogonal projection of $e_{j}$ onto $(T_{X} \, {\rm spt} \, \|V\|)^{\perp}$ and $C = C(n, q, \a) \in (0, \infty)$, $\widetilde{C} = \widetilde{C}(n,q, \a, \m) \in (0, \infty)$. (In particular, $C$, $\widetilde{C}$ do not depend on $\t$).
\end{theorem}

\begin{proof} We first establish conclusion (a). Let $\lambda_{1} = \lambda_{1}^{\prime}> \lambda_{2}^{\prime}> \ldots >\lambda_{p_{1}}^{\prime} = \lambda_{q}$ be the distinct elements of the set
$\{\lambda_{1}, \ldots, \lambda_{q}\}$ and $\mu_{1} = \mu_{1}^{\prime} < \mu_{2}^{\prime} < \ldots < \mu_{p_{2}}^{\prime} = \mu_{q}$ be the distinct 
elements of $\{\mu_{1}, \ldots, \mu_{q}\}$, so that $p_{1}, p_{2} \leq q$ and $p_{1} + p_{2} = p.$  By (\ref{slope-bounds-2}), 
provided $\e = \e(n, q), \g = \g(n, q) \in (0, 1)$ are sufficiently small, we have that $p_{1}, p_{2} \geq 2.$ By Remark (1) at the end of Section~\ref{step2}, Remark (3) of Section~\ref{outline} and Theorem~\ref{SS}, it follows that if $\e = \e(n, q, \a, \t), \g = \g(n, q, \a, \t) \in (0, 1)$ are sufficiently small, then  
\begin{equation}\label{L2-est-1-0}
V \res ({\mathbf R} \times (B_{3/4} \setminus \{|x^{2}| < \t\})) = \sum_{j=1}^{q} |{\rm graph} \, \widetilde{u}_{j}| + |{\rm graph} \, \widetilde{w}_{j}|
\end{equation}
where $\widetilde{u}_{j} \in C^{2} \, (B_{3/4} \setminus \{x^{2} > -\t\})$, $\widetilde{w}_{j} \in C^{2} \, (B_{3/4} \setminus \{x^{2} < \t\})$
are functions with small gradient solving the minimal surface equation and with $\widetilde{u}_{1} \leq \widetilde{u}_{2} \leq \ldots \leq \widetilde{u}_{q}$ and $\widetilde{w}_{1} \leq \widetilde{w}_{2} \leq \ldots \leq \widetilde{w}_{q}.$

If $p=4$, then $p_{1} = p_{2} = 2$  and by (\ref{slope-bounds-2}), provided 
$\e = \e(n, q), \g = \g(n, q) \in (0, 1)$ are sufficiently small, 
\begin{eqnarray*}\label{L2-est-1-1}
c {\hat E}_{V}  \leq {\rm max} \, \{|\lambda_{1}^{\prime}|, |\lambda_{2}^{\prime}|\} \leq c_{1}{\hat E}_{V}, \;\;\;\; c{\hat E}_{V}  \leq {\rm max} \, \{|\mu_{1}^{\prime}|, |\mu_{2}^{\prime}|\} \leq c_{1}{\hat E}_{V} \;\;\; {\rm and}&&\nonumber\\
&&\hspace{-3.5in} {\rm min} \, \{|\lambda_{1}^{\prime} - \lambda_{2}^{\prime}|, |\mu_{1}^{\prime} - \mu_{2}^{\prime}|\} \geq 2c{\hat E}_{V}
\end{eqnarray*}
where $c_{1}= c_{1}(n), c = c(n, q)  \in (0, \infty)$ are as in (\ref{slope-bounds-1}) and (\ref{slope-bounds-2}). Conclusion~(a) follows in this case from Hypothesis~\ref{hyp}(5) and elliptic estimates. 
Now suppose ${\mathbf C} \in {\mathcal C}_{q}(p)$ for some $p \in \{5,6, \ldots, 2q\}$ and assume by induction the following:
\begin{itemize}
\item[(A$_{1}$)] {There exist $\widetilde{\e} = \widetilde{\e}(n, q, \a, \t),$ $\widetilde{\g}= \widetilde{\g}(n, q, \a, \t)$ and 
$\widetilde{\b} = \widetilde{\b}(n, q, \a, \t) \in (0, 1)$ such that if Hypotheses~\ref{hyp}, Hypothesis~($\star$) and Hypothesis~($\star\star$) are satisfied with $M = M_{0}^{4},$ $\widetilde{\e}$, $\widetilde{\g}$, $\widetilde{\b}$ in place of $\e$, $\g$, $\b$ respectively, and with $V \in {\mathcal S}_{\a}$ and any cone $\widetilde{\mathbf C} \in \cup_{k=4}^{p-1}{\mathcal C}_{q}(k)$ in place of ${\mathbf C}$}, and if the induction hypotheses $(H1)$, $(H2)$ hold, then conclusion (a) with $\widetilde{\mathbf C}$ in place of ${\mathbf C}$ holds.   
\end{itemize}

By (\ref{separation}), 
\begin{equation}\label{L2-est-1-2}
|\lambda_{i+1}^{\prime} - \lambda_{i}^{\prime}| \geq 2c^{\prime}  Q_{V}^{\star}(p-1), \;\;\; |\mu_{j+1}^{\prime} - \mu_{j}^{\prime}| \geq 2c^{\prime} Q_{V}^{\star}(p-1)
\end{equation}
for some constant $c^{\prime} = c^{\prime}(n, q) \in (0, \infty)$ and all $i=1, 2, \ldots, p_{1}-1$ and $j=1, 2, \ldots, p_{2}-1.$
So if 
$$\left(Q_{V}^{\star}(p-1)\right)^{2} \geq \left(\frac{2}{3}\widetilde{\b}\right)^{2q}\widetilde{\g}{\hat E}^{2}_{V},$$
then it follows from (\ref{slope-bounds-1}), (\ref{L2-est-1-0}), (\ref{L2-est-1-2}) and elliptic estimates that conclusion (a) holds provided $\e = \e(n, q, \a, \t),$ $\g = \g(n, q, \a, \t) \in (0, 1)$ are sufficiently small. If on the other hand 
\begin{equation}\label{L2-est-1-3}
\left(Q_{V}^{\star}(p-1)\right)^{2} < \left(\frac{2}{3}\widetilde{\b}\right)^{2q}\widetilde{\g}{\hat E}^{2}_{V},
\end{equation}
then we argue as follows: Choose ${\mathbf C}_{1} \in \cup_{k=4}^{p-1} {\mathcal C}_{q}(k)$ such that 
\begin{equation}\label{L2-est-1-1-1}
Q_{V}^{2}({\mathbf C}_{1}) \leq \frac{3}{2} \left(Q_{V}^{\star}(p-1)\right)^{2}.
\end{equation} 
If  Hypothesis~($\star\star$) is satisfied with ${\mathbf C}_{1}$ in place of ${\mathbf C}$ and $\widetilde{\b}$ in place of $\b$, then it follows from assumption (A$_{1}$) (taken with $\widetilde{\mathbf C} = {\mathbf C}_{1}$),
(\ref{L2-est-1-1-1}), Hypothesis~($\star\star$), (\ref{L2-est-1-2}) and elliptic estimates that conclusion (a) holds  provided  $\e = \e(n, q, \a, \t),$ $\b = \b(n, q, \a, \t) \in (0, 1)$ are sufficiently small; if on the other hand 
Hypothesis~($\star\star$) is not satisfied with ${\mathbf C}_{1}$ in place of ${\mathbf C}$ and $\widetilde{\b}$ in place of $\b$, then $q \geq 3$, 
$p \geq 6$, ${\mathbf C}_{1} \in {\mathcal C}_{q}(k_{1})$ for some $k_{1} \in \{5, \ldots, p-1\}$, and 
\begin{equation}\label{L2-est-1-4}
Q_{V}^{2}({\mathbf C}_{1}) \geq \widetilde{\b} \left(Q_{V}^{\star}(k_{1}-1)\right)^{2}.
\end{equation}
Choose, in this case, a cone ${\mathbf C}_{2} \in \cup_{k=4}^{k_{1} - 1}{\mathcal C}_{q}(k)$ such that 
\begin{equation}\label{L2-est-1-6}
Q_{V}^{2}({\mathbf C}_{2}) \leq \frac{3}{2} \left(Q_{V}^{\star}(k_{1}-1)\right)^{2}
\end{equation}
and note that by (\ref{L2-est-1-1-1}), (\ref{L2-est-1-4}) and (\ref{L2-est-1-3}), we have that  
\begin{equation}\label{L2-est-1-7}
Q_{V}^{2}({\mathbf C}_{2}) \leq \widetilde{\g}{\hat E}_{V}^{2};
\end{equation}
by (\ref{L2-est-1-1-1}), (\ref{L2-est-1-2})  and (\ref{L2-est-1-4}), we have that   
\begin{equation}\label{L2-est-1-8}
|\lambda_{i+1}^{\prime} - \lambda_{i}^{\prime}| \geq \frac{4}{3}c^{\prime}\widetilde{\b}  Q_{V}^{\star}(k_{1}-1), \;\;\; |\mu_{j+1}^{\prime} - \mu_{j}^{\prime}| \geq \frac{4}{3}c^{\prime}\widetilde{\b} Q_{V}^{\star}(k_{1}-1)
\end{equation} 
for each $i  = 1, 2, \ldots, p_{1} - 1$ and $j  =1, 2, \ldots, p_{2} - 1;$ and since $Q_{V}^{\star}(p-1) \leq Q_{V}^{\star}(k_{1}-1)$, 
Hypothesis~($\star\star$) implies that
\begin{equation}\label{L2-est-1-9}
Q_{V}^{2}({\mathbf C}) \leq \b \left(Q_{V}^{\star}(k_{1}-1)\right)^{2}.
\end{equation}

So again, if Hypothesis~($\star\star$) is satisfied with ${\mathbf C}_{2}$ in place of ${\mathbf C}$ and $\widetilde{\b}$ in place of $\b$, it follows from (A$_{1}$) (taken with $\widetilde{\mathbf C} = {\mathbf C}_{2}$), (\ref{L2-est-1-8}), (\ref{L2-est-1-9}) and elliptic estimates that conclusion (a) holds  provided  $\e = \e(n, q, \a, \t),$ $\b = \b(n, q, \a, \t) \in (0, 1)$ are sufficiently small; if on the other hand  Hypothesis~($\star\star$) is not satisfied with ${\mathbf C}_{2}$ in place of ${\mathbf C}$ and $\widetilde{\b}$ in place of $\b$, then we may repeat the above argument in the obvious way. It is clear that at most $p$ repetitions of the argument are necessary to reach conclusion (a).

Now we prove conclusions (b) and (c). Let $\psi \, : \, {\mathbf R} \to [0, 1]$ be a decreasing $C^{2}$ function with $\psi(t) \equiv 1$ for $t \leq 13/16,$ $\psi(t) \equiv 0$ for $t \geq  29/32$, 
$|\psi^{\prime}(t)| \leq 32$ and $|\psi^{\prime\prime}(t)| \leq 1025$. For $\widetilde{X} = (\widetilde{x}^{1}, \widetilde{x}^{2}, \widetilde{y}) \in{\mathbf R} \times {\mathbf R} \times {\mathbf R}^{n-1}$, let 
$\widetilde{R}(\widetilde{X}) = |\widetilde{X}|$ and 
$\widetilde{r}(\widetilde{X})  = |(\widetilde{x}^{1}, \widetilde{x}^{2}, 0)|.$ We then have by the inequalities (2), (3) of the proof of Lemma~3.4 of \cite{S} that 
\begin{equation}\label{L2-est-S1}
\int_{B_{5/8}^{n+1}(0)} \frac{|\widetilde{X}^{\perp}|^{2}}{\widetilde{R}^{n+2}} \, d\|V\|(\widetilde{X}) \leq C\left(\int_{B_{1}^{n+1}(0)} \psi^{2}(\widetilde{R}) \, d\|V\|(\widetilde{X}) - \int_{B_{1}^{n+1}(0)} \psi^{2}(\widetilde{R}) \, d\|{\mathbf C}\|(\widetilde{X})\right) 
\end{equation}
and
 \begin{eqnarray}\label{L2-est-S2}
 &&\hspace{-.3in}\int_{B^{n+1}_{1}(0)} \left(1 + \sum_{j=3}^{n+1}|e_{j}^{\perp}|^{2}\right) \psi^{2}(\widetilde{R}) \, d\|V\|(\widetilde{X})\leq C\int_{B_{1}^{n+1}(0)} |(\widetilde{x}^{1}, \widetilde{x}^{2}, 0)^{\perp}|^{2}(\psi^{2}(\widetilde{R}) + 
 (\psi^{\prime}(\widetilde{R}))^{2})\,d\|V\|(\widetilde{X})  \nonumber\\
 &&\hspace{3.6in}-2\int_{B_{1}^{n+1}(0)} \widetilde{r}^{2}\widetilde{R}^{-1}\psi(\widetilde{R})\psi^{\prime}(\widetilde{R}) \, d\|V\|(\widetilde{X})
 \end{eqnarray}
 where $C = C(n) \in (0, \infty)$ and for $\|V\|$-a.e. $\widetilde{X} \in {\rm spt} \, \|V\|$,  the expression $(\widetilde{x}^{1}, \widetilde{x}^{2}, 0)^{\perp}$ denotes the orthogonal projection of $(\widetilde{x}^{1}, \widetilde{x}^{2}, 0)$ onto 
$(T_{\widetilde{X}} \, {\rm spt} \, \|V\|)^{\perp}.$ Also by the identity (6) of the same proof in \cite{S} we have that 
 \begin{equation}\label{L2-est-S3}
 \int_{B_{1}^{n+1}(0)} \psi^{2}(\widetilde{R}) \, d\|{\mathbf C}\|(\widetilde{X}) = -2\int_{B_{1}^{n+1}(0)} \widetilde{r}^{2}\widetilde{R}^{-1}\psi(\widetilde{R})\psi^{\prime}(\widetilde{R}) \, d\|{\mathbf C}\|(\widetilde{X}).
 \end{equation}
 Let $\d$ be a small positive constant to be chosen depending only on $n$, $q$ and $\a$;  let $\pi \, : \, {\mathbf R}^{n+1} \to \{0\} \times {\mathbf R}^{n}$ be the orthogonal projection and let ${\mathcal Y} = B_{15/16} \cap \{|\widetilde{x}^{2}| < 1/28\} \cap \pi \, {\rm spt} \, \|V\| \setminus \left(\{0\} \times {\mathbf R}^{n-1}\right).$ Denote by $(x,y)$ a general point in ${\mathbf R}^{n} = \{\widetilde{x}^{1} = 0\}$ where 
$x  \in {\mathbf R}$ and $y \in {\mathbf R}^{n-1}$. Write 
$${\mathcal Y}  = {\mathcal U} \, \cup \, {\mathcal W}$$ 
where ${\mathcal U}$ is the set of points $(x, y) \in {\mathcal Y}$ such that 
$$(15|x|/16)^{-n-2}\int_{{\mathbf R} \times B_{15|x|/16}(x, y)} {\rm dist}^{2} \, (\widetilde{X}, {\rm spt} \, \|{\mathbf C}\|) d\|V\|(\widetilde{X}) < \d$$ and
${\mathcal W}$ is the set of points $(x, y) \in {\mathcal Y}$ such that 
$$(15|x|/16)^{-n-2}\int_{{\mathbf R} \times B_{15|x|/16}(x, y)} {\rm dist}^{2} \, (\widetilde{X}, {\rm spt} \, \|{\mathbf C}\|) d\|V\|(\widetilde{X}) \geq \d.$$ 
Note that if $(x, y) \in {\mathcal Y}$ then $\pi^{-1}(x, y) \cap {\rm spt} \, \|V\| \neq \emptyset$, so it follows from  monotonicity of mass ratio that $\|V\|({\mathbf R} \times B_{|x|/16}(x, y)) \geq \omega_{n}(|x|/16)^{n}$. Consequently, for each point $(x, y) \in {\mathcal U}$, there is a point $Z^{(x, y)} \in {\rm spt} \, \|V\| \cap ({\mathbf R} \times B_{|x|/16}(x, y))$ with ${\rm dist} \, (Z^{(x, y)}, {\rm spt} \, \|{\mathbf C}\|) \leq 
\sqrt{2^{4n+1}\omega_{n}^{-1}\d}|x|$ and satisfying, by (\ref{slope-bounds-1}),  
$${\rm dist}_{\mathcal H} \, (\eta_{Z^{(x,y)}, 7|x|/8} \, {\rm spt} \, \|{\mathbf C}\| \cap ({\mathbf R} \times B_{1}), \{0\}\times B_{1})    
< C\sqrt{\d}$$ 
provided $\e_{0} = \e_{0}(\d)$ is sufficiently small. Here $C = C(n) \in (0, \infty).$ It also follows from Remark (1) at the end of Section~\ref{step2}, (\ref{slope-bounds-2}) and monotonicity of mass ratio that for any $\t^{\prime} \in (0, 1)$, we may ensure, by choosing $\e_{0} = \e_{0}(n, q, \a, \t^{\prime}), \g_{0} = \g_{0}(n, q, \a, \t^{\prime}) \in (0, 1)$ sufficiently small, that $\{Z \in {\rm spt} \, \|V\| \cap ({\mathbf R} \times B_{15/16}) \, : \,\Theta \, (\|V\|, Z) \geq q\} \subset \{(\widetilde{x}^{1}, \widetilde{x}^{2}, \widetilde{y}) \in {\mathbf R}^{n+1} \, : \, |\widetilde{x}^{2}| < \t^{\prime}\}$ and $\|V\|(({\mathbf R} \times B_{15/16}) \cap \{(\widetilde{x}^{1}, \widetilde{x}^{2}, \widetilde{y}) \in {\mathbf R}^{n+1} \, : \, |\widetilde{x}^{2}| < \t^{\prime}\}) < C\t^{\prime}$ where $C = C(n, q) \in (0, \infty).$  Using these facts with sufficiently small $\t^{\prime} = \t^{\prime}(n, q) \in (0, 1)$ together with Remark (3) of Section~\ref{outline} and Theorem~\ref{SS}, we find that 
$\omega_{n}^{-1}(1/16)^{-n}\|V\|(B_{1/16}^{n+1}(Z)) < q + 1/4$ for any $Z \in {\mathbf R} \times B_{14/16}$, and hence in particular that 
$$\omega_{n}^{-1}(7|x|/4)^{-n}\|V\|(B_{7|x|/4}^{n+1}(Z^{(x, y)})) < q + 1/4$$ 
for each $(x, y) \in {\mathcal U}.$ Furthermore, writing $Z^{(x, y)}_{1} = e_{1} \cdot Z^{(x, y)}$, we have for sufficiently small $\d = \d(n) \in (0, 1)$ and any $(x, y) \in {\mathcal U}$ that 
\begin{eqnarray*}
&&\hspace{-.2in}(7|x|/8)^{-n-2}\int_{{\mathbf R} \times 
B_{7|x|/8}(\pi \, Z^{(x, y)})} {\rm dist}^{2} \, (\widetilde{X}, Z^{(x, y)} + \{0\} \times {\mathbf R}^{n}) \, d\|V\|(\widetilde{X})\nonumber\\
&&\leq (7|x|/8)^{-n-2}\int_{({\mathbf R} \times 
B_{7|x|/8}(\pi \, Z^{(x, y)})) \cap \{|\widetilde{x}^{1} - Z^{(z, y)}_{1}| < \frac{3}{4}|x|\}} {\rm dist}^{2} \, (\widetilde{X}, Z^{(x, y)} + \{0\} \times {\mathbf R}^{n}) \, d\|V\|(\widetilde{X}) +\nonumber\\ 
&&\hspace{.5in}(7|x|/8)^{-n-2}\int_{({\mathbf R} \times 
B_{7|x|/8}(\pi \, Z^{(x, y)})) \cap  \{|\widetilde{x}^{1} - Z^{(x, y)}_{1}| \geq \frac{3}{4}|x|\}} {\rm dist}^{2} \, (\widetilde{X}, Z^{(x, y)} + \{0\} \times {\mathbf R}^{n}) \, d\|V\|(\widetilde{X})\nonumber\\
&&\leq c|x|^{-n-2}\int_{{\mathbf R} \times B_{|x|}(x, y)} {\rm dist}^{2} \, (\widetilde{X}, {\rm spt} \, \|{\mathbf C}\|) \, d\|V\|(\widetilde{X}) + \nonumber\\
&&\hspace{.1in}c|x|^{-n-2}\|V\|(B^{n+1}_{5|x|/4}(Z^{(x, y)})) \, {\rm dist}_{\mathcal H}^{2} \, (Z^{(x, y)} + \{0\} \times B_{7|x|/8}(0), {\rm spt} \, \|{\mathbf C}\| \cap ({\mathbf R} \times B_{7|x|/8}(\p \, Z^{(x, y)}))\nonumber
\end{eqnarray*}
where $c = c(n) \in (0, \infty)$ and we have used the pointwise inequality 
${\rm dist} \, (\widetilde{X}, Z^{(x,y)} + \{0\} \times {\mathbf R}^{n}) \leq 2 \, {\rm dist} \, (\widetilde{X}, {\rm spt} \, \|{\mathbf C}\|)$ for $\widetilde{X} \in ({\mathbf R} \times B_{7|x|/8}(\pi \, Z^{(x, y)}))  \cap \{|\widetilde{x}^{1} - Z^{(x, y)}_{1}| \geq \frac{3}{4}|x|\},$ valid if $\d = \d(n) \in (0, 1)$ and 
$\e_{0} = \e_{0}(n, q, \a) \in (0, 1)$ are sufficiently small. 
Thus provided $\e_{0} = \e_{0}(n, q, \d) \in (0, 1)$ is sufficiently small,   
\begin{equation}\label{L2-est-a}
(7|x|/8)^{-n-2}\int_{{\mathbf R} \times 
B_{7|x|/8}(\pi \, Z^{(x, y)})} {\rm dist}^{2} \, (\widetilde{X}, Z^{(x, y)} + \{0\} \times {\mathbf R}^{n}) \, d\|V\|(\widetilde{X}) < C\d
\end{equation} 
where $C = C(n, q) \in (0, \infty).$ In particular, $\|V\|(({\mathbf R} \times B_{7|x|/8}(\pi \, Z^{(x, y)})) \cap \{\widetilde{X} \, : \, {\rm dist} \, (\widetilde{X}, Z^{(x, y)} + \{0\} \times {\mathbf R}^{n}) \geq \d^{1/4}|x|\}) \leq C\sqrt{\d}|x|^{n}$ where $C = C(n, q) \in (0, \infty),$  and consequently, 
\begin{eqnarray*}
\omega_{n}^{-1}(7|x|/8)^{-n}\|V\|({\mathbf R} \times B_{7|x|/8}(\pi \, Z^{(x, y)})) &\leq& C\sqrt{\d} + \omega_{n}^{-1}(7|x|/8)^{-n}\|V\|(B^{n+1}_{(7/8 + \d^{1/4})|x|}(Z^{(x, y)}))\nonumber\\ 
&<& q + 1/2 
\end{eqnarray*}
provided $\d = \d(n, q, \a) \in (0, 1)$ is sufficiently small. Note also that (\ref{L2-est-a}) implies 
that ${\rm spt} \, \|V\| \cap ({\mathbf R} \times B_{3|x|/4}(\pi \, Z^{(x, y)})) \subset \{ \widetilde{X} \in {\mathbf R}^{n+1}\, : \, {\rm dist} \, (\widetilde{X}, Z^{(x, y)} + \{0\} \times {\mathbf R}^{n}) < |x|/2\}$ provided $\d =\d(n, q, \a) \in (0, 1)$ is sufficiently small. By applying Remark (3) of 
Section~\ref{step2} (with $\eta_{Z^{(x, y)}, 7|x|/8 \, \#} \, V$, 
$\eta_{Z^{(x, y)}, 7|x|/8} \,{\rm spt} \,  \|{\mathbf C}\|$ in place of $V$,${\mathbf P}$) we deduce that for each $(x, y) \in {\mathcal U}$, there exists a hyperplane $H_{(x,y)}$ with 
$H_{(x,y)} \cap \{\widetilde{x}^{2} >0\} \in \{G_{1}, \ldots, G_{q}\}$ (in case $x > 0$) or 
$H_{(x, y)} \cap \{\widetilde{x}^{2} < 0\} \in \{H_{1}, \ldots, H_{q}\}$ (in case $x<0$),  and an ${\mathcal H}^{n}$-measurable subset $\Sigma_{(x, y)} \subset H_{(x,y)} \cap {\rm spt} \, \|{\mathbf C}\| \cap ({\mathbf R} \times B_{|x|/4}(x, y))$ (where $\Sigma_{(x, y)} = \emptyset$ if Remark (3)(a) applies, and  
$\Sigma_{(x, y)}$ corresponds to the set $\Sigma$ as in Remark (3)(b) otherwise) such that 
\begin{eqnarray*}
&&\hspace{-.2in}\int_{({\mathbf R} \times (B_{|x|/4}(x, y)) \cap \{|\widetilde{x}^{1}| \leq |x|\} \setminus {\mathcal C}_{H_{(x,y)}}(\Sigma_{(x, y)})} |(\widetilde{x}^{1}, \widetilde{x}^{2}, 0)^{\perp}|^{2} \, d\|V\|(\widetilde{X}) \nonumber\\
&&\hspace{.1in}+\int_{({\mathbf R} \times B_{|x|/4}(x,y)) \cap {\mathcal C}_{H_{(x,y)}}(\Sigma_{(x, y)})} |\widetilde{x}^{2}|^{2} \, d\|V\|(\widetilde{X}) 
\leq C\int_{{\mathbf R} \times B_{15|x|/16}(x, y)} {\rm dist}^{2} (\widetilde{X}, {\rm spt} \, \|{\mathbf C}\|) \, d\|V\|(\widetilde{X}) 
\end{eqnarray*}
where $C = C(n, q, \a) \in (0, \infty)$ and ${\mathcal C}_{H}(A) = \{X \in {\mathbf R}^{n+1} \, : \, \pi_{H}(X) \in A\}.$   Since the pointwise inequality $|\widetilde{x}^{1}| \leq 2{\rm dist} \, (\widetilde{X}, {\rm spt} \, \|{\mathbf C}\|)$ holds whenever 
$\widetilde{X}  = (\widetilde{x}^{1}, \widetilde{x}^{2}, \widetilde{y}) \in ({\mathbf R} \times B_{|x|/4}(x, y)) \cap \{|\widetilde{x}^{1}| > |x|\}$, we also have that 
\begin{eqnarray*}
&&\hspace{-.2in}\int_{({\mathbf R} \times (B_{|x|/4}(x, y)) \cap \{|\widetilde{x}^{1}| > |x|\} \setminus {\mathcal C}_{H_{(x,y)}}(\Sigma_{(x, y)})} |(\widetilde{x}^{1}, \widetilde{x}^{2}, 0)^{\perp}|^{2} \, d\|V\|(\widetilde{X})\nonumber\\ 
&&\hspace{.3in}+\int_{({\mathbf R} \times B_{|x|/4}(x, y)) \cap \{|\widetilde{x}^{1}| > |x|\}} |\widetilde{x}^{1}|^{2} \, d\|V\|(\widetilde{X}) \leq C\int_{{\mathbf R} \times B_{15|x|/16}(x, y)} {\rm dist}^{2} (\widetilde{X}, {\rm spt} \, \|{\mathbf C}\|) \, d\|V\|(\widetilde{X}).
\end{eqnarray*}
Combining the two preceding integral estimates, we conclude that  for each $(x, y) \in {\mathcal U}$, 
\begin{eqnarray}\label{L2-est-basic}
&&\hspace{-.2in}\int_{({\mathbf R} \times B_{|x|/4}(x,y)) \cap {\mathcal C}_{H_{(x,y)}}(\Sigma_{(x, y)})} \widetilde{r}^{2} \, d\|V\|(\widetilde{X}) + 
\int_{({\mathbf R} \times B_{|x|/4}(x, y)) \setminus {\mathcal C}_{H_{(x,y)}}(\Sigma_{(x, y)})} |(\widetilde{x}^{1}, \widetilde{x}^{2}, 0)^{\perp}|^{2} \, d\|V\|(\widetilde{X})\nonumber\\ 
&&\hspace{3in}\leq C\int_{{\mathbf R} \times B_{15|x|/16}(x, y)} {\rm dist}^{2} (\widetilde{X}, {\rm spt} \, \|{\mathbf C}\|) \, d\|V\|(\widetilde{X}) 
\end{eqnarray}
where $C = C(n, q, \a) \in (0, \infty).$ We claim that (\ref{L2-est-basic}) also holds trivially (by taking $\Sigma_{(x, y)}$ to be equal to any component of ${\rm spt} \, \|{\mathbf C}\| \cap ({\mathbf R} \times B_{|x|/4}(x,y))$) whenever $(x, y) \in {\mathcal W}.$ Indeed, 
\begin{eqnarray*}
&&\hspace{-.3in}\int_{{\mathbf R} \times B_{|x|/4}(x, y)} \widetilde{r}^{2} \, d\|V\|(\widetilde{X}) = 
\int_{({\mathbf R} \times B_{|x|/4}(x, y)) \cap \{|\widetilde{x}^{1}| < |x|\}} \widetilde{r}^{2} \, 
d\|V\|(\widetilde{X}) +\nonumber\\ 
&&\hspace{3.85in}\int_{({\mathbf R} \times B_{|x|/4}(x, y)) \cap \{|\widetilde{x}^{1}| \geq |x|\}} \widetilde{r}^{2} \, d\|V\|(\widetilde{X})\nonumber\\
&&\hspace{.1in}\leq\frac{81}{16}|x|^{2}\|V\|(({\mathbf R} \times B_{|x|/4}(x, y)) \cap \{|\widetilde{x}^{1}| < |x|\})
+ 50\int_{{\mathbf R} \times B_{|x|/4}(x, y)} {\rm dist}^{2} \, (\widetilde{X}, {\rm spt} \, \|{\mathbf C}\|) \, d\|V\|(\widetilde{X})\nonumber\\
&&\hspace{.1in}\leq C|x|^{n+2} +  C\int_{{\mathbf R} \times B_{|x|/4}(x, y)} {\rm dist}^{2} \, (\widetilde{X}, {\rm spt} \, \|{\mathbf C}\|) \, d\|V\|(\widetilde{X})\nonumber\\
&&\hspace{.1in}\leq C\int_{{\mathbf R} \times B_{15|x|/16}(x, y)} {\rm dist}^{2} \, (\widetilde{X}, {\rm spt} \, \|{\mathbf C}\|) \, d\|V\|(\widetilde{X})
\end{eqnarray*}
whenever $(x, y) \in {\mathcal W},$ where $C = C(n, q, \a) \in (0, \infty).$ Thus (\ref{L2-est-basic}) holds for each $(x, y) \in {\mathcal Y}$ and some ${\mathcal H}^{n}$-measurable subset 
$\Sigma_{(x, y)} \subset H_{(x,y)} \cap {\rm spt} \, \|{\mathbf C}\| \cap ({\mathbf R} \times B_{|x|/4}(x, y))$. 

Choose now  a countable collection ${\mathcal I}$ of points $(x, y) \in {\mathcal Y}$ such that ${\mathcal Y} \subset \cup_{(x, y) \in {\mathcal I}} B_{|x|/8}(x, y)$ and the collection $\{B_{15|x|/16}(x, y)\}_{(x, y) \in {\mathcal I}}$ can be decomposed into at most $N = N(n)$ pairwise disjoint 
sub-collections. (This can be achieved e.g.\ as follows: Use the ``5-times covering lemma'' (\cite{S1}, Theorem 3.3) to extract a countable collection ${\mathcal I}$ of points 
$(x,y) \in {\mathcal Y}$ such that the collection of closed balls $\{\overline{B}_{|x|/41}(x,y)\}_{(x,y) \in {\mathcal I}}$ is pairwise disjoint and ${\mathcal Y} \subset \cup_{(x,y) \in {\mathcal I}} B_{|x|/8}(x,y).$ 
Then the collection ${\mathcal B} = \{B_{15|x|/16}(x,y)\}_{(x,y) \in {\mathcal I}}$ automatically will have the property that for each $(x_{0}, y_{0}) \in {\mathcal I}$, 
$${\rm card} \, \{(x,y) \in {\mathcal I} \, : \, B_{15|x|/16}(x,y) \cap B_{15|x_{0}|/16}(x_{0},y_{0}) \neq \emptyset\} \leq N \hspace{1.5in} (\dag)$$ 
for some fixed constant $N = N(n),$ from which it follows as required that $\cup {\mathcal B} = \cup_{j=1}^{N} \cup {\mathcal B}_{j}$ where ${\mathcal B}_{1}, \ldots, {\mathcal B}_{N} \subset {\mathcal B}$ and each ${\mathcal B}_{j}$ consists of pairwise disjoint balls. To see ($\dag$), note that 
$B_{15|x|/16}(x,y) \cap B_{15|x_{0}|/16}(x_{0},y_{0}) \neq \emptyset$ $\implies$ $|(x,y) - (x_{0}, y_{0})| \leq 15|x_{0}|/16 + 15|x|/16$ whence 
$|x| \leq 31|x_{0}| \leq 31 \times 31|x|$ and $|(x,y) - (x_{0}, y_{0})| \leq c|x_{0}| -|x|/41$ where 
$c = 15/16 + (31 \times 15)/16 + 31/41,$ which say that $B_{|x_{0}|/(31 \times 41)}(x,y) \subset B_{|x|/41}(x,y) \subset 
B_{c|x_{0}|}(x_{0},y_{0})$; since $B_{|x|/41}(x,y),$ $(x,y) \in {\mathcal I}$ are pairwise disjoint, the assertion ($\dag$) follows.) Let
$${\mathcal G} =  \bigcup_{(x, y) \in {\mathcal I}} \, \left(({\mathbf R} \times B_{|x|/8}(x,y)) \setminus  {\mathcal C}_{H_{(x,y)}}(\Sigma_{(x, y)})\right).$$  
We deduce from (\ref{L2-est-basic}) that
\begin{equation}\label{L2-est-basic-1}
\int_{({\mathbf R} \times {\mathcal Y}) \setminus {\mathcal G}} |\widetilde{r}|^{2}  \, d\|V\|(\widetilde{X}) +
\int_{({\mathbf R} \times {\mathcal Y}) \cap {\mathcal G}} |(\widetilde{x}^{1}, \widetilde{x}^{2}, 0)^{\perp}|^{2} \, d\|V\|(\widetilde{X})
\leq C \int_{{\mathbf R} \times B_{1}} {\rm dist}^{2} \, (\widetilde{X}, {\rm spt} \, \|{\mathbf C}\|) \, d\|V\|(\widetilde{X}) 
\end{equation}
where $C = C(n, q, \a) \in (0, \infty).$

Now let ${\mathcal J}$ be a collection of $J = J(n)$ points $w \in B_{15/16} \setminus \{|\widetilde{x}^{2}|< 1/28\}$ such that $B_{15/16} \setminus \{|\widetilde{x}^{2}|< 1/28\} \subset \cup_{w \in {\mathcal J}} B_{1/64}(w).$   For $z \in {\mathbf R}^{n}$ and $\r>0$, let $T_{\r}(z) = \{(\widetilde{x}\sin \th, \widetilde{x}\cos \th, \widetilde{y}) \, : \, (\widetilde{x}, \widetilde{y}) \in B_{\r}(z), \;\; \th \in [0, 2\pi)\}.$  Note that if $\e_{0}= \e_{0}(n, q, \a, 1/32) \in (0, 1),$ $\g_{0} = \g_{0}(n,q,\a, 1/32) \in (0, 1)$, $\b_{0}= \b_{0}(n,q,\a, 1/32) \in (0, 1)$ are sufficiently small, then for each $(x, y) \in {\mathcal I}$, 
\begin{eqnarray}\label{L2-est-S4}
\left(({\mathbf R} \times B_{|x|/8}(x,y)) \setminus {\mathcal C}_{H_{(x,y)}}(\Sigma_{(x, y)})\right) \cap {\rm spt} \, \|V\| &\subseteq& \left(T_{9|x|/64}^{\pm}(x, y) \setminus {\mathcal C}_{H_{(x,y)}}(\Sigma_{(x,y)})\right) \cap {\rm spt} \, \|V\|\nonumber\\ 
&\subseteq& \left(T_{3|x|/16}^{\pm}(x, y) \setminus {\mathcal C}_{H_{(x,y)}}(\Sigma_{(x,y)})\right) \cap {\rm spt} \, \|V\|\nonumber\\ 
&\subseteq& \left(({\mathbf R} \times B_{|x|/4}(x,y))  \setminus {\mathcal C}_{H_{(x,y)}}(\Sigma_{(x,y)})\right) \cap {\rm spt} \, \|V\|
\end{eqnarray}
and for each $w \in {\mathcal J},$
\begin{equation}\label{L2-est-S4-1}
{\mathbf R} \times B_{1/64}(w) \cap {\rm spt} \, \|V\| \subseteq T_{9/512}^{\pm}(w) \cap {\rm spt} \, \|V\| \subseteq T_{3/128}^{\pm}(w) \cap {\rm spt} \, \|V\|
\end{equation}
 where $T_{\r}^{+}(z) = T_{\r}(z) \cap \{|\widetilde{x}^{1}| < |\widetilde{x}^{2}|\} \cap \{\widetilde{x}^{2} > 0\}$; $T_{\r}^{-}(z) = T_{\r}(z) \cap \{|\widetilde{x}^{1}| < |\widetilde{x}^{2}|\} \cap \{\widetilde{x}^{2} < 0\};$ in  (\ref{L2-est-S4}) we choose the $+$ sign if $x > 0$ and the $-$ sign if $x < 0;$ in (\ref{L2-est-S4-1}) we choose the $+$ sign if $e_{2} \cdot w > 0$ and the $-$ sign if $e_{2} \cdot w < 0.$ 
 
 Now, applying [\cite{F}, 3.1.13] with 
 $\Phi = \{B_{3|x|/16}(x,y)\}_{(x,y) \in {\mathcal I}} \cup \{B_{3/128}(w)\}_{w \in {\mathcal J}}$, 
 and letting $h(p) = \frac{1}{20}\sup \{\inf\{1, {\rm dist} \, (p, {\mathbf R}^{n} \setminus B)\} \, : \, B \in \Phi\}$ for $p \in \cup \Phi$, we obtain a smooth partition of unity $\{\varphi_{s}\}_{s \in {\mathcal S}}$ having the following properties:
\begin{itemize}
\item[(i)] ${\mathcal S}$ is a countable subset of $\cup \Phi$ and $\varphi_{s} \, : \, \cup \Phi \to [0, 1]$
$\forall$ $s \in {\mathcal S}$.
\item[(ii)] $\{B_{h(s)}(s)\}_{s \in {\mathcal S}}$ is pairwise disjoint and for each $s \in {\mathcal S}$, 
$B_{h(s)}(s)\subset {\rm spt} \, \varphi_{s} \subset B_{10h(s)}(s) \subset B$ for some $B \in \Phi$. 
\item[(iii)] $\sum_{s \in {\mathcal S}} \varphi_{s}(p) = 1$ for each $p \in \cup \Phi$.
\item[(iv)] $|D\varphi_{s}(p)| \leq Ch(p)^{-1}$
for each $s \in {\mathcal S}$ and each $p \in \cup \Phi,$
where $C = C(n) \in (0, \infty).$ 
\end{itemize}
 Note in particular that it follows from (iv) and the definition of $h(\cdot)$ that for each $s \in {\mathcal S}$, 
\begin{equation}\label{L2-est-S5}
|D\varphi_{s}(\widetilde{x}, \widetilde{y})| \leq C \, |\widetilde{x}|^{-1}
\end{equation}
whenever $(\widetilde{x}, \widetilde{y}) \in \cup_{(x,y) \in {\mathcal I}} \, B_{5|x|/32}(x,y) \cup \cup_{w \in {\mathcal J}} B_{5/256}(w),$ where $C = C(n) \in (0, \infty).$ For each $s \in {\mathcal S}$, extend $\varphi_{s}$ to ${\mathbf R}^{n}$ by setting 
$\varphi_{s}(x) = 0$ for $x \in {\mathbf R}^{n} \setminus \cup\Phi$, and let $\widetilde{\varphi}_{s}$ be the (smooth) extension of $\varphi_{s}$  to 
$\{\widetilde{X} = (\widetilde{x}^{1}, \widetilde{x}^{2}, \widetilde{y}) \in {\mathbf R}^{n+1} \, : \, 
|\widetilde{x}^{1}| < |\widetilde{x}^{2}|\}$ defined by  
$\widetilde{\varphi}_{s}(\widetilde{x}^{1}, \widetilde{x}^{2}, \widetilde{y}) = \varphi_{s}(\pm\sqrt{|\widetilde{x}^{1}|^{2} + 
|\widetilde{x}^{2}|^{2}}, \widetilde{y})$ where the $+$ sign is chosen if $\widetilde{x}^{2} >0$ and 
the $-$ sign if $\widetilde{x}^{2} < 0.$ 

Let $\widetilde{\mathcal G} = {\mathcal G} \cup \left({\mathbf R} \times (B_{15/16} \setminus \{|\widetilde{x}^{2}| < 1/28\})\right).$ We claim that there exists a fixed constant $M = M(n)$ such that 
for each $(x,y) \in {\mathcal I}$, 
\begin{equation}\label{L2-est-S5-1}
{\rm card} \, \{s \in {\mathcal S} \, : \, {\rm spt} \, \widetilde{\varphi}_{s} \subset T_{3|x|/16}(x,y) \;\; {\rm and} \;\; {\rm spt} \, \widetilde{\varphi}_{s} \cap \widetilde{\mathcal G} \cap {\rm spt} \, \|V\| \neq \emptyset\} \leq M
\end{equation}
and for each $w \in {\mathcal J}$, 
\begin{equation}\label{L2-est-S5-2}
{\rm card} \, \{s \in {\mathcal S} \, : \, {\rm spt} \, \widetilde{\varphi}_{s} \subset T_{3/128}(w) \;\; {\rm and} \;\; {\rm spt} \, \widetilde{\varphi}_{s} \cap \widetilde{\mathcal G} \cap {\rm spt} \, \|V\| \neq \emptyset\} \leq M.
\end{equation}
To see (\ref{L2-est-S5-1}), fix $(x,y) \in {\mathcal I}$ and let 
${\mathcal S}_{(x,y)} =  \{s \in {\mathcal S} \, : \, {\rm spt} \, \widetilde{\varphi}_{s} \subset T_{3|x|/16}(x,y) \;\; {\rm and} \;\; {\rm spt} \, \widetilde{\varphi}_{s} \cap \widetilde{\mathcal G} \cap {\rm spt} \, \|V\| \neq \emptyset\}.$
Note that ${\rm spt} \, \widetilde{\varphi}_{s} \subset T_{3|x|/16}(x,y) \iff {\rm spt} \, \varphi_{s} \subset B_{3|x|/16}(x,y),$ and since 
$\widetilde{\mathcal G} \cap {\rm spt} \, \|V\| \subset \cup_{(x,y) \in {\mathcal I}} \left(({\mathbf R} \times B_{|x|/8}(x,y)) \setminus {\mathcal C}_{H_{(x,y)}}(\Sigma_{(x,y)})\right) \cup \cup_{w \in {\mathcal J}} ({\mathbf R} \times B_{1/64}(w)) \cap {\rm spt} \, \|V\|$, it follows from (\ref{L2-est-S4}), (\ref{L2-est-S4-1})  and (ii) above that 
if $s \in {\mathcal S}_{(x,y)}$ then either 
\begin{itemize}
\item[($\star$)] $B_{3|x|/16}(x,y) \cap B_{9|x^{\prime}|/64}(x^{\prime}, y^{\prime}) \neq \emptyset$ and $B_{10h(s)}(s) \cap B_{9|x^{\prime}|/64}(x^{\prime},y^{\prime}) \neq \emptyset$ for some 
$(x^{\prime}, y^{\prime}) \in {\mathcal I}$ or 
\item[($\star\star$)] $B_{10h(s)}(s) \cap B_{9/512}(w^{\prime}) \neq \emptyset$ for some $w^{\prime} \in {\mathcal J}.$ 
\end{itemize} 
If ($\star$) holds then $|x - x^{\prime}| < 3|x|/16 + 9|x^{\prime}|/64$ whence $|x^{\prime}|> 52|x|/73,$ and  $|s - (x^{\prime}, y^{\prime})| < 10h(s) +  9|x^{\prime}|/64;$ so if $h(s) < |x^{\prime}|/640$ then  $s \in B_{5|x^{\prime}|/32}(x^{\prime}, y^{\prime})$ and hence, since $B_{3|x^{\prime}|/16}(x^{\prime}, y^{\prime}) \in \Phi$, it follows from the definition of $h(s)$ that $h(s) \geq |x^{\prime}|/640$ contrary to our assumption. Hence in case ($\star$) holds, we must have that $h(s) \geq 52|x|/(640 \times 73).$ In case ($\star\star$) holds, similar reasoning shows that $h(s) \geq 1/5120.$  Thus  
for any fixed $(x,y) \in {\mathcal I}$, we have established that 
$s \in {\mathcal S}_{(x,y)}$ $\implies$ $h(s) \geq \min\{52|x|/(640\times73), 1/5120\}$ and (by (ii) above) $B_{h(s)}(s) \subset B_{3|x|/16}(x,y).$ Since $B_{h(s)}(s),$ $s \in {\mathcal S}$ are pairwise 
disjoint, this establishes (\ref{L2-est-S5-1}) for some fixed $M = M(n).$ Identical reasoning (using the fact that $|e_{2} \cdot w| > 1/28$ for each $w \in {\mathcal J}$) establishes 
(\ref{L2-est-S5-2}).

 Noting, by (\ref{L2-est-S4}) and the definition of 
$\Sigma_{(x,y)},$ that the set $\left(T_{3|x|/16}^{\pm}(x, y) \setminus {\mathcal C}_{H_{(x,y)}}(\Sigma_{(x,y)})\right) \cap {\rm spt} \, \|V\|,$ if non-empty, can be written as the union of normal graphs of Lipschitz functions defined over subsets of  a sub-collection of the half-hyperplanes $G_{1}, \ldots, G_{q}$ (if $x > 0$) or of the half-hyperplanes $H_{1}, \ldots, H_{q}$ (if $x < 0$), we see from the area formula and Remark (3) of Section~\ref{step2} that for any given $(x,y) \in {\mathcal I}$ and any $s \in {\mathcal S}$ with ${\rm spt} \, \widetilde{\varphi}_{s} \subset T_{3|x|/16}(x,y)$,
\begin{eqnarray}\label{L2-est-S6}
&&\hspace{-.2in}\int_{{\mathcal G} \cup ({\mathbf R} \times (B_{15/16}\setminus \{|\widetilde{x}^{2}| < 1/28\}))}\widetilde{\varphi}_{s}(\widetilde{X})\widetilde{r}^{2}\widetilde{R}^{-1}\psi(\widetilde{R})\psi^{\prime}(\widetilde{R}) \, d\|V\|(\widetilde{X})\nonumber\\
&& =\sum_{k=1}^{\ell(x,y)}\sum_{i=1}^{q_{k}(x,y)}\int_{\Omega_{k}(x,y)} \varphi_{s}\left(\pm\sqrt{\widetilde{r}^{2} + |u_{k}^{i}(\widetilde{X})|^{2}}, \widetilde{y}\right)\widetilde{r}_{u_{k}^{i}}^{2}\widetilde{R}_{u_{k}^{i}}^{-1}\psi(\widetilde{R}_{u_{k}^{i}})\psi^{\prime}(\widetilde{R}_{u_{k}^{i}})
\sqrt{1 + |\nabla \, u_{k}^{i}|^{2}} \, d{\mathcal H}^{n}(\widetilde{X})\nonumber\\
&& =\sum_{k=1}^{\ell(x,y)}q_{k}(x,y)\int_{\Omega_{k}(x,y)}\varphi_{s}(\pm\widetilde{r},\widetilde{y})\widetilde{r}^{2}\widetilde{R}^{-1}\psi(\widetilde{R})\psi^{\prime}(\widetilde{R}) \, d{\mathcal H}^{n}(\widetilde{X}) +\nonumber\\
&&\hspace{.2in}\sum_{k=1}^{\ell(x,y)}\sum_{i=1}^{q_{k}(x,y)}\int_{\Omega_{k}(x,y)} \left(\varphi_{s}\left(\pm\sqrt{\widetilde{r}^{2} + |u_{k}^{i}(\widetilde{X})|^{2}}, \widetilde{y}\right) - \varphi_{s}(\pm\widetilde{r}, 
\widetilde{y})\right) \widetilde{r}^{2}\widetilde{R}^{-1}\psi(\widetilde{R})\psi^{\prime}(\widetilde{R}) \, d{\mathcal H}^{n}(\widetilde{X})  + E
\end{eqnarray}
where we choose the $+$ sign if $x > 0$ and the $-$ sign if $x < 0$; $\ell(x,y)$ is a positive integer $\leq q$; $q_{k}(x,y)$ are positive integers with 
\begin{equation}\label{L2-est-S7}
\sum_{k=1}^{\ell(x,y)} q_{k}(x,y) \leq q;
\end{equation}
$\Omega_{k}(x,y)$ is, by (\ref{L2-est-S4}) and (\ref{L2-est-S4-1}),   a measurable subset of 
$$\left(\cup_{(x^{\prime}, y^{\prime}) \in {\mathcal I}}T_{19|x^{\prime}|/128}(x^{\prime}, y^{\prime}) \cup \cup_{w^{\prime} \in {\mathcal J}}T_{19/1024}(w^{\prime})\right) \cap ({\mathbf R} \times B_{|x|/4}(x,y)) \cap G_{j_{k}(x,y)}$$ 
(if $x >0$) or of  $$\left(\cup_{(x^{\prime}, y^{\prime}) \in {\mathcal I}}T_{19|x^{\prime}|/128}(x^{\prime}, y^{\prime}) \cup \cup_{w^{\prime} \in {\mathcal J}}T_{19/1024}(w^{\prime})\right) \cap({\mathbf R} \times B_{|x|/4}(x,y)) \cap H_{j_{k}(x,y)}$$ (if $x < 0$) for some integer $j_{k}(x,y) \in \{1, 2, \ldots, q\}$; $u_{k}^{i}$ 
are the Lipschitz functions as in Remark (3) of Section~\ref{step2} (applied with $\eta_{Z^{(x, y)}, 7|x|/8 \, \#} \, V$, 
$\eta_{Z^{(x, y)}, 7|x|/8} \,{\rm spt} \,  \|{\mathbf C}\|$ in place of $V$,${\mathbf P}$); 
$\widetilde{r}_{u_{k}^{i}} = \sqrt{\widetilde{r}^{2} + |u_{k}^{i}|^{2}};$ $\widetilde{R}_{u_{k}^{i}} = 
\sqrt{\widetilde{R}^{2} + |u_{k}^{i}|^{2}}$ and, by the estimates of Remark (3) of Section~\ref{step2},  
$$|E| \leq C\int_{{\mathbf R} \times B_{15|x|/16}(x,y)} {\rm dist}^{2} \, (\widetilde{X}, {\rm spt} \, \|{\mathbf C}\|) \, d\|V\|(\widetilde{X})$$
for some constant $C = (n, q) \in (0, \infty)$. Still assuming ${\rm spt} \, \widetilde{\varphi}_{s} \subset T_{3|x|/16}(x,y)$, we also see in view of  (\ref{L2-est-S7}) that 
\begin{eqnarray}\label{L2-est-S7-1}
&&\hspace{-.5in}\sum_{k=1}^{\ell(x,y)}q_{k}(x,y)\int_{\Omega_{k}(x,y)}\varphi_{s}(\pm\widetilde{r},\widetilde{y})\widetilde{r}^{2}\widetilde{R}^{-1}\psi(\widetilde{R})\psi^{\prime}(\widetilde{R}) \, d{\mathcal H}^{n}(\widetilde{X})\nonumber\\
&&\hspace{.2in}\geq\sum_{k=1}^{\ell(x,y)}q_{k}(x,y)\int_{P_{k} \cap T_{3|x|/16}(x,y)}\varphi_{s}(\pm\widetilde{r},\widetilde{y})\widetilde{r}^{2}\widetilde{R}^{-1}\psi(\widetilde{R})\psi^{\prime}(\widetilde{R}) \, d{\mathcal H}^{n}(\widetilde{X})\nonumber\\
&&\hspace{.2in}=\left(\sum_{k=1}^{\ell(x,y)}q_{k}(x,y)\right)\int_{B_{3|x|/16}(x,y)}\varphi_{s}(\pm\widetilde{r},\widetilde{y})\widetilde{r}^{2}\widetilde{R}^{-1}\psi(\widetilde{R})\psi^{\prime}(\widetilde{R}) \, d{\mathcal H}^{n}(\widetilde{X})\nonumber\\
&&\hspace{.2in}\geq q\int_{B_{3|x|/16}(x,y)}\varphi_{s}(\pm\widetilde{r},\widetilde{y})\widetilde{r}^{2}\widetilde{R}^{-1}\psi(\widetilde{R})\psi^{\prime}(\widetilde{R}) \, d{\mathcal H}^{n}(\widetilde{X})\nonumber\\
&&\hspace{.2in}=\int_{{\mathbf R} \times B_{15/16}}\varphi_{s}(\pm\widetilde{r},\widetilde{y})\widetilde{r}^{2}\widetilde{R}^{-1}\psi(\widetilde{R})\psi^{\prime}(\widetilde{R}) \, d\|{\mathbf C}\|(\widetilde{X})
\end{eqnarray}
where $P_{k} = G_{j_{k}}$ if $x > 0$ and $P_{k} = H_{j_{k}}$ if $x < 0$. Since we may bound, using the ${\rm sup}$ estimate of Remark (3)(b) of Section~\ref{step2} and (\ref{L2-est-S5}) (keeping in mind that 
$\Omega_{k}(x,y)$ $\subset$ $\left(\cup_{(x^{\prime}, y^{\prime}) \in {\mathcal I}}T_{19|x^{\prime}|/128}(x^{\prime}, y^{\prime})\right.$ $\cup$ $\left.\cup_{w^{\prime} \in {\mathcal J}}T_{19/1024}(w^{\prime})\right) \cap {\rm spt} \, \|{\mathbf C}\|$), the absolute value of the middle term of the last line of (\ref{L2-est-S6}) by a constant times $\int_{{\mathbf R} \times B_{15|x|/16}(x,y)}{\rm dist}^{2} \, (\widetilde{X}, {\rm spt}\,\|{\mathbf C}\|) \, d\|V\|(\widetilde{X}),$ 
we conclude from (\ref{L2-est-S6}) and (\ref{L2-est-S7-1}) that for each $(x,y) \in {\mathcal I}$ and each $s \in {\mathcal S}$ with 
${\rm spt} \, \widetilde{\varphi}_{s} \subset T_{3|x|/16}(x,y)$,
\begin{eqnarray}\label{L2-est-S8}
&&\hspace{-.2in}\int_{\mathcal G \cup ({\mathbf R} \times (B_{15/16} \setminus \{|\widetilde{x}^{2}| < 1/28\}))}\widetilde{\varphi}_{s}(\widetilde{X})\widetilde{r}^{2}\widetilde{R}^{-1}\psi(\widetilde{R})\psi^{\prime}(\widetilde{R}) \, d\|V\|(\widetilde{X}) \geq\nonumber\\
&&\hspace{1.5in}\int_{{\mathbf R} \times B_{15/16}}\varphi_{s}(\pm\widetilde{r},\widetilde{y})\widetilde{r}^{2}\widetilde{R}^{-1}\psi(\widetilde{R})\psi^{\prime}(\widetilde{R}) \, d\|{\mathbf C}\|(\widetilde{X})\nonumber\\
&&\hspace{3in} - C\int_{{\mathbf R} \times B_{15|x|/16}(x,y)}{\rm dist}^{2} \, (\widetilde{X}, {\rm spt}\,\|{\mathbf C}\|) \, d\|V\|(\widetilde{X})
\end{eqnarray}
where $C = C(n, q) \in (0, \infty);$ the $+$ sign is chosen if $x > 0$ and the $-$ sign if $x < 0.$ By a similar argument using part (a) and elliptic estimates, we see also that
for each $w \in {\mathcal J}$ and each $s \in {\mathcal S}$ with 
${\rm spt} \, \widetilde{\varphi}_{s} \subset T_{3/128}(w)$,
\begin{eqnarray}\label{L2-est-S8-1}
&&\hspace{-.2in}\int_{\mathcal G \cup ({\mathbf R} \times (B_{15/16} \setminus \{|\widetilde{x}^{2}| < 1/28\}))}\widetilde{\varphi}_{s}(\widetilde{X})\widetilde{r}^{2}\widetilde{R}^{-1}\psi(\widetilde{R})\psi^{\prime}(\widetilde{R}) \, d\|V\|(\widetilde{X}) \geq \nonumber\\
&&\hspace{1.5in}\int_{{\mathbf R} \times B_{15/16}}\varphi_{s}(\pm\widetilde{r},\widetilde{y})\widetilde{r}^{2}\widetilde{R}^{-1}\psi(\widetilde{R})\psi^{\prime}(\widetilde{R}) \, d\|{\mathbf C}\|(\widetilde{X})\nonumber\\
&&\hspace{3in} - C\int_{{\mathbf R} \times B_{1/32}(w)}{\rm dist}^{2} \, (X, {\rm spt}\,\|{\mathbf C}\|) \, d\|V\|(\widetilde{X})
\end{eqnarray}
where $C = C(n, q) \in (0, \infty)$; the $+$ sign is chosen if $e_{2}\cdot w> 0$ and the $-$ sign if $e_{2}\cdot w < 0.$ 

Now choose enumerations ${\mathcal J} = \{w_{j}\}_{j=1}^{J}$ and 
${\mathcal I} = \{(x_{J+j}, y_{J+j})\}_{j=1}^{\infty},$ let 
$${\mathcal S}_{j} = \left\{s \in {\mathcal S} \, : \, 
{\rm spt} \, \widetilde{\varphi}_{s} \subset T_{3/128}(w_{j}), \;\; {\rm spt} \, \widetilde{\varphi}_{s} \cap \left({\mathcal G} \cup \left({\mathbf R} \times (B_{15/16} \setminus \{|\widetilde{x}^{2}| < 1/28\})\right)\right) \cap {\rm spt} \, \|V\| \neq \emptyset\right\}$$ 
for $1 \leq j \leq J$ and 
$${\mathcal S}_{j} = \{s \in {\mathcal S} \, : \, 
{\rm spt} \, \widetilde{\varphi}_{s} \subset T_{3|x_{j}|/16}(x_{j},y_{j}), \;{\rm spt} \, \widetilde{\varphi}_{s} \cap \left({\mathcal G} \cup \left({\mathbf R} \times (B_{15/16} \setminus \{|\widetilde{x}^{2}| < 1/28\})\right)\right) \cap {\rm spt} \, \|V\| \neq \emptyset\}$$ 
for $j \geq J+1$, write 
$$\left\{s \in {\mathcal S} \, : \, {\rm spt} \, \widetilde{\varphi}_{s} \cap \left({\mathcal G} \cup \left({\mathbf R} \times (B_{15/16} \setminus \{|\widetilde{x}^{2}| < 1/28\})\right)\right) \cap {\rm spt} \, \|V\| \neq \emptyset\right\} = \cup_{j=1}^{\infty} {\mathcal S}^{\prime}_{j}$$
where ${\mathcal S}^{\prime}_{1} = {\mathcal S}_{1}$ and ${\mathcal S}^{\prime}_{j} = {\mathcal S}_{j} \setminus \cup_{i=1}^{j-1} {\mathcal S}_{i}^{\prime}$ for $j \geq 2,$ and note that ${\mathcal S}^{\prime}_{j}$ are pairwise disjoint and, by 
(\ref{L2-est-S5-1}), (\ref{L2-est-S5-2}), that ${\rm card} \, ({\mathcal S}^{\prime}_{j}) \leq M = M(n).$
Summing in (\ref{L2-est-S8}), (\ref{L2-est-S8-1}) first over $s \in {\mathcal S}_{j}^{\prime}$ for fixed $j$, and then over $j$ (where $j\in \{1, 2, \ldots, J\}$ in (\ref{L2-est-S8-1}) and $j \geq J+1$ in (\ref{L2-est-S8})) keeping in mind that the collection of balls $\{B_{15|x|/16}(x,y)\}_{(x,y) \in {\mathcal I}} = 
\{B_{15|x_{j}|/16}(x_{j}, y_{j})\}_{j=J+1}^{\infty}$ can be subdivided into at most $N = N(n)$ sub-collections of pairwise disjoint balls,  and adding the two resulting inequalities (and using the fact that $\sum_{s \in {\mathcal S}} \widetilde{\varphi}_{s}(\widetilde{X}) = 1$ for each point $\widetilde{X} \in {\mathcal G} \cup \left({\mathbf R} \times (B_{15/16} \setminus \{|\widetilde{x}^{2}| < 1/28\})\right) \cap {\rm spt} \, \|V\|$ and 
$\sum_{s \in {\mathcal S}} \varphi_{s}(\pm\widetilde{r}, \widetilde{y}) \leq 1$ for each point 
$\widetilde{X} = (\widetilde{x}^{1}, \widetilde{x}^{2}, \widetilde{y}) \in {\rm spt } \, \|{\mathbf C}\|$) we conclude that 
\begin{eqnarray}\label{L2-est-S9}
&&\hspace{-.2in}\int_{{\mathcal G} \cup ({\mathbf R} \times (B_{15/16} \setminus \{|\widetilde{x}^{2}| < 1/28\}))}\widetilde{r}^{2}\widetilde{R}^{-1}\psi(\widetilde{R})\psi^{\prime}(\widetilde{R}) \, d\|V\|(\widetilde{X}) - 
\int_{{\mathbf R} \times B_{15/16}}\widetilde{r}^{2}\widetilde{R}^{-1}\psi(\widetilde{R})\psi^{\prime}(\widetilde{R}) \, d\|{\mathbf C}\|(\widetilde{X})\nonumber\\ 
&&\hspace{3in}\geq -C\int_{{\mathbf R} \times B_{1}}{\rm dist}^{2} \, (\widetilde{X}, {\rm spt} \, \|{\mathbf C}\|) \, d\|V\|(\widetilde{X})
\end{eqnarray}
where $C = C(n, q) \in (0, \infty)$. In view of (\ref{L2-est-S1}), (\ref{L2-est-S2}), (\ref{L2-est-S3}), conclusions (b) and (c) now follow from the estimates (\ref{L2-est-basic-1}), (\ref{L2-est-S9}) and  conclusion (a).  Conclusion (d) follows from conclusion (b) by exactly the same argument as for the corresponding estimate in Lemma 3.4 of \cite{S}.
\end{proof}

For the proof of Corollary~\ref{L2-est-2}  below and subsequently,  we shall need the following elementary fact: If ${\mathbf C} \in {\mathcal C}_{q}$ is as in Hypothesis~\ref{hyp}(2) and if
$Z = (\z^{1}, \z^{2}, \eta) \in {\mathbf R} \times {\mathbf R} \times {\mathbf R}^{n-1} \equiv {\mathbf R}^{n+1},$ then for any $X \in {\mathbf R}^{n+1}$,
\begin{equation}\label{cone-dist}
\left|{\rm dist} \, (X, {\rm spt} \, \|{\mathbf C}\|) - {\rm dist} \, (X, {\rm spt} \, \|T_{Z \, \#} \, {\mathbf C}\|)\right| \leq |\z^{1}| + \nu|\z^{2}|
\end{equation}
where $T_{Z} \, : \, {\mathbf R}^{n+1} \to {\mathbf R}^{n+1}$ is the translation $X \mapsto X + Z$ and 
$\nu  = {\rm max} \, \{|\lambda_{1}|, \ldots, |\lambda_{q}|, |\m_{1}|, \ldots, |\m_{q}|\}.$

Indeed, by the triangle inequality 
$$\left|{\rm dist} \, (X, {\rm spt} \, \|{\mathbf C}\|) - {\rm dist} \, (X, {\rm spt} \, \|T_{Z \, \#} \, {\mathbf C}\|)\right| \leq {\rm dist}_{\mathcal H} \, ({\rm spt} \, \|{\mathbf C}\|, {\rm spt} \, \|T_{Z \, \#} \, {\mathbf C}\|)$$ 
and by translation invariance of ${\mathbf C}$ along $\{0\} \times {\mathbf R}^{n-1},$
\begin{eqnarray*}
&&\hspace{-.5in}{\rm dist}_{\mathcal H} \, ({\rm spt} \, \|{\mathbf C}\|, {\rm spt} \, \|T_{Z \, \#} \, {\mathbf C}\|)
 = {\rm dist}_{\mathcal H} \, ({\rm spt} \, \|{\mathbf C}\|, {\rm spt} \, \|T_{(\z^{1}, \z^{2}, 0) \, \#} \, {\mathbf C}\|)\nonumber\\
&& \leq {\rm dist}_{\mathcal H} \, ({\rm spt} \, \|{\mathbf C}\|, {\rm spt} \, \|T_{(\z^{1}, 0, 0) \, \#} \, {\mathbf C}\|) + {\rm dist}_{\mathcal H} \, ({\rm spt} \, \|T_{(\z^{1}, 0, 0) \, \#} \, {\mathbf C}\|, {\rm spt} \, \|T_{(\z^{1},\z^{2}, 0)\, \#} \, {\mathbf C}\|)\nonumber\\
&&={\rm dist}_{\mathcal H} \, ({\rm spt} \, \|{\mathbf C}\|, {\rm spt} \, \|T_{(\z^{1}, 0, 0) \, \#} \, {\mathbf C}\|) + {\rm dist}_{\mathcal H} \, ({\rm spt} \, \|{\mathbf C}\|, {\rm spt} \, \|T_{(0, \z^{2}, 0)\, \#} \, {\mathbf C}\|)\leq |\z^{1}| + \nu|\z^{2}|.\\
\end{eqnarray*}

\begin{corollary}\label{L2-est-2} 
Let $q$ be an integer $\geq 2$ and $\a \in (0, 1)$. 
For each $\r \in (0, 1/4]$, there exist numbers $\e = \e(n , q, \a, \r) \in (0, 1),$ 
$\g = \g(n, q, \a, \r) \in (0, 1)$ and $\b = \b(n, q, \a,  \r) \in (0, 1)$ such that the following is true:
If $V \in {\mathcal S}_{\a}$, ${\mathbf C} \in {\mathcal C}_{q}$ satisfy
Hypotheses~\ref{hyp}, Hypothesis~($\star$) with $M = \frac{3}{2}M_{0}^{3}$ and Hypothesis~($\star\star$), and if the induction hypotheses $(H1)$, $(H2)$ hold, then for each $Z = (\z^{1}, \z^{2}, \eta) \in {\rm spt} \, \|V\| \cap ({\mathbf R} \times B_{3/8})$ with $\Theta \, (\|V\|, Z) \geq q$ and each $\m \in (0, 1)$ we have the following:
\begin{eqnarray*}
&&\hspace{-.3in}({\rm a})\;\;\; |\z^{1}|^{2} + {\hat E}_{V}^{2}|\z^{2}|^{2} \leq C \int_{{\mathbf R} \times B_{1}} {\rm dist}^{2} \, (X, {\rm spt} \, \|{\mathbf C}\|) \, d\|V\|(X);\nonumber\\
&&\hspace{-.3in}({\rm b})\;\;\;\;\int_{B_{5\r/8}^{n+1}(Z)} \frac{{\rm dist}^{2} \, (X, {\rm spt} \, \|T_{Z \, \#} \, {\mathbf C}\|)}{|X - Z|^{n+2 - \m}} \, d\|V\|(X)\nonumber\\ 
&&\hspace{2in}\leq \widetilde{C}\r^{-n-2+\m}\int_{{\mathbf R} \times B_{\r}(\z^{2}, \eta)} {\rm dist}^{2} \, (X, {\rm spt} \,\|T_{Z\, \#} \, {\mathbf C}\|) \, d\|V\|(X).\nonumber\\
\end{eqnarray*}
Here $T_{Z} \, : \, {\mathbf R}^{n+1} \to {\mathbf R}^{n+1}$ is the translation $X \mapsto X + Z$;
$C = C(n, q, \a) \in (0, \infty)$ and $\widetilde{C} = \widetilde{C}(n, q, \a, \m) \in (0, \infty).$ (In particular, $C$, $\widetilde{C}$ do not depend on $\r$.) 
\end{corollary}

Our proof of this corollary will be based on several preliminary results, given below as Lemma~\ref{L2-est-2-L1}, Lemma~\ref{L2-est-2-L2}, Proposition~\ref{L2-est-2-P1}, Lemma~\ref{L2-est-2-L3} and Proposition~\ref{L2-est-2-P2}.  

\begin{lemma}\label{L2-est-2-L1}
For any given $\d \in (0, 1)$, there exist $\e^{\prime} = \e^{\prime}(n, q, \a, \d),$ $\g^{\prime} = \g^{\prime}(n, q, \a, \d) \in (0, 1)$ such that if Hypotheses~\ref{hyp} with $\e^{\prime}$, $\g^{\prime}$ in place of $\e$, $\g$ is satisfied by $V \in {\mathcal S}_{\a}$ and ${\mathbf C} \in {\mathcal K}$, and also Hypothesis~($\star$) with $M = \frac{3}{2}M_{0}^{3}$ is satisfied by $V$, then 

$$|\z^{1}|^{2} + {\hat E}_{V}^{2}|\z^{2}|^{2} < \d{\hat E}_{V}^{2}$$

for each $Z=(\z^{1}, \z^{2}, \eta)  \in {\rm spt} \, \|V\| \cap ({\mathbf R} \times B_{3/8})$ 
with $\Theta \, (\|V\|, Z) \geq q$.
\end{lemma}

\begin{proof}
The lemma follows by arguing by contradiction, using Remark (3) of Section~\ref{outline}, Theorem~\ref{SS}, the Remark at the end of Section~\ref{step2}  and the bounds (\ref{slope-bounds-1}), (\ref{slope-bounds-2}).
\end{proof}

\begin{lemma}\label{L2-est-2-L2}
Let $q$ be an integer $\geq 3.$ For any given $\d \in (0, 1)$, there exist $\e = \e(n, q, \a, \d),$ $\g = \g(n, q, \a, \d)$ and $\b = \b(n, q, \a, \d) \in (0,1)$ such that if 
\begin{itemize}
\item[(a)] $p^{\prime} \in \{5, \ldots, 2q\}$ and
\item[(b)] Hypotheses~\ref{hyp}, Hypothesis~($\star$) with $M = \frac{3}{2}M_{0}^{3}$ and  Hypothesis~($\star\star$) are satisfied with $V \in {\mathcal S}_{\a}$, ${\mathbf C} \in {\mathcal C}_{q}(p^{\prime})$, 
\end{itemize}
then

$$|\z^{1}|^{2} + {\hat E}_{V}^{2}|\z^{2}|^{2} \leq \d\left(Q_{V}^{\star}(p^{\prime}-1)\right)^{2}$$

for each $Z = (\z^{1}, \z^{2}, \eta) \in {\rm spt} \, \|V\| \cap ({\mathbf R} \times B_{3/8})$ with $\Theta \, (\|V\|, Z) \geq q.$
\end{lemma} 

We shall eventually prove this lemma by induction on $p^{\prime},$ but first we need to establish the following:

\begin{proposition}\label{L2-est-2-P1}
Let $q$ be an integer $\geq 2,$ $p \in \{4, \ldots, 2q\}$ and suppose that either
\begin{itemize}
\item[(i)] $p =4$ or
\item[(ii)] $q \geq 3$, $p \geq 5$ and Lemma~\ref{L2-est-2-L2} holds whenever 
$p^{\prime} \in \{5, \ldots, p\}.$
\end{itemize} 
Then Corollary~\ref{L2-est-2} holds whenever ${\mathbf C} \in \cup_{k=4}^{p} {\mathcal C}_{q}(k).$
\end{proposition}

\begin{proof} Let $\e_{0}$, $\g_{0}$ and $\b_{0}$ be the constants given by Theorem~\ref{L2-est-1} taken with $\t = 1/16$ (say). Suppose that the hypotheses of the Proposition are satisfied. Let 
$\r \in (0, 1/4]$ and suppose that the hypotheses of Corollary~\ref{L2-est-2}, for suitably small $\e, \g,\b$ to be determined depending only on $n$, $q$, $\a$  and $\r,$ are satisfied by a varifold $V \in {\mathcal S}_{\a}$ and  a cone ${\mathbf C} \in \cup_{k=4}^{p}{\mathcal C}_{q}(k)$ . 

To show that the conclusions of Corollary~\ref{L2-est-2} follow,  we need to apply Theorem~\ref{L2-est-1} with $\t = 1/16$ and $\eta_{Z, \r \, \#} \, V$ in place of $V$ for any $Z  = (\z^{1}, \z^{2}, \eta) \in {\rm spt} \, \|V\| \cap ({\mathbf R} \times B_{3/8})$  with $\Theta \,(\|V\|, Z) \geq q.$ Thus we need to show that it is possible to choose $\e$, $\g$, $\b$ depending only on $n$, $q$, $\a$, $\r$ such that

Hypotheses~\ref{hyp}, Hypothesis~($\star$)  and Hypothesis~($\star\star$) are satisfied with the varifold $\widetilde{V} = \eta_{Z, \r \, \#} \, V$ in place of $V$, with $\e_{0}$ , $\g_{0}$, $\b_{0}$ in place of $\e$, $\g$, $\b$ respectively and with $M = \frac{3}{2}M_{0}^{4}$. If this is so, then  part (b) of Corollary~\ref{L2-est-2} follows as the result of a direct application of Theorem~\ref{L2-est-1}(d) with $\widetilde{V}$ in place of $V$, and part (a) of Corollary~\ref{L2-est-2} follows from the argument of [\cite{W1}, Lemma~6.21] (which in turn is a minor modification of the corresponding argument of \cite{S}), which also requires application of Theorem~\ref{L2-est-1}(d) with $\widetilde{V}$ in place of $V.$

Hypothesis~\ref{hyp}(1) with $\widetilde{V}$ in place of $V$ follows from Theorem~\ref{flat-varifolds}; Hypothesis~\ref{hyp}(3) with $\widetilde{V}$ in place of $V$ and $\e_{0}$ in place of $\e$ holds if $\e < \r^{n+2}\e_{0}.$ Hypothesis~\ref{hyp}(4) 
with $\widetilde{V}$ in place of $V$ is satisfied since by the Remark at the end of Section~\ref{step2}, we may 
choose $\e=\e(n, q,\a, \r)$, $\g = \g(n, q, \a, \r)$ sufficiently small to ensure that $\{Z \, : \, \Theta \, (\|V\|, Z) \geq q\} \cap ({\mathbf R} \times B_{1/2}) \subset \{|x^{2}| < \r/128\}.$ 

 To verify that Hypotheses~\ref{hyp}(5) is satisfied with $\widetilde{V}$ in place of $V$ and $\g_{0}$ in place of $\g$, we proceed as follows: First, using Theorem~\ref{L2-est-1}(a) with $\t = \r/32$, we note that for $\e = \e(n, q, \a, \r), \g = \g(n, q, \a, \r), \b = \b(n, q, \a, \r)  \in (0, 1)$ sufficiently small,   
\begin{eqnarray}\label{L2-est-2-1}
\r^{-n-2}\int_{{\mathbf R} \times B_{\r}(\z^{2}, \eta)} |x^{1}- \z^{1}|^{2} \, d\|V\|(X) \geq&&\nonumber\\ &&\hspace{-2.5in}\r^{-n-2}\sum_{j=1}^{q}\left(\int_{B_{\r}(\z^{2}, \eta) \cap \{x^{2} < -\r/16\}} |h_{j} + u_{j} - \z^{1}|^{2}  + \int_{B_{\r}(\z^{2}, \eta) \cap  \{x^{2} > \r/16\}} |g_{j} + w_{j} - \z^{1}|^{2}\right)\nonumber\\
&&\hspace{-2.5in}\geq \frac{1}{2}\r^{-n-2}\sum_{j=1}^{q}\left(\int_{B_{\r/2}\cap \{x^{2} < -\r/16\}} |h_{j}|^{2} + \int_{B_{\r/2} \cap \{x^{2} > \r/16\}}|g_{j}|^{2}\right)\nonumber\\ 
&&\hspace{-2in}-\r^{-n-2}\sum_{j=1}^{q}\left(\int_{B_{\r}(\z^{2}, \eta) \cap \{x^{2} < -\r/16\}} |u_{j}|^{2}  + \int_{B_{\r}(\z^{2}, \eta) \cap \{x^{2} > \r/16\}} |w_{j}|^{2}\right) - C\r^{-2}|\z^{1}|^{2} \nonumber\\
&&\hspace{-2in} \geq 2^{-n-3}\overline{C}_{1}\left(\sum_{j=1}^{q} |\lambda_{j}|^{2} + |\mu_{j}|^{2}\right)  -\r^{-n-2}E_{V}^{2} 
-C\r^{-2}|\z^{1}|^{2} 
\end{eqnarray}
where $E_{V}^{2} = \int_{{\mathbf R} \times B_{1}} {\rm dist}^{2} \, (X, {\rm spt} \, \|{\mathbf C}\|) \, d\|V\|(X),$  $\overline{C}_{1}   = \overline{C}_{1}(n) \equiv \int_{B_{1/2} \setminus \{x^{2} > 1/16\}}|x^{2}|^{2} \, d{\mathcal H}^{n}(x^{2}, y),$ $C = C(n, q) \in (0, 1)$ and the rest of the notation is as in Theorem~\ref{L2-est-1}(a). If $\e = \e(n, q, \a, \r),$ $\g = \g(n, q, \a, \r) \in (0, 1)$ are sufficiently small, it follows from (\ref{L2-est-2-1}), (\ref{slope-bounds-2}) and Lemma~\ref{L2-est-2-L1} that 
\begin{equation}\label{L2-est-2-1-0}
{\hat E}_{\widetilde{V}} \geq C{\hat E}_{V} 
\end{equation}
where $C = C(n, q) \in (0, \infty).$ On the other hand, by (\ref{cone-dist}) and (\ref{slope-bounds-1}) we have that 
\begin{eqnarray*}
&&\int_{{\mathbf R} \times B_{1}} {\rm dist}^{2} \, (X, {\rm spt} \, \|{\mathbf C}\|) \, d\|\widetilde{V}\|(X) \leq 2\r^{-n-2}\int_{{\mathbf R} \times B_{1}}{\rm dist}^{2} \, (X, {\rm spt} \, \|{\mathbf C}\|)\, d\|V\|(X)\nonumber\\
&&\hspace{4in} +\, C\r^{-2}\left(|\z^{1}|^{2} + {\hat E}_{V}^{2}|\z^{2}|^{2}\right)
\end{eqnarray*}
where $C = C(n, q) \in (0, \infty)$ and, provided $\e = \e(n, q, \a, \r),$ $\g = \g(n,q,\a,\r)$, $\b=\b(n,q,\a,\r)$ are sufficiently small, 
\begin{eqnarray*}
&&\hspace{-1in}\int_{{\mathbf R} \times \left(B_{1/2} \setminus \{|x^{2}| < 1/16\}\right)} {\rm dist}^{2} \, (X, {\rm spt} \, \|\widetilde{V}\|) \, d\|{\mathbf C}\|(X)\nonumber\\
&&\hspace{-.75in}= \r^{-n-2}\int_{{\mathbf R} \times \left(B_{\r}(Z) \setminus \{|x^{2} - \z^{2}|< \r/16\}\right)} {\rm dist}^{2} (X, {\rm spt} \, \|V\|) \, d\|T_{Z \, \#} \, {\mathbf C}\|(X)\nonumber\\
&&\hspace{-.75in}\leq \r^{-n-2}\int_{{\mathbf R} \times \left(B_{17\r/16}(0,\eta) \setminus \{|x^{2}|< \r/32\}\right)} {\rm dist}^{2} (X, {\rm spt} \, \|V\|) \, d\|T_{Z \, \#} \, {\mathbf C}\|(X)\nonumber\\
&&\hspace{-.75in}\leq \r^{-n-2}\int_{{\mathbf R} \times \left(B_{5/8}(0) \setminus \{|x^{2}|< \r/32\}\right)} {\rm dist}^{2} (X, {\rm spt} \, \|V\|) \, d\|{\mathbf C}\|(X) + C\r^{-2}\left(|\z^{1}|^{2} + {\hat E}_{V}^{2}|\z^{2}|^{2}\right)\nonumber\\
&&\hspace{-.75in}\leq C\r^{-n-2}\int_{{\mathbf R} \times B_{1}} {\rm dist}^{2} \, (X, {\rm spt} \, \|{\mathbf C}\|) \, d\|V\|(X) + C\r^{-2}\left(|\z^{1}|^{2} + {\hat E}_{V}^{2}|\z^{2}|^{2}\right)
\end{eqnarray*}
where $C = C(n, q) \in (0, \infty),$ the second inequality follows from the area formula and (\ref{slope-bounds-1}), and  the last inequality follows from Theorem~\ref{L2-est-1}(a) applied with $\t = \r/64$. Thus 
\begin{equation}\label{L2-est-2-1-1}
Q_{\widetilde{V}}^{2}({\mathbf C}) \leq C\left(\r^{-n-2}Q_{V}^{2}({\mathbf C}) +  \r^{-2}(|\z^{1}|^{2} + {\hat E}_{V}^{2}|\z^{2}|^{2})\right)
\end{equation}
which, in view of (\ref{L2-est-2-1-0}) and Lemma~\ref{L2-est-2-L1} applied with sufficiently small $\d = \d(n, q, \a, \r) \in (0, 1),$ implies that Hypothesis~\ref{hyp}(5) is satisfied with $\widetilde{V}$ in place of $V$ and $\g_{0}$ in place of $\g$.

To verify that Hypothesis~($\star$) is satisfied with $\widetilde{V}$ in place of $V$ and $M = \frac{3}{2}M_{0}^{4}$, reasoning again as in (\ref{L2-est-2-1}), we see first that for any hyperplane $P$ of the form $P = \{x^{1}  = \lambda x^{2}\}$ with $|\lambda| < 1$, 
\begin{eqnarray}\label{L2-est-2-2}
\r^{-n-2}\int_{{\mathbf R} \times B_{\r}(\z^{2}, \eta)} {\rm dist}^{2} \, (X - Z, P) \, d\|V\|(X) \geq 
 2^{-n-4}\overline{C}_{1}\left(\sum_{j=1}^{q} |\lambda_{j} - \lambda|^{2} + |\mu_{j}-\lambda|^{2}\right)&&\nonumber\\
&&\hspace{-2in}  -\frac{1}{2}\r^{-n-2}E_{V}^{2} - C\r^{-2}|\z^{1} - \lambda\z^{2}|^{2}\nonumber\\
&&\hspace{-5.5in} \geq   2^{-n-4}\overline{C}_{1}{\rm dist}_{\mathcal H}^{2} \, ({\rm spt} \, \|{\mathbf C}\| \cap ({\mathbf R} \times B_{1}), P \cap ({\mathbf R} \times B_{1})) - \frac{1}{2}\r^{-n-2}E_{V}^{2} - C\r^{-2}|\z^{1} - \lambda\z^{2}|^{2}\nonumber\\
&&\hspace{-5.5in} \geq 2^{-n-4}\omega_{n}^{-1}(2q + 1)^{-1}\overline{C}_{1} \int_{{\mathbf R} \times B_{1}}{\rm dist}^{2} \, (X, P) \, d\|V\|(X)\nonumber\\
&&\hspace{-4in} - \left(2^{-n-3}\omega_{n}^{-1}(2q +1)^{-1}\overline{C}_{1} + 2^{-1}\r^{-n-2}\right)E_{V}^{2} - C\r^{-2}|\z^{1} - \lambda\z^{2}|^{2}
\end{eqnarray}
where $C = C(n, q) \in (0, \infty)$ and we have used the triangle inequality in the last step. On the other hand, noting, by the Constancy Theorem, that
$\left(\omega_{n}(2\r)^{n}\right)^{-1}\|V\|\left({\mathbf R} \times B_{2\r}(0, \eta)\right) \leq q + 1/2$ provided $\e = \e(n, q, \r) \in (0, 1)$ is sufficiently small, we see by Lemma~\ref{L2-est-2-L1} and the triangle inequality again that 
\begin{eqnarray}\label{L2-est-2-3}
\r^{-n-2}\int_{{\mathbf R} \times B_{\r}(\z^{2}, \eta)}|x^{1} - \z^{1}|^{2} \, d\|V\|(X)&&\nonumber\\ 
&&\hspace{-2.5in}\leq 2\r^{-n-2}\|V\|({\mathbf R} \times B_{2\r}(0, \eta)){\rm dist}^{2}_{\mathcal H} \,({\rm spt} \, \|{\mathbf C}\| \cap ({\mathbf R} \times B_{2\r}), \{0\} \times B_{2\r}) +  2\r^{-n-2}E_{V}^{2} + C\r^{-2}\d{\hat E}_{V}^{2}\nonumber\\
&&\hspace{-1in} \leq \left(2^{n+2}\omega_{n}(2q + 1)c_{1}^{2} + C\r^{-2}\d\right){\hat E}_{V}^{2} +  2 \r^{-n-2}E_{V}^{2},  
\end{eqnarray}
where $c_{1} = c_{1}(n)$ is as in (\ref{slope-bounds-1}). Since ${\hat E}_{V}^{2} \leq \frac{3}{2}M^{3}_{0}\inf_{\{P = \{x^{1} = \lambda x^{2}\}\}} \, \int_{{\mathbf R} \times B_{1}} {\rm dist}^{2} \, (X, P) \, d\|V\|(X)$ by hypothesis (of Corollary~\ref{L2-est-2}), in view of the fact that
$$\inf_{\{P = \{x^{1}  = \lambda x^{2}\}\}}\int_{{\mathbf R} \times B_{1}}{\rm dist}^{2} \, (X, P) \, d\|V\|(X) = \inf_{\{P = \{x^{1}  = \lambda x^{2}\}, \, |\lambda| < C{\hat E}_{V}\}}\int_{{\mathbf R} \times B_{1}}{\rm dist}^{2} \, (X, P) \,d\|V\|(X)$$ 
where $C = C(n) \in (0, \infty),$ we deduce from Lemma~\ref{L2-est-2-L1}, (\ref{L2-est-2-2}) and (\ref{L2-est-2-3}) that 
\begin{eqnarray*}
{\hat E}_{\widetilde{V}}^{2} \leq \frac{\left(2^{2n+6}\omega_{n}^{2}(2q+1)^{2}c_{1}^{2} + 2^{n+4}\omega_{n}(2q+1)(2\r^{-n-2}\g + C\r^{-2}\d)\right)\frac{3}{2}M^{3}_{0}}{\overline{C}_{1} - 2^{n+4}\omega_{n}(2q+1)\left((2^{-1}\r^{-n-2} + 2^{-n-3}\omega_{n}^{-1}(2q+1)^{-1}\overline{C}_{1})\g + C\r^{-2}\d\right) \frac{3}{2}M^{3}_{0}} \times\nonumber\\
&&\hspace{-2.5in}\times \; \inf_{\{P = \{x^{1} = \lambda x^{2}\}\}}\int_{{\mathbf R} \times B_{1}}{\rm dist}^{2} \, (X, P) \, d\|\widetilde{V}\|(X)\nonumber\\
&&\hspace{-4in} \leq \frac{3}{2}M^{4}_{0} \inf_{\{P = \{x^{1}  = \lambda x^{2}\}\}}\int_{{\mathbf R} \times B_{1}}{\rm dist}^{2} \, (X, P) \, d\|\widetilde{V}\|(X)
\end{eqnarray*}
provided $\e = \e(n, q, \a, \r),$ $\g = \g(n, q, \a, \r) \in (0, 1)$ are sufficiently small. 

It only remains to verify that Hypothesis~($\star\star$) with $\widetilde{V}$ in place of $V$ and $\b_{0}$ in place of $\b$ is satisfied whenever ${\mathbf C} \in \cup_{k=4}^{p}{\mathcal C}_{q}(k)$. If $p=4$ then ${\mathbf C} \in {\mathcal C}_{q}(4)$ and there is nothing further to verify, so assume that $q \geq 3$ and ${\mathbf C} \in {\mathcal C}_{q}(p^{\prime})$ for some $p^{\prime} \in \{5,\ldots, p\}$. Then for any 
${\mathbf C}^{\prime} \in \cup_{k=4}^{p^{\prime}-1} {\mathcal C}_{q}(k)$, we have by the definition of $Q_{V}^{\star}(p^{\prime}-1)$, the triangle inequality and Hypothesis~($\star\star$) (for $V$ and ${\mathbf C}$ with sufficiently small $\b$) that ${\rm dist}^{2}_{\mathcal H} \, \left({\rm spt} \, \|{\mathbf C}^{\prime}\| \cap ({\mathbf R} \times B_{1}), {\rm spt} \, \|{\mathbf C}\| \cap ({\mathbf R} \times B_{1})\right) \geq C\left(Q_{V}^{\star}(p^{\prime}-1)\right)^{2}$
where $C = C(n, q) \in (0, \infty),$ and hence by Theorem~\ref{L2-est-1}(a), for sufficiently small $\e = \e(n, q, \a, \r),$ $\g = \g(n, q, \a, \r),$ $\b = \b(n, q, \a, \r) \in (0, 1)$, that 
\begin{eqnarray}\label{L2-est-2-4}
\int_{{\mathbf R} \times B_{1}} {\rm dist}^{2} \, (X, {\rm spt}\, \|{\mathbf C}^{\prime}\|) \, d\|\widetilde{V}\|(X)&&\nonumber\\
&&\hspace{-2.5in}\geq \sum_{j=1}^{q}\r^{-n-2}\int_{B_{\r}(\z^{2}, \eta) \cap \{x^{2} < -\r/16\}} {\rm dist}^{2} \, ((h^{j}(x^{2},y)+ u^{j}(x^{2}, y), x^{2}, y) - Z, {\rm spt}\, \|{\mathbf C}^{\prime}\|) \, d{\mathcal H}^{n}(x^{2}, y)\nonumber\\
&&\hspace{-2in} + \,\sum_{j=1}^{q}\r^{-n-2}\int_{B_{\r}(\z^{2}, \eta) \cap \{x^{2} > \r/16\}} {\rm dist}^{2} \, ((g^{j}(x^{2},y) + w^{j}(x^{2},y), x^{2}, y) - Z, {\rm spt}\, \|{\mathbf C}^{\prime}\|) \, d{\mathcal H}^{n}(x^{2}, y)\nonumber\\
&&\hspace{-2.5in}\geq \sum_{j=1}^{q}\r^{-n-2}\int_{B_{\r/2}(0, \eta) \cap\{x^{2} < -\r/16\}} {\rm dist}^{2} \, ((h^{j}(x^{2},y)+ u^{j}(x^{2}, y), x^{2}, y), {\rm spt}\, \|{\mathbf C}^{\prime}\|) \, d{\mathcal H}^{n}(x^{2}, y)\nonumber\\
&&\hspace{-2in} + \,\sum_{j=1}^{q}\r^{-n-2}\int_{B_{\r/2}(0, \eta) \cap \{x^{2} > \r/16\}} {\rm dist}^{2} \, ((g^{j}(x^{2},y) + w^{j}(x^{2},y), x^{2}, y), {\rm spt}\, \|{\mathbf C}^{\prime}\|) \, d{\mathcal H}^{n}(x^{2}, y)\nonumber\\ 
&&\hspace{2in}-\, C^{\prime}\r^{-2}(|\z^{1}|^{2} + \d({\mathbf C}^{\prime})|\z^{2}|^{2})\nonumber\\
&&\hspace{-1.5in} \geq C\left(Q_{V}^{\star}(p-1)\right)^{2} - \r^{-n-2}E_{V}^{2} -C^{\prime}\r^{-2}(|\z^{1}|^{2} + \d({\mathbf C}^{\prime})|\z^{2}|^{2})\nonumber\\ 
&&\hspace{-1in}\geq \frac{1}{2}C\left(Q_{V}^{\star}(p^{\prime}-1)\right)^{2} - C^{\prime}\r^{-2}(|\z^{1}|^{2} + \d({\mathbf C}^{\prime})|\z^{2}|^{2})
\end{eqnarray}
where $C = C(n, q),$ $C^{\prime} = C^{\prime}(n, q) \in (0, \infty)$ and $\d({\mathbf C}^{\prime}) = {\rm dist}^{2}_{\mathcal H} \, \left({\rm spt} \, \|{\mathbf C}^{\prime}\| \cap ({\mathbf R} \times B_{1}), \{0\} \times B_{1}\right).$ Since ${\hat E}^{2}_{\widetilde{V}} \leq C\r^{-n-2}{\hat E}_{V}^{2}$ where $C = C(n, q) \in (0, \infty)$, we have that 
\begin{equation*}
Q_{\widetilde{V}}^{\star}(p^{\prime}-1) = \inf_{\{{\mathbf C}^{\prime} \in \cup_{k=4}^{p^{\prime}-1}{\mathcal C}_{q}(k): \d({\mathbf C}^{\prime}) < C\r^{-n-2}{\hat E}_{V}^{2}\}} \, Q_{\widetilde{V}}({\mathbf C}^{\prime}),
\end{equation*} 
so it follows from  (\ref{L2-est-2-4}) that 
\begin{eqnarray}\label{L2-est-2-5}
\left(Q_{\widetilde{V}}^{\star}(p^{\prime}-1)\right)^{2}\geq C\left(Q_{V}^{\star}(p^{\prime}-1)\right)^{2} - C^{\prime}\r^{-n-4}(|\z^{1}|^{2} + {\hat E}_{V}^{2}|\z^{2}|^{2})
\end{eqnarray}
where $C = C(n, q),$ $C^{\prime} = C^{\prime}(n, q) \in (0, \infty)$. On the other hand, by (\ref{L2-est-2-1-1}), 
\begin{eqnarray}\label{L2-est-2-6}
Q_{\widetilde{V}}^{2}({\mathbf C}) \leq C_{1}\left(\r^{-n-2}Q_{V}^{2}({\mathbf C}) +  \r^{-2}(|\z^{1}|^{2} + {\hat E}_{V}^{2}|\z^{2}|^{2})\right)&&\nonumber\\
&&\hspace{-2.5in} \leq C_{1}\b\r^{-n-2}\left(Q_{V}^{\star}(p^{\prime}-1)\right)^{2} + C_{1}\r^{-2}(|\z^{1}|^{2} + {\hat E}_{V}^{2}|\z^{2}|^{2})
\end{eqnarray}
where $C_{1} = C_{1}(n, q) \in (0, \infty)$. Since by assumption  Lemma~\ref{L2-est-2-L2} holds whenever $p^{\prime} \in \{5, \ldots, p\},$ we may apply Lemma~\ref{L2-est-2-L2} with any $\d = \d(n, q, \a, \r) \in (0, 1)$ satisfying $\max \{C^{\prime}, \b_{0}^{-1}C_{1}\}\r^{-n-4}\d < C/2$, where $C$, $C^{\prime}$, $C_{1}$ are as in (\ref{L2-est-2-5}) and (\ref{L2-est-2-6}), to conclude that 
Hypothesis~($\star\star$) with $\widetilde{V},$ $\b_{0}$ in place of $V$, $\b$ is satisfied.
The proof of the proposition is thus complete.
\end{proof}

\begin{lemma}\label{L2-est-2-L3}
Let $q \geq 3$ and $\d \in (0, 1)$. There exist $\widetilde{\e}_{1} = \widetilde{\e}_{1}(n, q, \a, \d), 
\widetilde{\b}_{1} = \widetilde{\b}_{1}(n, q, \a, \d), {\g}_{1} = {\g}_{1}(n, q, \a, \d)$ and ${\b}_{1} = {\b}_{1}(n, q, \a, \d) \in (0, 1)$ such that if 
\begin{itemize}
\item[(a)] $p^{\prime} \in \{5, \ldots, 2q\};$ 
\item[(b)] Hypotheses~\ref{hyp}(1)-(4), Hypothesis~($\star$), Hypothesis~($\star\star$) are satisfied with $V \in {\mathcal S}_{\a},$ ${\mathbf C} \in {\mathcal C}_{q}(p^{\prime})$, $M = \frac{3}{2}M_{0}^{3}$ and with $\widetilde{\e}_{1}, \widetilde{\b}_{1}$ in place of $\e, \b$ respectively; 
\item[(c)] either
\begin{itemize}
\item[(i)] $\left(Q_{V}^{\star}(4)\right)^{2} \leq \g_{1}{\hat E}_{V}^{2}$ or
\item[(ii)] $p^{\prime} \in \{6, \ldots, 2q\}$,  
$\left(Q_{V}^{\star}(p^{\prime}-j^{\prime})\right)^{2} \leq {\b}_{1}\left(Q_{V}^{\star}(p^{\prime}-j^{\prime}-1)\right)^{2}$
and $\left(Q_{V}^{\star}(p^{\prime} - j^{\prime})\right)^{2} \leq \g_{1}{\hat E}_{V}^{2}$ for some 
$j^{\prime} \in \{1, \ldots, p^{\prime} - 5\},$
\end{itemize}
\end{itemize}
then for each $Z  = (\z^{1}, \z^{2}, \eta) \in {\rm spt} \, \|V\| \cap ({\mathbf R} \times B_{3/8})$ with 
$\Theta \, (\|V\|, Z) \geq q,$
$$|\z_{1}|^{2} + {\hat E}_{V}^{2}|\z_{2}|^{2} < \d\left(Q_{V}^{\star}(4)\right)^{2}$$ 
in case {\rm (c)(i)} holds, and 
$$|\z_{1}|^{2} + {\hat E}_{V}^{2}|\z_{2}|^{2} < \d\left(Q_{V}^{\star}(p^{\prime}-j^{\prime})\right)^{2}$$ 
in case {\rm (c)(ii)} holds.
\end{lemma}

This lemma will follow, in view of the following proposition, from our inductive proof of Lemma~\ref{L2-est-2-L2} given below.

\begin{proposition}\label{L2-est-2-P2}
Let $q$ be an integer $\geq 3,$ $p \in \{5, \ldots, 2q\}$ and suppose that  either
\begin{itemize}
\item[(i)] $p =5$ or
\item[(ii)] $p \in \{6, \ldots, 2q\}$ and Lemma~\ref{L2-est-2-L2} holds whenever 
$p^{\prime} \in \{5, \ldots, p-1\}.$
\end{itemize} 
Then Lemma~\ref{L2-est-2-L3} holds whenever $p^{\prime}=p.$ 
\end{proposition}

\begin{proof} We argue by contradiction. Fix $p \in \{5, \ldots, 2q\}$ and suppose that the hypotheses of the proposition are satisfied. 

Note that if Lemma~\ref{L2-est-2-L3} with $p^{\prime} = p$ does not hold, then there exist a number $\d \in (0, 1),$ an integer $j^{\prime} \in \{1, \ldots, p-5\}$ in case $p \in \{6, \ldots, 2q\}$ and, for each $k=1, 2, \ldots$, a varifold $V_{k} \in {\mathcal S}_{\a}$, a point $Z_{k} = (\z^{1}_{k}, z^{2}_{k}, \eta_{k}) \in {\rm spt} \, \|V_{k}\| \cap ({\mathbf R} \times B_{3/8})$ with $\Theta \, (\|V_{k}\|, Z_{k}) \geq q$, a cone ${\mathbf C}_{k} \in {\mathcal C}_{q}(p)$ such that Hypotheses~\ref{hyp}(1), Hypotheses~\ref{hyp}(2), Hypotheses~\ref{hyp}(4), Hypothesis~($\star$) are satisfied with $V_{k}$  in place of $V$, ${\mathbf C}_{k}$ in place of ${\mathbf C}$, $M = \frac{3}{2}M_{0}^{3};$ 
$${\hat E}_{k} \to 0;$$
\begin{equation}\label{L2-est-2-9-1}
\left(Q_{k}^{\star}(p-1)\right)^{-1} Q_{V_{k}}({\mathbf C}_{k}) \to 0;
\end{equation}
\begin{equation}\label{L2-est-2-9'}
\mbox{either} \;\;  {\hat E}_{k}^{-1} Q_{k}^{\star}(4) \to 0 \;\; \mbox{or}
\end{equation}
\begin{equation}\label{L2-est-2-9''}
p \in \{6, \ldots, 2q\}, \;\; \left(Q_{k}^{\star}(p -j^{\prime}-1)\right)^{-1}Q_{k}^{\star}(p-j^{\prime})  \to 0 \;\; \mbox{and} \;\; {\hat E}_{k}^{-1}Q_{k}^{\star}(p - j^{\prime}) \to 0
\end{equation}
(or both) and yet

\begin{equation}\label{L2-est-2-10'}
|\z^{1}_{k}|^{2} + {\hat E}_{k}^{2}|\z^{2}_{k}|^{2} \geq \d\left(Q_{k}^{\star}(4)\right)^{2} \mbox{in case (\ref{L2-est-2-9'}) holds and}
\end{equation}
\begin{equation}\label{L2-est-2-10}
|\z^{1}_{k}|^{2} + {\hat E}_{k}^{2}|\z^{2}_{k}|^{2} \geq \d\left(Q_{k}^{\star}(p- j^{\prime})\right)^{2} \mbox{in case (\ref{L2-est-2-9''}) holds},
\end{equation}
where we have used the notation ${\hat E}_{k} = {\hat E}_{V_{k}}$ and $Q_{k}^{\star}(\cdot) = Q_{V_{k}}^{\star}(\cdot).$

For each $k=1, 2, \ldots,$ let $\overline{\mathbf C}_{k} \in {\mathcal K}$ be chosen as follows: in case (\ref{L2-est-2-9'}) holds, $\overline{\mathbf C}_{k} \in {\mathcal C}_{q}(4)$ is such that $\left(Q_{V_{k}}(\overline{\mathbf C}_{k})\right)^{2} < \frac{3}{2}\left(Q_{V_{k}}^{\star}(4)\right)^{2};$ 
in case (\ref{L2-est-2-9''}) holds, $\overline{\mathbf C}_{k} \in {\mathcal C}_{q}(p-j^{\prime})$ is such that $\left(Q_{V_{k}}(\overline{\mathbf C}_{k})\right)^{2} < \frac{3}{2}\left(Q_{V_{k}}^{\star}(p-j^{\prime})\right)^{2}.$ Note that since the rest of our argument is the same for either case, we  
use the same notation $\overline{\mathbf C}_{k}$ for either case.     
Let $\t_{k} \in (0, 1/8)$ be such that $\t_{k} \searrow 0^{+}.$ By passing to appropriate subsequences without changing notation, we have by Proposition~\ref{L2-est-2-P1} and Corollary~\ref{L2-est-2} that for each $k=1, 2, \ldots,$
\begin{equation}\label{L2-est-2-12}
|\z_{k}^{1}|^{2} + {\hat E}_{k}^{2}|\z_{k}^{2}|^{2} \leq CE_{k}^{2}
\end{equation}
where $C = C(n, q, \a) \in (0, \infty)$, and for each $\m \in (0, 1)$, 
\begin{eqnarray}\label{L2-est-2-13}
\sum_{j=1}^{q} \int_{B_{1/8}(\z_{k}^{2}, \eta_{k}) \cap \{x^{2} < -\t_{k}/4\}} \frac{|u_{j}^{k}(x^{2}, y) - (\z_{k}^{1} - \overline{\lambda}^{k}_{j}\z_{k}^{2})|^{2}}{|(h_{j}^{k}(x^{2}, y) + u^{k}_{j}(x^{2}, y), x^{2}, y) - (\z_{k}^{1}, \z_{k}^{2}, \eta_{k})|^{n+2-\m}} \, dx^{2} \, dy&&\nonumber\\ 
&&\hspace{-5.2in}+ \;\;\sum_{j=1}^{q}\int_{B_{1/8}(\z_{k}^{2}, \eta_{k}) \cap \{x^{2} > \t_{k}/4\}}\frac{|w^{k}_{j}(x^{2}, y) - (\z_{k}^{1} - \overline{\mu}^{k}_{j}\z^{2})|^{2}}{|(g_{j}^{k}(x^{2}, y) + w^{k}_{j}(x^{2}, y), x^{2}, y) - (\z_{k}^{1}, \z_{k}^{2}, \eta_{k})|^{n+2-\m}} \, dx^{2}dy \leq \widetilde{C}E_{k}^{2}
\end{eqnarray}
where $\widetilde{C} = \widetilde{C}(n, q, \a, \m) \in (0, \infty).$ Here, $E_{k}^{2} = \int_{{\mathbf R} \times B_{1}} {\rm dist}^{2} \, (X,{\rm spt} \, \|\overline{\mathbf C}_{k}\|) \, d\|V_{k}\|(X)$; for each $j \in \{1, 2, \ldots, q\}$ and 
$k=1, 2, \ldots$, the functions $u_{j}^{k}$, $w_{j}^{k}$ correspond to 
$u_{j}$, $w_{j}$ of Theorem~\ref{L2-est-1}(a) when $V_{k},$ $\overline{\mathbf C}_{k}$ are taken in place of $V,$ ${\mathbf C}$,
and the numbers $\overline{\lambda}^{k}_{j}$, $\overline{\mu}^{k}_{j}$ correspond to $\lambda_{j}, \mu_{j}$ of Hypothesis~\ref{hyp}(2) when  $\overline{\mathbf C}_{k}$ is taken in place of ${\mathbf C}$. Note then that  
$\overline{\lambda}^{k}_{1} \geq \overline{\lambda}^{k}_{2}\geq \ldots \geq \overline{\lambda}^{k}_{q}$, 
$\overline{\m}^{k}_{1} \leq \overline{\m}^{k}_{2}\leq \ldots \leq \overline{\m}^{k}_{q}$ and by (\ref{slope-bounds-1}) and (\ref{slope-bounds-2}), 
\begin{equation}\label{L2-est-2-14}
c{\hat E}_{k} \leq \max \, \{|\overline{\lambda}_{1}^{k}|, |\overline{\lambda}_{q}^{k}|\} \leq c_{1}{\hat E}_{k}, \;\; c{\hat E}_{k} \leq \max \, \{|\overline{\mu}_{1}^{k}|, |\overline{\mu}_{q}^{k}|\} \leq c_{1}{\hat E}_{k}, \;\; \min \, \{|\overline{\lambda}_{1}^{k} - \overline{\lambda}_{q}^{k}|, |\overline{\m}_{1}^{k} - \overline{\m}_{q}^{k}|\} \geq  2c{\hat E}_{k}
\end{equation}
where $c_{1} = c_{1}(n), c = c(n, q) \in (0, \infty)$ are as in (\ref{slope-bounds-1}) and (\ref{slope-bounds-2}). 

Writing $Q_{k} = Q_{V_{k}}(\overline{\mathbf C}_{k})$, we see 
by Theorem~\ref{L2-est-1}(a) and elliptic estimates that for each $j \in \{1, 2, \ldots, q\}$, there exist harmonic 
functions $\varphi_{j}  \, : \, B_{3/4} \cap \{x^{2} < 0\} \to {\mathbf R}$, $\psi_{j} \, : \, B_{3/4} \cap \{x^{2} > 0\} \to {\mathbf R}$ such that 
$Q_{k}^{-1}u_{j}^{k} \to \varphi_{j}$, $Q_{k}^{-1}w^{k}_{j} \to \psi_{j}$ where the convergence is in $C^{2}(K)$ for each compact subset $K$ of the respective domains of $\varphi_{j}$, $\psi_{j}$. By (\ref{L2-est-2-9-1}), $Q_{k}^{-1}Q_{V_{k}}({\mathbf C}_{k}) \to 0$, which implies that for each $j \in \{1, 2, \ldots, q\}$, 
there exist constants $\overline{\lambda}_{j}$, $\overline{\m}_{j}$ such that $\varphi_{j}(x^{2}, y) = \overline{\lambda}_{j}x^{2}$ for $(x^{2}, y) \in B_{1/2} \cap \{x^{2} <0\}$ and 
$\psi_{j}(x^{2}, y) = \overline{\m}_{j}x^{2}$ for $(x^{2}, y) \in B_{1/2} \cap \{x^{2} > 0\}.$ We find a point $\eta \in \{0\}\times {\mathbf R}^{n-1} \cap \overline{B_{3/8}(0)}$ and, by (\ref{L2-est-2-12}) and (\ref{L2-est-2-14}), numbers $\k_{1}, \k_{2}$, $\ell_{1}, \ldots, \ell_{q}, m_{1}, \ldots, m_{q}$  such that, passing to further subsequences without changing notation, $\eta_{k} \to \eta$, $Q_{k}^{-1}\z^{1}_{k} \to \k_{1}$, $Q_{k}^{-1}{\hat E}_{k}\z^{2}_{k} \to \k_{2}$,  ${\hat E}_{k}^{-1}\lambda_{j}^{k} \to \ell_{j}$ and ${\hat E}_{k}^{-1}\m_{j}^{k} \to m_{j}.$ 
We deduce from (\ref{L2-est-2-13}) that 
\begin{eqnarray*}
\sum_{j=1}^{q} \int_{B_{1/8}(0,\eta) \cap \{x^{2} <0\}} \frac{|\overline{\lambda}_{j}x^{2} - (\k_{1} - \ell_{j}\k_{2})|^{2}}{\left(|x^{2}|^{2} + |y -\eta|^{2}\right)^{\frac{n+2-\m}{2}}} \, dx^{2} dy&&\nonumber\\
&&\hspace{-1in} + \; \sum_{j=1}^{q}\int_{B_{1/8}(0,\eta) \cap \{x^{2} >0\}}\frac{|\overline{\m}_{j}x^{2} - (\k_{1} - m_{j}\k_{2})|^{2}}{\left(|x^{2}|^{2}+|y - \eta|^{2}\right)^{\frac{n+2-\m}{2}}} \, dx^{2}dy < \infty
\end{eqnarray*}
which readily implies that $\k_{1} - \ell_{j}\k_{2}=0$ and $\k_{1} - m_{j}\k_{2} = 0$ for each $j=1, 2, \ldots, q.$ Since by 
(\ref{L2-est-2-14}) not all $\ell_{1}, \ldots, \ell_{q}$ are equal, we must have that $\k_{1}= \k_{2} = 0.$ This contradicts (\ref{L2-est-2-10'}) in case (\ref{L2-est-2-9'}) holds and (\ref{L2-est-2-10}) in case 
(\ref{L2-est-2-9''}) holds. The proposition is thus proved.
\end{proof}

\noindent
\begin{proof}[Proof of Lemma~\ref{L2-est-2-L2}] We prove the lemma by induction on $p^{\prime}.$ Let $\d \in (0, 1)$ and consider first the case $p^{\prime} = 5$. Noting, in view of Proposition~\ref{L2-est-2-P2}, the validity of Lemma~\ref{L2-est-2-L3} with $p^{\prime} = 5,$ let $\widetilde{\e}_{1} = \widetilde{\e}_{1}(n, q, \a, \d)$, $\widetilde{\b}_{1} = \widetilde{\b}_{1}(n, q, \a, \d),$ $\g_{1} = \g_{1}(n, q, \a, \d)$, $\b_{1} = \b_{1}(n, q, \a, \d)$ be as in Lemma~\ref{L2-est-2-L3} with $p^{\prime} = 5$, and  
suppose that the hypotheses of Lemma~\ref{L2-est-2-L2} with $p^{\prime} = 5$ are satisfied by some $V \in {\mathcal S}_{\a}$ and ${\mathbf C} \in {\mathcal C}_{q}(5),$ with 
$\e = \min\{\widetilde{\e}, \e^{\prime}(n, q, \a, \d\g_{1})\}$, $\b = \widetilde{\b}_{1}$ and 
$\g = \min\{\g_{1}, \g^{\prime}(n, q, \a, \d\g_{1})\}$, where $\e^{\prime} = \e^{\prime}(n, q, \a, \cdot)$, $\g^{\prime} = \g^{\prime}(n, q, \a, \cdot)$ are as in Lemma~\ref{L2-est-2-L1}. Then hypotheses (a) and (b) of Lemma~\ref{L2-est-2-L3} with $p^{\prime} = 5$ are satisfied by $V$ and ${\mathbf C}.$ If also $\left(Q_{V}^{\star}(4)\right)^{2} \leq \g_{1}{\hat E}_{V}^{2}$, then by Lemma~\ref{L2-est-2-L3} we have the desired conclusion. If on the other hand $(Q_{V}^{\star}(4))^{2} > \g_{1}{\hat E}_{V}^{2},$ then applying Lemma~\ref{L2-est-2-L1} with $\d\g_{1}$ in place of $\d$, we again have the desired conclusion. So Lemma~\ref{L2-est-2-L2} is established in case $p^{\prime} = 5.$   

Fix now  $p \in \{6, \ldots, 2q\}$ and suppose by induction that Lemma~\ref{L2-est-2-L2} holds whenever $p^{\prime} \in \{5, \ldots, p-1\}.$  Then by Proposition~\ref{L2-est-2-P2}, Lemma~\ref{L2-est-2-L3}  with $p^{\prime} = p$ holds. Let $\d \in (0, 1),$ and let 
$\widetilde{\e}_{1}(n, q, \a, \cdot)$, $\widetilde{\b}_{1}(n, q, \a, \cdot),$ $\g_{1}(n, q, \a, \cdot)$, 
$\b_{1} = \b_{1}(n, q, \a, \cdot)$ be as in Lemma~\ref{L2-est-2-L3} with $p^{\prime} = p.$ Set $\b_{1}^{(0)} = 1$ and for $j = 1, 2, 3, \ldots p-5,$ set $\b_{1}^{(j)}  = \b_{1}(n, q, \a, \d\Pi_{k=0}^{j-1}\b_{1}^{(k)}),$ $\widetilde{\e}_{1}^{(j)} = \widetilde{\e}_{1}(n, q, \a, \d\Pi_{k=1}^{j}\b_{1}^{(k)})$, 
$\widetilde{\b}_{1}^{(j)} = \widetilde{\b}_{1}(n, q, \a, \d\Pi_{k=1}^{j}\b_{1}^{(k)})$ and 
$\g_{1}^{(j)} = \g_{1}(n, q, \a, \d\Pi_{k=1}^{j}\b_{1}^{(k)}).$ Let again $\e^{\prime} = \e^{\prime}(n, q, \a, \cdot)$, $\g^{\prime} = \g^{\prime}(n, q, \a, \cdot)$ be as in Lemma~\ref{L2-est-2-L1}, 
let $\d^{\prime} = \Pi_{j=1}^{p-5}\g_{1}^{(j)}\b_{1}^{(j)}\widetilde{\b}_{1}^{(j)}$ and let 
$\overline{\e} = \e^{\prime}(n, q, \a, \d\d^{\prime})$, 
$\overline{\g} = \g^{\prime}(n,q,\a,\d\d^{\prime}).$

Let ${\mathbf C} \in {\mathcal C}_{q}(p)$, $V \in {\mathcal S}_{\a}$ and suppose that 
the hypotheses of Lemma~\ref{L2-est-2-L2} are satisfied with 
$\e = \min\{\overline{\e}, \widetilde{\e}_{1}^{(j)} \, : \, 1 \leq j \leq p-5\},$ $\b = \min\{\widetilde{\b}_{1}^{(j)} \, : \, 1 \leq j \leq p-5\}$ and 
$\g = \min\{\overline{\g}, \g_{1}^{(j)} \, : \, 1 \leq j \leq p-5\}.$ Consider the following exhaustive list of alternatives:
\noindent
\begin{itemize}
\item[(a)] $\left(Q_{V}^{\star}(p-1)\right)^{2}~>~\d^{\prime}{\hat E}_{V}^{2}$. 
%\noindent
\item[(${\rm b}_{1}$)]  $\left(Q_{V}^{\star}(p-1)\right)^{2}\leq \d^{\prime}{\hat E}_{V}^{2}$ and $\left(Q_{V}^{\star}(p-1)\right)^{2}\leq \b_{1}^{(1)}\left(Q_{V}^{\star}(p-2)\right)^{2}$.
%\noindent
\item[(${\rm b}_{2}$)] $\left(Q_{V}^{\star}(p-1)\right)^{2}\leq\d^{\prime}{\hat E}_{V}^{2},\;\;$  $\left(Q_{V}^{\star}(p-1)\right)^{2}>\b_{1}^{(1)}\left(Q_{V}^{\star}(p-2)\right)^{2}\;\;\;$ and $\;\;\left(Q_{V}^{\star}(p-2)\right)^{2}\leq\b_{1}^{(2)}\left(Q_{V}^{\star}(p-3)\right)^{2}.$
%\noindent
\item[(${\rm b}_{3}$)] $\left(Q_{V}^{\star}(p-1)\right)^{2}\leq \d^{\prime}{\hat E}_{V}^{2},$  $\left(Q_{V}^{\star}(p-1)\right)^{2}>\b_{1}^{(1)}\left(Q_{V}^{\star}(p-2)\right)^{2}$, $\left(Q_{V}^{\star}(p-2)\right)^{2}>\b_{1}^{(2)}\left(Q_{V}^{\star}(p-3)\right)^{2}$ and 
$\left(Q_{V}^{\star}(p-3)\right)^{2}~\leq~\b_{1}^{(3)}\left(Q_{V}^{\star}(p-4)\right)^{2}.$
%\noindent
\item[] $\ldots$
%\noindent
\item[(${\rm b}_{p-5}$)] $\left(Q_{V}^{\star}(p-1)\right)^{2}\leq \d^{\prime}{\hat E}_{V}^{2},$  $\left(Q_{V}^{\star}(p-1)\right)^{2} > \b_{1}^{(1)}\left(Q_{V}^{\star}(p-2)\right)^{2}$, $\left(Q_{V}^{\star}(p-2)\right)^{2} > \b_{1}^{(2)}\left(Q_{V}^{\star}(p-3)\right)^{2},$  
$\left(Q_{V}^{\star}(p-3)\right)^{2} > \b_{1}^{(3)}\left(Q_{V}^{\star}(p-4)\right)^{2},$ $\ldots,$ $\left(Q_{V}^{\star}(6)\right)^{2} > \b_{1}^{(p-6)} \left(Q_{V}^{\star}(5)\right)^{2}\;$ and $\;\left(Q_{V}^{\star}(5)\right)^{2} \leq \b_{1}^{(p-5)} \left(Q_{V}^{\star}(4)\right)^{2}.$
%\noindent
\item[(c)] $\left(Q_{V}^{\star}(p-1)\right)^{2} \leq \d^{\prime}{\hat E}_{V}^{2},$  $\left(Q_{V}^{\star}(p-1)\right)^{2} > \b_{1}^{(1)}\left(Q_{V}^{\star}(p-2)\right)^{2}$, $\left(Q_{V}^{\star}(p-2)\right)^{2} >\b_{1}^{(2)}\left(Q_{V}^{\star}(p-3)\right)^{2},$  
$\left(Q_{V}^{\star}(p-3)\right)^{2} > \b_{1}^{(3)}\left(Q_{V}^{\star}(p-4)\right)^{2}$, $\ldots,$ $\left(Q_{V}^{\star}(6)\right)^{2} > \b_{1}^{(p-6)} \left(Q_{V}^{\star}(5)\right)^{2}\;$ and $\;\left(Q_{V}^{\star}(5)\right)^{2} > \b_{1}^{(p-5)} \left(Q_{V}^{\star}(4)\right)^{2}.$
\end{itemize}
\noindent
The conclusion of Lemma~\ref{L2-est-2-L2}  in case of alternative (a) follows from Lemma~\ref{L2-est-2-L1} applied with $\d\d^{\prime}$  in place of $\d$; the conclusion of Lemma~\ref{L2-est-2-L2}  in case of alternative (b$_{1}$)  follows from Lemma~\ref{L2-est-2-L3} applied with $p^{\prime} = p$ and $j=1$; the conclusion of Lemma~\ref{L2-est-2-L2} in case of alternative (b$_{2}$) follows from Lemma~\ref{L2-est-2-L3} applied with $p^{\prime} = p$, $j=2$ and 
$\d\b_{1}^{(1)}$ in place of $\d$; similarly, the conclusion of Lemma~\ref{L2-est-2-L2} in case of any of the alternatives 
(b$_{3}$)-(b$_{p-5}$) follows from an application of Lemma~\ref{L2-est-2-L3} with $p^{\prime} = p$ and appropriate value of $j$ and $\d$; the conclusion of Lemma~\ref{L2-est-2-L2} in case of alternative (c)  follows from Lemma~\ref{L2-est-2-L3} applied with $p^{\prime} = 5$ and $\d\Pi_{k=1}^{p-5}\b_{1}^{(k)}$ in place of $\d$. Thus the inductive poof of Lemma~\ref{L2-est-2-L2} is complete.
\end{proof}

\begin{proof}[Proof of Lemma~\ref{L2-est-2-L3}] 
Since we have now established Lemma~\ref{L2-est-2-L2} for all values of $p^{\prime} \in \{5, \ldots, 2q\}$, Lemma~\ref{L2-est-2-L3} follows from Proposition~\ref{L2-est-2-P2}.
\end{proof}

\begin{proof}[Proof of Corollary~\ref{L2-est-2}] 
Again, since Lemma~\ref{L2-est-2-L2} holds for all values of $p^{\prime} \in \{5, \ldots, 2q\}$, Corollary~\ref{L2-est-2} follows from Proposition~\ref{L2-est-2-P1}.
\end{proof}

\noindent
{\bf Remark:} Note that the proof of Corollary~\ref{L2-est-2} establishes that corresponding to each $\e, \g, \b \in (0, 1/2)$ and $\r \in (0, 1/2)$, there exist $\widetilde{\e} = \widetilde{\e}(n, q, \a, \r, \e) \in (0, 1/2)$, $\widetilde{\g} = \widetilde{\g}(n, q, \a, \r, \g) \in (0, 1/2)$, $\widetilde{\b} =\widetilde{\b}(n, q, \a, \r, \b) \in (0, 1/2)$ such that 
the following is true: Let $V \in {\mathcal S}_{\a}$ and ${\mathbf C} \in {\mathcal C}_{q}.$ If Hypotheses~\ref{hyp} are satisfied with $\widetilde{\e}$, $\widetilde{\g}$ in place of $\e$, $\g$ respectively, Hypothesis~($\star$) is satisfied with $M = \frac{3}{2}M_{0}^{2}$ and Hypothesis~($\star\star$) is satisfied with $\widetilde{\b}$ in place of $\b$, and if the induction hypotheses $(H1)$, $(H2)$ hold, then, for each $Z \in {\rm spt} \, \|V\| \cap ({\mathbf R} \times B_{3/8})$, 
Hypotheses~\ref{hyp}, Hypothesis~($\star$) with $M = \frac{3}{2}M_{0}^{3}$ and Hypothesis~($\star\star$) are satisfied with $\eta_{Z, \r \, \#} \, V$ in place of $V$.
  
\begin{lemma}\label{no-gaps}
Let $q$ be an integer $\geq 2$, $\a \in (0, 1)$, $\d \in (0, 1/8)$ and $\m \in (0, 1)$. There exist numbers $\e_{1} = \e_{1}(n, q, \a, \d) \in (0, 1)$, $\g_{1} = \g_{1}(n, q, \a, \d) \in (0, 1)$ and $\b_{1} = \b_{1}(n, q, \a) \in (0, 1)$ such that the following is true:  If $V \in {\mathcal S}_{\a}$, ${\mathbf C} \in {\mathcal C}_{q}$ satisfy Hypotheses~\ref{hyp}, Hypothesis~($\star$)  with $\e_{1}$, $\g_{1}$ in place of $\e$, $\g$ respectively and with $M = \frac{3}{2}M_{0}^{3}$, and if the induction hypotheses $(H1)$, $(H2)$ hold, then 
\begin{itemize}
\item[(a)] $B^{n+1}_{\d}(0,y) \cap \{Z \, : \, \Theta \, (\|V\|, Z) \geq q\} \neq \emptyset$
for each point $(0, y) \in \{0\} \times {\mathbf R}^{n-1} \cap B_{1/2}.$\\
\item[(b)] If additionally $V$, ${\mathbf C}$ satisfy Hypothesis~($\star\star$) with $\b_{1}$ in place of $\b$, then
\begin{eqnarray*}
\int_{B_{1/2}^{n+1}(0) \cap \{|(x^{1}, x^{2})| < \s\}} {\rm dist}^{2} \, (X, {\rm spt} \, \|{\mathbf C}\|) \, d\|V\|(X)&&\nonumber\\ 
&&\hspace{-1in}\leq C_{1}\s^{1-\m}\int_{{\mathbf R} \times B_{1}} {\rm dist}^{2} \, (X, {\rm spt} \, \|{\mathbf C}\|) \, d\|V\|(X)
\end{eqnarray*}
for each $\s \in [\d, 1/4)$, where $C_{1} = C_{1}(n, q, \a, \m) \in (0, \infty).$ (In particular $C_{1}$ is independent of $\d$ and $\s$.)
\end{itemize}
\end{lemma}

\begin{proof}
If part (a) were false, then there exist a number $\d \in (0, 1/2)$ and a sequence of varifolds $\{V_{k}\} \subset {\mathcal S}_{\a}$; a sequence of cones ${\mathbf C}_{k} = \sum_{j=1}^{q} |H_{j}^{k}|+|G_{j}^{k}| \in {\mathcal C}_{q}$
where, for each $k$, 
$H_{j}^{k} = \{(x^{1}, x^{2}, y) \in {\mathbf R}^{n+1} \, : \, x^{2} < 0 \;\; \mbox{and} \;\; x^{1} = \lambda_{j}^{k}x^{2}\}$, $G_{j}^{k} = \{(x^{1}, x^{2}, y) \in {\mathbf R}^{n+1} \, : \, x^{2} > 0 \;\; \mbox{and} \;\; x^{1} = \mu_{j}^{k}x^{2}\}$, 
with $\lambda_{1}^{k} \geq \lambda_{2}^{k} \geq \ldots \geq \lambda_{q}^{k}$ and  
$\mu_{1}^{k} \leq \mu_{2}^{k} \leq \ldots \leq \mu_{q}^{k}$; and a sequence of points $(0, y_{k}) \in \{0\} \times {\mathbf R}^{n-1} \cap B_{1/2}$ 
with $B^{n+1}_{\d}(0,y_{k}) \cap \{Z \, : \, \Theta \, (\|V_{k}\|, Z) \geq q\} = \emptyset$
such that Hypotheses~\ref{hyp} (1), (2), (4) are satisfied with $V_{k}$, ${\mathbf C}_{k}$ in place of $V$, ${\mathbf C}$; 
Hypothesis~($\star$) is satisfied with $M = \frac{3}{2}M_{0}^{3}$ and $V_{k}$ in place of $V$; ${\hat E}_{k} = {\hat E}_{V_{k}} \equiv \sqrt{\int_{{\mathbf R} \times B_{1}} |x^{1}|^{2} d\|V_{k}\|(X)} \to 0$ and
\begin{eqnarray}\label{no-gaps-1}
{\hat E}_{k}^{-2}\int_{{\mathbf R} \times (B_{1/2} \setminus \{|x^{2}| < 1/16\})} {\rm dist}^{2}(X, {\rm spt} \, \|V_{k}\|) \,d\|{\mathbf C}_{k}\|(X)&&\nonumber\\ 
&&\hspace{-2.5in}+\;\; {\hat E}_{k}^{-2}\int_{{\mathbf R} \times B_{1}} {\rm dist}^{2} \, (X, {\rm spt}\, \|{\mathbf C}_{k}\|) \, d\|V_{k}\|(X) \to 0.
\end{eqnarray}

After passing to a subsequence without changing notation, $(0, y_{k}) \to (0, y)$ for some point 
$(0, y) \in \{0\} \times {\mathbf R}^{n-1} \cap {\overline B_{1/2}},$ and hence
$$B^{n+1}_{\d/2}(0,y) \cap \{Z \, : \, \Theta \, (\|V_{k}\|, Z) \geq q\} = \emptyset$$ 
for all sufficiently large $k$.
This implies, by Remark 3 of Section~\ref{outline}, that for all sufficiently large $k$,
${\mathcal H}^{n-7+\g} \, ({\rm sing} \, V_{k}\res(B_{\d/2}^{n+1}(0, y)) = 0$ for each $\g >0$ if 
$n \geq 7$, and ${\rm sing} \, V_{k}\res (B_{\d/2}^{n+1}(0,y)) = \emptyset$ if $2 \leq n \leq 6$, so by 
Theorem~\ref{SS} and elliptic theory,  
\begin{equation}\label{no-gaps-2}
V_{k} \res \left({\mathbf R} \times B_{\d/4}((0, y))\right)= \sum_{j=1}^{q} |{\rm graph} \, u_{j}^{k}|
\end{equation}
for all sufficiently large $k$, where $v_{j}^{k} \in C^{\infty}\left( B_{\d/4}(0, y)\right)$, $u_{1}^{k} \leq u_{2}^{k} \leq \ldots \leq u_{q}^{k}$ on $B_{\d/4}((0, y))$ and $u_{j}^{k}$ are solutions of the minimal surface equation on 
$ B_{\d/4}((0, y))$ satisfying, by 
standard elliptic estimates,
$$\sup_{B_{\d/16}(0,y)} \, |D^{\ell} \, u_{j}^{k}| \leq C{\hat E}_{k}$$ 
for $\ell =0, 1, 2, 3$ and 
$j=1, 2, \ldots, q,$ where $C = C(n, \d).$ Passing to a further subsequence without changing notation, 
we deduce that for each $j=1, 2, \ldots, q$, 
${\hat E}_{k}^{-1} u_{j}^{k}  \to v_{j}$ in $C^{2} \, (B_{\d/16}(0, y))$ where 
$v_{j}$ are harmonic in $B_{\d/16}(0, y)$ with $v_{1} \leq v_{2} \leq \ldots  v_{q}$ on $B_{\d/16}(0, y).$ By  (\ref{no-gaps-1}), we see that 
$$\left.v_{j} \right|_{B_{\d/16}(0, y) \cap \{x^{2} < 0\}} = \left.{\widetilde h_{j}}\right|_{B_{\d/16}(0, y) \cap \{x^{2} < 0\}} \;\;\; \mbox{and}$$ 
$$\left.v_{j} \right|_{B_{\d/16}(0, y) \cap \{x^{2} > 0\}} = \left.{\widetilde g_{j}}\right|_{B_{\d/16}(0, y) \cap \{x^{2} > 0\}}$$ 
where ${\widetilde h}_{j}$ and ${\widetilde g}_{j}$ are linear functions of the form 
${\widetilde h}_{j}(x^{2}, y) = \widetilde{\lambda}_{j} x^{2}$, $\widetilde{g}_{j}(x^{2}, y) = \widetilde{\mu}_{j}x^{2}$, 
with $\widetilde{\lambda}_{j}, \widetilde{\mu}_{j} \in {\mathbf R}$, 
$\widetilde{\lambda}_{1} \geq \widetilde{\lambda}_{2} \geq \ldots \geq \widetilde{\lambda}_{q}$ and 
$\widetilde{\mu}_{1} \leq \widetilde{\mu}_{2} \leq \ldots \leq \widetilde{\mu}_{q}.$ 
By the maximum principle, we conclude that $\widetilde{\lambda}_{j} =  \widetilde{\mu}_{j} = \lambda$ for some 
$\lambda \in {\mathbf R}$ and
all $j=1, 2, \ldots, q.$ Therefore, by (\ref{no-gaps-1}) again, we see that the coarse blow-up (in the sense of Section~\ref{blow-up}) of 
$\{V_{k}\}$ and that of $\{{\mathbf C}_{k}\}$, are both equal to the hyperplane $x^{1} = \lambda x^{2}.$ But this is impossible in view of (\ref{slope-bounds-2}) so the assertion of part (a)
must hold.

To see part (b), argue as in [\cite{S}, Corollary 3.2 (ii)] noting that by (\ref{cone-dist}), (\ref{slope-bounds-1}) and Corollary~\ref{L2-est-2}(a), we have that for each $Z \in {\rm spt} \, \|V\| \cap ({\mathbf R} \times B_{3/8})$ with $\Theta \, (\|V\|, Z) \geq q$ and any $X \in {\mathbf R}^{n+1}$, 
$$|{\rm dist}\, (X, {\rm spt} \, \|{\mathbf C}\|) - {\rm dist}\, (X, {\rm spt} \, \|T_{Z \, \#} \, {\mathbf C}\|)|^{2} \leq C\int_{{\mathbf R} \times B_{1}}{\rm dist}^{2} \, (\widetilde{X}, {\rm spt} \,\|{\mathbf C}\|) \, d\|V\|(\widetilde{X})$$ 
where $C = C(n, q, \a) \in (0, \infty).$
\end{proof}

\section{Blowing up by  fine excess}\label{fineexcessblowup}
\setcounter{equation}{0}

Let $\{\e_{k}\}$,$\{\g_{k}\}$ and $\{\b_{k}\}$ be sequences of positive numbers such that $\e_{k}, \g_{k}, \b_{k} \to 0.$ Consider  sequences of varifolds $V_{k} \in {\mathcal S}_{\a}$ and cones ${\mathbf C}_{k} \in {\mathcal C}_{q}$ such that, for each $k=1, 2, \ldots,$ with $V_{k}$, ${\mathbf C}_{k}$  in place of $V$, ${\mathbf C}$ respectively, Hypotheses~\ref{hyp} hold with $\e_{k}$, $\g_{k}$ in place of $\e$, $\g$; Hypothesis ($\star$) holds with $M = \frac{3}{2}M_{0}^{3}$ and Hypothesis~($\star\star$) holds with $\b_{k}$ in place of $\b$. Thus, for each $k=1, 2, \ldots,$ we suppose: 
\begin{itemize}
\item[$(1_{k})$]  $\Theta \, (\|V_{k}\|, 0) \geq q$, \; $(2\omega_{n})^{-1}\|V_{k}\|(B_{2}^{n+1}(0)) < q + 1/2,$ \; $\omega_{n}^{-1}\|V_{k}\|({\mathbf R} \times B_{1}) < q + 1/2.$
\item[$(2_{k})$] ${\mathbf C}_{k} = \sum_{j=1}^{q} |H_{j}^{k}| + |G_{j}^{k}|$ where for each $j \in \{1, 2, \ldots, q\}$, $H_{j}^{k}$ is the half-space defined by
$H_{j}^{k} = \{(x^{1}, x^{2}, y) \in {\mathbf R}^{n+1} \, : \, x^{2} < 0 \;\; \mbox{and} \;\; x^{1} = \lambda_{j}^{k}x^{2}\},$ 
$G_{j}^{k}$ the half-space defined by $G_{j}^{k} = \{(x^{1},x^{2}, y) \in {\mathbf R}^{n+1}\, : \, x^{2} > 0 \;\; \mbox{and} \;\; x^{1} = \mu_{j}^{k}x^{2}\},$ with $\lambda_{j}^{k}, \mu_{j}^{k}$ constants,  
$\lambda_{1}^{k} \geq \lambda_{2}^{k} \geq \ldots \geq \lambda_{q}^{k}$ and $\mu_{1}^{k} \leq \mu_{2}^{k} \leq \ldots \leq \mu_{q}^{k}.$
\item[$(3_{k})$] ${\hat E}_{k}^{2} = {\hat E}_{V_{k}}^{2} \equiv \int_{{\mathbf R} \times B_{1}} |x^{1}|^{2} d\|V_{k}\|(X) < \e_{k}.$
\item[$(4_{k})$] $\{Z \, : \ \Theta \, (\|V_{k}\|, Z) \geq q\} \cap \left({\mathbf R} \times (B_{1/2} \setminus \{|x^{2}| < 1/16\})\right) = \emptyset$.
\item[$(5_{k})$] ${\hat E}_{k}^{-2}\left(Q_{k}({\mathbf C}_{k})\right)^{2} < \g_{k}$ where
\begin{eqnarray*}
\left(Q_{k}({\mathbf C}_{k})\right)^{2} = \left(Q_{V_{k}}({\mathbf C}_{k})\right)^{2} = \left(\int_{{\mathbf R} \times (B_{1/2} \setminus \{|x^{2}| < 1/16\})} {\rm dist}^{2}(X, {\rm spt} \, \|V_{k}\|) \,d\|{\mathbf C}_{k}\|(X)\right.&&\nonumber\\
&&\hspace{-2in} +\; \left.\int_{{\mathbf R} \times B_{1}} {\rm dist}^{2} \, (X, {\rm spt}\, \|{\mathbf C}_{k}\|) \, d\|V_{k}\|(X)\right).
\end{eqnarray*}
\item[$(6_{k})$] $${\hat E}^{2}_{k} < \frac{3}{2}M^{3}_{0}\inf_{\{P = \{x^{1} = \lambda x^{2}\}\}} \,\int_{{\mathbf R} \times B_{1}} {\rm dist}^{2} \, (X, P) \, d\|V_{k}\|(X).$$  
\item[$(7_{k})$] Either (i) or (ii) below holds: 
\begin{itemize}
\item[(i)] ${\mathbf C}_{k} \in {\mathcal C}_{q}(4).$ 
\item[(ii)] $q \geq 3,$ ${\mathbf C}_{k} \in {\mathcal C}_{q}(p_{k})$ for some $p_{k} \in \{5, 6, \ldots, 2q\}$ and $\left(Q_{k}^{\star}\right)^{-2}\left(Q_{k}({\mathbf C}_{k})\right)^{2} < \b_{k}$ where
\begin{eqnarray*}
\left(Q_{k}^{\star}\right)^{2} = \left(Q_{V_{k}}^{\star}(p_{k}-1)\right)^{2} = \inf_{\widetilde{\mathbf C}\in \cup_{j=4}^{p_{k}-1}{\mathcal C}_{q}(j)} \, \left(\int_{{\mathbf R} \times (B_{1/2} \setminus \{|x^{2}| < 1/16\})} {\rm dist}^{2}(X, {\rm spt} \, \|V_{k}\|) \,d\|\widetilde{\mathbf C}\|(X)\right.&&\\
&&\hspace{-2.5in} + \left.\int_{{\mathbf R} \times B_{1}} {\rm dist}^{2} \, (X, {\rm spt}\, \|\widetilde{\mathbf C}\|) \, d\|V_{k}\|(X)\right).
\end{eqnarray*}
\end{itemize}
\end{itemize}
Let $E_{k} = \sqrt{\int_{{\mathbf R} \times B_{1}} {\rm dist}^{2} \, (X, {\rm spt} \, \|{\mathbf C}_{k}\|) \,d\|V_{k}\|(X)},$ so that by $(5_{k}),$ 
\begin{equation}\label{excess-0}
{\hat E}_{k}^{-1}E_{k} \to 0.
\end{equation}

Note also that in case ${\mathbf C}_{k} \not\in {\mathcal C}_{q}(4)$ except for finitely many $k$, 
we have by $(3_{k})$ and $(5_{k})$ that $Q_{k}^{\star} \to 0.$

Let $\{\d_{k}\}, \{\t_{k}\}$ be sequences of decreasing positive numbers converging to 0. By passing to appropriate 
subsequences of $\{V_{k}\}$, $\{{\mathbf C}_{k}\}$, and possibly replacing ${\mathbf C}_{k}$ with a cone 
${\mathbf C}_{k}^{\prime} \in {\mathcal C}_{q}$ with ${\rm spt} \, \|{\mathbf C}_{k}^{\prime}\| = {\rm spt} \, \|{\mathbf C}_{k}\|$ without changing notation (see Remark (2) following the statement of Hypothesis~($\star\star$)), we deduce that, for each $k=1, 2, \ldots,$ the assertions $(A_{k})$-$(D_{k})$ below hold:
 \begin{itemize}
\item[(A$_{k}$)] By Lemma~\ref{no-gaps},
\begin{equation}\label{excess-1}
B^{n+1}_{\d_{k}}(0,y) \cap \{Z \, : \, \Theta \, (\|V_{k}\|, Z) \geq q\} \neq \emptyset
\end{equation}
for each point $(0, y) \in \{0\} \times {\mathbf R}^{n-1} \cap B_{1/2}$ and
\begin{equation}\label{excess-2}
\int_{B_{1/2}^{n+1}(0) \cap \{|(x^{1}, x^{2})| < \s\}} {\rm dist}^{2} \, (X, {\rm spt} \, \|{\mathbf C}_{k}\|) \, d\|V_{k}\|(X) 
\leq C\s^{1/2} E_{k}^{2}
\end{equation}
for each $\s \in [\d_{k}, 1/4),$ where $C= C(n, q, \a) \in (0, \infty).$
\item[(B$_{k}$)] By Theorem~\ref{L2-est-1}(a), 
\begin{equation}\label{excess-3}
V_{k} \res ({\mathbf R} \times (B_{3/4} \setminus \{|x^{2}| < \t_{k}\})) = \sum_{j=1}^{q} |{\rm graph} \, (h_{j}^{k} + u_{j}^{k})| + |{\rm graph} \, (g_{j}^{k} + w_{j}^{k})|
\end{equation}
where $h_{j}^{k}, g_{j}^{k}$ are the linear functions on ${\mathbf R}^{n}$ given by 
$h_{j}^{k}(x^{2},y) = \lambda_{j}^{k}x^{2}$, $g_{j}^{k}(x^{2}, y) = \mu_{j}^{k}x^{2},$ 
$u_{j}^{k} \in C^{2} \, (B_{3/4} \,\cap \,\{x^{2} < -\t_{k}\}),$ $w_{j}^{k} \in C^{2} \, (B_{3/4} \, \cap \, \{x^{2} > \t_{k}\})$  
with $h_{j}^{k} + u_{j}^{k}$ and $g_{j}^{k} + w_{j}^{k}$ solving the minimal surface equation on their respective domains and satisfying $h_{1}^{k} + u_{1}^{k} \leq h_{2}^{k} + u_{2}^{k} \leq \ldots \leq h_{q}^{k} + u_{q}^{k},$  
$g_{1}^{k} + w_{1}^{k} \leq g_{2}^{k} + w_{2}^{k} \leq \ldots \leq g_{q}^{k} + w_{q}^{k},$
${\rm dist} \, ((h_{j}^{k}(x^{2}, y) + u_{j}^{k}(x^{2},y), x^{2}, y), {\rm spt} \,\|{\mathbf C}_{k}\|) = 
(1 + (\lambda_{j}^{k})^{2})^{-1/2}|u_{j}^{k}(x^{2}, y)|$ for $(x^{2}, y) \in B_{3/4} \cap \{x^{2} < -\t_{k}\}$ and 
${\rm dist} \, ((g_{j}^{k}(x^{2}, y) + w_{j}^{k}(x^{2},y), x^{2}, y), {\rm spt} \,\|{\mathbf C}_{k}\|) = 
(1 + (\mu_{j}^{k})^{2})^{-1/2}|w_{j}^{k}(x^{2}, y)|$ for $(x^{2}, y) \in B_{3/4} \cap \{x^{2} > \t_{k}\}.$ 
\item[(C$_{k}$)] For each  point $Z = (\z^{1}, \z^{2}, \eta) \in {\rm spt} \, \|V_{k}\| \cap ({\mathbf R} \times B_{3/8})$ with $\Theta \, (\|V_{k}\|, Z) \geq q,$ by Corollary~\ref{L2-est-2}(a) (taken with $\r = 1/4$, say), 
\begin{equation}\label{excess-4}
|\z^{1}|^{2} + {\hat E}_{k}^{2}|\z^{2}|^{2} \leq C E_{k}^{2}
\end{equation}
where $C = C(n, q, \a) \in (0, \infty).$ 
\item[(D$_{k}$)] By (\ref{slope-bounds-1}) and (\ref{slope-bounds-2}), 
\begin{eqnarray}\label{excess-5-1}
c {\hat E}_{k} \leq {\rm max} \, \{|\lambda_{1}^{k}|, |\lambda_{q}^{k}|\} \leq c_{1} {\hat E}_{k}, \;\;\; c{\hat E}_{k}  \leq {\rm max} \, \{|\mu_{1}^{k}|, |\mu_{q}^{k}|\} \leq c_{1} {\hat E}_{k} \;\;\; {\rm and}\nonumber\\
&&\hspace{-3.5in} {\rm min} \, \{|\lambda^{k}_{1} - \lambda^{k}_{q}|, |\mu^{k}_{1} - \mu^{k}_{q}|\} \geq 2c {\hat E}_{k}
\end{eqnarray}
where $c_{1} = c_{1}(n) \in (0, \infty)$ and $c = c(n, q) \in (0, \infty)$. 
\end{itemize}
Furthermore, by Corollary~\ref{L2-est-2}(b), (\ref{excess-3}) and the area formula,
there exists, for each $\r \in (0, 1/4]$, an integer $K = K(\r) \geq 1$ such that the following assertion holds for each $k \geq K:$
\begin{itemize}
\item[(E$_{k}$)]  For each point $Z = (\z^{1}, \z^{2}, \eta) \in {\rm spt} \, \|V_{k}\| \cap ({\mathbf R} \times B_{3/8})$ with $\Theta \, (\|V_{k}\|, Z) \geq q$ and each $\m \in (0, 1)$, 
\begin{eqnarray}\label{excess-5}
&&\hspace{.5in}\sum_{j=1}^{q} \int_{B_{\r/2}(\z^{2}, \eta) \cap \{x^{2} < -\t_{k}\}} \frac{|u_{j}^{k}(x^{2}, y) - (\z^{1} - \lambda_{j}^{k}\z^{2})|^{2}}{|(h_{j}^{k}(x^{2},y)+u_{j}^{k}(x^{2}, y), x^{2}, y) - (\z^{1}, \z^{2}, \eta)|^{n+2-\m}} \, dx^{2} \, dy\nonumber\\ 
&&\hspace{.75in}+ \;\;\sum_{j=1}^{q}\int_{B_{\r/2}(\z^{2}, \eta) \cap \{x^{2} > \t_{k}\}}\frac{|w_{j}^{k}(x^{2}, y) - (\z^{1} - \mu_{j}^{k}\z^{2})|^{2}}{|(g_{j}^{k}(x^{2},y)+w_{j}(x^{2}, y), x^{2}, y) - (\z^{1}, \z^{2}, \eta)|^{n+2-\m}} \, dx^{2}dy\nonumber\\
&&\hspace{2.5in} \leq C_{1} \r^{-n-2+ \m}\int_{{\mathbf R} \times B_{\r}(\z^{2}, \eta)} {\rm dist}^{2} \, (X, {\rm spt} \, \|T_{Z \, \#} \, {\mathbf C}_{k}\|) \, d\|V_{k}\|(X)
\end{eqnarray}
where $C_{1}  = C_{1}(n, q, \a, \m) \in (0, \infty).$
\end{itemize}

Extend $u_{j}^{k},$ $w_{j}^{k}$ to all of $B_{3/4} \cap \{x^{2} <0\}$ and $B_{3/4} \cap \{x^{2} > 0\}$ respectively by defining values to be zero in 
$B_{3/4} \cap \{0 > x^{2} \geq -\t_{k}\}$ and $B_{3/4} \cap \{0 < x^{2} \leq \t_{k}\}$ respectively.

By (\ref{excess-5-1}), there exist numbers $\ell_{j}, m_{j}$ for each $j=1, 2, \ldots, q$ with 
\begin{eqnarray}\label{excess-6}
c \leq {\rm max} \, \{|\ell_{1}|, |\ell_{q}|\} \leq c_{1}, \;\;\; c \leq {\rm max} \,\{|m_{1}|, |m_{q}|\} 
\leq c_{1} \;\;\; \mbox{and}&&\nonumber\\
&&\hspace{-3in}{\rm min} \, \{|\ell_{1} - \ell_{q}|, |m_{1} - m_{q}|\} \geq 2c 
\end{eqnarray}
such that after passing to appropriate subsequences without changing notation, 
\begin{equation}\label{excess-7}
{\hat E}_{k}^{-1} \lambda_{j}^{k} \to \ell_{j} \;\;\; \mbox{and} \;\;\; {\hat E}_{k}^{-1}\mu_{j}^{k} \to  m_{j}
\end{equation}
for each $j=1, 2, \ldots, q.$ By (\ref{excess-3}) and elliptic estimates, there exist harmonic functions 
$\varphi_{j} \, : \, B_{3/4} \cap \{x^{2} < 0\} \to {\mathbf R}$ and $\psi_{j} \, : \, B_{3/4} \cap \{x^{2} > 0\} \to {\mathbf R}$ such that
\begin{equation}\label{excess-8}
E_{k}^{-1}u_{j}^{k} \to \varphi_{j} \;\;\; \mbox{and} \;\;\; E_{k}^{-1}w_{j}^{k} \to \psi_{j}
\end{equation}
where the convergence is in $C^{2}(K)$ for each compact subset $K$ of the respective domains. From (\ref{excess-2}), it follows that 
$$\int_{B_{1/2} \cap \{0 > x^{2} >- \s\}}|\varphi|^{2} \leq  C\s^{1/2}, \;\;\;\; 
\int_{B_{1/2} \cap \{0 < x^{2} < \s\}}|\psi|^{2} \leq C\s^{1/2}$$
for each $\s \in (0, 1/4)$,  where $C = C(n, q, \a) \in (0, \infty)$, and hence that the convergence in (\ref{excess-8}) is, respectively, also in $L^{2}\,(B_{1/2} \cap \{x^{2} < 0\})$ and 
$L^{2}\,(B_{1/2} \cap \{x^{2} > 0\}).$

Set $\varphi = (\varphi_{1}, \varphi_{2}, \ldots, \varphi_{q})$ and $\psi = (\psi_{1}, \psi_{2}, \ldots, \psi_{q}).$

\noindent
{\bf Definitions: (1) Fine blow-ups.}  Let 
$\varphi \, : \, B_{3/4} \cap \{x^{2} < 0\}  \to {\mathbf R}^{q}$ and $\psi \, : \, B_{3/4} \cap \{x^{2} > 0\} \to {\mathbf R}^{q}$ be a pair of functions arising, in the manner described above, corresponding to: (i) a  sequence 
of varifolds $\{V_{k}\} \subset {\mathcal S}_{\a}$ and a sequence of cones $\{{\mathbf C}_{k}\} \subset {\mathcal C}_{q}$ 
satisfying the hypotheses $(1_{k})-(7_{k})$ for some sequences of numbers $\{\e_{k}\},$ $\{\g_{k}\}$, $\{\b_{k}\}$ with 
$\e_{k}, \g_{k}, \b_{k} \to 0^{+}$ and (ii) sequences $\{\d_{k}\}$, $\{\t_{k}\}$ of decreasing positive numbers converging to zero such that (\ref{excess-1}), (\ref{excess-2}), (\ref{excess-3}) and (\ref{excess-5}) hold. We call 
the pair $(\varphi, \psi)$ a  {\em fine blow-up} of the sequence $\{V_{k}\}$ relative to $\{{\mathbf C}_{k}\}$.

\noindent
{\bf (2) The Class ${\mathcal B}^{F}$.} Let ${\mathcal B}^{F}$ be the collection of all fine blow-ups $(\varphi, \psi)$ such that the corresponding sequences of varifolds $V_{k} \in {\mathcal S}_{\a}$ satisfies condition $(6_{k})$ with $M_{0}^{2}$ in place of $M_{0}^{3};$ thus we assume the stronger condition 
$${\hat E}_{V_{k}}^{2} < \frac{3}{2}M^{2}_{0}\inf_{\{P = \{x^{1} = \lambda x^{2}\}\}}\int_{{\mathbf R} \times B_{1}} {\rm dist}^{2} \, (X, P) \, d\|V_{k}\|(X), \;\;\; k=1, 2, 3, \ldots$$ in place of $(6_{k})$ for any sequence $\{V_{k}\} \subset {\mathcal S}_{\a}$ giving rise to a fine blow-up belonging to ${\mathcal B}^{F}$.

\section{Continuity estimates for the fine blow-ups and their derivatives}\label{c0alpha}
\setcounter{equation}{0}

Here we first use estimates (\ref{excess-4}) and (\ref{excess-5}) to prove a continuity estimate (Lemma~\ref{continuity} below) for any $(\varphi, \psi)  \in {\mathcal B}^{F}.$ We then use it to establish the main result of this section (Theorem~\ref{c1alpha}), namely, the continuity estimate for the first derivatives of $(\varphi, \psi) \in {\mathcal B}^{F}.$
  
\begin{lemma}\label{continuity}
If $(\varphi, \psi) \in {\mathcal B}^{F}$, then
$$\varphi \in C^{0, \b} \, (\overline{B_{5/16} \cap \{x^{2} < 0\}}; {\mathbf R}^{q}), \;\;\; 
\psi \in C^{0, \b} \, (\overline{B_{5/16} \cap \{x^{2} > 0\}}; {\mathbf R}^{q})$$ for some $\b = \b(n, q, \a) \in (0,1)$ 
and the following estimates hold:
\begin{eqnarray*}
\sup_{\overline{B_{5/16} \cap \{x^{2} < 0\}}} \, |\varphi|^{2} \; + \; \sup_{x, z \in \overline{B_{5/16} \cap  \,\{x^{2} < 0\}}, x \neq z} \, \frac{|\varphi(x) - \varphi(z)|^{2}}{|x - z|^{2\b}}&&\nonumber\\
&&\hspace{-2in} \leq \; 
C\left(\int_{B_{1/2} \cap \{x^{2} <0\}} |\varphi|^{2} +  \int_{B_{1/2} \cap \{x^{2} > 0\}} |\psi|^{2}\right); 
\end{eqnarray*}
\begin{eqnarray*}
\sup_{\overline{B_{5/16} \cap \{x^{2} > 0\}}} \, |\psi|^{2} \;+  \sup_{x, z\in \overline{B_{5/16} \cap  \,\{x^{2} > 0\}}, x \neq z} \, \frac{|\psi(x) - \psi(z)|^{2}}{|x - z|^{2\b}}&&\nonumber\\ 
&&\hspace{-2in}\leq 
C\left(\int_{B_{1/2} \cap \{x^{2} <0\}} |\varphi|^{2} +  \int_{B_{1/2} \cap \{x^{2} > 0\}} |\psi|^{2}\right).
\end{eqnarray*}
Here $C = C(n,q, \a) \in (0, \infty).$
\end{lemma}

\begin{proof} By the definition of fine blow-up, there are sequences $\{V_{k}\} \subset {\mathcal S}_{\a}$, $\{{\mathbf C}_{k}\} \subset {\mathcal C}_{q}$ and sequences of decreasing positive numbers $\{\e_{k}\}$, $\{\g_{k}\}$, $\{\b_{k}\}$, $\{\d_{k}\}$, $\{\t_{k}\}$ converging to zero for which all of the assertions of Section~\ref{fineexcessblowup} hold, with $M^{2}_{0}$ in place of $M_{0}^{3}$ in $(6_{k})$. 

Let $Y \in \{0\} \times {\mathbf R}^{n-1} \cap B_{5/16}$ be arbitrary. By (\ref{excess-1}), for each $k=1, 2, 3, \ldots,$ there exist
$Z_{k} = (\z_{1}^{k}, \z_{2}^{k}, \eta^{k}) \in {\rm spt} \, \|V_{k}\|$ with $\Theta \, (\|V_{k}\|, Z_{k}) \geq q$ such that $Z_{k} \to Y.$ Using (\ref{excess-2}), (\ref{excess-4}), (\ref{excess-5}) (with $\z_{1}^{k}, \z_{2}^{k}, \eta^{k}$ in place of $\z_{1}, \z_{2}, \eta$ and  $\m = 1/2$) and (\ref{excess-7}), we deduce that for each $\r \in (0, 1/8]$, 
\begin{eqnarray}\label{continuity-1}
\sum_{j=1}^{q} \int_{B_{\r/2}(Y) \cap \{x^{2} <0\}} \frac{|\varphi_{j}(x) - (\k_{1}(Y) - \ell_{j}\k_{2}(Y))|^{2}}
{|x- Y|^{n+3/2}} \, dx &&\nonumber\\
&&\hspace{-3in}+ \;\;\sum_{j=1}^{q}\int_{B_{\r/2}(Y) \cap \{x^{2}> 0\}}\frac{|\psi_{j}(x) - (\k_{1}(Y)- 
m_{j}\k_{2}(Y))|^{2}}{|x - Y|^{n+3/2}} \, dx\nonumber\\
&&\hspace{-2in} \leq C_{1} \r^{-n-3/2}\sum_{j=1}^{q}\int_{B_{\r}(Y) \cap \{x^{2} <0\}} |\varphi_{j} - (\k_{1}(Y) - \ell_{j}\k_{2}(Y))|^{2}\nonumber\\
&&\hspace{-1.5in} +\, C_{1}\r^{-n-3/2}\sum_{j=1}^{q} \int_{B_{\r}(Y) \cap 
\{x^{2} > 0\}} |\psi_{j} - (\k_{1}(Y) - m_{j}\k_{2}(Y))|^{2} 
\end{eqnarray}
where $C_{1} = C_{1}(n, q, \a) \in (0, \infty)$ and we have set 
\begin{equation}\label{continuity-1-0}
\k_{1}(Y) = \lim_{k \to \infty} \, {E}_{k}^{-1}\z_{1}^{k}, \;\; \k_{2}(Y) = \lim_{k \to \infty} \, E_{k}^{-1} {\hat E}_{k} \z_{2}^{k},
\end{equation}
both of which limits exist after passing to 
a subsequence of the original sequence $\{k\}$. Note that by (\ref{excess-4}), 
\begin{equation}\label{continuity-1-1}
|\k_{1}(Y)|, |\k_{2}(Y)| \leq C, \;\; C = C(n, q,\a) \in (0, \infty).
\end{equation}

We remark also that our notation here is appropriate, and the limits in (\ref{continuity-1-0}) indeed depend only on $Y$ and are independent of the sequence of points $Z_{k}$ converging to $Y;$ this follows directly from the finiteness of the integrals on the left hand side of (\ref{continuity-1}) and the fact that, by Lemma~\ref{no-overlap},  at least two of the $\ell_{j}$'s and two of the $m_{j}$'s are distinct.

For $Y \in \{0\} \times {\mathbf R}^{n-1} \cap B_{5/16}$ and each $j=1, 2, \ldots, q$, 
define 
\begin{equation}\label{continuity-1-1-1}
\varphi_{j}(Y) = \k_{1}(Y)  - \ell_{j}\k_{2}(Y) \;\;\; {\rm and} \;\;\; \psi_{j}(Y) = \k_{1}(Y) - m_{j}\k_{2}(Y).
\end{equation}
Then by (\ref{continuity-1}), 
\begin{eqnarray}\label{continuity-2-1}
\s^{-n}\left(\int_{B_{\s}(Y) \cap \{x^{2} <0\}} |\varphi(x) - \varphi(Y)|^{2}\, dx \; + 
\; \int_{B_{\s}(Y) \cap \{x^{2}> 0\}}|\psi(x) - \psi(Y)|^{2}\, dx \right)&&\nonumber\\
&&\hspace{-4.75in}\leq C_{1} \left(\frac{\s}{\r}\right)^{3/2}\r^{-n}\left(\int_{B_{\r}(Y) \cap \{x^{2} <0\}} |\varphi(x) - \varphi(Y)|^{2} \, dx 
+  \int_{B_{\r}(Y) \cap \{x^{2} > 0\}} |\psi(x) - \psi(Y)|^{2} \, dx \right)
\end{eqnarray}
for each $0 < \s \leq \r/2 \leq 1/32$ and for the same constant $C_{1} = C_{1}(n, q, \a) \in (0, \infty)$ as in (\ref{continuity-1}).

To complete the proof of the lemma, we follow the argument of Lemma~\ref{general}. Consider an arbitrary point $z^{+} \in B_{5/16} \cap \{x^{2} > 0\}$ and let
$\r \in (0, 1/16].$  Denote by $z^{-}$ the image of $z^{+}$ under reflection across $\{0\} \times {\mathbf R}^{n-1}.$ 
Letting $Y \in \{0\} \times {\mathbf R}^{n-1}$ 
be the point such that $|z^{-} - Y| = |z^{+} - Y| = {\rm dist} \, (z^{+}, \{0\} \times {\mathbf R}^{n-1})$, and with $\g  =\g(n, q, \a) \in (0, 1/16]$ to be chosen, if ${\rm dist} \, (z^{+}, \{0\} \times {\mathbf R}^{n-1}) < \g\r$, then
\begin{eqnarray}
(\g\r)^{-n}\left(\int_{B_{\g\r}(z^{-}) \cap \{x^{2} < 0\}} |\varphi - \varphi(Y)|^{2} + \int_{B_{\g\r}(z^{+}) \cap \{x^{2} > 0\}} |\psi - \psi(Y)|^{2}\right)&&\nonumber\\
&&\hspace{-4.75in}\leq2^{n} (\g\r + |z^{-}-Y|)^{-n} \left(\int_{B_{\g\r +|z^{-}-Y|}(Y) \cap \{x^{2} < 0\}} |\varphi - \varphi(Y)|^{2}
+ \int_{B_{\g\r + |z^{+}-Y|}(Y) \cap \{x^{2} > 0\}} |\psi - \psi(Y)|^{2}\right)\nonumber\\ 
&&\hspace{-4in}\leq 2^{n}C_{1}\left(\frac{\g\r + |z^{-}-Y|}{\r - |z^{-}-Y|}\right)^{3/2} (\r - |z^{-} - Y|)^{-n}\left(\int_{B_{\r - |z^{-}-Y|}(Y) \cap \{x^{2} < 0\}} |\varphi - \varphi(Y)|^{2}\right.\nonumber\\
&&\hspace{-1in} + \left.\int_{B_{\r - |z^{+}-Y|}(Y) \cap \{x^{2} > 0\}} |\psi - \psi(Y)|^{2}\right)\nonumber\\ 
&&\hspace{-4in}\leq4^{n}C_{1} \left(\frac{2\g}{1 - \g}\right)^{3/2} \r^{-n}\left(\int_{B_{\r}(z^{-}) \cap \{x^{2} < 0\}} |\varphi - \varphi(Y)|^{2} + \int_{B_{\r}(z^{+}) \cap \{x^{2} > 0\}} |\psi - \psi(Y)|^{2}\right).
\end{eqnarray}
Choosing $\g = \g(n, q, \a) \in (0, 1/16]$ such that $4^{n}C_{1} \left(\frac{2\g}{1 - \g}\right)^{3/2} < 1/4$, 
we deduce that 
\begin{eqnarray}\label{continuity-3}
(\g\r)^{-n}\left(\int_{B_{\g\r}(z^{-}) \cap \{x^{2} < 0\}} |\varphi - \varphi(Y)|^{2} + \int_{B_{\g\r}(z^{+}) \cap \{x^{2} > 0\}} |\psi - \psi(Y)|^{2}\right)&&\nonumber\\
&&\hspace{-4.5in}\leq4^{-1}\r^{-n}\left(\int_{B_{\r}(z^{-}) \cap \{x^{2} < 0\}} |\varphi - \varphi(Y)|^{2} + \int_{B_{\r}(z^{+}) \cap \{x^{2} > 0\}} |\psi - \psi(Y)|^{2}\right)
\end{eqnarray} 
for any $z^{+} \in B_{5/16} \cap \{x^{2} >0\}$ and $\r \in (0, 1/16]$ provided $\g\r > |z^{+} - Y| = |z^{-} - Y| = {\rm dist} \, (z^{+}, \{0\} \times {\mathbf R}^{n-1}).$ If on the other hand $\g\r \leq {\rm dist} \, (z^{+}, \{0\} \times {\mathbf R}^{n-1})$, since $\varphi$ and $\psi$ are harmonic in $B_{1/2} \cap \{x^{2} < 0\}$ and $B_{1/2} \cap \{x^{2} > 0\}$ respectively, we have for each $\s \in (0, 1/2]$ and any constant vectors $b^{+}, b^{-} \in {\mathbf R}^{q}$, 
\begin{eqnarray}\label{continuity-4}
(\s\g\r)^{-n}\left(\int_{B_{\s\g\r}(z^{-})}|\varphi - \varphi(z^{-})|^{2} + \int_{B_{\s\g\r}(z^{+})} |\psi - \psi(z^{+})|^{2}\right)&&\nonumber\\
&&\hspace{-3in} \leq  C\s^{2}(\g\r)^{-n}
\left(\int_{B_{\g\r}(z^{-})}|\varphi - b^{-}|^{2} + \int_{B_{\g\r}(z^{+})} |\psi - b^{+}|^{2}\right) 
\end{eqnarray}
where $C = C(n) \in (0, \infty).$

Given any $z^{+} \in B_{5/16} \cap \{x^{2} >0\}$, let $j_{\star} \in \{0, 1, 2, \ldots \}$ be such that 
$\g^{j_{\star} +1} < {\rm dist} \, (z^{+}, \{0\} \times {\mathbf R}^{n-1}) \leq \g^{j_{\star}}.$ Then, 
with $Y \in \{0\} \times {\mathbf R}^{n-1}$ such that $|z^{+} - Y| = {\rm dist} \, (z^{+}, \{0\} \times {\mathbf R}^{n-1})$,  
by (\ref{continuity-4}),
\begin{eqnarray}\label{continuity-5}
(\s\g^{j_{\star} +1})^{-n}\left(\int_{B_{\s\g^{j_{\star} +1}}(z^{-})}|\varphi - \varphi(z^{-})|^{2} + \int_{B_{\s\g^{j_{\star} + 1}}(z^{+})} |\psi - \psi(z^{+})|^{2}\right)&&\nonumber\\
&&\hspace{-4in} \leq  C\s^{2}(\g^{j_{\star} + 1})^{-n}
\left(\int_{B_{\g^{j_{\star} + 1}}(z^{-})}|\varphi - \varphi(Y)|^{2} + \int_{B_{\g^{j_{\star} + 1}}(z^{+})} |\psi - \psi(Y)|^{2}\right) 
\end{eqnarray}
for any $\s \in (0, 1/2],$ and  if $j_{\star} \geq 1$, by (\ref{continuity-3}), 
\begin{eqnarray}\label{continuity-6}
(\g^{j})^{-n}\left(\int_{B_{\g^{j}}(z^{-}) \cap \{x^{2} < 0\}} |\varphi - \varphi(Y)|^{2} + \int_{B_{\g^{j}}(z^{+}) \cap \{x^{2} > 0\}} |\psi - \psi(Y)|^{2}\right)&&\nonumber\\
&&\hspace{-4.5in}\leq4^{-1}(\g^{j-1})^{-n}\left(\int_{B_{\g^{j-1}}(z^{-}) \cap \{x^{2} < 0\}} |\varphi - \varphi(Y)|^{2} + \int_{B_{\g^{j-1}}(z^{+}) \cap \{x^{2} > 0\}} |\psi - \psi(Y)|^{2}\right)\nonumber\\
&&\hspace{-4.5in}\leq 4^{-(j-1)}\g^{-n}\left(\int_{B_{\g}(z^{-}) \cap \{x^{2} < 0\}} |\varphi - \varphi(Y)|^{2} + \int_{B_{\g}(z^{+}) \cap \{x^{2} > 0\}} |\psi - \psi(Y)|^{2}\right)
\end{eqnarray}  
for $j=1, 2, \ldots, j_{\star}.$  If $j_{\star} \geq 1$, taking $j = j_{\star}$ in (\ref{continuity-6}) and $\s =1/2$ in (\ref{continuity-5}), we see by the triangle inequality that 
\begin{eqnarray*}
|\varphi(z^{-}) - \varphi(Y)|^{2} + |\psi(z^{+}) - \psi(Y)|^{2}&&\\
&&\hspace{-2.5in} \leq C4^{-(j_{\star}-1)}\g^{-n}\left(\int_{B_{\g}(z^{-}) \cap \{x^{2} < 0\}} |\varphi - \varphi(Y)|^{2} + \int_{B_{\g}(z^{+}) \cap \{x^{2} > 0\}} |\psi - \psi(Y)|^{2}\right)
\end{eqnarray*}
where $C = C(n, q, \a) \in (0, \infty),$ and hence by (\ref{continuity-6}) and the triangle inequality again that 
\begin{eqnarray}\label{continuity-7}
(\g^{j})^{-n}\left(\int_{B_{\g^{j}}(z^{-}) \cap \{x^{2} < 0\}} |\varphi - \varphi(z^{-})|^{2} + \int_{B_{\g^{j}}(z^{+}) \cap \{x^{2} > 0\}} |\psi - \psi(z^{+})|^{2}\right)&&\nonumber\\
&&\hspace{-4.5in}\leq C4^{-(j-1)}\g^{-n}\left(\int_{B_{\g}(z^{-}) \cap \{x^{2} < 0\}} |\varphi - \varphi(Y)|^{2} + \int_{B_{\g}(z^{+}) \cap \{x^{2} > 0\}} |\psi - \psi(Y)|^{2}\right)
\end{eqnarray}  
for $j=1, 2, \ldots, j_{\star}$, where $C = C(n, q, \a) \in (0, \infty).$ 

By applying (\ref{excess-4}) with 
$\widetilde{V}_{k} \equiv \eta_{0, 1/2 \, \#} \, V_{k}$ in place of $V_{k},$ and noting (e.g.\ by the argument establishing (\ref{L2-est-2-1-0})) that 
${\hat E}_{{\widetilde V}_{k}} \geq C{\hat E}_{V_{k}}$ where $C = C(n, q) \in (0, \infty)$, we deduce using also (\ref{excess-2}) and (\ref{excess-6}) that for each $Y \in \{0\} \times {\mathbf R}^{n-1} \cap B_{5/16}$ that 
\begin{equation}\label{continuity-7'}
|\varphi(Y)|^{2} + |\psi(Y)|^{2} \leq C\left(\int_{B_{1/2} \cap \{x^{2} < 0\}}|\varphi|^{2} + \int_{B_{1/2} \cap \{x^{2} > 0\}}|\psi|^{2}\right)
\end{equation}
where $C = C(n, q, \a) \in (0, \infty).$ With the help of (\ref{continuity-5}), (\ref{continuity-6}), (\ref{continuity-7}) and (\ref{continuity-7'}), we deduce that for any given $z^{+} \in B_{5/16} \cap \{x^{2} > 0\},$
\begin{eqnarray}\label{continuity-8}
\r^{-n}\left(\int_{B_{\r}(z^{-}) \cap \{x^{2} < 0\}} |\varphi - \varphi(z^{-})|^{2} + \int_{B_{\r}(z^{+}) \cap \{x^{2} > 0\}} |\psi - \psi(z^{+})|^{2}\right)&&\nonumber\\
&&\hspace{-2.5in}\leq C\r^{2\b}\left(\int_{B_{1/2} \cap \{x^{2} < 0\}} |\varphi|^{2} + \int_{B_{1/2}\cap \{x^{2} > 0\}} |\psi|^{2}\right)
\end{eqnarray}  
for all $\r \in (0, \g],$ where $C = C(n, q, \a) \in (0, \infty)$ and $\b = \b(n, q, \a) \in (0, 1),$ by considering, for any given 
$\r \in (0, \g]$, the alternatives $2\r \leq \g^{j_{\star} +1}$, in which case $\r = \s\g^{j_{\star} +1}$ for some $\s \in (0, 1/2]$ 
and we use (\ref{continuity-5}) and (\ref{continuity-6}) with $j=j_{\star}$, or $\g^{j + 1} < 2\r \leq \g^{j}$ for some $j \in \{1, 2, \ldots, j_{\star}\}$, in which case we use (\ref{continuity-7}). The conclusions of the lemma follow readily from (\ref{continuity-8}). 
\end{proof}

\begin{theorem}\label{c1alpha}
If $(\varphi, \psi) \in {\mathcal B}^{F}$, then 
$$\varphi \in C^{2} \, (\overline{B_{1/4} \cap \{x^{2} < 0\}}; {\mathbf R}^{q}), \;\;\; \psi \in C^{2} \, (\overline{B_{1/4} \cap \{x^{2} >0\}}; {\mathbf R}^{q})$$ and the following estimates hold:
\begin{eqnarray*}
\sup_{\overline{B_{1/4} \cap \{x^{2} < 0\}}} \, |D\varphi|^{2} + \sup_{x, z \in \overline{B_{1/4} \cap \{x^{2} < 0\}}, \, x\neq z} \, 
\frac{|D\varphi(x) - D\varphi(z)|^{2}}{|x-z|^{2}}&&\\
&&\hspace{-2.2in} \leq C\left(\int_{B_{1/2} \cap \{x^{2} < 0\}}|\varphi|^{2} + \int_{B_{1/2} \cap \{x^{2} > 0\}} |\psi|^{2}\right);
\end{eqnarray*}
\begin{eqnarray*}
\sup_{\overline{B_{1/4} \cap \{x^{2} < 0\}}} \, |D\psi|^{2} + \sup_{x, z \in \overline{B_{1/4} \cap \{x^{2} < 0\}}, \, x\neq z} \, 
\frac{|D\psi(x) - D\psi(z)|^{2}}{|x-z|^{2}}&&\\ 
&&\hspace{-2in}\leq C\left(\int_{B_{1/2} \cap \{x^{2} < 0\}}|\varphi|^{2} + \int_{B_{1/2} \cap \{x^{2} > 0\}} |\psi|^{2}\right).
\end{eqnarray*}
Here $C = C(n, q, \a) \in (0, \infty).$
\end{theorem}

\begin{proof}  By the definition of ${\mathcal B}^{F}$, there are sequences $\{V_{k}\} \subset {\mathcal S}_{\a}$, $\{{\mathbf C}_{k}\} \subset {\mathcal C}_{q}$ and sequences of decreasing positive numbers $\{\e_{k}\}$, $\{\g_{k}\}$, $\{\b_{k}\}$, $\{\d_{k}\}$, $\{\t_{k}\}$ for which all of the assertions of Section~\ref{fineexcessblowup} hold, with $M^{2}_{0}$ in place of $M_{0}^{3}$ in condition $(6_{k})$.

By (\ref{firstvariation}), 
\begin{equation}\label{c1alpha-1}
\int_{{\mathbf R} \times B_{1}} \nabla^{V_{k}} \, x^{1} \cdot \nabla^{V_{k}} \, \widetilde{\z} \, d\|V_{k}\|(X) = 0
\end{equation}
for each $k = 1, 2, \ldots$ and any $\widetilde{\z} \in C^{1}_{c}({\mathbf R} \times B_{1}).$ Let $\t \in (0, 1/32)$ be arbitrary. Choose any 
$\z \in C^{2}_{c}(B_{3/8})$ with $\frac{\partial \, \z}{\partial \, x^{2}} \equiv 0$ in $\{|x^{2}| < 2\t\},$ and set 
$\z_{1}(x^{1}, x^{2}, y) = \z(x^{2}, y)$ for $(x^{1}, x^{2}, y) \in {\mathbf R} \times B_{1/2}.$ Let $\widetilde{\z} \in C^{1}_{c}({\mathbf R} \times B_{3/8})$ be such that 
$\widetilde{\z} \equiv \z_{1}$ in a neighborhood of ${\rm spt} \, \|V_{k}\| \cap ({\mathbf R} \times B_{3/8})$ for all 
$k=1, 2, \ldots.$  By (\ref{c1alpha-1}) and (\ref{excess-3}), for all sufficiently large $k$, 
\begin{eqnarray}\label{c1alpha-2}
\int_{{\mathbf R} \times (B_{3/8} \cap \{|x^{2}| < 2\t\})} \nabla^{V_{k}} \, x^{1} \cdot \nabla^{V_{k}} \, \widetilde{\z} \, d\|V_{k}\|(X)&&\nonumber\\ 
&&\hspace{-2.5in}+ \;\sum_{j=1}^{q} \int_{B_{3/8} \cap \{x^{2} \leq -2\t\}} \left(1 + |D(h_{j}^{k} + u_{j}^{k})|^{2}\right)^{-1/2} D(h_{j}^{k} + u_{j}^{k}) \cdot D\z \,dx\nonumber\\
&&\hspace{-2in}+\; \sum_{j=1}^{q} \int_{B_{3/8} \cap \{x^{2}  \geq 2\t\}} \left(1 + |D(g_{j}^{k} + w_{j}^{k})|^{2}\right)^{-1/2} D(g_{j}^{k} + w_{j}^{k}) \cdot D\z \,dx = 0.
\end{eqnarray}   

Since $\frac{\partial \, \widetilde{\z}}{\partial \, x^{1}} = 0$ in a neighborhood of ${\rm spt} \, \|V_{k}\| \cap 
({\mathbf R} \times B_{1/2})$ and $\frac{\partial \, \widetilde{\z}}{\partial \, x^{2}} = 0$ in $\{|x^{2}| < 2\t\},$ 
it follows that 
\begin{eqnarray}\label{c1alpha-3}
\left|\int_{{\mathbf R} \times (B_{3/8} \cap \{|x^{2}| < 2\t\})} \nabla^{V_{k}} \, x^{1} \cdot \nabla^{V_{k}} \, \widetilde{\z} \, d\|V_{k}\|(X)\right|  = \left|\int_{{\mathbf R} \times (B_{3/8} \cap \{|x^{2}| < 2\t\})} e_{1} \cdot \nabla^{V_{k}} \, \widetilde{\z} \, d\|V_{k}\|(X)\right|&&\nonumber\\
&&\hspace{-5in} \leq \sup \, |D\z| \sum_{j=3}^{n+1}\int_{{\mathbf R} \times (B_{3/8} \cap \{|x^{2}| < 2\t\})} |e_{j}^{\perp_{k}}| \, d\|V_{k}\|(X) \nonumber\\
&&\hspace{-5.75in}\leq\sup \, |D\z| \; (\|V_{k}\|({\mathbf R} \times (B_{3/8} \cap \{|x^{2}| < 2\t\}))^{1/2}\left(\sum_{j=3}^{n+1}\int_{{\mathbf R} \times B_{3/8}}|e_{j}^{\perp_{k}}|^{2} \, d\|V_{k}\|(X)\right)^{1/2}\nonumber\\
&&\hspace{-5in} \leq C\sup \, |D\z| \sqrt{\t} E_{k}
\end{eqnarray}
where $C = C(n, q, \a) \in (0, \infty)$ and the last inequality is a consequence of Theorem~\ref{L2-est-1}(c) and the fact that $\|V_{k}\|({\mathbf R} \times (B_{3/8} \cap \{|x^{2}| < 2\t\})) \leq C\t,$ $C = C(n, q, \a) \in (0, \infty),$ for all sufficiently large $k$.

Since $h_{j}^{k}(x) = \lambda_{j}^{k}x^{2}$,  we have for each $j = 1, 2, \ldots, q$ and $k = 1, 2, \ldots,$
\begin{eqnarray}\label{c1alpha-4}
\int_{B_{3/8} \cap \{x^{2} \leq -2\t\}}\left(1 + |D(h_{j}^{k} + u_{j}^{k})|^{2}\right)^{-1/2}D(h_{j}^{k} + u_{j}^{k}) \cdot D\z \,dx &&\nonumber\\
&&\hspace{-4in}=\int_{B_{3/8} \cap \{x^{2} \leq -2\t\}} \left(1 + |D(h_{j}^{k} + u_{j}^{k})|^{2}\right)^{-1/2}Du_{j}^{k}\cdot D\z \,dx\nonumber\\
&&\hspace{-3.5in}-\; \lambda_{j}^{k}\left(1 + (\lambda_{j}^{k})^{2}\right)^{-1/2}\int_{B_{3/8} \cap \{x^{2} \leq -2\t\}} \frac{\left(1 + |D(h_{j}^{k} + u_{j}^{k})|^{2}\right)^{-1/2}D(2h_{j}^{k} + u_{j}^{k}) \cdot Du_{j}^{k}}
{\sqrt{1 + |D(h_{j}^{k} + u_{j}^{k})|^{2}} + \sqrt{1 + (\lambda_{j}^{k})^{2}}}\frac{\partial \, \z}{\partial \, x^{2}} \, dx\nonumber\\ 
&&\hspace{-3in}+  \; \lambda_{j}^{k}\left(1 + (\lambda_{j}^{k})^{2}\right)^{-1/2}\int_{B_{3/8} \cap \{x^{2} \leq -2\t\}}
\frac{\partial \, \z}{\partial \, x^{2}} \, dx.
\end{eqnarray} 

By the Cauchy-Schwarz inequality and elliptic estimates 
\begin{eqnarray}\label{c1alpha-5}
\left|\int_{B_{3/8} \cap \{x^{2} \leq -2\t\}} \frac{\left(1 + |D(h_{j}^{k} + u_{j}^{k})|^{2}\right)^{-1/2}D(2h_{j}^{k} + u_{j}^{k}) \cdot Du_{j}^{k}}{\sqrt{1 + |D(h_{j}^{k} + u_{j}^{k})|^{2}} + \sqrt{1 + (\lambda_{j}^{k})^{2}}} \frac{\partial \, \z}{\partial \, x^{2}}\, dx\right|&&\nonumber\\
&&\hspace{-3in}\leq C(\t) \, \sup \, |D\z|\, \sqrt{|\lambda_{j}^{k}|^{2} +\int_{B_{1/2} \cap \{x^{2} \leq -\t\}} |u_{j}^{k}|^{2} \, dx}
\sqrt{ \int_{B_{1/2} \cap \{x^{2} \leq -\t\}} |u_{j}^{k}|^{2} \, dx}\nonumber\\
&&\hspace{-2in}\leq C(\t) \sup \, |D\z| \, \sqrt{|\lambda_{j}^{k}|^{2} + E_{k}^{2}} \; E_{k}.
\end{eqnarray}

If $\z$ also satisfies 
\begin{equation}\label{c1alpha-5-1}
\int_{B_{3/8} \cap (\{0\} \times {\mathbf R}^{n-1})} \z \, dy = 0
\end{equation}
then, since $\int_{B_{3/8} \cap \{x^{2} \leq -2\t\}} \frac{\partial \, \z}{\partial \, x^{2}} \, dx = \int_{B_{3/8} \cap \{x^{2} \leq 0\}} 
\frac{\partial \, \z}{\partial \, x^{2}} \, dx = -\int_{B_{3/8} \cap (\{0\} \times {\mathbf R}^{n-1})} \z \, dy$, the last term on the right hand side of (\ref{c1alpha-4}) will be zero. Thus, for each fixed $\t >0$ and 
each $\z \in C^{1}_{c}(B_{3/8})$ with $\frac{\partial \, \z}{\partial \,x^{2}} = 0$ in $\{|x^{2}| < 2\t\}$ and satisfying 
(\ref{c1alpha-5-1}), we have  

\begin{eqnarray}\label{c1alpha-6}
\sum_{j=1}^{q}\int_{B_{3/8} \cap \{x^{2} \leq -2\t\}}\left(1 + |D(h_{j}^{k} + u_{j}^{k})|^{2}\right)^{-1/2} D(h_{j}^{k} + u_{j}^{k}) \cdot D\z \,dx &&\nonumber\\
&&\hspace{-2.5in}=\sum_{j=1}^{q}\int_{B_{3/8} \cap \{x^{2} \leq -2\t\}} \left(1 + |D(h_{j}^{k} + u_{j}^{k})|^{2}\right)^{-1/2}Du_{j}^{k}\cdot D\z \,dx + \e^{-}_{k}
\end{eqnarray}
and,  by a similar argument, 
\begin{eqnarray}\label{c1alpha-7}
\sum_{j=1}^{q}\int_{B_{3/8} \cap \{x^{2} \geq 2\t\}}\left(1 + |D(g_{j}^{k} + w_{j}^{k})|^{2}\right)^{-1/2}D(g_{j}^{k} + w_{j}^{k}) \cdot D\z \,dx &&\nonumber\\
&&\hspace{-2.5in}=\sum_{j=1}^{q}\int_{B_{3/8} \cap \{x^{2} \geq 2\t\}} \left(1 + |D(g_{j}^{k} + w_{j}^{k})|^{2}\right)^{-1/2}Dw_{j}^{k}\cdot D\z \,dx + \e^{+}_{k}
\end{eqnarray}
where $\lim_{k\to \infty} \, E_{k}^{-1}|\e_{k}^{-}| =\lim_{k \to \infty} \, E_{k}^{-1}|\e_{k}^{+}|= 0.$ We may divide (\ref{c1alpha-2}) by $E_{k}$ and let $k \to \infty$ to deduce, by 
(\ref{c1alpha-3}), (\ref{c1alpha-6}), (\ref{c1alpha-7}) and (\ref{excess-8}),  that for each $\t \in (0, 1/16)$, 
\begin{equation}\label{c1alpha-8}
\sum_{j=1}^{q}\int_{B_{3/8} \cap \{x^{2} \leq -2\t\}} D\varphi_{j} \cdot D\z + 
\sum_{j=1}^{q} \int_{B_{3/8} \cap \{x^{2} \geq 2\t\}} D\psi_{j} \cdot D\z + \e(\t) = 0
\end{equation}
where $\e(\t) \to 0$ as $\t \to 0.$ Upon integration by parts (in view of the fact that $\frac{\partial \, \z}{\partial \, x^{2}} = 0$ in $\{|x^{2}| < 2\t\}$), this gives
\begin{equation}\label{c1alpha-9}
\sum_{j=1}^{q}\int_{B_{3/8} \cap \{x^{2} \leq -2\t\}} \varphi_{j} \Delta \,\z + 
\sum_{j=1}^{q} \int_{B_{3/8} \cap \{x^{2} \geq 2\t\}} \psi_{j}  \Delta \, \z - \e(\t) = 0.
\end{equation}
Since $\varphi_{j} \in L^{1}(B_{1/2} \cap \{x^{2} \leq 0\})$ and $\psi_{j} \in L^{1}(B_{1/2} \cap \{x^{2} \geq 0\})$ for 
each $j=1, 2, \ldots, q$, we may let $\t \to 0$ in (\ref{c1alpha-9}) to conclude that 
\begin{equation}\label{c1alpha-10}
\sum_{j=1}^{q}\int_{B_{3/8} \cap \{x^{2} \leq 0\}} \varphi_{j} \Delta \,\z + 
\sum_{j=1}^{q} \int_{B_{3/8} \cap \{x^{2} \geq 0\}} \psi_{j}  \Delta \, \z  = 0
\end{equation}
for any  $\z \in C^{2}_{c}(B_{3/8})$ with 
$\frac{\partial \, \z}{\partial \, x^{2}} = 0$ in a neighborhood of $\{x^{2} = 0\}$ and satisfying 
(\ref{c1alpha-5-1}).

Now for any $\ell \in \{1, 2, \ldots, n-1\},$ $h \in (-1/16, 1/16)$ and any 
$\z \in C^{2}_{c}(B_{5/16})$ with $\frac{\partial \, \z}{\partial \, x^{2}} = 0$ in a neighborhood 
of $\{x^{2} = 0\}$, we have that $\d_{\ell, h} \, \z \in C^{1}_{c}(B_{3/8})$, 
$\frac{\partial}{\partial \, x^{2}} \, \d_{\ell, h} \, \z= 0$ in a neighborhood of $\{x^{2} = 0\}$ and $\d_{\ell, h} \, \z$ satisfies (\ref{c1alpha-5-1}), 
where $\d_{\ell, h} \, \z (x^{2}, y) = \z(x^{2}, y^{1}, \ldots, y^{\ell} + h, y^{\ell +1}, \ldots, y^{n-1}) - \z(x^{2}, y).$
Thus, by (\ref{c1alpha-10}), 
\begin{equation*}\label{c1alpha-11}
\sum_{j=1}^{q}\int_{B_{3/8} \cap \{x^{2} < 0\}} \varphi_{j} \,\Delta \, \d_{\ell, h} \, \z + 
\sum_{j=1}^{q} \int_{B_{3/8} \cap \{x^{2} > 0\}} \psi_{j} \, \Delta \, \d_{\ell, h} \, \z  = 0
\end{equation*}
and consequently,
\begin{equation}\label{c1alpha-12}
\sum_{j=1}^{q}\int_{B_{5/16} \cap \{x^{2} < 0\}} \d_{\ell, h} \, \varphi_{j} \,\Delta \,  \z + 
\sum_{j=1}^{q} \int_{B_{5/16} \cap \{x^{2} > 0\}} \d_{\ell, h} \, \psi_{j} \, \Delta \, \z  = 0
\end{equation}
for  any $\z \in C^{2}_{c}(B_{5/16})$ with $\frac{\partial \, \z}{\partial \, x^{2}} = 0$ in a neighborhood 
of $\{x^{2} = 0\},$ any $\ell \in \{1, 2, \ldots, n-1\}$ and $h \in (-1/16, 1/16).$ Since any 
$\z \in C^{2}_{c}(B_{5/16})$ which is even in the $x^{2}$ variable can be approximate in $C^{2}(B_{5/16})$ by 
a sequence $\z_{i} \in C^{2}_{c}(B_{5/16})$ satisfying, for each $i = 1, 2, 3, \ldots$, 
$\frac{\partial \, \z_{i}}{\partial \, x^{2}} = 0$ in a neighborhood of $\{x^{2} = 0\},$ we see that 
(\ref{c1alpha-12}) holds for any $\z \in C^{2}_{c}(B_{5/16})$ which is even in the $x^{2}$ variable and 
for each $\ell \in \{1, 2, \ldots, n-1\}$ and $h \in (-1/16, 1/16).$ Thus 
\begin{equation}\label{c1alpha-13}
\int_{B_{5/16}} \Phi_{\ell, h} \Delta \, \z = 0
\end{equation}
for any $\z \in C^{2}_{c}(B_{5/16})$ which is even in the $x^{2}$ variable, any 
$\ell \in \{1, 2, \ldots, n-1\}$ and $h \in (-1/16, 1/16)$, where $\Phi_{\ell, h} \, : \, B_{3/8} \to {\mathbf R}$ is the function defined by
$\Phi_{\ell, h}(x^{2}, y) = \sum_{j=1}^{q} \d_{\ell, h} \, \varphi_{j}(-x^{2}, y) + \d_{\ell, h} \, \psi_{j}(x^{2}, y)$ if $x^{2} \geq 0$ and 
$\Phi_{\ell, h}(x^{2}, y) = \sum_{j=1}^{q} \d_{\ell, h} \, \varphi_{j}(x^{2}, y) + \d_{\ell, h} \, \psi_{j}(-x^{2}, y)$ if $x^{2} <0.$ Since $\Phi$ is even in the $x^{2}$ variable, (\ref{c1alpha-13}) holds also for any $\z$ which is odd in the $x^{2}$ variable.
Thus (\ref{c1alpha-13}) holds for every $\z \in C^{2}_{c}(B_{5/16})$ and hence $\Phi_{\ell, h}$ is a smooth harmonic function in $B_{5/16}.$ Since we have directly from the definition of $\Phi_{\ell, h}$ and Lemma~\ref{continuity} that
\begin{equation}\label{c1alpha-14}
\left|\int_{B_{5/16}}h^{-1}\Phi_{\ell,h}\right| \leq C\left(\int_{B_{1/2} \cap \{x^{2} \leq 0\}} |\varphi|^{2} + \int_{B_{1/2} \cap \{x^{2} \geq 0\}} |\psi|^{2}\right)^{1/2}
\end{equation}
for all $h \in (-1/16,1/16) \setminus \{0\}$, where $C = C(n, q, \a) \in (0, \infty)$, it follows from standard estimates for harmonic functions that there exists a harmonic function $\Phi_{\ell} \, : \, B_{9/32} \to {\mathbf R}$ such that  
$h^{-1}\Phi_{\ell, h} \to \Phi_{\ell}$ in $C^{2}(B_{9/32})$ as $h \to 0,$ and 
\begin{equation}\label{c1alpha-15}
\sup_{B_{9/32}} \, |\Phi_{\ell}|^{2} + |D\Phi_{\ell}|^{2} + |D^{2}\Phi_{\ell}|^{2} \leq C\left(\int_{B_{1/2} \cap \{x^{2} \leq 0\}} |\varphi|^{2} + \int_{B_{1/2} \cap \{x^{2} \geq 0\}} |\psi|^{2}\right).
\end{equation}

Let $\Phi \, : \, B_{1/2} \to {\mathbf R}$ be the function defined by 
$\Phi(x^{2}, y) = \sum_{j=1}^{q} \varphi_{j}(x^{2}, y) + \psi_{j}(-x^{2}, y)$ if $x^{2} < 0$ and 
$\Phi(x^{2}, y) = \sum_{j=1}^{q} \varphi_{j}(-x^{2}, y) + \psi_{j}(x^{2}, y)$ if $x^{2}  \geq 0.$ Since $\Phi_{\ell} =  \frac{\partial}{\partial \, y^{\ell}} \, \Phi$ on $B_{1/2} \setminus (\{0\} \times {\mathbf R}^{n-1})$, it follows that for $(x^{2}, y) \in B_{9/32} \setminus (\{0\} \times {\mathbf R}^{n-1})$, 
\begin{equation*}
\Phi(x^{2}, y) = \Phi(x^{2},y^{1},  \ldots, y^{\ell - 1}, 0, y^{\ell+1}, \ldots, y^{n-1}) +  \int_{0}^{y^{\ell}} \Phi_{\ell}(x^{2}, y^{1}, \ldots, y^{\ell-1}, t, y^{\ell+1}, \ldots, y^{n-1}) \, dt,
\end{equation*}
so we may let $x^{2} \to 0$ on both sides of this and use Lemma~\ref{continuity},  (\ref{continuity-1-1-1}) and the arbitrariness of the index $\ell \in \{1, 2, \ldots, n-1\}$ to conclude that, with $Y = (0, y)$,
\begin{equation}\label{c1alpha-15-1}
\Phi(Y) = 2q\k_{1}(Y) - \left(\sum_{j=1}^{q}  (\ell_{j} + m_{j})\right)\k_{2}(Y)
\end{equation}
is a $C^{\infty}$ function of 
$Y \in B_{9/32} \cap (\{0\} \times {\mathbf R}^{n-1})$ (with $\frac{\partial}{\partial \, y^{\ell}} \Phi (Y) = \Phi_{\ell}(Y),$ $\frac{\partial^{2}}{\partial \, y^{m} \, \partial \, y^{\ell}} \Phi(Y) = \frac{\partial}{\partial \, y^{m}} \, \Phi_{\ell}(Y)$, $\frac{\partial^{3}}{\partial \, y^{k} \, \partial \, y^{m} \, \partial \, y^{\ell}} \Phi(Y) = \frac{\partial^{2}}{\partial \, y^{k} \, \partial \, y^{m}} \, \Phi_{\ell}(Y)$
for each $\ell, m, k \in \{1, 2, \ldots, n-1\}$) satisfying, by (\ref{c1alpha-15}) and Lemma~\ref{continuity}, the estimate
\begin{equation}\label{c1alpha-16}
\sup_{B_{9/32} \cap (\{0\} \times {\mathbf R}^{n-1})}\, |\Phi|^{2} + |D_{Y} \,\Phi|^{2} + |D_{Y}^{2} \, \Phi|^{2} + |D_{Y}^{3} \, \Phi|^{2} \leq C\left(\int_{B_{1/2} \cap \{x^{2} \leq 0\}} |\varphi|^{2} + \int_{B_{1/2} \cap \{x^{2} \geq 0\}} |\psi|^{2}\right)
\end{equation}
where $C = C(n, q, \a) \in (0, \infty).$

Next we derive regularity estimates for a different linear combination of $\k_{1}$ and $\k_{2}.$ For this, we note 
that by (\ref{firstvariation}) again, 
\begin{equation}\label{c1alpha-17}
\int_{{\mathbf R} \times B_{1}} \nabla^{V_{k}} \, x^{2} \cdot \nabla^{V_{k}} \, \widetilde{\z} \, d\|V_{k}\|(X) = 0
\end{equation}
for each $k=1, 2, \ldots$ and each $\widetilde{\z} \in C^{1}_{c}({\mathbf R} \times B_{1}).$ Let $\t \in (0, 1/16),$  
$\z \in C^{2}_{c}(B_{3/8})$ and $\widetilde{\z}$ be as before, so that in particular $\frac{\partial \, \z}{\partial \, x^{2}} = 0$ in $\{|x^{2}| < 2\t\}.$ Note that the unit normal $\nu_{j}^{k}$ to $(\mathcal{M}_{j}^{k})^{-} \equiv {\rm graph} \, (h_{j}^{k} + u_{j}^{k})$  is given by 
$\nu_{j}^{k} = \left((1 + |D(h_{j}^{k} + u_{j}^{k})|^{2}\right)^{-1/2}\left(1, -\lambda_{j}^{k} -\frac{\partial \, u_{j}^{k}}{\partial \, x^{2}}, -D_{y} \, u_{j}^{k}\right)$ so that, on $(\mathcal{M}_{j}^{k})^{-},$
\begin{eqnarray*}
\nabla^{V_{k}} \, x^{2} \cdot \nabla^{V_{k}} \, \widetilde{\z}  = e_{2} \cdot\left(D\widetilde{\z} - 
(D\widetilde{\z} \cdot \nu_{j}^{k})\nu_{j}^{k}\right)&&\nonumber\\
 &&\hspace{-2.5in}= \frac{\partial \, \z}{\partial \, x^{2}} - \left( 1 + |D(h_{j}^{k} + u_{j}^{k})|^{2}\right)^{-1}\left(\lambda_{j}^{k} + \frac{\partial \, u_{j}^{k}}{\partial \, x^{2}}\right)\left(\left(\lambda_{j}^{k} + \frac{\partial \, u_{j}^{k}}{\partial \, x^{2}}\right)\frac{\partial \, \z}{\partial \, x^{2}} + D_{y}u_{j}^{k} \cdot D_{y} \z\right)\nonumber\\
&&\hspace{-2.5in}= \left( 1 + |D(h_{j}^{k} + u_{j}^{k})|^{2}\right)^{-1}
\left(\left(1 + |D_{y} \, u_{j}^{k}|^{2}\right) \frac{\partial \, \z}{\partial \, x^{2}} - \left(\lambda_{j}^{k} + \frac{\partial \, u_{j}^{k}}{\partial \, x^{2}}\right)D_{y} \, u_{j}^{k} \cdot D_{y} \, \z\right). 
\end{eqnarray*}
Using this and similar expressions for $\nabla^{V_{k}} \, x^{2} \cdot \nabla^{V_{k}} \, \widetilde{\z}$ on $(\mathcal{M}_{j}^{k})^{+} \equiv {\rm graph} \, (g_{j}^{k} + w_{j}^{k}),$ we deduce from (\ref{c1alpha-17}) and
Theorem~\ref{L2-est-1}(a) that 
\begin{eqnarray}\label{c1alpha-18}
\int_{{\mathbf R} \times (B_{3/8} \cap \{|x^{2}| < 2\t\})} \nabla^{V_{k}} \, x^{2} \cdot \nabla^{V_{k}} \, \widetilde{\z} \, d\|V_{k}\|(X)&&\nonumber\\ 
&&\hspace{-3in}+ \;\sum_{j=1}^{q} \int_{B_{3/8} \cap \{x^{2} \leq -2\t\}} \frac{(1 + |D_{y} \, u_{j}^{k}|^{2}) \frac{\partial \, \z}{\partial \, x^{2}} - (\lambda_{j}^{k} + \frac{\partial \, u_{j}^{k}}{\partial \, x^{2}})D_{y} \, u_{j}^{k} \cdot D_{y} \, \z}
{\sqrt{ 1 + |D(h_{j}^{k} + u_{j}^{k})|^{2}}}\,dx\nonumber\\
&&\hspace{-2in}+\; \sum_{j=1}^{q} \int_{B_{3/8} \cap \{x^{2}  \geq 2\t\}} \frac{(1 + |D_{y} \, w_{j}^{k}|^{2}) \frac{\partial \, \z}{\partial \, x^{2}} - (\mu_{j}^{k} + \frac{\partial \, w_{j}^{k}}{\partial \, x^{2}})D_{y} \,w_{j}^{k} \cdot D_{y} \, \z}{\sqrt{ 1 + |D(g_{j}^{k} + w_{j}^{k})|^{2}}}
\,dx = 0.
\end{eqnarray}  

Since $\frac{\partial \, \widetilde{\z}}{\partial \, x^{1}} = 0$ in a neighborhood of ${\rm spt} \, \|V_{k}\| \cap 
({\mathbf R} \times B_{1/2})$ and $\frac{\partial \, \widetilde{\z}}{\partial \, x^{2}} = 0$ in $\{|x^{2}| < 2\t\},$ 
it follows that
\begin{eqnarray}\label{c1alpha-19}
\left|\int_{{\mathbf R} \times (B_{3/8} \cap \{|x^{2}| < 2\t\})} \nabla^{V_{k}} \, x^{2} \cdot \nabla^{V_{k}} \, \widetilde{\z} \, d\|V_{k}\|(X)\right|  =\left|\int_{{\mathbf R} \times (B_{3/8} \cap \{|x^{2}| < 2\t\})} e_{2} \cdot \nabla^{V_{k}} \, \widetilde{\z} \, d\|V_{k}\|(X)\right|&&\nonumber\\
&&\hspace{-6in}\leq \sup \, |D\z| \sum_{j=3}^{n+1}\int_{{\mathbf R} \times (B_{3/8} \cap \{|x^{2}| < 2\t\})}|e_{2}^{\perp_{k}}||e_{j}^{\perp_{k}}| \, d\|V_{k}\|(X)\nonumber\\
&&\hspace{-6in}\leq\sup \,|D\z| \left(\int_{{\mathbf R} \times (B_{3/8} \cap \{|x^{2}| < 2\t\})}|e_{2}^{\perp_{k}}|^{2} \, d\|V_{k}\|(X)\right)^{1/2}  \left(\sum_{j=3}^{n+1}\int_{{\mathbf R} \times B_{3/8}}|e_{j}^{\perp_{k}}|^{2} \, d\|V_{k}\|(X)\right)^{1/2}\nonumber\\
&&\hspace{-6in}\leq\sup \, |D\z|\left(\int_{{\mathbf R} \times (B_{3/8} \cap \{|x^{2}| < 2\t\})}1 - |e_{1}^{\perp_{k}}|^{2} \, d\|V_{k}\|(X)\right)^{1/2} \left(\sum_{j=3}^{n+1}\int_{{\mathbf R} \times B_{3/8}}|e_{j}^{\perp_{k}}|^{2} \, d\|V_{k}\|(X)\right)^{1/2}\nonumber\\
&&\hspace{-5.5in}\leq C \, \sup \, |D\z| \,\t^{1/4}{\hat E}_{k}E_{k}
\end{eqnarray}
for all sufficiently large $k$, where $C = C(n, q, \a) \in (0, \infty)$ and the last inequality follows from Theorem~\ref{L2-est-1}(c), Theorem~\ref{non-concentration}(b) and Lemma~\ref{no-gaps}. Since
\begin{eqnarray}\label{c1alpha-20}
\int_{B_{3/8} \cap \{x^{2} \leq -2\t\}} \frac{(1 + |D_{y} \, u_{j}^{k}|^{2}) \frac{\partial \, \z}{\partial \, x^{2}} - (\lambda_{j}^{k} + \frac{\partial \, u_{j}^{k}}{\partial \, x^{2}})D_{y} \, u_{j}^{k} \cdot D_{y} \, \z}
{\sqrt{ 1 + |D(h_{j}^{k} + u_{j}^{k})|^{2}}}\,dx\nonumber\\
&&\hspace{-4.5in}=-\int_{B_{3/8} \cap \{x^{2} \leq -2\t\}} \frac{\lambda_{j}^{k}D_{y} \, u_{j}^{k} \cdot D_{y} \, \z}
{\sqrt{ 1 + |D(h_{j}^{k} + u_{j}^{k})|^{2}}}\,dx +\int_{B_{3/8} \cap \{x^{2} \leq -2\t\}} \frac{|D_{y} \, u_{j}^{k}|^{2}\frac{\partial \, \z}{\partial \, x^{2}} -  \frac{\partial \, u_{j}^{k}}{\partial \, x^{2}}D_{y} \, u_{j}^{k} \cdot D_{y} \, \z}
{\sqrt{ 1 + |D(h_{j}^{k} + u_{j}^{k})|^{2}}}\,dx\nonumber\\
&&\hspace{-4.25in}- \int_{B_{3/8} \cap \{x^{2} \leq -2\t\}} \frac{(2\lambda_{j}^{k}\frac{\partial \, u_{j}^{k}}{\partial \, x^{2}} + |D \, u_{j}^{k}|^{2}) \frac{\partial \, \z}{\partial \, x^{2}}}
{\sqrt{1 + |\lambda_{j}^{k}|^{2}}\sqrt{ 1 + |D(h_{j}^{k} + u_{j}^{k})|^{2}} \left(\sqrt{1 + |\lambda_{j}^{k}|^{2}} + \sqrt{ 1 + |D(h_{j}^{k} + u_{j}^{k})|^{2}}\right)}\,dx\nonumber\\
&&\hspace{-3.5in}+\; \frac{1}{\sqrt{ 1 + |\lambda_{j}^{k}|^{2}}}\int_{B_{3/8} \cap \{x^{2} \leq -2\t\}}  \frac{\partial \, \z}{\partial \, x^{2}} \,dx,
\end{eqnarray}
it follows that if $\z$ satisfies also (\ref{c1alpha-5-1}), then 
\begin{eqnarray}\label{c1alpha-21}
\int_{B_{3/8} \cap \{x^{2} \leq -2\t\}} \frac{(1 + |D_{y} \, u_{j}^{k}|^{2}) \frac{\partial \, \z}{\partial \, x^{2}} - (\lambda_{j}^{k} + \frac{\partial \, u_{j}^{k}}{\partial \, x^{2}})D_{y} \, u_{j}^{k} \cdot D_{y} \, \z}
{\sqrt{ 1 + |D(h_{j}^{k} + u_{j}^{k})|^{2}}}\,dx\nonumber\\
&&\hspace{-4.5in}=-\int_{B_{3/8} \cap \{x^{2} \leq -2\t\}} \frac{2\lambda_{j}^{k}\frac{\partial \, u_{j}^{k}}{\partial \, x^{2}} \frac{\partial \, \z}{\partial \, x^{2}}}
{\sqrt{1 + |\lambda_{j}^{k}|^{2}}\sqrt{ 1 + |D(h_{j}^{k} + u_{j}^{k})|^{2}} \left(\sqrt{1 + |\lambda_{j}^{k}|^{2}} + \sqrt{ 1 + |D(h_{j}^{k} + u_{j}^{k})|^{2}}\right)}\,dx\nonumber\\
&&\hspace{-3.5in}- \; \int_{B_{3/8} \cap \{x^{2} \leq -2\t\}} \frac{\lambda_{j}^{k}D_{y} \, u_{j}^{k} \cdot D_{y} \, \z}
{\sqrt{ 1 + |D(h_{j}^{k} + u_{j}^{k})|^{2}}}\,dx + \eta_{k}^{-}\end{eqnarray}
where, by elliptic estimates, 
\begin{equation}\label{c1alpha-22}
|\eta_{k}^{-}| \leq C \sup |D\z| \, E_{k}^{2}, \;\;\;\;  C = C(n, q, \t) \in (0, \infty).
\end{equation} 
By the same argument, 
\begin{eqnarray}\label{c1alpha-23}
\int_{B_{3/8} \cap \{x^{2} \geq 2\t\}} \frac{(1 + |D_{y} \, w_{j}^{k}|^{2}) \frac{\partial \, \z}{\partial \, x^{2}} - (\mu_{j}^{k} + \frac{\partial \, w_{j}^{k}}{\partial \, x^{2}})D_{y} \, w_{j}^{k} \cdot D_{y} \, \z}
{\sqrt{ 1 + |D(g_{j}^{k} + w_{j}^{k})|^{2}}}\,dx\nonumber\\
&&\hspace{-4.5in}=-\int_{B_{3/8} \cap \{x^{2} \geq 2\t\}} \frac{2\mu_{j}^{k}\frac{\partial \, w_{j}^{k}}{\partial \, x^{2}} \frac{\partial \, \z}{\partial \, x^{2}}}
{\sqrt{1 + |\mu_{j}^{k}|^{2}}\sqrt{ 1 + |Dg_{j}^{k} + w_{j}^{k})|^{2}} \left(\sqrt{1 + |\mu_{j}^{k}|^{2}} + \sqrt{ 1 + |D(g_{j}^{k} + w_{j}^{k})|^{2}}\right)}\,dx\nonumber\\
&&\hspace{-3.5in}- \; \int_{B_{3/8} \cap \{x^{2} \geq 2\t\}} \frac{\mu_{j}^{k}D_{y} \, w_{j}^{k} \cdot D_{y} \, \z}
{\sqrt{ 1 + |D(g_{j}^{k} + u_{j}^{k})|^{2}}}\,dx + \eta_{k}^{+}\end{eqnarray}
where, again by elliptic estimates, 
\begin{equation}\label{c1alpha-24}
|\eta_{k}^{+}| \leq C \sup |D\z| \, E_{k}^{2}, \;\;\;\;  C = C(n, q, \t) \in (0, \infty).
\end{equation} 

Dividing (\ref{c1alpha-18}) by ${\hat E}_{k}E_{k}$ and letting $k \to \infty$, we conclude with the help of 
(\ref{c1alpha-19}), (\ref{c1alpha-21}), (\ref{c1alpha-22}), (\ref{c1alpha-23}), (\ref{c1alpha-24}), (\ref{excess-0}) and 
(\ref{excess-7}) that 
\begin{equation}\label{c1alpha-25}
\sum_{j=1}^{q}\ell_{j}\int_{B_{3/8} \cap \{x^{2} \leq -2\t\}} D\varphi_{j} \cdot D\z + 
\sum_{j=1}^{q} m_{j}\int_{B_{3/8} \cap \{x^{2} \geq 2\t\}} D\psi_{j} \cdot D\z + \eta(\t) = 0
\end{equation}
for any $\z \in C^{2}_{c}(B_{3/8})$ with $\frac{\partial \, \z}{\partial \, x^{2}} = 0$ in $\{|x^{2}| < 2\t\}$ and satisfying 
(\ref{c1alpha-5-1}), where $\eta(\t) \to 0$ as $\t \to 0.$ 

It follows from (\ref{c1alpha-25}) in the same way that (\ref{c1alpha-16}) follows from (\ref{c1alpha-8}) that if we let, for $Y \in B_{3/8} \cap (\{0\} \times {\mathbf R}^{n-1}),$ 
\begin{equation}\label{c1alpha-26}
\Psi(Y) = \left(\sum_{j=1}^{q}(\ell_{j} + m_{j})\right)\k_{1}(Y)  - \left(\sum_{j=1}^{q} (\ell_{j}^{2} + m_{j}^{2})\right)\k_{2}(Y),
\end{equation}
then $\Psi$ is a $C^{\infty}$ function on $B_{9/32} \cap (\{0\} \times {\mathbf R}^{n-1})$ satisfying the estimate 
\begin{equation}\label{c1alpha-27}
\sup_{B_{9/32} \cap (\{0\} \times {\mathbf R}^{n-1})}\, |\Psi|^{2} + |D_{Y} \, \Psi|^{2} + |D_{Y}^{2} \, \Psi|^{2} + |D_{Y}^{3} \, \Psi|^{2} \leq C\left(\int_{B_{1/2} \cap \{x^{2} \leq 0\}} |\varphi|^{2} + \int_{B_{1/2} \cap \{x^{2} \geq 0\}} |\psi|^{2}\right)
\end{equation}
where $C = C(n, q, \a) \in (0, \infty).$ 

Note that 
$$J \equiv 2q\sum_{j=1}^{q}(\ell_{j}^{2} + m_{j}^{2}) - \left(\sum_{j=1}^{q}(\ell_{j} + m_{j})\right)^{2} = \frac{1}{2}\sum_{i=1}^{q}\sum_{j=1}^{q} \left((m_{i} - m_{j})^{2} + (\ell_{i} - \ell_{j})^{2} + 2(\ell_{i} - m_{j})^{2}\right),$$
and so it follows from (\ref{excess-6}) that $\widetilde{C} \geq J \geq C >0$, where $\widetilde{C} =\widetilde{C}(n, q) \in (0, \infty)$ and $C = C(n, q) \in (0, \infty);$ thus, by (\ref{c1alpha-15-1}) and (\ref{c1alpha-26}), we may express each of $\k_{1}$ and $\k_{2}$  as a linear combination of 
$\Phi$ and $\Psi$ with coefficients, in absolute value, $\leq C = C(n, q) \in (0, \infty).$ Consequently, $\k_{1}, \k_{2}$ are in $C^{\infty}(B_{9/32} \cap (\{0\} \times {\mathbf R}^{n-1}))$ and, by (\ref{c1alpha-16}) and 
(\ref{c1alpha-27}), satisfy  the estimates

\begin{equation}\label{c1alpha-28}
\sup_{B_{9/32} \cap (\{0\} \times {\mathbf R}^{n-1})}\, |\k_{i}|^{2} + |D_{y} \, \k_{i}|^{2} + |D^{2}_{y} \, \k_{i}|^{2} + |D^{3}_{y} \, \k_{i}|^{2}\leq C\left(\int_{B_{1/2} \cap \{x^{2} \leq 0\}} |\varphi|^{2} + \int_{B_{1/2} \cap \{x^{2} \geq 0\}} |\psi|^{2}\right)
\end{equation}
for $i=1, 2$, where $C = C(n, q, \a) \in (0, \infty).$ This in turn implies that for each $j=1, 2, \ldots, q$, 
the functions $\left.\varphi_{j}\right|_{B_{9/32} \cap (\{0\} \times {\mathbf R}^{n-1})} \left(=\k_{1} - \ell_{j}\k_{2}\right)$ and $\left.\psi_{j}\right|_{B_{9/32} \cap (\{0\} \times {\mathbf R}^{n-1})} \left(=\k_{1} - m_{j}\k_{2}\right)$ belong to $C^{\infty}\left(B_{9/32} \cap (\{0\} \times {\mathbf R}^{n-1})\right)$ and satisfy the estimates
\begin{equation}\label{c1alpha-29}
\sup_{B_{9/32} \cap (\{0\} \times {\mathbf R}^{n-1})}\, |\varphi_{j}|^{2} + |D_{y} \, \varphi_{j}|^{2} + |D_{y}^{2} \,\varphi_{j}|^{2} + |D^{3}_{y} \, \varphi_{j}|^{2}\leq C\left(\int_{B_{1/2} \cap \{x^{2} \leq 0\}} |\varphi|^{2} + \int_{B_{1/2} \cap \{x^{2} \geq 0\}} |\psi|^{2}\right),
\end{equation}
\begin{equation}\label{c1alpha-30}
\sup_{B_{9/32} \cap (\{0\} \times {\mathbf R}^{n-1})}\, |\psi_{j}|^{2} + |D_{y} \, \psi_{j}|^{2} + |D_{y}^{2}\, \psi_{j}|^{2} + |D^{3}_{y} \, \psi_{j}|^{2}\leq C\left(\int_{B_{1/2} \cap \{x^{2} \leq 0\}} |\varphi|^{2} + \int_{B_{1/2} \cap \{x^{2} \geq 0\}} |\psi|^{2}\right)
\end{equation}
where $C = C(n, q, \a) \in (0, \infty).$ By Lemma~\ref{continuity} and the standard $C^{2, \a}$ boundary regularity theory for harmonic functions (\cite{M}), the desired conclusions 
of the present lemma in particular follow.
\end{proof}

\section{Improvement of fine excess}\label{finedecay}
\setcounter{equation}{0}
Let $q$ be an integer $\geq 2$, $\a \in (0, 1)$ and suppose that the induction hypotheses $(H1)$, $(H2)$ hold. The main result of this section (Lemma~\ref{multi-scale-final} below) establishes that there are fixed constants $\e = \e(n, q, \a) \in (0, 1)$, $\g = \g(n, q, \a) \in (0, 1)$ such that  whenever $V \in {\mathcal S}_{\a},$ ${\mathbf C} \in {\mathcal C}_{q}$ satisfy Hypotheses~\ref{hyp}, Hypothesis~($\star$) (of Section~\ref{fineblowup}) with a suitable constant $M$ depending only on $n$ and $q$,  
the fine excess of $V$ relative to a new cone ${\mathbf C}^{\prime} \in {\mathcal C}_{q}$ decays by a fixed factor at 
one of several fixed smaller scales.   

\begin{lemma}\label{excess-improvement} 
Let $q$ be an integer $\geq 2$,  $\a \in (0, 1)$ and $\th \in (0, 1/4)$. 
There exist numbers $\overline\e = \overline\e(n, q, \a, \th) \in (0, 1/2)$, $\overline\g = \overline\g(n, q, \a, \th) \in (0, 1/2)$ and 
$\overline\b = \overline\b(n, q, \a,\th) \in (0, 1/2)$ such that the following 
is true: If $V \in {\mathcal S}_{\a},$ ${\mathbf C} \in {\mathcal C}_{q}$ satisfy Hypotheses~\ref{hyp}, Hypothesis~($\star$), 
Hypothesis~($\star\star$) (of Section~\ref{fineblowup}) with $\e  = \overline\e$, $\g  = \overline\g,$ $M = \frac{3}{2}M_{0}$, $\b = \overline{\b}$,  and if the induction hypotheses $(H1)$, $(H2)$ hold, then there exist an orthogonal rotation $\G$ of ${\mathbf R}^{n+1}$ and a cone ${\mathbf C}^{\prime} \in {\mathcal C}_{q}$ 
such that, with 
$${\hat E}_{V}^{2} = \int_{{\mathbf R} \times B_{1}} |x^{1}|^{2} \, d\|V\|(X) \;\;\; {\rm and} \;\;\; E_{V}^{2} =  \int_{{\mathbf R} \times B_{1}} {\rm dist}^{2} \, (X, {\rm spt} \, \|{\mathbf C})\|) \, d\|V\|(X),$$ 
the following hold:
\begin{eqnarray*}
&&({\rm a}) \;\;\;\; |e_{1} - \G(e_{1})| \leq \overline{\k}E_{V} \;\;{\rm and} \;\; |e_{j} - \G(e_{j})| \leq \overline{\k}{\hat E}_{V}^{-1}E_{V}\;\; \mbox{for each} \;\; j=2, 3, \ldots, n+1;\nonumber\\   
&&({\rm b})\;\;\;\; {\rm dist}_{\mathcal H}^{2} \, ({\rm spt} \, \|{\mathbf C}^{\prime}\| \cap ({\mathbf R} \times B_{1}), {\rm spt} \, \|{\mathbf C}\| \cap ({\mathbf R} \times B_{1})) \leq \overline{C}_{0}E_{V}^{2};\nonumber\\
&&({\rm c})\;\;\;\;\th^{-n-2}\int_{\G\left({\mathbf R} \times \left(B_{\th/2} \setminus \{|x^{2}| \leq \th/16\}\right)\right)} {\rm dist}^{2} \, (X, {\rm spt} \, \|V\|) \, d\|\G_{\#} \, {\mathbf C}^{\prime}\|(X)\nonumber\\
&&\hspace{1.5in}+ \;\th^{-n-2}\int_{\G({\mathbf R} \times B_{\th})} {\rm dist}^{2} \, (X, {\rm spt} \, \|\G_{\#} \, {\mathbf C}^{\prime}\|) \, d\|V\|(X)  \leq \overline{\nu}\th^{2}E_{V}^{2};\nonumber\\
&&({\rm d}) \;\;\;\; \left(\th^{-n-2}\int_{{\mathbf R} \times B_{\th}} {\rm dist}^{2}\, (X, P)\, d\|\G^{-1}_{\#} \,V\|(X)\right)^{1/2}\nonumber\\
&&\hspace{1in} \geq 2^{-\frac{n+4}{2}}\sqrt{\overline{C}_{1}} \,{\rm dist}_{\mathcal H}\, ({\rm spt} \, \|{\mathbf C}\|\cap ({\mathbf R} \times B_{1}), P \cap ({\mathbf R} \times B_{1})) - \overline{C}_{2}E_{V}\nonumber\\ 
&&\;\;\;\;\mbox{for any $P \in G_{n}$ of the form $P = \{x^{1} = \lambda x^{2}\}$ for some $\lambda \in (-1, 1)$};\nonumber\\ 
&&({\rm e})\;\;\;\; \{Z \, : \, \Theta \, (\|\G_{\#}^{-1} \, V\|, Z) \geq q\} \cap \left({\mathbf R} \times (B_{\th/2} \cap \{|x^{2}| < \th/16\})\right) = \emptyset;\nonumber\\
&&({\rm f})\;\;\;\; \left(\omega_{n}\th^{n}\right)^{-1}\|\G_{\#}^{-1} \, V\|({\mathbf R} \times B_{\th}) < q + 1/2. 
\end{eqnarray*}
Here the constants $\overline{\k}, \overline{C}_{0}, \overline{\nu}, \overline{C}_{2} \in (0, \infty)$ each depends only on $n$, $q$, $\a$ and $\overline{C}_{1} = \overline{C}_{1}(n) = \int_{B_{1/2} \cap \{x^{2} > 1/16\}}|x^{2}|^{2} \, d{\mathcal H}^{n}(x^{2}, y)$. 
\end{lemma}

\begin{proof}
Consider any sequence of varifolds $\{V_{k}\} \subset {\mathcal S}_{\a}$ and any sequence of cones $\{{\mathbf C}_{k}\} \subset {\mathcal C}_{q}$ satisfying, for each $k=1, 2, \ldots,$ 
hypotheses $(1_{k})-(7_{k})$ of Section~\ref{fineexcessblowup} for some sequences $\{\e_{k}\}$, $\{\g_{k}\}$, $\{\b_{k}\}$ of numbers with $\e_{k}, \g_{k}, \b_{k} \to 0^{+}$ and with $M_{0}$ in place of $M_{0}^{3}$ (in hypothesis $(6_{k})$). The lemma will be established by showing that for each of infinitely many $k$, there exist an orthogonal rotation $\G_{k}$ of ${\mathbf R}^{n+1}$ 
and a cone ${\mathbf C}^{\prime}_{k} \in {\mathcal C}_{q}$ such that the conclusions of the lemma hold with 
$V_{k}$, ${\mathbf C}_{k}$, ${\mathbf C}_{k}^{\prime}$  in place of $V$, ${\mathbf C}$, ${\mathbf C}^{\prime}$ respectively, for fixed constants $\overline{\k}, \overline{C}_{0}, \overline{\g}_{0}, \overline{\nu}, \overline{C}_{1}, \overline{C}_{2} \in (0, \infty)$ depending only on $n$, $q$ and $\a$.     

Let ${\hat E}_{k} =  {\hat E}_{V_{k}}$ and $E_{k} = E_{V_{k}}.$ For $i=1, 2, \ldots,( n-1),$ let $Y_{i} = \frac{1}{2}\th \, e_{2+i} \in \{0\} \times {\mathbf R}^{n-1}.$ We infer from (\ref{excess-1}) that passing to a subsequence of $\{k\}$ without changing notation, for each $k = 1, 2, 3, \ldots,$ there exist points 
$Z_{i, \, k} = (\z_{1}^{i, \, k}, \z_{2}^{i, \, k}, \eta_{i, \, k}) \in {\rm spt} \, \|V_{k}\| \cap ({\mathbf R} \times B_{1}),$ $i=1, 2, \ldots, (n-1),$ such that $\Theta \, (\|V_{k}\|, Z_{i, \, k}) \geq q$ and $|Z_{i, \, k} - Y_{i}| \to 0$ as $k \to \infty;$ also, we may find orthogonal rotations $\G^{\prime}_{k}$ of ${\mathbf R}^{n+1}$ such that 
$$\G^{\prime}_{k}(\Sigma_{k}) = \{0\} \times {\mathbf R}^{n-1} \;\;{\rm and} \;\; \G_{k}^{\prime}\left(\frac{Z_{i,k}}{|Z_{i, k}|}\right) \to e_{2+i}
 \;\;\;\mbox{for each $i=1, 2, \ldots, (n-1)$},$$
where $\Sigma_{k}$ is the $(n-1)$-dimensional subspace spanned by $\{Z_{i, \, k} \, : \, i=1, 2, \ldots, (n-1)\}.$ 
Let $\G^{\prime\prime}_{k}$ be the orthogonal rotation of ${\mathbf R}^{n+1}$ such that $\G^{\prime\prime}_{k}(Y) = Y$
for each $Y \in \{0\}\times {\mathbf R}^{n-1}$ and $\G^{\prime\prime}_{k}\left(\frac{\pi_{12} \, \G_{k}^{\prime}(e_{1})}{|\pi_{12} \, \G_{k}^{\prime}(e_{1})|}\right) = e_{1}$, where 
$\pi_{12} \, : \, {\mathbf R}^{n+1} \to {\mathbf R}^{2} \times \{0\}$ is the orthogonal projection onto the $x^{1}x^{2}$-plane, and let 
\begin{equation*}\label{improvement-1-1}
\G_{k} = \G^{\prime\prime}_{k} \circ \G^{\prime}_{k}, \;\;\mbox{so that}
\end{equation*}
\begin{equation}\label{improvement-4}
\G_{k}(\Sigma_{k}) = \{0\}\times {\mathbf R}^{n-1}, \;\;\; \G_{k}\left(\frac{Z_{i, k}}{|Z_{i,k}|}\right) \to e_{2+i} \;\;\; \mbox{for each $i=1, 2, \ldots, (n-1)$}.
\end{equation}

Let $(\varphi, \psi) \in {\mathcal B}^{F}$ be the fine blow-up of a subsequence of $\{V_{k}\}$ relative to the corresponding subsequence of 
$\{{\mathbf C}_{k}\}.$ Since $\Theta \, (\|V_{k}\|, 0) \geq q$, it follows from (\ref{excess-5}) that $\varphi(0) = \psi(0) = 0$, and consequently, from (\ref{continuity-1-0}) and (\ref{c1alpha-28}) that after passing to further subsequences without changing notation, 
\begin{equation}\label{improvement-1-1-0}
|\z_{1}^{i, \, k}| + {\hat E}_{k}|\z_{2}^{i, \, k}| \leq C\th E_{k}
\end{equation}
for each $i = 1, 2, \ldots, n-1$ and $k=1, 2, \ldots,$ where $C = C(n, q,\a) \in (0, \infty);$ with the help of (\ref{improvement-1-1-0}), the following can then be verified:
\begin{equation}\label{improvement-5}
|e_{1} - \G_{k}(e_{1})| \leq CE_{k} \;\;\; {\rm and} \;\;\; |e_{j} - \G_{k}(e_{j})| \leq C{\hat E}_{k}^{-1}E_{k}, \;\;\; j=2, 3, \ldots, n+1,
\end{equation}
where $C = C(n, q, \a) \in (0, \infty).$ Note in particular that $C$ here is independent of $\th.$ Consequently, 
letting $\widetilde{V}_{k} = \eta_{0, 7/8 \, \#} (\G_{k \, \#} \, V_{k})$ and passing to a further subsequence without changing notation, we have for each $k=1, 2, 3, \ldots$ that
\begin{equation}\label{improvement-6}
d_{\mathcal H} \, (\G_{k}^{-1}(\{0\} \times {\mathbf R}^{n}) \cap ({\mathbf R} \times B_{1}), \{0\} \times B_{1}) \leq CE_{k} \;\;\; {\rm and}
\end{equation}
\begin{equation}\label{improvement-7}
E_{\widetilde{V}_{k}}^{2}  \equiv \int_{{\mathbf R} \times B_{1}} {\rm dist}^{2} \, (X, {\rm spt} \, \|{\mathbf C}_{k}\|) \, d\|\widetilde{V}_{k}\|(X) \leq CE_{k}^{2} 
\end{equation}
where $C = C(n, q, \a) \in (0, \infty).$ Furthermore, we claim that 
\begin{equation}\label{improvement-8}
\widetilde{C}{\hat E}_{k} \leq {\hat E}_{\widetilde{V}_{k}} \leq C{\hat E}_{k}
\end{equation}
for constants $\widetilde{C} = \widetilde{C}(n, q, \a) \in (0, \infty)$ and $C = C(n, q, \a) \in (0, \infty).$ The second of these inequalities follows directly from the definition of $\widetilde{V}_{k}$ and inequality (\ref{improvement-6}); to see the first, note first that since the coarse blow-up $v_{\star}$ of $\{V_{k}\}$ (by the excess ${\hat E}_{k}$) is homogeneous of degree 1 (in fact its graph is a union of half-hyperplanes meeting along $\{0\} \times {\mathbf R}^{n-1}$) and satisfies, by (\ref{excess-5-1}), $\int_{B_{1}}|v_{\star}|^{2} \geq \widetilde{c}$ where $\widetilde{c} = \widetilde{c}(n, q) \in (0, 1),$ we have for each $\s \in (0, 1)$ that $\s^{-n-2}\int_{B_{\s}}|v_{\star}|^{2} = \int_{B_{1}}|v_{\star}|^{2} \geq \widetilde{c}$, so that 
\begin{eqnarray*}
\int_{{\mathbf R} \times B_{\s}} |x^{1}|^{2} \, d\|V_{k}\|(X) = \sum_{j=1}^{q}\int_{B_{\s}} \sqrt{1 + |Du_{k}^{j}|^{2}}|u_{k}^{j}|^{2}  - \sum_{j=1}^{q}\int_{\Sigma_{k}} \sqrt{1 + |Du_{k}^{j}|^{2}}|u_{k}^{j}|^{2}&&\nonumber\\
&&\hspace{-1in} + \int_{{\mathbf R} \times \Sigma_{k}} |x^{1}|^{2} \, d\|V_{k}\|(X)\nonumber\\
&&\hspace{-5in}\geq \int_{B_{\s}}|u_{k}|^{2} - 2C_{\s}\left(\sup_{B_{\s}} \,|u_{k}|^{2}\right){\hat E}_{k}^{2} \geq \left(\frac{1}{2}\widetilde{c} \, \s^{n+2} - 2C_{\s}\left(\sup_{B_{\s}} \,|u_{k}|^{2}\right)\right){\hat E}_{k}^{2} 
\end{eqnarray*}
for sufficiently large $k$, where $u_{k}$, $\Sigma_{k}$ correspond to $u$, $\Sigma$ of Theorem~\ref{flat-varifolds} taken with $V_{k}$ in place of $V$ and the constant $C_{\s}$ is the same as the constant $C$ of Theorem~\ref{flat-varifolds}(a). Thus for sufficiently large $k$ depending on $\s$,  
\begin{equation}\label{improvement-8-1}
\int_{{\mathbf R} \times B_{\s}}|x^{1}|^{2} \, d\|V_{k}\|(X) \geq c{\hat E}_{k}^{2}
\end{equation}
where $c = c(n, q, \s) \in (0, 1)$, which, taken with a suitable choice of $\s \in (0, 1)$,  readily implies the first of the inequalities of (\ref{improvement-8}). 
 
We can now verify using Theorem~\ref{flat-varifolds}, (\ref{excess-4}) and inequalities (\ref{improvement-5})-(\ref{improvement-8}) that after passing to another subsequence without changing notation, for each $k=1, 2, \ldots,$ the hypotheses (1$_{k}$)-(7$_{k}$) of Section~\ref{fineexcessblowup} are satisfied with $\widetilde{V}_{k}$ in place of $V_{k},$ suitable numbers $\widetilde{\e}_{k}$, $\widetilde{\g}_{k}$, $\widetilde{\b}_{k} \to 0^{+}$ in place of $\e_{k}$, $\g_{k}$, $\b_{k}$ respectively and with $M_{0}^{2}$ in place of $M_{0}^{3}$ (in (6$_{k}$)); of these, verification of (1$_{k}$)-(5$_{k}$) is straightforward; to verify that (6$_{k}$) is satisfied with 
$\widetilde{V}_{k}$ in place of $V$ and $M_{0}^{2}$ in place of $M_{0}^{3}$, we proceed as follows: We note first that by (\ref{improvement-8}), 
\begin{equation*}
\inf_{P = \{x^{1} = \lambda x^{2}\}}\int_{{\mathbf R} \times B_{1}} {\rm dist}^{2} \, (X, P) \, d\|\widetilde{V}_{k}\|(X) =
\inf_{P = \{x^{1} = \lambda x^{2}\}; |\lambda| \leq C{\hat E}_{k}}\int_{{\mathbf R} \times B_{1}} {\rm dist}^{2} \, (X, P) \, d\|\widetilde{V}_{k}\|(X)
\end{equation*}
where $C = C(n, q, \a) \in (0, \infty)$, and that for any hyperplane 
$P = \{x^{1} = \lambda x^{2}\}$ with $|\lambda| <C{\hat E}_{k}$ and for sufficiently large $k$,
\begin{eqnarray*}
\int_{{\mathbf R} \times B_{1}} {\rm dist}^{2} \, (X, P) \, d\|\widetilde{V}_{k}\|(X) \geq \left(\frac{8}{7}\right)^{n+2} \int_{{\mathbf R} \times B_{1/2}} {\rm dist}^{2} \, (X, \G_{k}^{-1}(P)) \, d\|V_{k}\|(X)\nonumber\\
&&\hspace{-5in}\geq \frac{1}{2}\left(\frac{8}{7}\right)^{n+2}\int_{{\mathbf R} \times B_{1/2}}{\rm dist}^{2} \, (X, P) \, d\|V_{k}\|(X) - C\, {\rm dist}_{\mathcal H}^{2}(\G_{k}^{-1}(P) \cap ({\mathbf R} \times B_{1/2}), P \cap ({\mathbf R} \times B_{1/2}))\nonumber\\
&&\hspace{-5in}\geq 7^{-n-2}2^{n-1}\omega_{n}^{-1}(2q+1)^{-1}\overline{C}_{1} \int_{{\mathbf R} \times B_{1}}{\rm dist}^{2} \, (X, P) \, d\|V_{k}\|(X) - CE_{k}^{2}\nonumber\\
&&\hspace{-5in}\geq 7^{-n-2}2^{n-1}\omega_{n}^{-1}(2q+1)^{-1}\overline{C}_{1}\left(\frac{3}{2}M_{0}\right)^{-1} \int_{{\mathbf R} \times B_{1}} |x^{1}|^{2}\, d\|V_{k}\|(X) - CE_{k}^{2}\nonumber\\
\end{eqnarray*}
where $C = C(n, q) \in (0, \infty)$, the third inequality follows from  (\ref{L2-est-2-2}) with $\r = 1/2$ and $Z = 0$, and the last inequality holds 
by hypothesis of the present lemma. On the other hand, 
\begin{eqnarray*}
&&\int_{{\mathbf R} \times B_{1}} |x^{1}|^{2} \, d\|\widetilde{V}_{k}\|(X) \leq 2\left(\frac{8}{7}\right)^{n+2} \int_{{\mathbf R} \times B_{1}} |x^{1}|^{2} \, d\|V_{k}\|(X) +\nonumber\\
&&\hspace{1.8in}  \left(\frac{8}{7}\right)^{n+2} \omega_{n}(2q+1) {\rm dist}^{2} \, (\G_{k}^{-1}(\{0\} \times {\mathbf R}^{n}) \cap ({\mathbf R} \times B_{1}), \{0\} \times B_{1})\nonumber\\
&&\hspace{1.5in} \leq 2\left(\frac{8}{7}\right)^{n+2} \int_{{\mathbf R} \times B_{1}} |x^{1}|^{2} \, d\|V_{k}\|(X) + CE_{k}^{2}
\end{eqnarray*}
where $C = C(n, q, \a) \in (0, \infty).$ Hence 
\begin{equation*}
{\hat E}_{\widetilde{V}_{k}}^{2} \leq \frac{3M_{0}}{2\left(2^{-2n-7}\omega_{n}^{-1}(2q+1)^{-1}\overline{C}_{1} - C\g_{k}\right)}\int_{{\mathbf R} \times B_{1}}{\rm dist}^{2} \, (X, P) \, d\|\widetilde{V}_{k}\|(X)
\end{equation*}
where $C = C(n, q, \a) \in (0, \infty)$, and it follows from this that hypothesis (6$_{k}$) with $\widetilde{V}_{k}$ in place of $V_{k}$ and $M_{0}^{2}$ in place of $M_{0}^{3}$ is satisfied for all sufficiently large $k$; hypothesis (7$_{k}$) with $\widetilde{V}_{k}$ in place of $V_{k}$ can easily be verified using the estimate $Q^{\star}_{\widetilde{V}_{k}}(p_{k} - 1) \geq CQ^{\star}_{V_{k}}(p_{k}-1)$, where $C = C(n, q) \in (0, \infty)$, which  follows from (\ref{improvement-5}) and the fact that, for any ${\mathbf C} \in \cup_{j=4}^{p_{k}-1}{\mathcal C}_{q}(j)$,  
\begin{eqnarray*}
&&\int_{{\mathbf R} \times \left(B_{7/16} \setminus \{|x^{2}| \leq 7/(8 \cdot 16)\}\right)} {\rm dist}^{2} \, (X, {\rm spt} \, \|V_{k}\|) \, d\|{\mathbf C}\|(X)\nonumber\\
&&\hspace{1in} + \int_{{\mathbf R} \times B_{7/8}} {\rm dist}^{2} \, (X, {\rm spt} \, \|{\mathbf C}\|) \, d\|V_{k}\|(X) \geq \widetilde{c}_{1}\left(Q^{\star}_{V_{k}}(p_{k}-1)\right)^{2} 
\end{eqnarray*}
where $\widetilde{c}_{1} = \widetilde{c}_{1}(n, q) \in (0, 1),$ the validity of which can be seen by reasoning as in the proof of (\ref{improvement-8-1}), using the fact that the blow-up 
of $\{V_{k}\}$ by $Q^{\star}_{V_{k}}(p_{k}-1)$ is homogeneous of degree 1 (by hypothesis~(7$_{k}$)) and has, by (\ref{separation}), $L^{2}(B_{1})$ norm $\geq c$, $c = c(n, q) \in (0, 1).$ 

Thus, the fine blow-up $(\widetilde{\varphi}, \widetilde{\psi})$ of 
$\{\widetilde{V}_{k}\}$  relative to $\{\mathbf{C}_{k}\}$ belongs to ${\mathcal B}^{F}.$ Furthermore, it follows from (\ref{improvement-4}) and (\ref{excess-5}) (applied with $\widetilde{V}_{k}$ in place of $V_{k}$ and 
$\frac{8}{7}\G_{k} \, Z_{i, k}$, $i=1, 2, \ldots, n-1,$ in place of $Z$) that for each $i=1, 2, \ldots, (n-1)$, 
$\widetilde{\varphi}(Y_{i}) = \widetilde{\psi}(Y_{i}) = 0$ and consequently, since $Y_{i} = \frac{1}{2}\th \, e_{i+2},$  that 
there exist points $S_{j, i}, T_{j, i} \in B_{\th/2} \cap (\{0\} \times {\mathbf R}^{n-1})$ such that 
\begin{equation*}
\frac{\partial \, \widetilde{\varphi}_{j}}{\partial \, y^{i}}(S_{j,i}) = 0 \;\;\; \mbox{and} \;\;\; \frac{\partial \, \widetilde{\psi}_{j}}{\partial \, y^{i}}(T_{j, i}) = 0
\end{equation*}
for each $i=1, 2, \ldots, n-1$ and $j=1, 2, \ldots, q.$ By the estimate of Theorem~\ref{c1alpha}, this readily implies that
\begin{equation}\label{improvement-10}
|D_{y}\widetilde{\varphi}(0)|^{2} + |D_{y}\widetilde{\psi}(0)|^{2} \leq C\th^{2}\left(\int_{B_{1/2} \cap \{x^{2} < 0\}}|\widetilde{\varphi}|^{2} + \int_{B_{1/2} \cap \{x^{2} > 0\}} |\widetilde{\psi}|^{2}\right)
\end{equation}
where  $C = C(n, q, \a) \in (0, \infty).$ Letting, for $j=1, 2, \ldots, q$ and $x  = (x^{2}, y) \in {\mathbf R}^{n}$,
$$L_{\widetilde{\varphi}}^{j}(x) = D\widetilde{\varphi}_{j}(0) \cdot x, \;\;\; L_{\widetilde{\psi}}^{j}(x) = D\widetilde{\psi}_{j}(0) \cdot x, \;\;\; 
P_{\widetilde{\varphi}}^{j}(x^{2}, y) = \frac{\partial \, \widetilde{\varphi}_{j}}{\partial \, x^{2}}(0)x^{2} \;\;\; {\rm and} \;\;\; 
P_{\widetilde{\psi}}^{j}(x) = \frac{\partial \, \widetilde{\psi}_{j}}{\partial \, x^{2}}(0)x^{2},$$ 
it follows from (\ref{improvement-10}) that for each $(x^{2}, y) \in {\mathbf R}^{n}$, 

\begin{eqnarray*}\label{improvement-11}
|P_{\widetilde{\varphi}}^{j}(x^{2}, y) - L_{\widetilde{\varphi}}^{j}(x^{2}, y)|^{2} + |P_{\widetilde{\psi}}^{j}(x^{2}, y) - L_{\widetilde{\psi}}^{j}(x^{2}, y)|^{2}&&\nonumber\\
&&\hspace{-2in}\leq C\th^{2}|y|^{2} \left(\int_{B_{1/2} \cap \{x^{2} < 0\}}|\widetilde{\varphi}|^{2} + \int_{B_{1/2} \cap \{x^{2} > 0\}} |\widetilde{\psi}|^{2}\right)
\end{eqnarray*}
and consequently from Theorem~\ref{c1alpha}  that 
\begin{equation}\label{improvement-13}
\th^{-n-2}\left(\int_{B_{2\th} \cap \{x^{2} \leq 0\}} |\widetilde{\varphi} - P_{\widetilde{\varphi}}|^{2} + \int_{B_{2\th} \cap \{x^{2} \geq 0\}} |\widetilde{\psi} - 
P_{\widetilde{\psi}}|^{2} \right) \leq C\th^{2}, \;\; C = C(n, q, \a) \in (0, \infty).
\end{equation}

For $j=1, 2, \ldots, q$ and $k=1, 2, \ldots,$ let 
\begin{equation}\label{improvement-14}
\lambda_{j}^{\prime \, k} = \lambda_{j}^{k} + E_{\widetilde{V}_{k}}\frac{\partial \, \widetilde{\varphi}_{j}}{\partial \, 
x^{2}}(0), \;\;\; \mu_{j}^{\prime \, k} = \mu_{j}^{k} + E_{\widetilde{V}_{k}}\frac{\partial \, \widetilde{\psi}_{j}}{\partial \, x^{2}}(0),
\end{equation}
$$H_{j}^{\prime \, k} = \{(x^{1}, x^{2} , y) \, : \, x^{1} = \lambda_{j}^{\prime \, k}x^{2}, \;\; x^{2} \leq 0\}, \;\; 
G_{j}^{\prime \, k} = \{(x^{1}, x^{2} , y) \, : \, x^{1} = \mu_{j}^{\prime \, k}x^{2}, \;\; x^{2} \geq 0\} \;\; {\rm and}$$ 
$${\mathbf C}_{k}^{\prime}  =\sum_{j=1}^{q} |H_{j}^{\prime \, k}| + |G_{j}^{\prime \, k}|.$$
With the help of (\ref{improvement-4}), (\ref{excess-2}), (\ref{improvement-7}),  it is straightforward to verified that
\begin{equation}\label{improvement-15}
\th^{-n-2}\int_{\G_{k}^{-1}({\mathbf R} \times B_{\th})} {\rm dist}^{2} \, (X, {\rm spt} \, \|(\G_{k}^{-1})_{\#} \, {\mathbf C}_{k}^{\prime}\|) 
\, d\|V_{k}\|(X) \leq C\th^{2} E_{k}^{2}
\end{equation}
for all sufficiently large $k$, where ${\mathbf C}_{k}^{\prime}$ is as above and $\G_{k}$ is as in (\ref{improvement-4}), and $C = C(n, q, \a) \in (0, \infty).$ Furthermore, it follows from (\ref{improvement-7}), (\ref{improvement-14}) and  Theorem~\ref{c1alpha} that
\begin{equation}\label{improvement-16}
{\rm dist}_{\mathcal H}^{2} \, ({\rm spt} \, \|{\mathbf C}_{k}^{\prime}\| \cap ({\mathbf R} \times B_{1}), {\rm spt} \, \|{\mathbf C}_{k}\| \cap ({\mathbf R} \times B_{1})) \leq CE_{k}^{2}, \;\;C = C(n, q, \a) \in (0, \infty),
\end{equation}
and from (\ref{excess-3}) (applied with $\widetilde{V}_{k}$ in place of $V_{k}$) that 
\begin{eqnarray}\label{improvement-17}
\th^{-n-2}\int_{\G_{k}^{-1}\left({\mathbf R} \times \left(B_{\th/2} \setminus \{|x^{2}| \leq \th/16\}\right)\right)} {\rm dist}^{2} \, (X, {\rm spt} \, \|V_{k}\|) \, d\|\G^{-1}_{k \, \#} \, {\mathbf C}_{k}^{\prime}\|(X)&&\nonumber\\
&&\hspace{-4in} \leq C\th^{-n-2}\int_{\G_{k}^{-1}\left({\mathbf R} \times B_{\th}\right)} {\rm dist}^{2} \, (X, {\rm spt} \, \|\G^{-1}_{k \, \#} \, {\mathbf C}_{k}^{\prime}\|) \, d\|V_{k}\|(X), \;\; C = C(n, q, \a) \in (0, \infty). 
\end{eqnarray}
Again by (\ref{excess-3}) (applied with $\widetilde{V}_{k}$ in place of $V_{k}$), (\ref{improvement-7}),  (\ref{improvement-15}) and (\ref{improvement-16}), we have that for any hyperplane $P$ of the form $P = \{x^{1} = \lambda x^{2}\}$, $|\lambda| < 1$, writing $\widetilde{\th} = \frac{8}{7}\th$, 
\begin{eqnarray*}
\widetilde{\th}^{-n-2}\int_{{\mathbf R} \times B_{\widetilde{\th}}} {\rm dist}^{2}\,(X, P) \, d\|\widetilde{V}_{k}\|(X)&&\nonumber\\ 
&&\hspace{-2.25in}\geq \frac{1}{2}\widetilde{\th}^{-n-2}\sum_{j=1}^{q}\left(\int_{B_{\widetilde{\th}/2}\cap \{x^{2} < -\widetilde{\th}/16\}} |h_{j}^{k}  - \lambda x^{2}+ \widetilde{u}_{j}^{k}|^{2}  + \int_{B_{\widetilde{\th}/2} \cap  
\{x^{2} > \widetilde{\th}/16\}} |g_{j}^{k} - \lambda x^{2}  + \widetilde{w}_{j}^{k}|^{2}\right) \nonumber\\
&&\hspace{-2.25in}\geq \frac{1}{4}\widetilde{\th}^{-n-2}\sum_{j=1}^{q}\left(\int_{B_{\widetilde{\th}/2} \cap \{x^{2} < -\widetilde{\th}/16\}} |h_{j}^{k} - \lambda x^{2}|^{2} + \int_{B_{\widetilde{\th}/2} \cap \{x^{2} > \widetilde{\th}/16\}}|g_{j}^{k} - \lambda x^{2}|^{2}\right)\nonumber\\ 
&&\hspace{-1in}-\frac{1}{2}\widetilde{\th}^{-n-2}\sum_{j=1}^{q}\left(\int_{B_{\widetilde{\th}/2}\cap \{x^{2} < -\widetilde{\th}/16\}} |\widetilde{u}_{j}^{k}|^{2}  + \int_{B_{\widetilde{\th}/2} \cap \{x^{2} > \widetilde{\th}/16\}} |\widetilde{w}_{j}^{k}|^{2}\right)\nonumber\\
&&\hspace{-2.25in} \geq 2^{-n-4}\overline{C}_{1}{\rm dist}_{\mathcal H}^{2} \, ({\rm spt} \, \|{\mathbf C}_{k} \| \cap ({\mathbf R} \times B_{1}), P \cap ({\mathbf R} \times  B_{1}))\nonumber\\
&&- \,\frac{1}{2}\widetilde{\th}^{-n-2}\int_{{\mathbf R} \times B_{\widetilde{\th}}}{\rm dist}^{2} \, (X, {\rm spt} \, \|{\mathbf C}_{k}\|) \, d\|\widetilde{V}_{k}\|(X)\nonumber\\
&&\hspace{-2.25in} \geq 2^{-n-4}\overline{C}_{1}{\rm dist}_{\mathcal H}^{2} \, ({\rm spt} \, \|{\mathbf C}_{k} \| \cap ({\mathbf R} \times B_{1}), P \cap ({\mathbf R} \times B_{1}))\nonumber\\
&&- \,\widetilde{\th}^{-n-2}\int_{{\mathbf R} \times B_{\widetilde{\th}}}{\rm dist}^{2} \, (X, {\rm spt} \, \|{\mathbf C}_{k}^{\prime}\|) \, d\|\widetilde{V}_{k}\|(X) - CE_{k}^{2}\nonumber\\
&&\hspace{-2.25in} \geq 2^{-n-4}\overline{C}_{1}{\rm dist}_{\mathcal H}^{2} \, ({\rm spt} \, \|{\mathbf C}_{k} \| \cap ({\mathbf R} \times B_{1}), P \cap ({\mathbf R} \times B_{1})) - CE_{k}^{2}\nonumber\\
\end{eqnarray*}
where $\overline{C}_{1} = \int_{B_{1/2} \cap \{x^{2} > 1/16\}}|x^{2}|^{2} \, d{\mathcal H}^{n}(x^{2}, y)$, $C_{2} = C_{2}(n, q, \a) \in (0, \infty)$ and the notation is as in Theorem~\ref{L2-est-1} taken with $\widetilde{V}_{k}$ in place of $V$ (in particular with $\widetilde{u}^{j}_{k}$, $\widetilde{w}^{j}_{k}$ corresponding to $u^{j}$, $w^{j}$). This readily implies  that
\begin{eqnarray}\label{improvement-18}
&&\left(\th^{-n-2}\int_{{\mathbf R} \times B_{\th}} {\rm dist}^{2} \, (X, P) \, d\|\G_{k \, \#} \,V_{k}\|(X)\right)^{1/2}\nonumber\\
&&\hspace{1in} \geq \, 2^{-\frac{n+4}{2}}\sqrt{\overline{C}_{1}} \,{\rm dist}_{\mathcal H} \, ({\rm spt} \, \|{\mathbf C}_{k}\|\cap ({\mathbf R} \times B_{1}), P \cap ({\mathbf R} \times \times B_{1})) - CE_{k}
\end{eqnarray}
for each hyperplane $P = \{x^{1} = \lambda x^{2}\}$ with $|\lambda| < 1$ and all sufficiently large $k$, where $C= C(n, q, \a) \in (0, \infty).$

The inequalities (\ref{improvement-5}) and (\ref{improvement-15})-(\ref{improvement-18}) say that the conclusions (a)-(d) of the lemma, with $V_{k}$, ${\mathbf C}_{k}$,
${\mathbf C}_{k}^{\prime}$, $\G_{k}^{-1}$ in place of $V$, ${\mathbf C}$, 
${\mathbf C}^{\prime}$, $\G$, hold for all sufficiently large $k.$ Conclusion (e) with $V_{k}$ in place of $V$ and 
$\G_{k}^{-1}$ in place of $\G$
is clear, for all sufficiently large $k$, by (\ref{excess-4}) applied with $\widetilde{V}_{k}$ in place of $V_{k}$. Conclusion (f) with $V_{k}$ in place of $V$ and $\G_{k}^{-1}$ in place of $\G$ follows, for sufficiently large $k$, from the Constancy Theorem for stationary integral varifolds and the fact that $q \leq \Theta \, (\|\widetilde{V}_{k}\|, 0) \leq \left(\omega_{n}2^{n}\right)^{-1}\|\widetilde{V_{k}}\|(B_{2}^{n+1}(0)) < q + 1/2$ for each $k$.
\end{proof}

\begin{lemma}\label{multi-scale} 
Let $q \geq 2$ be an integer, $\a \in (0, 1)$  and $p \in \{4, 5, \ldots, 2q\}.$  
For $j=1, 2, \ldots, p-3,$ let $\th_{j} \in (0, 1/4)$ be such that $\th_{1} > 8\th_{2} > 64\th_{3} > \ldots >8^{p - 4}\th_{p-3}$. 
There exist numbers $\e^{(p)} = \e^{(p)}(n, q, \a, \th_{1}, \th_{2}, \ldots, \th_{p-3}) \in (0, 1/2)$, $\g^{(p)}= 
\g^{(p)}(n, q, \a, \th_{1}, \th_{2} \, \ldots, \th_{p-3}) \in (0, 1/2)$ such that if $V \in {\mathcal S}_{\a},$ ${\mathbf C} \in {\mathcal C}_{q}(p)$ satisfy Hypotheses~\ref{hyp}, Hypothesis~($\star$) (of Section~\ref{fineblowup}) with $\e = \e^{(p)}$, $\g = \g^{(p)}$, $M = \frac{3}{2}M_{0}$ and if the induction hypotheses $(H1)$, $(H2)$ hold,  
then there exist an orthogonal rotation $\G$ of ${\mathbf R}^{n+1}$ and a cone ${\mathbf C}^{\prime} \in {\mathcal C}_{q}$ 
such that, with 
$${\hat E}_{V}^{2} = \int_{{\mathbf R} \times B_{1}}|x^{1}|^{2} \, d\|V\|(X) \;\;\; {\rm and}$$ 
$$Q_{V}^{2}({\mathbf C}) = \int_{{\mathbf R} \times (B_{1/2} \setminus \{|x^{2}| < 1/16\})} {\rm dist}^{2}(X, {\rm spt} \, \|V\|) \,
d\|{\mathbf C}\|(X) + \int_{{\mathbf R} \times B_{1}} {\rm dist}^{2} \, (X, {\rm spt} \, \|{\mathbf C})\|) \, d\|V\|(X),$$ 
we have the following:
\begin{eqnarray*}
&&({\rm a})\;\;\;\; |e_{1} - \G(e_{1})| \leq \k^{(p)}Q_{V}({\mathbf C}) \;\;{\rm and} \;\; |e_{j} - \G(e_{j})| \leq \k^{(p)}{\hat E}_{V}^{-1}Q_{V}({\mathbf C}) \;\;\mbox{for each} \;\; j=2, 3, \ldots, n+1;\nonumber\\ 
&&({\rm b})\;\;\;\; {\rm dist}_{\mathcal H}^{2} \, ({\rm spt} \, \|{\mathbf C}^{\prime}\| \cap ({\mathbf R} \times B_{1}), {\rm spt} \, \|{\mathbf C}\| \cap ({\mathbf R} \times B_{1})) \leq C_{0}^{(p)}Q_{V}^{2}({\mathbf C});\nonumber\\
&&\hspace{-.25in}\mbox{and for some $j \in \{1, 2, \ldots, p-3\}$,}\nonumber\\ 
&&({\rm c})\;\;\;\;
\th_{j}^{-n-2}\int_{\G\left({\mathbf R} \times \left(B_{\th_{j}/2} \setminus \{|x^{2}| \leq \th_{j}/16\}\right)\right)} {\rm dist}^{2} \, (X, {\rm spt} \, \|V\|) \, d\|\G_{\#} \, {\mathbf C}^{\prime}\|(X)\nonumber\\
&&\hspace{1.5in} + \; \th_{j}^{-n-2}\int_{\G({\mathbf R} \times B_{\th_{j}})} {\rm dist}^{2} \, (X, {\rm spt} \, \|\G_{\#} \, {\mathbf C}^{\prime}\|) \, d\|V\|(X) \leq \nu_{j}^{(p)}\th_{j}^{2}Q_{V}^{2}({\mathbf C});\nonumber\\
&&({\rm d})\;\;\;\;\left(\th_{j}^{-n-2}\int_{{\mathbf R} \times B_{\th_{j}}} {\rm dist}^{2}\,(X, P) \, d\|\G^{-1}_{\#} \,V\|(X)\right)^{1/2}\nonumber\\ 
&&\hspace{1in}\geq 2^{-\frac{n+4}{2}}\sqrt{\overline{C}_{1}} \,{\rm dist}_{\mathcal H}\, ({\rm spt} \, \|{\mathbf C}\| \cap ({\mathbf R} \times B_{1}), P \cap ({\mathbf R} \times B_{1})) - C_{2}^{(p)}Q_{V}({\mathbf C})\nonumber\\
&&\;\;\;\;\mbox{for any $P \in G_{n}$ of the form $P = \{x^{1} = \lambda x^{2}\}$ for some $\lambda \in (-1, 1)$};\nonumber\\
&&({\rm e})\;\;\;\;\{Z \, : \, \Theta \, (\|\G_{\#}^{-1} \, V\|, Z) \geq q\} \cap \left({\mathbf R} \times (B_{\th_{j}/2} \cap \{|x^{2}| < \th_{j}/16\})\right) = \emptyset;\nonumber\\
&&({\rm f})\;\;\;\;\left(\omega_{n}\th_{j}^{n}\right)^{-1}\|\G_{\#}^{-1} \, V\|({\mathbf R} \times B_{\th_{j}}) < q + 1/2.
\end{eqnarray*}
Here the dependence of the various constants on the parameters is as follows: 
$$\k^{(p)} = \k^{(p)}(n, q, \a,\th_{1}, \ldots, \th_{p-4}), \; C_{0}^{(p)} = C_{0}^{(p)}(n, q, \a, \th_{1}, \ldots, \th_{p-4}), \; C_{2}^{(p)} = C_{2}^{(p)}(n, q, \a, \th_{1}, \ldots, \th_{p-4})$$ 
in case $q \geq 3$ and $p \in \{5, 6, \ldots, 2q\};$  $\k^{(4)} = \overline{\k}$, $C_{0}^{(4)} = \overline{C}_{0}$, $C_{2}^{(4)} = \overline{C}_{2},$ where $\overline{\k} = \overline{\k}(n, q, \a)$, $\overline{C}_{0} = \overline{C}_{0}(n, q, \a)$, $\overline{C}_{2} = \overline{C}_{2}(n, q, \a)$ are as  in Lemma~\ref{excess-improvement}; $\nu_{1}^{(p)} = \overline{\nu},$ where $\overline{\nu} = \overline{\nu}(n, q, \a)$ is as in Lemma~\ref{excess-improvement}; and, in case $q \geq 3,$ for each $j=2, 3, \ldots, p-3,$ $\nu_{j}^{(p)} = \nu_{j}^{(p)}(n, q, \a, \th_{1}, \ldots, \th_{j-1}).$ In particular, $\nu_{j}^{(p)}$ is independent of $\th_{j}, \th_{j+1}, \ldots, \th_{p-3}$ for $j=1, 2, \ldots, p-3.$ 
\end{lemma}

\begin{proof}
If $p=4$ then we may simply set $\e^{(4)}(n,q, \a, \th_{1}) = \overline\e(n, q, \a, \th_{1})$ and $\g^{(4)} (n, q, \a, \th_{1})= \overline\g(n, q, \a, \th_{1})$, where $\overline\e,$ $\overline\g$ are as in Lemma~\ref{excess-improvement}, and deduce from Lemma~\ref{excess-improvement} with $\th = \th_{1}$  that  there exist a cone 
${\mathbf C}^{\prime} \in {\mathcal C}_{q}$ and an orthogonal rotation $\G$ of ${\mathbf R}^{n+1}$ such that the conclusions of the lemma hold with $j=1$ in (c)-(f); with $\k^{(4)}=\overline{\k}$, $C_{0}^{(4)} = \overline{C}_{0}$,
$C_{2}^{(4)} = \overline{C}_{2}$ and $\nu_{1}^{(4)} = \overline{\nu}$, where $\overline{\k}_{j}$, $\overline{C}_{0}$, $\overline{C}_{2}$, $\overline{\nu}$ are as in Lemma~\ref{excess-improvement}. Thus the lemma holds if $p=4.$

Else $q \geq 3$ and $p \in \{5, 6, \ldots, 2q\}.$ Assume by induction the validity of  the lemma with any $p^{\prime} \in \{4, 5, \ldots, p-1\}$ in place of $p$.  Let $\th_{j} \in (0, 1/4),$ $j=1, 2, \ldots, p-3$ be given such that $\th_{1} > 8\th_{2} > 64\th_{3} > \ldots >8^{p - 4}\th_{p-3}$. To prove the lemma as stated,  it suffices to show that for arbitrary sequences $\{V_{k}\} \subset {\mathcal S}_{\a}$, $\{{\mathbf C}_{k}\} \subset {\mathcal C}_{q}(p)$ that satisfy hypotheses $(1_{k})-(5_{k})$ of Section~\ref{fineexcessblowup} as well as hypothesis $(6_{k})$ of Section~\ref{fineexcessblowup} with $M_{0}$ in place of $M_{0}^{2},$ there exist 
a subsequence $\{k^{\prime}\}$ of $\{k\}$ and, for each $k^{\prime},$ a cone 
${\mathbf C}_{k^{\prime}}^{\prime} \in {\mathcal C}_{q}$ and an orthogonal rotation $\G_{k^{\prime}}$ of ${\mathbf R}^{n+1}$ such that the conclusions of the lemma hold with $V_{k^{\prime}}$, ${\mathbf C}_{k^{\prime}}$, ${\mathbf C}_{k^{\prime}}^{\prime},$ $\G_{k}^{\prime}$ in place of $V$, ${\mathbf C}$, ${\mathbf C}^{\prime},$ $\G$ respectively, and with 
suitable constants $\k^{(p)}$, $C_{0}^{(p)}$, $C_{2}^{(p)}$  and $\nu_{1}^{(p)}, \ldots, \nu_{p-3}^{(p)}$ depending only on the parameters as specified in the statement of the lemma. So suppose, for $k=1, 2, \ldots,$  $V_{k} \in {\mathcal S}_{\a}$, ${\mathbf C}_{k} \in {\mathcal C}_{q}(p)$ satisfy hypotheses $(1_{k})-(6_{k})$ of Section~\ref{fineexcessblowup} with $M_{0}$ in place of $M_{0}^{2}.$ For each $k$, choose a cone $\widetilde{\mathbf C}_{k} \in \cup_{j=4}^{p-1} \, {\mathcal C}_{q}(j)$ such that 
\begin{eqnarray}\label{multi-1}
\left(\widetilde{Q}_{k}\right)^{2} \equiv \left(\int_{{\mathbf R} \times (B_{1/2} \setminus \{|x^{2}| < 1/16\})} {\rm dist}^{2}(X, {\rm spt} \, \|V_{k}\|) \,d\|\widetilde{\mathbf C}_{k}\|(X)\right.&&\nonumber\\
&&\hspace{-2in} + \left.\int_{{\mathbf R} \times B_{1}} {\rm dist}^{2} \, (X, {\rm spt}\, \|\widetilde{\mathbf C}_{k}\|) \, d\|V_{k}\|(X)\right) \leq \frac{3}{2}\left(Q_{k}^{\star}\right)^{2}
\end{eqnarray}
where 
\begin{eqnarray*}
\left(Q_{k}^{\star}\right)^{2} = \inf_{\widetilde{\mathbf C}\in \cup_{j=4}^{p-1}{\mathcal C}_{q}(j)} \, \left(\int_{{\mathbf R} \times (B_{1/2} \setminus \{|x^{2}| < 1/16\})} {\rm dist}^{2}(X, {\rm spt} \, \|V_{k}\|) \,d\|\widetilde{\mathbf C}\|(X)\right.&&\\
&&\hspace{-2in} + \left.\int_{{\mathbf R} \times B_{1}} {\rm dist}^{2} \, (X, {\rm spt}\, \|\widetilde{\mathbf C}\|) \, d\|V_{k}\|(X)\right).
\end{eqnarray*}

Let $\b = \overline\b(n, q, \a, \th_{1})$ where $\overline\b$ is as in Lemma~\ref{excess-improvement}, and consider the following two alternatives:
\begin{itemize}
\item[$(A)$] for infinitely many $k$, 
\begin{eqnarray*}
\int_{{\mathbf R} \times (B_{1/2} \setminus \{|x^{2}| < 1/16\})} {\rm dist}^{2}(X, {\rm spt} \, \|V_{k}\|) \,d\|{\mathbf C}_{k}\|(X)&&\nonumber\\
&&\hspace{-2in} +\int_{{\mathbf R} \times B_{1}} {\rm dist}^{2} \, (X, {\rm spt}\, \|{\mathbf C}_{k}\|) \, d\|V_{k}\|(X) < \b\left(Q_{k}^{\star}\right)^{2}
\end{eqnarray*}
\item[$(B)$] for all sufficiently large $k$, 
\begin{eqnarray*}
\int_{{\mathbf R} \times (B_{1/2} \setminus \{|x^{2}| < 1/16\})} {\rm dist}^{2}(X, {\rm spt} \, \|V_{k}\|) \,d\|{\mathbf C}_{k}\|(X)&&\nonumber\\
&&\hspace{-2in} +\int_{{\mathbf R} \times B_{1}} {\rm dist}^{2} \, (X, {\rm spt}\, \|{\mathbf C}_{k}\|) \, d\|V_{k}\|(X) \geq \b\left(Q_{k}^{\star}\right)^{2}.
\end{eqnarray*}
\end{itemize}

If alternative $(A)$ holds, we deduce directly from Lemma~\ref{excess-improvement}, applied with $\th = \th_{1},$ that 
for infinitely many $k$, there exist a cone ${\mathbf C}_{k}^{\prime} \in {\mathcal C}_{q}$ and an orthogonal rotation $\G_{k}$ of 
${\mathbf R}^{n+1}$  such that the conclusions of the present lemma 
hold with $V_{k}$, ${\mathbf C}_{k},$  ${\mathbf C}_{k}^{\prime}$, $\G_{k}$ in place of $V$, ${\mathbf C}$, ${\mathbf C}^{\prime}$, $\G;$ with $j=1$ in the conclusions (c)-(f); and with $\overline{\k}$, $\overline{C}_{0}$, $\overline{C}_{2}$, $\overline{\nu}$ (as in Lemma~\ref{excess-improvement}) in place of 
$\k^{(p)}$, $C_{2}^{(p)}$, $\nu_{1}^{(p)}$.

If alternative $(B)$ holds, we have by hypothesis $(5_{k})$ and (\ref{multi-1}) that for all sufficiently large $k$, 
\begin{eqnarray}\label{multi-2}
\left(\int_{{\mathbf R} \times (B_{1/2} \setminus \{|x^{2}| < 1/16\})} {\rm dist}^{2}(X, {\rm spt} \, \|V_{k}\|) \,d\|\widetilde{\mathbf C}_{k}\|(X)\right.&&\nonumber\\
&&\hspace{-2in} + \left.\int_{{\mathbf R} \times B_{1}} {\rm dist}^{2} \, (X, {\rm spt}\, \|\widetilde{\mathbf C}_{k}\|) \, d\|V_{k}\|(X)\right) \leq \frac{3\g_{k}}{2\b} {\hat E_{k}}^{2}.
\end{eqnarray}
Since $\widetilde{\mathbf C}_{k} \in {\mathcal C}_{q}(p^{\prime})$ for some $p^{\prime} \in \{4, 5, \ldots, p-1\}$ and infinitely many $k$, we may, by the induction hypothesis, apply the lemma with $p^{\prime}$ in place of $p$ and 
$\th_{2}, \th_{3}, \ldots, \th_{p^{\prime} -2}$ in place of $\th_{1}, \th_{2}, \ldots, \th_{p^{\prime}-3}$ to deduce that for infinitely many $k$, there exist a cone ${\mathbf C}_{k}^{\prime} \in {\mathcal C}_{q}$ and an orthogonal rotation $\G_{k}$ of ${\mathbf R}^{n+1}$ such that the conclusions 
(a)-(f) hold with $V_{k}$, $\widetilde{\mathbf C}_{k},$  ${\mathbf C}_{k}^{\prime}$, $\G_{k}$ in place of $V$, ${\mathbf C}$, ${\mathbf C}^{\prime}$, $\G$---in particular with $\widetilde{Q}_{k}$ in place of $Q_{V}({\mathbf C})$---and such that:

 \noindent 
(i) in case $p^{\prime} = 4$ (which must be the case if $p = 5$), with $\k^{(4)} = \overline{\k}$, $C_{0}^{(4)} = \overline{C}_{0}$, $C_{2}^{(4)} = \overline{C}_{2}$ and $\nu_{1}^{(4)} = \overline{\nu}$ 
(where $\overline{\k}$, $\overline{C}_{0}$, $\overline{C}_{1}$, $\overline{C}_{2}$, $\overline{\nu}$ are as in Lemma~\ref{excess-improvement} and $C$ is as in Theorem~\ref{L2-est-2}(a)), and

\noindent
(ii) in case $p^{\prime} \in \{5, 6, \ldots, p-1\}$ (possible, of course, only if $p \geq 6$)  with 
$$\k^{(p^{\prime})} = \k^{(p^{\prime})}(n, q, \a, \th_{2}, \ldots, \th_{p^{\prime} - 3}), \;\; C_{0}^{(p^{\prime})} = C_{0}^{(p^{\prime})}(n, q, \a, \th_{2}, \ldots, \th_{p^{\prime} - 3}),$$ 
$$C_{2}^{(p^{\prime})} = C_{2}^{(p^{\prime})}(n, q, \a, \th_{2}, \ldots, \th_{p^{\prime} - 3})$$  
in place of $\k^{(p)}$,  $C_{0}^{(p)}$, $C_{2}^{(p)}$ respectively; with $\nu_{1}^{(p^{\prime})}(n, q, \a) = \overline{\nu}$ (where $\overline{\nu}$ is as in Lemma~\ref{excess-improvement}) in place of $\nu_{1}^{(p)}$; and with $\nu_{j-1}^{(p^{\prime})}(n, q, \a, \th_{2}, \ldots, \th_{j-1})$
in place of $\nu_{j-1}^{(p)}$ for each $j=3, \ldots, p^{\prime} - 2.$ 

Since by (\ref{multi-1}) and the definition of alternative $(B)$ we have that 
\begin{eqnarray*}
\widetilde{Q}_{k}^{2} \leq \frac{3}{2\b}\left(\int_{{\mathbf R} \times (B_{1/2} \setminus \{|x^{2}| < 1/16\})} {\rm dist}^{2}(X, {\rm spt} \, \|V_{k}\|) \,d\|{\mathbf C}_{k}\|(X)\right.&&\nonumber\\
&&\hspace{-2in} +\; \left.\int_{{\mathbf R} \times B_{1}} {\rm dist}^{2} \, (X, {\rm spt}\, \|{\mathbf C}_{k}\|) \, d\|V_{k}\|(X)\right),
\end{eqnarray*}
and ${\rm dist}_{\mathcal H}^{2} \, ({\rm spt} \, \|{\mathbf C}_{k}\| \cap ({\mathbf R} \times B_{1}), {\rm spt} \, \|\widetilde{\mathbf C}_{k}\| \cap ({\mathbf R} \times B_{1})) \leq {\overline C}(\widetilde{Q}_{k}^{2} + Q_{k}^{2}),$ where $\overline{C} = \overline{C}(n, q) \in (0, \infty),$ setting 
$$\k^{(5)}(n, q , \a, \th_{1}) = \frac{3\overline{\k}}{2\overline{\b}(n, q, \a, \th_{1})},\;\; C_{0}^{(5)}(n, q , \a, \th_{1}) = 2\overline{C} + \frac{3(\overline{C} + \overline{C}_{0})}{\overline{\b}(n, q, \a, \th_{1})},$$
$$C_{2}^{(5)}(n, q , \a, \th_{1}) = 2^{-\frac{n+4}{2}}\sqrt{\overline{C}_{1}\overline{C}} + \left(2^{-\frac{n+4}{2}}\sqrt{\overline{C}_{1}\overline{C}}+ \overline{C}_{2}\right)\sqrt{\frac{3}{2\overline{\b}(n, q, \a, \th_{1})}},$$
$$\nu_{1}^{(5)}(n, q, \a) = \overline{\nu},\;\;\; \nu_{2}^{(5)}(n, q, \a, \th_{1}) = \frac{3\overline{\nu}}{2\overline{\b}(n, q, \a, \th_{1})},$$
and, for $p \geq 6$, 
\begin{eqnarray*}
&&\k^{(p)}(n, q, \a, \th_{1} , \ldots, \th_{p-4}) =\\ 
&&\hspace{.1in}{\rm max} \, \left\{\overline{\k}, \; \frac{3}{2\overline{\b}(n, q, \a, \th_{1})}\k^{(p^{\prime})}(n, q, \a, \th_{2}, \ldots, \th_{p^{\prime} - 3}) \, : \, p^{\prime} = 5, \ldots, p-1\right\},
\end{eqnarray*}
\begin{eqnarray*}
&&C_{0}^{(p)}(n, q, \a, \th_{1}, \ldots, \th_{p-4}) =\\ 
&&\hspace{.5in} {\rm max} \, \left\{\overline{C}_{0}, \; 2\overline{C} + \frac{3}{\overline{\b}(n, q, \a, \th_{1})}\left(\overline{C} + C_{0}^{(p^{\prime})}(n, q, \a, \th_{2}, \ldots, \th_{p^{\prime} - 3})\right) \, : \, p^{\prime}=5, \ldots, p-1\right\},
\end{eqnarray*}
\begin{eqnarray*}
&&C_{2}^{(p)}(n, q, \a, \th_{1}, \ldots, \th_{p-4}) =\\ 
&&\hspace{.5in}{\rm max} \, \left\{\overline{C}_{2}, \; a + \sqrt{\frac{3}{2\overline{\b}(n, q, \a, \th_{1})}}\left(a + C_{2}^{(p^{\prime})}(n, q, \a, \th_{2}, \ldots, \th_{p^{\prime} - 3})\right) \, : \, p^{\prime}=5, \ldots, p-1\right\},
\end{eqnarray*}
where $a = 2^{-\frac{n+4}{2}}\sqrt{\overline{C}_{1}\overline{C}}$,
\begin{eqnarray*}
&&\nu_{1}^{(p)}(n, q, \a) = \overline{\nu},\;\;\; \nu_{2}^{(p)}(n, q, \a, \th_{1}) = \frac{3\overline{\nu}}{2\overline{\b}(n, q, \a, \th_{1})} \;\;\; \mbox{and, for} \;\;\; j=3, \ldots, p-3,\\
&&\nu_{j}^{(p)}(n, q, \a, \th_{1}, \ldots, \th_{j-1}) =\\
&&\hspace{.5in} {\rm max} \, \left\{\frac{3}{2\overline{\b}(n, q, \a, \th_{1})}\nu_{j-1}^{(p^{\prime})}(n, q, \a, \th_{2}, \ldots, \th_{j-1}) \, : \, p^{\prime} = j+2,  \ldots, p-1\right\}, 
\end{eqnarray*}
we see that if alternative $(B)$ holds, the conclusions (a)-(f) of the lemma 
follow with $V_{k}$, ${\mathbf C}_{k},$  ${\mathbf C}_{k}^{\prime}$, $\G_{k}$ in place of $V$, ${\mathbf C}$, ${\mathbf C}^{\prime}$, $\G;$ with constants $\k^{(p)}, C_{0}^{(p)},  C_{2}^{(p)}$ depending only on $n$, $q$, $\a$, $\th_{1}, \th_{2}, \ldots, \th_{p - 4};$ with $\nu_{1}^{(p)}$ depending only on $n$, $q$, $\a$  and 
for each $j=2, 3, \ldots, p - 3$, with $\nu_{j}^{(p)}$ depending only on 
$n$, $q,$ $\a$ and $\th_{1}, \th_{2}, \th_{3}, \ldots, \th_{j-1}.$ Note that in checking that conclusion (d) holds with 
$V_{k}$, ${\mathbf C}_{k}$  in place of $V$, ${\mathbf C},$ we have used the fact that 
${\rm dist}_{\mathcal H} \, ({\rm spt} \, \|\widetilde{\mathbf C}_{k}\| \cap ({\mathbf R} \times B_{1}),\{0\} \times B_{1}) \geq  
{\rm dist}_{\mathcal H} \, ({\rm spt} \, \|{\mathbf C}_{k}\| \cap ({\mathbf R} \times B_{1}),\{0\} \times B_{1}) - 
{\rm dist}_{\mathcal H} \, ({\rm spt} \, \|{\mathbf C}_{k}\| \cap ({\mathbf R} \times B_{1}), {\rm spt} \, \|\widetilde{\mathbf C}_{k}\| \cap ({\mathbf R} \times B_{1})) \geq {\rm dist}_{\mathcal H} \, ({\rm spt} \, \|{\mathbf C}_{k}\| \cap ({\mathbf R} \times B_{1}),\{0\} \times B_{1}) - 
\sqrt{\overline{C}}(\widetilde{Q}_{k} + Q_{k}).$  Similar reasoning applies in checking conclusion (b). This completes the proof. 
\end{proof}

\begin{lemma}\label{multi-scale-final} 
Let $q \geq 2$ be an integer and $\a \in (0, 1)$. For $j=1, 2, \ldots, 2q-3,$ let $\th_{j} \in (0, 1/4)$ be such that $\th_{1} > 8\th_{2} > 64\th_{3} > \ldots >8^{2q - 4}\th_{2q-3}$. 
There exist numbers $\e = \e(n, q, \a, \th_{1}, \th_{2}, \ldots, \th_{2q-3}) \in (0, 1/2)$, $\g= 
\g(n, q, \a, \th_{1}, \th_{2} \, \ldots, \th_{2q-3}) \in (0, 1/2)$ such that the following 
is true: If $V \in {\mathcal S}_{\a},$ ${\mathbf C} \in {\mathcal C}_{q}$ satisfy Hypotheses~\ref{hyp}, Hypothesis~($\star$) (of Section~\ref{fineblowup}) with $M = \frac{3}{2}M_{0}$ and if the induction hypotheses $(H1)$, $(H2)$ hold, 
then there exist an orthogonal rotation $\G$ of ${\mathbf R}^{n+1}$ and a cone ${\mathbf C}^{\prime} \in {\mathcal C}_{q}$ 
such that, with ${\hat E}_{V}$ and $Q_{V}({\mathbf C})$ as defined in Lemma~\ref{multi-scale}, we have the following:
\begin{eqnarray*}
&&({\rm a})\;\;\;\;|e_{1} - \G(e_{1})| \leq \k Q_{V}({\mathbf C}) \;\;{\rm and} \;\; |e_{j} - \G(e_{j})| \leq \k{\hat E}_{V}^{-1}Q_{V}({\mathbf C}) \;\;\mbox{for each} \;\; j=2, 3, \ldots, n+1;\nonumber\\
&&({\rm b})\;\;\;\; {\rm dist}_{\mathcal H}^{2} \, ({\rm spt} \, \|{\mathbf C}^{\prime}\| \cap ({\mathbf R} \times B_{1}), {\rm spt} \, \|{\mathbf C}\| \cap ({\mathbf R} \times B_{1})) \leq C_{0}Q_{V}^{2}({\mathbf C});\nonumber\\
&&\hspace{-.37in}\mbox{and for some $j \in \{1, 2, \ldots, 2q-3\}$},\nonumber\\ 
&&({\rm c})\;\;\;\; 
\th_{j}^{-n-2}\int_{\G\left({\mathbf R} \times \left(B_{\th_{j}/2} \setminus \{|x^{2}| \leq \th_{j}/16\}\right)\right)} {\rm dist}^
{2} \, (X, {\rm spt} \, \|V\|) \, d\|\G_{\#} \, {\mathbf C}^{\prime}\|(X)\nonumber\\
&& \hspace{1.5in} + \; \th_{j}^{-n-2}\int_{\G({\mathbf R} \times B_{\th_{j}})} {\rm dist}^{2} \, (X, {\rm spt} \, \|\G_{\#} \, {\mathbf C}^{\prime}\|) \, d\|V\|(X) \leq \nu_{j}\th_{j}^{2}Q_{V}^{2}({\mathbf C});\nonumber\\
&&({\rm d})\;\;\;\;\left(\th_{j}^{-n-2}\int_{{\mathbf R} \times B_{\th_{j}}} {\rm dist}^{2}(X, P) \, d\|\G^{-1}_{\#} \,V\|(X)\right)^{1/2}\nonumber\\ 
&&\hspace{1in}\geq 2^{-\frac{n+4}{2}}\sqrt{\overline{C}_{1}} \,{\rm dist}_{\mathcal H}\, ({\rm spt} \, \|{\mathbf C}\| \cap ({\mathbf R} \times B_{1}),P \cap ({\mathbf R} \times  B_{1})) - C_{2}Q_{V}({\mathbf C})\nonumber\\
&&\;\;\;\;\mbox{for any $P \in G_{n}$ of the form $P = \{x^{1} = \lambda x^{2}\}$ for some $\lambda \in (-1, 1)$};\nonumber\\
&&({\rm e})\;\;\;\; \{Z \, : \, \Theta \, (\|\G_{\#}^{-1} \, V\|, Z) \geq q\} \cap \left({\mathbf R} \times (B_{\th_{j}/2} \cap \{|x^{2}| < \th_{j}/16\})\right) = \emptyset;\nonumber\\
&&({\rm f})\;\;\;\;\left(\omega_{n}\th_{j}^{n}\right)^{-1}\|\G_{\#}^{-1} \, V\|({\mathbf R} \times B_{\th_{j}}) < q + 1/2. 
\end{eqnarray*}
Here the constants $\k, C_{0}, C_{2} \in (0, \infty)$ depend only on $n, \a$ in case $q=2$ and only on $n$, $q,$ $\a$ and $\th_{1}, \th_{2}, \ldots, \th_{2q-4}$ in case $q \geq 3$; $\nu_{1} = \nu_{1}(n, q,\a);$ and, 
in case $q \geq 3,$ for each $j=2, 3, \ldots, 2q-3,$ $\nu_{j} = \nu_{j}(n, q, \a, \th_{1}, \ldots, \th_{j-1}).$ (In particular, $\nu_{j}$ is independent of $\th_{j}, \th_{j+1}, \ldots, \th_{2q-3}$ for each $j = 1, 2, \ldots, 2q-3$.) 
\end{lemma}

\begin{proof}
Set $\e = {\rm min} \,  \left\{\e^{(4)}, \e^{(5)}, \ldots, \e^{(2q)}\right\}$ and $\g = {\rm min} \, \left\{\g^{(4)}, \g^{(5)}, \ldots, \g^{(2q)}\right\}$, where 
$$\e^{(p)} = \e^{(p)}(n, q, \a, \th_{1}, \ldots, \th_{p-3}), \;\;\; \g^{(p)} = \g^{(p)}(n, q, \a, \th_{1}, \ldots, \th_{p-3}), \;\;\; 4 \leq p \leq 2q,$$ 
are as in Lemma~\ref{multi-scale}. Set $\nu_{1} = \overline{\nu}$ and for each $j=2, \ldots, 2q-3$, 
$$\nu_{j} = {\rm max} \, \left\{\nu^{(j+3)}_{j}, \nu^{(j+4)}_{j}, \ldots, \nu^{(2q)}_{j}\right\} \; \left( = \nu_{j}^{(2q)}\right)$$ 
where $\overline{\nu}$ is as in Lemma~\ref{excess-improvement}, and for each $p \in \{5, \ldots, 2q\}$, the numbers $\nu_{j}^{(p)}$  are as in Lemma~\ref{multi-scale} taken with scales 
$\th_{1}, \ldots, \th_{p-3}$. Note that then, $\nu_{1} = \nu_{1}(n, q, \a)$ and in case $q \geq 3$, 
$$\nu_{j} = \nu_{j}(n, q, \a, \th_{1}, \ldots, \th_{j-1}) \;\; {\rm for} \;\;  2 \leq j \leq 2q-3.$$ 
$${\rm Set} \;\; \k = {\rm max}\, \left\{\k^{(4)}, \k^{(5)}, \ldots, \k^{(2q)}\right\} \; \left(= \k^{(2q)}\right), \;\; C_{0} = {\rm max} \, \left\{C_{0}^{(4)}, C_{0}^{(5)}, \ldots, C_{0}^{(2q)}\right\}\; \left(= C_{0}^{(2q)}\right),$$
$$C_{2} = {\rm max} \, \left\{C_{2}^{(4)}, C_{2}^{(5)}, \ldots, C_{2}^{(2q)}\right\} \; \left(=C_{2}^{(2q)}\right),$$ 
where for each $p \in \{4, 5, \ldots, 2q\}$, the numbers $\k^{(p)}$,  $C_{0}^{(p)}$, $C_{2}^{(p)}$ are as in Lemma~\ref{multi-scale} taken with scales $\th_{1}, \ldots, \th_{p-3}.$ Since ${\mathbf C} \in {\mathcal C}_{q}$ implies that ${\mathbf C} \in {\mathcal C}_{q}(p)$ for some $p \in \{4, 5, \ldots, 2q\}$,  the conclusions of the present lemma follow directly from Lemma~\ref{multi-scale}. 
\end{proof}

\section{Properties of coarse blow-ups: Part III}\label{propertiesIII}
\setcounter{equation}{0}

Subject to the induction hypotheses $(H1)$, $(H2)$, we complete in this section the proof that ${\mathcal B}_{q}$ is a proper blow-up class, by showing that ${\mathcal B}_{q}$ satisfies property (${\mathcal B{\emph 7}}$). Recall that in order to do this, it only remains to rule out the possibility that ${\mathcal B}_{q}$ contains an element whose graph is the union of $q$ half-hyperplanes in the half-space $\{x^{2} \leq 0\}$ and $q$ half-hyperplanes in $\{x^{2} \geq 0\},$ with all half-hyperplanes meeting along $\{0\} \times {\mathbf R}^{n-1},$ and with at least two of the half-hyperplanes on each side distinct. (This is {\bf Case 2} stated at the beginning of Section~\ref{step3}.) 

\begin{lemma}\label{finedistgap}
Let $q \geq 2$ be an integer and $\a \in (0, 1)$. There exist constants $\e_{1} = \e_{1}(n, q, \a) \in (0, 1)$ and $\g_{1} = \g_{1}(n, q, \a) \in (0, 1)$ such that 
if the induction hypotheses $(H1)$, $(H2)$ hold, $V \in {\mathcal S}_{\a}$, 
$$\Theta \, (\|V\|, 0) \geq q, \;\;\; (\omega_{n}2^{n})^{-1}\|V\|(B_{2}^{n+1}(0)) < q + 1/2, \;\;\; \omega_{n}^{-1}\|V\|({\mathbf R} \times B_{1}) < q + 1/2,$$  
$$\{Z \, : \, \Theta \, (\|V\|, Z) \geq q\} \cap \left({\mathbf R} \times (B_{1/2} \setminus \{|x^{2}| < 1/16\})\right) = \emptyset,$$ 
$${\hat E}_{V}^{2} \equiv \int_{{\mathbf R} \times B_{1}} |x^{1}|^{2} d\|V\|(X) < \e_{1}\;\;\; {\rm and}$$ 
$${\hat E}_{V}^{2} < \frac{3}{2}\inf_{\{P = \{x^{1} = \lambda x^{2}\}\}} \int_{{\mathbf R} \times B_{1}} {\rm dist}^{2} \, (X, P) \, d\|V\|(X) \;\;\; {\rm then}$$
$$\int_{{\mathbf R} \times (B_{1/2} \setminus \{|x^2| < 1/16\})} {\rm dist}^{2}(X, {\rm spt} \, \|V\|) \,d\|{\mathbf C}\|(X) + \int_{{\mathbf R} \times B_{1}} {\rm dist}^{2} \, (X, {\rm spt}\, \|{\mathbf C}\|) \, d\|V\|(X) \geq \g_{1} {\hat E}_{V}^{2}$$ 
for any cone ${\mathbf C} \in {\mathcal C}_{q}.$ 
\end{lemma}

\begin{proof}
For $j=1, 2, \ldots, 2q-3,$ choose numbers $\th_{j}  = \th_{j}(n, q, \a) \in (0, 1/2)$ as follows:  First choose 
$\th_{1} = \th_{1}(n , q, \a) \in (0, 1/2)$ such that $\nu_{1}\th_{1}^{2(1- \a)} < 1,$ where $\nu_{1} = \nu_{1}(n, q, \a)$ is as in Lemma~\ref{multi-scale-final}. Having chosen $\th_{1}, \th_{2}, \ldots, \th_{j}$, 
$1 \leq j \leq 2q-4$, choose $\th_{j+1} = \th_{j+1}(n, q, \a)$ such that $\th_{j+1} < 8^{-1}\th_{j}$ and 
$\nu_{j+1}\th_{j+1}^{2(1 - \a)}< 1,$ where $\nu_{j+1} = \nu_{j+1}(n, q, \a, \th_{1}, \th_{2}, \ldots, \th_{j})$ is as in Lemma~\ref{multi-scale-final}. 

Let $\e_{1}  \in (0, \e),$ $\g_{1} \in (0, \g)$ be constants to be eventually chosen  depending only on $n,$ $q$ and $\a$, where $\e = \e(n, q, \a, \th_{1}, \ldots, \th_{2q-3}),$  $\g = \g(n, q, \a, \th_{1}, \ldots, \th_{2q-3})$ are as in Lemma~\ref{multi-scale-final}. Suppose that the hypotheses of the present lemma are satisfied with $V \in {\mathcal S}_{\a}$ but the conclusion fails, i.e.\ there exists ${\mathbf C} \in {\mathcal C}_{q}$ such that 
\begin{eqnarray}\label{finedistgap-0}
\int_{{\mathbf R} \times (B_{1/2} \setminus \{|x^2| < 1/16\})} {\rm dist}^{2}(X, {\rm spt} \, \|V\|) \,d\|{\mathbf C}\|(X)&&\nonumber\\
&&\hspace{-1.5in}  + \; \int_{{\mathbf R} \times B_{1}} {\rm dist}^{2} \, (X, {\rm spt}\, \|{\mathbf C}\|) \, d\|V\|(X) < \g_{1} {\hat E}_{V}^{2}.
\end{eqnarray}
In particular, $V$, ${\mathbf C}$ then satisfy the hypotheses of Lemma~\ref{multi-scale-final}. In what follows, for ${\mathbf C}^{\prime} \in {\mathcal C}_{q}$, $\G$ an orthogonal rotation of ${\mathbf R}^{n+1}$ and $\r \in (0, 1]$, we shall use the notation 
\begin{eqnarray*}
Q_{V}({\mathbf C}^{\prime}, \G, \r) = \left(\r^{-n-2}\int_{\G\left({\mathbf R} \times (B_{\r/2} \setminus \{|x^2| < \r/16\})\right)} {\rm dist}^{2}(X, {\rm spt} \, \|V\|) \,d\|\G_{\#} \, {\mathbf C}^{\prime}\|(X)\right.&&\\ 
&& \hspace{-2.5in}+ \; \left.\r^{-n-2}\int_{\G\left({\mathbf R} \times B_{\r}\right)} {\rm dist}^{2} \, (X, {\rm spt}\, \|\G_{\#} \, {\mathbf C}^{\prime}\|) \, d\|V\|(X)\right)^{1/2}. 
\end{eqnarray*}
We claim that we may apply Lemma~\ref{multi-scale-final} iteratively to obtain, for each $k=0,1, 2, 3, \ldots$,  an orthogonal rotation $\G_{k}$ of ${\mathbf R}^{n+1},$ with $\G_{0} =$ Identity, and a cone ${\mathbf C}_{k} \in {\mathcal C}_{q}$  with ${\mathbf C}_{0} = {\mathbf C}$, satisfying, for $k \geq 1$, 
\begin{equation}\label{finedistgap-a}
|\G_{k}(e_{1}) - \G_{k-1}(e_{1})|^{2} \leq C\d_{k}Q^{2}_{V};
\end{equation}
\begin{equation}\label{finedistgap-0-0}
|\G_{k}(e_{j}) - \G_{k-1}(e_{j})|^{2} \leq C\d_{k}{\hat E}_{V}^{-2}Q^{2}_{V};
\end{equation}
\begin{equation}\label{finedistgap-1}
{\rm dist}_{\mathcal H}^{2} \, ({\rm spt} \, \|{\mathbf C}_{k}\| \cap ({\mathbf R} \times B_{1}), {\rm spt} \, \|{\mathbf C}_{k-1}\| \cap ({\mathbf R} \times B_{1})) \leq  C\d_{k}Q^{2}_{V};
\end{equation}
\begin{equation}\label{finedistgap-2}
Q_{V}^{2}({\mathbf C}_{k}, \G_{k}, \s_{k})  \leq  \nu_{j_{k}}\th_{j_{k}}^{2} Q_{V}^{2}({\mathbf C}_{k-1}, \G_{k-1}, \s_{k-1}) \leq \ldots \leq \d_{k} Q^{2}_{V}
\end{equation}
for some $j_{k} \in \{1, 2, \ldots, 2q-3\};$
\begin{eqnarray}\label{finedistgap-4}
&&\left(\s_{k}^{-n-2}\int_{{\mathbf R} \times B_{\s_{k}}} {\rm dist}^{2} \, (X, P) \, d\|\G_{k \, \#} \,V\|(X)\right)^{1/2} \geq\nonumber\\ 
&&\hspace{.5in}2^{-\frac{n+4}{2}}\sqrt{\overline{C}_{1}} \,{\rm dist}_{\mathcal H} \, ({\rm spt} \, \|{\mathbf C}_{k-1}\|\cap ({\mathbf R} \times B_{1}), P \cap ({\mathbf R} \times   B_{1})) - C_{2} \,Q_{V}({\mathbf C}_{k-1}, \G_{k-1}, \s_{k-1}) 
\end{eqnarray}
for each $P \in G_{n}$ of the form $P = \{x^{1} = \lambda x^{2}\}$ for some $\lambda \in (-1, 1);$
\begin{equation}\label{finedistgap-4-a}
\{Z \, : \, \Theta \, (\|\G_{k \, \#} \, V\|, Z) \geq q\} \cap \left({\mathbf R} \times (B_{\s_{k}} \setminus \{|x^{2}| < \s_{k}/16\})\right) = \emptyset; \;\; {\rm and}
\end{equation}
\begin{equation}\label{finedistgap-4-b}
\left(\omega_{n}\s_{k}^{n}\right)^{-1}\|\G_{k \, \#}^{-1} \, V\|({\mathbf R} \times B_{\s_{k}}) < q + 1/2
\end{equation}
where $Q_{V} = Q_{V}({\mathbf C}, \G_{0}, 1)$, $C = C(n, q, \a) \in (0, \infty)$, $C_{2} = C_{2}(n, q, \a) \in (0, \infty)$ and, 
for each $k=1, 2, 3, \ldots,$ 
$$\s_{k} = \th_{j_{k}}\s_{k-1}, \;\;\;\; \d_{k} = \nu_{j_{k}}\th_{j_{k}}^{2}\d_{k-1}$$ for some $j_{k} \in \{1, 2, \ldots, 2q-3\}$ where $\s_{0} = \d_{0} = 1.$ Thus 
$$\s_{k} = \prod_{j=1}^{2q-3} \th_{j}^{k_{j}} \;\;\;\;\;\; {\rm and} \;\;\;\;\;\; \d_{k} = \prod_{j=1}^{2q-3} \left(\nu_{j}\th_{j}^{2}\right)^{k_{j}}$$
for some non-negative integers $k_{1}, k_{2}, \ldots, k_{2q-3}$ such that $\sum_{j=1}^{2q-3} k_{j} = k.$ Note in particular that 
$$\th_{2q-3}^{k} \leq \s_{k} \leq \th_{1}^{k}, \;\;\;\;\;\; \d_{k} < \s_{k}^{2\a} < 4^{-k\a} \;\;\;\;\;\; {\rm and} \;\;\;\;\;\; \sum_{j=k}^{\infty} \d_{j} < c\d_{k}$$ 
for $k=1, 2, \ldots,$ where $c = c(\a) \in (0, \infty).$

To verify these assertions inductively, note that (\ref{finedistgap-1})-(\ref{finedistgap-4-a}) with $k=1$ follow directly from Lemma~\ref{multi-scale-final}. Suppose $k \geq 2$ and that (\ref{finedistgap-1})-(\ref{finedistgap-4-a}) hold with $1, 2, 3, \ldots, k-1$ in place of $k$. We wish to apply Lemma~\ref{multi-scale-final} with 
$\eta_{\s_{k-1} \, \#} \, \G_{k-1 \, \#}^{-1} \, V$ in place of $V$ and ${\mathbf C}_{k-1}$ in place of ${\mathbf C}.$ Note first that by the triangle inequality and (\ref{finedistgap-4-b}) with $k-1$ in place of $k$,
\begin{eqnarray*}
{\hat E}^{2}_{\eta_{\s_{k-1} \, \#} \, \G_{k-1 \, \#}^{-1} \, V} = \s_{k-1}^{-n-2}\int_{{\mathbf R} \times B_{\s_{k-1}}}
|x^{1}|^{2} \, d\|\G^{-1}_{k-1 \, \#} \, V\|(X)&&\nonumber\\
&&\hspace{-2.5in} \leq 2 \,\s_{k-1}^{-n-2}\int_{{\mathbf R} \times B_{\s_{k-1}}}
{\rm dist}^{2} \, (X, {\rm spt} \, \|{\mathbf C}_{k-1}\|) \, d\|\G^{-1}_{k-1 \, \#} \, V\|(X)\nonumber\\ 
&&\hspace{-1in} + \;  \omega_{n}(2q+1) \,{\rm dist}^{2}_{\mathcal H} \, ({\rm spt} \, \|{\mathbf C}_{k-1}\| \cap ({\mathbf R} \times B_{1}), \{0\} \times B_{1}) 
\end{eqnarray*}
and by applying (\ref{finedistgap-1}) with $1, 2, \ldots, k-1$ in place of $k$, summing over $k$, and using the fact that 
$\sum_{k=1}^{\infty} \d_{k}^{1/2} < 2^{-\a}(1 - 2^{-\a})^{-1}$, 
\begin{eqnarray*}
&&{\rm dist}_{\mathcal H} \, ({\rm spt} \, \|{\mathbf C}_{k-1}\| \cap ({\mathbf R} \times B_{1}), \{0\} \times B_{1}) \leq\nonumber\\
&&\hspace{1in}{\rm dist}_{\mathcal H} \, ({\rm spt} \, \|{\mathbf C}\| \cap ({\mathbf R} \times B_{1}), \{0\} \times B_{1}) + CQ_{V}, \; C = C(n, q,\a) \in (0, \infty); \;\; {\rm thus}, 
\end{eqnarray*}
\begin{equation*}
{\hat E}^{2}_{\eta_{\s_{k-1} \, \#} \, \G_{k-1 \, \#}^{-1} \, V} \leq 2\omega_{n}(2q+1) \,{\rm dist}^{2}_{\mathcal H} \, ({\rm spt} \, \|{\mathbf C}\| \cap ({\mathbf R} \times B_{1}), \{0\} \times B_{1}) + CQ^{2}_{V}, \;\; C = C(n, q, \a) \in (0, \infty), 
\end{equation*}
so that, by (\ref{slope-bounds-1}),
\begin{equation}\label{finedistgap-5-0}
{\hat E}^{2}_{\eta_{\s_{k-1} \, \#} \, \G_{k-1 \, \#}^{-1} \, V} \leq 2(2q+1)\omega_{n} c_{1}^{2}{\hat E}_{V}^{2}+ CQ^{2}_{V}, \; C = C(n, q, \a) \in (0, \infty), 
\end{equation}
where $c_{1} = c_{1}(n) \in (0, \infty)$ is as in (\ref{slope-bounds-1}); in particular,
\begin{equation}\label{finedistgap-5}
{\hat E}^{2}_{\eta_{\s_{k-1} \, \#} \, \G_{k-1 \, \#}^{-1} \, V} \leq C{\hat E}^{2}_{V}, \;\;\; C=C(n, q, \a) \in (0, \infty).
\end{equation}
Again by (\ref{finedistgap-1}), 
\begin{eqnarray*}
{\rm dist}_{\mathcal H} \, ({\rm spt} \, \|{\mathbf C}_{k-2}\| \cap ({\mathbf R} \times B_{1}), \{0\} \times B_{1}) \geq 
 {\rm dist}_{\mathcal H} \, ({\rm spt} \, \|{\mathbf C}\| \cap {\mathbf R} \times B_{1}, \{0\} \times B_{1})&&\nonumber\\
 &&\hspace{-4in} - \sum_{j=1}^{k-2} {\rm dist}_{\mathcal H} \, ({\rm spt} \, \|{\mathbf C}_{j-1}\| \cap ({\mathbf R} \times B_{1}), {\rm spt} \, \|{\mathbf C}_{j}\| \cap ({\mathbf R} \times B_{1}))\nonumber\\
 &&\hspace{-3.5in}\geq {\rm dist}_{\mathcal H} \, ({\rm spt} \, \|{\mathbf C}\| \cap ({\mathbf R} \times B_{1}), \{0\} \times B_{1})  -CQ_{V}\sum_{j=1}^{k-2}\d_{j}^{1/2}
\end{eqnarray*}
which implies by (\ref{finedistgap-4}) and (\ref{finedistgap-2}) that 
\begin{equation*}
{\hat E}_{\eta_{\s_{k-1} \, \#} \, \G_{k-1 \, \#}^{-1} \, V} \geq 2^{-\frac{n+4}{2}}\sqrt{\overline{C}_{1}}\,{\rm dist}_{\mathcal H} \, ({\rm spt} \, 
\|{\mathbf C}\|\cap ({\mathbf R} \times B_{1}), \{0\} \times B_{1}) - CQ_{V}
\end{equation*}
where $C = C(n, q, \a) \in (0, \infty);$ hence by (\ref{slope-bounds-2}) and (\ref{finedistgap-0}), we see that  
\begin{equation}\label{finedistgap-7}
{\hat E}_{\eta_{\s_{k-1} \, \#} \, \G_{k-1 \, \#}^{-1} \, V} \geq (C_{1} - C\g_{1}){\hat E}_{V} 
\end{equation}
where $C_{1} = C_{1}(n, q)$, $C = C(n, q) \in (0, \infty).$ Thus if $2C\g_{1} < C_{1}$, it follows from (\ref{finedistgap-0}), (\ref{finedistgap-2}) and (\ref{finedistgap-7}) that 
\begin{eqnarray}\label{finedistgap-8}
\int_{{\mathbf R} \times \left(B_{1/2} \setminus \{|x^{2}| < 1/16\}\right)} {\rm dist}^{2} \, (X, {\rm spt} \, 
\|\eta_{\s_{k-1} \, \#} \, \G_{k-1 \, \#}^{-1} \, V\|) \, d\|{\mathbf C}_{k-1}\|(X)&&\nonumber\\
&&\hspace{-4in} + \int_{{\mathbf R} \times B_{1}} {\rm dist} \, ^{2}(X, {\rm spt} \, \|{\mathbf C}_{k-1}\|) \, 
d\|\eta_{\s_{k-1} \, \#} \, \G_{k-1 \, \#}^{-1} \, V\|(X)\leq C\g_{1}{\hat E}^{2}_{\eta_{\s_{k-1} \, \#} \, \G_{k-1 \, \#}^{-1} \, V}
\end{eqnarray}
and from (\ref{finedistgap-5}), that
\begin{equation*}
{\hat E}^{2}_{\eta_{\s_{k-1} \, \#} \, \G_{k-1 \, \#}^{-1} \, V} \leq C\e_{1}
\end{equation*}
where $C = C(n, q, \a) \in (0, \infty).$ By (\ref{finedistgap-4}) again with $k-1$ in place of $k$ and (\ref{finedistgap-1}) with $1, 2, \ldots, k-1$ in place of $k$, 
\begin{eqnarray*}
&&\left(\s_{k-1}^{-n-2}\int_{{\mathbf R} \times B_{\s_{k-1}}} {\rm dist}^{2} \, (X, P) \, d\|\G_{k-1 \, \#} \,V\|(X)\right)^{1/2} \geq\nonumber\\ 
&&\hspace{.5in}2^{-\frac{n+4}{2}}\sqrt{\overline{C}_{1}} \,{\rm dist}_{\mathcal H} \, ({\rm spt} \, \|{\mathbf C}\|\cap ({\mathbf R} \times B_{1}), P \cap ({\mathbf R} \times  B_{1})) - C Q_{V}\nonumber\\
\end{eqnarray*}
so that 
\begin{eqnarray*}
&&\int_{{\mathbf R} \times B_{1}} {\rm dist}^{2} \, (X, P) \, d\|\eta_{\s_{k-1} \, \#}\G_{k-1 \, \#} \,V\|(X) \geq\nonumber\\ 
&&\hspace{1in}2^{-n-5}\overline{C}_{1} \,{\rm dist}^{2}_{\mathcal H} \, ({\rm spt} \, \|{\mathbf C}\|\cap ({\mathbf R} \times B_{1}), P \cap ({\mathbf R} \times  B_{1})) - CQ^{2}_{V}\nonumber\\
&&\hspace{.5in}\geq 2^{-n-5}\overline{C}_{1}\omega_{n}^{-1}(2q+1)^{-1}\int_{{\mathbf R} \times B_{1}}{\rm dist}^{2} \, (X, P) \, d\|V\|(X) -CQ_{V}^{2}\nonumber\\
&&\hspace{.5in}\geq 2^{-n-5}\overline{C}_{1}\omega_{n}^{-1}(2q+1)^{-1}\left(\frac{3}{2}\right)^{-1}{\hat E}_{V}^{2} - CQ_{V}^{2}\nonumber\\
&&\hspace{.5in}\geq 2^{-n-6}\overline{C}_{1} \omega_{n}^{-2}(2q+1)^{-2}c_{1}^{-2}\left(\frac{3}{2}\right)^{-1}{\hat E}^{2}_{\eta_{\s_{k-1} \, \#} \G_{k-1 \, \#} \, V} - CQ_{V}^{2}\nonumber\\
&&\hspace{.5in}\geq  \left(2^{-n-6}\overline{C}_{1} \omega_{n}^{-2}(2q+1)^{-2}c_{1}^{-2}\left(\frac{3}{2}\right)^{-1} - C\g_{1}\right){\hat E}^{2}_{\eta_{\s_{k-1} \, \#} \G_{k-1 \, \#} \, V} 
\end{eqnarray*}
where $C = C(n, q, \a) \in (0, \infty)$ and we have used our hypothesis that
$${\hat E}_{V}^{2} < \frac{3}{2} \inf_{P = \{x^{1} = \lambda x^{2}\}} \, \int_{{\mathbf R} \times B_{1}}{\rm dist}^{2} \, (X, P) \, d\|V\|(X).$$ 
This readily implies that if we choose $\g_{1} = \g_{1}(n, q, \a) \in (0, 1)$ sufficiently small, then 
\begin{equation*}
{\hat E}^{2}_{\eta_{\s_{k-1} \, \#} \G_{k-1 \, \#} \, V} \leq \frac{3}{2}M_{0}\int_{{\mathbf R} \times B_{1}} {\rm dist}^{2} \, (X, P) \, d\|\eta_{\s_{k-1} \, \#}\G_{k-1 \, \#} \,V\|(X)
\end{equation*}
for any hyperplane $P$ of the form $P = \{x^{1} = \lambda x^{2}\}$. So  if we choose $\g_{1} = \g_{1}(n, q, \a)$ and $\e_{1} = \e_{1}(n, q, \a)$ sufficiently small, we can apply Lemma~\ref{multi-scale-final} with  $\eta_{\s_{k-1} \, \#} \, \G_{k-1 \, \#}^{-1} \, V$ in place of $V$ and ${\mathbf C}_{k-1}$ in place of ${\mathbf C}$ to obtain an orthogonal rotation $\G$ of ${\mathbf R}^{n+1}$ and a cone ${\mathbf C}_{k} \in {\mathcal C}_{q}$ such that, with $\G_{k} = \G_{k-1} \circ \G$, (\ref{finedistgap-a})-(\ref{finedistgap-4-b}) hold. This completes the inductive proof that (\ref{finedistgap-a})-(\ref{finedistgap-4-b}) hold for 
all $k=1, 2, 3, \ldots.$ Writing 
$${\mathbf C}_{k} = \sum_{j=1}^{q} |H_{j}^{k}| + |G_{j}^{k}|$$ 
where for each $j \in \{1, 2, \ldots, q\}$, $H_{j}^{k}$ is the half-space defined by
$H_{j}^{k} = \{(x^{1}, x^{2}, y) \in {\mathbf R}^{n+1} \, : \, x^{2} < 0 \;\; \mbox{and} \;\; x^{1} = \lambda_{j}^{k}x^{2}\},$ 
$G_{j}^{k}$ the half-space defined by $G_{j}^{k} = \{(x^{1},x^{2}, y) \in {\mathbf R}^{n+1}\, : \, x^{2} > 0 \;\; \mbox{and} \;\; x^{1} = \mu_{j}^{k}x^{2}\},$ with $\lambda_{j}^{k}, \mu_{j}^{k}$ constants,  
$\lambda_{1}^{k} \geq \lambda_{2}^{k} \geq \ldots \geq \lambda_{q}^{k}$ and $\mu_{1}^{k} \leq \mu_{2}^{k} \leq \ldots \leq \mu_{q}^{k}$, note that by (\ref{slope-bounds-2}) (applied with $\eta_{\s_{k} \, \#} \, \G_{k \, \#} \,V$ in place of $V$ and ${\mathbf C}_{k}$ in place of ${\mathbf C}$), (\ref{finedistgap-7}) and (\ref{finedistgap-8}), we also have that
\begin{equation}\label{finedistgap-9}
|\lambda_{1}^{k} - \lambda_{q}^{k}| \geq C{\hat E}_{V} \;\;\; {\rm and} \;\;\; |\mu_{1}^{k} - \mu_{q}^{k}| \geq C{\hat E}_{V}, \;\;\; C = C(n, q, \a) \in (0, \infty),
\end{equation}
for all $k =1, 2, 3, \ldots.$

By (\ref{finedistgap-1}) $\{{\rm spt} \, \|{\mathbf C}_{k}\|\cap ({\mathbf R} \times B_{1})\}$ is a Cauchy sequence (in Hausdorff distance) and hence, since $\Theta \, (\|{\mathbf C}_{k}\|, 0) = q$ for each $k=1,2,\ldots$, there is a 
varifold ${\mathbf H} \in {\mathcal C}_{q}$ such that passing to a subsequence $\{k^{\prime}\}$ of $\{k\}$, 
${\mathbf C}_{k^{\prime}} \res B_{2}^{n+1}(0) \to {\mathbf H} \res B_{2}^{n+1}(0)$ and  
\begin{equation}\label{finedistgap-9-1-0}
{\rm dist}_{\mathcal H}^{2} \, ({\rm spt} \, \|{\mathbf H}\| \cap ({\mathbf R} \times B_{1}), {\rm spt} \, \| {\mathbf C}_{k^{\prime}}\| \cap ({\mathbf R} \times B_{1})) \leq  C\d_{k^{\prime}}Q^{2}_{V}
\end{equation}
for each $k^{\prime}$, where $C = C(n, q, \a) \in (0, \infty).$ 
By (\ref{finedistgap-9}), ${\rm spt} \, \|{\mathbf H}\|$ is not a hyperplane. Since $\d_{k} \leq \s_{k}^{2\a}$, it follows from (\ref{finedistgap-2}), (\ref{finedistgap-4-b}) and (\ref{finedistgap-9-1-0}) that 
\begin{equation}\label{finedistgap-9-1-1}
\int_{{\mathbf R} \times B_{1}} {\rm dist}^{2} \, (X, {\rm spt} \,\|{\mathbf H}\|) \, d\|\eta_{\s_{k^{\prime}} \, \#} \G^{-1}_{k^{\prime} \, \#} \, V\|(X) \leq C\s_{k^{\prime}}^{2\a} Q_{V}^{2} \;\; {\rm and}
\end{equation}
\begin{equation}\label{finedistgap-9-1-1-1}
\int_{{\mathbf R} \times \left(B_{1/2} \setminus \{|x^{2}| < 1/16\}\right)} {\rm dist}^{2} \, (X, {\rm spt} \, \|\eta_{\s_{k^{\prime}} \, \#} \G^{-1}_{k^{\prime} \, \#} \, V\|) \, d\|{\mathbf C}_{k^{\prime}}\|(X)\leq C \s_{k^{\prime}}^{2\a} Q_{V}^{2}
\end{equation}
for all $k^{\prime},$ where $C = C(n, q, \a) \in (0, \infty).$ Now, since $q \leq \Theta \, (\|V\|, 0) \leq (\omega_{n}2^{n})^{-1}\|V\|(B_{2}^{n+1}(0)) < q + 1/2$, it follows from the monotonicity formula that $$q \leq \Theta \, (\|\eta_{\s_{k^{\prime}} \, \#} \G^{-1}_{k^{\prime} \, \#} \, V\|, 0) \leq (\omega_{n}2^{n})^{-1}\|\eta_{\s_{k^{\prime}} \, \#} \G^{-1}_{k^{\prime} \, \#} \, V\|(B_{2}^{n+1}(0))  \leq (\omega_{n}2^{n})^{-1}\|V\|(B_{2}^{n+1}(0)) < q + 1/2.$$ 
Hence, there is a stationary integral varifold $W$ on $B_{2}^{n+1}(0)$ with 
$$q \leq \Theta \, (\|W\|, 0) \leq (\omega_{n}2^{n})^{-1}\|W\|(B_{2}^{n+1}(0)) < q + 1/2$$ 
such that, passing to a further subsequence without changing notation,  
\begin{equation}\label{finedistgap-10}
\eta_{\s_{k^{\prime}} \, \#} \G^{-1}_{k^{\prime} \, \#} \,V \to W
\end{equation}
as varifolds on $B_{2}^{n+1}(0)$. The estimate (\ref{finedistgap-9-1-1}) implies that ${\rm spt} \, \|W\|\cap ({\mathbf R} \times B_{1}) \subseteq {\rm spt} \,\|{\mathbf H}\| \cap ({\mathbf R} \times B_{1}),$  and since 
${\rm dist}_{\mathcal H} \, ({\rm spt} \, \|\eta_{\s_{k^{\prime}} \, \#} \G_{k^{\prime} \, \#} \, V\| \cap ({\mathbf R} \times B_{1}), {\rm spt} \, \|W\| \cap ({\mathbf R} \times B_{1})) \to 0$, it follows from (\ref{finedistgap-9-1-1-1}), the triangle inequality and the weak convergence $\|{\mathbf C}_{k^{\prime}}\| \to \|{\mathbf H}\|$ on ${\mathbf R} \times B_{1/2}$ that 
${\rm spt} \, \|{\mathbf H}\| \cap ({\mathbf R} \times (B_{1/2} \setminus \{|x^{2}| < 1/16\})) \subseteq {\rm spt} \, \|W\| \cap ({\mathbf R} \times (B_{1/2} \setminus \{|x^{2}| < 1/16\})).$ Hence 
${\rm spt} \, \|W\| \cap ({\mathbf R} \times B_{1}) = {\rm spt} \, \|{\mathbf H}\| \cap ({\mathbf R} \times B_{1}),$ from which it also follows that $\Theta \, (\|W\|, 0) = q.$ Thus  (\ref{finedistgap-10}) contradicts case $\Theta \, (\|{\mathbf C}_{0}\|, 0) = q$ of the induction hypothesis $(H2)$, proving the Lemma.
\end{proof}

\begin{corollary}\label{F7}
Let $q$ be an integer $\geq 2$, and suppose that the induction hypotheses $(H1)$, $(H2)$ hold. Then 
the class ${\mathcal B}_{q}$ (defined in Section~\ref{blow-up}) satisfies property $({\mathcal B{\emph 7}})$ of the definition of proper blow-up classes (given in Section~\ref{proper-blow-up}). 
\end{corollary}

\begin{proof} If not, in view of Lemma~\ref{no-overlap}, there exists an element $v_{\star} \in {\mathcal B}_{q}$ 
such that, for $j=1, 2, \ldots, q$, $v_{\star}^{j}(x^{2}, y) = L_{1}^{j}(x^{2}, y)$ if $x^{2} < 0$ and $v_{\star}^{j}(x^{2}, y) = L_{2}^{j}(x^{2}, y)$ if $x^{2} \geq 0$  where $L_{1}^{j}, L_{2}^{j} \, : \, {\mathbf R}^{n} \to {\mathbf R}$ are linear functions with $L_{1}^{j}(0,y) = L_{2}^{j}(0, y) = 0$ and 
\begin{equation}\label{partIII-0-0}
L_{1}^{j_{1}} \neq L_{1}^{j_{1}+1} \;\;\; {\rm and} \;\;\;L_{2}^{j_{2}} \neq L_{2}^{j_{2} +1} \;\;\; \mbox{for some} \;\;\; j_{1}, j_{2} \in \{1, 2, \ldots, q-1\}.
\end{equation}

Since  the average $(v_{\star})_{a} = q^{-1}\sum_{j=1}^{q} v_{\star}^{j}$ is linear (by property (${\mathcal B{\emph 3}}$)) and 
$\|v_{\star} - (v_{\star})_{a}\|_{L^{2}(B_{1})}^{-1}(v_{\star} - (v_{\star})_{a}) \in {\mathcal B}_{q}$ (by property (${\mathcal B{\emph 5I}}$)), where $v_{\star} - (v_{\star})_{a} = (v_{\star}^{1} - (v_{\star})_{a}, \ldots, v_{\star}^{q} - (v_{\star})_{a}),$ we may assume without loss of generality that $(v_{\star})_{a} = 0$ and that
\begin{equation}\label{partIII-0}
\|v_{\star}\|_{L^{2}(B_{1})} = 1.
\end{equation}

By the definition of ${\mathcal B}_{q}$, for each $k = 1, 2, 3, \ldots$, there exists a stationary integral 
varifold $V_{k}$ of $B_{2}^{n+1}(0)$ with
\begin{equation}\label{partIII-0-1}
\left(\omega_{n}2^{n}\right)^{-1}\|V_{k}\|(B_{2}^{n+1}(0)) < q + 1/2, \;\;\;\; q-1/2\leq \omega_{n}^{-1}\|V_{k}\|({\mathbf R} \times B_{1}) < q+1/2 \;\; {\rm and}
\end{equation}
\begin{equation}\label{partIII-0-2}
{\hat E}_{k}^{2} \equiv \int_{{\mathbf R} \times B_{1}} |x^{1}|^{2}d\|V_{k}\|(X) \to 0
\end{equation}
as $k \to \infty$, such that the following hold: If  
$\s \in (0, 1)$, $k$ is sufficiently large depending on $\s$, $\Sigma_{k} \subset B_{\s}$ is the measurable set corresponding to $\Sigma$ and $v_{k}^{j} \, : \, B_{\s} 
\to {\mathbf R},$ $j=1, 2, \ldots, q$, are the Lipschitz functions corresponding to $u^{j}$ in Theorem~\ref{flat-varifolds} applied with $V_{k}$ in place of $V$, then by Theorem~\ref{flat-varifolds}, 
$v_{k}^{1} \leq v_{k}^{2} \leq \ldots \leq v_{k}^{q},$ 
\begin{equation}\label{partIII-1}
{\rm Lip} \, v_{k}^{j} \leq 1/2 \;\;\; \mbox{for each} \;\;\;  j \in \{1, 2, \ldots, q\},
\end{equation}
\begin{equation}\label{partIII-2}
\|V_{k}\|({\mathbf R} \times \Sigma_{k}) + {\mathcal H}^{n} \, (\Sigma_{k}) \leq C_{\s}{\hat E}_{k}^{2}
\end{equation}
where $C_{\s} \in (0, \infty)$ is a constant depending only on $n$, $q$ and $\s$,  
\begin{equation}\label{partIII-2-0}
{\rm spt} \, \|V_{k}\| \cap ({\mathbf R} \times (B_{\s} \setminus \Sigma_{k})) = 
\cup_{j=1}^{q} {\rm graph} \, v_{k}^{j} \cap ({\mathbf R} \times (B_{\s} \setminus \Sigma_{k})), \;\; {\rm and}
\end{equation}
\begin{equation}\label{partIII-2-1}
{\hat E}_{k}^{-1}v_{k}^{j} \to v_{\star}^{j}
\end{equation}
where the convergence is in $L^{2}(B_{\s})$ and 
weakly in $W^{1, 2} \, (B_{\s}).$ Note that by (\ref{partIII-0-1}), after passing to a subsequence without changing notation, there exists a stationary integral varifold $V$ of $B_{2}^{n+1}(0)$ such that $V_{k} \to V,$  
and by (\ref{partIII-0-2}), ${\rm spt} \, \|V \res ({\mathbf R}  \times B_{1})\| \subset \{0\} \times \overline{B}_{1}$, so by (\ref{partIII-0-1}) and the Constancy Theorem for stationary integral varifolds, 
$V\res ({\mathbf R} \times B_{1}) = q|\{0\} \times B_{1}|.$ Hence by replacing $V_{k}$ with 
$\eta_{0, 1/2 \, \#} \, V_{k},$  and noting that by homogeneity of $v_{\star}$, the coarse blow-up of $\{\eta_{0, 1/2 \, \#} \, V_{k}\}$ is still $v_{\star}$, we may assume that for all sufficiently large $k$,
\begin{equation}\label{partIII-3}
q - 1/4 \leq (\omega_{n}2^{n})^{-1}\|V_{k}\|(B_{2}^{n+1}(0)) < q + 1/4.
\end{equation}
By using the argument justifying the assertion (\ref{no-overlap-4-1-0}), we may pass to a subsequence without changing notation and find points $Z_{k} \in {\rm spt} \,\|V_{k}\|\cap B_{1}^{n+1}(0)$ with $\Theta \, (\|V_{k}\|, Z_{k}) \geq q$ and $Z_{k} \to 0.$ Replacing $V_{k}$ with $\eta_{Z_{k}, 1 - |Z_{k}| \, \#} \,V_{k}$, we may thus assume that 
\begin{equation}\label{partIII-4}
\Theta \, (\|V_{k}\|, 0) \geq q
\end{equation}
for each $k=1,2,3, \ldots,$ and in view of (\ref{partIII-3}), the monotonicity formula implies that the new $V_{k}$ satisfy (\ref{partIII-0-1}). We now argue that for each sufficiently large $k$, we must have that  
\begin{equation}\label{partIII-4-0-1}
\int_{{\mathbf R} \times B_{1}}|x^{1}|^{2} \, d\|V_{k}\|(X)  <\frac{3}{2} \inf_{\{P = \{x^{1} = \lambda x^{2}\}\}} \int_{{\mathbf R} \times B_{1}}{\rm dist}^{2} \, (X, P) \, d\|V_{k}\|(X).
\end{equation}
If this were false, then there is a subsequence $\{k^{\prime}\}$ of $\{k\}$ and corresponding to each $k^{\prime}$,
a number $\lambda_{k^{\prime}} \in {\mathbf R}$ such that, with $P_{k^{\prime}} = \{x^{1} = \lambda_{k^{\prime}} x^{2}\}$,
we have 
$$\int_{{\mathbf R} \times B_{1}}{\rm dist}^{2} \, (X, P_{k^{\prime}}) \, d\|V_{k^{\prime}}\|(X) \leq \frac{5}{6}{\hat E}^{2}_{k^{\prime}}$$
for all $k^{\prime}.$ In particular, for each $\s \in (1/2, 1)$ and sufficiently large $k^{\prime}$, 
\begin{equation}\label{partIII-4-0-2}
(1 + \lambda_{k^{\prime}}^{2})^{-1}\sum_{j=1}^{q}\int_{B_{\s} \setminus \Sigma_{k^{\prime}}} (v^{j}_{k^{\prime}}(x^{2},y) - \lambda_{k^{\prime}}x^{2})^{2} \, dx^{2}dy \leq \frac{5}{6}{\hat E}_{k^{\prime}}^{2},
\end{equation}
whence $(1 + \lambda_{k^{\prime}}^{2})^{-1}\lambda_{k^{\prime}}^{2}\int_{B_{1/2} \setminus \Sigma_{k^{\prime}}}|x^{2}|^{2} \, dx^{2}dy \leq \frac{11}{3}{\hat E}_{k^{\prime}}^{2};$ thus, $|\lambda_{k^{\prime}}| \leq C{\hat E}_{k^{\prime}}$ for all sufficiently large $k^{\prime}$, where 
$C = C(n) \in (0, \infty),$ and hence, passing to a further subsequence without changing notation, ${\hat E}_{k^{\prime}}^{-1}\lambda_{k^{\prime}} \to \ell$ for some $\ell \in {\mathbf R}.$ It follows from (\ref{partIII-4-0-2}) and (\ref{partIII-2}) that 
$$\sum_{j=1}^{q}\int_{B_{\s}}(v_{k^{\prime}}^{j} - \lambda_{k^{\prime}} x^{2})^{2} \,dx^{2}dy\leq \frac{5}{6}(1 + \lambda_{k^{\prime}}^{2}){\hat E}_{k^{\prime}}^{2} + 2C_{\s}\sup_{B_{\s}} \, (|v_{k^{\prime}}|^{2} + \lambda_{k^{\prime}}^{2}|x^{2}|^{2}) \, {\hat E}_{k^{\prime}}^{2}.$$
First dividing this by ${\hat E}_{k^{\prime}}^{2}$ and letting $k^{\prime} \to \infty$, and then letting $\s \to 1$, we see that $\sum_{j=1}^{q}\int_{B_{1}}(v_{\star}^{j} - \ell x^{2})^{2} \leq 5/6.$ Since $v_{\star}^{j}(x^{2}, y) = \ell_{j}x^{2}$ if 
$x^{2} < 0$ and $v_{\star}^{j}(x^{2}, y) = m_{j} x^{2}$ if $x^{2} > 0$ for some $\ell_{j}, m_{j} \in {\mathbf R}$, 
this implies that $\int_{B_{1}}|v_{\star}|^{2} - 2\ell\sum_{j=1}^{q} (\ell_{j} + m_{j}) \int_{B_{1} \cap \{x^{2} > 0\}} |x^{2}|^{2} + \ell^{2} \int_{B_{1}}|x^{2}|^{2} \leq 5/6,$ which is a contradiction since $(v_{\star})_{a} \equiv 0$  (so that $\sum_{j=1}^{q} \ell_{j} = \sum_{j=1}^{q} m_{j} = 0$) and $\int_{B_{1}} |v_{\star}|^{2} =1.$ Thus (\ref{partIII-4-0-1}) must hold for all sufficiently large $k$.

For $j=1, 2, 3, \ldots, q$ and $k=1,2, 3, \ldots,$ let $h_{j}^{k} = {\hat E}_{k}L_{1}^{j}$, $g_{j}^{k} = {\hat E}_{k} L_{2}^{j}$, $H_{j}^{k} = {\rm graph} \, h_{j}^{k}$, $G_{j}^{k} = {\rm graph} \, g_{j}^{k}$ and 
${\mathbf C}_{k} = \sum_{j=1}^{q} |H_{j}^{k}| + |G_{j}^{k}|.$ By (\ref{partIII-1}), (\ref{partIII-2}) and (\ref{partIII-2-0}), 
\begin{equation}\label{partIII-4-1}
\int_{{\mathbf R} \times B_{\s}} {\rm dist}^{2} \, (X, {\rm spt} \, \|{\mathbf C}_{k}\|) \, d\|V_{k}\|(X) \leq 2\int_{B_{\s}} |v_{k} - {\hat E}_{k}v_{\star}|^{2} + C_{\s} \sup_{X \in {\rm spt}  \|V_{k}\| 
\cap ({\mathbf R} \times B_{\s})} \, {\rm dist}^{2} \,(X, {\rm spt} \, \|{\mathbf C}_{k}\|){\hat E}_{k}^{2}.
\end{equation}
By (\ref{partIII-0}) and homogeneity of $v_{\star}$, $\int_{B_{\s}} |v^{\star}|^{2} = \s^{n+2}$, so by (\ref{partIII-2-1}), 
for each $\th \in (0, 1/8)$ and $\s \in (0, 1)$, $\int_{B_{\s}}|v_{k}|^{2} \geq (1-\th)\s^{n+2}{\hat E}_{k}^{2}$ for sufficiently large $k$. Since 
\begin{eqnarray*}
\int_{{\mathbf R} \times B_{\s}} |x^{1}|^{2} \, d\|V_{k}\|(X) = \sum_{j=1}^{q}\int_{B_{\s}} \sqrt{1 + |Dv_{k}^{j}|^{2}}|v_{k}^{j}|^{2}  - \sum_{j=1}^{q}\int_{\Sigma_{k}} \sqrt{1 + |Dv_{k}^{j}|^{2}}|v_{k}^{j}|^{2}&&\nonumber\\
&&\hspace{-1in} + \int_{{\mathbf R} \times \Sigma_{k}} |x^{1}|^{2} \, d\|V_{k}\|(X)\nonumber\\
&&\hspace{-4in}\geq \int_{B_{\s}}|v_{k}|^{2} - 2C_{\s}\left(\sup_{B_{\s}} \,|v_{k}|^{2}\right){\hat E}_{k}^{2},
\end{eqnarray*}
it follows that 
\begin{equation}\label{partIII-4-2}
\int_{{\mathbf R} \times (B_{1} \setminus B_{\s})} |x^{1}|^{2} \, d\|V_{k}\|(X) \leq \left(1 - (1-\th)\s^{n+2} + 2C_{\s}\left(\sup_{B_{\s}} \,|v_{k}|^{2}\right)\right){\hat E}_{k}^{2}
\end{equation}
for all sufficiently large $k$. By the triangle inequality, 
\begin{eqnarray}\label{partIII-4-3}
\int_{{\mathbf R} \times (B_{1} \setminus B_{\s})} {\rm dist}^{2} \, (X, {\rm spt} \, \|{\mathbf C}_{k}\|) \, d\|V_{k}\|(X)&&\nonumber\\ 
&&\hspace{-2.75in}\leq 2\int_{{\mathbf R} \times (B_{1} \setminus B_{\s})} |x^{1}|^{2} \, d\|V_{k}\|(X) + 3\,{\rm dist}_{\mathcal H}^{2}({\rm spt} \, \|{\mathbf C}_{k}\| \cap ({\mathbf R} \times B_{1}), \{0\} \times B_{1})\|V_{k}\|({\mathbf R} \times (B_{1} \setminus B_{\s}))\nonumber\\
&&\hspace{-2in}\leq 2\int_{{\mathbf R} \times (B_{1} \setminus B_{\s})} |x^{1}|^{2} \, d\|V_{k}\|(X) + C{\mathcal H}^{n}(B_{1} \setminus B_{\s}){\hat E}_{k}^{2}
\end{eqnarray}
for all sufficiently large $k$, where $C = C(n, q) \in (0, \infty).$ Here we have used the fact that 
$V_{k} \res ({\mathbf R} \times B_{1}) \to q|\{0\} \times B_{1}|.$ Thus, if $\g_{1} = \g_{1}(n, q, \a) \in (0, 1/2)$ is the constant as in 
Lemma~\ref{finedistgap}, then we may fix $\th  = \th(n, q, \a) \in (0, 1/8)$ sufficiently small and 
$\s = \s(n, q, \a) \in (0, 1)$ sufficiently close to 1 in order to conclude from (\ref{partIII-2-1}), (\ref{partIII-4-1}), (\ref{partIII-4-2}) and (\ref{partIII-4-3}) that for all sufficiently large $k$, 
\begin{equation}\label{partIII-4-4}
\int_{{\mathbf R} \times B_{1}} {\rm dist}^{2} \, (X, {\rm spt} \, \|{\mathbf C}_{k}\|) \, d\|V_{k}\|(X) \leq \frac{\g_{1}}{4}{\hat E}_{k}^{2}.
\end{equation}

In view of (\ref{partIII-0-0}), we have by the argument leading to (\ref{no-overlap-3}) that for all sufficiently large $k$,   
$\Sigma_{k} \subset B_{\s} \cap \{|x^{2}| < 1/64\}$ and that
\begin{equation}\label{partIII-5}
 V_{k} \res(({\mathbf R} \times B_{\s}) \cap \{x^{2} \leq -1/64\}) =   \sum_{j=1}^{q} |{\rm graph} \, u_{k}^{j}| \res 
(({\mathbf R} \times B_{\s}) \cap \{x^{2} \leq -1/64\}) \;\;\; {\rm and}
\end{equation}
\begin{equation}\label{partIII-6}
 V_{k} \res(({\mathbf R} \times B_{\s}) \cap \{x^{2} \geq 1/64\}) =   \sum_{j=1}^{q} |{\rm graph} \, w_{k}^{j}| \res 
(({\mathbf R} \times B_{\s}) \cap \{x^{2} \geq 1/64\}) 
\end{equation}
where $u_{k}^{1} \leq u_{k}^{2} \leq \ldots u_{k}^{q}$ and $w_{k}^{1} \leq w_{k}^{2} \leq \ldots w_{k}^{q}$ (thus, $\left.v_{k}\right|_{B_{\s} \cap \{x^{2} \leq -1/64\}} \equiv u_{k}$ and $\left.v_{k}\right|_{B_{\s} \cap \{x^{2} \geq 1/64\}} \equiv w_{k}$), $u_{k}^{j}, w_{k}^{j}$ are $C^{2}$ functions on $B_{\s} \cap \{x^{2} \leq -1/64\},$  
$B_{\s} \cap \{x^{2} \geq 1/64\}$ respectively, solving the minimal surface equation there, and satisfying, by elliptic theory, 
\begin{equation}\label{partIII-7}
\sup_{B_{\k\s}\cap \{x^{2} \leq -1/64\}} |D\, u_{k}^{j}|^{2} + 
\sup_{B_{\k\s}\cap \{x^{2} \geq 1/64\}} |D\, w_{k}^{j}|^{2} \leq C(\k, \s){\hat E}_{k}^{2}
\end{equation}
for each $\k \in (0, 1)$, $j=1,2, \ldots, q$ where $C(\k, \s) \in (0, \infty)$ is a constant depending only on $n$, $\k$ and $\s$. We see from (\ref{partIII-5}), (\ref{partIII-6}), (\ref{partIII-7}) and (\ref{partIII-2-1}) that 
\begin{eqnarray}\label{partIII-8}
\int_{{\mathbf R} \times (B_{1/2} \setminus \{|x^{2}| < 1/64\})} {\rm dist}^{2} \, (X, {\rm spt} \, \|V_{k}\|) \, d\|{\mathbf C}_{k}\|(X)&&\nonumber\\
&&\hspace{-3in} \leq 2\sum_{j=1}^{q}\left(\int_{B_{1/2} \cap \{x^{2} < -1/64\}} |{\hat E}_{k}L_{1}^{j} - u_{k}^{j}|^{2} 
+ \int_{B_{1/2} \cap \{x^{2} > 1/64\}} |{\hat E}_{k} L_{2}^{j} - w_{k}^{j}|^{2}\right) \leq \eta_{k} {\hat E}_{k}^{2}
\end{eqnarray} 
where $\eta_{k} \to 0.$ By (\ref{partIII-4-4}) and (\ref{partIII-8}), 
\begin{equation*}
\int_{{\mathbf R} \times (B_{1/2} \setminus \{|x^{2}| < 1/64\})} {\rm dist}^{2} \, (X, {\rm spt} \, \|V_{k}\|) \, d\|{\mathbf C}_{k}\|(X) + \int_{{\mathbf R} \times B_{1}} {\rm dist}^{2} \, (X, {\rm spt} \, \|{\mathbf C}_{k}\|) \, d\|V_{k}\|(X) \leq \frac{\g_{1}}{2}{\hat E}_{k}^{2}
\end{equation*} 
for sufficiently large $k$, which in view of (\ref{partIII-0-1}), (\ref{partIII-0-2}) and (\ref{partIII-4}) contradicts 
Lemma~\ref{finedistgap}.\end{proof}

\begin{theorem}\label{blowups}
Let $q$ be an integer $\geq 2$, $\a \in (0, 1)$  and suppose that the induction hypotheses $(H1),$ $(H2)$ hold. Let ${\mathcal B}_{q}$ be the class of functions defined in Section~\ref{blow-up}. (Thus, each $v \in {\mathcal B}_{q}$ is  a coarse blow-up, in the sense described in Section~\ref{blow-up}, of a sequences of varifolds in ${\mathcal S}_{\a}$ converging weakly, in ${\mathbf R} \times B_{1}$, to $q|\{0\} \times B_{1}|$). If $v = (v^{1},v^{2}, \ldots, v^{q}) \in {\mathcal B}_{q}$, then $v^{j}$ is harmonic in $B_{1}$ for each $j=1, 2, \ldots, q;$ furthermore, if 
$\{V_{k}\} \subset {\mathcal S}_{\a}$ is a sequence whose coarse blow-up is $v$, and if 
for each of infinitely many values of $k$, there is a point $Z_{k} \in {\rm spt} \,\|V_{k}\| \cap (B_{3/4} \times {\mathbf R})$ with $\Theta \, (\|V_{k}\|, Z_{k}) \geq q$, then $v^{1} = v^{2} = \ldots = v^{q}.$
\end{theorem}

\begin{proof}
By the discussion of Section~\ref{step2} and Corollary~\ref{F7}, ${\mathcal B}_{q}$  is a proper blow-up class for a constant $C = C(n, q) \in (0, \infty).$ The present theorem follows from Theorem~\ref{blowup-reg} and the remark at the end of Section~\ref{step2}.
\end{proof}

\section{The Sheeting Theorem}\label{sheeting}
\setcounter{equation}{0}
This section is devoted to the proof of the Sheeting Theorem (Theorem $3.3^{\prime}$) subject to the induction 
hypotheses $(H1),$ $(H2)$.

\begin{lemma}\label{excess-s}
Let $q$ be an integer $\geq 2$, $\a \in (0, 1)$ and $\th \in (0, 1/4).$ Suppose that the induction hypotheses $(H1)$, $(H2)$ hold. There exists a number $\b_{0} = \b_{0}(n, q, \a, \th) \in (0, 1/2)$ such that if $V \in {\mathcal S}_{\a}$, $(\omega_{n}2^{n})^{-1}\|V\|(B_{2}^{n+1}(0)) < q + 1/2$, $q-1/2 \leq (\omega_{n})^{-1}\|V\|(B_{1} \times {\mathbf R}) < q + 1/2$, and $\int_{{\mathbf R} \times B_{1}} {\rm dist}^{2} \, (X, P) \, d\|V\|(X) < \b_{0}$
for some affine hyperplane $P$ of ${\mathbf R}^{n+1}$ with ${\rm dist}^{2}_{\mathcal H} \,(P \cap (B_{1} \times {\mathbf R}), B_{1} \times \{0\}) < \b_{0},$ then the following hold:
\begin{itemize}
\item[(a)] Either  $V \res(B_{1/2} \times {\mathbf R}) = \sum_{j=1}^{q} |{\rm graph} \, u_{j}|$
where $u_{j} \in C^{2} \,(B_{1/2}; {\mathbf R})$ for $j=1, 2, \ldots, q;$ 
$u_{1} \leq u_{2} \leq \ldots \leq u_{q}$ on $B_{1/2};$ $u_{j_{0}} < u_{j_{0}+1}$ on $B_{1/2}$ for some $j_{0} \in \{1,2, \ldots, q-1\}$ and, for each $j \in \{1,2, \ldots, q\}$, 
$${\rm sup}_{B_{1/2}} \, |u_{j} - p|^{2} + |D \, u_{j} - D \, p|^{2} + |D^{2} \, u_{j}|^{2} \leq C\int_{{\mathbf R} \times B_{1}} {\rm dist}^{2} \, (X, P) \, d\|V\|(X)$$
where $C = C(n, q) \in (0, \infty)$ and $p \, : \, {\mathbf R}^{n} \to {\mathbf R}$ is the affine function such that ${\rm graph} \,  p = P;$ 
or, there exists an  affine hyperplane $P^{\prime}$ with 
$${\rm dist}_{\mathcal H}^{2} \, \left( P^{\prime} \cap ({\mathbf R} \times B_{1}), P \cap ({\mathbf R} \times B_{1})\right) \leq C_{1}\int_{{\mathbf R} \times B_{1}} {\rm dist}^{2} \,(X, P) \,  d\|V\|(X) \;\;\; and$$
\begin{eqnarray*}
\th^{-n-2}\int_{{\mathbf R} \times B_{\th}} {\rm dist}^{2} \, (X, P^{\prime}) \, d\|V\|(X) \leq C_{2}\th^{2}\int_{{\mathbf R} \times B_{1}} {\rm dist}^{2} \,(X, P) \,  d\|V\|(X)
 \end{eqnarray*}
where $C_{1} = C_{1}(n, q) \in (0, \infty)$ and $C_{2} = C_{2}(n, q) \in (0, \infty).$
\item[(b)] $\left(\omega_{n}(2\th)^{n}\right)^{-1}\|V\|(B_{2\th}^{n+1}(0)) < q + 1/2$ and $q-1/2 \leq (\omega_{n}\th^{n})^{-1} \|V\|({\mathbf R} \times B_{\th}) < q + 1/2.$
\end{itemize}
\end{lemma}

\begin{proof}
For each $k=1,2,3,\ldots,$ let $V_{k} \in {\mathcal S}_{\a}$ be such that 
\begin{equation}\label{excess-s-0}
(\omega_{n}2^{n})^{-1}\|V_{k}\|(B_{2}^{n+1}(0)) < q + 1/2 \;\; {\rm and} \;\; q-1/2 \leq (\omega_{n})^{-1} 
\|V_{k}\|({\mathbf R} \times B_{1}) < q + 1/2,
\end{equation}
and let $P_{k}$ be an affine hyperplane of ${\mathbf R}^{n+1}$ such that 
\begin{equation}\label{excess-s-1}
{\rm dist}_{\mathcal H} \, (P_{k} \cap ({\mathbf R} \times B_{1}), \{0\} \times B_{1}) \to 0 \;\;\; {\rm and}
\end{equation}
\begin{equation}\label{excess-s-2}
\int_{{\mathbf R} \times B_{1}} {\rm dist}^{2} \, (X, P_{k}) \,  d\|V_{k}\|(X) \to 0.
\end{equation}
The lemma will be established by proving that for each of infinitely many $k$, the conclusions hold with $V_{k}$ in place of $V$, $P_{k}$ in place of $P$ and with fixed constants $C  = C(n, q)$, $C_{1} = C_{1}(n, q)$, $C_{2} = C_{2}(n, q)  \in (0, \infty).$

By (\ref{excess-s-1}), (\ref{excess-s-2}) and the triangle inequality, 
${\hat E}_{k} \equiv \sqrt{\int_{{\mathbf R} \times B_{1}}|x^{1}|^{2} \, d\|V_{k}\|(X)} \to 0$, and hence, by (\ref{excess-s-0}) and the Constancy Theorem, $V_{k} \res ({\mathbf R} \times B_{1}) \to q|\{0\} \times B_{1}|$, so that 
$$q - 1/2 \leq \left(\omega_{n}\th^{n}\right)^{-1}\|V_{k}\|({\mathbf R} \times B_{\th}) < q + 1/2$$ 
for sufficiently large $k$. Furthermore, by monotonicity of mass ratio, $$\left(\omega_{n}(2\th)^{n}\right)^{-1}\|V_{k}\|(B_{2\th}^{n+1}(0)) \leq 
\left(\omega_{n}2^{n}\right)^{-1}\|V_{k}\|(B_{2}^{n+1}(0)) < q + 1/2.$$  
Thus, conclusion (b) with $V_{k}$ in place of $V$ holds for sufficiently large $k$.

For each $k=1,2, 3, \ldots, $ there exists, by  (\ref{excess-s-1}), a rigid motion  $\G_{k}$ of 
${\mathbf R}^{n+1}$ with $\G_{k} \to {\rm Identity}$ such that $\G_{k}(P_{k}) = \{0\} \times {\mathbf R}^{n}.$
Let $\widetilde{V}_{k} = \eta_{9/10 \, \#} \, \G_{k \, \#} \, V_{k}$. Then by (\ref{excess-s-0}),  $\left(\omega_{n}2^{n}\right)^{-1}\|\widetilde{V}_{k}\|(B_{2}^{n+1}(0)) < q+ 1/2$, and by (\ref{excess-s-2}),
\begin{equation}\label{excess-s-4}
\int_{{\mathbf R} \times B_{19/18}} |x^{1}|^{2} \,d\|\widetilde{V}_{k}\|(X) \leq \left(\frac{9}{10}\right)^{-n-2} \int_{{\mathbf R} \times B_{1}} {\rm dist}^{2} \,(X, P_{k}) \, d\|V_{k}\|(X) \to 0.
\end{equation}
It follows again by the Constancy Theorem, for all sufficiently large $k,$
$$q-1/2 \leq (\omega_{n})^{-1} \|\widetilde{V}_{k}\|({\mathbf R} \times B_{1}) < q + 1/2.$$

Let $\widetilde{v} = (\widetilde{v}^{1}, \widetilde{v}^{2}, \ldots, \widetilde{v}^{q})  \in {\mathcal B}_{q}$ be the coarse blow-up of (a subsequence) of $\widetilde{V}_{k}$ by the coarse excess ${\hat E}_{\widetilde{V}_{k}} \equiv \sqrt{\int_{{\mathbf R} \times B_{1}} |x^{1}|^{2} \, d\|\widetilde{V}_{k}\|(X)}$. Suppose first that the $\widetilde{v}^{j}$'s are not all identical to one another. Then by Theorem~\ref{blowups}, for all sufficiently large $k$,
$$Z \in {\rm spt} \, \|\widetilde{V}_{k}\| \cap ({\mathbf R} \times B_{3/4}) \;\;\; \implies \;\;\; \Theta \, (\|\widetilde{V}_{k}\|, Z) < q;$$  
hence, by Remark 3 of Section~\ref{outline}, we may apply Theorem~\ref{SS} followed by elliptic theory to $\widetilde{V}_{k}$ and conclude, after transforming by $\G_{k}^{-1} \circ \eta_{9/10}^{-1}$, that 
$$V_{k} \res ({\mathbf R} \times B_{1/2}) = \sum_{j=1}^{q} |{\rm graph} \, u_{k}^{j}|$$
for all sufficiently large $k,$  where $u_{k}^{j} \in C^{2}(B_{1/2}; {\mathbf R}),$ $u_{k}^{1} \leq u_{k}^{2} \leq \ldots u_{k}^{q};$ $u_{k}^{j_{0}} < u_{k}^{j_{0}+1}$ on $B_{1/2}$ for some $j_{0} \in \{1, 2, \ldots, q-1\}$ and, for each $j \in \{1, 2, \ldots, q\}$,  
$${\rm sup}_{B_{1/2}} \, |u_{k}^{j} - p_{k}|^{2} + |D \, u_{k}^{j} - D \, p_{k}|^{2} + |D^{2} \, u_{k}^{j}|^{2} \leq C\int_{{\mathbf R} \times B_{1}} {\rm dist}^{2} \, (X, P_{k}) \, d\|V_{k}\|(X)$$
where $C = C(n, q) \in (0, \infty)$ and $p_{k} \, : \, {\mathbf R}^{n} \to {\mathbf R}$ is the affine function such that ${\rm graph} \,  p_{k} = P_{k}.$

If on the other hand $\widetilde{v}^{1}=\widetilde{v}^{2}= \ldots =\widetilde{v}^{q}$ ($=\widetilde{v}$, say) on $B_{1}$, then letting 
$\widetilde{p}(x) = \widetilde{v}(0) + D\widetilde{v}(0) \cdot x$ and $\widetilde{P}_{k} = {\rm graph} \, {\hat E}_{\widetilde{V}_{k}}\widetilde{p}$, it follows from Theorem~\ref{flat-varifolds} and the standard estimates for harmonic functions that 
\begin{equation}\label{excess-s-5}
{\rm dist}_{\mathcal H} \, (\widetilde{P}_{k} \cap ({\mathbf R} \times B_{1}), \{0\} \times B_{1}) \leq C{\hat E}_{\widetilde{V}_{k}} \;\;\; {\rm and}
\end{equation}
\begin{equation}\label{excess-s-6}
\th^{-n-2}\int_{{\mathbf R} \times B_{2\th}} {\rm dist}^{2} \, (X, \widetilde{P}_{k}) \, d\|\widetilde{V}_{k}\|(X) 
\leq C\th^{2}{ \hat E}_{\widetilde{V}_{k}}^{2}
\end{equation}
for all sufficiently large $k$, where $C = C(n, q) \in (0, \infty).$ Setting $P^{\prime}_{k} =  \eta_{9/10}^{-1}\G_{k}^{-1} \, \widetilde{P}_{k}$, it follows readily from (\ref{excess-s-5}), (\ref{excess-s-6}) and (\ref{excess-s-4}) that 
$${\rm dist}_{\mathcal H} \, (P^{\prime}_{k} \cap ({\mathbf R} \times B_{1}), P_{k} \cap ({\mathbf R} \times B_{1})) \leq C\int_{{\mathbf R} \times B_{1}} {\rm dist}^{2} \, (X, P_{k}) \, d\|V_{k}\|(X) \;\; {\rm and}$$
$$\th^{-n-2}\int_{{\mathbf R} \times B_{\th}} {\rm dist}^{2} \, (X, P^{\prime}_{k}) \, d\|V_{k}\|(X) 
\leq C\th^{2}\int_{{\mathbf R} \times B_{1}} {\rm dist}^{2} \, (X, P_{k}) \, d\|V_{k}\|(X)$$
for all sufficiently large $k,$ where $C = C(n, q) \in (0, \infty).$ Thus, conclusion (a)  with $V_{k}$, $P_{k}$, $P^{\prime}_{k}$ in place of $V$, $P$,  $P^{\prime}$ and with a fixed constant $C = C(n, q) \in(0, \infty)$ holds for infinitely many $k.$ 
\end{proof}

\begin{theorem}\label{sheetingthm}
Let $q$ be an integer $\geq 2,$ $\a \in (0, 1),$ $\g \in (0, 1)$ and suppose that the induction hypotheses $(H1)$, $(H2)$ hold. There exists a number $\e = \e(n, q, \a, \g) \in (0, 1)$ such that the following is true: If $V \in {\mathcal S}_{\a}$, 
$(\omega_{n}2^{n})^{-1}\|V\|(B_{2}^{n+1}(0)) < q + 1/2$, $q-1/2 \leq \omega_{n}^{-1}\|V\| (B_{1} \times {\mathbf R}) < q + 1/2$ and 
${\hat E}_{V}^{2} \equiv \int_{{\mathbf R} \times B_{1}}|x^{1}|^{2} \, d\|V\|(X) < \e,$ then 
$$V \res (B_{\g/2} \times {\mathbf R}) = \sum_{j=1}^{q} \, |{\rm graph} \, u_{j}|$$
where $u_{j} \in C^{1, \lambda}(B_{\g/2})$ for each $j=1, 2, \ldots, q$,  $u_{1} \leq u_{2} \leq  \ldots \leq  u_{q}$
and 
$$\sup_{B_{\g/2}} \, \left(|u_{j}| + |Du_{j}|\right) + \sup_{Y_{1}, Y_{2} \in B_{\g/2}, \, Y_{1} \neq Y_{2}} \, \frac{|Du_{j}(Y_{1}) - u_{j}(Y_{2})|}{|Y_{1} - Y_{2}|^{\lambda}} \leq C\left(\int_{{\mathbf R} \times B_{1}} |x^{1}|^{2} \, d\|V\|(X)\right)^{1/2}.$$
Here $C=C(n, q, \a, \g) \in (0, \infty)$ and $\lambda  = \lambda(n, q, \a, \g) \in (0, 1).$ Furthermore, we have in fact that $u_{j} \in C^{\infty}(B_{\g/2})$ and $u_{j}$ solves the minimal surface equation on $B_{\g/2}$ for each $j=1, 2, \ldots, q.$ 
\end{theorem}

\begin{proof}
Let $\widetilde{\g} = (1-\g)/4.$ Let $C = C(n, q)$, $C_{1} = C_{1}(n, q)$ and $C_{2} = C_{2}(n, q)$ be the constants as in the conclusion of Lemma~\ref{excess-s}. Choose $\th = \th(n, q) \in (0, 1/4)$ such that $C_{2}\th^{2} < 1/4$ and  
$\e = \e(n, q, \a, \g) \in (0, 1)$ such that 
$\e < (1 + C_{1})^{-1}\widetilde{\g}^{n+2}\b_{0}/8$ where $\b_{0} = \b_{0}(n, q, \a, \th)$ is as in Lemma~\ref{excess-s}. Additional restrictions on $\e$  will be imposed during the course of the proof, but we will choose  $\e$ depending only on $n$, $q$, $\a$ and $\g$. Suppose that ${\hat E}_{V}^{2} \equiv \int_{{\mathbf R} \times B_{1}} |x^{1}|^{2} \, d\|V\|(X) \leq \e,$ and let 
\begin{eqnarray*}
\b = \min \, \left\{4^{-1}\left(1 + 2C_{1}\right)^{-1}\b_{0}, \;\; 4^{-1}\widetilde{\g}^{n}\omega^{-1}_{n}(2q+1)^{-1}\b_{0}, \;\; 4^{-1}\left(2 + \omega_{n}(2q+1)C_{1}\right)^{-1}\left(\frac{2\th}{3}\right)^{n+2}\overline{\e}_{0},\right.&&\nonumber\\ 
&&\hspace{-4in}\left.8^{-1}\omega_{n}4^{-n}\th^{n}\left(256\th^{-2} + (q+1)\overline{C}\right)^{-1}\left(2 + \omega_{n}(2q+1)C_{1}\right)^{-1}\right\}
\end{eqnarray*}
Here $\overline{\e}_{0} = \e_{0}(n, q, 5/6)$, $\overline{C} = C(n, q, 5/6)$, where $\e_{0} = \e_{0}(n, q, \cdot)$ is as in Theorem~\ref{flat-varifolds} and $C = C(n, q, \cdot)$ is as in Theorem~\ref{flat-varifolds}(a). Note that 
$\b$ depends only on $n$, $q$, $\a$ and $\g$. Let $P_{0}$ be any affine hyperplane such that 
\begin{equation}\label{sheeting-0}
{\rm dist}^{2}_{\mathcal H} \, (P_{0} \,\cap \,({\mathbf R} \times B_{1}), \{0\} \times B_{1}) < \b.
\end{equation}
Fix any point $Y \in B_{\g}(0)$ and let $\widetilde{V} = \eta_{Y, \widetilde{\g}\, \#} \, V.$ Note then that
\begin{eqnarray}\label{sheeting-1}
{\hat E}_{\widetilde{V}, \, P_{0}}^{2} \equiv \int_{{\mathbf R} \times B_{1}} {\rm dist}^{2} \, (X, P_{0}) \, d\|\widetilde{V}\|(X) = 
\widetilde{\g}^{-n-2} \int_{{\mathbf R} \times B_{\widetilde{\g}}(Y)} {\rm dist}^{2} \, (X, Y + \widetilde{\g}P_{0}) \,d\|V\|(X)&&\nonumber\\
&&\hspace{-5.25in}\leq 2\widetilde{\g}^{-n-2} \int_{{\mathbf R} \times B_{\widetilde{\g}}(Y)} |x^{1}|^{2}\,d\|V\|(X)\nonumber\\
&& \hspace{-4.25in}+ 2\widetilde{\g}^{-n-2}\|V\|({\mathbf R} \times B_{\widetilde{\g}}(Y)){\rm dist}^{2}_{\mathcal H} \, (Y+ \widetilde{\g}P_{0} \, \cap \, ({\mathbf R} \times B_{\widetilde{\g}}(Y)), 
\{0\} \times B_{\widetilde{\g}}(Y))\nonumber\\
&&\hspace{-5.25in} \leq 2\widetilde{\g}^{-n-2}{\hat E}_{V}^{2} + \widetilde{\g}^{-n}\omega_{n}(2q + 1){\rm dist}^{2}_{\mathcal H} \ (P_{0} \cap ({\mathbf R} \times B_{1}), \{0\} \times B_{1})\leq 2\widetilde{\g}^{-n-2}\e +\b_{0}/4 < \b_{0}.
\end{eqnarray}
Furthermore, assuming $\e < \e_{0}\left(n, q, \frac{3 +\g}{4}\right)$ where $\e_{0}$ is as in Theorem~\ref{flat-varifolds} and applying Theorem~\ref{flat-varifolds}  with $\s = (3 + \g)/4$, we have that 
\begin{eqnarray}\label{sheeting-1-0}
\left|\left(\omega_{n}\widetilde{\g}^{n}\right)^{-1}\|V\|({\mathbf R} \times B_{\widetilde{\g}}(Y)) - q\right| =\left| (\omega_{n}\widetilde{\g}^{n})^{-1}\sum_{j=1}^{q}
\int_{B_{\widetilde{\g}} \setminus \Sigma} \left(\sqrt{1 + |Du^{j}|^{2}} - 1\right) \, dx\right.&&\nonumber\\
&&\hspace{-3in}  - \; \left.\left(\omega_{n}\widetilde{\g}^{n}\right)^{-1}\left(q{\mathcal H}^{n}(B_{\widetilde{\g}} \cap \Sigma) - \|V\|({\mathbf R} \times (B_{\widetilde{\g}} \cap \Sigma))\right)\right|\nonumber\\
&&\hspace{-5in} \leq \left(\omega_{n}\widetilde{\g}^{n}\right)^{-1}\sum_{j=1}^{q}\int_{B_{\widetilde{\g}} \setminus \Sigma}\frac{|Du^{j}|^{2}}{\sqrt{1 + |Du^{j}|^{2}}} \, dx + \left(\omega_{n}\widetilde{\g}^{n}\right)^{-1}(q + 1)\widetilde{C}{\hat E}_{V}^{2}\nonumber\\ 
&&\hspace{-5in} \leq \left(\omega_{n}\widetilde{\g}^{n}\right)^{-1} \int_{{\mathbf R} \times B_{\widetilde{\g}}(Y)} |\nabla^{V} \, x^{1}|^{2}  \, d\|V\|(X) + \left(\omega_{n}\widetilde{\g}^{n}\right)^{-1}(q + 1)\widetilde{C}{\hat E}_{V}^{2}\nonumber\\
&&\hspace{-5in}\leq 16\widetilde{\g}^{-2}\left(\omega_{n}\widetilde{\g}^{n}\right)^{-1}\int_{{\mathbf R} \times B_{\widetilde{2\g}}(Y)}|x^{1}|^{2} \, d\|V\|(X) + \left(\omega_{n}\widetilde{\g}^{n}\right)^{-1}(q + 1)\widetilde{C}{\hat E}_{V}^{2}\nonumber\\
&&\hspace{-5in}\leq \left(\omega_{n}\widetilde{\g}^{n}\right)^{-1}(16\widetilde{\g}^{-2} + (q+1)\widetilde{C}){\hat E}_{V}^{2}.
\end{eqnarray}
Here $\widetilde{C} = C\left(n, q, \frac{3 + \g}{4}\right)$ where $C = C(n, q, \cdot)$ is as in Theorem~\ref{flat-varifolds}(a), and $u^{j}$, $\Sigma$ are as in Theorem~\ref{flat-varifolds}; we have also used the fact that $\int_{{\mathbf R} \times B_{\widetilde{\g}}(Y)}|\nabla^{V} \, x^{1}|^{2} \, d\|V\|(X) \leq 16\widetilde{\g}^{-2} \int_{{\mathbf R} \times B_{\widetilde{\g}}(Y)}|x^{1}|^{2} \, d\|V\|(X)$, which follows from (\ref{tilt-ht-0}). Thus if $\e = \e(n, q, \a, \g) \in (0, 1)$ is sufficiently small, this says that  
\begin{equation}\label{sheeting-2}
q -1/2 \leq (\omega_{n})^{-1}\|\widetilde{V}\|({\mathbf R} \times B_{1}) < q + 1/2.
\end{equation}
Since $(\omega_{n}2^{n})^{-1}\|\widetilde{V}\|(B_{2}^{n+1}(0))  = \left(\omega_{n}(2\widetilde{\g})^{n}\right)^{-1}\|V\|(B^{n+1}_{2\widetilde{\g}}(Y)) \leq \left(\omega_{n}(2\widetilde{\g})^{n}\right)^{-1}\|V\|({\mathbf R} \times B_{2\widetilde{\g}}(Y))$, the same estimate 
with $2\widetilde{\g}$ in place of $\widetilde{\g}$ shows that 
\begin{equation}\label{sheeting-2-1}
(\omega_{n}2^{n})^{-1}\|\widetilde{V}\|(B_{2}^{n+1}(0)) < q + 1/2
\end{equation}
provided $\e = \e(n, q, \a, \g) \in (0, 1)$ is sufficiently small.
 
We claim that either (I) or (II) below must hold: 
\begin{itemize}
\item[(I)] for each $k \in \{0, 1, 2, \ldots \}$, 
\begin{equation}\label{sheeting-3}
\left(\omega_{n}(2\th^{k})^{n}\right)^{-1}\|\widetilde{V}\|(B_{2\th^{k}}^{n+1}(0)) < q + 1/2; \; q-1/2 \leq (\omega_{n}(\th^{k})^{n})^{-1} \|\widetilde{V}\|({\mathbf R} \times B_{\th^{k}}) < q + 1/2;
\end{equation}
there exists an affine hyperplane $P_{k}$ such that, if $k \geq 1$,  
\begin{equation}\label{sheeting-4}
{\rm dist}_{\mathcal H}^{2} \, \left( P_{k} \cap ({\mathbf R} \times B_{1}), P_{k-1} \cap ({\mathbf R} \times B_{1})\right) \leq C_{1}4^{-k}{\hat E}_{\widetilde{V}, P_{0}}^{2} \;\; {\rm and}
\end{equation}
\begin{eqnarray}\label{sheeting-5}
(\th^{k})^{-n-2}\int_{{\mathbf R} \times B_{\th^{k}}} {\rm dist}^{2} \, (X, P_{k}) \, d\|\widetilde{V}\|(X)&&\nonumber\\
&&\hspace{-2in} \leq 4^{-1}(\th^{k-1})^{-n-2}\int_{{\mathbf R} \times B_{\th^{k-1}}} {\rm dist}^{2} \, (X, P_{k-1}) \, d\|\widetilde{V}\|(X) \leq \ldots \leq 
4^{-k}{\hat E}_{\widetilde{V}, P_{0}}^{2}. 
\end{eqnarray}
\item[(II)] there exists $\r_{0} \in (0, 1)$ such that $\widetilde{V} \res ({\mathbf R} \times B_{\r_{0}}) = 
\sum_{j=1}^{q} |{\rm graph} \, u_{j}|$ for 
functions $u_{j} \in C^{2} \, (B_{\r_{0}}; {\mathbf R}),$ $j= 1, 2, \ldots, q,$
satisfying $u_{1} \leq u_{2} \leq \ldots \leq u_{q}$ on $B_{\r_{0}};$ 
$u_{j_{0}} < u_{j_{0}+1}$ on $B_{\r_{0}}$ for some $j_{0} \in \{1, 2, \ldots, q-1\}$ and 
\begin{eqnarray}\label{sheeting-6}
{\rm sup}_{B_{\r_{0}/2}} \, \r_{0}^{-2}|u_{j}|^{2} + |D \, u_{j}|^{2} + \r_{0}^{2}|D^{2} \, u_{j}|^{2}&&\nonumber\\
&&\hspace{-2in} \leq (C+2C_{1}){\hat E}^{2}_{\widetilde{V}, P_{0}} + 4{\rm dist}^{2}_{\mathcal H} \, (P_{0} \cap ({\mathbf R} \times B_{1}),\{0\} \times B_{1})
\end{eqnarray}
for each $j \in \{1, 2, \ldots, q\};$ moreover, 
\begin{eqnarray}\label{sheeting-6-1}
\r^{-n}\int_{{\mathbf R} \times B_{\r}} |\nabla^{V} \, x^{1}|^{2} \, d\|V\|(X)&&\nonumber\\
&&\hspace{-1.5in} \leq C^{\prime}\left({\hat E}^{2}_{\widetilde{V}, P_{0}} + {\rm dist}^{2}_{\mathcal H} \, (P_{0} \cap ({\mathbf R} \times B_{1}),\{0\} \times B_{1})\right) \;\; {\rm for} \;\; \r_{0} < \r < \theta
\end{eqnarray}
where $C^{\prime} = C^{\prime}(n, q) \in (0, \infty)$. 
\end{itemize}
To see this, let $k_{0}$ be the smallest integer ($\geq 1$) such that alternative (I) fails to hold. If $k_{0}=1$, in view of (\ref{sheeting-1}) and (\ref{sheeting-2}), it follows directly from Lemma~\ref{excess-s} applied with $\widetilde{V}$ in place of $V$ and with $P = P_{0}$ 
that (II) must hold with $\r_{0} = 1/2$. Suppose $k_{0} \geq 2.$ Then by assumption, the inequalities (\ref{sheeting-3}), (\ref{sheeting-4}) and (\ref{sheeting-5}) hold for each $k = 1, 2, \ldots, k_{0}-1$ and consequently, by (\ref{sheeting-0}), (\ref{sheeting-1}),
(\ref{sheeting-4}) and the triangle inequality, 
\begin{eqnarray}\label{sheeting-7}
{\rm dist}^{2}_{\mathcal H} \, (P_{k_{0} - 1} \cap ({\mathbf R} \times B_{1}), \{0\} \times B_{1}) \leq  \left(\sqrt{C_{1}}{\hat E}_{\widetilde{V}, P_{0}} + {\rm dist}_{\mathcal H} \, (P_{0} \cap ({\mathbf R} \times B_{1}), \{0\} \times B_{1})\right)^{2}&&\nonumber\\
&&\hspace{-4in} \leq 4C_{1}\widetilde{\g}^{-n-2}\e +\b_{0}/2 < \b_{0}.
\end{eqnarray}
Applying Lemma~\ref{excess-s} with $\eta_{\th^{k_{0}-1} \, \#} \, \widetilde{V}$ in place of $V$ and $P_{k_{0}-1}$ in place of $P$, we see  by the defining property of $k_{0}$ that  $\widetilde{V} \res({\mathbf R} \times B_{\th^{k_{0}-1}/2}) = \sum_{j=1}^{q} |{\rm graph} \, u_{j}|$
where $u_{j} \in C^{2} \,(B_{\th^{k_{0}-1}/2}; {\mathbf R})$ for $j=1, 2, \ldots, q;$ 
$u_{1} \leq u_{2} \leq \ldots \leq u_{q}$ on $B_{\th^{k_{0}-1}/2};$ $u_{j} < u_{j+1}$ on $B_{\th^{k_{0}-1}/2}$ for some $j \in \{1,2, \ldots, q-1\}$ and
\begin{eqnarray}\label{sheeting-8}
{\rm sup}_{B_{\th^{k_{0}-1}/2}} \, (\th^{k_{0}-1})^{-2}|u_{j} - p|^{2} + |D \, u_{j} - D \, p|^{2} + (\th^{k_{0}-1})^{2}|D^{2} \, u_{j}|^{2}&&\nonumber\\
&&\hspace{-4in} \leq 
C(\th^{k_{0}-1})^{-n-2}\int_{{\mathbf R} \times B_{\th^{k_{0}-1}}} {\rm dist}^{2} \, (X, P_{k_{0}-1}) \, d\|\widetilde{V}\|(X)
\leq C4^{-(k_{0}-1)}{\hat E}_{\widetilde{V}, P_{0}}^{2}
\end{eqnarray}
for each $j \in \{1,2, \ldots, q\};$  here $p \, : \, {\mathbf R}^{n} \to {\mathbf R}$ is the affine function such that ${\rm graph} \,  p = P_{k_{0}-1}.$ In view of (\ref{sheeting-1}) and (\ref{sheeting-7}), this evidently implies alternative (II) with $\r_{0} = \th^{k_{0}-1}$. Thus either (I) or (II) holds as claimed.

Suppose that (I) holds for some $P_{0}$ satisfying (\ref{sheeting-0}). It is standard then that there exists a hyperplane $\widetilde{P}$ with 
\begin{equation}\label{sheeting-9}
{\rm dist}_{\mathcal H} \, (\widetilde{P} \cap ({\mathbf R} \times B_{1}), P_{0} \cap ({\mathbf R} \times B_{1})) \leq C_{1}{\hat E}_{\widetilde{V}, P_{0}}
\end{equation}
such that
\begin{equation}\label{sheeting-10}
\r^{-n-2}\int_{{\mathbf R} \times B_{\r}} {\rm dist}^{2} \, (X,\widetilde{P}) \, d\|\widetilde{V}\|(X) \leq C_{3}\r^{2\m}{\hat E}_{\widetilde{V},P_{0}}^{2}
\end{equation}
for each $\r \in (0, 1)$, where $C_{3} = C_{3}(n, q) \in (0, \infty)$ and $\m = \m(n, q) \in (0, 1).$ Note that $\widetilde{P}$ does not depend on $P_{0}$, nor do the constants $C_{3}$ and $\m$. Moreover, in this case, we claim that we have 
for each $\r \in (0, 1/4)$ that 
\begin{equation}\label{sheeting-10-3}
\left(\omega_{n}(2\r)^{n}\right)^{-1}\|\widetilde{V}\|(B_{2\r}^{n+1}(0)) < q + 1/2 \;\; {\rm and} \;\; q-1/2 \leq (\omega_{n}\r^{n})^{-1}\|\widetilde{V}\|({\mathbf R} \times B_{\r}) < q + 1/2.
\end{equation}
To see this, given $\r \in (0, 1/4)$, choose $k$ such that $\th^{k+1} \leq 4\r < \th^{k}$ and note by (\ref{sheeting-3}), (\ref{sheeting-4}), (\ref{sheeting-5}) and the triangle inequality that 
\begin{eqnarray}\label{sheeting-10-3-1}
&&\left(\th^{k}\right)^{-n-2}\int_{{\mathbf R} \times B_{\th^{k}}}|x^{1}|^{2} \, d\|\widetilde{V}\|(X)\nonumber\\ 
&&\hspace{.5in} \leq \; 2(2 + \omega_{n}(2q+1)C_{1})(\widetilde{\g}^{-n-2}{\hat E}_{V}^{2} + {\rm dist}^{2} \, (P_{0} \cap ({\mathbf R} \times B_{1}), \{0\} \times B_{1}))\nonumber\\ 
&&\hspace{1in}\leq 2(2 + \omega_{n}(2q+1)C_{1})(\widetilde{\g}^{-n-2}\e +\b)
\end{eqnarray}
so provided 
$\e  = \e(n, q, \a, \g) \in (0, 1)$ is sufficiently small, we may, in view of (\ref{sheeting-3}), apply Theorem~\ref{flat-varifolds} with $\eta_{\th^{k} \, \#} \, \widetilde{V}$ in place of $V$, $\s = 5/6$ and estimate exactly 
as in (\ref{sheeting-1-0}) (with $\eta_{\th^{k} \, \#} \, \widetilde{V}$ in place of $V$, $Y = 0$ and $\th^{-k}\r$ in place of $\widetilde{\g}$) to deduce that 
\begin{eqnarray*}
\left|\left(\omega_{n}\left(\th^{-k}\r\right)^{n}\right)^{-1}\|\eta_{\th^{k} \, \#} \, \widetilde{V}\|\left({\mathbf R} \times B_{\th^{-k}\r}(Y)\right) - q\right|&&\nonumber\\ 
&&\hspace{-3in}\leq 2\left(\omega_{n}\left(\th^{-k}\r\right)^{n}\right)^{-1}\left(16\left(\th^{-k}\r\right)^{-2} + (q+1)\overline{C}\right)\left(2 + \omega_{n}(2q+1)C_{1}\right)\left(\widetilde{\g}^{-n-2}\e +\b\right)\nonumber\\
&&\hspace{-3in}\leq 2\omega_{n}^{-1}4^n\th^{-n}\left(256\th^{-2} + (q+1)\overline{C}\right)\left(2 + \omega_{n}(2q+1)C_{1}\right)\left(\widetilde{\g}^{-n-2}\e +\b\right).
\end{eqnarray*}
(Recall that $\overline{C} = C(n, q,5/6)$ where $C = C(n, q, \cdot)$ is as in Theorem~\ref{flat-varifolds}(a)). From this, (\ref{sheeting-2-1}) and the monotonicity formula, we deduce that (\ref{sheeting-10-3}) holds provided 
$\e = \e(n, q, \a, \g) \in (0, 1)$ is sufficiently small. It then follows, if alternative (I) holds for some $P_{0}$ satisfying (\ref{sheeting-0}), that ${\rm spt} \, \|\widetilde{V}\| \cap \pi^{-1}(0)$ consists of a single point ($= \widetilde{P} \cap \pi^{-1}(0)$); to see this, first note that ${\rm spt} \, \|\widetilde{V}\|\cap \pi^{-1}(0) \neq \emptyset$ by the second inequality in (\ref{sheeting-10-3}). Let $Z \in {\rm spt} \, \|\widetilde{V}\| \cap \pi^{-1}(0)$ and ${\mathbf C}_{Z}$ be a tangent cone to $\widetilde{V}$ at $Z$. Thus 
 $\eta_{Z, \s_{j} \, \#} \, \widetilde{V} \to {\mathbf C}_{Z} \neq 0$ for some sequence of numbers $\s_{j} \to 0^{+},$  and by (\ref{sheeting-10}),  ${\rm dist}_{\mathcal H} \, ({\rm spt} \, \|\eta_{Z, \s_{j} \, \#} \, \widetilde{V}\| \cap ({\mathbf R} \times B_{1/2}), \s_{j}^{-1}(\widetilde{P} - Z) \cap ({\mathbf R} \times B_{1/2})) \to 0$, which can only be true if $Z \in \widetilde{P}.$ But by (\ref{sheeting-9}), $\widetilde{P} \cap \pi^{-1}(0)$ consists of a single point. 
 
 Since alternative (II) implies that ${\rm spt} \, \|\widetilde{V}\| \cap \pi^{-1}(0)$ has at least two distinct points, we see that if alternative (I) holds for some $P_{0}$ satisfying (\ref{sheeting-0}), then (I) must hold for all $P_{0}$ satisfying (\ref{sheeting-0}). Taking $P_{0} = {\mathbf R}^{n} \times \{0\},$ we deduce from (\ref{sheeting-9}) and 
(\ref{sheeting-10}) that
 \begin{equation}\label{sheeting-10-1}
{\rm dist}_{\mathcal H} \, (\widetilde{P} \cap ({\mathbf R} \times B_{1}),\{0\} \times  B_{1})) \leq C_{1}\widetilde{\g}^{-n-2}{\hat E}_{V}^{2}\leq C_{1}\widetilde{\g}^{-n-2}\e \;\;\;{\rm and}
\end{equation}
\begin{equation}\label{sheeting-10-2}
\r^{-n-2}\int_{{\mathbf R} \times B_{\r}} {\rm dist}^{2} \, (X,\widetilde{P}) \, d\|\widetilde{V}\|(X) \leq C_{3}\widetilde{\g}^{-n-2}{\hat E}_{V}^{2} \leq C_{3}\widetilde{\g}^{-n-2}\e
\end{equation}
for $\r \in (0, 1).$ So if we choose $\e = \e(n, q, \a, \g)$ such that $C_{1}\widetilde{\g}^{-n-2}\e < \b$ we may in particular take $P_{0} = \widetilde{P}$ in (\ref{sheeting-10}). 

Thus we have so far established the following: Given $q$, $\a$, $\g$ as in the theorem and that the induction hypotheses $(H1)$, $(H2)$ hold, there exists $\e = \e(n, q, \a, \g) \in (0, 1)$ such that if $V \in {\mathcal S}_{\a}$ satisfies the hypotheses of the theorem, $Y \in B_{\g}$, $\widetilde{V} = \eta_{Y, \widetilde{\g} \, \#} \, V$ where $\widetilde{\g} = (1 - \g)/4$, then either alternative (I) above holds for all affine hyperplanes $P_{0}$ satisfying (\ref{sheeting-0}) or 
alternative (II) above holds for all such $P_{0}$; furthermore, if alternative (I) holds, then the bounds (\ref{sheeting-10-3}) are satisfied for each $\r \in (0, 1/4),$  the estimates (\ref{sheeting-10-1}), (\ref{sheeting-10-2}) are satisfied and 
\begin{equation}\label{sheeting-10-5}
\r^{-n-2}\int_{{\mathbf R} \times B_{\r}} {\rm dist}^{2} \, (X,\widetilde{P}) \, d\|\widetilde{V}\|(X) \leq C_{3}\r^{2\m} \int_{{\mathbf R} \times B_{1}}{\rm dist}^{2} \, (X, \widetilde{P}) \, d\|\widetilde{V}\|(X) 
\end{equation}
for $\r \in (0, 1)$, where $\widetilde{P}$ is a (uniquely determined) affine hyperplane, $C_{3} = C_{3}(n, q) \in (0, \infty)$ and $\m = \m(n, q) \in (0, 1).$  

Now suppose the hypotheses of the theorem are satisfied with a number $\e^{\prime} = \e^{\prime}(n, q, \a, \g) \in (0, \e)$ in place of $\e$, and that alternative (I) (with $Y \in B_{\g}$, $\widetilde{V}  = \eta_{Y, \widetilde{\g} \, \#} \, V$ as above) still holds. Then for any $\r \in (0, 1/4)$ we have by (\ref{sheeting-10-1}), (\ref{sheeting-10-2}), (\ref{sheeting-10-3}) and the triangle inequality that 
\begin{equation*}
\r^{-n-2}\int_{{\mathbf R} \times B_{\r}} |x^{1}|^{2} \, d\|\widetilde{V}\|(X)\leq(2C_{3} + \omega_{n}(2q + 1)C_{1})\widetilde{\g}^{-n-2}{\hat E}_{V}^{2} < 2(C_{3} + C_{1})\e^{\prime}.
\end{equation*}
Thus if we choose $\e^{\prime} = \e^{\prime}(n, q, \a, \g)$ sufficiently small, we may, in view of this and (\ref{sheeting-10-3}), repeat the argument leading to (\ref{sheeting-10-5}) with $\eta_{\r \, \#} \, \widetilde{V}$ (for which alternative (I) must hold) in place of $\widetilde{V}$; by applying (\ref{sheeting-10-5}) with $\eta_{\r \, \#} \, \widetilde{V}$ in place of $\widetilde{V}$ and $\r^{-1}\s$ in place of $\r$, we deduce that if alternative (I) holds for some $P_{0}$ satisfying (\ref{sheeting-0}), then there exists a unique affine hyperplane $\widetilde{P}$ satisfying (\ref{sheeting-10-1}) and
\begin{equation}\label{sheeting-10-4}
\s^{-n-2}\int_{{\mathbf R} \times B_{\s}} {\rm dist}^{2} \, (X, \widetilde{P}) \, d\|\widetilde{V}\|(X) \leq C_{3}\left(\frac{\s}{\r}\right)^{2\m}
\r^{-n-2}\int_{{\mathbf R} \times B_{\r}}{\rm dist}^{2} \, (X, \widetilde{P}) \, d\|\widetilde{V}\|(X) 
\end{equation}
for each $0 < \s < \r < 1/4.$ If on the other hand (I) fails for some $P_{0}$ satisfying 
(\ref{sheeting-0}), then it fails with $P_{0} = \{0\}\times {\mathbf R}^{n}$ in which case (by (II)) there exists
$\r_{0} \in (0, 1)$ such that $\widetilde{V} \res ({\mathbf R} \times B_{\r_{0}}) = 
\sum_{j=1}^{q} |{\rm graph} \, u_{j}|$ for 
functions $u_{j} \in C^{2} \, (B_{\r_{0}}; {\mathbf R}),$ $j= 1, 2, \ldots, q$
satisfying $u_{1} \leq u_{2} \leq \ldots \leq u_{q}$ on $B_{\r_{0}};$ 
$u_{j} < u_{j+1}$ on $B_{\r_{0}}$ for some $j \in \{1, 2, \ldots, q-1\}$ and 
\begin{equation*}
{\rm sup}_{B_{\r_{0}/2}} \, \r_{0}^{-2}|u_{j}|^{2} + |D \, u_{j}|^{2} + \r_{0}^{2}|D^{2} \, u_{j}|\leq (C+2C_{1})\widetilde{\g}^{-n-2}{\hat E}_{V}^{2}\nonumber\\
\end{equation*}
for each $j \in \{1, 2, \ldots, q\}.$

Thus we have shown that if the hypotheses of the theorem are satisfied with sufficiently small $\e^{\prime} = \e^{\prime}(n, q, \g) \in (0, 1)$ in place of $\e$, then for each point $Y \in B_{\g},$ precisely one of the following alternatives 
(I$_{Y}$) and (II$_{Y}$) must hold:
\begin{itemize}
\item[(I$_{Y}$)] there exists an affine hyperplane $P_{Y}$ with
\begin{equation}\label{sheeting-11}
{\rm dist}_{\mathcal H}^{2} \, (P_{Y} \cap ({\mathbf R} \times B_{1}(Y)), \{0\} \times B_{1}(Y)) \leq C_{1}\widetilde{\g}^{-n-2}{\hat E}_{V}^{2}
\end{equation}
such that
\begin{equation}\label{sheeting-12}
\s^{-n-2}\int_{{\mathbf R} \times B_{\s}(Y)} {\rm dist}^{2} \, (X, P_{Y}) \, d\|V\|(X) \leq C_{3}\left(\frac{\s}{\r}\right)^{2\m}\r^{-n-2}\int_{{\mathbf R} \times B_{\r}(Y)} {\rm dist}^{2}(X, P_{Y}) \, d\|V\|(X) 
\end{equation}
for each $0 < \s < \r < \widetilde{\g}/4$, where $C_{3} = C_{3}(n, q) \in (0, \infty)$ and $\m= \m(n, q) \in (0, 1),$ or 
\item[(II$_{Y}$)] there exists 
$\r_{Y} \in (0, 1/2]$ such that $V \res ({\mathbf R} \times B_{\r_{Y}}(Y)) = 
\sum_{j=1}^{q} |{\rm graph} \, u_{j}^{Y}|$ for 
functions $u_{j}^{Y} \in C^{2} \, (B_{\r_{Y}}(Y); {\mathbf R}),$ $j= 1, 2, \ldots, q,$
satisfying $u_{1}^{Y} \leq u_{2}^{Y} \leq \ldots \leq u_{q}^{Y}$ on $B_{\r_{Y}}(Y);$ 
$u_{j_{0}}^{Y} < u_{j_{0}+1}^{Y}$ on $B_{\r_{Y}}(Y)$ for some $j_{0} \in \{1, 2, \ldots, q-1\}$ and 
\begin{equation}\label{sheeting-13}
{\rm sup}_{B_{\r{Y}}(Y)} \, \r_{Y}^{-2}|u_{j}^{Y}|^{2} + |D \, u_{j}^{Y}|^{2} + \r_{Y}^{2}|D^{2} \, u_{j}^{Y}|^{2} \leq 
(C+ 2C_{1}) \widetilde{\g}^{-n-2}{\hat E}_{V}^{2}
\end{equation}
for each $j \in \{1, 2, \ldots, q\},$ where $C = C(n, q) \in (0, \infty).$
\end{itemize}

Let $\Omega = \{Y \in B_{\g} \,: \, ({\rm I}_{Y}) \;\; {\rm fails}\}.$ Since $u_{j}^{Y} < u_{j+1}^{Y}$ on $B_{\r_{Y}}(Y)$ for some $j$ 
whenever $({\rm II}_{Y})$ holds, it follows that $\Omega$ is an open set. Hence, since for every 
$Y \in \Omega,$ each of the functions $u_{j}^{Y}$ as in $({\rm II}_{Y})$ solves the minimal surface equation on $B_{\r_{Y}}(Y)$, by unique continuation of solutions to the minimal surface equation, we see that
\begin{equation}\label{sheeting-14}
V \res ({\mathbf R} \times \Omega) = \sum_{j=1}^{q} |{\rm graph} \, u_{j}|
\end{equation}
for functions $u_{j} \in C^{\infty} \, (\Omega; {\mathbf R}),$ solving the minimal surface equation on $\Omega$ and 
satisfying $u_{1} \leq u_{2} \leq \ldots \leq u_{q}$ on $\Omega;$  
$u_{j}< u_{j+1}$ for some $j \in \{1, 2, \ldots, q-1\}$ in each connected component of $\Omega$ (by the maximum principle) and 
\begin{equation}\label{sheeting-14-1}
{\rm sup}_{\Omega} |u_{j}|^{2} + |D \, u_{j}|^{2}  \leq (C + 2C_{1})\widetilde{\g}^{-n-2}{\hat E}_{V}^{2} \leq (C + 2C_{1})\widetilde{\g}^{-n-2}\e^{\prime}
\end{equation}
for each $j \in \{1, 2, \ldots, q\}.$ This implies that for each affine function $p \, : \, {\mathbf R}^{n} \to {\mathbf R}$ with $\sup_{B_{1}} |p|^{2} \leq C_{1}\widetilde{\g}^{-n-2}\e^{\prime}$ and each $j = 1, 2, \ldots, q$, the function $w_{j} = u_{j} - p \in C^{\infty} \, (\Omega)$ solves on $\Omega$ a uniformly elliptic equation of the type $a_{\ell k} D_{\ell}D_{\k} \, w_{j} + b_{\ell}D_{\ell} \, w_{j}= 0$ with smooth coefficients $a_{\ell k}$, $b_{\ell}$ satisfying $\sup_{\Omega} \, |a_{\ell k}| +|b_{\ell}| \leq \k$, $\k = \k(n, q, \g) \in (0, \infty).$ By using the standard second derivative estimates for solutions to such equations, we conclude that for each $Y \in \Omega,$ each $j = 1, 2, \ldots, q$ and each affine function $p_{j} \, : \, {\mathbf R}^{n} \to {\mathbf R}$ with $\sup_{B_{1}} \, |p_{j}|^{2} \leq C_{1}\widetilde{\g}^{-n-2}\e^{\prime}$, 
\begin{equation}\label{sheeting-15}
\s^{-n-2}\int_{B_{\s}(Y)} |u_{j} - p_{j}^{Y}|^{2} \leq C_{4}\left(\frac{\s}{\r}\right)^{2} \r^{-n-2}\int_{B_{\r}(Y)} |u_{j} - p_{j}|^{2}
\end{equation}
for $0 < \s \leq \r/2 < \frac{1}{2}{\rm dist} \, (Y, B_{\g} \setminus \Omega)$,  where $C_{4} = C_{4}(n, q, \g) \in (0, \infty)$ and 
$p_{j}^{Y}(X) = u_{j}(Y) + Du_{j}(Y) \cdot (X -Y).$ Since for each $Y \in B_{\g} \setminus \Omega$, 
 ${\rm spt} \, \|V\| \cap \pi^{-1}(Y)$ consists of a single point $(z_{Y}, Y)$ ($= P_{Y} \cap \pi^{-1}(Y)),$ for each $j=1,2, \ldots, q$, we may extend $u_{j}$ to all of $B_{\g}$ by setting $u_{j}(Y) = z_{Y}$ for $Y \in B_{\g} \setminus \Omega.$ Then by (\ref{sheeting-14}),
\begin{equation}\label{sheeting-15-1}
{\rm spt} \, \|V\| \cap ({\mathbf R} \times B_{\g})  = \cup_{j=1}^{q} {\rm graph} \, u_{j}.
\end{equation}
Now let $\widetilde\Sigma_{1}$, $\widetilde\Sigma_{2}$, $\widetilde\Sigma_{3},$ $\Sigma^{\prime}$ be the sets as in Theorem~\ref{flat-varifolds} taken with $\s = \g$. We claim that these sets are all empty if $\e^{\prime} = \e^{\prime}(n, q, \a, \g)$ is sufficiently small. Indeed, it is clear from (\ref{sheeting-14}), (\ref{sheeting-14-1}) and the definitions of $\widetilde\Sigma_{j}$, $\Sigma^{\prime}$ that $\widetilde\Sigma_{j} \cap ({\mathbf R} \times \Omega) = \emptyset$ for $j=1, 2, 3$ and that $\Sigma^{\prime} \cap \Omega = \emptyset.$ For each $Y \in B_{\g} \setminus \Omega$, we see by applying (\ref{tilt-ht-0}) with $\Gamma_{Y \, \#} \, V$ in place of $V$ where $\G_{Y}$ is a rigid motion of ${\mathbf R}^{n+1}$ that takes 
$(z_{Y}, Y) \in P_{Y}$ to the origin and $P_{Y}$ to $\{0\} \times {\mathbf R}^{n},$ and using the estimate (\ref{sheeting-12}), that $Y \not\in \pi \, \widetilde\Sigma_{1}$ provided $\e^{\prime} = \e^{\prime}(n, q, \a, \g)$ is sufficiently small. Since (\ref{sheeting-12}) implies that for each $Y \in B_{\g} \setminus \Omega$, the varifold $V$ has a unique tangent cone at $(z_{Y}, Y)$ with support equal to 
$P_{Y} - (z_{Y}, Y)$, it follows from the constancy theorem that $\Theta \, (\|V\|, (z_{Y}, Y))$ is a positive integer and furthermore, from the fact that varifold convergence implies Hausdorff convergence of supports, that ${\rm Tan} \, ({\rm spt} \, \|V\|, Y) = P_{Y} - (z_{Y}, Y);$ consequently, we see that $Y \not\in \pi \, \widetilde\Sigma_{2}$ and by (\ref{sheeting-11}), that $Y \not\in \pi \, \widetilde\Sigma_{3}.$ Finally, we argue that $\Theta(\|V\|, (z_{Y}, Y)) \geq q$ for each $Y \in B_{\g} \setminus \Omega$,  from which it follows that $\Sigma^{\prime} \cap (B_{\g} \setminus \Omega)  = \emptyset.$
If $\Theta \, (\|V\|, (z_{Y_{0}}, Y_{0})) < q$ for some $Y_{0} \in B_{\g} \setminus \Omega$, there is, by upper semi-continuity of density, some $\s_{0} > 0$ such that $\Theta \, (\|V\|, X) < q$ for each 
$X \in {\rm spt} \, \|V\| \cap ({\mathbf R} \times B_{\s_{0}}(Y_{0})),$ and hence, by Remark 3 of 
Section~\ref{outline},  the estimate (\ref{sheeting-12}) taken with $\s = \s_{0}$, $\rho = \widetilde\gamma/8$ and the estimate (\ref{sheeting-10-3}) taken with $\rho = \widetilde\gamma^{-1}\s_{0}$, 
we may, provided $\e^{\prime} = \e^{\prime}(n, q, \a, \g)$ is sufficiently small,  apply Theorem~\ref{SS} to conclude that 
$$V \res ({\mathbf R} \times B_{\s_{0}/2}) = \sum_{j=1}^{q} |{\rm graph} \, w^{j}|$$
for smooth functions $w_{1} \leq w_{2} \leq \ldots \leq w_{q}$ on $B_{\s_{0}/2}(Y_{0})$ solving the minimal surface equation.
Since $\pi^{-1}(Y_{0}) \cap {\rm spt} \, \|V\|$ consists of a single point, we must have by the maximum principle that $w_{1}= w_{2} = \ldots = w_{q}$ on $B_{\s_{0}/2}(Y_{0}),$ contrary to the assumption that $\Theta \, (\|V\|, (z_{Y_{0}}, Y_{0})) < q.$ This concludes the proof of the claim that the sets
$\widetilde\Sigma_{j}, \Sigma^{\prime}$  are all empty. By Theorem~\ref{flat-varifolds} and (\ref{sheeting-15-1}) then, 
for each $j=1, 2, \ldots, q$, the function $u_{j} \, : \, 
 B_{\g}  \to {\mathbf R}$ is Lipschitz with Lipschitz constant $\leq 1/2$, so that by (\ref{sheeting-12}),(\ref{sheeting-14-1}) and the area formula, it follows that 
\begin{equation}\label{sheeting-16}
\s^{-n-2}\int_{B_{\s}(Y)} |u_{j} - p^{Y}|^{2}  = 2C_{3}\left(\frac{\s}{\r}\right)^{2\a} \r^{-n-2}\int_{B_{\r}(Y)}|u^{j} - p^{Y}|^{2}
\end{equation}
for each $Y \in B_{\g} \setminus \Omega$ and each $\s$, $\r$ with $0 < \s < \r < \widetilde{\g}/4$, where $p^{Y} \, : \, {\mathbf R}^{n} \to {\mathbf R}$ is the affine function such that ${\rm graph} \, p^{Y} = P_{Y}.$

In view of (\ref{sheeting-15}) and (\ref{sheeting-16}), we conclude from Lemma~\ref{general} that
$u_{j} \in C^{1, \lambda} \, (B_{\g/2})$ with 
$$\sup_{B_{\g/2}} \, |u_{j}|^{2} + |Du_{j}|^{2} + \sup_{Y_{1}, Y_{2} \in B_{\g/2}, \, Y_{1} \neq Y_{2}} \, \frac{|Du_{j}(Y_{1}) - Du_{j}(Y_{2})|^{2}}{|Y_{1} - Y_{2}|^{2\lambda}}  \leq C_{5} {\hat E}_{V}^{2}$$
for each $j = 1,2, \ldots, q$, where $C_{5} = C_{5}(n, q, \g) \in (0, \infty)$ and $\lambda = \lambda(n, q, \g) \in (0, 1).$ 

To show that for each $j =1, 2, \ldots, q,$ the function $u_{j} \in C^{\infty}(B_{\g/2})$ and solves the minimal surface equation on $B_{\g/2}$, we argue as follows: We know that on the open set $\Omega \subseteq B_{\g}$, each $u_{j} \in C^{2}$  and solves the minimal surface equation (and hence is smooth), and on $B_{\g} \setminus \Omega$, the functions $u_{j}$ all agree, so if $B_{\g/2} \subseteq \Omega$ or $B_{\g/2} \cap \Omega = \emptyset$, there is nothing further to prove. Else, for any connected component $\Omega^{\prime}$ of $\Omega$ such that $B_{\g/2} \cap \Omega^{\prime} \neq \emptyset$, we must have that $B_{\g/2} \setminus \Omega^{\prime} \neq \emptyset$  whence $\partial \, \Omega^{\prime} \cap B_{\g/2} \neq \emptyset.$ Fix any such $\Omega^{\prime},$ and let $B \subset \Omega^{\prime}$ be an open ball such that ${\overline B} \cap \partial \, \Omega^{\prime} \cap B_{\g/2} \neq \emptyset.$ (To find such $B$, pick any point $p \in \Omega^{\prime}$ closer to $\partial \, \Omega^{\prime}$ than to $\partial \, B_{\g/2}$ and 
let $B = B_{R}(p)$ where $R = \sup \, \{r  \, : \, B_{r}(p) \subset \Omega^{\prime}\}$.) Let $x_{0} \in \partial \, B \cap \partial \, \Omega^{\prime} \cap B_{\g/2}.$ Pick any $j \in \{1, 2, \ldots, q-1\}$ and let 
$w_{j} = u_{j+1} - u_{j}$.  Then  $w_{j}$ solves in $B$ a uniformly elliptic equation with smooth coefficients. Since 
$w_{j} \in C^{1}(B_{\g/2}),$ $w_{j} \geq 0$ and  $w_{j}(x_{0}) = 0,$ it follows that $Dw_{j}(x_{0}) = 0,$ and hence by the Hopf boundary point lemma, $w_{j} \equiv 0$ in $B.$ This implies by unique continuation of solutions to the minimal surface equation that $w_{j} \equiv 0$ in $\Omega^{\prime}$  whence all of the $u_{j}$'s agree on $\Omega^{\prime}$ which is impossible by the definition of $\Omega$ (see the line preceding (\ref{sheeting-14-1})). Thus we must have either $B_{\g/2} \subseteq \Omega$ or $B_{\g/2} \cap \Omega  = \emptyset$, and the proof of the theorem is complete. 
\end{proof}

\section{The Minimum Distance Theorem}\label{MD}
\setcounter{equation}{0}

Let $q$ be an integer $\geq 2$ and let ${\mathbf C}_{0}$ be a stationary integral hypercone in ${\mathbf R}^{n+1}$ such that ${\rm spt} \, \|{\mathbf C}_{0}\|$ consists of three or more distinct half-hyperplanes of ${\mathbf R}^{n+1}$ meeting along a common $(n-1)$-dimensional 
subspace $L_{{\mathbf C}_{0}}$ of ${\mathbf R}^{n+1}.$ In this section we will use the multiplicity $q$ case of the Sheeting Theorem (i.e.\ Theorem~\ref{sheetingthm}) to establish, subject to the induction hypotheses $(H1)$, $(H2),$ the validity of Theorem~\ref{no-transverse}  whenever
\begin{equation}\label{MD-0}
\Theta \, (\|{\mathbf C}_{0}\|, 0)  \in \{q+1/2, q+1\};
\end{equation}
our argument will also establish Theorem~\ref{no-transverse} in case $\Theta \, (\|{\mathbf C}_{0}\|, 0) \in \{3/2, 2\}$; see the remark at the end of this section.

Suppose that ${\mathbf C}_{0}$ satisfies (\ref{MD-0}), and without loss of generality assume that $L_{{\mathbf C}_{0}} = \{0\} \times {\mathbf R}^{n-1}.$ Thus, 
${\rm spt} \, \|{\mathbf C}_{0}\| = {\rm spt} \, \|\Delta_{0}\| \times {\mathbf R}^{n-1}$, where $\Delta_{0}$ is a 1-dimensional stationary cone in ${\mathbf R}^{2}$, whence 
$\Delta_{0} = \sum_{j=1}^{m_{0}} q_{j}^{(0)}|R_{j}^{(0)}|$ and 
\begin{equation}\label{MD-1-0}
{\mathbf C}_{0} = \sum_{j=1}^{m_{0}}q_{j}^{(0)}|H_{j}^{(0)}|
\end{equation}
where $m_{0}$ is an integer $\geq 3$,  
$q_{j}^{(0)}$ is a positive integer for each $j=1, 2, \ldots, m_{0},$ $R_{j}^{(0)} = \{t{\mathbf w}_{j}^{(0)} \, : \, t > 0\}$
for some unit vector ${\mathbf w}_{j}^{(0)} \in {\mathbf S}^{1} \subset {\mathbf R}^{2}$ with ${\mathbf w}_{j}^{(0)} \neq {\mathbf w}_{k}^{(0)}$ for $j \neq k,$ and $H_{j}^{(0)} = R_{j}^{(0)} \times {\mathbf R}^{n-1}.$ Stationarity of ${\mathbf C}_{0}$ is equivalent to the requirement   
\begin{equation}\label{MD-0-0}
\sum_{j=1}^{m_{0}} q_{j}^{(0)}{\mathbf w}_{j}^{(0)} = 0.
\end{equation}
Since, by (\ref{MD-0}),  
\begin{equation}\label{MD-1}
\sum_{j=1}^{m_{0}} q_{j}^{(0)} \in  \{2q+1, 2q + 2\}, 
\end{equation}
we see readily from (\ref{MD-0-0}) that 
\begin{equation}\label{MD-2}
q_{j}^{(0)} \leq q\;\;\; \mbox{for each} \;\; j=1, 2, \ldots, m_{0}.
\end{equation}

The theorem we wish to prove is the following:
 
\begin{theorem}\label{no-transverse-q}
Let $q$ be an integer $\geq 2$, $\a \in (0, 1)$ and suppose that the induction hypotheses $(H1)$, $(H2)$ hold. Let ${\mathbf C}_{0}$  be the stationary cone as in (\ref{MD-1-0}), where $m_{0} \geq 3$ and $H_{j}^{(0)} \neq H_{k}^{(0)}$ for $j \neq k$, and suppose that ${\mathbf C}_{0}$ satisfies (\ref{MD-0}). For each $\g \in (0, 1/2)$, there exists a number 
$\e_{0} = \e_{0}(n,q, \a, \g, {\mathbf C}_{0}) \in (0, 1)$ such that if  $V \in {\mathcal S}_{\a},$ $\Theta \, (\|V\|, 0) \geq \Theta \, (\|{\mathbf C}_{0}\|, 0)$ and
$(\omega_{n}2^{n})^{-1}\|V\|(B_{2}^{n+1}(0)) < \Theta \, (\|{\mathbf C}_{0}\|, 0) + \g$, then 
$${\rm dist}_{\mathcal H} \, ({\rm spt} \, \|V\| \cap B_{1}^{n+1}(0), {\rm spt} \, \|{\mathbf C}_{0}\| \cap B_{1}^{n+1}(0)) \geq \e_{0}.$$
\end{theorem}

\noindent
{\bf Notation:} Let ${\mathbf C}_{0}$ be as in (\ref{MD-1-0}), with the associated unit vectors ${\mathbf w}_{j}^{(0)} \in {\mathbf R}^{2}$, $j =1, 2, \ldots, m_{0},$ as described above. We shall use the following notation in connection with ${\mathbf C}_{0}$:

$\s_{0}  = \max \, \{{\mathbf w}_{j}^{(0)} \cdot {\mathbf w}_{k}^{(0)} \, : \, j, k = 1, 2, \ldots, m_{0}, \; j \neq k\}.$ 

$N(H_{j}^{(0)})$ is the conical neighborhood of $H_{j}^{(0)}$ defined by 
$$N(H_{j}^{(0)}) = \left\{(x, y) \in {\mathbf R}^{2} \times {\mathbf R}^{n-1} \, : \, x \cdot {\mathbf w}_{j}^{(0)} > \sqrt{\frac{1 + \s_{0}}{2}}|x|\right\}.$$ 

Given ${\mathbf C}_{0}$ as above, ${\mathcal K}$ denotes the family of hypercones ${\mathbf C}$ of ${\mathbf R}^{n+1}$ of the form 
\begin{equation}\label{MD-3}
{\mathbf C} = \sum_{j=1}^{m_{0}}\sum_{\ell =1}^{q_{j}^{(0)}} |H_{j, \, \ell}|,
\end{equation}
where $H_{j, \, \ell}$ are half-hyperplanes of ${\mathbf R}^{n+1}$ meeting along $\{0\} \times {\mathbf R}^{n-1}$ with  
$H_{j, \, \ell} \in N(H_{j}^{(0)})$ for each $j \in \{1, 2, \ldots, m_{0}\}$,   $\ell \in \{1, 2, \ldots, q_{j}^{(0)}\},$ and $H_{j, \, 1}, H_{j, \, 2}, \ldots, H_{j, \, q_{j}^{(0)}}$ not necessarily distinct for each $j \in \{1, 2, \ldots, m_{0}\}.$ Note that unless otherwise specified, we do {\em not} assume a cone 
${\mathbf C} \in {\mathcal K}$ is stationary in ${\mathbf R}^{n+1}.$

For $p \in \{m_{0}, m_{0}+1, \ldots, 2q\}$, ${\mathcal K}(p)$ denotes the set of cones ${\mathbf C} \in {\mathcal K}$ as in (\ref{MD-3}) such that the number of \emph{distinct} elements
in the set $\{H_{j, \, \ell} \, : \, j=1, 2, \ldots, m_{0}, \;\; \ell = 1, 2, \ldots, q_{j}^{(0)}\}$ is $p$. Thus 
$${\mathcal K} = \cup_{p=m_{0}}^{2\Theta \, (\|{\mathbf C}_{0}\|,0)} {\mathcal K}(p).$$

Also, for $X \in {\mathbf R}^{n+1},$ let $r(X) = {\rm dist} \, (X, \{0\} \times {\mathbf R}^{n-1}).$  

For the rest of this section, we shall fix ${\mathbf C}_{0}$ as above, with fixed labelling of the elements of the set $\{H_{j}^{(0)} \, : \, j=1, \ldots, m_{0}\}$ of constituent half-hyperplanes of ${\rm spt} \, \|{\mathbf C}_{0}\|$ and with $q_{j}^{(0)}$, $1 \leq j \leq m_{0},$ denoting the multiplicity on $H_{j}^{(0)}$.
 
For $\a \in (0, 1)$, $\g \in (0, 1/2)$ and appropriate $\e \in (0, 1/2)$, consider the following:

\begin{hypotheses}\label{MD-hyp}
\begin{itemize}
\item[]
\item[(1)] $V \in {\mathcal S}_{\a}$, $0 \in {\rm spt} \, \|V\|$, $\Theta \, (\|V\|, 0) \geq \Theta \, (\|{\mathbf C}_{0}\|, 0)$, 
$(\omega_{n}2^{n})^{-1} \|V\|(B_{2}^{n+1}(0)) < \Theta \, (\|{\mathbf C}_{0}\|, 0) + \g.$ 
\item[(2)] ${\mathbf C}  = \sum_{j=1}^{m_{0}} \sum_{\ell=1}^{q_{j}^{(0)}} |H_{j, \, \ell}| \in {\mathcal K},$ where $H_{j, \, \ell}$ are half-hyperplanes of ${\mathbf R}^{n+1}$ meeting along $\{0\} \times {\mathbf R}^{n-1}$ with $H_{j, \, \ell} \in N(H_{j}^{(0)})$ for each $j \in \{1, 2, \ldots, m_{0}\}$ and $\ell \in\{1, 2, \ldots, q_{j}^{(0)}\}.$   
\item[(3)] ${\rm dist}_{\mathcal H} \, ({\rm spt} \, \|{\mathbf C}\| \cap B_{1}^{n+1}(0), {\rm spt} \, \|{\mathbf C}_{0}\| \cap B_{1}^{n+1}(0)) < \e.$
\item[(4)] $$\int_{B_{1}^{n+1}(0)} {\rm dist}^{2} \, (X, {\rm spt} \,\|{\mathbf C}\|) \, d\|V\|(X)< \e.$$
\item[(5)] For each $j=1,2, \ldots, m_{0},$ 
$$\|V\|((B_{1/2}^{n+1}(0) \setminus \{r(X) < 1/8\}) \cap N(H^{(0)}_{j})) \geq \left(q_{j}^{(0)} - \frac{1}{4}\right){\mathcal H}^{n}((B_{1/2}^{n+1}(0) \setminus \{r(X) < 1/8\}) \cap H^{(0)}_{j}).$$ 
\end{itemize}
\end{hypotheses}

Fix a number $s = s(n,q) \in (0, 1/16)$ such that 
\begin{equation}\label{MD-3-0}
{\mathcal H}^{n}\left(B_{\frac{1}{2}-s}^{n}(0) \setminus \{r(X) < 1/8+s\}\right) \geq \left(1 - \frac{1}{4q}\right){\mathcal H}^{n}\left(B_{1/2}^{n}(0) \setminus \{r(X) < 1/8\}\right)
\end{equation}
and note that by (\ref{MD-2}), 
\begin{eqnarray}\label{MD-3-1}
&&q_{j}^{(0)}{\mathcal H}^{n}((B_{\frac{1}{2}-s}^{n+1}(0) \setminus \{r(X) < 1/8+s\}) \cap H^{(0)}_{j})\nonumber\\
&&\hspace{2in}\geq \left(q_{j}^{(0)} - \frac{1}{4}\right){\mathcal H}^{n}((B_{1/2}^{n+1}(0) \setminus \{r(X) < 1/8\}) \cap H^{(0)}_{j})
\end{eqnarray}
for each $j=1, 2, \ldots, m_{0}.$ 

\noindent{\bf Remarks:} {\bf (1)} \emph{For each $\g \in (0, 1/2)$ and $\t \in (0, 1/8)$, there 
exists $\e = \e(n, q, \t, \g, {\mathbf C}_{0}) \in (0, 1)$ such that if the induction hypotheses (H1), (H2) and Hypotheses~\ref{MD-hyp} hold then
\begin{itemize}
\item[(a)] $\{Z \in {\rm spt} \, \|V\| \cap B_{15/16}^{n+1}(0) \, : \, \Theta \, (\|V\|, Z) \geq q + 1/2\} \subset \{X\in {\mathbf R}^{n+1} \, : \, r(X) < \t/2\}$ and
\item[(b)] for each $j \in \{1, 2, \ldots, m_{0}\}$ and $\ell \in \{1, 2, \ldots, q_{j}^{(0)}\}$, there exists a function
$$\widetilde{u}_{j, \, \ell}  \in C^{2}\,\left(\left(B_{15/16}^{n+1}(0) \cap H_{j}^{(0)} \setminus \{r(X) < \t\}\right); \left(H_{j}^{(0)}\right)^{\perp}\right)$$ 
with small $C^{2}$ norm such that $\widetilde{u}_{j, \, \ell}$ solves the minimal surface equation on its domain and $$V \res \left(B_{15/16}^{n+1}(0) \setminus \{r(X) < \t\}\right) = \sum_{j=1}^{m_{0}} \sum_{\ell=1}^{q_{j}^{(0)}} |{\rm graph} \, \widetilde{u}_{j, \, \ell}|.$$
\end{itemize}
} 

To see this, argue by contradiction, using the Constancy Theorem, upper semi-continuity of the density function $\Theta \, (\cdot, \cdot),$ (\ref{MD-2}), induction hypothesis ($H1$) and Theorem~\ref{sheetingthm}.
 
\noindent
{\bf (2)}  \emph{For each $\g \in (0, 1/2)$ and $\t \in (0, 1/8)$, there 
exists $\e = \e(n, q, \t, \g, {\mathbf C}_{0}) \in (0, 1)$ such that if Hypotheses~\ref{MD-hyp}(1)-(4)  hold and if (in place of Hypothesis~\ref{MD-hyp}(5)) 
$$\;\; \int_{B_{1/2}^{n+1}(0) \setminus \{r(X) < 1/8\}} {\rm dist}^{2} \, (X, {\rm spt} \,\|V\|) \, d\|{\mathbf C}\|(X)< \e,$$ then $\{Z \in {\rm spt} \, \|V\| \cap B_{15/16}^{n+1}(0) \, : \, \Theta \, (\|V\|, Z) \geq q + 1/2\} \subset \{X\in {\mathbf R}^{n+1} \, : \, r(X) < \t\}.$} 

Again, this is easily seen by arguing by contradiction using the Constancy Theorem, upper semi-continuity of density and (\ref{MD-2}).

\noindent
{\bf (3)} \emph{Let $q$ be an integer $\geq 2$. If the induction hypotheses $(H1)$, $(H2)$ hold, $V \in {\mathcal S}_{\a},$ $\Omega \subseteq B_{2}^{n+1}(0)$ is open and $\Theta \, (\|V\|, Z) < q +1/2$ for each $Z \in {\rm spt} \, \|V\| \cap \Omega$, then 
${\mathcal H}^{n-7 + \g} \, ({\rm sing} \, V \, \res \, \Omega) = 0$ for each $\g > 0$ if $n \geq 7$, ${\rm sing} \, V \, \res \, \Omega$ is discrete if $n=7$ and 
${\rm sing} \, V \, \res \, \Omega = \emptyset$ if $2 \leq n \leq 6.$}

This can be seen by reasoning exactly as in Remarks (2) and (3) of Section~\ref{outline}, with the additional help of Theorem~\ref{sheetingthm}.

\noindent
{\bf (4)} \emph{Let $\g \in (0,1/2)$, $\rho \in (0, 1/2]$ and $\e^{\prime} \in (0, 1/2)$. There exists a number $\e = \e(\rho, \e^{\prime}, \a, \g,{\mathbf C}_{0}) \in (0, 1/2)$ such that  if Hypotheses~\ref{MD-hyp} are satisfied, then for each $Z \in {\rm spt} \, \|V\| \cap B_{1/8}^{n+1}(0)$ with 
$\Theta \, (\|V\|, Z) \geq q + 1/2,$ Hypotheses~\ref{MD-hyp} are also satisfied with 
$\eta_{Z, \rho \, \#} \, V$ in place of $V$ and $\e^{\prime}$ in place of $\e$.} 

Indeed, given any $\r \in (0, 1/2]$, if $V$, ${\mathbf C}$ are as in Hypotheses~\ref{MD-hyp} with sufficiently small $\e  = \e(\r, \a, {\mathbf C}_{0}) \in (0, 1/2),$ then it follows from Remark (1) applied with suitably small $\t  = \t(\r, \g) \in (0, 1/16)$  and the fact that $\|V\|(B_{1}^{n+1}(0) \cap \{X \, : \, r(X) < \t\}) \leq C\t$ where $C = C(n, q) \in (0, \infty)$ that for any $Z \in {\rm spt} \, \|V\| \cap B_{1/8}^{n+1}(0)$ with $\Theta(\|V\|, Z) \geq q + 1/2$ , Hypothesis~\ref{MD-hyp}(1) is satisfied with $\eta_{Z, \r \, \#} \, V$ in place of $V$; also, since by the triangle inequality 
\begin{eqnarray*}
&&\int_{B_{1}^{n+1}(0)} {\rm dist}^{2} \, (X, {\rm spt} \, \|{\mathbf C}\|) d\|\eta_{Z, \r \, \#} \, V\|(X) \leq
2\r^{-n-2}\int_{B_{\r}^{n+1}(Z)}{\rm dist}^{2} \, (X, {\rm spt} \, \|{\mathbf C}\|) \, d\|V\|(X)\nonumber\\
&&\hspace{3.5in} +  C\r^{-2}{\rm dist}^{2} \, (Z, \{0\} \times {\mathbf R}^{n-1}),
\end{eqnarray*}
where $C = C(n, q, \g) \in (0, \infty)$, it follows again by Remark (1) (taken with $\t  = \r\sqrt{(2C)^{-1}\e^{\prime}}$) that if $\e = \e(\r, \e^{\prime}, \a, \g, {\mathbf C}_{0})$ is sufficiently small, then Hypothesis~\ref{MD-hyp}(4) is satisfied with 
$\eta_{Z, \r \, \#} \, V$ in place of $V$ and $\e^{\prime}$ in place of $\e$;  finally, applying Remark (1) once again with $\t = \r s,$ where $s = s(n, q)$ is as in (\ref{MD-3-0}), we deduce using (\ref{MD-3-1}),  the area formula and the inclusion 
${\rm spt} \, \|V\| \cap \left(B_{\r - \t}^{n+1}(0, \eta) \setminus B^{2}_{\frac{\r}{8} + \t}(0) \times {\mathbf R}^{n-1}\right) \subset  {\rm spt} \, \|V\| \cap \left(B_{\r}^{n+1}(Z) \setminus B^{2}_{\t}(0) \times {\mathbf R}^{n-1}\right)$
where $(0, \eta)$ is the orthogonal projection of $Z$ onto $\{0\} \times {\mathbf R}^{n-1}$ that if $\e = \e(\r, \a, \g, {\mathbf C}_{0})$ is sufficiently small, then Hypothesis~\ref{MD-hyp}(5) is satisfied with $\eta_{Z, \r \, \#} \, V$ in place of $V.$

\medskip

With the notation as above, for $V \in {\mathcal S}_{\a}$, ${\mathbf C} \in {\mathcal K}$ as in Hypotheses~\ref{MD-hyp} and appropriate $\b \in (0, 1/2)$, we will also need to consider the following:

\noindent
{\bf Hypothesis $(\dag)$}: \emph{Either {\rm(i)} or {\rm (ii)} below holds:
\begin{eqnarray*}
&&{\rm (i)} \;\;\;\;{\mathbf C} \in {\mathcal K}(m_{0}).\nonumber\\ 
&&{\rm (ii)} \;\;\;\;2\Theta \, (\|{\mathbf C}_{0}\|, 0) \geq m_{0}+1, \;\; {\mathbf C} \in {\mathcal K}(p) \;\; \mbox{for some $p \in \{m_{0}+1, m_{0}+2, \ldots, 2\Theta \, (\|{\mathbf C}_{0}\|, 0)\}$ and}\nonumber\\ 
&&\int_{B_{1}^{n+1}(0)} {\rm dist}^{2} \, (X, {\rm spt} \,\|{\mathbf C}\|) \, d\|V\|(X) +  \int_{B_{1}^{n+1}(0) \setminus \{r(X) < 1/16\}} {\rm dist}^{2} \, (X, {\rm spt} \,\|V\|) \, d\|{\mathbf C}\|(X)\nonumber\\
&&\hspace{1in} \leq \; \b\; {\rm inf}_{\widetilde{\mathbf C} \in \cup_{j=m_{0}}^{p-1} {\mathcal K}(j)} 
\left(\int_{B_{1}^{n+1}(0)} {\rm dist}^{2} \, (X, {\rm spt} \,\|\widetilde{\mathbf C}\|) \, d\|V\|(X)\right.\nonumber\\
&&\hspace{2in}+ \left.\int_{B_{1}^{n+1}(0) \setminus \{r(X) < 1/16\}} {\rm dist}^{2} \, (X, {\rm spt} \,\|V\|) \, d\|\widetilde{\mathbf C}\|(X)   \right). 
\end{eqnarray*}} 

\noindent
{\bf Remark:} If Hypotheses~\ref{MD-hyp}, 
Hypothesis ($\dag$) for some $\b \in (0, 1/2)$ are satisfied, and if ${\mathbf C}^{\prime} \in {\mathcal K}$ is such that ${\rm spt} \, \|{\mathbf C}^{\prime}\| = {\rm spt} \, \|{\mathbf C}\|,$ then, Hypotheses~\ref{MD-hyp}, Hypothesis ($\dag$) taken with ${\mathbf C}^{\prime}$ in place of ${\mathbf C}$ and $2q\b$ in place of $\b$
will be satisfied.

\medskip

\noindent
{\bf Case $\Theta \, (\|{\mathbf C}_{0}\|,0) = q + 1/2$:} From now on until we state otherwise, we shall assume that $\Theta \, (\|{\mathbf C}_{0}\|, 0) = q + 1/2.$

\medskip

The basic $L^{2}$-estimates of [\cite{S}, Theorem 3.1]  hold under our assumptions, namely, the induction hypotheses $(H1)$, $(H2)$, Hypotheses~\ref{MD-hyp} and Hypothesis~($\dag$), and are given in Theorem~\ref{MD-L2-est-1} and Corollary~\ref{MD-L2-est-2} below:

\begin{theorem}\label{MD-L2-est-1}
Let $q$ be an integer $\geq 2$, $\a \in (0,1)$, $\g \in (0, 1/2),$ $\m \in (0, 1)$ and $\t \in (0, 1/8).$ Suppose that the induction hypotheses $(H1)$, $(H2)$ hold. Let ${\mathbf C}_{0}$ be a stationary cone as above, with $\Theta \, (\|{\mathbf C}_{0}\|, 0)  = q+1/2$. There exist numbers $\e_{0} = \e_{0}(n, q, \a, \g,\t,  {\mathbf C}_{0}) \in (0, 1/2)$, $\b_{0}  = \b_{0}(n, q, \a, \g, \t,{\mathbf C}_{0}) \in (0, 1/2)$ such that if $V \in {\mathcal S}_{\a},$ ${\mathbf C} \in {\mathcal K}$ satisfy Hypotheses~\ref{MD-hyp} with $\e_{0}$ in place of $\e$ and Hypothesis {\rm (}$\dag${\rm )} with $\b_{0}$ in place of $\b$, then, after taking appropriate ${\mathbf C}^{\prime}  \in {\mathcal K}$ with 
${\rm spt} \, \|{\mathbf C}^{\prime}\| = {\rm spt} \, \|{\mathbf C}\|$ in place of ${\mathbf C}$, relabelling ${\mathbf C}^{\prime}$ as ${\mathbf C}$ (see the preceding Remark) and writing 
${\mathbf C}  = \sum_{j=1}^{m_{0}} \sum_{\ell=1}^{q_{j}^{(0)}} |H_{j, \, \ell}|$ where $H_{j, \, \ell}$ are half-hyperplanes of ${\mathbf R}^{n+1}$ meeting along $\{0\} \times {\mathbf R}^{n-1}$ with $H_{j, \, \ell} \in N(H_{j}^{(0)})$ for each $j \in \{1, 2, \ldots, m_{0}\}$ and $\ell \in\{1, 2, \ldots, q_{j}^{(0)}\},$ the following hold: 
\begin{eqnarray*}
&&({\rm a})\;\;\;\; V \res (B^{n+1}_{7/8}(0) \setminus \{r(X) < \t\}) = \sum_{j=1}^{m_{0}} \sum_{\ell=1}^{q_{j}^{(0)}}|{\rm graph} \, u_{j, \, \ell} \cap B^{n+1}_{7/8}(0)|\;\; \mbox{where}\nonumber\\
&&u_{j,\,  \ell} \in C^{2} \, (B_{7/8}^{n+1}(0) \,\cap \, H_{j, \, \ell} \setminus \{r(X) < \t\}; H_{j, \, \ell}^{\perp}) \;\;\mbox{for $1 \leq j \leq m_{0}$, $1 \leq \ell  \leq q_{j}^{(0)},$}\nonumber\\ 
&&\mbox{$u_{j, \, \ell}$ solves the minimal surface equation on $B_{7/8}^{n+1}(0) \cap H_{j, \, \ell} \setminus \{r(X) < \t\}$,}\nonumber\\ 
&&{\rm dist} \, (X + u_{j, \, \ell}(X), {\rm spt} \,\|{\mathbf C}\|) = |u_{j, \, \ell}(X)|\;\; \mbox {for}\;\; X \in B_{7/8}^{n+1}(0) 
\cap H_{j, \, \ell} \setminus \{r(X)  < \t\} \;\; \mbox{and for each}\nonumber\\
&&\mbox{$j \in \{1, \ldots, m_{0}\}$ and $\ell_{1}, \ell_{2} \in \{1,  \ldots, q_{j}^{(0)}\}$, either 
${\rm graph} \, u_{j, \, \ell_{1}} \cap B^{n+1}_{7/8}(0) \equiv {\rm graph} \, u_{j, \ell_{2}} \cap B^{n+1}_{7/8}(0)$}\nonumber\\
&&\mbox{or  ${\rm graph} \, u_{j, \, \ell_{1}} \cap {\rm graph} \, u_{j, \ell_{2}} \cap B^{n+1}_{7/8}(0) = \emptyset$};\nonumber\\
&&({\rm b})\;\;\;\;\int_{B_{3/4}^{n+1}(0)} \frac{|X^{\perp}|^{2}}{|X|^{n+2}} \, d\|V\|(X) \leq C \int_{B_{1}^{n+1}(0)} {\rm dist}^{2} \, (X, {\rm spt} \,\|{\mathbf C}\|) \, d\|V\|(X);\nonumber\\
&&({\rm c})\;\;\;\; \sum_{j=3}^{n+1}\int_{B_{3/4}^{n+1}(0)} |e_{j}^{\perp}|^{2} \, d\|V\|(X) \leq C \int_{B_{1}^{n+1}(0)} {\rm dist}^{2} \, (X, {\rm spt} \,\|{\mathbf C}\|) \, d\|V\|(X);\nonumber\\
&&({\rm d})\;\;\;\;\int_{B_{3/4}^{n+1}(0)} \frac{{\rm dist}^{2} \, (X, {\rm spt} \, \|{\mathbf C}\|)}{|X|^{n+2 - \m}} \, d\|V\|(X) 
\leq C_{1}\int_{B_{1}^{n+1}(0)} {\rm dist}^{2} \, (X, {\rm spt} \,\|{\mathbf C}\|) \, d\|V\|(X).\nonumber\\
\end{eqnarray*}
Here $C = C(n, \a, \g, {\mathbf C}_{0}) \in (0, \infty)$ and $C_{1}  = C_{1}(n, \a, \g, \m, {\mathbf C}_{0}) \in (0, \infty).$ In particular, $C$, $C_{1}$ do not depend on $\t.$
\end{theorem}

\begin{proof} Note first that by Remark (1) following Hypotheses~\ref{MD-hyp}, provided the hypotheses of the theorem are satisfied with $\e_{0} = \e_{0}(n, q, \a, \t, {\mathbf C}_{0})$  sufficiently small, we have that $\Theta \, (\|V\|, Z) < q + 1/2$ for each $Z \in {\rm spt} \, \|V\| \cap B_{15/16}^{n+1}(0) \setminus \{r(X) < \t/2\},$ and 
$$V \res \left(B_{7/8}^{n+1}(0) \setminus \{r(X) < \t\}\right) = \sum_{j=1}^{m_{0}} \sum_{\ell=1}^{q_{j}^{(0)}} |{\rm graph} \, \widetilde{u}_{j, \, \ell}|,$$ 
where for each $j \in \{1, 2, \ldots, m_{0}\}$ and $\ell \in \{1, 2, \ldots, q_{j}^{(0)}\}$, 
$$\widetilde{u}_{j, \, \ell}  \in C^{2}\,\left(\left(B_{7/8}^{n+1}(0) \cap H_{j}^{(0)} \setminus \{r(X) < \t\}\right); \left(H_{j}^{(0)}\right)^{\perp}\right)$$ 
and $\widetilde{u}_{j, \, \ell}$ are solutions to the minimal surface equation over $H_{j}^{(0)} \cap \left(B_{7/8}^{n+1}(0) \setminus \{r(X) < \t\}\right)$ with small $C^{2}$ norm. So if ${\mathbf C} \in {\mathcal K}(m_{0})$, then the desired conclusions in part (a) readily follow because 
then ${\mathbf C} = \sum_{j=1}^{m_{0}} q_{j}^{(0)}|H_{j}^{\prime}|$ for 
distinct half-hyperplanes $H_{j}^{\prime}$ meeting along $\{0\} \times {\mathbf R}^{n-1}$, which, by Hypotheses~\ref{MD-hyp}(3), satisfy ${\rm dist}_{\mathcal H} \, (H_{j}^{\prime} \cap B_{1}^{n+1}(0), H_{j}^{(0)} \cap B_{1}^{n+1}(0)) < \e_{0}$ for each $j \in \{1, 2, \ldots, m_{0}\}.$ Else we must have that $2\Theta \,(\|{\mathbf C}_{0}\|,0) \geq m_{0} +1$ and that ${\mathbf C} \in {\mathcal K}(p)$ for some $p \in \{m_{0}+1, m_{0}+2, \ldots, 2\Theta \, (\|{\mathbf C}_{0}\|, 0)\}.$ For each $j \in \{1, 2, \ldots, m_{0}\}$, let $q_{j}^{\prime} \in \{1, 2, \ldots, q_{j}^{(0)}\}$ be the number of distinct elements in the set $\{H_{j, \, 1}, H_{j, \, 2} , \ldots, H_{j, \, q_{j}^{(0)}}\}$ and label them $H^{\prime}_{j, \, \ell^{\prime}}$, $\ell^{\prime} = 1, 2, \ldots, q_{j}^{\prime}.$  Let ${\mathbf w}_{j, \, \ell^{\prime}}^{\prime} \in {\mathbf R}^{2}$ be the unit vector such that 
$H_{j, \, \ell^{\prime}}^{\prime} = \{(t{\mathbf w}_{j, \, \ell^{\prime}}^{\prime}, y) \, : \, t >0, \;\; y \in {\mathbf R}^{n-1}\}.$ 
Provided that $\b_{0} = \b_{0}(\a, \g, {\mathbf C}_{0}) \in (0, 1/2)$ is sufficiently small,  it follows from the definition of ${\mathcal K}$, Hypotheses ~\ref{MD-hyp}(3) and Hypothesis ($\dag$) (ii) that for each $j \in \{1, 2, \ldots, m_{0}\}$ and $\ell^{\prime}_{1}, \ell^{\prime}_{2} \in \{1, 2, \ldots, q_{j}^{\prime}\}$, 
$$|{\mathbf w}_{j, \, \ell^{\prime}_{1}}^{\prime} - {\mathbf w}_{j, \, \ell^{\prime}_{2} }^{\prime}| \geq c \, {\mathcal Q}_{V}^{\prime}$$
for some constant $c = c(\a, \g, {\mathbf C}_{0}) \in (0, \infty),$  where
\begin{eqnarray*}
{\mathcal Q}^{\prime}_{V} = {\rm inf}_{\widetilde{\mathbf C} \in \cup_{j=m_{0}}^{p-1} {\mathcal K}(j)} 
\left(\int_{B_{1}^{n+1}(0)} {\rm dist}^{2} \, (X, {\rm spt} \,\|\widetilde{\mathbf C}\|) \, d\|V\|(X)\right.&&\\
&&\hspace{-2in}+ \left.\int_{B_{1}^{n+1}(0) \setminus \{r(X) < 1/16\}} {\rm dist}^{2} \, (X, {\rm spt} \,\|V\|) \, d\|\widetilde{\mathbf C}\|(X)   \right).
\end{eqnarray*}
By exactly the same inductive proof of Theorem~\ref{L2-est-1}(a),  conclusion (a) now follows from this provided $\e_{0} = \e_{0}(n,\a, \g,\t, {\mathbf C}_{0}), \b_{0} = \b_{0}(n, \a, \g,\t, {\mathbf C}_{0}) \in (0, 1/2)$ are  sufficiently small.

The rest of the theorem is proved by arguing exactly as in \cite{S}, Theorem~3.4; the key point which enables us to use the argument of \cite{S}, Theorem~3.4 is having at our disposal the appropriate regularity theorem, namely, Theorem~\ref{sheetingthm}. Specifically, letting
$$T_{\r, \k}(\z) = \{(x, y) \in {\mathbf R}^{2} \times {\mathbf R}^{n-1} \, : \, (|x| - \r)^{2} + |y - \z|^{2} < \k^{2}\r^{2}/64\}$$
for $\k \in (0, 1]$, $\r \in (0, 1/2)$ and $\z \in {\mathbf R}^{n-1},$ we have the following for any given 
 $\b \in (0, 1)$:

\noindent
{\bf Claim:}  There exists a constant $\d = \d(n, q, \a, \g, \b, {\mathbf C}_{0}) \in (0, 1/2)$ such that if $V$, ${\mathbf C}$ are as in the theorem, 
$(\xi, \z) \in {\rm spt} \, \|{\mathbf C}\| \cap B_{13/16}^{n+1}(0) \cap \{r(X) < 1/16\}$ where $\z \in {\mathbf R}^{n-1},$  
\begin{equation}\label{MD-L2-a}
{\rm spt} \, \|V\| \cap T_{|\xi|, 1/16}(\z) \neq \emptyset \;\;\; {\rm and} 
\end{equation}
\begin{equation}\label{MD-L2-b}
|\xi|^{-n-2}\int_{T_{|\xi|, 1}(\z)} {\rm dist}^{2} \, (X, {\rm spt} \, \|{\mathbf C}\|) \, d\|V\|(X) < \d,
\end{equation}
then there exist distinct integers $j_{1}, j_{2}, 
\ldots, j_{p} \in \{1, 2, \ldots, m_{0}\}$ and, for each $k \in \{1, 2, \ldots, p\}$, functions 
$u^{(|\xi|, \z)}_{j_{k}, k_{\ell}} \in C^{2}(T_{|\xi|, 3/4}(\z) \cap H_{j_{k}, k_{\ell}}; H_{j_{k}, k_{\ell}}^{\perp})$
with $\ell = 1, 2, \ldots, n_{k}$  for some $n_{k} \leq q$  such that 
\begin{equation}\label{MD-L2-c}
V \res T_{|\xi|, 1/2}(\z) = \sum_{k=1}^{p} \sum_{\ell=1}^{n_{k}} |{\rm graph} \, u^{(|\xi|, \z)}_{j_{k}, k_{\ell}} \cap T_{|\xi|, 1/2}(\z)|
\end{equation}
and for each $k \in \{1, \ldots, p\}$, $\ell \in \{1, \ldots, n_{k}\}$, 
$$|\xi|^{-1} \sup_{T_{|\xi|, 3/4}(\z) \cap H_{j_{k}, k_{\ell}}}\, |u_{j_{k}, k_{\ell}}^{(|\xi|, \z)}| + \sup_{T_{|\xi|, 3/4}(\z) \cap H_{j_{k}, k_{\ell}}} \, |\nabla \, u_{j_{k}, k_{\ell}}^{(|\xi|, \z)}| \leq \b/2$$

To verify this claim, observe first that by using the monotonicity of mass ratio and a covering argument, we see that under the hypotheses of the theorem,  
$\|V\|(B_{1}^{n+1}(0) \cap \{r(X) < \t\}) \leq C\t$ for each $\t \in (0, 1/4)$ where $C = C(n, q) \in (0, \infty).$ Using this with sufficiently small $\t = \t(n, q) \in (0, 1/2)$ and conclusion (a), we deduce that 
$\omega_{n}^{-1}\r^{-n}\|V\|(B^{n+1}_{\r}(Z)) \leq \omega_{n}^{-1}(16)^{n}\|V\|(B^{n+1}_{1/16}(Z)) < q + 3/4$ for each $Z \in B^{n+1}_{13/16}(0)$ and $\r \in (0, 1/16)$ provided $\e_{0} = \e_{0}(n, q, \a, \g,  {\mathbf C}_{0}) \in (0, 1)$ is sufficiently small. Since (\ref{MD-L2-b}) for sufficiently small 
$\d = \d(n, q, \a, \g, {\mathbf C}_{0}) \in (0, 1/2)$ implies that 
$V \res T_{|\xi|, 7/8}(\z) = \sum_{j=1}^{m_{0}}V_{j}$ where ${\rm spt} \, \|V_{j}\| \subset N(H_{j}^{(0)}) \cap T_{|\xi|, 7/8}(\z)$  (allowing for the possibility  that $V_{j} = 0$ for some values of $j$) we see by applying Theorem~\ref{sheetingthm} and Remark (3) at the end of Section~\ref{step2} that (\ref{MD-L2-c}) follows from (\ref{MD-L2-a}) and (\ref{MD-L2-b}) as claimed. 

Now let, as in Lemma~2.6 of \cite{S}, $U$ be the union of all $T_{|\xi|, 1/2}(\z) \cap {\rm spt} \, \|{\mathbf C}\|$ over all $(\xi, \z) \in {\rm spt} \, \|{\mathbf C}\| \cap B_{7/8}^{n+1}(0)$ such that for each $j \in \{1, \ldots, m_{0}\}$ and each $\ell \in \{1, \ldots, q_{j}^{(0)}\}$, there exists $u_{j, \ell}^{(|\xi|, \z)} \in C^{2}(T_{|\xi|, 3/4}(\z) \cap H_{j, \ell}; H_{j, \ell}^{\perp})$ 
with 
$$|\xi|^{-1} \sup_{T_{|\xi|, 3/4}(\z) \cap H_{j, \ell}}\, |u_{j, \ell}^{(|\xi|, \z)}| + \sup_{T_{|\xi|, 3/4}(\z) \cap H_{j, \ell}} \, |\nabla \, u_{j, \ell}^{(|\xi|, \z)}| \leq \b/2,$$
$${\rm dist} \, (X + u_{j, \ell}^{(|\xi|, \z)}(X), {\rm spt} \,\|{\mathbf C}\|) = |u_{j, \ell}^{(|\xi|, \z)}(X)| \;\;\;\; \mbox{for each $X \in T_{|\xi|, 1/2}(\z) \cap H_{j, \ell}$ and}$$  
$$V \res T_{|\xi|, 1/2}(\z) = \sum_{j=1}^{m_{0}} \sum_{\ell=1}^{q_{j}^{(0)}} |{\rm graph} \, u^{(|\xi|, \z)}_{j, \ell} \cap T_{|\xi|, 1/2}(\z)|.$$
For each $j \in \{1, \ldots, m_{0}\},$ $\ell \in \{1, \ldots, q_{j}^{(0)}\}$, define $u_{j, \ell} \in C^{2}(U \cap H_{j, \ell}; H_{j, \ell}^{\perp})$ by 
$\left.u_{j, \ell} \right|_{T_{|\xi|, 1/2}(\z) \cap H_{j, \ell}} = u_{j, \ell}^{(|\xi|, \z)}.$ With the help of the above claim and unique continuation of solutions to the minimal surface equation, we may now verify the validity of Lemma~2.6 of \cite{S} (by following the same proof), with the conclusion that for each $j \in \{1, \ldots, m_{0}\}$ and $\ell \in \{1, \ldots, q_{j}^{(0)}\}$, 
$$H_{j, \ell} \cap B_{7/8}^{n+1}(0) \setminus \{r(X) < \t\} \subset U;$$ 
there exists $u_{j, \ell} \in C^{2}(U \cap H_{j, \ell}; H_{j, \ell}^{\perp})$ 
such that 
$$\sup_{U \cap H_{j, \ell}} r^{-1}|u_{j, \ell}| + |\nabla \, u_{j, \ell}| \leq \b;$$
and 
\begin{eqnarray*}
&&\int_{B_{7/8}^{n+1}(0) \setminus G} r^{2}(X) \, d\|V\|(X) + \sum_{j=1}^{m_{0}}\sum_{\ell=1}^{q_{j}^{(0)}}\int_{U \cap H_{j, \ell}}r^{2}(X) |\nabla \, u_{j, \ell}(X)|^{2} \, d{\mathcal H}^{n}(X)\nonumber\\
&&\hspace{3in} \leq C \int_{B_{1}^{n+1}(0)} {\rm dist}^{2} \, (X, {\rm spt} \, \|{\mathbf C}\|) \, d\|V\|(X)
\end{eqnarray*}
where $G = \cup_{j=1}^{m_{0}} \cup_{\ell=1}^{q_{j}^{(0)}} {\rm graph} \, u_{j, \ell}$ and 
$C = C(n, \a, \g, {\mathbf C}_{0}) \in (0, \infty).$
Consequently, the argument of Lemma~3.4 of \cite{S} carries over 
to give conclusions (b)-(d) of the present theorem.
\end{proof}

%\medskip

\begin{corollary}\label{MD-L2-est-2}
Let $q$ be an integer $\geq 2$, $\a \in (0,1)$, $\g \in (0, 1/2)$ and $\m \in (0, 1).$  Suppose that the induction hypotheses $(H1)$, $(H2)$ hold and let ${\mathbf C}_{0}$ be the stationary cone as in \ref{MD-1-0}, with $\Theta \, (\|{\mathbf C}_{0}\|, 0) = q + 1/2$. For each $\r \in (0, 1/4]$, there exist numbers $\e = \e(n, \a, \g,\t, \r, {\mathbf C}_{0}) \in (0, 1/2)$, $\b  = \b(n, \a, \g, \t, \r,{\mathbf C}_{0}) \in (0, 1/2)$ such that if $V \in {\mathcal S}_{\a},$ ${\mathbf C} \in {\mathcal K}$ satisfy Hypotheses~\ref{MD-hyp} and Hypothesis {\rm (}$\dag${\rm )}, then for each $Z = (\z^{1}, \z^{2}, \eta) \in {\rm spt} \, \|V\| \cap (B_{3/8}^{n+1}(0))$ with $\Theta \, (\|V\|, Z) \geq \Theta \, (\|{\mathbf C}_{0}\|, 0),$ we have the following: 
\begin{eqnarray*}
&&\hspace{-.2in}({\rm a})\;\;\;\; |\z^{1}|^{2} + |\z^{2}|^{2} \leq C \int_{B_{1}^{n+1}(0)} {\rm dist}^{2} \, (X, {\rm spt} \,\|{\mathbf C}\|) \, d\|V\|(X).\nonumber\\
&&\hspace{-.2in}({\rm b})\int_{B_{\r/2}^{n+1}(Z)} \frac{{\rm dist}^{2} \, (X, {\rm spt} \, \|T_{Z \, \#} \, {\mathbf C}\|)}{|X - Z|^{n+2 - \m}} \, d\|V\|(X)\nonumber\\ 
&&\hspace{2in}\leq C_{1}\r^{-n-2+\m}\int_{B_{\r}^{n+1}(Z)} {\rm dist}^{2} \, (X, {\rm spt} \,\|T_{Z \, \#} \, {\mathbf C}\|) \, d\|V\|(X)\nonumber\\
&&\hspace{-.2in}\mbox{where $T_{Z} \, : \, {\mathbf R}^{n+1} \to {\mathbf R}^{n+1}$ is the translation $X \mapsto X + Z.$}
\end{eqnarray*}
Here $C = C(n, \a, \g, {\mathbf C}_{0}) \in (0, \infty)$ and $C_{1}  = C_{1}(n, \a, \g, \m, {\mathbf C}_{0}) \in (0, \infty).$ (In particular, $C$, $C_{1}$ do not depend on $\r$.)
\end{corollary}

\begin{proof} The proof requires application of Theorem~\ref{MD-L2-est-1} with $\eta_{Z, \r \, \#} \, V$ in place of $V$, 
where $Z \in {\rm spt} \, \|V\| \cap B_{3/8}^{n+1}(0)$ is any point such that $\Theta \, (\|V\|, Z) \geq \Theta \, (\|{\mathbf C}_{0}\|, 0).$ It follows from Remark (4) above that whenever Hypotheses~\ref{MD-hyp} are satisfied with $\e = \e(n, q, \a, \t, \r, {\mathbf C}_{0})$ sufficiently small, they  are also satisfied with $\eta_{Z, \r \, \#} \, V$ in place of $V$ and $\e_{0}$ (as in Theorem~\ref{MD-L2-est-1}) in place of $\e$; to verify that Hypothesis~($\dag$) is satisfied with $\eta_{Z, \r \, \#} \, V$ in place of $V$ and $\b_{0}$ (as in Theorem~\ref{MD-L2-est-1}) in place of $\b$, and complete at the same time the proof of the corollary inductively, we may follow the steps of the 
proof of Corollary~\ref{L2-est-2} (i.e.\ Lemmas~\ref{L2-est-2-L1}, \ref{L2-est-2-L2}, \ref{L2-est-2-L3} and Propositions~\ref{L2-est-2-P1}, \ref{L2-est-2-P2}) in conjunction with the argument of Lemma 3.9 of \cite{S} (with modifications as in \cite{W1}). 
\end{proof}

We shall need the following easy consequence of the preceding corollary for the proof of Theorem~\ref{no-transverse-q} at the end of this section.

\begin{corollary}\label{MD-height-est}
Let $q$ be an integer $\geq 2$, $\a \in (0, 1),$ $\g \in (0, 1/2),$ $\e^{\prime} \in (0, 1/2)$ and suppose that the induction hypotheses $(H1)$, $(H2)$ hold.  Let ${\mathbf C}_{0}$ be the stationary cone as in \ref{MD-1-0}, with $\Theta \, (\|{\mathbf C}_{0}\|, 0) = q + 1/2$. There exists a number 
$\e_{1} = \e_{1}(n, \a, \g, \e^{\prime}, {\mathbf C}_{0}) \in (0, 1/2)$ such that if $V \in {\mathcal S}_{\a},$
${\mathbf C} \in {\mathcal K}$ satisfy Hypotheses~\ref{MD-hyp} with $\e_{1}$ in place of $\e,$ then 
$$({\rm a}) \;\;\;\; \int_{B_{1}^{n+1}(0)} {\rm dist}^{2} \, (X, {\rm spt} \,\|{\mathbf C}\|) \, d\|V\|(X)
+  \int_{B_{1}^{n+1}(0) \setminus \{r(X) < 1/16\}} {\rm dist}^{2} \, (X, {\rm spt} \,\|V\|) \, d\|{\mathbf C}\|(X) < \e^{\prime}$$
and for each $Z = (\z^{1}, \z^{2}, \eta) \in {\rm spt} \, \|V\| \cap (B_{3/8}^{n+1}(0))$ with $\Theta \, (\|V\|, Z) \geq \Theta \, (\|{\mathbf C}_{0}\|, 0),$ we have that 
\begin{eqnarray*}
&&({\rm b}) \;\;\;\; |\z^{1}|^{2} + |\z^{2}|^{2} \leq C\left(\int_{B_{1}^{n+1}(0)} {\rm dist}^{2} \, (X, {\rm spt} \,\|{\mathbf C}\|) \, d\|V\|(X)\right.\nonumber\\ 
&&\hspace{2in} +  \left.\int_{B_{1}^{n+1}(0) \setminus \{r(X) < 1/16\}} {\rm dist}^{2} \, (X, {\rm spt} \,\|V\|) \, d\|{\mathbf C}\|(X)\right)
\end{eqnarray*}
\begin{eqnarray*}
&&({\rm c}) \;\;\int_{B_{1}^{n+1}(0)} {\rm dist}^{2} \, (X, {\rm spt} \,\|{\mathbf C}\|) \, d\|V^{Z}\|(X)
+ \int_{B_{1}^{n+1}(0) \setminus \{r(X) < 1/16\}} {\rm dist}^{2} \, (X, {\rm spt} \,\|V^{Z}\|) \, d\|{\mathbf C}\|(X)\nonumber\\
&&\hspace{1in} \leq C\left(\int_{B_{1}^{n+1}(0)} {\rm dist}^{2} \, (X, {\rm spt} \,\|{\mathbf C}\|) \, d\|V\|(X)\right.\nonumber\\ 
&&\hspace{2.5in} +  \left.\int_{B_{1}^{n+1}(0) \setminus \{r(X) < 1/64\}} {\rm dist}^{2} \, (X, {\rm spt} \,\|V\|) \, d\|{\mathbf C}\|(X)\right)
\end{eqnarray*}
where $V^{Z} = \eta_{Z, 1/2 \, \#} \, V$ and $C = C(n, \a, \g, {\mathbf C}_{0}) \in (0, \infty).$
\end{corollary}

\begin{proof}
Conclusion (a) is easily seen by arguing by contradiction using Allard's integral varifold compactness theorem (\cite{AW}; also \cite{S1}, Section 42.8). Conclusion (b) in case ${\mathbf C} \in {\mathcal K}(m_{0})$ follows directly from Corollary~\ref{MD-L2-est-2}. So suppose that ${\mathbf C} \not\in {\mathcal K}(m_{0}).$ Noting in this case  
that $2\Theta(\|{\mathbf C}_{0}\|, 0) \geq m_{0}+1$, fix $p \in \{m_{0}+1, m_{0}+2, \ldots, 2\Theta(\|{\mathbf C}_{0}\|, 0)\}$ and assume by induction that conclusion (b) of the corollary (taken with $\e^{\prime} = 1/4$, say) holds whenever ${\mathbf C} \in \cup_{j=m_{0}}^{p-1}{\mathcal K}(j),$ with $\overline{\e}$ denoting the required value of $\e_{1}$. Choose a cone $\widetilde{\mathbf C}_{1} \in \cup_{j=m_{0}}^{p-1} {\mathcal K}(j)$ such that 
\begin{eqnarray*}
&&\int_{B_{1}^{n+1}(0)} {\rm dist}^{2} \, (X, {\rm spt} \,\|\widetilde{\mathbf C}_{1}\|) \, d\|V\|(X) +  \int_{B_{1}^{n+1}(0) \setminus \{r(X) < 1/16\}} {\rm dist}^{2} \, (X, {\rm spt} \,\|V\|) \, d\|\widetilde{\mathbf C}_{1}\|(X)\nonumber\\
&&\hspace{1in} \leq \; \frac{3}{2}\; {\rm inf}_{\widetilde{\mathbf C} \in \cup_{j=m_{0}}^{p-1} {\mathcal K}(j)} 
\left(\int_{B_{1}^{n+1}(0)} {\rm dist}^{2} \, (X, {\rm spt} \,\|\widetilde{\mathbf C}\|) \, d\|V\|(X)\right.\nonumber\\
&&\hspace{2in}+ \left.\int_{B_{1}^{n+1}(0) \setminus \{r(X) < 1/16\}} {\rm dist}^{2} \, (X, {\rm spt} \,\|V\|) \, d\|\widetilde{\mathbf C}\|(X)   \right) 
\end{eqnarray*}
and let $\b_{1} = \frac{2}{3}\b(n, \a, \g, 1/4, 1/4, {\mathbf C}_{0})$ where $\b$ is as in Corollary~\ref{MD-L2-est-2}. Suppose ${\mathbf C} \in {\mathcal K}(p)$ and that Hypotheses~\ref{MD-hyp} hold with 
the value of $\e$ equal to $\e(n, \a, \g, 1/4, 1/4, {\mathbf C}_{0})$ where $\e(n, \a, \g, \cdot, \cdot,{\mathbf C}_{0})$ is as in Corollary~\ref{MD-L2-est-2}. If  
\begin{eqnarray*}
&&\int_{B_{1}^{n+1}(0)} {\rm dist}^{2} \, (X, {\rm spt} \,\|{\mathbf C}\|) \, d\|V\|(X) +  \int_{B_{1}^{n+1}(0) \setminus \{r(X) < 1/16\}} {\rm dist}^{2} \, (X, {\rm spt} \,\|V\|) \, d\|{\mathbf C}\|(X)\nonumber\\
&&\hspace{1in} \leq \; \b_{1}\; \left(\int_{B_{1}^{n+1}(0)} {\rm dist}^{2} \, (X, {\rm spt} \,\|\widetilde{\mathbf C}_{1}\|) \, d\|V\|(X)\right.\nonumber\\
&&\hspace{2in}+ \left.\int_{B_{1}^{n+1}(0) \setminus \{r(X) < 1/16\}} {\rm dist}^{2} \, (X, {\rm spt} \,\|V\|) \, d\|\widetilde{\mathbf C}_{1}\|(X)   \right), 
\end{eqnarray*}
then conclusion (b) follows directly from Corollary~\ref{MD-L2-est-2}; if on the other hand the reverse inequality holds, then by taking 
$\e^{\prime} = \e^{\prime}(n, \a, \g, {\mathbf C}_{0})$  sufficiently small in conclusion (a), 
we can ensure that Hypotheses~\ref{MD-hyp} are satisfied with $\widetilde{\mathbf C}_{1}$ in place of ${\mathbf C}$ and $\overline{\e}$ in place of $\e,$ so conclusion (b) in this case follows by the induction hypothesis. Thus conclusion (b) holds whenever ${\mathbf C} \in {\mathcal K}(p)$, and since ${\mathbf C} \in {\mathcal K} \implies {\mathbf C} \in {\mathcal K}(j)$ for some $j \in \{m_{0}, \ldots, 2\Theta(\|{\mathbf C}_{0}\|, 0)\}$, the inductive proof of conclusion (b) is complete.

To see conclusion (c), note that  
\begin{eqnarray*}
&&\int_{B_{1}^{n+1}(0)} {\rm dist}^{2} \, (X, {\rm spt} \,\|{\mathbf C}\|) \, d\|V^{Z}\|(X) =
{2}^{n+2}\int_{B_{1/2}^{n+1}(Z)} {\rm dist}^{2} \, (X, T_{Z}\, {\rm spt} \, \|{\mathbf C}\|) \, d\|V\|(X)\nonumber\\
&&\hspace{2in}\leq 2^{n+2}\int_{B_{1}^{n+1}(0)}{\rm dist}^{2} \, (X, {\rm spt} \,\|{\mathbf C}\|) \, d\|V\|(X) + C\left(|\z^{1}|^{2} + |\z^{2}|^{2}\right)
\end{eqnarray*}
and 
\begin{eqnarray*}
&&\int_{B_{1}^{n+1}(0) \setminus \{r(X) < 1/16\}} {\rm dist}^{2} \, (X, {\rm spt} \,\|V^{Z}\|) \, d\|{\mathbf C}\|(X)\nonumber\\
&&\hspace{1in}=2^{n+2}\int_{B_{1/2}^{n+1}(Z) \setminus\{r(X - Z) < 1/32\}} {\rm dist}^{2} \, (X, {\rm spt} \, \|V\|) \, d\|T_{Z \, \#} \, {\mathbf C}\|(X)\nonumber\\
&&\hspace{1in}\leq2^{n+2}\int_{B_{1/2}^{n+1}(Z) \setminus\{r(X) < 1/64\}} {\rm dist}^{2} \, (X, {\rm spt} \, \|V\|) \, d\|T_{Z \, \#} \, {\mathbf C}\|(X)\nonumber\\
&&\hspace{1in}\leq2^{n+2}\int_{B_{1}^{n+1}(0) \setminus \{r(X) < 1/64\}}{\rm dist}^{2} \, (X, {\rm spt} \, \|V\|)d\|{\mathbf C}\|(X) + C\left(|\z^{1}|^{2} + |\z^{2}|^{2}\right)
\end{eqnarray*}
where $C = C(n, q) \in (0, \infty)$, $T_{Z} \, : \, {\mathbf R}^{n+1} \to {\mathbf R}^{n+1}$ is the translation $X \mapsto X + Z$ and we have used the fact that ${\mathbf C}$ is translation invariant along $\{0\}\times {\mathbf R}^{n-1}$ and assumed that $\e = \e(n, \a, \g, {\mathbf C}_{0})$ is sufficiently small to ensure that ${\rm dist} \, (Z, \{0\} \times {\mathbf R}^{n-1}) < 1/64.$ In view of conclusion (b), the validity of conclusion (c) readily follows from these two inequalities.
\end{proof}

%\medskip

\begin{lemma}\label{MD-no-gaps}
Let $q$ be an integer $\geq 2$, $\a \in (0, 1)$, $\d \in (0, 1/8)$, $\g \in (0, 1/2)$ and ${\mathbf C}_{0}$ be as above. Suppose that the induction hypotheses $(H1)$, $(H2)$ hold and that $\Theta \, (\|{\mathbf C}_{0}\|,0) = q + 1/2.$ There exist numbers $\e_{1} = \e_{1}(n, \a, \g,\d,  {\mathbf C}_{0}) \in (0, 1/2)$ and $\b_{1} = \b_{1}(n, \a, \g, {\mathbf C}_{0}) \in (0, 1/2)$ such that if $V \in {\mathcal S}_{\a},$ ${\mathbf C} \in {\mathcal K}$ satisfy Hypotheses~\ref{MD-hyp} with $\e_{1}$ in place of $\e$ then 
\begin{itemize}
\item[(a)] $B^{n+1}_{\d}(0,y) \cap \{Z \, : \, \Theta \, (\|V\|, Z) \geq q + 1/2\} \neq \emptyset$
for each point $(0, y) \in \{0\} \times {\mathbf R}^{n-1} \cap B_{1/2}^{n+1}(0),$ and\\
\item[(b)] if additionally Hypothesis {\rm (}$\dag${\rm )} holds with $\b_{1}$ in place of $\b$ and if $\m \in (0, 1)$, then 
\begin{eqnarray*}
\int_{B_{1/2}^{n+1}(0) \cap \{r(X) < \s\}} {\rm dist}^{2} \, (X, {\rm spt} \, \|{\mathbf C}\|) \, d\|V\|(X)&&\nonumber\\ 
&&\hspace{-1in}\leq C_{1}\s^{1-\m}\int_{B_{1}^{n+1}(0)} {\rm dist}^{2} \, (X, {\rm spt} \, \|{\mathbf C}\|) \, d\|V\|(X)
\end{eqnarray*}
for each $\s \in [\d, 1/4)$, where $C_{1} = C_{1}(n, q, \a, \m, {\mathbf C}_{0}) \in (0, \infty).$ (In particular $C_{1}$ is independent of $\d$ and $\s$.)
\end{itemize}
\end{lemma}

\begin{proof}
Suppose for some number $\d \in (0, 1/8)$ and some point $(0, y) \in \{0\} \times {\mathbf R}^{n-1} \cap B_{1/2}^{n+1}(0)$, 
$B^{n+1}_{\d}(0,y) \cap \{Z \, : \, \Theta \, (\|V\|, Z) \geq q + 1/2\} = \emptyset.$ Then by Remark (3) following 
Hypotheses~\ref{MD-hyp},  it follows that 
\begin{equation}\label{MD-no-gaps-0}
{\mathcal H}^{n-7+\g} \, ( {\rm sing} \, V \res (B_{\d}^{n+1}(0, y))) = 0 \;\; \mbox{if $n \geq 7$ and} \;\; 
{\rm sing} \, V \res(B_{\d}^{n+1}(0, y)) = \emptyset \;\; \mbox{if $2 \leq n \leq 6.$}
\end{equation}

From this and  hypothesis $({\mathcal S{\emph 2}})$ we deduce (with the help of  an elementary covering argument in case $n \geq 7$) that 
\begin{equation*}
\int_{{\rm spt} \, \|V\| \cap B_{\d}^{n+1}(0, y)} |A|^{2} \z^{2} \, d{\mathcal H}^{n} \leq \int_{{\rm spt} \, \|V\| \cap B_{\d}^{n+1}(0, y)} |\nabla^{V} \, \z|^{2}\, d{\mathcal H}^{n}
\end{equation*}
for any  $\z \in C^{1}_{c}(B_{\d}^{n+1}(0, y))$, where $A$ denotes 
the second fundamental form of ${\rm reg} \, V.$ Choosing $\z \in C^{1}_{c}(B_{\d}^{n+1}(0, y))$ such that 
$\z \equiv 1 $ in $B_{\d/2}^{n+1}(0, y)$ and $|D \, \z| \leq 4\d^{-1}$, we conclude from the preceding inequality that 
\begin{equation}\label{MD-no-gaps-1}
\int_{{\rm spt} \, \|V\| \cap B_{\d/2}^{n+1}(0, y)} |A|^{2} \, d{\mathcal H}^{n}\leq C\d^{n-2}
\end{equation}
where $C = C(n, q) \in (0, \infty).$ Now let $\t \in (0, \d/4)$ be arbitrary for the moment and assume that $\e \in (0, \e_{0})$, where $\e_{0} = \e_{0}(\a, \b, \t, {\mathbf C}_{0})$ is as in Theorem~\ref{MD-L2-est-1}. Using Theorem~\ref{MD-L2-est-1} (a), (\ref{MD-no-gaps-1}) and the argument leading to the inequality (6.12) of \cite{SS} (with $\s = \t$), we deduce, provided $\e = \e(\a, \b, \t, {\mathbf C}_{0})$ is sufficiently small and positive, that
$$C \leq \t^{1/2} \d^{-1/2},$$ 
where $C = C(\b, {\mathbf C}_{0}) \in (0, \infty)$. This however is a contradiction if we choose e.g.\ $\t = C^{2}\d^{2},$ and we conclude that part (a) must hold provided $\e = \e(\a, \b, \d, {\mathbf C}_{0}) \in (0, 1/2)$ is sufficiently small. To prove the estimate of part (b), first note that in view of Corollary~\ref{MD-L2-est-2} (a),(b) (with $\t = 1/16$, say), it follows from the argument leading to the estimate (3) of \cite{S}, p.619 that for each $Z \in {\rm spt} \, \|V\| \cap B_{3/8}^{n+1}(0)$ with $\Theta \, (\|V\|, Z) \geq q$, 
$$\int_{B_{1/4}^{n+1}(Z)}  \frac{{\rm dist}^{2} \, (X, {\rm spt} \, \|{\mathbf C}\|)}{|X - Z|^{n-\a}} \, d\|V\|(X) \leq C\int_{B_{1}^{n+1}(0)} {\rm dist}^{2} \, (X, {\rm spt} \, \|{\mathbf C}\|) \, d\|V\|(X)$$
where $C = C(\b, \a, {\mathbf C}_{0}) \in (0, \infty).$ By the argument of \cite{S}, Corollary 3.2 (ii) (cf. proof of Lemma~\ref{no-gaps}(b)), the required estimate follows from this and part (a). 
\end{proof}

\noindent
{\bf Remark:} Note that Theorem~\ref{MD-L2-est-1}, Corollary~\ref{MD-L2-est-2}, Corollary~\ref{MD-height-est} and Lemma~\ref{MD-no-gaps} all continue to hold in case $\Theta \, (\|{\mathbf C}_{0}\|, 0) = q+1$ provided that Theorem~\ref{sheetingthm} holds with $q+1$ in place of $q$.\\ 

Let $\g \in (0, 1/2)$ and consider a sequence of varifolds $\{V_{k}\} \subset {\mathcal S}_{\a}$ and a sequence of cones 
$\{{\mathbf C}_{k}\}$ satisfying, for each $k=1, 2, \ldots,$ Hypotheses~\ref{MD-hyp} and Hypothesis ($\dag$) with $V_{k}$, ${\mathbf C}_{k}$ in place of $V$, ${\mathbf C}$ and $\e_{k}$, $\b_{k}$ in place of $\e$, $\b$, where 
$\e_{k}, \b_{k} \to 0^{+}.$  Thus, we suppose, for each $k=1, 2, \ldots,$ 
\begin{itemize}
\item[($1_{k}$)] $V_{k} \in {\mathcal S}_{\a}$, $0 \in {\rm spt} \, \|V_{k}\|$, $\Theta \, (\|V_{k}\|, 0) \geq q+1/2$, 
$(\omega_{n}2^{n})^{-1} \|V_{k}\|(B_{2}^{n+1}(0)) < q +1/2 + \g;$
\item[($2_{k}$)] ${\mathbf C}_{k}  = \sum_{j=1}^{m_{0}} \sum_{\ell=1}^{q_{j}^{(0)}} |H^{k}_{j, \, \ell}| \in {\mathcal K},$ where $H^{k}_{j, \, \ell}$ are half-hyperplanes of ${\mathbf R}^{n+1}$ meeting along $\{0\} \times {\mathbf R}^{n-1}$ with $H^{k}_{j, \, \ell} \in N(H_{j}^{(0)})$ for each $j \in \{1, 2, \ldots, m_{0}\}$ and $\ell \in\{1, 2, \ldots, q_{j}^{(0)}\};$   
\item[($3_{k}$)] ${\rm dist}_{\mathcal H} \, ({\rm spt} \, \|{\mathbf C}_{k}\| \cap B_{1}^{n+1}(0), {\rm spt} \, \|{\mathbf C}_{0}\| \cap B_{1}^{n+1}(0)) < \e_{k};$
\item[($4_{k}$)] $$\int_{B_{1}^{n+1}(0)} {\rm dist}^{2} \, (X, {\rm spt} \,\|{\mathbf C}_{k}\|) \, d\|V_{k}\|(X)< \e_{k};$$
\item[($5_{k}$)] For each $j=1,2, \ldots, m_{0},$ 
$$\|V_{k}\|((B_{1/2}^{n+1}(0) \setminus \{r(X) < 1/8\}) \cap N(H^{(0}_{j})) \geq \left(q_{j}^{(0)} - \frac{1}{4}\right){\mathcal H}^{n}((B_{1/2}^{n+1}(0) \setminus \{r(X) < 1/8\}) \cap H^{(0)}_{j});$$
\item[($6_{k}$)] Either (i) or (ii) below holds:
\begin{itemize}
\item[(i)] ${\mathbf C}_{k} \in {\mathcal K}(m_{0});$ 
\item[(ii)] $2\Theta \, (\|{\mathbf C}_{0}\|, 0) \geq m_{0}+1$, ${\mathbf C}_{k} \in {\mathcal K}(p_{k})$ for some $p_{k} \in \{m_{0}+1, m_{0}+2, \ldots, 2\Theta \, (\|{\mathbf C}_{0}\|, 0)\}$ and 
\begin{eqnarray*}
\int_{B_{1}^{n+1}(0)} {\rm dist}^{2} \, (X, {\rm spt} \,\|{\mathbf C}_{k}\|) \, d\|V_{k}\|(X) +  \int_{B_{1}^{n+1}(0) \setminus \{r(X) < 1/16\}} {\rm dist}^{2} \, (X, {\rm spt} \,\|V_{k}\|) \, d\|{\mathbf C}_{k}\|(X)&&\\
&&\hspace{-5in} \leq \; \b_{k}\; {\rm inf}_{\widetilde{\mathbf C} \in \cup_{j=m_{0}}^{p_{k}-1} {\mathcal K}(j)} 
\left(\int_{B_{1}^{n+1}(0)} {\rm dist}^{2} \, (X, {\rm spt} \,\|\widetilde{\mathbf C}\|) \, d\|V_{k}\|(X)\right.\\
&&\hspace{-3.5in}+\;  \left.\int_{B_{1}^{n+1}(0) \setminus \{r(X) < 1/16\}} {\rm dist}^{2} \, (X, {\rm spt} \,\|V_{k}\|) \, d\|\widetilde{\mathbf C}\|(X)   \right). 
\end{eqnarray*} 
\end{itemize} 
\end{itemize}
Note that it follows from ($2_{k}$) and ($3_{k}$) that 
$H^{k}_{j, \, \ell} \to H_{j}^{(0)}$ for each $j \in \{1, \ldots, m_{0}\}$ and $\ell \in \{1, \ldots,q_{j}^{(0)}\}$.

Let ${\E}_{k} = \sqrt{\int_{B_{1}^{n+1}(0)} {\rm dist}^{2} \, (X, {\rm spt} \,\|{\mathbf C}_{k}\|) \, d\|V_{k}\|(X)}$.

Let $\{\d_{k}\}, \{\t_{k}\}$ be sequences of decreasing positive numbers converging to 0. By passing to appropriate 
subsequences of $\{V_{k}\}$, $\{{\mathbf C}_{k}\}$ without changing notation, we have the following:
 \begin{itemize}
\item[(A$_{k}$)] By Lemma~\ref{MD-no-gaps},
\begin{equation}\label{MD-excess-1}
B^{n+1}_{\d_{k}}(0,y) \cap \{Z \, : \, \Theta \, (\|V_{k}\|, Z) \geq q + 1/2\} \neq \emptyset
\end{equation}
for each point $(0, y) \in \{0\} \times {\mathbf R}^{n-1} \cap B_{1/2}^{n+1}(0)$ and
\begin{equation}\label{MD-excess-2}
\int_{B_{1/2}^{n+1}(0) \cap \{r(X) < \s\}} {\rm dist}^{2} \, (X, {\rm spt} \, \|{\mathbf C}_{k}\|) \, d\|V_{k}\|(X) 
\leq C\s^{1/2} {\E}_{k}^{2}
\end{equation}
for each $\s \in [\d_{k}, 1/4),$ where $C= C(n,q, \a, \g, {\mathbf C}_{0}) \in (0, \infty).$ 
\item[(B$_{k}$)] By Theorem~\ref{MD-L2-est-1} (a), 
\begin{equation}\label{MD-excess-3}
V_{k} \res (B^{n+1}_{7/8}(0) \setminus \{r(X) < \t_{k}\}) = \sum_{j=1}^{m_{0}} \sum_{\ell=1}^{q_{j}^{(0)}}
|{\rm graph} \, u_{j, \, \ell}^{k}|
\end{equation}
where, for each $k=1, 2, \ldots$,  $j \in \{1, 2, \ldots, m_{0}\}$ and $\ell \in \{1, 2, \ldots, q_{j}^{(0)}\},$ 
$u_{j,\,  \ell}^{k} \in C^{2} \, (B_{7/8}^{n+1}(0) \,\cap \, H_{j, \, \ell}^{k} \setminus \{r(X) < \t_{k}\}; (H_{j, \, \ell}^{k})^{\perp}),$  
$u_{j, \, \ell}^{k}$ solves the minimal surface equation on $B_{7/8}^{n+1}(0) \cap H_{j, \, \ell}^{k} \setminus \{r(X) < \t_{k}\}$ and 
satisfies ${\rm dist} \, (X + u_{j, \, \ell}^{k}(X), {\rm spt} \,\|{\mathbf C}\|) = |u_{j, \, \ell}^{k}(X)|$ for $X \in B_{7/8}^{n+1}(0) \cap H_{j, \, \ell}^{k} \setminus \{r(X)  < \t_{k}\}.$ 
\item[(C$_{k}$)] For each  point $Z = (\z^{1}, \z^{2}, \eta) \in {\rm spt} \, \|V_{k}\| \cap B_{3/8}^{n+1}(0)$ with $\Theta \, (\|V_{k}\|, Z) \geq q +1/2$, by Corollary~\ref{MD-L2-est-2} (a), 
\begin{equation}\label{MD-excess-4}
|\z^{1}|^{2} + |\z^{2}|^{2} \leq C {\E}_{k}^{2}
\end{equation}
where $C = C(n, q, \a, \g, {\mathbf C}_{0}) \in (0, \infty).$  
\item[(D$_{k}$)]  For each fixed $\m \in (0, 1)$, $\r \in (0, 1/4)$, each sufficiently large $k$  and each point $Z = (\z^{1}, \z^{2}, \eta) \in {\rm spt} \, \|V_{k}\| \cap B_{3/8}^{n+1}(0)$ with $\Theta \, (\|V_{k}\|, Z) \geq q + 1/2$, by Corollary~\ref{MD-L2-est-2}(b), 
\begin{eqnarray}\label{MD-excess-5}
\sum_{j=1}^{m_{0}} \sum_{\ell=1}^{q_{j}^{(0)}}\int_{B_{1/4}^{n+1}(Z) \cap H_{j, \, \ell}^{k} \setminus 
\{r(X) < \t_{k}\}} \frac{|u_{j, \, \ell}^{k}(X) - (\z^{1}, \z^{2}, 0)^{\perp_{H_{j}^{k}}}|^{2}}{|X+ u_{j,\, \ell}^{k}(X) - Z|^{n+2-\m}} \, dX &&\nonumber\\
&& \hspace{-2in}\leq C_{1} \r^{-n-2+ \m}\int_{B_{\r}^{n+1}(Z)} {\rm dist}^{2} \, (X, {\rm spt} \, \|T_{Z \, \#} \, {\mathbf C}_{k}\|) \, d\|V_{k}\|(X),
\end{eqnarray}
where $C_{1}  = C_{1}(\a, \g, \m, {\mathbf C}_{0}) \in (0, \infty).$ 
\end{itemize}

Extend $u_{j, \ell}^{k}$ to all of $B_{7/8}^{n+1}(0) \cap H_{j, \, \ell}^{k}$  by defining values to be zero in 
$B_{7/8}^{n+1}(0) \cap H_{j, \, \ell}^{k} \cap \{r(X) < \t_{k}\}.$ For each $j \in \{1, 2, \ldots, m_{0}\}$ and 
$\ell \in \{1, 2, \ldots, q_{j}^{(0)}\}$, let $h_{j, \, \ell} \, : \, H_{j}^{(0)} \to (H_{j}^{(0)})^{\perp}$
be the linear functions such that $\{X + h_{j, \, \ell}(X) \, : \, X \in H_{j}^{(0)}\} = H_{j, \, \ell}$ and 
let $\widetilde{u}_{j, \, \ell}^{k}(X) = u_{j, \, \ell}^{k}(X + h_{j, \, \ell}(X)).$ By (\ref{MD-excess-3}) and elliptic estimates, there exist, for each $j=1, 2, \ldots, m_{0}$ and 
 $\ell =1, 2, \ldots, q_{j}^{(0)}$ (for any manner in which the labelling is chosen for the elements of the 
 sets $\{u_{j, \, 1}^{k}, u_{j, \, 2}^{k}, \ldots, u_{j, \, q_{j}^{(0)}}^{k}\}, \;\;\; k=1, 2, 3, \ldots)$ harmonic functions 
$v_{j, \, \ell} \, : \, B_{3/4} \cap H_{j}^{(0)} \to (H_{j}^{(0)})^{\perp}$  such that, after passing to a subsequence,
\begin{equation}\label{MD-excess-8}
{\E}_{k}^{-1}{\widetilde u}_{j, \, \ell}^{k} \to v_{j, \, \ell} 
\end{equation}
where the convergence is in $C^{2}(K)$ for each compact subset $K$ of $B_{3/4} \cap H_{j}^{(0)}.$ From (\ref{MD-excess-2}), it follows that for each $\s \in (0, 1/4)$, 
$$\sum_{j=1}^{m_{0}}\sum_{\ell=1}^{q_{j}^{(0)}}
\int_{B_{3/4}^{n+1}(0) \cap H_{j}^{(0)} \cap \{r(X) < \s\}}|v_{j, \, \ell}|^{2} \leq  C\s^{1/2}, \;\;\; C = C(\a, \g, {\mathbf C}_{0})$$
and hence that the convergence in (\ref{MD-excess-8}) is also in $L^{2}\,(B_{3/4} \cap H_{j}^{(0)}).$

\begin{lemma}\label{MD-continuity}
For each $j \in \{1, 2, \ldots, m_{0}\}$ and $\ell \in \{1, 2, \ldots, q_{j}^{(0)}\}$, we have that 
$$v_{j, \, \ell} \in C^{0, \m} \,\left(\overline{B_{5/16}^{n+1}(0) \cap H_{j}^{(0)}}; \left(H_{j}^{(0)}\right)^{\perp}\right)$$ 
\noindent
with the estimate 
\begin{eqnarray*}
\sup_{\overline{B_{5/16}^{n+1}(0) \cap H_{j}^{(0)}}} \, |v_{j, \, \ell}|^{2}  \; + \;
\sup_{X_{1}, X_{2} \in \overline {B_{5/16}^{n+1}(0) \cap H_{j}^{(0)}}, \, X_{1} \neq X_{2}} \, 
\frac{|v_{j, \, \ell}(X_{1}) - v_{j,  \, \ell}(X_{2})|^{2}}{|X_{1} - X_{2}|^{2\m}}&&\nonumber\\
&&\hspace{-2in} \leq C\sum_{j=1}^{m_{0}} \sum_{\ell=1}^{q_{j}^{(0)}}\int_{B_{3/4}^{n+1}(0) \cap H_{j}^{(0)}}|v_{j, \, \ell}|^{2}
\end{eqnarray*}
where $\m  = \m(n, q, \a, \g, {\mathbf C}_{0}) \in (0, 1)$ and $C = C(n, q, \a, \g, {\mathbf C}_{0}) \in (0, \infty).$
\end{lemma}

\begin{proof} Note first that for each given $Y \in B_{5/16}^{n+1}(0) \cap \{0\} \times {\mathbf R}^{n-1}$, there exists, by (\ref{MD-excess-1}), a sequence of points $Z_{k} = (\z_{1}^{k}, \z_{2}^{k}, \eta_{k}) \in {\rm spt} \, \|V_{k}\| \cap B_{3/4}^{n+1}(0)$ with $\Theta \, (\|V_{k}\|, Z_{k}) \geq q$ 
such that $Z_{k} \to Y.$ Passing  to a subsequence without changing notation, the limits 
$\lim_{k \to \infty} \, {\E}_{k}^{-1}\z_{1}^{k}$ and  $\lim_{k \to \infty} \, {\E}_{k}^{-1}\z_{2}^{k}$ exist by 
(\ref{MD-excess-4}). Write
\begin{equation}\label{MD-continuity-1} 
\k(Y) = \left(\lim_{k \to \infty} \, {\E}_{k}^{-1}\z_{1}^{k}, \lim_{k \to \infty} \, {\E}_{k}^{-1}\z_{2}^{k}, 0 \right)
\end{equation}
and note by (\ref{MD-excess-4}) that 
\begin{equation}\label{MD-continuity-1'}
|\k(Y)| \leq C
\end{equation}
where $C = C(n, q, \a, \g, {\mathbf C}_{0}) \in (0, \infty).$ It follows from (\ref{MD-excess-2}), (\ref{MD-excess-3}) and (\ref{MD-excess-5}) that for each $\m \in (0, 1)$, 
\begin{eqnarray}\label{MD-continuity-2}
\sum_{j=1}^{m_{0}} \sum_{\ell=1}^{q_{j}^{(0)}}\int_{B_{1/4}^{n+1}(Y) \cap H_{j}^{(0)}} \frac{|v_{j, \, \ell}(X) - \k(Y)^{\perp_{H_{j}^{(0)}}}|^{2}}{|X  - Y|^{n+2-\m}} \, dX &&\nonumber\\
&& \hspace{-2in}\leq C_{1} \r^{-n-2+ \m}\sum_{j=1}^{m_{0}}\sum_{\ell = 1}^{q_{j}^{(0)}}\int_{B_{\r}^{n+1}(Y) \cap H_{j}^{(0)}}|v_{j, \, \ell} -\k(Y)^{\perp_{H_{j}^{(0)}}}|^{2} 
\end{eqnarray}
for $\r \in (0, 1/8]$, where $C_{1}  = C_{1}(\a,\g, \m, {\mathbf C}_{0}) \in (0, \infty).$ In view of 
(\ref{MD-continuity-1'}), this in particular implies that for each $j=1, 2, \ldots, m_{0}$, $\k(Y)^{\perp_{H_{j}^{(0)}}}$ is uniquely defined (depending only on $Y$ and independent of the sequence $\{Z_{k}\}$ 
tending to $Y$), and hence, since the set of normal directions to $H_{j}^{(0)}$, $j=1, 2, \ldots, m_{0},$ spans ${\mathbf R}^{2} \times \{0\}$, the vector $\k(Y)$ is also uniquely defined. For $Y \in B_{1/4}^{n+1}(0) \cap \{0\} \times {\mathbf R}^{n-1}$, $j \in \{1, 2, \ldots, m_{0}\}$ and $\ell \in \{1, 2, \ldots, q_{j}^{(0)}\}$, define $v_{j, \ell}(Y) = \k(Y)^{\perp_{H_{j}^{(0)}}}.$ The proof of the lemma can now be completed by modifying the proof of Lemma~\ref{continuity} in an obvious way.
\end{proof}

\begin{theorem}\label{MD-c1alpha}
For each $j \in \{1, 2, \ldots, m_{0}\}$, $\ell \in \{1, 2, \ldots, q_{j}^{(0)}\}$, we have that 
$$v_{j, \, \ell} \in C^{2} \,\left(\overline{B_{1/4}^{n+1}(0) \cap H_{j}^{(0)}}; \left(H_{j}^{(0)}\right)^{\perp}\right)$$
\noindent
with the estimate 
\begin{eqnarray*}
\sup_{\overline{B_{1/4}^{n+1}(0) \cap H_{j}^{(0)}}} \, |Dv_{j, \, \ell}|^{2}  \; + \;
\sup_{X_{1}, X_{2} \in \overline {B_{1/4}^{n+1}(0) \cap H_{j}^{(0)}}, \, X_{1} \neq X_{2}} \, 
\frac{|Dv_{j, \, \ell}(X_{1}) - Dv_{j,  \, \ell}(X_{2})|^{2}}{|X_{1} - X_{2}|^{2}}&&\nonumber\\
&&\hspace{-2in} \leq C\sum_{j=1}^{m_{0}} \sum_{\ell=1}^{q_{j}^{(0)}}\int_{B_{3/4}^{n+1}(0) \cap H_{j}^{(0)}}|v_{j, \, \ell}|^{2}
\end{eqnarray*}
where $C = C(n, q, \a, \g, {\mathbf C}_{0}) \in (0, \infty).$
\end{theorem}

\begin{proof} For $Y \in B_{1/2}^{n+1}(0) \cap \{0\} \times {\mathbf R}^{n-1}$, let $\widetilde{\k}(Y)  = \sum_{j=1}^{m_{0}} \k(Y)^{\perp_{H_{j}^{(0)}}}$ 
where $\k$ is the function defined by (\ref{MD-continuity-1}). By modifying the argument leading to the estimate (\ref{c1alpha-16})  in obvious ways, 
it can be seen that $\widetilde{\k} \in C^{\infty} \,(B_{1/2}^{n+1}(0) \cap \{0\} \times {\mathbf R}^{n-1}; {\mathbf R}^{n+1})$ with 
\begin{equation*}
\sup_{B_{1/2}^{n+1}(0) \cap (\{0\} \times {\mathbf R}^{n-1})}\, |\widetilde{\k}|^{2} + |D_{Y} \,\widetilde{\k}|^{2} + |D_{Y}^{2} \, \widetilde{\k}|^{2} +  |D_{Y}^{3} \, \widetilde{\k}|^{2} \leq C\sum_{j=1}^{m_{0}}\sum_{\ell=1}^{q_{j}^{(0)}}\int_{B_{3/4}^{n+1}(0) \cap H_{j}^{(0)}}|v_{j, \, \ell}|^{2} 
\end{equation*}
where $C = C(\a, \g, {\mathbf C}_{0}) \in (0, \infty).$ Since 
the set of normal directions to $H_{j}^{(0)}$, $j=1, 2, \ldots, m_{0},$ span ${\mathbf R}^{2} \times \{0\}$, it follows that for each $j=1, 2, \ldots, m_{0},$
$\k^{\perp_{H_{j}^{(0)}}} \in C^{\infty} \,(B_{1/2}^{n+1}(0) \cap \{0\} \times {\mathbf R}^{n-1}; {\mathbf R}^{n+1})$
with 
\begin{eqnarray}\label{MD-c1alpha-1}
&\hspace{-1in}\sup_{B_{1/2}^{n+1}(0) \cap (\{0\} \times {\mathbf R}^{n-1})}\, |\k^{\perp_{H_{j}^{(0)}}}|^{2} + 
|D_{Y} \,\k^{\perp_{H_{j}^{(0)}}}|^{2} + |D_{Y}^{2} \, \k^{\perp_{H_{j}^{(0)}}}|^{2}  + 
|D_{Y}^{3} \, \k^{\perp_{H_{j}^{(0)}}}|^{2} \nonumber\\
&\hspace{4in}\leq C\sum_{j=1}^{m_{0}}\sum_{\ell=1}^{q_{j}^{(0)}}\int_{B_{3/4}^{n+1}(0) \cap H_{j}^{(0)}}|v_{j, \, \ell}|^{2} 
\end{eqnarray}
where $C = C(\a, \g, {\mathbf C}_{0}) \in (0, \infty).$ Since by Lemma~\ref{MD-continuity}, for each $j=1, 2, \ldots, m_{0}$ and $\ell =1, 2, \ldots, q_{j}^{(0)}$, $v_{j, \, \ell}$ is continuous in $\overline{B_{1/2}^{n+1}(0) \cap H_{j}^{(0)}}$ with boundary values 
$\left.v_{j, \, \ell} \right|_{B_{1/2}^{n+1}(0) \cap \{0\} \times {\mathbf R}^{n-1}} \equiv \k^{\perp_{H_{j}^{(0)}}}$, and 
$v_{j, \, \ell}$ is harmonic in $B_{3/4}^{n+1}(0) \cap H_{j}^{(0)}$, the desired conclusions of the lemma  follow, in view of the estimate (\ref{MD-c1alpha-1}), from the standard boundary regularity theory for harmonic functions.
\end{proof}

\begin{lemma}\label{MD-excess-improvement} 
Let $q$ be an integer $\geq 2$, $\a \in (0, 1)$, $\g \in (0, 1/2)$ and $\th \in (0, 1/4)$. Let ${\mathbf C}_{0}$ be the stationary cone as in (\ref{MD-1-0}), with $\Theta \, (\|{\mathbf C}_{0}\|, 0) = q + 1/2.$
There exist numbers $\overline\e = \overline\e(\a, \g, \th, {\mathbf C}_{0}) \in (0, 1/2)$ and 
$\overline\b = \overline\b(\a,\g, \th, {\mathbf C}_{0}) \in (0, 1/2)$ such that if $V \in {\mathcal S}_{\a},$ ${\mathbf C} \in {\mathcal K}$ satisfy Hypotheses~\ref{MD-hyp} and Hypothesis {\rm(}$\dag${\rm)} 
with $\e = \overline{\e}$ and $\b= \overline{\b}$ and  if the induction hypotheses $(H1)$, $(H2)$ hold, then 
there exist an orthogonal rotation $\G$ of ${\mathbf R}^{n+1}$ and a cone ${\mathbf C}^{\prime} \in {\mathcal K}$ 
such that, with 
$${\E}_{V}^{2} =  \int_{B_{1} ^{n+1}(0)} {\rm dist}^{2} \, (X, {\rm spt} \, \|{\mathbf C})\|) \, d\|V\|(X),$$ 
\noindent
the following hold:
\begin{eqnarray*}
&&({\rm a})\;\;\;\; |e_{j} - \G(e_{j})| \leq \overline{\k}{\E}_{V},\;\;\;{\rm for} \;\; j=1, 2, 3, \ldots, n+1;\nonumber\\ 
&&({\rm b})\;\;\;\; {\rm dist}_{\mathcal H}^{2} \, ({\rm spt} \, \|{\mathbf C}^{\prime}\| \cap B_{1}^{n+1}(0), {\rm spt} \, \|{\mathbf C}\| \cap B_{1}^{n+1}(0)) \leq \overline{C}_{0}{\E}_{V}^{2};\nonumber\\
&&({\rm c})\;\;\;\; \th^{-n-2}\int_{\G\left(B_{\th}^{n+1}(0) \setminus \{r(X) \leq \th/16\}\right)} {\rm dist}^{2} \, (X, {\rm spt} \, \|V\|) \, d\|\G_{\#} \, {\mathbf C}^{\prime}\|(X)\nonumber\\
&&\hspace{1.5in}+ \;\th^{-n-2}\int_{B_{\th}^{n+1}(0)} {\rm dist}^{2} \, (X, {\rm spt} \, \|\G_{\#} \, {\mathbf C}^{\prime}\|) \, d\|V\|(X)  \leq \overline{\nu}\th^{2}{\E}_{V}^{2}.
\end{eqnarray*}
Here the constants $\overline{\k}, \overline{C}_{0},  \overline{\nu} \in (0, \infty)$ each depends only on $\a$, $\g$ and 
${\mathbf C}_{0}.$ 
\end{lemma}

\begin{proof}
If the lemma is false, then there exist a sequence of varifolds $\{V_{k}\} \subset {\mathcal S}_{\a}$ and a sequence of cones $\{{\mathbf C}_{k}\} \subset {\mathcal K}$ satisfying, for each $k=1, 2, 3, \ldots,$ the conditions ($1_{k}$)--($6_{k}$) 
above (listed  immediately following the proof of Lemma~\ref{MD-no-gaps}), but not satisfying, for any choice of orthogonal rotation $\G$ of ${\mathbf R}^{n+1}$ and ${\mathbf C}^{\prime} \in {\mathcal K}$,  the conclusion of the Lemma taken with $V_{k}$, ${\mathbf C}_{k}$ in place of $V$, ${\mathbf C}$. 
Choose any two sequences of decreasing positive numbers $\{\d_{k}\}$ and $\{\t_{k}\}$ with $\d_{k} \to 0$ and $\t_{k} \to 0$ and corresponding subsequences of $\{V_{k}\}$, $\{{\mathbf C}_{k}\}$ for which the assertions (\ref{MD-excess-1})--(\ref{MD-excess-5}) are valid, and let 
$\{v_{j, \, \ell}\}_{j=1, 2, \ldots, m_{0}; \; \ell = 1, 2, \ldots, q_{j}^{(0)}}$ be the blow-up of $\{V_{k}\}$ relative to $\{{\mathbf C}_{k}\}$. Thus, for each $j  = 1, 2, \ldots, m_{0}$ and $\ell = 1, 2, \ldots, q_{j}^{(0)},$ 
$$v_{j, \, \ell} \in L^{2} \, \left(B_{3/4}^{n+1}(0) \cap H_{j}^{(0)}; \left(H_{j}^{(0)}\right)^{\perp}\right) \cap 
C^{2} \, \left(B_{3/4}^{n+1}(0) \cap H_{j}^{(0)}; \left(H_{j}^{(0)}\right)^{\perp}\right)$$
\noindent
are the functions produced as in (\ref{MD-excess-8}). Note then that Theorem~\ref{MD-c1alpha} is applicable to the functions $v_{j, \, \ell}.$  By following exactly the corresponding steps in the proof of Lemma~\ref{excess-improvement}, and using Theorem~\ref{MD-c1alpha} where the proof of Lemma~\ref{excess-improvement} depended on Theorem~\ref{c1alpha}, we see that corresponding to 
infinitely many $k$, there are orthogonal rotations $\G_{k}$, cones ${\mathbf C}_{k}^{\prime} \in {\mathcal K}$ such that  
the conclusions of the present lemma hold with $V_{k},$ ${\mathbf C}_{k}$, ${\mathbf C}_{k}^{\prime}$ and $\G_{k}$ in place of $V$, ${\mathbf C}$, ${\mathbf C}^{\prime}$ and $\G$, and with constants  $\overline{\k}$, $\overline{C}_{0}$, $\overline{\nu}$ depending only on $\a$, $\g$ and ${\mathbf C}_{0}$. This contradicts our assumption, establishing the lemma. 
\end{proof}

\begin{lemma}\label{MD-multi-scale-final} 
Let $\a \in (0, 1)$, $q$ be an integer $\geq 2$ and $\g \in (0, 1/2).$ Let ${\mathbf C}_{0} = \sum_{j=1}^{m_{0}} \sum_{\ell = 1}^{q_{j}^{(0)}} |H_{j, \, \ell}|$ be the stationary cone 
as in (\ref{MD-1-0}) with $\Theta \, (\|{\mathbf C}_{0}\|, 0) = q + 1/2.$ 
For $j=1, 2, \ldots, 2q-m_{0} +1,$ let $\th_{j} \in (0, 1/4)$ be such that $\th_{1} > 8\th_{2} > 64\th_{3} > \ldots 
>8^{2q -m_{0}}\th_{2q-m_{0} +1}$. 
There exists a number $\e_{0} = \e_{0}(\a, \g, \th_{1}, \th_{2}, \ldots, \th_{2q-m_{0}+1}, {\mathbf C}_{0}) \in (0, 1/2)$ such that the following 
is true: If $V \in {\mathcal S}_{\a},$ ${\mathbf C} \in {\mathcal K}$ satisfy Hypotheses~\ref{MD-hyp} and if the induction hypotheses $(H1)$, $(H2)$ hold,
then there exist orthogonal rotations $\G, \Delta$ of ${\mathbf R}^{n+1}$ and cones ${\mathbf C}^{\prime}, {\mathbf C}^{\prime\prime} \in {\mathcal K}$ 
such that, with 
$${\mathcal Q}_{V}^{2}({\mathbf C}) = \int_{B_{1}^{n+1}(0) \setminus \{r(X) < 1/16\}} {\rm dist}^{2}(X, {\rm spt} \, \|V\|) \,
d\|{\mathbf C}\|(X) + \int_{B_{1}^{n+1}(0)} {\rm dist}^{2} \, (X, {\rm spt} \, \|{\mathbf C})\|) \, d\|V\|(X),$$ 
$${\mathcal R}_{V}^{2}({\mathbf C}) = \int_{B_{1}^{n+1}(0) \setminus \{r(X) < 1/64\}} {\rm dist}^{2}(X, {\rm spt} \, \|V\|) \,
d\|{\mathbf C}\|(X) + \int_{B_{1}^{n+1}(0)} {\rm dist}^{2} \, (X, {\rm spt} \, \|{\mathbf C})\|) \, d\|V\|(X),$$ 
we have the following:
\begin{eqnarray*}
&&({\rm a})\;\;\;\;|e_{j} - \G(e_{j})| \leq \k{\mathcal Q}_{V}({\mathbf C}) \;\;{\rm and} \;\; |e_{j} - \Delta(e_{j})| \leq \k{\mathcal R}_{V}({\mathbf C})\;\;\; {\rm for} \;\; j=1,2, 3, \ldots, n+1;\nonumber\\  
&&({\rm b})\;\;\;\;{\rm dist}_{\mathcal H}^{2} \, ({\rm spt} \, \|{\mathbf C}^{\prime}\| \cap B_{1}^{n+1}(0), {\rm spt} \, \|{\mathbf C}\| \cap B_{1}^{n+1}(0)) \leq C_{0}{\mathcal Q}_{V}^{2}({\mathbf C})\;\;{\rm and}\nonumber\\
&&\;\;\;\;\;\;\;\;\;{\rm dist}_{\mathcal H}^{2} \, ({\rm spt} \, \|{\mathbf C}^{\prime\prime}\| \cap B_{1}^{n+1}(0), {\rm spt} \, \|{\mathbf C}\| \cap B_{1}^{n+1}(0)) \leq C_{0}{\mathcal R}_{V}^{2}({\mathbf C});\nonumber\\
&&({\rm c})\;\;\;\;\mbox{for some $j^{\prime} \in \{1, 2, \ldots, 2q-m_{0} +1\}$},\nonumber\\ 
&&\hspace{.5in}\th_{j^{\prime}}^{-n-2}\int_{\G\left(B_{\th_{j^{\prime}}}^{n+1}(0) \setminus \{r(X) \leq \th_{j^{\prime}}/16\}\right)} {\rm dist}^{2} \, (X, {\rm spt} \, \|V\|) \, d\|\G_{\#} \, {\mathbf C}^{\prime}\|(X)\nonumber\\
&&\hspace{1.5in} + \; \th_{j^{\prime}}^{-n-2}\int_{B_{\th_{j^{\prime}}}^{n+1}(0)} {\rm dist}^{2} \, (X, {\rm spt} \, \|\G_{\#} \, {\mathbf C}^{\prime}\|) \, d\|V\|(X) \leq \nu_{j^{\prime}}\th_{j^{\prime}}^{2}{\mathcal Q}_{V}^{2}({\mathbf C}),  \;\; {\rm and}\nonumber\\
&&\;\;\;\;\;\;\;\;\;\mbox{for some $j^{\prime\prime} \in \{1, 2, \ldots, 2q-m_{0} +1\}$},\nonumber\\ 
&&\hspace{.5in}\th_{j^{\prime\prime}}^{-n-2}\int_{\D\left(B_{\th_{j^{\prime\prime}}}^{n+1}(0) \setminus \{r(X) \leq \th_{j^{\prime\prime}}/64\}\right)} {\rm dist}^{2} \, (X, {\rm spt} \, \|V\|) \, d\|\D_{\#} \, {\mathbf C}^{\prime\prime}\|(X)\nonumber\\
&&\hspace{1.5in} + \; \th_{j^{\prime\prime}}^{-n-2}\int_{B_{\th_{j^{\prime\prime}}}^{n+1}(0)} {\rm dist}^{2} \, (X, {\rm spt} \, \|\D_{\#} \, {\mathbf C}^{\prime\prime}\|) \, d\|V\|(X) \leq \nu_{j^{\prime\prime}}\th_{j^{\prime\prime}}^{2}{\mathcal R}_{V}^{2}({\mathbf C}).
\end{eqnarray*}
Here $\k$, $C_{0}$ depend only on $\a, \g, {\mathbf C}_{0}$ in case $2q = m_{0}$ and only on $\a, \g, \th_{1}, \ldots, \th_{2q-m_{0}}$ and  ${\mathbf C}_{0}$ in case $2q \geq m_{0}+1;$ 
$\nu_{1} = \nu_{1}(\a, \g, {\mathbf C}_{0})$ and, in case $2q \geq m_{0}+1,$ for each $j=2, 3, \ldots, 2q-m_{0}+1,$ $\nu_{j} = \nu_{j}(\a, \g, \th_{1}, \ldots, \th_{j-1}, {\mathbf C}_{0}).$ In particular, $\nu_{j}$ is independent of $\th_{j}, \th_{j+1}, \ldots, \th_{2q-m_{0}+1}$ for $j=1, 2, \ldots, 2q-m_{0}+1.$ 
\end{lemma}

\begin{proof} First use Lemma~\ref{MD-excess-improvement} and the argument of Lemmas~\ref{multi-scale} and \ref{multi-scale-final} to obtain each of those conclusions above in which ${\mathcal Q}_{V}({\mathbf C})$ appears on the right hand side, with a set of constants $\k_{1}$, $C_{0}^{(1)}$, $\nu_{j}^{(1)}$ in place of $\k$, $C_{0}$, $\nu_{j}$, $j=1, 2, \ldots, 2q-m_{0}+1,$  depending only on the allowed parameters stated in the conclusion. Then repeat the entire argument leading to these conclusions but with 
${\mathcal R}_{V}({\mathbf C})$ in place of ${\mathcal Q}_{V}({\mathbf C})$ (so in particular  part (ii) of Hypothesis~($\dag$) reads ${\mathcal R}_{V}({\mathbf C}) \leq \b \inf_{\widetilde{\mathbf C} \in \cup_{j=m_{0}}^{p-1}{\mathcal K}(j)}{\mathcal R}_{V}({\widetilde{\mathbf C}})$) to obtain those 
conclusions above where  ${\mathcal R}_{V}({\mathbf C})$ appears on the right hand side,  with a set of constants $\k_{2}$, $C_{0}^{(2)}$, $\nu_{j}^{(2)}$ in place of $\k$, $C_{0}$, $\nu_{j}$, $j=1, 2, \ldots, 2q-m_{0}+1,$ depending again only on the allowed parameters. Set 
$\k = \max \, \{\k_{1}, \k_{2}\}$, $C_{0} = \max \, \{C_{0}^{(1)}, C_{0}^{(2)}\}$ and 
$\nu_{j} = \max \, \{\nu_{j}^{(1)}, \nu_{j}^{(2)}\}$ for $j=1, 2, \ldots, 2q-m_{0}+1.$
\end{proof}

\begin{proof}[Proof of Theorem~\ref{no-transverse-q}]
\noindent
{\bf Case 1:} $\Theta \, (\|{\mathbf C}_{0}\|,0) = q + 1/2$, $q \geq 2$. If the theorem is false in this case, then there exist a number $\g \in (0, 1/2)$ and, for each $\ell=1, 2, 3, \ldots,$ a varifold $V_{\ell} \in {\mathcal S}_{\a}$ with $\Theta \, (\|V_{\ell}\|, 0) \geq q+1/2$ and 
$(\omega_{n}2^{n})^{-1}\|V_{\ell}\|(B_{2}^{n+1}(0)) \leq q + 1/2+\g$ such that ${\rm dist}_{\mathcal H} \, ({\rm spt} \, \|V_{\ell}\| \cap B_{1}^{n+1}(0), {\rm spt} \, \|{\mathbf C}_{0}\| \cap B_{1}^{n+1}(0)) \to 0$ as $\ell \to \infty.$ By Allard's integral varifold compactness 
theorem (\cite{AW}; see also \cite{S1}, Section 42.8)  and the constancy theorem (\cite{S1}, Section 41), it follows, after passing to a subsequence without changing notation, 
that $$V_{\ell} \res B_{1}^{n+1}(0) \to \left(\sum_{j=1}^{m_{0}} q_{j}^{(0)}|H_{j}^{(0)}|\right) \res B_{1}^{n+1}(0)$$
as varifolds, where $q_{j}^{(0)}$, $j=1, 2, \ldots, m_{0},$ are positive integers with $\sum_{j=1}^{m_{0}} q_{j}^{(0)} = 2q +1.$ We may assume, by redefining the multiplicities of the original cone ${\mathbf C}_{0}$ if necessary, that ${\mathbf C}_{0} = \sum_{j=1}^{m_{0}} q_{j}^{(0)}|H_{j}^{(0)}|.$ Thus
\begin{equation}\label{MDthm-0}
V_{\ell} \res B_{1}^{n+1}(0) \to {\mathbf C}_{0} \res B_{1}^{n+1}(0) \;\;\; \mbox{as varifolds.}
\end{equation}

For $j=1, 2, \ldots, 2q-m_{0}+1$, choose numbers $\th_{j} = \th_{j}(\a, \g, {\mathbf C}_{0}) \in (0, 1/8)$ as follows:
First choose $\th_{1}$ such that $\nu_{1}\th_{1}^{2(1-\a)} < 1/4,$ where $\nu_{1} = \nu_{1}(\a, \g, {\mathbf C}_{0})$ is as in 
Lemma~\ref{MD-multi-scale-final}. Having chosen $\th_{1}, \th_{2}, \ldots, \th_{j},$ $1 \leq j \leq 2q - m_{0}$, choose 
$\th_{j+1}$ such that $\th_{j+1} < 8^{-1}\th_{j}$ and $\nu_{j+1}\th_{j+1}^{2(1-\a)} < 1/4,$ where $\nu_{j+1} = \nu_{j+1}(\a, \g, \th_{1}, \ldots, \th_{j}, {\mathbf C}_{0})$ is as in Lemma~\ref{MD-multi-scale-final}.

Note that it is easily seen by arguing by contradiction that corresponding to any given $\e^{\prime} \in (0, 1/2)$,  there exist $\e = \e(\e^{\prime}, \a, \g, {\mathbf C}_{0}) \in (0, 1/2)$ such that if Hypotheses~\ref{MD-hyp} are satisfied, then 
$${\mathcal Q}_{V}({\mathbf C}) \leq {\mathcal R}_{V}({\mathbf C}) < \e^{\prime}$$
where ${\mathcal Q}_{V}({\mathbf C})$, ${\mathcal R}_{V}({\mathbf C})$ are defined as in Lemma~\ref{MD-multi-scale-final}.
By Remark (4) following the statement of Hypotheses~\ref{MD-hyp}, it then follows that 
if Hypotheses~\ref{MD-hyp} are satisfied with sufficiently small 
$\e = \e(\e^{\prime}, \a, \g, {\mathbf C}_{0}),$ then for each $Z \in {\rm spt} \, \|V\| \cap B_{1/8}^{n+1}(0)$ with $\Theta \, (\|V\|, Z) \geq q + 1/2$, 
\begin{equation}\label{MDthm-0-0-0}
{\mathcal Q}_{V^{Z}}({\mathbf C}) \leq {\mathcal R}_{V^{Z}}({\mathbf C}) < \e^{\prime}
\end{equation}
where $V^{Z} = \eta_{Z, 1/2 \, \#} \, V.$

Now fix $\ell$ sufficiently large, let $V = V_{\ell}$ and let $Z \in {\rm spt} \, \|V\| \cap B_{1/8}^{n+1}(0)$ with $\Theta \, (\|V\|, Z) \geq q+1/2$. We claim that we may apply Lemma~\ref{MD-multi-scale-final} iteratively to obtain, for each $k=0,1, 2, 3, \ldots$,  an orthogonal rotation $\G_{k}^{Z}$ of ${\mathbf R}^{n+1}$ with $\G_{0}^{Z} =$ Identity, and a cone ${\mathbf C}_{k}^{Z} \in {\mathcal K}$  with ${\mathbf C}_{0}^{Z} = {\mathbf C}_{0}$ satisfying, for $k \geq 1$, 
\begin{equation}\label{MDthm-0-0-1}
|\G_{k}^{Z}(e_{j}) - \G_{k-1}^{Z}(e_{j})| \leq \k\left(\s_{k}^{Z}\right)^{\a}{\mathcal Q}_{V^{Z}}({\mathbf C}_{0}), \;\;\;\;\; j=1,2, \ldots, n+1;
\end{equation}
\begin{equation}\label{MDthm-0-0-2}
{\rm dist}_{\mathcal H}\, ({\rm spt} \, \|{\mathbf C}_{k}^{Z}\| \cap B_{1}^{n+1}(0), {\rm spt} \, \|{\mathbf C}_{k-1}^{Z}\| \cap B_{1}^{n+1}(0)) \leq  C_{0}\left(\s_{k}^{Z}\right)^{\a}{\mathcal Q}_{V^{Z}}({\mathbf C}_{0});
\end{equation}
\begin{equation}\label{MDthm-0-0-3}
\left(\s_{k}^{Z}\right)^{-n-2}\int_{B_{\s_{k}^{Z}}^{n+1}(0)} {\rm dist}^{2} \, (X,  {\rm spt} \,\|\left(\G_{k}^{Z}\right)_{\#} \, {\mathbf C}_{k}^{Z}\|) \, d\|V^{Z}\|(X) \leq \left(\s_{k}^{Z}\right)^{2\a} {\mathcal Q}_{V^{Z}}^{2}({\mathbf C}_{0}), \;\;\; {\rm and}
\end{equation}
\begin{equation}\label{MDthm-0-0-4}
\left(\s_{k}^{Z}\right)^{-n-2}\int_{\G^{Z}_{k}(B^{n+1}_{\s_{k}^{Z}}(0) \setminus \{|r(X)| \leq \s_{k}^{Z}/16\})} {\rm dist}^{2} \, (X, {\rm spt} \, \|V^{Z}\|) \, d\|\left(\G_{k}^{Z}\right)_{\#} \, {\mathbf C}_{k}^{Z}\|(X)\leq \left(\s_{k}^{Z}\right)^{2\a} {\mathcal Q}_{V^{Z}}^{2}({\mathbf C}_{0})
\end{equation}
where $\k = \k(\a, \g, {\mathbf C}_{0})$, $C_{0} = C_{0}(\a, \g, {\mathbf C}_{0})$ are as in Lemma~\ref{MD-multi-scale-final} and $\{\s_{k}^{Z}\}$ is a sequence of positive numbers such that $\s_{0}^{Z}  = 1$ and for each $k=1, 2, \ldots,$ 
$\s_{k}^{Z} = \th_{j_{k}^{Z}}\s_{k-1}^{Z}$ for some $j_{k}^{Z} \in \{1, 2, \ldots, 2q-m_{0}+1\}.$
To see this, note first that it follows from Remark (4) following the statement of Hypotheses~\ref{MD-hyp} that if $V = V_{\ell}$ 
with $\ell$ fixed sufficiently large,  then for each $Z \in {\rm spt} \, \|V\| \cap B_{1/8}^{n+1}(0)$ 
with $\Theta \, (\|V\|, Z) \geq q + 1/2$, Hypotheses~\ref{MD-hyp} are satisfied with $V^{Z}$ in place of $V,$ ${\mathbf C}_{0}$ in place of ${\mathbf C}$ and with $\e = \e_{0}(\a, \g, {\mathbf C}_{0})$ where $\e_{0}$ is as in Lemma~\ref{MD-multi-scale-final}; hence by applying Lemma~\ref{MD-multi-scale-final} with $V^{Z}$ in place of $V$ and ${\mathbf C} = {\mathbf C}_{0},$ we deduce that (\ref{MDthm-0-0-1})-(\ref{MDthm-0-0-4}) hold in case $k=1.$ So let $k \geq 2$ and suppose by induction that (\ref{MDthm-0-0-1})-(\ref{MDthm-0-0-4}) are valid with $1, 2, \ldots, k-1$ in place of $k$. Then for any given $\e \in (0, 1/4)$, provided $V = V_{\ell}$ with $\ell$ sufficiently large, Hypotheses~\ref{MD-hyp} are satisfied with $\left(\G_{k-1}^{Z}\right)^{-1}_{\#} \, \eta_{0, \s_{k-1}^{Z} \, \#} \, V^{Z}$ in place of $V$ and with 
${\mathbf C} = {\mathbf C}_{k-1}^{Z};$ here, the validity of Hypotheses~\ref{MD-hyp}(1)-(4) with 
$\left(\G_{k-1}^{Z}\right)^{-1}_{\#} \, \eta_{0, \s_{k-1}^{Z} \, \#} \, V^{Z}$ in place of $V$ and ${\mathbf C}_{k-1}^{Z}$ in place of ${\mathbf C}$ is clear, and in verifying Hypothesis~\ref{MD-hyp}(5) with $\left(\G_{k-1}^{Z}\right)^{-1}_{\#} \, \eta_{0, \s_{k-1}^{Z} \, \#} \, V^{Z}$ in place of $V$, note first that by Remarks (1) and  (4) following the statement of Hypotheses~\ref{MD-hyp} (taken with $\r = \s_{1}^{Z}$ and $\t = \frac{1}{32}\min \{\th_{1}, \th_{2} \, \ldots, \th_{2q-m_{0}+1}\}  = \frac{1}{32}\th_{2q-m_{0}+1}$) we have that 
\begin{equation}\label{MDthm-0-0-5}
\eta_{0, \s_{1}^{Z} \, \#} \, V^{Z} \res \left(B_{1}^{n+1}(0) \setminus \left\{r(X) < \frac{1}{32}\th_{2q-m_{0}+1}\right\}\right) = \sum_{j=1}^{m_{0}} \sum_{i=1}^{q_{j}^{(0)}} |{\rm graph} \, \widetilde{u}_{j, \, i}|,
\end{equation}
where for each $j \in \{1, 2, \ldots, m_{0}\}$ and $i \in \{1, 2, \ldots, q_{j}^{(0)}\}$, 
$$\widetilde{u}_{j, \, i}  \in C^{2}\,\left(H_{j}^{(0)} \cap \left(B_{1}^{n+1}(0) \setminus \left\{r(X) < \frac{1}{32}\th_{2q-m_{0}+1}\right\}\right); \left(H_{j}^{(0)}\right)^{\perp}\right)$$ 
and $\widetilde{u}_{j, \, i}$ are solutions to the minimal surface equation over $H_{j}^{(0)} \cap \left(B_{1}^{n+1}(0) \setminus \{r(X) < \frac{1}{32}\th_{2q-m_{0}+1}\}\right)$ with small $C^{2}$ norm; so in particular, in view of (\ref{MDthm-0-0-1}),  Hypothesis~\ref{MD-hyp}(5) is satisfied with $\left(\G_{1}^{Z}\right)^{-1}_{\#} \, \eta_{0, \s_{1}^{Z} \, \#} \, V^{Z}$ in place of $V$. On the other hand, by (\ref{MDthm-0-0-2}), (\ref{MDthm-0-0-3}) and (\ref{MDthm-0-0-4}), we may apply Remarks (2) and (3) following the statement of Hypotheses~\ref{MD-hyp} with $\left(\G_{r}^{Z}\right)^{-1}_{\#} \, \eta_{0, \s_{r}^{Z} \, \#} \, V^{Z}$ in place of $V$ and $\t = \frac{1}{2}\th_{2q-m_{0}+1}$, followed by 
Theorem~\ref{SS}, to deduce that for each $r \in \{2, 3, \ldots,k-1\}$, 
\begin{equation}\label{MDthm-0-0-6}
\left(\G_{r}^{Z}\right)^{-1}_{\#}\eta_{0, \s_{r}^{Z} \, \#} \, V^{Z} \res \left(B_{1/2}^{n+1}(0) \setminus \left\{r(X) < \frac{1}{2}\th_{2q-m_{0}+1}\right\}\right) = \sum_{j=1}^{m_{0}} \sum_{i=1}^{p_{j}^{Z,r}} |{\rm graph} \, \widetilde{u}_{j, \, i}^{Z,r}|,
\end{equation}
for some integers $p_{j}^{Z,r} \geq 1$, where for each $j \in \{1, 2, \ldots, m_{0}\}$ and $i \in \{1, 2, \ldots, p_{j}^{Z,r}\}$, 
$$\widetilde{u}_{j, \, i}^{Z,r}  \in C^{2}\,\left(H_{j}^{(0)} \cap \left(B_{1/2}^{n+1}(0) \setminus \left\{r(X) < \frac{1}{2}\th_{2q-m_{0}+1}\right\}\right); \left(H_{j}^{(0)}\right)^{\perp}\right)$$ 
and $\widetilde{u}_{j, \, i}^{Z,r}$ are solutions to the minimal surface equation over $H_{j}^{(0)} \cap \left(B_{1/2}^{n+1}(0) \setminus \{r(X) < \frac{1}{2}\th_{2q-m_{0}+1}\}\right)$ with small $C^{2}$ norm.
Since $\s_{r}^{Z} \geq \th_{2q-m_{0}+1} \s_{r-1}^{Z}$ for each $r \geq 1$, it follows from (\ref{MDthm-0-0-5}), (\ref{MDthm-0-0-6}) and unique continuation of solutions to the minimal surface equation that 
\begin{equation}\label{MDthm-9-0}
p_{j}^{Z,r} = q_{j}^{(0)}
\end{equation}
for each $r \in \{2, 3, \ldots, k-1\}$ and $j  \in \{1, 2,\ldots, m_{0}\},$ whence, by (\ref{MDthm-0-0-6}), we see that 
Hypothesis~\ref{MD-hyp}(5) with $\left(\G_{k-1}^{Z}\right)^{-1}_{\#} \, \eta_{0, \s_{k-1}^{Z} \, \#} \, V^{Z}$ in place of $V$ is satisfied as  claimed. Hence we may apply Lemma~\ref{MD-multi-scale-final} with 
$\left(\G_{k-1}^{Z}\right)^{-1}_{\#} \, \eta_{0, \s_{k-1}^{Z} \, \#} \, V^{Z}$ in place of $V$ and ${\mathbf C}_{k-1}^{Z}$ in place of ${\mathbf C}$ to obtain an orthogonal rotation $\G_{k}^{Z}$ of ${\mathbf R}^{n+1}$ and a cone ${\mathbf C}_{k}^{Z} \in {\mathcal K}$ satisfying (\ref{MDthm-0-0-1})-(\ref{MDthm-0-0-4}). This establishes inductively the validity of (\ref{MDthm-0-0-1})-(\ref{MDthm-0-0-4}) for each $k=1, 2, 3, \ldots$. Using (\ref{MDthm-0-0-1})-(\ref{MDthm-0-0-4}) in a standard way, we reach the conclusion that if $V = V_{\ell}$ with $\ell$ fixed sufficiently large, then corresponding to each $Z \in {\rm spt} \, \|V\| \cap B_{1/8}^{n+1}(0)$ with $\Theta \, (\|V\|, Z) \geq q+1/2$, there exist a cone ${\mathbf C}^{Z} \in {\mathcal K}$ with 
\begin{equation}\label{MDthm-4}
{\rm dist}_{\mathcal H} \, ({\rm spt} \, \|{\mathbf C}^{Z}\| \cap B_{1}^{n+1}(0), {\rm spt} \, \|{\mathbf C}_{0}\| \cap B_{1}^{n+1}(0)) \leq C{\mathcal Q}_{V^{Z}}({\mathbf C}_{0});
\end{equation}
an orthogonal rotation $\G^{Z}$ of ${\mathbf R}^{n+1}$ satisfying, for each $k=0, 1, 2, \ldots$, 
\begin{equation}\label{MDthm-5}
|\G^{Z}(e_{j}) - \G_{k}^{Z}(e_{j})| \leq C\left(\s_{k}^{Z}\right)^{\a}{\mathcal Q}_{V^{Z}}({\mathbf C}_{0}), \;\;\; j=1, 2, \ldots, n+1; 
\end{equation}
such that
\begin{equation}\label{MDthm-6}
\left(\s^{Z}_{k}\right)^{-n-2}\int_{\G^{Z}_{k}(B^{n+1}_{\s^{Z}_{k}/2}(0) \setminus \{|r(X)| \leq \s^{Z}_{k}/16\})} {\rm dist}^{2} \, (X, {\rm spt} \, \|V^{Z}\|) \, d\|\G^{Z}_{\#} \, {\mathbf C}^{Z}\|(X)\leq C \left(\s_{k}^{Z}\right)^{2\a} {\mathcal Q}_{V^{Z}}^{2}({\mathbf C}_{0}) 
\end{equation}
for each $k=0, 1, 2, \ldots$ and 
\begin{equation}\label{MDthm-7}
\r^{-n-2}\int_{B_{\r}^{n+1}(0)} {\rm dist}^{2} \, (X, {\rm spt} \,\|\G^{Z}_{\#} \, {\mathbf C}^{Z}\|) \, d\|V^{Z}\|(X) \leq C\r^{2\a} {\mathcal Q}_{V^{Z}}^{2}({\mathbf C}_{0}),
\end{equation}
for all $\r \in (0, 1/4]$, where $C = C(\a, \g, {\mathbf C}_{0}) \in (0, \infty).$

Let $T_{V} = \{ Z \in {\rm spt} \, \|V\|\, : \, \Theta \,(\|V\|, Z) \geq q+1/2\} \cap B_{1}^{n+1}(0).$ We now use the estimates (\ref{MDthm-4})-(\ref{MDthm-7}), Lemmas~\ref{MD-no-gaps}(a), \ref{MD-multi-scale-final} and Corollary~\ref{MD-height-est} to establish that  $T_{V} \cap B_{1/32}^{n+1}(0)$ is an $(n-1)$-dimensional embedded $C^{1, \a}$ submanifold of $B_{1/32}^{n+1}(0)$ containing the origin. Indeed, note first that the estimates (\ref{MDthm-4}), (\ref{MDthm-6}) and (\ref{MDthm-7}) imply that for any given $\e \in (0, 1/4)$, if $V = V_{\ell}$  with fixed $\ell$ sufficiently large, then for each $Z \in T_{V} \cap B_{1/16}^{n+1}(0)$ and each $k \geq 1$, Hypotheses~\ref{MD-hyp} are satisfied with $\left(\G^{Z}\right)^{-1}_{\#} \, \eta_{Z, \frac{1}{2}\s_{k}^{Z} \, \#} \, V$ in place of $V$ and 
${\mathbf C} = {\mathbf C}^{Z}.$ (In verifying Hypothesis~\ref{MD-hyp}(5) with $\left(\G^{Z}\right)^{-1}_{\#} \, \eta_{Z, \frac{1}{2}\s_{k}^{Z} \, \#} \, V = \left(\G^{Z}\right)^{-1}_{\#} \, \eta_{0, \s_{k}^{Z} \, \#} \, V^{Z}$ in place of $V$, we argue exactly as we did in verifying Hypothesis~\ref{MD-hyp}(5) with $\left(\G_{k-1}^{Z}\right)^{-1}_{\#} \, \eta_{0, \s_{k-1}^{Z} \, \#} \, V^{Z}$ in place of $V$ as part of the inductive step described above.)

Consequently, we see that  for each point $(0, y) \in \{0\} \times {\mathbf R}^{n-1} \cap B_{1/16}^{n+1}(0),$
\begin{equation}\label{MDthm-9-1}
T_{V} \cap {\mathbf R}^{2} \times \{(0, y)\} \neq \emptyset;
\end{equation}
for if there is a point $(0, y) \in \{0\} \times {\mathbf R}^{n-1} \cap B_{1/16}^{n+1}(0)$ with $T_{V} \cap \left({\mathbf R}^{2} \times \{(0, y)\}\right) = \emptyset$, then, since $T_{V} \cap B_{1/16}^{n+1}(0)$ is a relatively closed 
subset of $B_{1/16}^{n+1}(0)$ and $0 \in T_{V}$, we can find $r \in (0, 1/16)$ such that $T_{V} \cap \left({\mathbf R}^{2} \times B_{r}^{n-1}(0,y)\right) = \emptyset$ but 
$T_{V} \cap \left({\mathbf R}^{2} \times \partial \, B_{r}^{n-1}(0,y)\right) \neq \emptyset,$ whence we may, in view of (\ref{MDthm-5}), (\ref{MDthm-6})  and (\ref{MDthm-7}), pick any point $Z \in T_{V} \cap \left({\mathbf R}^{2} \times \partial \, B_{r}^{n-1}(0,y)\right),$ choose $k$ such that $\s_{k}^{Z} < r/4$ and apply
Lemma~\ref{MD-no-gaps}(a) with $\left(\G^{Z}\right)^{-1}_{\#} \, \eta_{Z, \frac{1}{2}\s_{k}^{Z} \, \#} \, V$ in place of $V,$ ${\mathbf C} = {\mathbf C}^{Z}$ and $\d = 1/8$ to get a contradiction with the assumption $T_{V} \cap \left({\mathbf R}^{2} \times B_{r}^{n-1}(0,y)\right) = \emptyset$. 
   
For $Z \in T_{V}$, let $S_{Z} = Z + \G^{Z}\left(\{0\} \times {\mathbf R}^{n-1}\right)$ and note that for each $Z \in T_{V}$ and each $\r \in (0, 1/4]$, 
\begin{equation}\label{MDthm-10}
T_{V} \cap \left(B_{\r}^{n+1}(Z) \setminus \left\{X \in {\mathbf R}^{n+1} \, : \, {\rm dist} \, (X, S_{Z}) <  \frac{1}{8}\r\right\}\right) = \emptyset;
\end{equation}
this is easily seen by choosing, for given $Z \in T_{V}$ and $\r \in (0, 1/4],$ the unique integer $k$ such that $\frac{15}{32}\s_{k+1}^{Z} < \r \leq \frac{15}{32}\s_{k}^{Z}$, and applying Remark (2) following Hypotheses~\ref{MD-hyp} with $\t = \frac{1}{16}\th_{2q-m_{0}+1}$ and with $\left(\G^{Z}\right)^{-1}_{\#} \, \eta_{Z, \frac{1}{2}\s_{k}^{Z} \, \#} \, V$ in place of $V.$ This and (\ref{MDthm-9-1}) imply that for each $(0, y) \in \{0\} \times {\mathbf R}^{n-1} \cap B_{1/16}^{n+1}(0)$, the set $T_{V} \cap {\mathbf R}^{2} \times \{(0, y)\}$ consists of a unique point, so that 
\begin{equation}\label{MDthm-graph}
T_{V} \cap B_{1/16}^{n+1}(0) = {\rm graph} \, \varphi
\end{equation}
for a function $\varphi  = (\varphi_{1}, \varphi_{2})\, : \, B_{1/16}^{n-1}(0) \to {\mathbf R}^{2}.$ Moreover, (\ref{MDthm-10}) and the estimates  (\ref{MDthm-0-0-0}), (\ref{MDthm-5}) say that $\varphi$ is Lipschitz with ${\rm Lip} \, (\varphi) \leq 1$ and, writing $\widetilde{\varphi}(Z) = (\varphi_{1}(Z), \varphi_{2}(Z), Z)$ for $Z \in B_{1/16}^{n-1}(0)$,  that 
\begin{equation}\label{MDthm-slope}
D\widetilde{\varphi}(Z)\left(\{0\}\times {\mathbf R}^{n-1}\right)  = \G^{Z}(\{0\} \times {\mathbf R}^{n-1})
\end{equation}
for ${\mathcal H}^{n-1}$-a.e. $Z \in B_{1/16}^{n-1}(0).$

We now argue that $\left.\varphi\right|_{B_{1/32}^{n-1}(0)} $ must be of class $C^{1, \a}$. For this, first observe that by employing exactly the argument leading to (\ref{MDthm-4})-(\ref{MDthm-7}) but 
using those conclusions of Lemma~\ref{MD-multi-scale-final} involving ${\mathcal R}_{V}({\mathbf C})$ (in place of those involving ${\mathcal Q}_{V}({\mathbf C})$), we obtain for each $Z \in T_{V}$ orthogonal rotations $\D^{Z}, \D_{k}^{Z}$ of ${\mathbf R}^{n+1}$ for $k=1, 2, 3, \ldots;$ a cone ${\mathbf W}^{Z} \in {\mathcal K}$ and numbers $\t_{k}^{Z} \in (0, 1]$ for $k=1, 2, 3, \ldots,$ where for each $k$, $\t_{k}^{Z} = \th_{{\ell}_{k}^{Z}}\t_{k-1}^{Z}$ for some $\ell_{k}^{Z} \in \{1, 2, \ldots, 2q-m_{0}+1\},$ such that   
\begin{equation}\label{MDthm-4-alt}
{\rm dist}_{\mathcal H} \, ({\rm spt} \, \|{\mathbf W}^{Z}\| \cap B_{1}^{n+1}(0), {\rm spt} \, \|{\mathbf C}_{0}\| \cap B_{1}^{n+1}(0)) \leq C{\mathcal R}_{V^{Z}}({\mathbf C}_{0});
\end{equation}
\begin{equation}\label{MDthm-5-alt}
|\D^{Z}(e_{j}) - \D_{k}^{Z}(e_{j})| \leq C\left(\t_{k}^{Z}\right)^{\a}{\mathcal R}_{V^{Z}}({\mathbf C}_{0}), \;\;\; j=1, 2, \ldots, n+1; 
\end{equation}
\begin{equation}\label{MDthm-6-alt}
\left(\t^{Z}_{k}\right)^{-n-2}\int_{\D^{Z}_{k}(B^{n+1}_{\t^{Z}_{k}/2}(0) \setminus \{|r(X)| \leq \t^{Z}_{k}/64\})} {\rm dist}^{2} \, (X, {\rm spt} \, \|V^{Z}\|) \, d\|\D^{Z}_{\#} \, {\mathbf W}^{Z}\|(X)\leq C \left(\t_{k}^{Z}\right)^{2\a} {\mathcal R}_{V^{Z}}^{2}({\mathbf C}_{0}) 
\end{equation}
for each $k=0, 1, 2, \ldots$ and 
\begin{equation}\label{MDthm-7-alt}
\r^{-n-2}\int_{B_{\r}^{n+1}(0)} {\rm dist}^{2} \, (X, {\rm spt} \,\|\D^{Z}_{\#} \, {\mathbf W}^{Z}\|) \, d\|V^{Z}\|(X) \leq C\r^{2\a} {\mathcal R}_{V^{Z}}^{2}({\mathbf C}_{0}),
\end{equation}
for all $\r \in (0, 1/4]$, where $C = C(\a, \g, {\mathbf C}_{0}) \in (0, \infty).$ 

Since the sequence of varifolds $W_{k} = \eta_{Z, \s_{k}^{Z} \, \#} \, V,$ $k=1, 2, 3, \ldots$ has a subsequence $W_{k^{\prime}}$ which converges to a cone ${\mathbf P}$ satisfying, by (\ref{MDthm-6}) and
(\ref{MDthm-7}), ${\rm spt} \, \|{\mathbf P}\| = {\rm spt}\, \|\G^{Z}_{\#} \,{\mathbf C}^{Z}\|,$ it follows from 
(\ref{MDthm-7-alt}) taken with $\r = \s_{k^{\prime}}^{Z}$  that   ${\rm spt}\, \|\G^{Z}_{\#} \,{\mathbf C}^{Z}\| \subseteq {\rm spt} \, \|\D^{Z}_{\#} \, {\mathbf W}^{Z}\|.$ Same reasoning applied to the sequence $\eta_{Z, \t_{k}^{Z} \, \#} \, V$ establishes the reverse inclusion, so we have that 
${\rm spt}\, \|\G^{Z}_{\#} \,{\mathbf C}^{Z}\| = {\rm spt} \, \|\D^{Z}_{\#} \, {\mathbf W}^{Z}\|$ 
whence in particular that
\begin{equation}\label{MDthm-tilt}
\G^{Z}(\{0\} \times {\mathbf R}^{n-1})  = \D^{Z}(\{0\} \times {\mathbf R}^{n-1}).
\end{equation}  

Recall (cf. paragraph preceding (\ref{MDthm-9-1})) that given any $\e \in (0, 1/2)$, if $V = V_{\ell}$ with $\ell$ fixed sufficiently large depending on $\e$, then for each $Z \in T_{V} \cap B_{1/16}^{n+1}(0)$ and $k \geq 1$, Hypotheses~\ref{MD-hyp} are satisfied with $V_{k, Z} \equiv \left(\D^{Z}\right)^{-1}_{\#} \, \eta_{Z, \frac{1}{2}\t_{k}^{Z} \, \#} \, V$ in place of $V$ and ${\mathbf W}^{Z}$ in place of ${\mathbf C}.$ Consequently, by Remark (4) following the statement of Hypotheses~\ref{MD-hyp}, we see that given any $\e \in (0, 1/2)$,  if $V = V_{\ell}$ for $\ell$ fixed sufficiently large, then for any $Z \in T_{V}$, $k \geq 1$ and $\widetilde{Z} \in T_{V_{k, Z}},$ Hypotheses~\ref{MD-hyp} are satisfied with $\eta_{\widetilde{Z}, 1/2 \, \#} \, V_{k, Z}$  in place of $V$ and ${\mathbf W}^{Z}$ in place of ${\mathbf C}.$ Now take any two distinct points $Z_{1}, Z_{2} \in T_{V} \cap B_{1/32}^{n+1}(0)$, let $m$ be the unique integer satisfying $\t_{m+1}^{Z_{1}} < 2|Z_{1} - Z_{2}| \leq \t_{m}^{Z_{1}}$ and let 
$\widetilde{V} = V_{m, Z_{1}} = \left(\D^{Z_{1}}\right)^{-1}_{\#} \, \eta_{Z_{1}, \frac{1}{2}\t_{m}^{Z_{1}} \, \#} \, V.$ Letting $\widetilde{Z} = \left(\D^{Z_{1}}\right)^{-1} \left(\frac{2(Z_{2} - Z_{1})}{\t_{m}^{Z_{1}}}\right)$ and noting that $\widetilde{Z} \in {\rm spt} \, \|\widetilde{V}\| \cap B_{1/16}^{n+1}(0)$ with $\Theta \, (\|\widetilde{V}\|, \widetilde{Z}) \geq q + 1/2,$ we may apply Lemma~\ref{MD-multi-scale-final} iteratively (utilising its conclusions involving ${\mathcal Q}_{(\cdot)}(\cdot)$), starting with $\widetilde{V}^{\widetilde{Z}} = \eta_{\widetilde{Z}, 1/2 \, \#} \, \widetilde{V}$ in place of $V$ and ${\mathbf W}^{Z_{1}}$ in place of ${\mathbf C}$ (and with $\th_{1}, \th_{2}, \ldots \th_{2q-m_{0}+1}$ equal to the same fixed constants as chosen at the beginning of the proof of the present theorem), in the manner exactly as in the argument leading to (\ref{MDthm-4})-(\ref{MDthm-7}), to conclude that there exist a cone $\widetilde{\mathbf C} \in {\mathcal K}$ with 
\begin{equation}\label{MDthm-11}
{\rm dist}_{\mathcal H} \, ({\rm spt} \, \|\widetilde{\mathbf C}\| \cap B_{1}^{n+1}(0), {\rm spt} \, \|{\mathbf W}^{Z_{1}}\| \cap B_{1}^{n+1}(0)) \leq C{\mathcal Q}_{\widetilde{V}^{\widetilde{Z}}}({\mathbf W}^{Z_{1}});
\end{equation}
orthogonal rotations $\widetilde{\G}, \widetilde{\G}_{0}, \widetilde{\G}_{1}, \ldots$ of ${\mathbf R}^{n+1}$ with $\widetilde{\G}_{0}=$ Identity; and a sequence of positive numbers $\{\widetilde{\s}_{k}\}$ with $\widetilde{\s}_{0}  = 1$ and  
$\widetilde{\s}_{k} = \th_{\widetilde{j}_{k}}\widetilde{\s}_{k-1}$ for some $\widetilde{j}_{k}\in \{1, 2, \ldots, 2q-m_{0}+1\}$ and each $k \geq 1,$ satisfying, for each $k=0, 1, 2, \ldots$, 
\begin{equation}\label{MDthm-12}
|\widetilde{\G}(e_{j}) - \widetilde{\G}_{k}(e_{j})| \leq C\left(\widetilde{\s}_{k}\right)^{\a}{\mathcal Q}_{\widetilde{V}^{\widetilde{Z}}}({\mathbf W}^{Z_{1}}), \;\;\; j=1, 2, \ldots, n+1; 
\end{equation}
\begin{equation}\label{MDthm-13}
\left(\widetilde{\s}_{k}\right)^{-n-2}\int_{\widetilde{\G}_{k}(B^{n+1}_{\widetilde{\s}_{k}/2}(0) \setminus \{|r(X)| \leq \widetilde{\s}_{k}/16\})} {\rm dist}^{2} \, (X, {\rm spt} \, \|\widetilde{V}^{\widetilde{Z}}\|) \, d\|\widetilde{\G}_{\#} \, \widetilde{\mathbf C}\|(X)\leq C \left(\widetilde{\s}_{k}\right)^{2\a} {\mathcal Q}_{\widetilde{V}^{\widetilde{Z}}}^{2}({\mathbf W}^{Z_{1}}); 
\end{equation}
and for each $\r \in (0, 1/4],$
\begin{equation}\label{MDthm-14}
\r^{-n-2}\int_{B_{\r}^{n+1}(0)} {\rm dist}^{2} \, (X, {\rm spt} \,\|\widetilde{\G}_{\#} \, \widetilde{\mathbf C}\|) \, d\|\widetilde{V}^{\widetilde{Z}}\|(X) \leq C\r^{2\a} {\mathcal Q}_{\widetilde{V}^{\widetilde{Z}}}^{2}({\mathbf W}^{Z_{1}}),
\end{equation}
where $C = C(\a, \g, {\mathbf C}_{0}) \in (0, \infty)$ is as in (\ref{MDthm-4})-(\ref{MDthm-7}). Noting that $\widetilde{V}^{\widetilde{Z}} = \left(\D^{Z_{1}}\right)^{-1}_{\#} \, \eta_{Z_{2}, \frac{1}{2}\t_{m}^{Z_{1}} \, \#} \, V,$ we deduce from (\ref{MDthm-13}), (\ref{MDthm-14}) and the inequalities (\ref{MDthm-6}), (\ref{MDthm-7}) taken with $Z = Z_{2},$ and reasoning exactly as for (\ref{MDthm-tilt}), that 
\begin{equation}\label{MDthm-tilt-1}
\D^{Z_{1}}\circ\widetilde{\G}(\{0\} \times {\mathbf R}^{n-1})= \G^{Z_{2}} (\{0\} \times {\mathbf R}^{n-1}).
\end{equation}
This together with (\ref{MDthm-tilt}) taken with $Z = Z_{1}$ and (\ref{MDthm-12}) taken with $k=0$ implies that 
$${\rm dist}_{\mathcal H} \, \left(\G^{Z_{1}}\left(\{0\} \times {\mathbf R}^{n-1}\right) \cap B_{1}^{n+1}(0),
\G^{Z_{2}}\left(\{0\} \times {\mathbf R}^{n-1}\right) \cap B_{1}^{n+1}(0)\right) \leq C{\mathcal Q}_{\widetilde{V}^{\widetilde{Z}}}({\mathbf W}^{Z_{1}})$$
where $C = C(\a, \g, {\mathbf C}_{0}) \in (0, \infty).$ On the other hand, we see directly from Corollary~\ref{MD-height-est}(c) ( taken with $\widetilde{V}$ in place of $V$ and 
${\mathbf W}^{Z_{1}}$ in place of ${\mathbf C}$) that 
${\mathcal Q}_{\widetilde{V}^{\widetilde{Z}}}({\mathbf W}^{Z_{1}}) \leq C{\mathcal R}_{\widetilde{V}}({\mathbf W}^{Z_{1}}),$ and by (\ref{MDthm-6-alt}), (\ref{MDthm-7-alt}) and (\ref{MDthm-0-0-0}) that
${\mathcal R}_{\widetilde{V}}({\mathbf W}^{Z_{1}}) \leq C(\t_{m}^{Z_{1}})^{\a}{\mathcal R}_{V^{Z_{1}}}({\mathbf C}_{0}) \leq C|Z_{1} - Z_{2}|^{\a}$ where $C = C(n, \a, \g, {\mathbf C}_{0}) \in (0, \infty).$ We have thus established that for any pair of points $Z_{1}, Z_{2} \in T_{V} \cap B_{1/32}^{n+1}(0)$, 
 $${\rm dist}_{\mathcal H} \, \left(\G^{Z_{1}}\left(\{0\} \times {\mathbf R}^{n-1}\right) \cap B_{1}^{n+1}(0),
\G^{Z_{2}}\left(\{0\} \times {\mathbf R}^{n-1}\right) \cap B_{1}^{n+1}(0)\right) \leq C|Z_{1} - Z_{2}|^{\a},$$
which in view of (\ref{MDthm-slope}) says that 
\begin{equation}\label{MDthm-reg}
\left.\varphi\right|_{B_{1/32}^{n-1}(0)} \in C^{1, \a}(B_{1/32}^{n-1}(0)).
\end{equation}
Now fix $j \in \{1, 2, \ldots, m_{0}\}$ and assume, for notational convenience and without loss of generality, that $H_{j}^{(0)} = \{(0, x^{2}, y)  \in {\mathbf R}^{n+1}\, : \, x^{2} > 0, \, y \in {\mathbf R}^{n-1}\}.$ Let 
$T^{\prime}_{V}$ be the orthogonal projection of $T_{V} \cap B_{1/32}^{n+1}(0)$ onto 
the hyperplane $\{x^{1}= 0\} \equiv {\mathbf R}^{n},$ so that 
$T^{\prime}_{V} = \{(0, \varphi_{2}(y), y) \, : \, y \in B_{1/32}^{n-1}(0)\}.$ Assuming that $V = V_{\ell}$ with $\ell$ sufficiently large, note then that  $T^{\prime}_{V} \subset  \{|x^{2}| < 1/128\}$ and by (\ref{MDthm-graph}) and (\ref{MDthm-reg}) that $B_{1/64}^{n}(0) \setminus T_{V}^{\prime}$ has exactly two components. Let $\Omega^{\prime}$ be the component of 
$B_{1/64}^{n}(0) \setminus T_{V}^{\prime}$ containing $B_{1/64}^{n}(0) \cap \{x^{2} > 1/128\}.$  
Keeping in mind that (\ref{MDthm-0-0-6}) and (\ref{MDthm-9-0}) are valid for each $Z \in T_{V} \cap B_{1/32}^{n+1}(0)$ and each $r=1,2, 3, \ldots$, it follows from (\ref{MDthm-0-0-6}), (\ref{MDthm-9-0}) and unique continuation of solutions to the minimal surface equation that 
$$V \res \left(({\mathbf R} \times \Omega^{\prime}) \cap N_{j}\right) = \sum_{i=1}^{q_{j}^{0}} |{\rm graph} \, u_{i}|$$
where $N_{j} = \cup_{Z \in T_{V} \cap B_{1/32}^{n+1}(0)} \left(Z + \G^{Z} (N(H_{j}^{(0)})\right)$ and, for each $i=1, 2, \ldots, q_{j}^{(0)}$,  $u_{i}  \in C^{2}(\Omega^{\prime})$ with $u_{i}$ solving the minimal surface equation on $\Omega^{\prime}$, $|Du_{i}| < 1$,
$u_{1} \leq u_{2} \leq \ldots \leq u_{q_{j}^{(0)}}$ and, by the maximum principle, either $u_{i} \equiv u_{i+1}$ or $u_{i} < u_{i+1}$ for each $i=1, 2, \ldots, q_{j}^{(0)}-1.$ Since for each $i=1, 2, \ldots, q_{j}^{(0)},$ $u_{i}$ extends to $\overline{\Omega^{\prime}} \cap B_{1/64}^{n}(0)$ as a Lipschitz function 
with boundary values given by 
$\left.u_{j}\right|_{\partial \, \Omega^{\prime} \cap B_{1/64}^{n}(0)} (0, \varphi_{2}(y), y)= \varphi_{1}(y)$ for each point $(0, \varphi_{2}(y), y) \in \partial \, \Omega^{\prime} \cap B_{1/64}^{n}(0) = T^{\prime}_{V} \cap B_{1/64}^{n}(0),$ it follows from (\ref{MDthm-reg}) and standard $C^{1, \a}$ boundary regularity theory for uniformly elliptic equations (\cite{M}) that $u_{i} \in  C^{1, \a}(\overline{\Omega^{\prime}} \cap B_{1/64}^{n}(0)).$ 

We have thus established that $V \res B_{1/64}^{n+1}(0) = \sum_{j=1}^{2q+1} |M_{j}|$ where, for each $j \in \{1, 2, \ldots, 2q+1\}$, $M_{j}$ is an embedded $C^{1, \a}$ hypersurface-with-boundary with $\partial \, M_{j}  = T_{V} \cap B_{1/64}^{n+1}(0)$ and, for each $j, k \in \{1, 2, \ldots, 2q+1\},$ either $M_{j} \cap M_{k} = T_{V} \cap B_{1/64}^{n+1}(0)$
or $M_{j} = M_{k}.$ This directly contradicts hypothesis (${\mathcal S3}$) that $V$ is assumed to satisfy, completing the proof  of the theorem in case $\Theta \, (\|{\mathbf C}_{0}\|, 0) = q + 1/2$.\\   

\noindent
{\bf Case 2:  $\Theta \, (\|{\mathbf C}_{0}\|,0) = q + 1,$ $q \geq 2$}. Note that the validity of Theorem~\ref{no-transverse-q} in case $\Theta \, (\|{\mathbf C}_{0}\|,0) = q+1/2$ enables us to repeat the entire proof of Theorem~\ref{sheetingthm} with $q+1$ in place of $q$, yielding Theorem~\ref{sheetingthm} with $q+1$ in place of $q$. Consequently, the assertion of Remark (3) following the statement of Hypotheses~\ref{MD-hyp} holds with $q + 1$ in place of $q + 1/2$. Thus we may simply repeat (see Remark following the proof of Lemma~\ref{MD-no-gaps})  all of the steps of the above argument taking $q +1$ in place of $q + 1/2.$ This establishes Theorem~\ref{no-transverse-q} in case $\Theta \, (\|{\mathbf C}_{0}\|, 0) = q+1.$ 

The proof of Theorem~\ref{no-transverse-q} is now complete.
\end{proof}

%\medskip

\noindent
{\bf Remark:} The case $q=1$ of Theorem~$3.3^{\prime}$  is a special case of Allard's Regularity Theorem (which is reproduced by taking $q=1$ in our proofs of Lemma~\ref{excess-s} and Theorem~\ref{sheetingthm}). The validity of the case $\Theta \, (\|{\mathbf C}_{0}\|, 0) = 3/2$ of Theorem~\ref{no-transverse} follows from the validity of the case $q=1$ of Theorem~$3.3^{\prime}$; indeed, in this case, the same argument as for Theorem~\ref{no-transverse-q} carries over (with obvious simplifications) provided the induction hypothesis $(H1)$ is replaced by Theorem~$3.3^{\prime}$, case $q=1.$ In fact, when $\Theta \, (\|{\mathbf C}_{0}\|, 0) = 3/2$, Theorem~\ref{no-transverse} is true without the stability hypotheses (${\mathcal S{\emph2}}$) on $V$ (so $V$ only needs to be stationary and satisfy (${\mathcal S{\emph3}}$)); see [\cite{S}, Corollaries 2 and 3].
This in turn enables us to prove Theorem~\ref{no-transverse} in case $\Theta \, (\|{\mathbf C}_{0}\|, 0) = 2$, by
repeating the above proof of Theorem~\ref{no-transverse-q} (case $\Theta \,(\|{\mathbf C}_{0}\|, 0) = q +1$), taking 
$q= 1$ and, in place of induction hypotheses $(H1)$ and $(H2),$ case $q=1$ of Theorem~$3.3^{\prime}$ and  
case $\Theta \, (\|{\mathbf C}_{0}\|, 0) = 3/2$ of Theorem~\ref{no-transverse} respectively.

%\medskip

Theorem~\ref{sheetingthm} and Theorem \ref{no-transverse-q} together with the above remark and the remark preceding the statement of Theorem~$3.3^{\prime}$ complete the inductive proof of both Theorem~\ref{main} and Theorem~\ref{no-transverse}.

\section{The Regularity and Compactness Theorem}\label{reg-compactness}
\setcounter{equation}{0}

\noindent
\begin{proof}[Proof of Theorem~\ref{compactness}]  Note first that if $V \in {\mathcal S}_{\a},$ then it follows from Theorem~\ref{main}, Theorem~\ref{no-transverse} and Remark~3 of Section~\ref{outline} that ${\mathcal H}^{n-7+\g} \, ({\rm sing} \, V \cap (B_{2}^{n+1}(0)) = 0$ for each $\g >0$ if 
$n \geq 7$ and  ${\rm sing} \, V \cap B_{2}^{n+1}(0) = \emptyset$ if $2 \leq n\leq 6.$

Suppose, for each $k=1, 2, 3, \ldots,$ $V_{k} \in {\mathcal S}_{\a}$ and that
$$\Lambda = \limsup_{k \to \infty} \, \|V_{k}\|(B_{2}^{n+1}(0)) < \infty.$$ 
\noindent
By Allard's integer varifold compactness theorem, there exists a stationary integral varifold $V$ of $B_{2}^{n+1}(0),$ with 
$\|V\|(B_{2}^{n+1}(0)) < \Lambda + 1,$ such that, after passing to a subsequence, $V_{k} \to V$ as varifolds in $B_{2}^{n+1}(0).$ 
Set $K = {\rm sing} \, V \cap B_{2}^{n+1}(0).$ 

We argue that $V \in {\mathcal S}_{\a}$ as follows. By Theorem~\ref{main} and unique continuation of solutions to the minimal surface equation, if $M$ is a connected component of ${\rm reg} \, V$  and $0 < \r^{\prime} < \r < 2,$ there exists a number $\e = \e(M, \r, \r^{\prime}) \in (0, 1/2)$ such that for all sufficiently large $k$, ${\rm spt} \,\|V_{k}\| \cap \{X \in B_{\r}^{n+1}(0) \, : \, {\rm dist} \, (X, M \cap B_{\r}^{n+1}(0)) < \e\} \supset \cup_{j=1}^{q} {\rm graph} \, u_{j}^{k} \supset {\rm spt} \,\|V_{k}\| \cap \{X \in B_{\r^{\prime}}^{n+1}(0) \, : \, {\rm dist} \, (X, M \cap B_{\r^{\prime}}^{n+1}(0)) < \e\}$
for some integer $q \geq 1$ and functions $u_{j}^{k} \in C^{1, \a} \, (M \cap B_{\r}^{n+1}(0); M^{\perp})$ solving the minimal surface equation on $M \cap B_{\r}^{n+1}(0)$. It follows that $\int_{{\rm reg} \, V} |A|^{2}\z^{2} \leq \int_{{\rm reg} \,V}|\nabla \, \z|^{2}$ for each $\z \in C^{1}_{c}({\rm reg} \, V)$, where $A$ denotes the second fundamental form 
of ${\rm reg} \,V.$ It is also clear, from Theorem~\ref{no-transverse}, that $V$ satisfies the structural property $({\mathcal S 
{\emph 3}});$ for if not, there exists a point $Z \in {\rm spt} \, \|V\| \cap B_{2}^{n+1}(0)$ such that the (unique) tangent cone ${\mathbf C}_{Z}$ to $V$ at $Z$ is supported by the union of a finite number ($\geq 3$) of half-hyperplanes meeting along an $(n-1)$-dimensional subspace. By the definition of tangent cone and the fact that varifold convergence of stationary integral varifolds implies 
convergence in Hausdorff distance of the supports of the associated weight measures, for any given $\e_{1} >0$, 
there exists a number $\s \in (0, {\rm dist} \, (Z, \partial \, B_{2}^{n+1}(0)))$ such that for all sufficiently large $k$, 
${\rm dist} \, ({\rm spt} \, \|\eta_{Z, \s \, \#} \, V_{k}\| \cap B_{1}^{n+1}(0), {\rm spt} \, \|{\mathbf C}_{Z}\| \cap B_{1}^{n+1}(0)) < \e_{1}.$ This however contradicts Theorem~\ref{no-transverse} if we take $\e_{1} = \e(1/2, {\mathbf C}_{Z})$, where $\e$ is as in 
Theorem~\ref{no-transverse}. Thus $V \in {\mathcal S}_{\a},$ and hence ${\mathcal H}^{n-7+\g}(K) = 0$ for each $\g > 0$ if 
$n \geq 7$ and $K = \emptyset$ if  $2 \leq n \leq 6.$

Finally, suppose $n=7$ and consider any $V \in {\mathcal S}_{\a}$. To complete the proof of the theorem, it only remains to show that $K$ is discrete.  If this were false, there would 
exist points $Z, Z_{j} \in K,$ 
$j=1, 2, 3 \ldots,$ such that $Z_{j} \neq Z$ for each $j=1, 2, 3, \ldots$ and $Z_{j} \to Z$ as $j \to \infty.$ Letting $\s_{j} = |Z - Z_{j}|$, we obtain, passing to a subsequence without changing notation, a tangent cone ${\mathbf C}  = \lim_{j \to \infty} \, \eta_{Z, \s_{j} \, \#} \, V.$ By the discussion above, ${\mathbf C} \in {\mathcal S}_{\a}.$ Since $\s_{j}^{-1}(Z_{j} - Z) \in {\mathbf S}^{n-1} \cap 
{\rm sing} \, \eta_{Z, \s_{j} \, \#} \, V,$ it follows, passing to a further subsequence, that $\s_{j}^{-1}(Z_{j} - Z) \to Z^{\star} \in {\mathbf S}^{n-1}$ and by Hausdorff convergence and Theorem~\ref{main},  $Z^{\star} \in {\rm sing} \, {\mathbf C}.$ Since 
${\mathbf C}$ is a cone, it follows that $\{t Z^{\star} \, : \, t >0\} \subset {\rm sing} \, {\mathbf C}$, which is impossible 
since ${\mathbf C} \in {\mathcal S}_{\a}$ and we have established that for $n=7$, 
${\mathcal H}^{\g} \, (K) = 0$ for each $\g >0$ and any $V \in {\mathcal S}_{\a}.$
\end{proof}

\section{Generalization to Riemannian manifolds}\label{manifolds}
\setcounter{equation}{0}

Let $N$ be a smooth $(n+1)$-dimensional Riemannian manifold (without boundary) and for $X \in N$, let ${\rm exp}_{X}$ denote the exponential map at $X$.  For each $X \in N$, let
$R_{X} \in (0, \infty]$ be the injectivity radius at $X$.

Let ${\widetilde V}$ be a  stationary integral $n$-varifold on $N$.  Let $X_{0} \in {\rm spt} \, \|\widetilde{V}\|$,  ${\mathcal N}_{\r_{0}}(X_{0})$ be a normal coordinate ball of radius $\r_{0} \in (0, R_{X_{0}})$ around $X_{0}.$ 
Then $V = {\rm exp}_{X_{0} \,\#}^{-1} \, \widetilde{V} \res {\mathcal N}_{\r_{0}}(X_{0})$ is an integral $n$-varifold on $B_{\r_{0}}^{n+1}(0) \subset T_{X_{0}} \, N \approx {\mathbf R}^{n+1}$ which is stationary with respect to the functional  
\begin{equation}\label{manifolds-0}
{\mathcal F}_{X_{0}} (V) = \int_{B_{\r_{0}}^{n+1}(0) \times G_{n}} |\Lambda_{n} DF(X) \circ S| \, dV(X, S)
\end{equation}
where $F \equiv {\rm exp}_{X_{0}}.$ Let $\psi \in C^{1}_{c} \, (B_{\r_{0}}^{n+1}(0); {\mathbf R}^{n+1})$ and let $\varphi_{t}$, $t \in (-\e,\e)$ be the flow generated by $\psi.$ By computing directly the first variation $\d_{{\mathcal F}_{X_{0}}} \, V(\psi) \equiv\left.\frac{d}{dt}\right|_{t=0} {\mathcal F}_{X_{0}}(\varphi_{t \,\#} V)$ of $V$ with respect to ${\mathcal F}_{X_{0}}$ and setting 
$\d_{{\mathcal F}_{X_{0}}} \, V (\psi) = 0$, we see that the following bound holds (cf. \cite{SS} (1.7), (1.9) and (1.11)) for some constant $\m$ depending only on the metric on $N.$ (Such $\m \in (0, \infty)$ exists by replacing $N$ with a suitable open subset of $N$ if necessary):
\begin{itemize}
\item[(${\mathcal S^{\star}{\emph 1}}$)] 
$$\left|\int_{B_{\r_{0}}^{n+1}(0) \times G_{n}} {\rm div}_{S} \, \psi(X) \, dV(X, S) \right| \leq \m\int_{B_{\r_{0}}^{n+1}(0)} \left(|\psi(X)| + |X||\nabla\psi(X)|\right) \, d\|V\|(X)$$
\noindent
for all $\psi \in C_{c}^{1}(B_{\r_{0}}^{n+1}(0); {\mathbf R}^{n+1}).$ 
\end{itemize}

Furthermore, for $\psi \in C^{1}_{c} \, (B_{\r_{0}}^{n+1}(0) \setminus {\rm sing} \, V; {\mathbf R}^{n+1})$, the second variation $$\d^{2}_{{\mathcal F}_{X_{0}}} \, V(\psi) \equiv \left.\frac{d^{2}}{dt^{2}}\right|_{t=0} {\mathcal F}_{X_{0}}(\varphi_{t \, \#} V)$$ of $V$ with respect to ${\mathcal F}_{X_{0}}$ is given by (cf. \cite{SS} (1.8), (1.10), (1.12))
\begin{equation*}
\d^{2}_{{\mathcal F}_{X_{0}}} \, V(\psi) = \int_{{\rm reg} \, V} \left(\sum_{i=1}^{n} |(D_{\t_{i}} \, \psi)^{\perp}|^{2} + ({\rm div}_{\rm reg \, V} \, \psi)^{2}  - \sum_{i,j=1}^{n} (\t_{i} \cdot D_{\t{j}} \, \psi 
)\cdot(\t_{j} \cdot D_{\t_{i}} \, \psi)\right) \, d{\mathcal H}^{n} + R(\psi)
\end{equation*}
where $\{\t_{1}, \t_{2}, \ldots, \t_{n}\}$ is an orthonormal basis for the tangent space $T_{X} \, ({\rm reg} \, V)$ of ${\rm reg} \,V$ at $X$, $D_{\t} \, \psi$ denotes the directional derivative of $\psi$ in the direction $\t$ and   
$$|R(\psi)| \leq c\m  \int_{{\rm reg} \, V} \left(\widetilde{c}\m|\psi|^{2} + |\psi||\nabla \, \psi| + |X| |\nabla \, \psi|^{2}\right) \, d{\mathcal H}^{n}$$
with $c$, $\widetilde{c}$ absolute constants. If ${\rm reg} \, V$ is orientable and $\nu$ is a continuous choice of unit normal to ${\rm reg} \, V$, we may, for any  $\z \in C^{1}_{c} \, ({\rm reg} \, V)$,  extend $\z \nu$ to 
a vector field in $C^{1}_{c} \, (B_{\r_{0}}^{n+1}(0) \setminus {\rm sing} \,V; {\mathbf R}^{n+1})$ and take in the above 
$\psi = \z\nu$ to deduce that (cf. \cite{SS} (1.14), (1.15))
\begin{equation*}
\d^{2}_{{\mathcal F}_{X_{0}}} \, V(\psi) = \int_{{\rm reg} \, V} \left(|\nabla \, \z|^{2}  - |A|^{2}\z^{2} + H^{2}\z^{2}\right) \, d{\mathcal H}^{n} + R(\psi)
\end{equation*}
where $A$ denotes the second fundamental form of ${\rm reg} \, V$, $|A|$ the length of $A$, $H$ the mean curvature of ${\rm reg} \, V$  and
$$|R(\psi)| \leq c\m \int_{{\rm reg} \, V} \left(\widetilde{c}\m|\z|^{2} + |\z||\nabla \, \z| + \z^{2}|A| |X| |\nabla \, \z|^{2} + |X|\z^{2}|A|^{2}\right) \, d{\mathcal H}^{n}.$$ If $\d^{2}_{{\mathcal F}_{X_{0}}} \, (\psi) \geq 0$ for all $\psi  = \z \nu$, $\z \in C^{1}_{c}({\rm reg} \, V)$, then we have (cf. \cite{SS} (1.17))
\begin{itemize}
\item[(${\mathcal S^{\star}{\emph 2}}$)]
\begin{eqnarray*}
\int_{{\rm reg} \, V \cap B_{\r_{0}}^{n+1}(0)} |A|^{2} \z^{2} \, d{\mathcal H}^{n}\leq \int_{{\rm reg} \, V \cap B_{\r_{0}}^{n+1}(0)} |\nabla \, \z|^{2} \, d{\mathcal H}^{n}&&\\
&&\hspace{-3.5in} + c_{1}\m\int_{{\rm reg} \, V \cap B_{\r_{0}}^{n+1}(0)}
\left(c_{2}\m\z^{2} + \z|\nabla \, \z| + \z^{2}|A| + |X||\nabla \, \z|^{2} + |X| \z^{2}|A|^{2} + c_{2}\m|X|^{2} \z^{2} |A|^{2}\right) \, d{\mathcal H}^{n}
\end{eqnarray*}
\end{itemize}
for all $\z \in C^{1}_{c}({\rm reg} \, V)$ where $c_{1}$, $c_{2}$ are constants depending only on $n$.

For the rest of this discussion, we take $\m$, $c_{1}$, $c_{2}$ to be chosen as above and fixed.

\noindent
{\bf Definitions:} Let $\m$, $c_{1}$, $c_{2}$ be the positive numbers as above. 

\noindent
{\bf (1)} By a \emph{stable integral $n$-varifold $\widetilde{V}$ on} $N$ we mean a stationary integral $n$-varifold 
$\widetilde{V}$ on $N$ such that for each $X_{0} \in {\rm spt} \, \|\widetilde{V}\|$ and each normal ball ${\mathcal N}_{\r_{0}}(X_{0}) \subset N$ around $X_{0}$, the integral $n$-varifold $V = ({\rm exp}_{X_{0}}^{-1})_{\#} \, \widetilde{V} \res {\mathcal N}_{\r_{0}}(X_{0})$
on $B_{\r_{0}}^{n+1}(0) \subset {\mathbf R}^{n+1}$  satisfies (${\mathcal S^{\star}{\emph 2}}$).

\noindent
{\bf (2)} For $\a \in (0, 1)$, let $\widetilde{{\mathcal S}}_{\a}$ denote the collection of stable integral $n$-varifolds on $N$ satisfying the structural condition (${\mathcal S{\emph 3}}$) of Section~\ref{maintheorems} taken with normal ball ${\mathcal N}_{\r}(Z) \subset N$ in place of $B_{\r}^{n+1}(Z).$  

\noindent
{\bf (3)} For $\a \in (0, 1)$, let ${\mathcal S}_{\a}^{\star}$ denote the collection of integral $n$-varifolds $V$ on $B_{1}^{n+1}(0) \subset {\mathbf R}^{n+1}$ such that 
\begin{equation}\label{manifolds-1}
V = \eta_{0, \r \, \#} \, {\rm exp}_{X \, \#}^{-1} \, \widetilde{V} \res {\mathcal N}_{\r}(X)
\end{equation}
for some $\widetilde{V} \in \widetilde{{\mathcal S}}_{\a}$, $X \in {\rm spt} \, \|\widetilde{V}\|$ and $\r  \in (0, R_{X}).$

\noindent
{\bf (4)} For $\r \in (0, \infty)$ and $\a \in (0, 1)$, let ${\mathcal S}_{\a}^{\star}(\r)$
be the set of integral $n$-varifolds $V \in {\mathcal S}_{\a}^{\star}$ such that (\ref{manifolds-1}) holds for 
some $\widetilde{V} \in \widetilde{{\mathcal S}}_{\a}$ and $X \in {\rm spt} \, \|\widetilde{V}\|$ with 
$R_{X} \geq \r.$ 

\noindent
{\bf Remark:} Let $\r \in (0, \infty)$ and suppose that $V \in {\mathcal S}_{\a}^{\star}(\r)$. Then for each $Y \in {\rm spt} \, \|V\| \cap B_{1/2}^{n+1}(0)$, 
\begin{equation}\label{manifolds-2}
\eta_{0, \r/2 \, \#} \, \t_{Y \, \#} \, V \in {\mathcal S}_{\a}^{\star}(\r/2)
\end{equation}
where $\t_{Y} = {\rm exp}^{-1}_{{\rm exp}_{X}(\r Y)} \circ {\rm exp}_{X} \circ \eta_{0, \r^{-1}}.$ Note that 
$\t_{Y}(Y)  =0.$ 

We assert that the following direct analog of Theorem~\ref{compactness}  holds:

\begin{theorem}[\sc Regularity and Compactness Theorem---Manifold version]\label{compactness-N}
Let $N$ be a smooth $(n+1)$-dimensional Riemannian manifold, $X_{0} \in N$ and $\a \in (0, 1/2).$ Let $\{\widetilde{V}_{k}\} \subset \widetilde{{\mathcal S}}_{\a}$ be a sequence with $X_{0} \in {\rm spt} \, \|\widetilde{V}_{k}\|$ for each $k=1, 2, \ldots,$ and with $$\limsup_{k \to \infty} \, \|\widetilde{V}_{k}\|(N)< \infty.$$ 
Then there exist a subsequence $\{k^{\prime}\}$ of $\{k\}$ and a varifold $\widetilde{V} \in \widetilde{{\mathcal S}}_{\a}$  with
$X_{0} \in {\rm spt} \, \|\widetilde{V}\|$ and with
${\mathcal H}^{n-7+\g} \,({\rm sing} \ \widetilde{V} \cap N) = 0$ for each $\g >0$ if $n \geq 7$, ${\rm sing} \, \widetilde{V} \cap N$ discrete if $n= 7$ and ${\rm sing} \,\widetilde{V} \cap N = \emptyset$ if $2 \leq n \leq 6$  
such that $\widetilde{V}_{k^{\prime}} \to \widetilde{V}$ as varifolds of $N$ and smoothly (i.e. 
in the $C^{m}$ topology for every $m$) locally in $N \setminus {\rm sing} \, \widetilde{V}$.

In particular, if $\widetilde{W} \in \widetilde{{\mathcal S}}_{\a}$, then ${\mathcal H}^{n-7+\g} \,({\rm sing} \, \widetilde{W} \cap N) = 0$ for each $\g >0$ if $n \geq 7$, ${\rm sing} \, \widetilde{W} \cap N$ is discrete if $n= 7$ and ${\rm sing} \,\widetilde{W} \cap N = \emptyset$ if $2 \leq n \leq 6.$  
\end{theorem}

By the preceding discussion, this theorem is equivalent to the assertion obtained from it by replacing  $N$ with $B_{1}^{n+1}(0) \subset {\mathbf R}^{n+1},$ $X_{0}$ with $0$ and $\widetilde{{\mathcal S}}_{\a}$ with ${\mathcal S}_{\a}^{\star};$ the proof of the latter amounts to making minor modifications, as described below, to the proof of Theorem~\ref{compactness}.

\noindent
{\bf Step 1}: Let $V$ be an integral $n$-varifold of $B_{1}^{n+1}(0)$ such that (\ref{manifolds-1}) holds for some 
stationary integral $n$-varifold $\widetilde{V}$ of $N$, $X_{0} \in {\rm spt} \,\|\widetilde{V}\|$ in place of $X$ and $\r_{0} \in (0, R_{X_{0}})$ in place of $\r.$ By the discussion involving (5.3)--(5.9) of \cite{SS}, we have, for each $0 < \s < \d$, where $\d = \d(n, \m\r_{0}) \in (0, 1)$, the following facts:
\begin{equation}\label{monotonicity-N-0}
\t^{-n}\|V\|(B_{\t}^{n+1}(0)) \leq (1 + 12n\m\r_{0}\s)\s^{-n}\|V\|(B_{\s}^{n+1}(0))
\end{equation}
for all $\t$ with $0 < \t \leq \s$; the density 
$\Theta \, (\|V\|, 0) = \lim_{\t \to 0} \, \frac{\|V\|(B_{\t}^{n+1}(0))}{\omega_{n}\t^{n}}$ exists ( and is finite); the function 
$\Theta \, (\|\cdot\|, 0)$ is upper semi-continuous; 
\begin{equation}\label{monotonicity-N}
\int_{B^{n+1}_{\s}(0)} \frac{|X ^{\perp}|^{2}}{|X|^{n+2}} d\|V\|(X) \leq 
\frac{\|V\|(B_{\s}^{n+1}(0))}{\omega_{n}\s^{n}} - \Theta \, (\|V\|, 0) + C\s \frac{\|V\|(B_{\s}^{n+1}(0))}{\omega_{n}\s^{n}}
\end{equation} 
where $C = C(n, \m\r_{0}) \in (0, \infty)$; tangent cones to $V$ at  
$0 \in {\rm spt} \,\|V\|$ exist and are stationary integral hypercones of ${\mathbf R}^{n+1}$. 

Let ${\rm Var Tan} \, (V, 0)$ denote the set of tangent cones to $V$ at $0.$ For $Y \in {\rm spt} \, \|V\| \cap B_{1/2}(0)$, let 
$\Theta \, (\|V\|, Y) = \Theta \, (\|\eta_{0, \r_{0}/2 \, \#} \, \t_{Y \, \#} V\|, 0)$ (see \ref{manifolds-2}) and 
${\rm Var Tan} \, (V, Y) = {\rm Var Tan} \, (\eta_{0, \r_{0}/2 \, \#} \, \t_{Y \, \#} V, 0)$. Recalling the well known fact that if 
${\mathbf C}$ is a stationary cone in a Euclidean space ${\mathbf R}^{m}$, then the set $\{Z \in {\mathbf R}^{m} \, : \, \Theta \, (\|{\mathbf C}\|, Z) = \Theta \, (\|{\mathbf C}\|, 0)\}$ is a linear subspace of ${\mathbf R}^{m},$ we deduce by the argument of Almgren's generalized stratification of stationary integral varifolds (\cite{A}, Remark 2.28; see also \cite{S3}, Section 3.4) the following:

\noindent
\emph{Let $V$ be an integral $n$-varifold of $B_{1}^{n+1}(0)$ such that (\ref{manifolds-1}) holds for some 
stationary integral $n$-varifold $\widetilde{V}$ of $N$, $X \in {\rm spt} \,\|\widetilde{V}\|$ and $\r \in (0, R_{X}).$  For $k=0, 1, 2, \ldots, n,$  let $S_{k} = \{Y \in {\rm spt} \, \|V\| \cap B_{1/2}^{n+1}(0)\, : \, 
{\rm dim} \, \{Z \in {\mathbf R}^{n+1} \, : \, \Theta \, (\|{\mathbf C}\|, Z) = \Theta \, (\|{\mathbf C}\|, 0)\} \leq k 
\;\; \forall {\mathbf C} \in {\rm Var Tan} \, (V,Y)\}.$ Then ${\rm dim}_{\mathcal H} \, (S_{k}) \leq k.$} 

\noindent
{\bf Step 2:} We claim that the following analogs of Theorems~\ref{main} and \ref{no-transverse} hold. 

\begin{theorem}[\sc Sheeting Theorem---Manifold Version]\label{sheeting-N}
Let $\a \in (0, 1/2)$, $\r_{0}  \in (0, \infty)$ and $q$ be any integer $\geq 1.$ Let $\a^{\prime} = (2\a + 1)/4.$ There exists a number 
$\e_{0} = \e_{0}(n, q, \a, \m\r_{0}) \in (0, 1)$ such that if $V \in {\mathcal S}^{\star}_{\a}(\r_{0})$, $\omega_{n}^{-1}\|V\|(B_{1}^{n+1}(0)) < q + 1/2$, $\s \in (0, 1/2)$, 
$(q-1/2) \leq \left(\omega_{n}\s^{n}\right)^{-1}\|V\| (B_{\s}^{n+1}(0)) < (q + 1/2)$ and 
$\s^{-1}{\rm dist}_{\mathcal H} \, ({\rm spt} \, \|V\| \cap ({\mathbf R} \times B_{\s}), \{0\} \times B_{\s}) + \s^{2\a^{\prime}} < \e_{0}$ then 
$$V \res ({\mathbf R} \times B_{\s/2}) = \sum_{j=1}^{q} \, |{\rm graph} \, u_{j}|$$
where $u_{j} \in C^{1, \b}(B_{\s/2})$ for each $j=1, 2, \ldots, q$; $u_{1} \leq u_{2} \leq \ldots \leq u_{q}$ and 
\begin{eqnarray*}
\s^{-1}\sup_{B_{\s/2}} \, |u_{j}| + \sup_{B_{\s/2}} \, |Du_{j}| + \s^{\b}\sup_{X_{1}, X_{2} \in B_{\s/2},\, X_{1} \neq X_{2}} \, \frac{|Du_{j}(X_{1}) - Du_{j}(X_{2})|}{|X_{1} - X_{2}|^{\b}}&&\\
&&\hspace{-2in} \leq C\, \left(\s^{-n-2}\int_{{\mathbf R} \times B_{\s}}|x^{1}|^{2} \, d\|V\|(X) + \s^{2\a^{\prime}}\right)^{1/2}.
\end{eqnarray*}
Here $C = C(n, q, \a, \m\r_{0}) \in (0, \infty)$ and $\b = \b(n, q, \a, \m\r_{0}) \in (0, 1)$.
\end{theorem}

\noindent
{\bf Remark:} If the conclusions of Theorem~\ref{sheeting-N} hold, and $V$ corresponds, as in (\ref{manifolds-1}), to 
some $\widetilde{V} \in \widetilde{{\mathcal S}}_{\a},$  $X = X_{0} \in N \cap {\rm spt} \, \|\widetilde{V}\|$ and $\r = \r_{0} \in (0, R_{X_{0}})$,  then it follows that for each $j \in \{1, 2, \ldots, q\}$, 
$V_{j} \equiv |{\rm graph} \, \r_{0}u_{j}(\r_{0}^{-1}(\cdot))|$  is stationary with respect to the functional 
${\mathcal F} (\cdot) = {\mathcal F}_{X_{0}} \, \left((\cdot) \res {\mathbf R} \times B_{\s/2}\right)$ where ${\mathcal F}_{X_{0}}$ is as in (\ref{manifolds-0}); thus, by computing the associated Euler-Lagrange equation and applying elliptic regularity theory, we see that $u_{j} \in C^{\infty}(B_{\s/2})$ and satisfies an equation of the form 
\begin{equation}\label{manifolds-2-1}
\sum_{k, \ell=1}^{n} a_{k\ell}^{j} D_{k}D_{\ell} u_{j}  =f^{j}
\end{equation}
on $B_{\s/2}$, with $|f^{j}(x)| \leq \m\r_{0}$ and $a_{k\ell}(x) = \d_{k\ell} - \frac{D_{k}u_{j}(x) D_{\ell}u_{j}(x)}{\sqrt{1 + |Du_{j}(x)|^{2}}} + b_{k\ell}^{j}(x)$, where $|b_{k\ell}^{j}(x)| \leq \m\r_{0}\s,$ for $x \in B_{\s/2}.$ 

\begin{theorem}[\sc Minimum Distance Theorem---Manifold Version]\label{mindist-N}
Let $\a \in (0, 1/2),$ $\r_{0} \in (0,\infty)$ and $\g \in (0, 1/2).$ Let $\a^{\prime} = (2\a + 1)/4$. Suppose that ${\mathbf C}_{0}$ is an $n$-dimensional stationary cone in ${\mathbf R}^{n+1}$ such that 
${\rm spt} \, \|{\mathbf C}_{0}\|$ is equal to a finite union of at 
least three distinct $n$-dimensional half-hyperplanes of ${\mathbf R}^{n+1}$ meeting along an $(n-1)$-dimensional subspace. Then there exists $\e = \e(\a, \g, \m\r_{0}, {\mathbf C}_{0}) \in (0, 1)$  such that if 
$V \in {\mathcal S}^{\star}_{\a}(\r_{0}),$ $\s \in (0, 1/2)$, $\Theta \, (\|V\|, 0) \geq \Theta \, (\|{\mathbf C}_{0}\|, 0)$ and
$(\omega_{n})^{-1}\|V\|(B_{1}^{n+1}(0)) \leq \Theta_{{\mathbf C}_{0}}(0) +\g$ then 
$$\s^{\a^{\prime}} + \s^{-1}{\rm dist}_{\mathcal H} \, ({\rm spt} \, \|V\| \cap B_{\s}^{n+1}(0), {\rm spt} \, \|{\mathbf C}_{0}\| \cap B_{\s}^{n+1}(0)) \geq \e.$$
In particular  $\s^{-1}{\rm dist}_{\mathcal H} \, ({\rm spt} \, \|V\| \cap B_{\s}^{n+1}(0), {\rm spt} \, \|{\mathbf C}_{0}\| \cap B_{\s}^{n+1}(0)) \geq \e/2$ for sufficiently small $\s > 0.$
 \end{theorem}

The proof of Theorems~\ref{sheeting-N} and ~\ref{mindist-N}  amounts to an easy modification of the induction argument given above for Theorems~$3.3^{\prime}$ and \ref{no-transverse}, which is the ``Euclidean case,'' viz.  the case when $\m = 0$ (which corresponds to the case when 
$N$ is an open subset of ${\mathbf R}^{n+1}$ in Theorem~\ref{compactness-N}). We outline the proof as follows:

\noindent
(i) It follows from \cite{SS}, Theorem 1, that Theorem~\ref{sheeting-N} holds if $V$ satisfies, in place of the structural condition (${\mathcal S{\emph 3}}$), that
$${\rm dim}_{\mathcal H} \, ({\rm sing \, V}) \leq n-7 \;\; \mbox{in case} \;\; n \geq 7\;\; {\rm and} \;\; {\rm sing} \, V = \emptyset \;\; \mbox{in case} \;\; n \leq 6,$$ 
together with all other hypotheses as in Theorem~\ref{sheeting-N}.
 
 \noindent
(ii) Let $\r_{0} \in (0, \infty)$ and let $V$ be an integral $n$-varifold on $B_{1}^{n+1}(0)$ such that (\ref{manifolds-1}) holds with $\r = \r_{0}$ for some stationary integral $n$-varifold $\widetilde{V}$ on $N$ and $X_{0} \in {\rm spt} \, \|\widetilde{V}\|$ with $R_{X_{0}} \geq  \r_{0}.$ Let $\s \in (0, 1),$ $\Lambda \in [1, \infty)$ and suppose that $(\omega_{n}\s^{n})^{-1}\|V\|(B_{\s}^{n+1}(0)) \leq \Lambda$ and  $\s^{-n-2}\int_{{\mathbf R} \times B_{\s}} |x^{1}|^{2} \, d\|V\|(X) + \s < 1.$ 
By taking $\psi(X) = x^{1}\widetilde{\z}^{2}(X)e^{1}$ in (${\mathcal S^{\star}\emph{1}}$), where $\widetilde{\z} \in C^{1}_{c}({\mathbf R} \times B_{3/4})$, we deduce that 
\begin{equation}\label{manifolds-3}
\int_{{\mathbf R} \times B_{3/4}} |\nabla \, x^{1}|^{2}\widetilde{\z}^{2} \, d \, \|\eta_{0, \s \, \#} \, V\|(X)  \leq C\left(\int_{{\mathbf R} \times B_{3/4}}|x^{1}|^{2}|\nabla \, \widetilde{\z}|^{2} \, d \,\|\eta_{0, \s \, \#} \, V\|(X) + \s\right)
\end{equation}
for each $\widetilde{\z} \in C^{1}_{c}({\mathbf R} \times B_{3/4})$ where $C = C(n, \Lambda, M, \m\r_{0}) \in (0, \infty),$ and $M = \sup_{{\rm spt} \, \|\eta_{0, \s \, \#} \, V\| \cap ({\mathbf R} \times B_{3/4})} \, |\widetilde{\z}| + |D\widetilde{\z}|.$ Choosing $\widetilde{\z}$ such that $\widetilde{\z}(x^{1}, x^{\prime}) = \z(x^{\prime})$ in a neighborhood of ${\rm spt} \, \|\eta_{0, \s \, \#} \, V\| \cap ({\mathbf R} \times B_{3/4}),$ where $\z \in C^{1}_{c}(B_{3/4})$ is such that 
$\z \equiv 1$ on $B_{1/2}$, $0 \leq \z \leq 1$ and $|D\z| \leq 8$, we deduce from this that
\begin{equation}\label{manifolds-4}
\int_{{\mathbf R} \times B_{1/2}} |\nabla \, x^{1}|^{2} \, d \, \|\eta_{0, \s \, \#} \, V\|(X)  \leq C\left(\int_{{\mathbf R} \times B_{3/4}}|x^{1}|^{2}\, d \,\|\eta_{0, \s \, \#} \, V\|(X) + \s\right)
\end{equation}
where $C = C(n, \Lambda, \m\r_{0}).$ 

\noindent
(iii) Let $\r_{0}$, $V$ be as in (ii) and let $\s \in (0, 3/4)$. With $\eta_{0,\s \, \#} \, V$ in place of $V$,  $$\sqrt{\s^{-n-2}\int_{{\mathbf R} \times B_{\s}} |x^{1}|^{2} \, d\|V\|(X) + \s}$$ in place of ${\hat E}_{V}$ and with the constants $\e_{0}$, $C$ depending on $n$, $q$, $\m\r_{0},$ Theorem~\ref{flat-varifolds} holds; its proof amounts to modifying the argument of \cite{A}, Theorem 3.8 in obvious ways, making use of (\ref{monotonicity-N-0}), (\ref{monotonicity-N}) and (\ref{manifolds-4}).
  
\noindent
(iv) Consequently, the case $q=1$ of Theorem~\ref{sheeting-N}  follows by the excess improvement argument as in [\cite{AW}, Chapter 8].

\noindent
(v) From (iii) and the inequalities (\ref{monotonicity-N}), (\ref{manifolds-3}), we deduce that for $\r_{0}$, $V$ as in (ii) and $\s \in (0, 3/4)$, Theorem~\ref{non-concentration}  hold with $\eta_{0, \s \, \#} \, V$ in place of $V$ and 
$$\sqrt{\s^{-n-2}\int_{{\mathbf R} \times B_{\s}} |x^{1}|^{2} \, d\|V\|(X) + \s}$$ in place of ${\hat E}_{V}$, again with the constants $\e_{1}$, $C$ etc. depending also on $\m\r_{0}.$ 

\noindent
(vi) For what follows, fix $\a \in (0, 1/2)$, $\r_{0} \in (0, \infty),$ and let $\a^{\prime} = (2\a + 1)/4.$ For $V \in {\mathcal S}_{\a}^{\star}(\r_{0})$ and $\s \in (0, 3/4)$ let 
$${\hat E}^{\star}_{V}(\s) = \sqrt{\s^{-n-2}\int_{{\mathbf R} \times B_{\s}} |x^{1}|^{2} \, d\|V\|(X) + \s^{2\a^{\prime}}}.$$

Let $q$ be an integer $\geq 2$, and assume inductively the validity of Theorem~\ref{sheeting-N} with $1, 2, \ldots, q-1$ in place of $q$ and that of Theorem~\ref{mindist-N} if $\Theta \,(\|{\mathbf C}_{0}\|, 0) \in \{3/2,2, 5/2, \ldots, q - 1/2, q\}.$

\noindent
(vii) For each $k=1, 2, 3, \ldots,$ let $\s_{k} \in (0, 3/4)$, $V_{k} \in {\mathcal S}^{\star}_{\a}(\r_{0})$ be such that 
$\omega_{n}^{-1}\|V_k\|(B_{1}^{n+1}(0)) < q + 1/2$, $\s_{k} \to 0$ and $(q-1/2) \leq \left(\omega_{n}\s_{k}^{n}\right)^{-1}\|V_{k}\| (B_{\s_{k}}^{n+1}(0)) < (q + 1/2).$ If ${\hat E}^{\star}_{V_{k}}(\s_{k}) \to 0$, then as in the discussion following Theorem~\ref{flat-varifolds}, we may blow up the sequence $\{\eta_{0, \s_{k} \, \#} \, V_{k} \res B_{1}^{n+1}(0)\}$ by  ${\hat E}^{\star}_{V_{k}}(\s_{k})$. We shall continue to call a function $v  \in W^{1, 2}_{\rm loc} \, (B_{1} ; {\mathbf R}^{q}) \cap L^{2} \, (B_{1} ; {\mathbf R}^{q})$ produced this way a coarse blow-up.

\noindent
(viii) By the reasoning of Remarks 2, 3 of Section~\ref{outline} and Step 1 above, we have the following:

\noindent
 \emph{Let $q$ be an integer $\geq 2$ and suppose that the induction hypotheses as in {\rm (vi)} hold. If $V \in {\mathcal S}^{\star}_{\a}(\r_{0}),$ $\Omega \subseteq B_{1}^{n+1}(0)$ is open and $\Theta \, (\|V\|, Z) < q$ for each $Z \in {\rm spt} \, \|V\| \cap \Omega$, then 
${\mathcal H}^{n-7 + \g} \, ({\rm sing} \, V \, \res \, \Omega) = 0$ for each $\g > 0$ if $n \geq 7$ and 
${\rm sing} \, V \, \res \, \Omega = \emptyset$ if $2 \leq n \leq 6.$}

\noindent
(ix) The collection ${\mathcal B}_{q}^{\star}$ of all coarse blow-ups $v$ (as in (vii)) is a proper blow-up class, viz. ${\mathcal B}_{q}^{\star}$ satisfies properties (${\mathcal B{\emph 1}}$)--(${\mathcal B{\emph 7}}$) of Section~\ref{proper-blow-up}. Verification of properties
(${\mathcal B{\emph 1}}$)--(${\mathcal B{\emph 3}}$), (${\mathcal B{\emph 5}}$) and (${\mathcal B{\emph 6}}$) proceeds in the same way as for the Euclidean case described in Section~\ref{step2} above. In view of (i), property 
(${\mathcal B{\emph 4}}$) follows from the corresponding argument for the Euclidean case, also described in Section~\ref{step2}, with the inequality 
(\ref{monotonicity-N}) taking the place of the monotonicity identity (\ref{monotonicity-1}).

Property (${\mathcal B{\emph 7}}$) is verified by establishing separately the same two cases as {\bf Case 1} and {\bf Case 2} of Section~\ref{step3}.  With regard to {\bf Case 1}, 
note that by taking $\psi(X) = \widetilde{\z}(X)e^{2}$ in (${\mathcal S^{\star}{\emph 1}}$), where $\widetilde{\z} \in C^{1}_{c}({\mathbf R} \times B_{3/4}),$ it follows that for each $k=1,2, \ldots,$
\begin{eqnarray*}
\left|\int_{{\mathbf R} \times B_{3/4}} \nabla \, x^{2} \cdot \nabla \,\widetilde{\z} \, d\|\eta_{0, \s_{k} \, \#} \, V_{k}\|(X) \right| &\leq& C \sup_{{\rm spt} \, \|\eta_{0, \s_{k} \, \#} \, V_{k}\| \cap ({\mathbf R} \times B_{3/4})} \, \left(|\widetilde{\z}| + |D\widetilde{\z}|\right)\s_{k}\\ 
&\leq&  C \sup_{{\rm spt} \, \|\eta_{0, \s_{k} \, \#} \, V_{k}\| \cap ({\mathbf R} \times B_{3/4})} \, \left(|\widetilde{\z}| + |D\widetilde{\z}|\right)\s_{k}^{1 - 2\a^{\prime}} \left({\hat E}_{V_{k}}^{\star}(\s_{k})\right)^{2}
\end{eqnarray*}
where $C = C(n, q, \m\r_{0}).$ {\bf Case 1} is established by taking this in place of (\ref{no-overlap-main}) and (\ref{manifolds-3}) in place of (\ref{tilt-ht-0}) in the argument of Lemma~\ref{no-overlap}.
With regard to {\bf Case 2}, we note that the following analogue of Lemma~\ref{excess-improvement} holds. Here 
${\mathcal C}_{q}$, ${\mathcal C}_{q}(p)$ are as defined in Section~\ref{fineblowup}.
\begin{lemma}\label{excess-improvement-N}
Let $q$ be an integer $\geq 2$,  $\a \in (0, 1/2)$, $\th \in (0, 1/4)$ and $\r_{0} \in (0, \infty)$. 
There exist numbers $\overline\e = \overline\e(n, q, \a, \th, \m\r_{0}) \in (0, 1/2)$, $\overline\g = \overline\g(n, q, \a, \th, \m\r_{0}) \in (0, 1/2)$ and 
$\overline\b = \overline\b(n, q, \a,\th,\m\r_{0}) \in (0, 1/2]$ such that the following 
is true: Let $\s \in (0, 1)$ and suppose that the induction hypotheses as in (vi) and the following hold. 
\begin{itemize}
\item[(1)] $V \in {\mathcal S}_{\a}^{\star}(\r_{0})$, \; $\Theta \, (\|V\|, 0) \geq q$, \; $(\omega_{n}\s^{n})^{-1}\|V\|(B_{\s}^{n+1}(0)) < q + 1/2.$
\item[(2)] ${\mathbf C} = \sum_{j=1}^{q} |H_{j}| + |G_{j}| \in {\mathcal C}_{q},$ where for each $j \in \{1, 2, \ldots, q\}$, $H_{j}$ is the half-space defined by
$H_{j} = \{(x^{1}, x^{2}, y) \in {\mathbf R}^{n+1} \, : \, x^{2} < 0 \;\; \mbox{and} \;\; x^{1} = \lambda_{j}x^{2}\},$ 
$G_{j}$ the half-space defined by $G_{j} = \{(x^{1},x^{2}, y) \in {\mathbf R}^{n+1}\, : \, x^{2} > 0 \;\; \mbox{and} \;\; x^{1} = \mu_{j}x^{2}\},$ with $\lambda_{j}, \mu_{j}$ constants,  
$\lambda_{1} \geq \lambda_{2} \geq \ldots \geq \lambda_{q}$ and $\mu_{1} \leq \mu_{2} \leq \ldots \leq \mu_{q}$. 
\item[(3)] $\left({\hat E}^{\star}_{V}(\s)\right)^{2} \equiv \int_{{\mathbf R} \times B_{1}} |x^{1}|^{2} d\|\eta_{0, \s \, \#} \, V\|(X) + \s^{2\a^{\prime}}< \e,$ where $\a^{\prime} = (2\a + 1)/4.$
\item[(4)] $\{Z \, : \, \Theta \, (\|\eta_{0, \s \, \#} \, V\|, Z) \geq q\} \cap \left({\mathbf R} \times (B_{1/2} \setminus \{|x^{2}| < 1/16\}) \right) = \emptyset.$
\item[(5)] 
\begin{eqnarray*}
\int_{{\mathbf R} \times (B_{1/2} \setminus \{|x^{2}| < 1/16\})} {\rm dist}^{2}(X, {\rm spt} \, \|\eta_{0, \s \, \#} \, V\|) \,d\|{\mathbf C}\|(X)&&\\ 
&&\hspace{-2in}+ \; \int_{{\mathbf R} \times B_{1}} {\rm dist}^{2} \, (X, {\rm spt}\, \|{\mathbf C}\|) \, d\|\eta_{0, \s \, \#} \, V\|(X) \leq\g \left({\hat E}_{V}^{\star}(\s)\right)^{2}.
\end{eqnarray*}
\item[(6)] $\left({\hat E}^{\star}_{V}\right)^{2} < \frac{3}{2}M_{0}\inf_{\{P \in G_{n} \, : \, P \cap (\{0\} \times {\mathbf R}^{n}) = \{0\} \times {\mathbf R}^{n-1}\}} \, \int_{{\mathbf R} \times B_{1}} {\rm dist}^{2} \, (X, P) , d\|V\|(X) + \s^{2\a^{\prime}},$ where $M_{0} = M_{0}(n, q) \in (1, \infty)$
is the constant defined in Section~\ref{fineblowup}. 
\item[(7)] Either
\begin{itemize}
\item[(i)] ${\mathbf C} \in {\mathcal C}_{q}(4)$ or
\item[(ii)] $q \geq 3$, ${\mathbf C} \in {\mathcal C}_{q}(p)$ for some $p \in \{5, \ldots, 2q\}$ and 
\begin{eqnarray*}
\int_{{\mathbf R} \times (B_{1/2} \setminus \{|x^{2}| < 1/16\})} {\rm dist}^{2}(X, {\rm spt} \, \|\eta_{0, \s \, \#} \, V\|) \,d\|{\mathbf C}\|(X) + \int_{{\mathbf R} \times B_{1}} {\rm dist}^{2} \, (X, {\rm spt}\, \|{\mathbf C}\|) \, d\|\eta_{0, \s \, \#} \, V\|(X)&&\\ 
&&\hspace{-5in}\leq\b \inf_{\widetilde{\mathbf C}\in \cup_{k=4}^{p-1}{\mathcal C}_{q}(k)} \, \left(\int_{{\mathbf R} \times (B_{1/2} \setminus \{|x^{2}| < 1/16\})} {\rm dist}^{2}(X, {\rm spt} \, \|\eta_{0, \s \, \#} \, V\|) \,d\|\widetilde{\mathbf C}\|(X)\right.\\
&&\hspace{-3in} + \left.\int_{{\mathbf R} \times B_{1}} {\rm dist}^{2} \, (X, {\rm spt}\, \|\widetilde{\mathbf C}\|) \, d\|\eta_{0, \s \, \#} \, V\|(X)\right).
\end{eqnarray*}
\end{itemize}
\end{itemize}
Then there exists an orthogonal rotation $\G$ of ${\mathbf R}^{n+1}$ and a cone ${\mathbf C}^{\prime} \in {\mathcal C}_{q}$ 
such that the conclusions of Lemma~\ref{excess-improvement} hold with 
$\eta_{0, \s \, \#} \, V$ in place of $V$, ${\hat E}_{V}^{\star}(\s)$  in place of ${\hat E}_{V},$ 
$$E_{V}^{\star}({\mathbf C}, \s) \equiv \sqrt{\int_{{\mathbf R} \times B_{1}} {\rm dist}^{2} \, (X, {\rm spt} \, \|{\mathbf C}\|) \, d\|\eta_{0, \s \, \#} \, V\|(X) + \s^{2\a^{\prime}}}$$ in place of $E_{V}$ and with the constants $\overline{\k}$, $\overline{C}_{0}$, $\overline{\g}_{0}$, $\overline{\n}$, $\overline{C}_{1}$, $\overline{C}_{2} \in (0, \infty)$ depending only on $n$, $q$, $\a$ and $\m\r_{0}$.
\end{lemma}

In proving this, note first that if $\s^{2\a^{\prime}} > \int_{{\mathbf R} \times B_{1}} |x^{1}|^{2} \, d\|\eta_{0, \s \, \#} \, V\|(X)$, then, provided $\g < \th^{n+4}/2$, we trivially have that 
\begin{eqnarray*}
\th^{-n-2}\int_{{\mathbf R} \times B_{\th}} {\rm dist}^{2} \, (X,{\rm spt} \, \|{\mathbf C}\|) \, d\|\eta_{0, \s \, \#} \, V\|(X) \leq \th^{-n-2}\g \left({\hat E}_{V}^{\star}(\s)\right)^{2}&&\\
&&\hspace{-2in} \leq 2\th^{-n-2}\g \s^{2\a^{\prime}} \leq \th^{2}\s^{2\a^{\prime}} \leq \th^{2}\left(E_{V}^{\star}({\mathbf C}, \s)\right)^{2}
\end{eqnarray*}
and thus the conclusions (a)--(d) hold with ${\mathbf C}^{\prime} = {\mathbf C}$ and $\G$= Identity, and the conclusions (e) and (f) can be checked as in the proof of Lemma~\ref{excess-improvement}. Hence we may assume without loss of generality that $\left({\hat E}_{V}^{\star}(\s)\right)^{2} \leq 2\int_{{\mathbf R} \times B_{1}}|x^{1}|^{2} \, d\|\eta_{0, \s \, \#} \, V\|(X).$ With this additional assumption and with the help of inequality (\ref{monotonicity-N}), the obvious analogues of Theorem~\ref{L2-est-1} and Corollary~\ref{L2-est-2} can be established; consequently, Lemma~\ref{excess-improvement-N} can be proved by making obvious modifications to the entire argument leading to Lemma~\ref{excess-improvement}, as described in Sections~\ref{fineblowup}--\ref{finedecay}.

The obvious analog of Lemma~\ref{multi-scale-final} then follows; note in particular that in the conclusions of this modified lemma we must take
\begin{eqnarray*}
Q_{V}^{\star}({\mathbf C}, \s) \equiv \left(\int_{{\mathbf R} \times (B_{1/2} \setminus \{|x^{2}| < 1/16\})} {\rm dist}^{2}(X, {\rm spt} \, \|\eta_{0, \s \, \#} \, \eta_{0, \s \, \#} \, V\|) \,d\|{\mathbf C}\|(X)\right.&&\\
&&\hspace{-3in} +\; \left.\int_{{\mathbf R} \times B_{1}} {\rm dist}^{2} \, (X, {\rm spt}\, \|{\mathbf C}\|) \, d\|\eta_{0, \s \, \#} \, V\|(X) + \s^{2\a^{\prime}}\right)^{1/2}
\end{eqnarray*}
in place of $Q_{V}$, and that the modified lemma yields that for some $j \in \{1, 2, \ldots, 2q-3\}$, ${\mathbf C}^{\prime} \in {\mathcal C}_{q}$ and some orthogonal rotation $\G$ of ${\mathbf R}^{n+1}$,  
\begin{eqnarray*}
\int_{{\mathbf R} \times (B_{1/2} \setminus \{|x^{2}| < 1/16\})} {\rm dist}^{2}(X, {\rm spt} \, \|\eta_{0, \th_{j}\s \, \#} \, \eta_{0, \s \, \#} \, V\|) \,d\|{\G_{\#} \,{\mathbf  C}^{\prime}}\|(X)&&\\
&&\hspace{-3in} +\; \int_{{\mathbf R} \times B_{1}} {\rm dist}^{2} \, (X, {\rm spt}\, \|\G_{\#} \, {\mathbf C}^{\prime}\|) \, d\|\eta_{0, \th_{j}\s \, \#} \, V\|(X) \leq \nu_{j}\th_{j}^{2}\left(Q_{V}^{\star}(\s)\right)^{2}
\end{eqnarray*}
where the parameters $\th_{1}, \ldots, \th_{2q-3}$ and the constants $\nu_{1}, \ldots \nu_{2q-3}$ are analogous to the same quantities as in Lemma~\ref{multi-scale-final}, with $\nu_{1}$ depending only on $n$, $q$, $\a$, $\m\r_{0}$ and
for $j \in \{2, 3, \ldots, 2q-3\}$, $\nu_{j}$ depending only on $n$, $q$, $\a,$ $\th_{1}, \ldots, \th_{j-1}$,$\m\r_{0}$. By choosing $\th_{1}, \th_{2}, \ldots, \th_{2q-3}$ in that order, depending only on $n$, $q$, $\a$ and $\m\r_{0}$, to ensure that $\nu_{j}\th_{j}^{2} < \frac{1}{2}\th_{j}^{2\a}$ and $\th_{j}^{2\a^{\prime}} < \frac{1}{2}\th_{j}^{2\a}$ for each $j=1, 2, \ldots, 2q-3$, we deduce that under the hypotheses of the modified lemma, $$\left(Q_{V}^{\star}(\G_{\#} \, {\mathbf C}^{\prime}, \th_{j}\s)\right)^{2} \leq \th_{j}^{2\a} \left(Q_{V}^{\star}({\mathbf C}, \s)\right)^{2}$$ for some $j \in \{1, 2, \ldots, 2q-3\},$ ${\mathbf C}^{\prime} \in {\mathcal C}_{q}$ and an orthogonal rotation $\G$ of ${\mathbf R}^{n+1}.$ In view of the Remark preceding Theorem~\ref{compactness-N}, the iterative application of this  as in Lemma~\ref{finedistgap} gives the analog of Lemma~\ref{finedistgap}; arguing as in Corollary~\ref{F7} then 
establishes {\bf Case 2}, completing the proof that ${\mathcal B}_{q}^{\star}$ is a proper blow-up class.

\noindent
(x) In view of (i) and (\ref{manifolds-2-1}), the argument of Section~\ref{sheeting} carries over to yield Theorem~\ref{sheeting-N} for $q \geq 2$, subject to the induction hypotheses as in (vi).  Theorem~\ref{mindist-N} first in case $\Theta \, (\|{\mathbf C}_{0}\|, 0) = q + 1/2$ and then in case $\Theta \, ({\mathbf C}_{0}\|, 0) = q+1$ follows, again subject to the induction hypotheses as in (vi), from the argument, with obvious modifications, of  Section~\ref{MD}; note in particular that in view of the ``monotonicity inequality'' (\ref{monotonicity-N}) needed in the proof, and the need to use directly the first variation inequality (${\mathcal S^{\star}{\emph 1}}$) in establishing regularity of blow-ups as in Theorem~\ref{MD-c1alpha}, we must take 
$${\E}_{V}^{\star}({\mathbf C}, \s) \equiv \left(\int_{B_{1} ^{n+1}(0)} {\rm dist}^{2} \, (X, {\rm spt} \, \|{\mathbf C})\|) \, d\|\eta_{0, \s \, \#} \, V\|(X) + \s^{2\a^{\prime}}\right)^{1/2}$$ in place of the excess ${\E}$ used in Section~\ref{MD} (see Lemma~\ref{MD-excess-improvement}). Same modification applies to the excess ${\mathcal Q}$ used in Lemma~\ref{MD-multi-scale-final}.

\noindent
{\bf Step 3:} In view of {\bf Step 1}, {\bf Step 2} and the fact that Allard's integral varifold compactness theorem (\cite{AW}, Theorem 6.4) holds in Riemannian manifolds, Theorem~\ref{compactness-N} follows from the argument of Theorem~\ref{compactness} in Section~\ref{reg-compactness}.

\section{A sharp varifold maximum principle}\label{max}
\setcounter{equation}{0}

We conclude this paper by pointing out an immediate application of Theorem~\ref{compactness-N}, namely, 
the following optimal strong maximum principle for co-dimension 1 stationary integral varifolds:

\begin{theorem}\label{maxm}
Let $N$ be a smooth $(n+1)$-dimensional Riemannian manifold.
\begin{itemize}
\item[(a)] If $V_{1}$, $V_{2}$ are stationary integral $n$-varifolds on $N$ such that 
$${\mathcal H}^{n-1} \, ({\rm spt} \, \|V_{1}\| \cap {\rm spt} \, \|V_{2}\|) = 0,$$
then $\;\;\; {\rm spt} \, \|V_{1}\| \cap {\rm spt} \, \|V_{2}\|= \emptyset.$
\item[(b)] Let $\Omega_{1}$, $\Omega_{2}$ be open subsets of $N$ with $\Omega_{1} \subset \Omega_{2}$ and $M_{i} = \partial \, \Omega_{i}$, $i=1, 2.$ 
If for $i=1, 2$, $M_{i}$ is connected,  
$${\mathcal H}^{n-1} \, ({\rm sing} \, M_{i}) = 0 \, \;\;\; {\rm and}$$
$V_{i} \equiv |M_{i}|$ is stationary in $N,$ 
then either 
$${\rm spt} \, \|V_{1}\| = {\rm spt} \, \|V_{2}\| \;\;\; {\rm or} \;\;\; {\rm spt} \, \|V_{1}\| \cap {\rm spt} \, \|V_{2}\| = \emptyset.$$ 
Here ${\rm sing} \, M_{i} = M_{i} \setminus {\rm reg} \, M_{i},$ where ${\rm reg} \, M_{i}$ is the set of points $X \in M_{i}$ with the property that there exists a number 
$\s = \s(X) >0$ such that $M_{i} \cap B^{n+1}_{\s}(X)$ is a smooth, properly embedded hypersurface of $B_{\s}^{n+1}(X)$
with no boundary in $B_{\s}^{n+1}(X)$.
\end{itemize}
\end{theorem}

\noindent
{\bf Remark:} These results were established by T. Ilmanen (\cite{I}) under the stronger hypotheses that
$${\mathcal H}^{n-2} \, ({\rm spt} \, \|V_{1}\| \cap {\rm spt} \, \|V_{2}\|) = 0$$ 
\noindent
in part (a) and 
$${\mathcal H}^{n-2} \, ({\rm sing} \, M_{i}) = 0, \;\;\; i=1, 2,$$ 
\noindent
in part (b). Obvious examples show that for any $\g > 0$, neither of these hypotheses can be weakened to 
${\mathcal H}^{n-1+\g} \, (\cdot) = 0.$ 

\begin{proof}
The argument of \cite{I} carries over, with (2) of \cite{I} replaced by the hypothesis 
$${\mathcal H}^{n-1} \, ({\rm sing} \, M) = 0$$ 
\noindent
and  Theorems (8), (9) therein replaced by our Theorem~\ref{compactness-N}.
\end{proof}

\end{document}